\documentclass[12pt,a4paper]{book}
 \include{preamble}

\begin{document}

\begin{titlepage}

\thispagestyle{empty}

\newgeometry{left=7.4cm,bottom=2cm, top=1cm, right=1cm}

\tikz[remember picture,overlay] \node[opacity=1,inner sep=0pt] at (-21mm,-128mm){\includegraphics{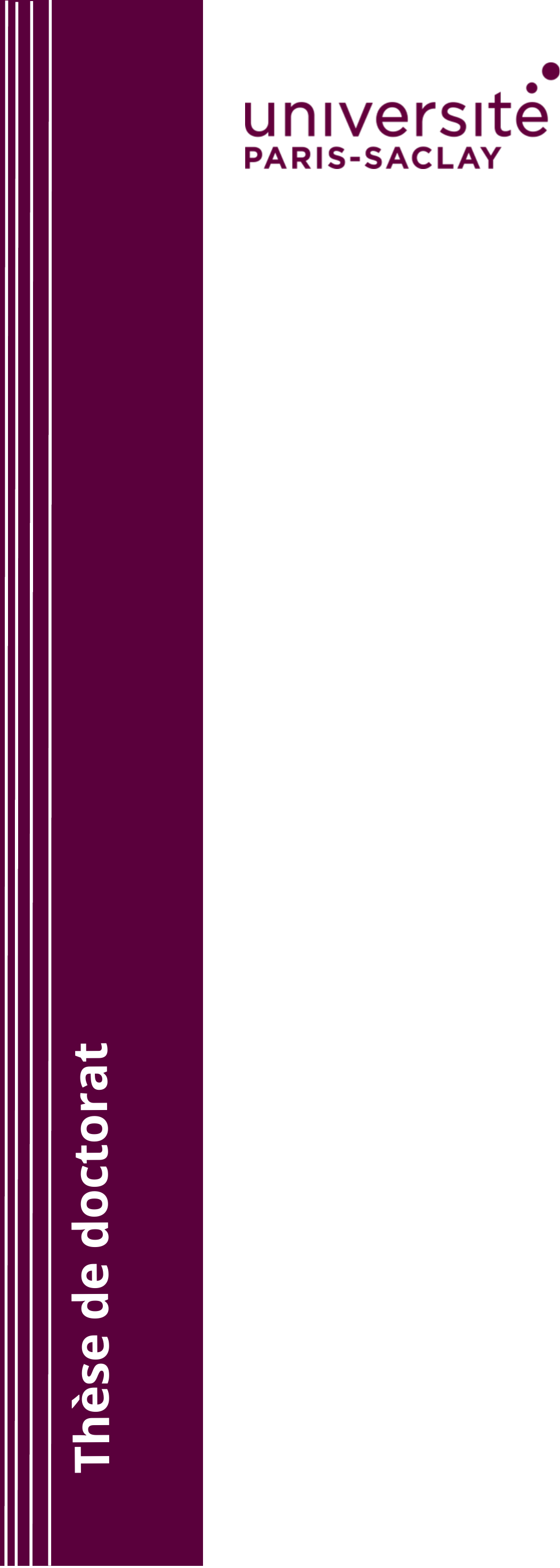}};

\fontfamily{fvs}\fontseries{m}\selectfont


\color{white}

\begin{picture}(0,0)

\put(-130,-720){\rotatebox{90}{NNT: 2021UPASM031}}
\end{picture}

\vspace{-21mm} 
\flushright \includegraphics[scale=0.015]{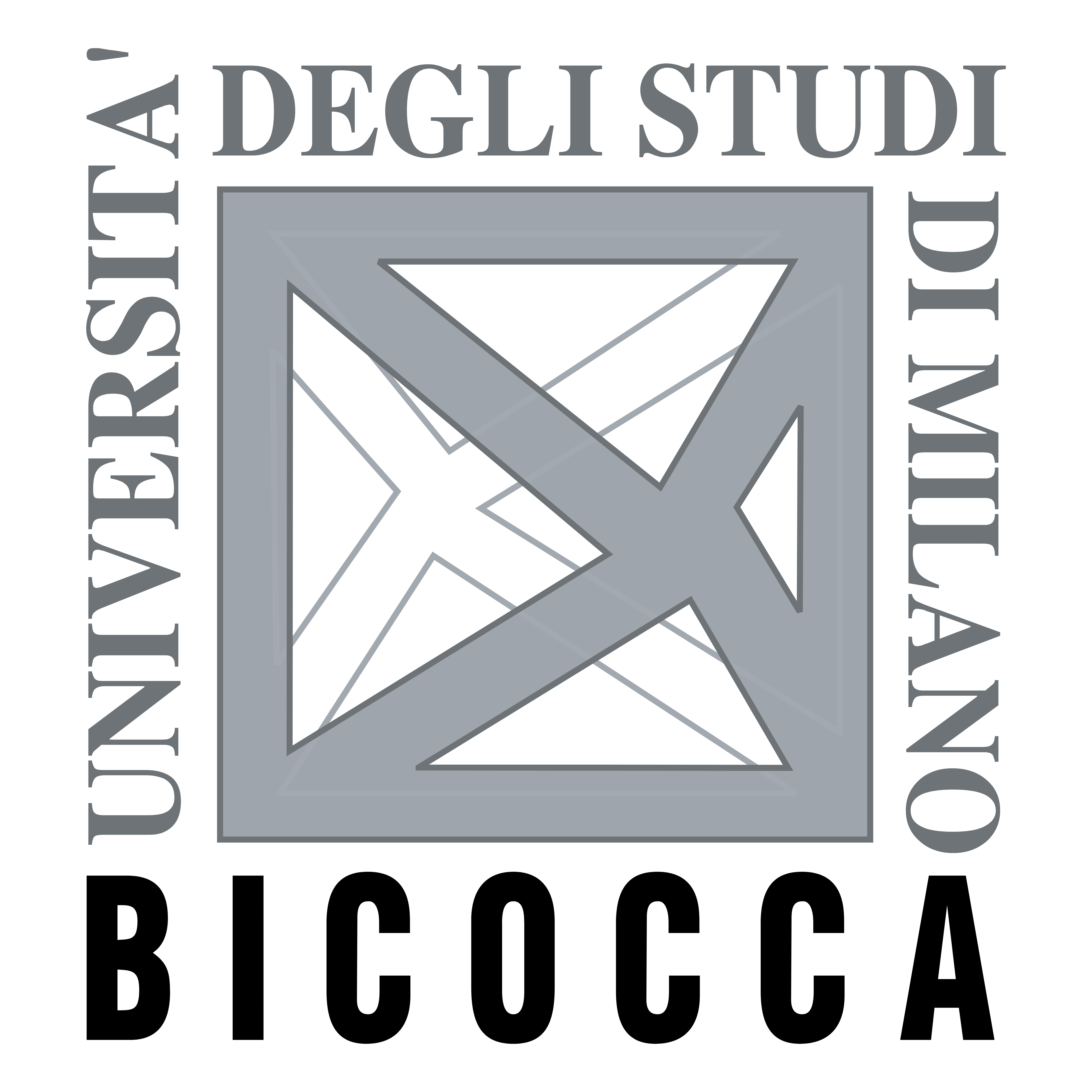}

\flushright
\vspace{10mm} 
\color{Prune}
\fontfamily{fvs}\fontseries{m}\fontsize{22}{26}\selectfont
  Weak regularization by degenerate Lévy noise 
  and its applications

\vspace{5mm} 
\color{black}
\fontfamily{fvs}\fontseries{m}\fontsize{18}{22}\selectfont
Régularisation faible par un bruit de Lévy dégénéré 
et applications

\normalsize
\vspace{1.5cm}

\color{black}
\textbf{Thèse de doctorat de l'Université Paris-Saclay}

\vspace{15mm}

École doctorale de Math\'ematique Hadamard (EDMH) n$^{\circ}$574\\
\small Spécialité de doctorat: Math\'ematiques appliqu\'ees\\
\footnotesize Unité de recherche: Laboratoire de math\'ematiques et mod\'elisation\\
d'Evry (UEVE), UMR 8071, CNRS-INRA \\
\footnotesize Référent: Universit\'e d'Evry-Val d'Essonne
\vspace{15mm}

\textbf{{Thèse présentée et soutenue à Évry, le 01 Octobre 2021, par}}\\
\bigskip
\Large {\color{Prune} \textbf{Lorenzo MARINO}}

\vspace{\fill} 

\flushleft \small \textbf{Composition du jury:}
\medskip

\scriptsize
\begin{tabular}{|p{8cm}l}
\arrayrulecolor{Prune}
\textbf{Pierre-Gilles Lemarié-Rieusset} &  Pr\'esident \\
Professeur, Universit\'e d'Evry Val d'Essonne   &   \\
\textbf{Arturo Kohatsu-Higa} &  Rapporteur {\&}  Examinateur \\
Professeur, Ritsumeikan University  &   \\
\textbf{Diego Pallara} &  Rapporteur {\&} Examinateur \\
Professeur, Universit\`a del Salento   &\\
\textbf{Franco Flandoli} &  Examinateur \\
Professeur, Scuola Normale Superiore di Pisa   &   \\
\end{tabular}

\bigskip

\flushleft \small \textbf{Direction de la th\`ese:}
\medskip

\scriptsize
\begin{tabular}{|p{8cm}l}
\arrayrulecolor{Prune}
\textbf{Stéphane Menozzi} &  Directeur \\
Professeur, Université d'Evry Val d'Essonne   &   \\
\textbf{Enrico Priola} &  Codirecteur \\
Professeur, Universit\`a degli studi di Pavia   &   \\
\end{tabular}

\end{titlepage}


\normalfont
\normalsize
\newgeometry{top=3.3cm, bottom=3.3 cm, left=2.7cm, right=2.7cm, headheight= 1cm,
heightrounded
}
\selectlanguage{english}

\tableofcontents

\selectlanguage{english}
\pagestyle{fancy}
\fancyhf{}
\fancyfoot[LE,RO]{\thepage}
\renewcommand{\headrulewidth}{0.5pt}
\renewcommand{\footrulewidth}{0.5pt}
\chapter{Introduction}
\fancyhead[LE]{Chapter $1$. Introduction}
\section{The model}
\fancyhead[RO]{Section \thesection. The model}
\label{Sec:INTRO:modello_considerato}
The present PhD thesis focuses on the study of the regularisation by degenerate Lévy noise phenomena for chains of differential equations that are possibly ill-posed. In particular, our main goal is to determine, under suitable conditions on the dynamics, which is the minimal Hölder regularity on the coefficients which ensures the well-posedness of the associated stochastic equation. More in details, given a ``large'' space $\R^N$ and a ``small'' one $\R^d$ ($d\le N$), we are interested in a stochastic differential equation (SDE) of the following form:
\begin{equation}\label{INTRO:SDE_Generale1}
   dX_t \, = \, F(t,X_t)dt +B\sigma(t,X_{t-})dZ_t, \quad t\ge 0
\end{equation}
where $\{Z_t\}_{t\ge 0}$ is a $d$-dimensional Lévy process on some probability space
$(\Omega,\mathcal{F},\mathbb{P})$ and the drift $F\colon [0,+\infty)\times\R^N\to \R^N$ and the diffusion matrix $\sigma\colon [0,+\infty)\times\R^N\to \R^d\otimes\R^d$ are Hölder continuous in space, uniformly in time. Above, the process $\{Z_t\}_{t\ge 0}$ is transmitted through the diffusion coefficient $\sigma$ on $\R^d$, thanks to a uniform ellipticity property, and is then immersed in the large space $\R^N$ through the matrix $B$ in $\R^N\otimes\R^d$ for which we assume, without loss of generality, that $\rank B=d$.

A detailed analysis of the well-posedness of degenerate stochastic chains of the form \eqref{INTRO:SDE_Generale1} is not only important for its intrinsic mathematical interest but especially for the different scientific contexts in which this type of dynamics appears as a model. To name a few, \emph{non-local degenerate kinetic diffusions} of the following form:
\begin{equation}\label{Intro:SDE_cinetica}
\begin{cases}
dX^1_t \, &= \,F_1(t,X^1_t,X^2_t)dt+ dZ_t; \\
dX^2_t \, &= \, F_2(t,X^1_t,X^2_t)dt,
\end{cases}
\end{equation}
corresponding to Equation \eqref{INTRO:SDE_Generale1} with only two macro-components ($N=2d$), $\sigma=1$, $B=(I_{d\times d},0_{d\times d })^t$, are often used in Hamiltonian mechanics in the analysis of turbulent regimes or in models that consider anomalous diffusion phenomena such as, for example, the study of heat variations between two different materials in contact (cf.\ \cite{Baeumer:Benson:Meerschaert01,Eckmann:Pillet:Rey-Bellet99}). Furthermore, SDEs like \eqref{Intro:SDE_cinetica} with $F_2(t,x_1,x_2)=x_1$ can be used to describe the position/velocity dynamics of a moving particle in dynamics where only the velocity component $X^1_t$ is randomly perturbed (for example, due to turbulent effects given by the air at high speed regimes) while the position variable $X^2_t$ feels the random noise only through its dependence on the first component. This kind of degenerate kinetic models also appears as the diffusive limit for linearized Boltzmann equations if the equilibrium function is assumed to be a Lévy-type distribution (\cite{Alexandre12, Mellet:Mischler:Mouhot11, Mellet16, Chen:Zhang18}). See also \cite{Cushman:Park:Kleinfelter:Moroni05} for an application to Richardson's law of turbulence. Another natural use of non-local degenerate kinetic diffusions appears in finance and, in particular, in the pricing of \emph{Asian options}, a particular type of ``exotic'' financial derivatives (cf.\ \cite{Barndorff-Nielsen:Shephard01,Brockwell01,Jeamblanc:Yor:Chesney09,Bally:Kohatsu-Higa10}). They can be described as options whose payoff, i.e. how much the stakeholder earns at a given expiry date, depends on the average of the values assumed by the underlying during the entire life of the contract and not only on the value of the underlying at its closure, as in the more common European options.

\begin{figure}
    \centering
    \label{fig:my_label}
\begin{tikzpicture}[scale=0.7]
\draw[gray, very thick] (2.7513,4.7619)--(3.4921,4.9735)--(3.1746,5.582)--(3.9683,5.582)--(3.6772,6.2963)--(4.3915,6.2169)--(4.4709,6.7196);
\draw[gray, very thick] (4.4709,6.746)--(5.1058,6.6402)--(5.0265,6.0317)--(5.6349,6.2434)--(5.582,5.6349)--(6.2434,5.7672);
\draw[gray, very thick] (6.2434,5.7407)--(6.8783,5.8201)--(6.6931,6.2963)--(7.2222,6.1376);
\draw[gray, very thick] (7.1164,6.7196)--(7.6984,6.6402)--(7.4603,7.0899)--(8.1481,7.037);
\draw[gray, very thick] (8.9683,7.7513)--(8.2804,7.7249)--(8.3862,7.3016);
\draw[black, very thick] (2.2487,8.254)--(2.7778,7.7249);
\draw[black, very thick] (2.7513,7.7513)--(2.2487,7.2487)--(2.7513,6.746)--(2.2487,6.2434)--(2.7513,5.7407)--(2.2487,5.2381)--(2.7513,4.7619)--(2.2487,4.2593)--(2.7513,3.7566);
\draw[black, very thick] (3.9947,8.254)--(4.4974,7.7513)--(3.9947,7.2487)--(4.4974,6.746)--(3.9947,6.2434)--(4.4974,5.7407)--(3.9947,5.2381)--(4.4974,4.7619)--(3.9947,4.2593)--(4.4974,3.7566);
\draw[black, very thick] (5.7407,8.254)--(6.2434,7.7513)--(5.7407,7.2487)--(6.2434,6.746)--(5.7407,6.2434)--(6.2434,5.7407)--(5.7407,5.2381)--(6.2434,4.7619)--(5.7407,4.2593)--(6.2698,3.7566);
\draw[black, very thick] (8.5185,8.254)--(8.9947,7.7513)--(8.4921,7.2487)--(8.9947,6.746)--(8.4921,6.2434)--(8.9947,5.7407)--(8.4921,5.2381)--(8.9947,4.7619)--(8.4921,4.2593)--(8.9947,3.7566)--(8.9947,3.7302);
\filldraw[black](2.7249,4.7619) circle (7pt);
\filldraw[black](4.4974,6.746) circle (7pt);
\filldraw[black](6.2434,5.7407) circle (7pt);
\filldraw[black](8.9683,7.7513) circle (7pt);

\draw[black, very thick] (1.746,8.0159)--(1.746,6.7196)--(1.4815,6.9974);
\draw[black, very thick] (1.5344,7.7778)--(1.746,7.9894)--(1.9841,7.7646);
\draw[black, very thick] (1.746,6.7328)--(1.9841,6.9841);
\draw[black, very thick] (3.5053,7.9762)--(3.5053,6.7593)--(3.4921,6.7857);
\draw[black, very thick] (3.2937,7.7381)--(3.5053,8.0026)--(3.7434,7.7646);
\draw[black, very thick] (3.254,6.9974)--(3.4921,6.7328)--(3.7302,6.9841);
\draw[black, very thick] (6.746,6.7593)--(6.746,8.0026);
\draw[black, very thick] (6.5079,6.9974)--(6.746,6.746)--(6.9974,6.9974);
\draw[black, very thick] (6.4947,7.7513)--(6.746,8.0026)--(6.9974,7.7513);
\draw[black, very thick] (-.9,5.7513) circle(36pt);
\node at (-.9,5.7407) { \bf Noise};
\node at (1,5.7407) {{\Large $\Rightarrow$}};
\end{tikzpicture}
\caption{Interconnected elasto-plastic oscillators from \cite{Delarue:Menozzi10}.}
\end{figure}
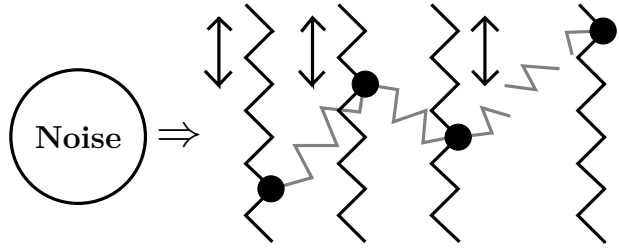

All the examples mentioned above focus only on the kinetic case, i.e.\ when Equation \eqref{INTRO:SDE_Generale} is composed by only two components ($N=2d$). In a more general framework ($N=nd$), models like:
\begin{equation}\label{INTRO:SDE_Generale}
\begin{cases}
   dX^1_t \, = \, F_1(t,X^1_t,\dots,X^n_t)dt +dZ_t \\
   dX^2_t \, = \, F_2(t,X^1_t,\dots,X^n_t)dt;\\
   dX^3_t \, = \, F_3(t,X^2_t,\dots,X^n_t)dt;\\
  \vdots \\
   dX^n_t \, = \, F_n(t,X^{n-1}_t,X^n_t)dt,
\end{cases}
\end{equation}
appear, for example, in seismology, in the study of shock waves propagation through different structures. Furthermore, SDEs like \eqref{INTRO:SDE_Generale} are often used to represent interconnected elasto-plastic oscillators, i.e.\ systems of connected springs, when a random perturbation is applied only to the first of them (cf.\ Figure \ref{fig:my_label}). In this regard, see for example \cite{Bensoussan:Turi08, Bensoussan:Mertz:Pironneau:Turi09} in a diffusive context. More generally, models that consider Lévy-type noises appear more versatile and realistic than their Brownian counter-parts since they allow for the presence of random jumps in the dynamics.

To highlight the more general context in which our specific problem fits, we begin the present thesis with a brief introduction to the regularisation by noise theory.

\section{Regularisation by noise}
\fancyhead[RO]{Section \thesection. Regularisation by noise}
\label{Sec:INTRO:Regolarizzazione_rumore}

Given an initial point $x$ in $\R^N$, the classical Cauchy-Lipschitz theory ensures that an ordinary differential equation of the following type:
\begin{equation}\label{INTRO:eq:ODE}
\begin{cases}
dX_t \, = \, F(t,X_t)dt, \quad t> 0; \\
X_0\, = \, x,
\end{cases}    
\end{equation}
admits a unique global solution if the drift $F$ is smooth enough, i.e.\ Lipschitz continuous and with at most linear growth.
This result was then extended by several authors under increasingly weak conditions on the regularity of $F$. In this regard, we recall the famous DiPerna-Lions theory in \cite{Diperna:Lions89} where weakly Lipschitz continuous drifts are considered and the subsequent generalisation to drifts of bounded variation (cf.\ \cite{Ambrosio04,Crippa:DeLellis08}). However, all the cited above works share some sort of boundedness condition on the derivatives of $F$ which, at least in the case of time independent drifts, can be stated as $\text{div} \, F \in L^\infty (\R^N)$.

Weakening the above divergence condition immediately allows the emergence of models that confute the well-posedness of Equation \eqref{INTRO:eq:ODE}. A classic counter-example is indeed the \emph{Peano} model obtained from \eqref{INTRO:eq:ODE} imposing $d=1$ and $F(s,x)=\text{sgn}( x)|x|^\beta$ for some $\beta$ in $(0,1)$. In fact, it is not difficult to show that this equation, when starting at $x=0$, admits an infinite number of solutions of the form:
\begin{equation}
\label{INTRO:eq:solution_Peano}
X_s\, =\, (1-\beta)^{\frac{1}{1-\beta}}\mathds{1}_{[T_0,+\infty)}(s)\left(s-T_0\right)^{\frac{1}{1-\beta}},\quad s\ge0; 
\end{equation}
as $T_0>0$ varies. Intuitively, the presence of the singularity at the origin, where the drift $F$ is not Lipschitz regular, allows the solutions to stay at that point for any length of time. This kind of phenomenon is often called a \emph{bifurcation} of flows. We focused here only on a counter-example with Hölder regular drift and only on the problem of the non-uniqueness of solutions because this is the framework we will be interested in. However, we mention that in the case of discontinuous drifts, different and more variegate phenomena can actually occur (cf.\ \cite{book:Flandoli11, Diperna:Lions89}), among which, for example, the \emph{coalescence} of the flows, i.e.\ the collision of solutions starting from different initial points, or even the non-existence of solutions.

\paragraph{Regularisation by non-degenerate Brownian noise.}

The situation presented in the previous paragraph drastically changes if a random noise is added to the dynamics, i.e. if one consider instead of the deterministic Equation \eqref{INTRO:eq:ODE} its stochastic counterpart given by
\begin{equation}\label{INTRO:eq:SDE:Brownian}
\begin{cases}
dX_t \, = \, F(t,X_t)dt +\sigma(t,X_t)dW_t,\\
X_0\, = \, x,
\end{cases}    
\end{equation}
where $\{W_s\}_{s\ge 0}$ is a Brownian motion on $\R^N$. Indeed, the presence of a sufficiently ``strong'' random noise may actually restore the well-posedness of the dynamics. For stochastic Peano-type models with evanescent noise (cf.\ Equation \eqref{INTRO:eq:SDE:Brownian} with $\sigma=\varepsilon$ small and $F(t,x)=\text{sgn}(x )|x|^\beta$), this phenomenon was described by Delarue and Flandoli in \cite{Delarue:Flandoli14} (see also the previous work \cite{Bafico:Baldi81}): in small time, the noise fluctuations dominate the dynamics so that the solution can escape the singularity at zero, while in a long time range, the drift $F$ prevails forcing the solution to fluctuate around one of the two extremal solutions (cf.\ Equation \eqref{INTRO:eq:solution_Peano} with $T_0=0$). It is therefore important to recall how, in a short time, there is a strong competition between the roughness of the drift $F$ and the mean fluctuations of the Brownian noise $W_t$.

We now recall that in the stochastic framework, some particular attention must be paid to what exactly is meant by uniqueness of solutions to Equation \eqref{INTRO:eq:SDE:Brownian}. A \emph{weak solution} of SDE \eqref{INTRO:eq:SDE:Brownian} is a pair of processes $\{(X_t,W_t)\}_{t\ge0}$ where $\{W_t \}_{t\ge 0}$ is a Brownian motion on a stochastic basis $(\Omega,\mathcal{F},\{\mathcal{F}_t\}_{t\ge 0},\mathbb{P})$ while $\{X_t\} _{t\ge 0}$ satisfies Equation \eqref{INTRO:eq:SDE:Brownian} with respect to the noise $W_t$. Intuitively, the Brownian motion $W_t$ in this case has to be constructed as well and it is itself part of the notion of solution. We will then say that \emph{weak uniqueness}, or uniqueness in law, holds for SDE \eqref{INTRO:eq:SDE:Brownian} when for each pair of weak solutions $\{(X^1_t,W^1_t)\} _{t\ge0}$ and $\{(X^2_t,W^2_t)\}_{t\ge0}$, the marginal laws of the processes $\{X^1_t\}_{t\ge0}$ and $\{X^2_t\}_{t\ge0}$ coincide. We emphasise in particular that this notion of uniqueness does not even require the two solutions to be defined on the same stochastic basis. On the other hand, a solution $\{X_t\}_{t\ge 0}$ of Equation \eqref{INTRO:eq:SDE:Brownian} is said to be \emph{strong} if, given a Brownian noise $\{ W_t\}_{t\ge0}$, $X_t$ solves the given equation with respect to $W_t$ and it is adapted to the natural filtration of $\{W_t\}_{t\ge 0}$. Intuitively, a strong solution $\{X_t\}_{t\ge0}$ can be seen, at any fixed time $t$, as a functional of the given noise $W_t$, seen as an input of the system. A usual method to establish the strong uniqueness of stochastic dynamics such as \eqref{INTRO:eq:SDE:Brownian}, i.e. when a strong solution is unique for some given $W_t$ noise, is to exploit the Yamada-Watanabe theorem (cf.\ \cite{Yamada:Watanabe71}) which ensures the well-posedness in the strong sense to the dynamics, starting from the existence of a weak solution and the path-wise uniqueness, i.e.\ when two solution processes are indistinguishable.

In the non-degenerate Brownian context, it has been shown that SDEs like \eqref{INTRO:eq:SDE:Brownian} are well-posed in a weak sense as soon as the drift $F$ is bounded measurable and the diffusion matrix $ \sigma$ is uniformly elliptic and continuous in space. Such a result was established by Stroock and Varadan in \cite{book:Stroock:Varadhan79}, exploiting the now well-known \emph{martingale problem} approach. In recent years, great attention has also been given to additive noise dynamics (i.e.\ $\sigma=1$) with distributional drifts in space, especially for their connection with physical models in materials theory (cf.\ \cite{Alberts:Khanin:Quastel14,Brox86}). A first result in this regard can be found in \cite{Flandoli:Issoglio:Russo17} where they considered a time-inhomogeneous drift $F$ belonging to a Hölder space with negative index but strictly larger than $-1/2$. For a more precise presentation of such negative Hölder spaces, the reader is referred to \cite{book:Friz:Hairer}, Section $13$. Subsequently, Delarue and Diel in \cite{Delarue:Diel16} proved that the same result still holds when one consider regularity indexes that are only greater than $-2/3$, even if only for a one-dimensional dynamics. Such a work was then extended to the multidimensional case in \cite{Cannizzaro:Chouk18}. We also mention the work of Bass and Chen \cite{Bass:Chen03} where an additive noise dynamics with time-independent drift and belonging to the Kato class is considered.

As far as the path-wise uniqueness for solutions to Equation \eqref{INTRO:eq:SDE:Brownian} is concerned, the first result was obtained by Zvonkin in \cite{Zvonkin74} for a one-dimensional dynamics with bounded measurable drift $F$ and a Hölder continuous diffusion matrix $\sigma$ with regularity index strictly larger than $1/2$. This result was then extended by Veretennikov \cite{Veretennikov80} to the multi-dimensional case for a Lipschitz continuous diffusion matrix. As previously mentioned, these works focus on the path-wise uniqueness but they could then be used to show as well the existence of strong solutions for the dynamics, thanks to the Yamada-Watanabe theorem. In this regard, we mention the more direct approach for constructing strong solutions under the same assumptions on the drift above, obtained in \cite{Meyer:Proske10,Menoukeu:Meyer:Nilssen:Proske:Zhang13} through Malliavin calculus. Subsequently, Krylov and Röckner in \cite{Krylov:Rockner05} and Zhang in \cite{Zhang10} showed the strong well-posedness for stochastic dynamics as \eqref{INTRO:eq:SDE:Brownian} when the drift is only integrable, i.e.\ when $F$ belongs to $L^p(0,T;L^q(\R^N))$ under the Prodi-Serrin condition: $N/q+2/p<1$, $p,q\ge 2$, respectively for an additive noise and a multiplicative one with a diffusion matrix $\sigma$ in Sobolev spaces. See also the recent works \cite{Krylov21} and \cite{Rockner:Zhao21,Nam20} in which the additive critical case (i.e. $\frac{N}{q}+\frac{2}{p }=1$) respectively for homogeneous and inhomogeneous in time drifts is studied. However, we remark that in \cite{Nam20} they only consider a Lorentz-type integrability, that is known to be stronger than the usual Lebesgue one. See Definition $2.1$ in the cited article, for more details. Finally, we mention the work \cite{Fedrizzi:Flandoli11} by Fedrizzi and Flandoli, where under the same conditions of Krylov and Röckner, the continuous dependence of the solutions $X_t$ on the initial condition $x$ is shown, and the subsequent works \cite{Zhang15,Xie:Zhang16} by Zhang \emph{et al.}, where instead the differentiability in a weak sense (i.e.\ in Sobolev spaces) still with respect to the initial condition $x$ is analysed. The works mentioned above intuitively illustrate what has been called, following Flandoli's terminology in \cite{book:Flandoli11}, a \emph{regularisation by noise} phenomenon: when a deterministic differential equation is ill posed (because the existence or the uniqueness of solutions fail) while the associated stochastic dynamics is well posed in a weak or strong sense. We suggest to the interested reader to see the monograph \cite{book:Flandoli11} where a more general analysis of the topic is presented. 

We have focused for the moment on \emph{non-degenerate} stochastic dynamics, i.e. when the noise has the same dimension as the underlying system it acts on (i.e.\ when $N=d$ in \eqref{INTRO:SDE_Generale1}). However, we remark that this condition is not actually always satisfied in many practical cases. One can thinks, for example, in Hamiltonian mechanics, of the Langevin equations with a perturbation on the velocity component, or in finance, in the analysis of Asian options. It has also been shown (\cite{ delCastillo10, Woyczynski01}) that in many practical applications, random fluctuations in real complex systems are indeed often of a non-Gaussian nature. The starting point and main motivation for this work was then to analyse the regularisation by noise phenomena introduced above for chains of ordinary differential equations when the random perturbation is no longer a Brownian motion but a more general Lévy process with suitable properties, which possibly acts only on some components of the system (i.e.\ when $d<N$ in \eqref{INTRO:SDE_Generale1}). However, due to some time constraints, we decided to focus only on the well-posedness of such dynamics in a weak sense. We now briefly present an overview of the results already known in this area.

\paragraph{Weak regularisation by non-degenerate stable noise.}
The problem of the weak well-posedness for stochastic dynamics of the following form:
\begin{equation}
\label{INTRO:eq:non-deg_SDE}
X_t \, = \, x +\int_0^tF(X_s) ds +Z_t, \quad t\ge 0,
\end{equation}
where $\{Z_t\}_{t\ge0}$ is a symmetric, $\alpha$-stable process on $\R^N$, has been extensively studied during the years. One of the first contributions in the one-dimensional case has been obtained by Tanaka \emph{et al.} in \cite{Tanaka:Tsuchiya:Watanabe74} where the uniqueness in law is proved for Equation \eqref{INTRO:eq:non-deg_SDE} when the drift $F$ is bounded and continuous and the Lévy symbol $\Phi$ associated to $\{Z_t\}_{t\ge 0}$ satisfies some natural conditions at infinity: $\Re \Phi (\xi)^{-1}\, = \, O(1/\vert \xi \vert)$ when $\vert \xi \vert \to \infty$. An extension to the multi-dimensional case was then obtained in \cite{Komatsu83} assuming that the drift $F$ is bounded and continuous and the law of $\{Z_t\}_{t\ge 0}$ has a density with respect to the Lebesgue measure on $\R^N$. Stochastic dynamics like \eqref{INTRO:eq:non-deg_SDE} for drifts $F$ in suitable $L^p$ spaces have also been considered in \cite{Jin18}. 
To the best of our knowledge, the first work dealing with distributional drifts $F$ in an $\alpha$-stable context has been \cite{Athreya:Butkovsky:Mytnik18} in the one-dimensional case, where the authors consider a time-homogeneous drift $F$ belonging to a (negative) Hölder space of index strictly larger than $(1-\alpha)/2$. We also mention the almost simultaneous works \cite{Ling:Zhao19} and \cite{Chaudru:Menozzi20} where drifts $F$ on general Besov spaces are analysed, under suitable conditions on the parameters, which are respectively homogeneous and inhomogeneous in time. The results shown above rely on a Young-type formulation to make sense of the dynamics. Beyond the Young regime, we mention instead \cite{kremp:Perkowski20} where SPDE techniques such as the paracontrolled products (cf.\ \cite{Gubinelli:Imkeller:Perkowski}) are exploited to obtain the weak well-posedness for dynamics driven by inhomogeneous in time drifts belonging to negative Hölder spaces with regularity index strictly larger than $(2-2\alpha)/3$.

Finally, we point out that the above papers all focus on the \emph{sub-critical} $\alpha$-stable case, i.e. when $\alpha>1$. In fact, we remark that if $\alpha\le 1$, stochastic dynamics like \eqref{INTRO:eq:non-deg_SDE} are much more difficult to deal with since the noise no longer dominates the drift in small time. In this regard, we cite the works \cite{Zhao19} and \cite{Chaudru:Menozzi:Priola20} where the authors consider drifts that are $(1-\alpha)$-Hölder continuous and continuous, respectively.

\paragraph{Weak regularisation by degenerate noise.}

In all the previously mentioned results, the noise played a fundamental role, allowing to regularise the coefficients at each component of the system. It is then clear that in a \emph{degenerate} noise context, i.e.\ when the random perturbation directly acts only on some parts of the dynamics, the problem immediately appears much more delicate. Let us start by considering the classical example of a chain of $n$ different ordinary deterministic equations on $\R^d$ where only the first one is perturbed by a Brownian diffusion:
\begin{equation}\label{INTRO:eq4}
\begin{cases}
  dX^1_t \, =\,F_1(t,X^1_t,\dots,X^n_t)dt+\sigma(t,X^1_t,\dots,X^n_t)dW_t;\\ 
  dX^2_t \, =\,F_2(t,X^1_t,\dots,X^n_t)dt;\\
  \vdots             \\
  dX^n_t \, =\,F_n(t,X^{n-1}_t,X^n_t)dt.
\end{cases}
\end{equation}
In order to obtain a regularisation by noise phenomenon even in this framework, it is then necessary that the noise $W_t$ acting only on the first component propagates through the system, thus reaching all the components of the model. Classical assumptions that ensure this kind of behaviour are the uniform ellipticity of the diffusion matrix $\sigma$ (cf.\ [\textbf{UE}] in Section \ref{Sec:INTRO:Buona_posizione_Levy}) which guarantees the preservation of the noise on $\R^d$, and the so-called Hörmander condition (cf.\ \cite{Hormander67}) for the hypoellipticity of the system. Under a \emph{strong} Hörmander condition, i.e. when the diffusive vector fields and their commutators generate the space, the main works in the Brownian context have been obtained by Kusuoka and Stroock \cite{Kusuoka:Stroock84, Kusuoka:Stroock85, Kusuoka:Stroock87}, exploiting Malliavin's calculus techniques.

In the literature concerning degenerate noise models like \eqref{INTRO:eq4}, several authors have instead supposed that each component of the drift $F$ is differentiable with respect to its first component and that the resulting gradient is non-singular, i.e. they assumed that $ D_{x_{i-1}}F_i$ (with $i=2,\dots ,n$) has maximum rank. This non-degeneracy assumption on the subdiagonal of the Jacobian matrix of $F$ is the reason why this kind of condition has often been called of \emph{weak Hörmander-type}. Roughly speaking, up to a mollification of the coefficients if necessary, the drift $F$ is actually necessary to generate the space $\R^N$ through the Lie commutators. One of the first works to consider this type of condition was \cite{Menozzi11} where the author showed that SDE \eqref{INTRO:eq4} is well-posed in a weak sense when the drift $F$ is Lipschitz continuous and the diffusion matrix $\sigma$ is Hölder continuous. This result was then extended in \cite{Menozzi18} to diffusion matrices that are only continuous in space. For the characterisation of the weak well-posedness in a kinetic framework, corresponding to Equation \eqref{INTRO:eq4} for $n=2$, we also mention Zhang's work \cite{Zhang18} where a semi-linear drift $F$ is considered, such that $F_2(x_1,x_2)=Ax_1$ and under conditions of local integrability for $F_1$ and the continuity of the diffusion coefficient $\sigma$. Still in the case of two oscillators under a weak Hörmander-type condition, Chaudru de Raynal showed in \cite{Chaudru18} the weak well-posedness of the dynamics as soon as the drift $F$ is Hölder continuous with respect to the degenerate variable $x_2$ with a regularity index strictly greater than $1/3$. In such a work, it is also shown through Peano-type counter-examples that the  Hölder regularity threshold of $1/3$ for the drift $F$ is actually ``almost'' optimal. This result was then extended in \cite{Chaudru:Menozzi17} to the more general case of $n$ oscillators (cf.\ Equation \eqref{INTRO:eq4}). Intuitively, this minimal threshold for the Hölder regularity on the drift can be explained as the price to pay to balance the noise degeneration. Referring to the previous discussion about Peano-type models in \cite{Delarue:Flandoli14} (cf.\ Equation \eqref{INTRO:eq:solution_Peano}), this threshold is indeed related to the fact that the noise fluctuations are not strong enough to push the solution away from the singularity at the origin if the drift is too rough.

To the best of our knowledge, there are not so many papers that deal instead with the weak well-posedness of stochastic dynamics with degenerate and jump noises, i.e. when in SDE \eqref{INTRO:eq4} the Brownian motion $\{W_t\}_ {t\ge0}$ is substituted with an  $\alpha$-stable process $\{Z_t\}_{t\ge0}$. Within a regularisation by noise perspective (i.e. when the associated deterministic dynamics is ill-posed), we can only mention \cite{Huang:Menozzi16} where the authors showed the weak well-posedness for a linearised version of the dynamics in \eqref{INTRO:eq4} with $F(t,x)=Ax$ and a Hölder continuous diffusion coefficient $\sigma$, under some dimensional constraints on $d$, $N$.

\section{Uniqueness in law for degenerate chains}
\fancyhead[RO]{Section \thesection. Uniqueness in law for degenerate chains}
\label{Sec:INTRO:Unicita_in_legge}

Starting from the work \cite{book:Stroock:Varadhan79} in the non-degenerate diffusive case, the connection between the solutions of the \emph{martingale problem} and the weak solutions (in a probabilistic sense) of SDE \eqref{INTRO:SDE_Generale1} is well-known. Roughly speaking, this method characterises the process $\{X_t\}_{t\ge0}$ through its infinitesimal generator $L_t$. To be able to introduce it properly in our jump context, let $\mathcal{D}([0,+\infty);\R^N)$ be the space of all the paths from $[0,+\infty)$ to $\R^N$ that are càdlàg, i.e.\ right continuous paths with finite left limits. Skorokhod showed in \cite{Skorohod56} that it is possible to equip such a set with a natural distance so that it becomes a separable metric space. We can then think of $\mathcal{D}([0,+\infty);\R^N)$ as a Borel measurable space and consider probability measures on it.

We now introduce the canonical, or evaluation, process $\{y_t\}_{t\ge0}$ associated with the space $\mathcal{D}([0,+\infty);\R^N)$ and given by
\[y_t(\omega) \, = \, \omega(t), \quad \omega \in \mathcal{D}([0,+\infty);\R^N).\]
For more details, see for example \cite{book:Ethier:Kurtz86} or \cite{book:Bass11}. Fixed a point $x$ in $\R^N$, a probability measure $\mathbb{P}$ on $\mathcal{D}([0,+\infty);\R^N)$ is a solution of the martingale problem (for the infinitesimal generator $L_t$ and with the initial point $x$) if:
\begin{itemize}
     \item $\mathbb{P}\left(y_0 =x\right)\, = \, 1$;
     \item for each function $\phi$ in the domain $\dom\bigl(\partial_t+L_t\bigr)$, the process
   \begin{equation}
\label{INTRO:def_martingala}\Bigl{\{}\phi(t,y_t)-\phi(0,x)- \int_{0}^{t}\bigl(\partial_s+L_s\bigr)\phi (s,y_s) \, ds\Bigr{\}}_{t \ge 0}
\end{equation}
   is a $\mathbb{P}$-martingale with respect to the natural filtration of the canonical process $\{y_t\}_{t\ge0}$.
   \end{itemize}
  
A third possible characterisation of the process $\{X_t\}_{t\ge0}$ can be obtained from its marginal law $\mu_t$:
\begin{equation}
\label{INTRO:Fokker_Plank}
\begin{cases}
\partial_t \mu_t \, = \, L^\ast_t \mu_t, \quad t>0; \\
\mu_0\, = \, x,
\end{cases}
\end{equation}
where $L^\ast_t$ is the (formal) adjoint operator of $L_t$. The above dynamics is usually called the \emph{forward Fokker-Plank equation}. We recall that a continuous family of measures $\{\mu_t\}_{t\ge0}$ is a solution to Equation \eqref{INTRO:Fokker_Plank} if for each test function $\phi$ in $C^\infty_c( \R^N)$,
\[\partial_t \int_{\R^N}\phi(y)\, d\mu_t(y) \,= \, \int_{\R^N}L_t\phi(y)\, d\mu_{ t}(y), \]
in a distributional in time sense while the initial condition requires that $\mu_t$ converge (in a suitable sense) to the Dirac mass $\delta_{x}$ centred in $x$. For a more exhaustive introduction on this topic, we refer the reader to the monograph \cite{book:Bogachev:krylov:Rockner:Shaposhikov15}. An extensive study of these phenomena over the years (cf.\ \cite{book:Stroock:Varadhan79,book:Ethier:Kurtz86, Kurtz98, Kurtz11}) has now  concluded that the three ways described above (stochastic differential equation, martingale problem and forward Fokker-Plank equation) are, under minimal conditions on the coefficients, indeed equivalent in specifying a jump diffusion, in the sense that existence and/or uniqueness in one case implies existence and/or uniqueness also in the other ones.

We now focus on the uniqueness of solutions to the martingale problem associated with the operator $L_t$ which, as we will see later, is the main motivation for the analysis carried out in the first two chapters of this thesis. Let $\mathbb{P}_1$, $\mathbb{P}_2$ be two measures on the Skorokhod space $\mathcal{D}([0,+\infty);\R^N)$ which are solutions of the martingale problem with the same starting point $x$ in $\R^N$. From the definition of the martingale problem in \eqref{INTRO:def_martingala}, it makes sense to introduce the Cauchy problem associated with $L_t$ with zero terminal condition. In particular, given a final time $T>0$ and a source $f\colon [0,T]\times \R^N \to \R$ in a sufficiently ``rich'' class $\mathcal{F}$ of functions, we consider the following partial differential equation with Cauchy terminal condition:
\begin{equation}
\label{INTRO:eq:sourceless_IPDE}
\begin{cases}
  \partial_t u(t,x) +L_t u(t,x) \, = \, f(t,x), & \mbox{on } [0,T)\times \R^N; \\
  u(T,x) \, = \, 0, & \mbox{on } \R^N.
\end{cases}    
\end{equation}
Let us assume for the moment that a solution $u\colon [0,T]\times \R^N\to \R$ of the above Cauchy problem exists. Furthermore, if $u$ is regular enough, it immediately follows from \eqref{INTRO:def_martingala} that
\begin{equation}
\label{INTRO:eq:Analytical}
\bigl{\{}u(t,y_t)-\int_0^t f(s,y_s) \, ds-u(0,x)\bigr{\}}_{t \in [0,T]}
\end{equation}
is a $\mathbb{P}_i$-martingale for each $i\in \{1,2\}$, where, we recall, we have indicated with $\{y_t\}_{t\ge 0}$ the canonical process on $\mathcal{D}([0,T];\R^N)$. Now, taking the expected value at the final time $T$ in the above equation, we can exploit the martingale property in \eqref{INTRO:eq:Analytical} and the system \eqref{INTRO:eq:sourceless_IPDE} satisfied by $u$ (in particular, $u(T,\cdot)=0$) to link the two solutions of the martingale problem as follows:
\[\mathbb{E}^{\mathbb{P}_1}\left[\int_0^T f(s,y_s) \, ds\right] \, = \,u(0,x) \, = \, \mathbb{E}^{\mathbb{P}_2}\left[\int_0^T f(s,y_s) \, ds\right],\]
where $\mathbb{E}^{\mathbb{P}_i}$ represents the expected value with respect to the probability measure $\mathbb{P}_i$. If the class of functions $\mathcal{F}$ we consider is rich enough, we can then conclude that the marginal law of the canonical process $\{y_t\}_{t\ge0}$ is the same under the two measures considered, for every fixed time $t$. By exploiting regular conditional probability techniques, it is then possible to show that the process $\{y_t\}_{t\ge0}$ has the same finite-dimensional distributions with respect to the two probability measures, i.e.\ $\mathbb{P}_1= \mathbb{P}_2$ (cf.\ Theorem $4.4.2$ in \cite{book:Ethier:Kurtz86}) and therefore, that the solution of the martingale problem for $L_t$ is unique.

The crucial difficulty in the above reasoning is the regularity assumption on the solution $u$ to Cauchy problem \eqref{INTRO:eq:sourceless_IPDE} that is indeed necessary to conclude in Equation \eqref{INTRO:eq:Analytical}. In fact, even when considering a smooth source $f$ in $C^\infty_c([0,T]\times \R^N)$, it would not be necessarily true that the solution $u$ is also smooth, because of the weak regularity (only Hölder continuous in space) of the coefficients $F$ and $\sigma$. When the noise is additive ($\sigma=1$) and the coefficients are regular enough, the required regularity of the solution $u$ can be obtained, for example, by exploiting stochastic flow techniques (cf.\ \cite{book:kunita19}). Furthermore, classical methods allows in this case to demonstrate as well the strong well-posedness of the considered stochastic dynamics. To apply the above reasoning in our multiplicative framework and under minimal regularity assumptions, it will instead be necessary to firstly approximate the coefficients appearing in the operator $L_t$ with a sequence of smooth functions, for example, through a mollification argument. 
It will then be possible to apply the method described above for such sufficiently regular coefficients and conclude that the two solutions to the ``mollified'' martingale problem are indeed equal. Finally, in order to go back to the original dynamics, we will need to establish a ``suitable'' analytic theory associated with the operator $L_t$ in order to control the convergence of the solutions with respect to the approximation parameter. An extensive literature on the subject (cf.\ \cite{Priola09,Chaudru17,Chaudru:Honore:Menozzi18_Strong, Flandoli:Gubinelli:Priola10}), has shown that a first step along this direction is to establish a particular type of estimates, called \emph {Schauder} ones, which allow to control the approximated solutions to Cauchy Problem \eqref{INTRO:eq:sourceless_IPDE} in terms of the mollified coefficients, on a suitable functional space. Furthermore, it is also possible to prove, by means of compactness arguments, that in fact such Schauder estimates also hold for the limit of the mollified solutions.

\subsection{Schauder estimates for degenerate systems}

As explained at the end of the previous section, we are now interested in a detailed analysis of a non-local equation of the Kolmogorov type which can be written as follows:
\begin{equation}
\label{INTRO:eq:IPDE_degenere}
\begin{cases}
  \partial_t u(t,x) +L_t u(t,x) \, = \, f(t,x), & \mbox{on } [0,T)\times \R^N; \\
  u(T,x) \, = \, u_T(x), & \mbox{on } \R^N,
\end{cases}    
\end{equation}
where the source $f\colon [0,T)\times \R^N\to \R$ and the final condition $u_T\colon \R^N\to\R$ are sufficiently smooth functions. We are going to study such equation in a finite time range, i.e.\ for any $t$ in $[0,T]$ for some fixed final time $T>0$. Above, $L_t$ represents the time dependent integro-partial differential operator which can be seen as the infinitesimal generator of the solution $\{X_t\}_{t\ge 0}$ to Equation \eqref{INTRO:SDE_Generale1}, i.e.\ such that
\begin{equation}
\label{INTRO:eq:def_L_t}
     L_t \, = \, \langle F(t,x),D_x\rangle +\mathcal{L}_t, \quad \text{ on } [0,T)\times \R^N,
\end{equation}
where $\mathcal{L}_t$ is the infinitesimal generator associated with the process $\{B\sigma Z_t\}_{t\ge0}$. We also remark that, similarly to the stochastic dynamics from which it is derived, Equation \eqref{INTRO:eq:IPDE_degenere} is degenerate in the sense that the principal component $\mathcal{L}_t$ of the operator $L_t$ only considers some directions of the space $\R^N$, the ones associated with the matrix $B$ in \eqref {INTRO:SDE_Generale1}.

In particular, the first two chapters of the current manuscript will focus on establishing the Schauder estimates for solutions to the degenerate system in \eqref{INTRO:eq:IPDE_degenere}. An important feature of such estimates is to highlight how much a solution $u$ to Cauchy problem \eqref{INTRO:eq:IPDE_degenere} gains, in terms of regularity, with respect to the source $f$, a notion usually called \emph{parabolic bootstrap} (or regularising effect, as referred to later) associated with the operator $L_t$. We firstly illustrate this phenomenon in a non-degenerate diffusive setting, i.e.\ when in \eqref{INTRO:eq:def_L_t} the operator $\mathcal{L}_t$ coincides with the classical Laplacian $\Delta_x$ acting on the whole space $\R^N$. In this framework and for bounded and suitably regular coefficients (i.e.\ $F$ in $C^{\frac{\beta}{2},\beta}$ in the notations below), Friedman in \cite{book:Friedman08} and Krylov in \cite{book:Krylov96} showed the following estimates:
\begin{equation}\label{INTRO:Stime_di_Schauder1}
    \Vert u \Vert_{C^{\frac{2+\beta}{2},2+\beta}} \, \le \, C \Vert f \Vert_{C^{\frac{\beta}{2},\beta}},
\end{equation}
where $\Vert \cdot \Vert_{C^{\gamma,\gamma'}}$ denotes the usual Hölder norm on the space-time $[0,T)\times \R^N$ with indexes $\gamma$ in time and $\gamma'$ in space. The Schauder estimates in \eqref{INTRO:Stime_di_Schauder1} then suggest that each solution $u$ is actually  Hölder regular of order $\frac{2+\beta}{2}$ in time and order $2+\beta$ in space if we suppose that the source $f$ is $\frac{\beta}{2}$-Hölder continuous in time and $\beta$-Hölder continuous in space. We point out that it is precisely the boundedness of the coefficients in space and their regularity in time that allows us to consider a parabolic bootstrap also with respect to the time variable. When considering unbounded in space coefficients but still of linear growth, Krylov and Priola in \cite{Krylov:Priola10} showed the following variation of the Schauder estimates:
\begin{equation}\label{INTRO:stime_di_Schauder2}
     \Vert u \Vert_{L^\infty(C^{2+\beta})} \, \le \, C \Vert f \Vert_{L^\infty(C^\beta)},
\end{equation}
where $\Vert \cdot \Vert_{L^\infty(C^\gamma)}$ represents the Hölder norm in space of order $\gamma$, uniformly in time. Instead, this kind of estimates only considers the regularity gain in space, of order $2$ in this case, of the solution $u$ with respect to the source $f$. We underline however that, at least in the uniformly elliptic case, it is still possible to deduce the regularity in time of the solution $u$ in a second moment.

An immediate consequence of the Schauder estimates is the uniqueness of solutions to Cauchy Problem \eqref{INTRO:eq:IPDE_degenere} in the considered functional space. Indeed, let $u_1$, $u_2$ be two solutions to Equation \eqref{INTRO:eq:IPDE_degenere}. Exploiting the linearity of the problem, we can see that their difference $u_1-u_2$ is then a solution of the following Cauchy problem:
\[
\begin{cases}
\left(\partial_t+L_t\right)(u_1-u_2)(t,x) \, = \, 0;\\
(u_1-u_2)(T,x) \, = \, 0.
\end{cases}\]
The Schauder estimates showed, for example, in \eqref{INTRO:stime_di_Schauder2} now imply the following control:
\[\Vert u_1-u_2 \Vert_{L^\infty(C^{2+\beta})} \, \le \, 0,\]
which immediately leads to the uniqueness of solutions in the considered functional space. The study of Cauchy problems such as \eqref{INTRO:eq:IPDE_degenere} is also fundamental in the analysis of the associated SDE \eqref{INTRO:SDE_Generale1}. Besides in the proof of the uniqueness of solutions for the martingale problem as explained at the end of the previous section, we recall that the Schauder estimates are often exploited in the Zvonkin transform in order to conclude the strong well-posedness for the SDE. For more details, see \cite{Zvonkin74,Veretennikov80} in a non-degenerate diffusive framework or \cite{Chaudru:Honore:Menozzi18_Strong, Hao:Wu:Zhang19} in a degenerate respectively diffusive and $\alpha$-stable one. Such estimates also appear in connection with some stochastic partial differential equations (SPDEs). For example, we refer to an application to the stochastic transport equation presented in \cite{Flandoli:Gubinelli:Priola10}, where Schauder estimates as in \eqref{INTRO:stime_di_Schauder2} are exploited to show the existence of a differentiable stochastic flow for the corresponding stochastic characteristics equation. Finally, we mention that a typical example of a degenerate system of the form \eqref{INTRO:eq:IPDE_degenere} is the following stable kinetic equation:
\[
 \partial_t u(t,x)  \,  = \, \mathcal{L}^\alpha u(t,x) -x_1 \cdot \nabla_{x_2}u(t,x) +f(t,x), \quad 
 \text{on } \R^{2d},
\]
where $x=(x_1,x_2)$  is in $\R^{2d}$ and $\mathcal{L}^\alpha$ is an $\alpha$-stable operator acting only on $x_1$, which naturally appears in the study of the linearized Boltzmann equation (cf.\ \cite{book:Villani02, Chen:Zhang18}).

We now briefly summarise the main results concerning Schauder estimates for degenerate parabolic systems. In a diffusive context, i.e.\ when in \eqref{INTRO:eq:IPDE_degenere} the operator $L_t$ can be rewritten as
\[L_t \, = \, \langle F(t,x),D_x\rangle +\frac{1}{2}{\rm Tr}\left(Ba(t,x)B^\ast D^2_x\right) \, =: \, \langle F(t,x),D_x\rangle+\mathcal{L}_t,\]
for some diffusion matrix $a(t,x)$ in $\R^d\otimes \R^d$, Lunardi in \cite{Lunardi97} was the first one to establish Schauder estimates for Kolmogorov equations of Ornstein-Uhlenbeck type, i.e.\ when $a(t,x)=1$ and $F(t,x)=Ax$ with a particular structure. This result has been obtained by exploiting anisotropic Hölder spaces, where the Hölder index depends on the considered spacial direction, precisely to control the different regularity orders caused by the degeneration of the system. See the next section for a more comprehensive explanation of this phenomenon and a precise definition of anisotropic spaces. Later, in \cite{Lorenzi05} and \cite{Priola09}, the authors obtained Schauder-type estimates for hypoelliptic Kolmogorov equations whose drift is partially non-linear, only along the components where the principal operator $\mathcal{ L}_t$ is non-degenerate, that is:
\[
F(x) \, = \, Ax +\begin{pmatrix}
    \tilde{F}_1(x)\\
    0_{(N-d),d}
\end{pmatrix}.
\]
Finally, the fully non-linear diffusive case with a diffusion matrix $a$ depending on time and space has been addressed for the first time in \cite{Chaudru:Honore:Menozzi18_Sharp}. Their result has been obtained under minimal regularity assumptions and the non-degeneracy on the coefficients $F$ and $a$, which leads to a weak Hörmander type condition for the system. The method of proof in \cite{Chaudru:Honore:Menozzi18_Sharp} is based on a perturbative approach by forward parametrix expansions that we will also exploit later on. Finally, we mention Di Francesco and Polidoro's work \cite{DiFrancesco:Polidoro06}, where local Schauder estimates are established for a linear degenerate system by assuming a different notion of continuity on the diffusion coefficient $a$ which, to a certain extent, also takes into account the linear transport associated with the drift.

In the last years, Schauder estimates for non-local operators, especially of the $\alpha$-stable type, have attracted great interest in the mathematical community (cf.\ \cite{Bae:Kassmann15, Bass09, Chaudru:Menozzi:Priola19, Dong:Kim13, Fernanadez:Ros-Oton17, Imbert:Jin:Shvydkoy18,Priola12,Ros-Oton:Serra16,Zhang:Zhao18}). However, almost all the known works focus only on the non-degenerate case. To the best of our knowledge, the only existing work dealing with the non-local degenerate case is \cite{Hao:Wu:Zhang19}, where the authors establish Schauder estimates for a stable kinetic system (i.e.\ when $N= 2d$ and $\rank B=d$ in \eqref{INTRO:eq:IPDE_degenere}). Their method is based on a Littlewood-Paley decomposition method already exploited in other works by the same author (cf.\ \cite{Zhang:Zhao18}), which is adapted to the natural anisotropic framework for degenerate systems.

In the current manuscript, we will analyse in detail two specific cases that generalise the previously cited results:
\begin{itemize}
\item a degenerate system with a non-linear drift in which the principal part of the operator is $\alpha$-stable. This model will be presented in detail in Section \ref{Sec:INTRO:Stime_Schauder_stabili};
\item a degenerate system driven by a Levy type Ornstein-Uhlenbeck operator. We refer to Section \ref{Sec:INTRO:Stime_Schauder_Levy} for a more detailed presentation of this model.
\end{itemize}

\subsection{Anisotropic geometry of degenerate systems}
In this section, we focus on understanding which is the most suitable functional space in which to establish our degenerate Schauder estimates. As mentioned earlier, Schauder estimates in a non-degenerate context are usually expressed in terms of the Hölder norm with respect to the ``usual'' Euclidean distance. Thus, we would like now to construct an Hölder space with respect to a new distance that is coherent with the degenerate structure of our system. In particular, at the end of this section, we will define a family of Hölder spaces with multi-indexes of regularity, depending on the considered coordinate.

To show to the reader how this anisotropic structure naturally appears when considering degenerate systems, we will present two different approaches: an analytical one, based on multi-scale dilation operators, and a more probabilistic one, which instead exploits the characteristic time scale of the solution process to the associated SDE. In order to expose them as clear as possible, we will also focus only on a linear example with only two components ($n=2$).

From an analytical point of view, we are interested in a Kolmogorov operator of $\alpha$-stable type for some $\alpha\in (0,2)$:
\[L^{\text{K}}\, := \, \Delta^{\frac{\alpha}{2}}_{x_1}+x_1\cdot \nabla_{x_2} \, \, \text { on } \R^{2d}\]
where $x=(x_1,x_2)$ is a point in $\R^{2d}$. Above, the operator $\Delta^{\frac{\alpha}{2}}_{x_1}$ represents the fractional Laplacian of order $\alpha/2$ with respect to the variable $x_1$, given by
\begin{equation}
\label{INTRO:def_laplaciano_fraz}
\Delta^{\frac{\alpha}{2}}_{x_1}\phi(x_1,x_2) \, := \, \text{p.v.}\int_{\R^d}\left[\phi(x_1+z,x_2)-\phi(x_1,x_2)\right]\frac{dz}{|z|^{d+\alpha}},
\end{equation}
for any sufficiently regular function $\phi\colon \R^{2d}\to \R$. In the diffusive case, i.e.\ when $\Delta^{\frac{\alpha}{2}}_{x_1}$ becomes the classical Laplacian $\Delta_{x_1}$ along $x_1$:
\[\Delta_{x_1}\phi(x_1,x_2)\,=\, \sum_{i=1}^d\partial^2_{x^i_1}\phi(x_1,x_2),\]
the parabolic equation involving such operator was indeed studied by Kolmogorov \cite{Kolmogorov34} and it was one of the first examples that inspired Hörmander's hypoelliptic theory \cite{Hormander67}. To understand how the different components of the system behave with respect to each other, we introduce a dilation operator $\delta\colon [0,\infty)\times \R^{2d}\to[0,\infty)\times \R^{2d} $ given by:
\[\delta(t,x) \, := \, (\delta_0t,\delta_1x_1,\delta_2x_2).\]
The correct values of $\delta_0,\delta_1,\delta_2$ are then to be determined so that the dilation $\delta$ is invariant for the equation:
\begin{equation}\label{INTRO:eq:dinamica_Kolmogorov}
   \partial_tu(t,x)+L^{\text{K}} u(t,x) \, = \, 0 \quad \text{ on }(0,\infty)\times\R^{2d}, 
\end{equation}
in the sense that it transforms solutions of the above equation into other solutions of the same one. The idea of a dilation $\delta$ which summarises the multi-scale, or anisotropic, behaviour of the considered degenerate system was initially introduced by Lanconelli and Polidoro (\cite{Lanconelli:Polidoro94}) for the study of the diffusive Kolmogorov equation. Since then, it has become a common tool for the analysis of the anisotropic geometry of degenerate equations, as evidenced by the rich literature in which it appears (\cite{Lunardi97,Huang:Menozzi:Priola19,Hao:Wu:Zhang19,DiFrancesco:Pascucci05,Metafune:Pallara:Priola02}). Exploiting now the natural scaling property of the fractional Laplacian, we can note that:
\[
\begin{split}
\left(\partial_t+L^{\text{K}}\right)(u\circ \delta) \, &= \, \delta_0(\partial_tu\circ \delta) + \delta^{\alpha}_1(\Delta^{\frac{\alpha}{2}}_{x_1}u\circ \delta)+\delta_2\left(x_1\cdot (\nabla_{x_2}u\circ \delta)\right) \\
&=\,\delta_0(\partial_tu\circ \delta) + \delta^{\alpha}_1(\Delta^{\frac{\alpha}{2}}_{x_1}u\circ \delta)+\delta^{-1}_1\delta_2\left([x_1\cdot \nabla_{x_2}u]\circ \delta\right)
= \, 0,
\end{split}
\]
where, we recall, we denoted $[x_1\cdot \nabla_{x_2}u]\circ \delta(t,x)  := \delta_1x_1\cdot \nabla_{x_2}u(\delta(t,x))$. In order to obtain homogeneous terms in the above equation, it is then natural to consider, for each fixed $\lambda>0$, the following dilation operator $\delta_\lambda$:
\begin{equation}
\label{INTRO:def_operatore_dilatazione}
\delta_\lambda(t,x_1,x_2) := (\lambda^\alpha t,\lambda x_1,\lambda^{1+\alpha}x_2).
\end{equation} 
In particular, we notice that, for our choice of $\delta_\lambda$, it holds that
\[\bigl(\partial_t +L^{\text{K}}\bigr) u = 0 \, \Longrightarrow \bigl(\partial_t +L^{\text{K}} \bigr)(u \circ
\delta_\lambda) = 0.\]
Summarising, the appearance of this multi-scale phenomenon is essentially due to the particular structure of the considered system, composed of a main part $\Delta^{\alpha/2}_{x_1}$ which provides a regularising effect of order $\alpha$ but only on the first component and of a transport term $x_1\cdot\nabla_{x_2}$ which allows this effect to be transmitted to the second component, even if with a lower intensity (of order $\alpha/(1+ \alpha)$). As we will see afterwards, such indexes will appear again in the parabolic bootstrap associated with the Schauder estimates for degenerate systems.

From a more probabilistic point of view, the scale indexes appearing in the dilatation operator $\delta$ can be found in the exponents of the characteristic times of an $\alpha$-stable process and its integral in time. In fact, the characteristic time of a multi-dimensional stochastic process can help explain the speed rate between the different components of the process itself. We start by considering the stochastic counterpart of Equation \eqref{INTRO:eq:dinamica_Kolmogorov}:
\begin{equation}
\label{INTRO_dinamica_Kolmogorov_stoc}
 \begin{cases}
dX^1_t \, = \, dZ_t,\\
dX^2_t \, = \, X^1_tdt,\quad t\ge0,
\end{cases}   
\end{equation}
where, for simplicity, we assumed that the solution process starts in the origin. Such stochastic dynamics is associated with the Kolmogorov operator introduced earlier in the sense that $\mathcal{L}^{\text{K}}$ is the infinitesimal generator of the process $X_t:=(X^1_t,X^2_t)$. SDE \eqref{INTRO_dinamica_Kolmogorov_stoc} can explicitly be solved through an integration in time:
\begin{equation}
\label{INTRO:dinamica_Kolmogorov_explic}
X_t \, = \, (X^1_t,X^2_t) \, = \, \left(Z_t, \int_0^t Z_s \, ds\right).   
\end{equation}
If we now check the characteristic times associated with the two components of the process $X_t$, we note that they are given by $(t^{\frac{1}{\alpha}},t^{1+\frac{1}{ \alpha}})$. Indeed, it is known that the stable process is $\alpha$-self-similar and its time integral simply adds one more order. We actually found the same scales shown in the dilation operator $\delta$, even if in this case, re-scaled with respect to the current time. In the diffusive case ($\alpha=2$), the multi-scale behaviour of the solution process $X_t$ appears even clearer. Indeed, thanks to the existence of finite second moments, it is possible to translate the above reasoning in terms of the covariance matrix of the process. If we replace the $\alpha$-stable process $Z_t$ with a Brownian motion $W_t$ in \eqref{INTRO_dinamica_Kolmogorov_stoc} and in \eqref{INTRO:dinamica_Kolmogorov_explic}, we immediately obtain that the solution $X_t=(X^1_t,X ^2_t)$ is a Gaussian process with zero mean and covariance $K_t$ in $\R^{2d}\otimes \R^{2d}$ given by
\[K_t \, = \, \begin{pmatrix}
               tI_{d\times d} & \frac{t^2}{2}I_{d\times d} \\
                 \frac{t^2}{2}I_{d\times d}     & \frac{t^3}{3}I_{d\times d}
             \end{pmatrix}.
\]
The equivalence, in terms of the associated quadratic forms, between the covariance matrix $K_t$ and a diagonal one was then shown in \cite{Konakov:Menozzi:Molchanov10}:
\[\sqrt{K_t} \, \asymp \, \begin{pmatrix}
               t^{\frac{1}{2}}I_{d\times d} & 0_{d\times d} \\
                 0_{d\times d}     & t^\frac{3}{2}I_{d\times d}
             \end{pmatrix},\]
where, for two matrixes $A,B$ in $\R^{nd}\otimes \R^{nd}$, the notation $A\asymp B$ means that there exists a constant $C\ge 1$ such that $C^{-1}|A\xi|^2\le |B\xi|^2 \le C|A\xi|^2$ for any $\xi$ in $\R^{nd}$. Such a property has been called by the authors a \emph{good scaling} property. See also Definition $3.2$ in \cite{Delarue:Menozzi10} for an extension to the general chain ($n\ge2$). It is now clear how one can find the same exponents shown previously on the diagonal entries of the scale matrix, even if only in the diffusive setting.

We now want to introduce a parabolic distance $\mathbf{d}_P$ on $[0,+\infty)\times \R^N$ which is homogeneous with respect to the multi-scale structure of the considered system, in the sense that:
\[\mathbf{d}_P\bigl(\delta_\lambda(t,x);\delta_\lambda(s,x')\bigr) = \lambda \mathbf{d}_P\bigl((t,x );(s,x')\bigr).\]
A natural choice is then given by:
\[\mathbf{d}_P\left((t,x),(s,x')\right) \, = \, |t-s|^{\frac{1}{\alpha}}+|x_1- x'_1|+|x_2-x'_2|^{\frac{1}{1+\alpha}}.\]
The newly introduced distance $d_P$ can be seen as a natural generalisation of the classical parabolic distance (cf.\ \cite{book:Krylov96,book:Friedman08}) to the multi-scale structure of our degenerate $\alpha$-stable dynamics. Since later on we will often use only the spacial part of the distance, let us introduce as well:
\begin{equation}
\label{INTRO:distanza_omogenea_d}
\mathbf{d}(x,x')\, := \, |x_1-x'_1|+|x_2-x'_2|^{\frac{1}{1+\alpha}}.
\end{equation}
For more details on homogeneous metrics, see also Stein's book \cite{book:Stein93}.

We are finally ready to introduce a functional space that is suitable for our purposes, an anisotropic Hölder space $C^\beta_{d}(\R^{2d})$ associated with the distance $\mathbf{d}_P$, in the sense that is homogeneous with respect to the dilation operators $\delta_\lambda$ defined in \eqref{INTRO:def_operatore_dilatazione}. In fact, the anisotropic semi-norm $\Vert \cdot \Vert_{C^\beta_d}$ will actually be considered component-wise. Indeed, fixed a coordinate, we will calculate the standard Hölder norm along that particular direction, but with a regularity index that is rescaled according to the dilatation operator $\delta_\lambda$ in \eqref{INTRO:def_operatore_dilatazione}, uniformly in time and with respect to the other coordinates. Finally, we will sum together all the contributions coming from the different components. More precisely, we can introduce for a function $\phi \colon \R^{2d} \to \R$ and a point $z$ in $\R^{2d}$, the function $\Pi^1_z\phi\colon\R^d\to \R$ given by
\[\Pi^1_z\phi(x_1) \, := \,  \phi(z_1+x_1,z_2).\]
A similar notation will be exploited for $\Pi^2_z\phi$ as well. Fixed a final time $T>0$, we define the \emph{homogeneous} space $L^\infty(0,T;C^\beta_{d}(\R^{2d}))$ as the family of all the Borel measurable functions $\phi\colon [0,T]\times \R^{2d}\to \R$ for which the following semi-norm is finite:
\[\Vert \phi \Vert_{L^\infty(C^\beta_{d})} \, := \, \sup_{t,z}\left([\Pi^1_z\phi(t,\cdot)]_{C^\beta(\R^d)}+[\Pi^2_z\phi(t,\cdot)]_{C^{\frac{\beta}{1+\alpha}}(\R^d)}\right)\, \asymp \, \sup_{t,x,x'}\frac{|\phi(t,x)-\phi(t,x')|}{\mathbf{d}(x,x')}.\]
Fixed $\alpha$ in $(0,2)$ such that $\alpha+\beta$ is in $(1,2)$, we can also define the anisotropic H\"older space $L^\infty(0,T;C^{\alpha+\beta}_{d}(\R^{2d}))$ of order $\alpha+\beta$ with respect to the following semi-norm:
\begin{equation}
\label{INTRO:def_holder_seminorm}
\begin{split}
    \Vert \phi \Vert_{L^\infty(C^{\alpha+\beta}_{d})} \, &:= \, \Vert D_{x_1}\phi \Vert_{L^\infty}\\
    &\qquad \qquad+\sup_{t,z}\left([\Pi^1_zD_{x_1}\phi(t,\cdot)]_{C^{\alpha+\beta-1}(\R^d)}+[\Pi^2_z\phi(t,\cdot)]_{C^{\frac{\alpha+\beta}{1+\alpha}}(\R^d)}\right),
\end{split}
\end{equation}
where, we recall, $\Vert \cdot \Vert_{L^\infty}$ represents the uniform norm on $[0,T]\times \R^{2d}$. If instead $\alpha+\beta>2$, a natural extension involving as well the second order derivatives is then necessary. Since the source $f$ and the terminal condition $u_T$ will be assumed to be bounded in our model, we also introduce the inhomogeneous version (with a $b$ at the bottom) of the anisotropic Hölder spaces just described, by simply adding the uniform norm of the function itself. For example,
\[\Vert \phi \Vert_{L^\infty(C^{\alpha+\beta}_{b,d})}\, := \, \Vert \phi \Vert_{L^\infty} + \Vert \phi \Vert_{L^\infty(C^{\alpha+\beta}_{d})}.\]
Finally, we will say that a function $\phi$ is in $C^{\alpha+\beta}_{b,d}(\R^{2d})$ if $\phi$ is independent from time and its anisotropic Hölder norm is finite. For example, this will be the case of the terminal condition $u_T$.

\setcounter{equation}{0}

\section{Schauder estimates for a degenerate $\alpha$-stable system with non-linear drift}
\fancyhead[RO]{Sezione \thesection. Schauder estimates for non-linear stable systems}
\label{Sec:INTRO:Stime_Schauder_stabili}

We summarise here the results presented in Chapter \ref{Chap:Schauder_estimates_Stable} of the present work, which have been published in \emph{Bulletin des Sciences Mathematiques}. Our goal is to establish optimal Schauder estimates, in terms of the minimal regularity assumptions on the coefficients, for the solutions to a degenerate parabolic integral-partial differential equation on $\R^{nd}$. The degeneration in this context comes from the fact that the main part of the operator, of $\alpha$-stable type, acts only on the first $d$ components of the system. More precisely, given a source $f\colon [0,T]\times \R^{nd} \to \R$ and a terminal condition $u_T\colon \R^{nd}\to \R$, we are interested in a Cauchy problem of the following form:
\begin{equation}
\label{INTRO:Degenerate_Stable_PDE}
\begin{cases}
   \partial_t u(t,x) + \langle Ax + F(t,x), D_{x}u(t,x)\rangle +
\mathcal{L}_\alpha u(t,x) \, = \, -f(t,x) & \mbox{on } [0,T)\times \R^{nd}; \\
    u(T,x) \, = \, u_T(x) & \mbox{on }\R^{nd}.
  \end{cases}
\end{equation}
where $x:=(x_1,\dots,x_n)$ is a point in $\R^{nd}$ with each $x_i$ in $\R^d$ and $\langle \cdot,\cdot \rangle$ represents the inner product on $\R^{nd}$. Above, $F\colon [0,T]\times \R^{nd}\to \R^{nd}$ is a sufficiently regular function and $A$ is a matrix in $\R^{{nd}} \otimes\R^{{nd}}$ on which we will impose suitable conditions. The operator $\mathcal{L}_\alpha$ is the infinitesimal generator of a symmetric, non-degenerate $\alpha$-stable process that acts only on the first component $x_1$ of the system. More precisely, the operator $\mathcal{L}_\alpha$ can be represented, for any smooth enough function $\phi\colon [0,T]\times\R^{nd}\to \R$, in the following form:
\begin{equation}
\label{INTRO:def_stable_operator}
\mathcal{L}_\alpha\phi(t,x)\, := \, \text{p.v.}\int_{\R^d}\bigl[\phi(t,x+Bz)-\phi(t,x) \bigr] \,\nu_\alpha(dz),
\end{equation}
where $B:= (I_{d\times d},0_{d\times d},\dots, I_{d\times d})^\ast$ is the immersion matrix from $\R^d$ to $\R^{nd}$ and $\nu_\alpha$ is the symmetric L\'evy measure on $\R^d$ associated with the $\alpha$-stable process. 

We now recall that the Lévy symbol $\Phi$ associated with the operator $\mathcal{L}_\alpha$ (or more exactly, with the process which the operator is the infinitesimal generator for) is usually defined using the Lévy-Khintchine formula which, in our symmetric and $\alpha$-stable case, can be stated as follows (cf.\ \cite{book:Sato99}):
\[
\Phi(p) \, = \, -\int_{\mathbb{S}^{d-1}}\vert p\cdot s \vert^\alpha \, \mu(ds),
\]
where ``$\cdot$'' denotes the inner product on the ``small'' space $\R^d$. Above, $\mu$ is a measure on the sphere $\mathbb{S}^{d-1}$ usually called the \emph{spectral} (or spherical) measure associated with $\nu_\alpha$, in the sense that a change of variables in polar coordinates $y=\rho s$, where $(\rho,s) \in (0,\infty)\times \mathbb{S}^{d-1}$, allows one to decompose the $\alpha$-stable Lévy measure $\nu_\alpha$ as
\begin{equation}
\label{INTRO:eq_decomposition_misura_stabile}
\nu_\alpha(dy) \, :=\, C_\alpha\mu(ds)\frac{d\rho}{\rho^{1+\alpha}}.
\end{equation}
For a proof of this fact, see for example Theorem $14.3$ in \cite{book:Sato99}. In particular, we will assume that the $\alpha$-stable Lévy measure $\nu_\alpha$ is symmetric and \emph{non-degenerate} in the sense that its spherical measure $\mu$ satisfies the following condition:
\begin{description}
  \item[{[ND]}] there exists a constant $\eta\ge 1$ such that for any $p$ in $\R^d$,
\begin{equation}
\eta^{-1}\vert p \vert^\alpha \, \le \, \int_{\mathbb{S}^{d-1}}\vert p\cdot s \vert^\alpha \, \mu(ds) \, \le\,\eta
\vert p \vert^\alpha.
\label{INTRO:def_non-deg}
\end{equation}
\end{description}
As we will see later, this condition implies in particular the existence of a fundamental solution for the operator $\mathcal{L}_\alpha$, since the Fourier transform of the $\alpha$-stable process $\{Z_t\}_{t\ge0}$ associated with $\mathcal{L}_\alpha$ is then integrable. It is also important to remark that the family of non-degenerate (in the above sense) spectral measures is very rich and variegate.
Indeed, condition [\textbf{ND}] is satisfied, for example, by the Lebesgue measure on the sphere, corresponding to the ``usual'' fractional Laplacian given by:
\[\mathcal{L}_\alpha\phi(t,x) \,:= \, \Delta^{\frac{\alpha}{2}}_{x_1}\phi(t,x) \, = \, \text{p.v.}\int_{\R^d}\left[\phi(t,x_1+z,x_2,\dots,x_n)-\phi(t,x)\right]\frac{dz}{|z|^{d+\alpha}},\]
but also by very singular (with respect to the Lebesgue measure) measures on the sphere, such as the sum of Dirac masses along the canonical coordinates which are associated with the following operator:
\begin{equation}
\label{INTRO:eq:def_stabile_cilindrico}
  \mathcal{L}_\alpha\phi(t,x) \,= \,  \sum_{i=1}^{d}\Delta^{\frac{\alpha}{2}}_{x^i_1}\phi(t,x) \,
\end{equation}
where, we recall, $x_1=(x_1^1,\dots,x_1^d)$ is a point in $\R^d$ and $\Delta^{\frac{\alpha}{2}}_{x^i_1}$ represents the scalar ``usual'' fractional Laplacian acting on the coordinate $x^i_1$. This type of operator is usually called a \emph{cylindrical} fractional Laplacian and in this case, the associated Lévy measure $\nu_\alpha$ is concentrated on the axes $\{x_1=0\}\cup\dots\cup\{x_d=0\}$. For more details, see e.g.\ Equation $(1.2)$ in \cite{Bass:Chen06}. Finally, we mention that in the literature the non-degeneracy condition [\textbf{ND}] on the Lévy measure $\nu_\alpha$ often appears as well in the following formulations:
\begin{itemize}
     \item (minimal support) the support of the spherical measure $\mu$ is not contained in any proper linear subspace of $\R^d$;
     \item (Picard condition) There exists a constant $C:=C(\alpha)$ such that for any $\rho>0$ and any $u$ in $\mathbb{S}^{d-1}$, it holds that
     \[\int_{\{|u\cdot y|\}\le \rho} |u\cdot y|^2 \nu_\alpha(dy)\, \ge \, C\rho^{2-\alpha }.\]
\end{itemize}
For more details and a proof of the equivalence between the above conditions, see for example \cite{Picard96, Sztonyk10_anisotropic, Priola12}.

As already explained in the introduction, the analysis of degenerate systems, where the regularising effect applies directly only on a subspace ($\R^d$) of the considered state space ($\R^{nd}$), requires some additional assumptions on the equation for the regularisation  to actually be transmitted throughout the components. In the diffusive case ($\alpha=2$), i.e. when $\mathcal{L}_\alpha=\Delta_{x_1}$, a natural condition is given by the Hörmander hypoellipticity (cf.\ \cite{Hormander67}) which imposes, at least formally, that the iterated $n-1$ Lie commutators associated with $\partial_{x_1}$ and $\langle Ax,D_x\rangle$ generate all the space. In our non-local case, although there does not seem to exist a general Hörmander theorem (cf.\ \cite{Komatsu:Takeuchi01}), some natural assumptions (conditions [\textbf{H}] and \textbf{[ND]}) ensure that the Markov semigroup generated by the operator
\[L^{\text{ou}} \, :=\, \mathcal{L}_\alpha +\langle Ax,D_x \rangle\]
admits a sufficiently smooth density \cite{Priola:Zabczyk09}. In this regard, we also refer to Cass' work \cite{Cass09} where a non-local  extension, even if still incomplete, of Hörmander result is presented, under very general conditions. More precisely, we will impose a particular structure for the matrix $A$ which ensures the hypoellipticity of the system:
\begin{description}
  \item[{[H]}] the matrix $A$ has the following sub-diagonal form:
\begin{equation}
\label{INTRO:eq:matrix_A_sub-diag}
A \, := \, \begin{pmatrix}
               0_{d\times d} & \dots         & \dots         & \dots     & 0_{d\times d} \\
               A_{2,1}       & 0_{d\times d} & \dots         & \dots     & 0_{d\times d} \\
               0_{d\times d} & A_{3,2}       & 0_{d\times d} & \dots     & 0_{d\times d} \\
               \vdots        & \ddots        & \ddots        & \ddots    & \vdots        \\
               0_{d\times d} & \dots         & 0_{d\times d} & A_{n,n-1} & 0_{d\times d}
             \end{pmatrix}
\end{equation}
  and each element $A_{i,i-1}$ in $\R^d\otimes\R^d$ has full rank $d$.
\end{description}
The specific structure of matrix $A$ appears natural (cf.\ \cite{Lanconelli:Polidoro94}) since it is invariant under the intrinsic dilations $\delta_\lambda$ (defined in \eqref{INTRO:def_operatore_dilatazione} for $n =2$) to our degenerate system, in the sense that
\begin{equation}
\label{INTRO:decomposition_exponential}
e^{tA} \, = \, \mathbb{M}_te^A\mathbb{M}^{-1}_t,   
\end{equation}
where $\mathbb{M}_t$ is a matrix in $\R^{nd}\otimes \R^{nd}$ given by
\begin{equation}
\label{INTRO:eq_matrix_Mt}
\bigl[\mathbb{M}_t\bigr]_{i,j} \, := \,
\begin{cases}
  t^{i-1}I_{d\times d}, & \mbox{if } i=j; \\
  0_{d\times d}, & \mbox{otherwise}.
\end{cases}    
\end{equation}
Intuitively, $\mathbb{M}_t$ takes into account the multi-scale structure associated with the distance $d$ but with respect to the characteristic time $t^{\frac{1}{\alpha}}$. The decomposition in \eqref{INTRO:decomposition_exponential} can easily be obtained from the definition of exponential matrix and from the identity $\mathbb{M}_tA\mathbb{M}^{-1}_t = tA$. Finally, we remark that our model only considers one particular specific structure among those possibly included in the general hypoellipticity theory developed by Hörmander. In fact, the non-degeneracy of the sub-diagonal elements in the matrix $A$ requires, at each level of the chain, to exploit a single additional Lie bracket to generate the corresponding direction. It is precisely this property that allows the transmission of the $\alpha$-stable regularising effects to each component of the chain, as explained in the previous section. Moreover, we highlight that in our non-linear context, the ``classical'' Hörmander condition (cf.\ \cite{Hormander67}) cannot be considered, due to the low regularity of the drift $F$ that, we will see, is only Hölder continuous in space. In particular, this prevents us from explicitly calculating the commutators needed in the usual Hörmander condition.

The anisotropic geometry associated to our degenerate system on $\R^{nd}$ can easily be understood as an extension on $n$ components of the one introduced in the previous section in the kinetic case ($n=2$), in the sense that, for example, the spacial distance $\mathbf{d}$ is now defined as
\begin{equation}
\label{INTRO:def_distanza_generica}
\mathbf{d}(x,x') \, = \, \sum_{i=1}^n|(x-x')_{i}|^{\frac{1}{1+\alpha(i-1)}}, \quad x,x' \in \R^{nd}.
\end{equation}
In particular, a function $\phi$ in $C^\beta_{d}(\R^{nd})$ is $\beta/(1+\alpha(i-1))$-Hölder continuous with respect to the coordinate $x_i$, uniformly in the other variables $x_j$ ($j\neq i$). The function $F=(F_1,\dots,F_n)$ can be understood as a (possibly non-linear) perturbation of the drift $Ax$ in the degenerate Ornstein-Uhlenbeck model $L^{\text{ou}}$. We also remark that, even further on, the non-linearity of the system we will mention refers to the drift shape. In fact, the perturbed operator $L^{\text{ou}}+\langle F(t,x), D_x\cdot\rangle$ is still a linear differential operator. We will impose a particular structure on this perturbation $F$ in such a way that it does not destroy the hypoellipticity of the system. We will also assume a certain degree of regularity for the drift $F$, which is necessary for our purposes.
\begin{description}
  \item[{[R]}] for any $i$ in $\llbracket 1,n\rrbracket$, $F_i$ only depends on time and on the last $n-(i-1)$ spacial variables, i.e.\ $F_i(t,x)=F_i(t,x_i,\dots,x_n)$. Moreover, $F_i$ belongs to $L^\infty(0,T;C^{\gamma_i+\beta}_{d}(\R^{nd}))$, where
    \begin{equation}\label{INTRO:Drift_assumptions}
    \gamma_i \,:= \,
    \begin{cases}
        1+ \alpha(i-2), & \mbox{if } i>1; \\
        0, & \mbox{if } i=1.
   \end{cases}
   \end{equation}
\end{description}
We point out that no boundedness condition in space has been imposed on the drift $F$ but only a Hölder regularity assumption with increasing multi-indexes. In fact, this will be one of the main difficulties we will have to deal with in our method of proof. We also remark that, unlike the non-degenerate case (cf.\ \cite{Chaudru:Menozzi:Priola19}), here it is necessary to impose an additional increasing regularity on the degenerate components ($i>1$) of the perturbation $F_i$, represented by the $\gamma_i$ parameter above. Indeed, this assumption appears natural if one thinks that, due to the degenerate structure of the system, the regularising effect of the $\alpha$-stable operator $\mathcal{L}_\alpha$, which only acts on the first coordinate, becomes weaker and weaker as it gets transmitted throughout the components of the system. In a sense, the additional regularity on $F$ is the price to pay for re-balancing the increasing time singularities appearing while descending along the components.

Finally, it is necessary to impose some conditions on the range of stability indexes $\alpha$ in $(0,2)$ and Hölder regularity indexes $\beta$ in $(0,1)$ we can consider. More precisely,
\begin{description}
   \item[{[P]}] $\alpha+\beta<2$ and if $\alpha<1$, it also holds that
   \[\beta<\alpha, \quad \alpha+\beta>1,\quad 1-\alpha <\frac{\alpha-\beta}{1+\alpha(n-1)}.\]
\end{description}
Some considerations on the additional conditions in the super-critical case ($\alpha<1$) are now necessary. The limitation $\beta<\alpha$ essentially reflects the weak integrability property (of order strictly lower than $\alpha$) for a possibly non-isotropic, $\alpha$-stable, process. The condition $\alpha+\beta>1$ naturally appears to give a point-wise meaning to the gradient of the solution $u$ with respect to the non-degenerate variable $x_1$. In this regard, we also mention Tanaka's work \cite{Tanaka:Tsuchiya:Watanabe74} where it is shown that the weak well-posedness, a property strictly interconnected with the Schauder estimates considered here, can already fail for a one-dimensional  SDE driven by an additive, non-degenerate $\alpha$-stable noise if $\alpha+\beta< 1$, where $\beta$ is the Hölder regularity of the deterministic drift. This counterexample can be understood as a stochastic generalisation of the Peano's example shown in \eqref{INTRO:eq:solution_Peano}. The last assumption is instead a technical condition and it seems necessary for our method of proof via forward parametrix expansions to work. Finally, we remark that in the sub-critical case, when $\alpha \ge 1$, these conditions are always satisfied.

Due to the low regularity assumed on the coefficients, it is possible to consider Equation \eqref{INTRO:Degenerate_Stable_PDE} only in a distributional sense. Indeed, the expected regularity (through the parabolic \emph{bootstrap} in space given by the Schauder estimates) for a ``classical'' solution $u$ of the problem is in the space $L^\infty(0,T;C^{\alpha+\beta}_{b,d}(\R^N))$, not enough for giving a point-wise meaning to the gradient $D_xu$ with respect to the degenerate variables and consequently, for giving a point-wise sense to Equation \eqref{INTRO:Degenerate_Stable_PDE}. Instead, we will consider here two other notions of solution, that are more suitable for our purposes. As mentioned before, by weak solution of Equation \eqref{INTRO:Degenerate_Stable_PDE} we essentially mean a solution in the sense of distributions, i.e.\ a function $u \colon [0,T]\times \R^{nd} \to \R $ such that for any test function $\phi$ (smooth function with compact support) on $(0,T]\times \R^{nd}\to \R$, it holds that
\begin{multline*}
     \int_{0}^{T}\int_{\R^{nd}}\Bigl(-\partial_t+(L_t)^*\Bigr)\phi(t,y)u(t,y) \, dy +\int_{\R^{nd}}u_T(y)\phi(T,y) \, dy \\
= \, - \int_{0}^{T}\int_{\R^{nd}}\phi(t,y)f(t,y) \, dy,
\end{multline*}
where the operator $(L_t)^\ast$ represents the (formal) adjoint of $L_t$ given by
\[L_t \, = \, \langle Ax +F(t,x) ,D_x \rangle + \mathcal{L}_\alpha \,\,\text{ on }\,\, \R^{nd} .\]
Instead, a \emph{mild} solution (in the sense of Stroock and Varadhan \cite{book:Stroock:Varadhan79}) to Equation \eqref{INTRO:Degenerate_Stable_PDE} is a function $u \colon [0,T]\times \ R^{nd} \to \R$ obtained as the limit, in a suitable functional space, of the sequence of the classical solutions to regularised versions of the considered Cauchy problem. For more details, see Definition \ref{definition:mild_sol} in Chapter \ref{Chap:Schauder_estimates_Stable} or the book \cite{book:Kolokoltsov11}.

The main results established in Chapter \ref{Chap:Schauder_estimates_Stable} can now be summarised in the following theorem:
\begin{theorem}
Under the above assumptions, let $f$ be in $L^\infty(0,T;C^\beta_{b,d}(\R^{nd}))$ and $u_T$ in $C^{\alpha+\beta}_{b,d}(\R^{nd})$. Then, there exists a unique mild and weak solution $u\colon [0,T]\times \R^{nd}\to \R$ to Cauchy problem \eqref{INTRO:Degenerate_Stable_PDE}. Moreover, $u$ belongs to $L^\infty(0,T;C^{\alpha+\beta}_{b,d}(\R^{nd}))$ and there exists a constant $C$, independent from $f$ and $u_T$, such that
\begin{equation}
\label{INTRO:eq:Schauder_estimates}
\Vert u \Vert_{L^\infty(C^{\alpha+\beta}_{b,d})} \, \le \, C \bigl[\Vert f \Vert_{L^\infty(C^{\beta}_{b,d})} + \Vert u_T
\Vert_{C^{\alpha+\beta}_{b,d}}\bigr].
\end{equation}
\end{theorem}

In conclusion, we remark that it is possible, with minor modifications to the arguments presented below, to consider as well degenerate systems with a completely non-linear deterministic drift or a space-time dependent diffusion coefficient in a suitable functional space, as explained at the end of Chapter \ref{Chap:Schauder_estimates_Stable}.

\subsection{Sketch of the proof}

We briefly present here the method of proof we follow to exhibit the Schauder estimates \eqref{INTRO:eq:Schauder_estimates} under the assumptions introduced in the previous section. The main difficulties to face in our context will be related to the degeneration of the $\alpha$-stable operator $\mathcal{L}_\alpha$  which only acts on the first component and to the unboundedness of the perturbation $F$. We also recall that we want to establish the Schauder estimates under the minimal Hölder regularity assumptions for the coefficients of the equation. In particular, we will not be able to rely on derivatives along the degenerate components but instead, we will exploit duality properties between Hölder and Besov spaces and in particular, establish delicate Besov norm controls.

Our approach is based on a perturbative method known as \emph{parametrix technique}, originally introduced by Levi \cite{Levi1907} for the analysis of linear elliptic PDEs of even order with variable coefficients. In the non-degenerate diffusive framework, we cite the works by Friedman \cite{book:Friedman08} and by McKean and Singer \cite{Mckean:Singer67} who exploit this technique to obtain Aronson-type estimates (cf.\ \cite{Aronson59, Aronson67}) for the fundamental solution of the system, respectively under an assumption of Hölder regularity in space-time and one of measurability in time and Hölder continuity in space, of the coefficients. In a degenerate diffusive context more similar to the present one, this method was exploited in \cite{Chaudru:Honore:Menozzi18_Sharp} in order to establish the Schauder estimates for a degenerate diffusive chain. Finally, we mention that Hadamard in \cite{book:Hadamard32, book:Hadamard64} extended this technique to the study of hyperbolic equations. Another, more classic, method for obtaining Schauder estimates is to proceed by exploiting a priori controls on the fundamental solution associated with the system. Existence and uniqueness of solutions to the considered PDE, in a suitable functional space, are in this case addressed only in a second moment. We refer to \cite{book:Friedman08} and \cite{book:Krylov96} for a clear presentation of this a priori approach or to \cite{Krylov:Priola10} for its extension to the non-degenerate diffusive case with unbounded coefficients.

\subsubsection{The frozen proxy operator}

The crucial element in the parametrix expansion method  consists in choosing a suitable \emph{proxy} operator for the equation of interest, such that an operator $\tilde{L}_t$ whose properties (existence and regularity of the density or of the associated Markov semigroup) are known and that is close, in a suitable sense, to the original operator $L_t$. We can then apply a first order expansion, such as a Duhamel formula, to solve PDE \eqref{INTRO:Degenerate_Stable_PDE} around the chosen proxy operator. By exploiting the known properties of the proxy operator, we could finally be able to obtain a suitable bound on the expansion error.

When dealing with bounded coefficients, a common choice for the proxy is given by the original operator with constant coefficients. When potentially unbounded perturbations $F$ are considered, as in our case, it is instead natural (cf.\ \cite{Krylov:Priola10, Chaudru:Menozzi:Priola19}) to use an operator involving a non-zero first order term, like the one given by the flow associated with the transport term $Ax+F$:
\begin{equation}
\label{INTRO:eq:flusso_associato}
 \begin{cases}
  d \theta_{\tau,s}(\xi) \, = \, \bigl[A\theta_{\tau,s}(\xi)+F(s,\theta_{\tau,s}(\xi))\bigr]ds \,\,\, \mbox{ on
  } [\tau,T];\\
  \theta_{\tau,\tau}(\xi) \, = \, \xi,
\end{cases}   
\end{equation}
where $(\tau,\xi)$ in $[0,T]\times \R^N$ are two fixed \emph{freezing} parameters whose exact value will be chosen afterwards. However, we immediately remark that a solution to the above dynamics may not be unique, since the drift $F$ is only Hölder continuous in space. For this reason, we will have to choose a particular flow, denoted by $\theta_{\tau,s}(\xi)$, and consider it fixed from now on. The suitable proxy operator is then obtained by freezing the original one $L_t$ along such fixed flow $\theta_{\tau,t}(\xi)$:
\[\tilde{L}^{\tau,\xi}_t \, = \, \mathcal{L}_\alpha+\langle A x + F(t,\theta_{\tau,t}(\xi)), D_{x}\rangle.\]

We can now introduce the ``frozen'' Cauchy problem
associated with the proxy operator:
\begin{equation}\label{INTRO:Frozen_PDE}
\begin{cases}
   \left(\partial_t+\tilde{L}^{\tau,\xi}_t\right) \tilde{u}^{\tau,\xi}(t,x) \, = \, -f(t,x) &
   \mbox{on }   [0,T)\times \R^{nd}; \\
    \tilde{u}^{\tau,\xi}(T,x) \, = \, u_T(x) & \mbox{on }\R^{nd}.
  \end{cases}
\end{equation}
To determine the properties of the frozen proxy operator $\tilde{L}^{\tau,\xi}$, we now switch to the corresponding stochastic dynamics. Given an initial point $(t,x)$ in $[0,T]\times \R^{nd}$ and a symmetric $\alpha$-stable process $\{Z_s\}_{s\ge t} $ on $\R^d$ with Lévy measure given by $\nu_\alpha$, we are interested in
\[
\begin{cases}
\tilde{X}_s^{\tau,\xi} \, &= \,  \left[A\tilde{X}_s^{\tau,\xi} + F(s, \theta_{\tau,s} (\xi))\right]ds + BdZ_s; \quad s>t\\
\tilde{X}_t^{\tau,\xi} \, &= \, x.
\end{cases}
\]
An integration in time against the exponential matrix then allows to obtain a more explicit representation of the solution process:
\begin{equation}
\label{INTRO:eq:rappresentazione_X_tilda}
\begin{split}
\tilde{X}_s^{\tau,\xi} \, &= \, e^{A(s-t)}x + \int_{t}^{s}e^{A(s-u)} F(u, \theta_{\tau,u}
(\xi)) \, du + \int_t^se^{A(s-u)}BdZ_u \\
&=:\, \tilde{m}^{\tau,\xi}_{t,s}(x) +  \int_t^se^{A(s-u)}BdZ_u.
\end{split}
\end{equation}
Thanks to the symmetry of $Z_t$, one can think of the transport term $\tilde{m}^{\tau,\xi}_{t,s}(x)$ as the ``mean'' of the frozen process $\tilde{X}_s^{\tau,\xi}$, affine to the initial point $x$ (although when $\alpha<1$, the expected value of the process is actually not defined) or at least as the value around which the process fluctuates. The above identity is now crucial to show that the stochastic convolution
\[\Lambda_{t,s}\, := \, \int_t^s e^{A(s-v)}BdZ_v\]
is again a symmetric, non-degenerate $\alpha$-stable process on $\R^{nd}$ but rescaled according to the anisotropic structure of the degenerate system. Indeed, through a reasoning in Fourier spaces, one can show that
\begin{equation}
    \label{INTRO:prova1}
    \begin{split}
\mathbb{E}\Bigl[\text{exp}\Bigl(i\langle p,\Lambda_{t,s}\rangle\Bigr)\Bigr] \, &= \,
\text{exp}\Bigl(-\int_{t}^{s}\int_{\mathbb{S}^{d-1}}\left\vert \langle p,e^{(s-u)A}B\varsigma\rangle\right\vert^\alpha \, \mu(d\varsigma)du\Bigr) \\
&=: \text{exp}\Bigl(\Phi_{\Lambda_{t,s}}(p)\Bigr),
\end{split}
\end{equation}
so that $\Phi_{\Lambda_{t,s}}$ is the Lévy symbol for the random variable $\Lambda_{t,s}$ at any fixed times $t<s$. By the change of variable $v=(s-u)/(s-t)$, we can then rewrite it in the following way:
\[
\begin{split}
\Phi_{\Lambda_{t,s}}(p) \, = \, (t-s)\int_{0}^{1}\int_{\mathbb{S}^{d-1}}\vert \langle  p, e^{(s-t)vA}B\varsigma \rangle
\vert^\alpha \, \mu(d\varsigma)dv\\
(t-s)\int_{0}^{1}\int_{\mathbb{S}^{d-1}}\vert \langle  \mathbb{M}_{s-t}p, e^{vA}B\varsigma \rangle
\vert^\alpha \, \mu(d\varsigma)dv,
\end{split}
\]
where, in the last step, we used the decomposition of $e^{A(s-t)v}$ given in \eqref{INTRO:decomposition_exponential}:
\[e^{(s-t)vA}B \, = \, \mathbb{M}_{s-t}e^{vA}B,\]
where the intrinsic scale matrix $\mathbb{M}_{s-t}$ has been defined in \eqref{INTRO:eq_matrix_Mt}. Furthermore, we can now re-normalise the term inside the inner product:
\[
\begin{split}
\Phi_{\Lambda_{t,s}}(p) \, &= \,
(t-s)\int_{0}^{1}\int_{\mathbb{S}^{d-1}}\vert \langle \mathbb{M}_{s-t}p,\frac{e^{vA}B\varsigma}{\vert e^{vA}B\varsigma\vert}\rangle\vert^\alpha\,\vert e^{vA}B\varsigma\vert^\alpha \mu(d\varsigma)dv\\
&=: \, (t-s)\int_{[0,1]\times \mathbb{S}^{d-1}}\vert \langle \mathbb{M}_{s-t}p,l(v,\varsigma)\rangle\vert^\alpha\,m_\alpha(d\varsigma,dv),
\end{split}\]
where $m_\alpha(d\varsigma,dv) := \vert e^{vA}B\varsigma \vert^\alpha \mu(d\varsigma)dv$ is a product measure on $[0,1]\times \mathbb{S}^{d-1}$ and $l\colon [0,1]\times \mathbb{S}^{d-1}\to \mathbb{S}^{nd-1}$ is the \emph{lift} function given by
\[l(v,\varsigma) \, := \, \frac{e^{vA}B\varsigma}{\vert e^{vA}B\varsigma \vert}.\]
In particular, we denote by $\mu_S\, := \, \text{Sym}(l_\ast(m_\alpha))$ the symmetrization of the pushforward measure from $m_\alpha$ through $l$.  Then, 
\begin{equation}
\label{INTRO:prova2}
\mathbb{E}\Bigl[\text{exp}\Bigl(i\langle p,\Lambda_{t,s}\rangle\Bigr)\Bigr] \, = \,
\text{exp}\Bigl(-(s-t)\int_{\mathbb{S}^{d-1}}\vert \langle \mathbb{M}_{s-t}p,\eta\rangle\vert^\alpha\,\mu_S(d\eta)\Bigr).    
\end{equation}
If we now denote by $\{S_u\}_{u\ge 0}$ the Lévy process on $\R^{nd}$ whose Lévy symbol $\Phi_{S}$ is given by
\[\Phi_{S}(p) \, = \, -\int_{S^{nd-1}}| \langle p, \eta \rangle |^\alpha\mu_S(d\eta),\]
we finally obtain the following identity:
\begin{equation}
\label{INTRO:decomposition}
\Lambda_{s,t} \, \overset{(\text{law})}{=} \, \mathbb{M}_{s-t}S_t.
\end{equation}
We would like now to point out that, even if the spectral measure $\mu_S$ is non-degenerate (in the sense of condition [\textbf{ND}]), such a measure is strongly singular with respect to the Lebesgue measure on $\mathbb{S}^{nd-1}$. Indeed, it is not difficult to check that its support is given by the image on  $\mathbb{S}^{nd-1}$ of the support of the measure $\vert e^{vA}B\varsigma\vert^\alpha \mu(d\varsigma)dv$ through the lift function $(v,\varsigma) \mapsto e^{vA}B\varsigma/\vert e^{vA}B\varsigma \vert$. 
In particular, even assuming that the support of the spectral measure $\mu$ of the driving process $\{Z_t\}_{t\ge 0}$ is the full sphere $\mathbb{S}^{d-1}$, 
the dimension of the support of $\mu_S$ will only be $d-1+1=d$.

Under the non-degeneracy condition [\textbf{ND}], we now know that the process $\{S_t\}_{t\ge 0}$ admits a regular density $p_S(t,z)$. Equations \eqref{INTRO:eq:rappresentazione_X_tilda} and \eqref{INTRO:decomposition} now imply the existence of a density $\tilde{p}^{\tau,\xi}$ associated with the frozen operator $\tilde{L}^{\tau,\xi}_t$ given by
\begin{equation}\label{INTRO:eq:definition_tilde_p}
\tilde{p}^{\tau,\xi}(t,s,x,y)\,  = \, \frac{1}{\det (\mathbb{M}_{s-t})} p_S\bigl( s-t,\mathbb{M}^{-1}_{s-t} (y -
\tilde{m}^{ \tau,\xi}_{t,s}(x))\bigr).
\end{equation}
In particular, the estimates on the density of the frozen process $\tilde{X}^{\tau,\xi}_s$ will be inferred from the corresponding ones for the density of $\{S_u\}_{u\ge 0}$. For completeness, let us introduce as well the ``frozen'' Markov semigroup $\bigl{\{}\tilde{P}^{\tau,\xi}_{t,s}\bigr{\}}_{t\le s}$ given by
\[\tilde{P}^{\tau,\xi}_{t,s}\phi(x) \, := \, \int_{\R^{nd}}\tilde{p}^{\tau,\xi}(t,s,x,y)\phi(y) \,dy,\]
for any regular enough function $\phi\colon \R^{nd}\to \R$. We now want to understand which regularity properties this frozen density $\tilde{p}^{\tau,\xi}(t,s,x,\cdot)$ possesses. Firstly, we will show that its derivatives in space can be controlled by another density but at the cost of additional singularities in time. More precisely, we will prove that for any  $i$ in $\llbracket 1, n\rrbracket$ and $k\in\{1,2\}$, it holds that
\begin{equation}\label{INTRO:eq:Smoothing_effects}
\bigl{\vert}
D^k_{x_i}
\tilde{p}^{ \tau,\xi} (t,s,x,y) \bigr{\vert} \, \le \, C  \frac{\bar{p}\bigl( s-t,\mathbb{M}^{-1}_{s-t} (y -
\tilde{m}^{ \tau,\xi}_{t,s}(x))\bigr)}{\det (\mathbb{M}_{s-t})}(s-t)^{-k \frac{1+\alpha(i-1)}{\alpha}}.
\end{equation}
where $\bar{p}(t,\cdot)$ is a density  on $\R^{nd}$ with suitable properties. Roughly speaking, a space derivative of $\tilde{p}^{\tau,\xi}(t,s,x,\cdot)$ induces a singularity in time whose intensity depends on the derivation direction, or said differently, on the anisotropic nature of our system. 

A fundamental and recurring element in our analysis will consist of re-balancing these singularities in time through the space regularity of the density $\bar{p}(t,\cdot)$. In particular, we will show that it posses a smoothing effect in the space of order $\alpha$, in the sense that, for every $\gamma$ in $[0,\alpha)$, it holds that
\begin{equation}
\label{INTRO:eq:Smoothing_effects2}
\int_{\R^{nd}} \frac{\bar{p}\bigl( s-t,\mathbb{M}^{-1}_{s-t} (y -
\tilde{m}^{ \tau,\xi}_{t,s}(x))\bigr)}{\det (\mathbb{M}_{s-t})} \mathbf{d}^\gamma(y,\tilde{m}^{ \tau,\xi}_{t,s}(x))\, dy \, \le \, C(s-t)^{\frac{\gamma}{\alpha}},   
\end{equation}
where, we recall, the homogeneous distance $\mathbf{d}$ has been defined in \eqref{INTRO:distanza_omogenea_d} when $n=2$ and then extended to our more general framework in \eqref{INTRO:def_distanza_generica}. 
We can now highlight one of the main differences between the diffusive case, considered for example in \cite{Chaudru:Honore:Menozzi18_Sharp}, and the $\alpha$-stable one studied here, characterised by a regularising effect strictly bounded by the stability order of the operator. From a stochastic point of view, this fact is essentially related to the weaker integrability properties associated with Lévy processes with respect to the Brownian motion.

The properties presented in \eqref{INTRO:eq:Smoothing_effects} and \eqref{INTRO:eq:Smoothing_effects2} can finally be summarised in the regularising effect for the proxy operator $\tilde{L}^{\tau,\xi}_t$ in the following fundamental bound for its frozen semigroup:
\begin{equation}\label{INTRO:eq:Control_of_semigroup}
\bigl{\vert} D^k_{x_i}\tilde{P}^{\tau,\xi}_{t,s}\phi(x) \bigr{\vert} \, \le \, C \Vert\phi\Vert_{C^\gamma_d}
(s-t)^{\frac{\gamma}{\alpha}-k\frac{1+\alpha(i-1)}{\alpha}}, \quad \forall \phi \in C^\gamma_d(\R^{nd}).
\end{equation}
Controls like \eqref{INTRO:eq:Control_of_semigroup} 
can be established from the smoothing properties of the density in \eqref{INTRO:eq:Smoothing_effects}, exploiting some \emph{cancellation} techniques. Noticing that
\begin{equation}
\label{INTRO:eq:cancellazione0}
    \int_{\R^{nd}}D_{x_i}\tilde{p}^{\tau,\xi}(t,s,x,y)\, dy\,=\,0,
\end{equation}
for any freezing parameter $(\tau,\xi)$, the underlying idea consists of introducing an additional term inside the integral in order to rely on the Hölder regularity of $\phi$, i.e.:
\begin{equation}
\label{INTRO:eq:cancellazione1}
\begin{split}
\bigl{\vert}D^k_{x_i}\tilde{P}^{\tau,\xi}_{t,s}\phi(x) \bigr{\vert}  \, &= \, \int_{\R^{nd}}    D^k_{x_i}\tilde{p}^{\tau,\xi}(t,s,x,y)\phi(y)\, dy \\
&= \, \int_{\R^{nd}}  D^k_{x_i}\tilde{p}^{\tau,\xi}(t,s,x,y)\left[\phi(y)- \phi(\tilde{m}^{\tau,\xi}_{t,s}(x))  \right]\, dy\\
&\le \, \Vert \phi \Vert_{C^\gamma_d}\int_{\R^{nd}}  \left|D^k_{x_i}\tilde{p}^{\tau,\xi}(t,s,x,y)\right|d^{\gamma}(y-\tilde{m}^{\tau,\xi}_{t,s}(x))\, dy \\
&\le \, C\Vert \phi \Vert_{C^\gamma_d} (s-t)^{\frac{\gamma}{\alpha}-k\frac{1+\alpha(i-1)}{\alpha}}.
\end{split}
\end{equation}

By carefully exploiting the smoothing effects in space in order to balance the singularities in time coming from the derivation of the density (cf.\ Control \eqref{INTRO:eq:Smoothing_effects}), it is now possible to show that the Schauder estimates hold for the solution $\tilde{u}^{\tau,\xi}$ to the frozen Cauchy problem \eqref{INTRO:Frozen_PDE}:
\begin{equation}\label{INTRO:eq:Schauder_Estimates_for_proxy}
\Vert \tilde{u}^{\tau,\xi} \Vert_{L^\infty(C^{\alpha+\beta}_{b,d})} \, \le \, C\bigl[\Vert f \Vert_{L^\infty(C^{\beta}_{b,d})}+\Vert u_T \Vert_{C^{\alpha+\beta}_{b,d}} \bigr].
\end{equation}
In particular, these estimates are valid for any \emph{fixed} value of the freezing parameters $(\tau,\xi)$ and the constant $C$ above is independent from them. Thanks to the stability given by the above bounds, it is also possible to represent the unique solution $\tilde{u}^{\tau,\xi}$ to the frozen Cauchy problem \eqref{INTRO:Frozen_PDE} in terms of the Markov semigroup $\tilde{P}^{\tau,\xi}_{t,s}$:
\begin{equation}\label{INTRO:Duhamel_representation_of_proxy}
\tilde{u}^{\tau,\xi}(t,x) \, = \, \tilde{P}^{\tau,\xi}_{t,T}u_T(x) + \int_{t}^{T}\tilde{P}^{\tau,\xi}_{t,s}f(s,x) \,ds.
\end{equation}
The next step consists of applying the actual perturbative method in order to infer the Schauder estimates for a solution $u$ to the original PDE from those, just obtained in \eqref{INTRO:eq:Schauder_Estimates_for_proxy}, for the solution $\tilde{u}^ {\tau,\xi}$ to the frozen Cauchy problem associated with the proxy operator.

\subsubsection{A Duhamel type expansion}

At least formally (in practice, through a mollification of the coefficients using the definition of mild solution), the original system \eqref{INTRO:Degenerate_Stable_PDE} can be rewritten around the proxy operator $\tilde{L}^{\tau,\xi} _{\alpha}$ as follows:
\[\begin{cases}
\left(\partial_t + \tilde{L}^{\tau,\xi}_t\right)u(t,x) \, = \, -f(t,x)- \left(L_t-\tilde{L}^{\tau,\xi}_t\right)u(t,x), &
   \mbox{on }   (0,T)\times \R^{nd}; \\
    u(T,x) \, = \, u_T(x) & \mbox{on }\R^{nd}.
\end{cases}
\]
The uniqueness of solutions for the frozen Cauchy problem \eqref{INTRO:Frozen_PDE} and the representation of $\tilde{u}^{\tau,\xi}$ in \eqref{INTRO:Duhamel_representation_of_proxy} now imply the following Duhamel formula, which corresponds to a first order parametrix expansion:
\begin{equation}\label{INTRO:eq:Expansion_along_proxy}
u(t,x) \, = \, \tilde{u}^{\tau,\xi}(t,x) + \int_{t}^{T}  \tilde{P}^{\tau,\xi}_{t,s}R^{\tau,\xi}(s,x)\,
ds,
\end{equation}
where $R^{\tau,\xi}$ is the remainder term given by
\begin{equation}\label{INTRO:eq:def_remainder_regul}
R^{\tau,\xi}(t,x)\,:=\, \langle F(t,x)-F(t,\theta_{\tau,t}(\xi)),D_{x}u(t,x)
\rangle.
\end{equation}
Since we already established the suitable control \eqref{INTRO:eq:Schauder_Estimates_for_proxy} for the frozen solution $\tilde{u}^{\tau,\xi}$, we can infer from the Duhamel formula \eqref{INTRO:eq:Expansion_along_proxy}  that in order to obtain the Schauder estimates for $u$, it only remains to investigate the term associated with the remainder:
\begin{equation}\label{INTRO:eq:qqq}
  \int_{t}^{T}\tilde{P}^{\tau,\xi}_{t,s}R^{\tau,\xi}(s,x) \, ds,
\end{equation}
that precisely represents the approximation error around the proxy.

Up until now, the freezing parameters $(\tau,\xi)$ have been considered fixed but free. They will now be chosen appropriately according to the control we want to establish. In particular, in this work we will follow a \emph{forward} parametrix approach, in the sense that we will impose $(\tau,\xi)=(t,x)$ and therefore the flow $\theta_{\tau,s}(\xi)$ given in \eqref{eq:_INtro_def_flusso} will move forwardly from the initial point $(t,x)$ towards $(s,y)$, where $y$ is the integration variable in the frozen density. This method has been widely used by Friedman \cite{book:Friedman08} and Il'in et al.\ \cite{ilcprimein:Kalavsnikov:Oleuinik62} to establish point-wise estimates on the derivatives of the fundamental solution for the heat equation or by Chaudru de Reynal in \cite{Chaudru17} to derive the strong well-posedness for a degenerate kinetic-type diffusive SDE. In particular, the forward approach allows to exploit more easily the cancellation techniques in \eqref{INTRO:eq:cancellazione1} which we have seen to be crucial in the control of the derivatives of the frozen density $\tilde{p}^{\tau,\ xi}(t,s,x,y)$.

For completeness, we mention that there is also a \emph{backward} parametrix approach (obtained fixing $(\tau,\xi)=(s,y)$), introduced by McKean and Singer in \cite{Mckean:Singer67}. However, we remark that in this case the density $\tilde{p}^{\tau,\xi}(t,s,x,y)$ frozen in $\xi=y$ is no longer a true density with respect to the variable $y$, since the freezing parameter also acts as an integration variable. This makes the smoothing effects presented in \eqref{INTRO:eq:Control_of_semigroup} for the Markov semigroup much more difficult to establish. For more details on the perturbative method by backward parametrix expansions, we refer to Section \ref{Sec:INTRO:Buona_posizione_Levy}. Finally, we mention that in the case of bounded coefficients, a natural choice is given by the trivial flow $\theta_{\tau,s}(\xi)=\xi$.

\subsubsection{Controls on the expansion error}

To conclude the perturbative method, we finally have to show that the expansion error in the Duhamel formula adds a negligible contribution to the Schauder estimates in \eqref{INTRO:eq:Schauder_Estimates_for_proxy} for the frozen solution $\tilde{u}^{\tau ,\xi}$. As mentioned at the beginning of this section, this control will be the most difficult to establish, due to the low regularity of the drift $F$ along the degenerate components $x_i$ ($i>1$). To briefly highlight the main additional difficulties, we will focus here on establishing a control in uniform norm  of the term in \eqref{INTRO:eq:qqq}. The other controls in uniform norm for the derivatives of the solution and those in Hölder semi-norm require even longer to be established although they still share the same approach.

We start by decomposing the remainder term with respect to the different components of the system:
\begin{equation}
\label{INTRO:eq:decomposition_error}
\Bigl{\vert}\int_{t}^{T}  \tilde{P}^{\tau,\xi}_{t,s}R^{\tau,\xi}(s,x)\, ds\Bigr{\vert}\, = \, \Bigl{\vert}
\sum_{j=1}^{n}\int_{t}^{T}\int_{\R^{nd}}\tilde{p}^{\tau,\xi}(t,s,x,y)\Delta^{\tau,\xi} F_j(s,y)\cdot D_{ y_j}u(s,y) \, dyds\Bigl{\vert},
\end{equation}
where, for simplicity, we denoted
\[
\Delta^{\tau,\xi} F_j(s,y) \, := \,F_j(s,y)- F_j(s,\theta_{\tau,s}
(\xi)), \quad j \in \llbracket 1,n\rrbracket.
\]
Recalling that the solution $u$ is assumed to be bounded differentiable in the non-degenerate component $(j=1)$, the first term in the above sum can be easily controlled exploiting the smoothing effect for the frozen density in \eqref{INTRO:eq:Smoothing_effects}:
\[
\begin{split}
\Bigl{\vert}\int_{t}^{T}\int_{\R^{nd}}\tilde{p}^{\tau,\xi}(t,s,x,y) \Delta^{\tau,\xi} &F_1(s,y)\cdot D_{ y_1}u(s,y) \, dyds\Bigl{\vert} \\
&\le \, C\Vert D_{ y_1}u\Vert_{\infty}\Vert F_1 \Vert_{C^\beta_d} \int_{t}^{T}\int_{\R^{nd}}\tilde{p}^{\tau,\xi}
(t,s,x,y)d^{\beta}\bigl(y,\theta_{\tau,s}(\xi)\bigr) \, dyds \\
&\le \, C\Vert D_{ y_1}u\Vert_{\infty} \Vert F_1 \Vert_{C^\beta_d} \int_{t}^{T} (s-t)^{\frac{\beta}{\alpha}} \, ds \\
&\le \,
C\Vert D_{ y_1}u\Vert_{\infty}\Vert F_1 \Vert_{C^\beta_d}(T-t)^{\frac{\alpha+\beta}{\alpha}}.
\end{split}
\]
In particular, we remark that the choice $(\tau,\xi)=(t,x)$ is natural if we want to balance the difference between the drifts  
\[|F_1(s,y)-F_1(s,\theta_{\tau,s}(\xi))|\]
and the anisotropic structure of the frozen density:
 \[\frac{1}{\det (\mathbb{M}_{s-t})} \bar{p}\bigl( s-t,\mathbb{M}^{-1}_{s-t} (y -
\tilde{m}^{ \tau,\xi}_{t,s}(x))\bigr)\]
and then exploit the multi-scale regularising effects associated with the frozen density in \eqref{INTRO:eq:Smoothing_effects2}. Indeed, it easy to control the dynamics in \eqref{INTRO:eq:flusso_associato} and \eqref{INTRO:eq:rappresentazione_X_tilda} 
since our choice of freezing parameters immediately implies that $\tilde{m}^{t,x}_{t,s}(x)=\theta_{t,s}(x)$.

We can now focus on the degenerate terms ($j>1$) appearing in the decomposition in Equation \eqref{INTRO:eq:decomposition_error}. Since a priori $u$ is not differentiable in $y_j$ if $j>1$ (indeed it is only $(\alpha+\beta)/(1+\alpha(j-1))$-Hölder continuous along that variable), we start by moving the derivative on the other terms through an integration by parts:
\begin{equation}\label{INTRO:besov_intro}
    \Bigl{\vert}\int_{t}^{T}\int_{\R^{nd}}D_{ y_j}\cdot\Bigl{\{}\tilde{p}^{\tau,\xi}(t,s,x,y)\Delta^{\tau,\xi} F_j(s,y)\Bigr{\}} u(s,y) \, dyds\Bigl{\vert}.
\end{equation}
The main idea will be then to bound $D_{ y_j}\cdot\bigl{\{}\tilde{p}^{\tau,\xi}(t,s,x,y)\Delta^{\tau,\xi} F_j(s,y)\bigr{\}}$ 
exploiting the regularity of the solution $u$, belonging to $L^\infty(0,T;C^{\alpha+\beta}_{b,d}(\R^{nd}))$. However, a control of this type do not seems straightforward to establish since the contribution $\tilde{p}^{\tau,\xi}(t,s,x,y)\Delta^{\tau,\xi} F_j(s,y)$ is still not differentiable in $y_j$, due to the low regularity of the drift $F$. Indeed, to obtain it, it will be necessary to consider a duality argument in Besov spaces. 

We recall that for $\tilde{\gamma}$ in $\R$, the following identification holds
\begin{equation}
\label{INTRO:eq:identific}
    C^{\tilde{\gamma}}_b(\R^d)  \, = \,  B^{\tilde{\gamma}}_{\infty,\infty}(\R^d),
\end{equation}
where for $p,q$ in $[1,\infty]$, $B^{\tilde{\gamma}}_{p,q}(\R^{nd})$ represents the Besov space on $\R^{nd}$ with indexes $(\tilde{\gamma},p,q)$. Furthermore, it is well-known (cf.\ Proposition $3.6$ in \cite{book:Lemarie-Rieusset02}) that
$B^{\tilde{\gamma}}_{\infty,\infty}(\R^d)$
and $B^{-\tilde{\gamma}}_{1,1}(\R^d)$ are dual spaces, such that
\begin{equation}\label{INTRO:Besov:duality_in_Besov}
\bigl{\vert}\int_{\R^d} fg \, dx \bigr{\vert}\, \le \, C\Vert f \Vert_{B^{\tilde{\gamma}}_{\infty,\infty}}\Vert g
\Vert_{B^{-\tilde{\gamma}}_{1,1}},
\end{equation}
for any $f$ in $B^{\tilde{\gamma}}_{\infty,\infty}(\R^d)$ and any $g$ in $B^{-\tilde{\gamma}}_{1,1}(\R^d)$.  We now recall that there are different ways to define the Besov spaces (through modulus of continuity, Littlewood-Paley decomposition, etc.) but the thermal characterisation, through convolution with a fractional heat kernel, seems to be the most natural for our purposes. For a more detailed discussion on the topic, we suggest Section $2.6.4$ in Triebel \cite{book:Triebel83}. We can then define the Besov space of indexes $(\tilde{\gamma},p,q)$ on $\R^d$ as:
\[B^{\tilde{\gamma}}_{p,q}(\R^d):= \{f \in \mathcal{S}'(\R^d)\colon \Vert f \Vert_{\mathcal{H}^{\tilde{\gamma}}_{p,q}} < + \infty\},\]
where $\mathcal{S}(\R^d)$ denotes the Schwartz space on $\R^d$. The norm $\Vert \cdot \Vert_{\mathcal{H}^{\tilde{\gamma}}_{p,q}}$ is then given by
\begin{equation}\label{INTRO:alpha-thermic_Characterization}
\Vert f \Vert_{\mathcal{H}^{\tilde{\gamma}}_{p,q}} \, := \, \Vert (\phi_0\hat{f})^\vee \Vert_{L^p}+ \Bigl(\int_{0}^{1}
v^{(1-\frac{\tilde{\gamma}}{\alpha})q}\Vert \partial_vp_h(v,\cdot)\ast f \Vert^q_{L^p} \, \frac{dv}{v}\Bigr)^{\frac{1}{q}},
\end{equation}
where $\phi_0$ is a test function in $C^\infty_c(\R^d)$ such that $\phi_0(0) \neq 0$ and $p_h$ is an isotropic, $\alpha$-stable heat kernel on $\R^d$, i.e.\ a Lévy density whose Lévy symbol is equivalent to $-\vert \lambda \vert^\alpha$. The main advantage in the thermal characterisation of Besov spaces is precisely to allow the passage of the derivative $D_{y_j}$ on the additional isotropic density $p_h(s-t,y)$ and thus, to exploit the Hölder regularity of the drift $F_j$.

We can now use the identification in \eqref{INTRO:eq:identific} and the duality property in \eqref{INTRO:Besov:duality_in_Besov} along the component $y_j$ to bound the term in \eqref{INTRO:besov_intro}: 
\begin{multline*}
\Bigl{\vert}\int_{t}^{T}\int_{\R^{nd}}D_{ y_j}\cdot\Bigl{\{}\tilde{p}^{\tau,\xi}(t,s,x,y) \Delta^{\tau,\xi}F_j(s,y)\Bigr{\}} u(s,y) \, dyds\Bigl{\vert} \\
\le \,\Vert u\Vert_{L^\infty(C^{\alpha+\beta}_{b,d})} \int_{t}^{T}\int_{\R^{(n-1)d}}\Bigl{\Vert}y_j\to D_{ y_j}\cdot\Bigl{\{}\tilde{p}^{\tau,\xi}(t,s,x,y) \Delta^{\tau,\xi}F_j (s,y)\Bigr{\}} \Bigr{\Vert}_{B^{-(\alpha_j+\beta_j)}_{1,1}} \,
dy_{\smallsetminus j}ds,
\end{multline*}
where, for simplicity, we have denoted by $y_{\smallsetminus j}$  the variable in $\R^{d(n-1)}$ without the component $y_j$. Now, it only remains to suitably control the above integral of the Besov norm. We will show in Section \ref{Sec:stabile:Second_Besov} of Chapter \ref{Chap:Schauder_estimates_Stable} the following sharp bound:
\begin{multline*}
\int_{\R^{(n-1)d}}\Bigl{\Vert}y_j\to D_{ y_j}\cdot\Bigl{\{}\tilde{p}^{\tau,\xi}(t,s,x,
y_{\smallsetminus j},\cdot) \Delta^{\tau,\xi}F_j (s,y)\Bigr{\}} \Bigr{\Vert}_{B^{-(\alpha_j+\beta_j)}_{1,1}} \, dy_{\smallsetminus j}\\
\le \, C\Vert F_j
\Vert_{L^\infty(C^{\gamma_j+\beta}_d)}(s-t)^{\frac{\beta}{\alpha}},
\end{multline*}
where, we recall, we have chosen $(\tau,\xi)=(t,x)$.

\subsubsection{Circular argument and conclusion of the proof}
Exploiting the various controls summarised above, we will finally be able to infer from the Duhamel formula \eqref{INTRO:eq:Expansion_along_proxy} the following estimate on the solution $u$ in $L^\infty(0,T ;C^{\alpha+\beta}_{b,d}(\R^{nd}))$ of the original system \eqref{INTRO:Degenerate_Stable_PDE}:
\begin{equation}\label{INTRO:eq:A_Priori_Estimates_Reg}
 \Vert u \Vert_{L^\infty(C^{\alpha+\beta}_{b,d})} \, \le \, C\bigl[\Vert u_T
\Vert_{C^{\alpha+\beta}_{b,d}}+\Vert f \Vert_{L^\infty(C^{\beta}_{b,d})}\bigr]+C\sup_{i}\Vert
F_i\Vert_{L^\infty(C^{\gamma_i+\beta}_d)}\Vert u \Vert_{L^\infty(C^{\alpha+\beta}_{b,d})},
\end{equation}
where the constant $C$ is independent from $f$, $u_T$ and $F$. In particular, we point out that when considering the associated complete norm one will lose, as in \eqref{INTRO:eq:A_Priori_Estimates_Reg}, the dependence on $(s-t)$ in short time. This in fact will not allow to directly apply a circular type argument to shift the term $\Vert u \Vert_{L^\infty(C^{\alpha+ \beta}_{b,d})}$ to the left of the above control. Instead, if one firstly suppose that the Hölder norm of the drift $F$ is small enough, i.e., for example, such that
\begin{equation}\label{INTRO:condition}
C\sup_{i}\Vert
F_i\Vert_{L^\infty(C^{\gamma_i+\beta}_d)} \, \le \, \frac{1}{2},
\end{equation}
it is then possible to use a circular argument to conclude that the Schauder estimates hold as well for $u$:
\begin{equation}\label{INTRO:eq:A_Priori_Estimates_circular}
 \Vert u \Vert_{L^\infty(C^{\alpha+\beta}_{b,d})} \, \le \, 2C\bigl[\Vert u_T
\Vert_{C^{\alpha+\beta}_{b,d}}+\Vert f \Vert_{L^\infty(C^{\beta}_{b,d})}\bigr].
\end{equation}
In the general case, it will be necessary to start the proof with an additional rescaling argument on the coefficients so that one can impose a condition like in \eqref{INTRO:condition}.

Finally, we remark that the procedure described in this section can actually be applied effectively only if the considered time interval is small enough. Intuitively, this appears natural since the expansion error to control, on which the perturbative method is based, requires that the original operator $\mathcal{L}_\alpha$ and the proxy one $\mathcal{L}^ {\tau,\xi}_\alpha$ are not too far apart. Then, to obtain the Schauder estimates for an arbitrary but finite final time, we will then have to iterate the above reasoning multiple times over each sufficiently small time interval.

\setcounter{equation}{0}

\section{Schauder estimates for a linear degenerate Lévy system}
\fancyhead[RO]{Section \thesection. Schauder estimates for linear Lévy systems}
\label{Sec:INTRO:Stime_Schauder_Levy}

We now briefly summarise the results presented in Chapter \ref{Chap:Schauder_Estimates_Levy}, which has been published in \emph{Journal of Mathematical Analysis and Applications}. Although we will focus here on a degenerate system with a linear drift, this chapter can be understood as an extension of the previous one from different points of view. Given a ``large'' space $\R^N$, we are interested in the following Ornstein-Uhlenbeck operator:
\begin{equation}
\label{INTRO:Degenerate_Levy_PDE}
  L^{\text{ou}} \, := \, \mathcal{L}+ \langle A x , D_x\rangle \quad \text{ on } \R^N,
\end{equation}
where  $\langle \cdot,\cdot \rangle$ denotes the Euclidean inner product on  $\R^N$, $A$ is a  matrix in $\R^N\otimes \R^N$ and $\mathcal{L}$ is a possibly degenerate Lévy operator, in the sense that it may acts in a non-degenerate way only on a subspace of $\R^N$. More precisely, fixed an integer $d\le N$ and a matrix $B$ in $\R^N\otimes \R^d$ such that $\text{rank}(B)=d$, the operator $\mathcal{L}$ can be represented, for any smooth enough function $\phi\colon \R^N\to \R$, as
\begin{multline}
\label{INTRO:eq:def_operator_L}
\mathcal{L}\phi(x)\, := \, \frac{1}{2}\text{Tr}\bigl(B\Sigma B^\ast D^2_x\phi(x)\bigr) +\langle Bb,D_x \phi(x)\rangle \\
+
\text{p.v.}\int_{\R^d_0}\bigl[\phi(x+Bz)-\phi(x)-\langle D_x\phi(x), Bz\rangle\mathds{1}_{B(0,1)}(z) \bigr] \,\nu(dz),
\end{multline}
where $b$ is a vector in $\R^d$, $\Sigma$ is a non-negative definite matrix in $\R^d\otimes \R^d$ and $\nu$ is a L\'evy measure on $\R^d$. We already point out that for this model we will no longer assume the \emph{symmetry} of the measure $\nu$ as we did in Chapter \ref{Chap:Schauder_estimates_Stable}. It is precisely for this reason that we cannot now delete the first order term $\langle D_x\phi(x),Bz\rangle$ in the definition of the operator $\mathcal{L}$ above (cf.\ Equation \eqref{INTRO:def_stable_operator} in the previous section).

We are interested here in establishing the well-posedness and the associated Schauder estimates for elliptic and parabolic equations involving the operator $L^{\text{ou}}$ under the \emph{minimal} Hölder regularity on the coefficients. More precisely, fixed $\lambda>0$, we will consider the following elliptic equation:
\begin{equation}\label{INTRO:eq:Elliptic_IPDE}
\lambda u(x) - L^{\text{ou}}u(x) \, = \, g(x), \quad x \in \R^N,
\end{equation}
and, for a fixed final time $T>0$, the following Cauchy problem:
\begin{equation}\label{INTRO:eq:Parabolic_IPDE}
\begin{cases}
  \partial_tu(t,x) \, = \, L^{\text{ou}}u(t,x)+f(t,x), \quad (t,x) \in (0,T)\times \R^N; \\
  u(0,x) \, = \, u_0(x), \quad x \in \R^N,
\end{cases}
\end{equation}
where $f$, $g$ and $u_0$ are given functions. Since our goal is to establish optimal regularity results, we will fix an exponent $\beta\in (0,1)$ and assume, in the elliptic problem \eqref{eq:Elliptic_IPDE}, that the source $g$ belongs to the anisotropic Hölder space $C^\beta_{b,d}(\R^N)$ while in the parabolic one \eqref{INTRO:eq:Parabolic_IPDE}, that $u_0$ is in $C^{\alpha+\beta}_ {b,d}(\R^N)$ and $f$ in $L^\infty(0,T;C^{\beta}_{b,d}(\R^N))$. A precise definition of the anisotropic Hölder spaces $C^{\gamma}_{b,d}(\R^N)$ in this context will be given later, when we will introduce the conditions on the system.

Another difference as compared to Chapter \ref{Chap:Schauder_estimates_Stable} is that for this linear model, we will only assume that $A$ and $B$ satisfy a weak Hörmander type condition, usually called the \emph{Kalman's rank} one, which ensures that the system is hypoelliptic and that for some $\alpha \in (0,2)$, the Lévy operator $\mathcal{L}$ is comparable, in a suitable sense, to a possibly truncated, non-degenerate $\alpha$-stable operator on the same subspace ($B\R^N\sim \R^d$) of $\R^N$. More in detail, the Lévy measure associated with the integral part of the operator $\mathcal{L}$ will be controlled from below by the Lévy measure of a possibly truncated $\alpha$-stable operator. Recalling from \eqref{INTRO:eq_decomposition_misura_stabile} that any $\alpha$-stable Lévy measure $\nu_\alpha$ can be decomposed into a spherical part $\mu$ and a radial one $r^{-(1+\alpha )}dr$, we will impose that the measure $\nu$ satisfies the following condition, usually called of \emph{stable domination}:
\begin{description}
\item[{[SD]}] there exist $r_0>0$, $\alpha$ in $(0,2)$ and a finite, non-degenerate (in the sense of \eqref{INTRO:def_non-deg}) measure $\mu$ on the sphere $\mathbb{S}^{d-1}$ such that
      \begin{equation}
      \label{INTRO:def_non-deg1}
          \nu(C) \, \ge \, \int_{0}^{r_0}\int_{\mathbb{S}^{d-1}} \mathds{1}_{C}(r\theta)\, \mu(d\theta)\frac{dr}{r^{1+\alpha}}, \quad
      C\in \mathcal{B}(\R^d_0).
      \end{equation}
\end{description}
Intuitively, Condition [\textbf{SD}] ensures the existence of a density associated with the operator $\mathcal{L}$ with a regularising effect of order (at least) $\alpha$. Indeed, it is precisely the small jumps (with a small radius $r_0$) associated with the Lévy measure that allow, if sufficiently intense, to generate a density for the corresponding process. Under Condition [\textbf{SD}], we know in particular that the contributions associated with the small jumps of $\nu$ are controlled from below by those of a non-degenerate $\alpha$-stable L\'evy measure, whose absolute continuity is well-known in this context. Clearly, the non-degeneracy condition [\textbf{ND}] assumed in the previous chapter can be now understood as a special case of [\textbf{SD}] when $r_0=+\infty$ and an equality holds in \eqref{INTRO:def_non-deg1}. In particular, we remark that the class of Lévy operators on $\R^d$ satisfying [\textbf{SD}] is very rich and variegated and includes some quasi-stable ones which do not previously appeared often in the literature. A possible example on $\R^2$ is given by the relativistic fractional Laplacian $\Delta^{\alpha/2}_{\text{rel}}$ acting only on the first component, which is defined by:
\[ \Delta^{\alpha/2}_{\text{rel}}\phi(x)\, := \, \text{p.v.}\int_{\R}\bigl[\phi(x_1+z ,x_2 )-\phi( x)\bigr]\frac{1+\vert z \vert^{\frac{d+\alpha -1}{2}}}{\vert z \vert^{d+\alpha} }e^{-\vert z \vert}\, dz,\]
where $x=(x_1,x_2)$ is a point in $\R^2$. Such operator, also called the relativistic Schrödinger operator, is often used for its connections with the study of the relativistic stability of matter. For more details, see for example \cite{Bogdan:Byczkowski99, Carmona:Masters:Simon90, Fefferman86, Lieb90} and the references therein.

As previously mentioned, conditions like [\textbf{SD}] ensure a regularising effect of (at least) order $\alpha$ to the operator $\mathcal{L}$ which however acts only on a subspace of $\R^N$, where the operator applies in a non-degenerate way. In order to obtain a global effect on the whole space $\R^N$, it is then necessary that this regularising effect spreads throughout the system. For this reason, we will assume that the matrices $A$, $B$ satisfy the following \emph{Kalman condition}:
\begin{description}
  \item[{[K]}] it holds that $N \, = \, \text{rank}\bigl[B,AB,\dots,A^{N-1}B\bigr]$,
\end{description}
where $\bigr[B,AB,\dots,A^{N-1}B\bigr]$ is the matrix in $\R^N\otimes\R^{dN}$ whose columns are given by $B,AB,\dots,A^{N-1}B$. It is well-known (see, for example, \cite{book:Zabczyk95}) that there exists an equivalence between Condition [\textbf{K}] and the following one:
\begin{equation}
\label{INTRO:eq_condition_K_equiv}
\det K_t \, := \, \det \int_0^te^{sA}BB^\ast e^{sA^\ast} \, ds \, > \, 0, \quad \forall t>0.
\end{equation}
At least in the diffusive case ($\alpha=2$ and $\mathcal{L}=BB^\ast\Delta_x$), $K_t$ is the covariance matrix for the solution process to the associated stochastic dynamics. In turn, Equation \eqref{INTRO:eq_condition_K_equiv} can be shown to be equivalent to the Hörmander hypoellipticity condition (cf.\ \cite{Hormander67}) for the Ornstein-Uhlenbeck operator $L^{\text{ou}}$, which ensures, in particular, the existence and regularity of a distributional solution to the following equation:
\[L^{\text{ou}}u(x) \, = \, \mathcal{L}u(x)+\langle Ax,Du(x)\rangle \, = \, \phi(x) , \quad \text{$x$ in }\R^N\]
for any sufficiently regular function $\phi\colon \R^N\to \R$. See also Ishikawa's book \cite{book:Ishikawa16}, Chapter $3.6$, for more details in the non-degenerate case. Finally, we highlight that Condition [\textbf{K}] is well-known in control theory. Indeed, it was introduced by Kalman (cf.\ \cite{Kalman60a,Kalman60b}) as a sufficient condition (actually equivalent) for the \emph{null controllability} of linear systems of the form:
\begin{equation}
\label{INTRO:controlli}
\dot{x}_t\, = \, Ax_t+Bu_t,    
\end{equation}
i.e.\ because at each final state $x$ in $\R^N$, there exists a control $t\mapsto u_t$ in $\R^d$ such that the corresponding solution $t\mapsto x_t$ to \eqref{INTRO:controlli} starting at $0$ reaches $x$ in finite time. For more details on this topic, see, for example, \cite{Kalman:Ho:Narendra63} or \cite{book:Zabczyk95}.

Thanks to Condition [\textbf{K}], we can now precisely introduce the anisotropic distance $\mathbf{d}$ and the associated Hölder spaces $C^{\beta}_{b,d}(\R^N)$ even in this context. We start by fixing $n$ as the smallest integer such that the Kalman condition [\textbf{K}] holds:
\[n \, = \, \min\left\{r \in \N \colon N=\text{rank}\left[B,AB,\dots,A^{r-1}B\right]\right\}. \]
From a more probabilistic point of view, i.e.\ considering the following stochastic dynamics:
\[dX_t \, = \, AX_tdt+BdZ_t, \quad t\ge0\]
where $\{Z_t\}_{t\ge0}$ is a Lévy process on $\R^d$ with a Lévy triplet $(b,\Sigma, \nu)$, the integer $n$ can be understood as the minimal number of applications of $A$ which allow the noise, located on $B\R^N$, to be transmitted throughout the whole space $\R^N$. More precisely, given $i$ in $\llbracket 1,n\rrbracket$, we define $V_i$ as the image space reached by the first $i-1$ iterated Lie commutators between $A$ and $B$:
\begin{equation}
\label{INTRO:eq:iterazioni_commutatori}
V_i \,  := \, \begin{cases}
            \text{Im} (B), & \mbox{if } i=1, \\
            \bigoplus_{k=1}^{i}\text{Im}(A^{k-1}B), & \mbox{otherwise}.
        \end{cases}    
\end{equation}
Since clearly $V_1\subset V_2\subset\dots V_n=\R^N$, it makes sense to denote: 
\[
W_i \,  := \, \begin{cases}
            V_1, & \mbox{if } i=1, \\
            (V_{i-1})^\perp \cap V_i, & \mbox{ptherwise}.
        \end{cases}\]
Intuitively, each space $W_i$ characterises how much a subsequent application of the commutator on $V_i$ adds in terms of the covered space. Finally, we can introduce the orthogonal projections $E_i\colon \R^N\to\R^N$ from $\R^N$ to $W_i$. Noting that $\dim E_1(\R^N)=\dim B\R^N=d$, it makes sense to fix $d_1:=d$ and write
\begin{equation}
\label{INTRO:eq_decomp}
d_i \,:= \, \dim E_i(\R^N) \quad \text{ for }i>1.    
\end{equation}
Applying a change of variables if necessary, we will assume from now on that the space $\R^N$ can be decomposed as $x=(x_1,\dots,x_n)$ such that $E_i x=x_i$ and $x_i$ is in $\R^{d_i}$, for any $i$ in $\llbracket 1, n\rrbracket$. This explicit decomposition now allows one to easily extend to this framework the anisotropic measure $\mathbf{d}$, defined as in \eqref{INTRO:def_distanza_generica}, and the Hölder spaces associated with the dilation operator 
\begin{equation}
\label{INTRO:def_operatore_dilatazione_generale}
    \delta_\lambda(t,x_1,\dots,x_n) \, = \, (\lambda^\alpha,\lambda x_1,\dots, \lambda^{1+\alpha(i-1)}x_n).
\end{equation}
In particular, a function $\phi$ in $C^{\gamma}_{d}(\R^N)$ will be such that $x_i\to \phi(x_1,\dots,x_i,\dots,x_n) $ is in $C^{\frac{\gamma}{1+\alpha(i-1)}}(\R^{d_i})$, uniformly in the other variables $x_j$, $(j\neq i)w$.

Finally, we mention that the linear structure considered here is more general than the one assumed (on the matrix $A$) in the previous chapter. In fact, it has been proved by Lanconelli and Polidoro in \cite{Lanconelli:Polidoro94} that the Kalman condition [\textbf{K}] on $A$ and $B$ imposes the following form on the matrices:
\begin{equation}\label{INTRO:eq:Lancon_Pol}
B \, = \,
    \begin{pmatrix}
        B_0    \\
        0      \\
        \vdots  \\
        0
    \end{pmatrix}
\,\, \text{ e } \,\, A \, = \,
    \begin{pmatrix}
        \ast   & \ast  & \dots  & \dots  & \ast   \\
         A_2   & \ast  & \ddots & \ddots  & \vdots   \\
        0      & A_3   & \ast  & \ddots & \vdots \\
        \vdots &\ddots & \ddots& \ddots & \ast \\
        0      & \dots & 0     & A_n    & \ast
    \end{pmatrix}
\end{equation}
where $B_0$ is a non-degenerate matrix in $\R^{d_1}\otimes \R^{d_1}$, $A_i$ is a matrix in $\R^{d_i}\otimes \R^{d_{i-1}}$ such that
$\text{rank}(A_i)=d_i$ for any $i$ in $\llbracket 2,n\rrbracket$ and the elements $\ast$ can be non-zero and different from each others. Moreover, it holds that $d_1\ge d_2\ge \dots\ge d_n\ge 1$. We also mention that the presence of non-zero elements $\ast$ in the matrix $A$ adds a further difficulty to our method of proof. In fact, it has been shown (cf.\ \cite{Lanconelli:Polidoro94}) that the matrix $A$ is invariant under the dilation $\delta_\lambda$ (defined in \eqref{INTRO:def_operatore_dilatazione_generale}) if and only if $A$ has non-zero elements only on the sub-diagonal. In particular, explicit decompositions like \eqref{INTRO:decomposition_exponential} are no longer available in the current context.

Conditions [\textbf{SD}] and [\textbf{K}] described above allow us to show the Schauder estimates associated with the Ornstein-Uhlenbeck operator $L^{\text{ou}}$. For coherence, we also mention that, similarly to the previous chapter, we will only consider \emph{weak} solutions in the sense of distributions. We can now summarise the main results obtained both in the elliptical and in the parabolic framework.

\begin{teorema}[Elliptic case]
Fixed $\lambda>0$, let $g$ be in $C^{\alpha+\beta}_{b,d}(\R^N)$. Then, there exists a unique weak solution $u$ to elliptic equation \eqref{INTRO:eq:Elliptic_IPDE}. Moreover, $u$ belongs to $C^{\beta}_{b,d}(\R^N)$ and there exists a constant $C>0$ such that
\begin{equation}
\label{INTRO:eq:Schauder_ellittico}
\Vert u \Vert_{C^{\alpha+\beta}_{b,d}} \, \le \, C\bigl(1+\frac{1}{\lambda}\bigr)\Vert g \Vert_{C^{\beta}_{b,d}}.
\end{equation}
\end{teorema}

As in the parabolic case, Schauder estimates in an elliptical context are often used in connection to the associated stochastic dynamics. For example in \cite{Cerrai:Clement07}, such estimates for degenerate diffusive operators on non-smooth domains are a fundamental tool to obtain the uniqueness of the associated martingale problem. In a non-degenerate $\alpha$-stable ($\alpha\ge1$) case, we also mention the global Schauder estimates obtained by Priola in \cite{Priola12,priola18} and their applications to the well-posedness of the associated SDE, respectively in a strong sense and in the sense of Davie (cf.\ \cite{Davie07}).

\begin{teorema}[Parabolic case]
Fixed $T>0$ and $\beta$ in $(0,1)$, let $u_0$ be in $C^{\alpha+\beta}_{b,d}(\R^N)$ and $f$ in
$L^\infty\bigl(0,T;C^{\beta}_{b,d}(\R^N)\bigr)$. Then, there exists a unique weak solution $u$ to Cauchy problem \eqref{INTRO:eq:Parabolic_IPDE}. Moreover, $u$ belongs to $L^\infty\bigl(0,T;C^{\alpha+\beta}_{b,d}(\R^N)\bigr)$ and there exists a constant $C:=C(T)>0$ such that
\begin{equation}\label{INTRO:eq:Schauder_parabolico}
\Vert u \Vert_{L^\infty(C^{\alpha+\beta}_{b,d})} \, \le \, C\bigl[\Vert u_0 \Vert_{C^{\alpha+\beta}_{b,d}}+\Vert f
\Vert_{L^\infty(C^{\beta}_{b,d})} \bigr].    
\end{equation}
\end{teorema}
Finally, we mention that the results presented here for the parabolic case can naturally be extended to Cauchy problems with time dependent coefficients $A_t$, $B_t$, under some additional assumptions, such as the boundedness of the matrix $A_t$ and the uniform ellipticity of $B_t$ on the small space $\R^d$. We refer to Section \ref{Sec:Levy:Extensions_time_dependent} of Chapter \ref{Chap:Schauder_Estimates_Levy} for more details.

\subsection{Sketch of the proof}

Inspired by Priola's work \cite{Priola09}, where analogous Schauder estimates are proved in the degenerate diffusive case, we decided to follow a semigroup approach in the proof of our result. This method, originally introduced by Da Prato and Lunardi in \cite{Daprato:Lunardi95}, consists in establishing a priori controls on the semigroup $P_t$ associated with the Ornstein-Uhlenbeck operator $L^{\text{ou}}$ on suitable functional spaces and to derive from them the parabolic Schauder estimates thanks to a mild (or Duhamel) representation of the solution $u$ (cf.\ Equation \eqref{INTRO:Duhamel_representation_of_proxy}). We remark already that this approach in the parabolic context focuses only on the regularity in space for the solutions $u$. This appears clear, for example, in the definition of the anisotropic Hölder spaces $L^\infty(0,T;C^\gamma_{b,d}(\R^N))$, only uniform in time. In particular, the parabolic Schauder estimates in \eqref{INTRO:eq:Schauder_parabolico} show no bootstrap effect in time with respect to the initial condition $u_0$. Finally, we mention that another possible method for the analysis of Ornstein-Uhlenbeck operators such as $L^{\text{ou}}$ was introduced by Manfredini in \cite{Manfredini97} and exploits a more abstract reasoning in terms of Lie groups associated with the differential dynamics. In particular, Schauder estimates in this context are constructed with respect to intrinsic Hölder spaces, i.e.\ spaces that take into account the joint regularity in space-time of the functions involved. For more details on the subject, see for example \cite{Pascucci03}.

As we have already said, the crucial element in our method consists in an a priori analysis of the properties of the Markovian semigroup $\{P_t\}_{t\ge0}$ (formally) associated with the Ornstein-Uhlenbeck operator $L ^{\text{ou}}$. Then, the Laplace transform formula will allow to write any solution $u$ to the elliptic problem \eqref{eq:Elliptic_IPDE} in terms of $P_t$:
\begin{equation}
\label{INTRO:eq:representation_elliptic_sol}
u(x) \, = \, \int_{0}^{\infty}e^{-\lambda t}\bigl[P_tg\bigr](x) \, dt \, =: \, \int_{0}^{\infty}e^{-\lambda t}P_tg(x) \, dt.
\end{equation}
In the parabolic case instead, it will be possible to use the variation of constants formula to show an analogous representation for a solution $u$ to Cauchy problem \eqref{INTRO:eq:Parabolic_IPDE}:
\begin{equation}
\label{INTRO:eq:representation_parabolic_sol}
u(t,x) \, = \, P_tu_0(x) +\int_{0}^{t}\bigl[P_{t-s}f(s,\cdot)\bigr](x)\, ds \, =: \, P_tu_0(x) +\int_{0}^{t}P_{t-s}f(s,x)\, ds.
\end{equation}
To establish the Schauder estimates for a solution $u$ (both in the elliptic and in the parabolic case), it is then clear that it is firstly necessary to obtain analogous controls on the associated semigroup $P_t$. More importantly, we are interested in understanding how the operator $P_t$ behaves on the anisotropic Hölder spaces $C^\gamma_{b,d}(\R^N)$ we are considering. We also highlight that the usual techniques for obtaining this type of estimates in the diffusive framework are, however, difficult to extend to our non-local context. For example, methods of proof as in \cite{Lunardi97}, through explicit formulas on the density of the semigroup, in \cite{Lorenzi05} or \cite{Saintier07} in the case $n=2$, through a priori Bernstein type estimates combined with interpolation methods or in \cite{Priola09}, via the Malliavin calculus for probabilistic representations of the semigroup $P_t$, can no longer be used here, mainly due to the non-local nature or the low integrability associated with the operator $\mathcal{L}$. To overcome this difficulty, we will instead exploit a \emph{perturbative} approach allowing us to consider the Lévy operator $\mathcal{L}$ as a perturbation, in a suitable sense, of an $\alpha$-stable operator whose properties are well-known. Finally, we mention that these decomposition techniques were originally introduced in \cite{Schilling:Sztonyk:Wang12} for the study of the coupling properties for Lévy processes and in \cite{Schilling:Wang12}, in relation to the generalisation of some Liouville type theorems for non-local Ornstein-Uhlenbeck operators.

\subsubsection{Smoothing effect associated with the Ornstein-Uhlenbeck operator}
As already seen in the previous section, if we want to determine the properties of the semigroup associated with the operator $L^{\text{ou}}$, it is convenient to initially consider its stochastic counterpart. Given a probability space $(\Omega,\mathcal{F},\mathbb{P})$, we then introduce the Lévy process $\{Z_t\}_{t\ge0}$ uniquely (in law) determined by the following Levy symbol:
\[\Phi(p) \, = \, ib\cdot p -\frac{1}{2}p\cdot \Sigma p + \int_{\R^d_0}\bigl(e^{ip \cdot z}-1-ip\cdot z\mathds{1}_{B(0,1)}(z)\bigr) \, \nu(dz), \quad p \in \R^d,\]
where, we recall, the triplet $(b,\Sigma,\nu)$ is the same one appearing in the definition of the operator $\mathcal{L}$ and, in particular, the Lévy measure $\nu$ satisfies the stable domination condition  [\textbf{SD}]. We also remark that the process $\{Z_t\}_{t\ge 0}$ is associated with the operator $\mathcal{L}$ in the sense that the infinitesimal generator of the degenerate process $\{BZ_t\}_{t\ge 0}$ on $\R^N$ is given by $\mathcal{L}$. Given a point $x$ in $\R^N$, we are then interested in the $N$-dimensional Ornstein-Uhlenbeck process $\{X_t\}_{t\ge 0}$ driven by $BZ_t$, which is the unique solution (in a strong sense) of the following SDE:
\[\begin{cases}
  dX^x_t \, = \, AX^x_tdt +BZ_t, \quad t>0;\\
  X^x_0 \, =\, x.
\end{cases}\]
By integrating in time the above equation against the exponential matrix $e^{(t-s)A}$, it is possible to show the following explicit representation of the solution process $\{X_t\}_{t\ge 0}$:
\[
X^x_t \, = \, e^{tA}x + \int_{0}^{t}e^{(t-s)A}B\, dZ_s, \quad \quad t\ge0.\]
The Markow semigroup associates with $\{X^x_t\}_{t\ge0}$ is now defined as the family of contraction operators
$\{P_t\colon t \ge 0\}$ on $C_b(\R^N)$, the space of all the real valued, bounded, continuous functions on $\R^N$ such that
\begin{equation}\label{INTRO:eq:Def_transition_semi-group}
P_t\phi(x) \, = \, \mathbb{E}\bigl[\phi(X^x_t)\bigr], \quad x \in \R^N.
\end{equation}
Finally, we mention that the semigroup $P_t$ is generated by the operator $L^{\text{ou}}$ in the sense that its infinitesimal generator coincides with $L^{\text{ou}}$ on the test functions space $C^\infty_c(\R^N)$.

Now, an additional reasoning on Fourier spaces, similar to that carried out in the previous section in \eqref{INTRO:prova1}-\eqref{INTRO:prova2}, allows to show that the random part of $X_t$ again satisfies the stable domination assumption [\textbf{SD}], now on the ``bigger'' space $\R^N$, even if re-scaled according to the anisotropic structure of the dynamics (cf.\ matrix $\mathbb{M}_t$ in \eqref{INTRO:eq_matrix_Mt}). Similarly to \eqref{INTRO:decomposition}, it is in fact possible to obtain that
\begin{equation}
X_t \, \overset{(\text{law})}{=} \, e^{tA}x +\mathbb{M}_t S^t_t,
\label{INTRO:decomposition_Xt_levy}
\end{equation}
where, for each fixed parameter $t$, $\{S^t_u\}_{u\ge 0}$ is a Lévy process on $\R^N$ with suitable properties. In particular, its Lévy measure $\tilde{\nu}^t$ again satisfies the stable domination condition [\textbf{SD}] extended over $\R^N$:
    \begin{equation}
       \label{INTRO:def_non-deg2}
        \tilde{\nu}^t(C) \, \ge \, \int_{0}^{R_0}\int_{\mathbb{S}^{N-1}} \mathds{1}_{C}(r\theta)\, \tilde{\mu}^t(d\theta)\frac{dr}{r^{1+\alpha}} \, =:\, \tilde{\nu}^t_\alpha(C), \quad C\in \mathcal{B}(\R^N_0),
    \end{equation}
for some $R_0>0$ and a family $\{\tilde{\mu}^t\colon t\ge0\}$ of finite, non-degenerate measures on the sphere $\mathbb{S}^{N-1}$. We now highlight that the dependence on the parameter $t$ in the process $\{S^t_u\}_{u\ge0}$ precisely appears because of the possibly non-zero elements $\ast$ in the representation \eqref{INTRO:eq:Lancon_Pol} of the matrix $A$. As already remarked before, the matrix $A$ in this case is no longer invariant under the anisotropic dilation operators $\delta_\lambda$ given in \eqref{INTRO:def_operatore_dilatazione_generale}. Indeed, decompositions like \eqref{INTRO:decomposition_exponential}, exploited in the previous section precisely for this kind of result, are now no longer valid but must be modified with ``approximate'' versions of the form:
\[e^{tA} \, = \, \mathbb{M}_tR_t\mathbb{M}^{-1}_t,\quad t \ll 1,\]
where $R_t$ is a (time dependent) locally bounded, non-degenerate matrix in $\R^N\otimes \R^N$.  We also mention that this is one of the main reasons why the a priori estimates we want to establish will only be valid in a small time interval.

From a more analytical point of view, the identity in \eqref{INTRO:decomposition_Xt_levy} suggests that the semigroup generated by the Ornstein-Uhlenbeck operator $L^{\text{ou}}$ coincides with a non-degenerate semigroup even if ``multiplied'' by the matrix $\mathbb{M}_t$ which takes into account the original degeneracy of the operator considered. More precisely, starting from the identity in law in \eqref{INTRO:decomposition_Xt_levy} we can establish a first representation of the Markov semigroup $P_t$:
\begin{equation}\label{INTRO:eq:Representation_semi-group1}
P_t\phi(x) \, = \, \int_{\R^N}\phi(e^{tA}x+\mathbb{M}_ty) \, \mathbb{P}_{S^t_t}(dy), \quad \phi \in C_b(\R^N),
\end{equation}
where, for any real random variable $X$, $\mathbb{P}_{X}$ denotes the law of $X$. 

To determine the smoothing effect associated with the semigroup $P_t$, it is now clear that a more in-depth analysis of the measure $\mathbb{P}_{S^t_t}$ is needed. First of all, we notice that Condition [\textbf{SD}] allows us to understand the Lévy process $\{S^t_u\}_{u\ge0}$ as a \emph{perturbation} of a possibly truncated, $\alpha$-stable one $\{Y^t_u\}_{t\ge 0}$, associated to the Lévy measure $\tilde{\nu}^t_\alpha$ given in \eqref{INTRO:def_non-deg2}. Properties of the measure $\mathbb{P}_{S^t_t}$ can then be inferred from the better known properties of the stable one. Indeed, we will push this perturbative method even further, considering instead of the \emph{full} stable process only those contributions associated with its small jumps which, as it is known, ensure the existence of a density for the process. This additional step will in fact allow us to obtain even more precise controls on this density and its space derivatives. More specifically, let $(\tilde{\Sigma}^t,\tilde{b}^t,\tilde{\nu}^t)$ be the Lévy triplet associated with the process $\{S^t_u\}_ {u\ge 0}$ at any fixed time $t$, which, we recall, characterises the Lévy symbol $\Phi_{S^t}$ through the Lévy-Khintchine formula on $\R^N$:
\[\Phi_{S^t}(\xi) \, := \, i\langle\tilde{b}^t,\xi\rangle -\frac{1}{2}\langle \xi, \tilde{\Sigma}^t\xi\rangle + \int_{\R^d_0}\bigl(e^{i\langle\xi, z\rangle}-1-i\langle \xi, z\rangle\mathds{1}_{B(0,1)}(z)\bigr) \, \tilde{\nu}^t(dz).\]
In a sufficiently small time interval, so that we can assume that $t^{\frac{1}{\alpha}}<R_0\wedge 1$, we truncate the Lévy measure $\tilde{\nu}^t_\alpha$ at the typical characteristic time associated with an $\alpha$-stable process at time $t$, i.e. we choose time $t^{1/\alpha}$. In particular, we now introduce the Lévy symbol $\Phi^{\text{tr}}_t$ associated with the small jumps of the $\alpha$-stable process $\{Y^t_u\}_{u\ge0}$, defined by
\[\Phi^{\text{tr}}_t(\xi) \, := \, \int_{\vert z\vert \le t^{\frac{1}{\alpha}}} \bigl[e^{i\langle \xi,z\rangle}-1-i\langle
\xi,z\rangle\bigr]\, \tilde{\nu}^t_\alpha(dz).\]
It will also be necessary to consider the remainder term which intuitively appears as the error from having considered the process $\{S^t_u\}_{u\ge0}$ as a perturbation of the $\alpha$-stable one $ \{Y^t_u\}_{u\ge0}$. More precisely, we introduce the Lévy symbol $\Phi^{\text{err}}_t$ defined by
\[\Phi^{\text{err}}_t(\xi) :=  \Phi_{S^t}(\xi)- \Phi^{\text{tr}}_t(\xi), \quad \xi \in \R^N,\] 
and associated with the Lévy triplet $(\tilde{\Sigma}^t,\tilde{b}^t,\tilde{\nu}^t-\mathds{1}_{B(0,t^{1/\alpha})}\tilde{\nu}^t_\alpha)$. We remark that it is precisely the stable domination Condition [\textbf{SD}] on $\{S^t_u\}_{u\ge 0}$ which allows us to conclude that $\Phi^{\text{err}}_t $ is actually a Lévy symbol, since it ensures the non-negativity of the measure $\tilde{\nu}^t-\mathds{1}_{B(0,t^{1/\alpha})}\tilde{\nu}^t_\alpha$ associated with the remainder term. We can now denote by $\{\mathbb{P}^{\text{tr}}_t\}_{t\ge 0}$ and $\{\pi_t\}_{t\ge 0}$ the families of infinitely divisible probabilities associated with the Lévy symbols $\Phi^{\text{tr}}_t$ and $\Phi^{\text{err}}_t$. From the representation through characteristic functions of the probability measures, we can then disintegrate the probability $\mathbb{P}_{S^t_t}$ in the following way:
\begin{equation}
\label{INTRO:convolution_measures}
 \mathbb{P}_{S^t_t} \, = \, \mathbb{P}^{\text{tr}}_t \ast \pi_t,\quad t>0,
\end{equation}
where $\ast$ represents the convolution between probability measures.

From the above arguments, we can now focus on the family of measures $\{\mathbb{P}^{\text{tr}}_t\}_{ t\ge 0}$ which, we recall, are associated with the small jumps of an $\alpha$-stable process. Exploiting known results on the Lévy symbol associated with a non-degenerate $\alpha$-stable process, such as the Hartman-Wintner condition and some controllability assumptions on Fourier spaces, we can obtain the existence of a regular (in space) density for the probability measure $\mathbb{P}^{\text{tr}}_t$ and suitable controls on its derivatives, at least in a sufficiently small time interval. More precisely, we will show that there exists a final time $T_0:=T_0(N)>0$ such that on $(0,T_0]$, the probability $\mathbb{P}_t^{\text{tr}}$ admits a density $p^{\text{tr}}(t,\cdot)$ which is $3$-times differentiable with continuous derivative on $\R^N$. Moreover, for each $k$ in $\llbracket 0 , 3\rrbracket$, it holds that
\begin{equation}
\label{INTRO:smoothing_effect_Levy}
\vert \partial^k_y p^{\text{tr}}(t,y) \vert \, \le \, Ct^{-\frac{N+k}{\alpha}} \Bigl(1+\frac{\vert y
\vert}{t^{1/\alpha}}\Bigr)^{-(N+3)} \, =:\, Ct^{-\frac{k}{\alpha}}\bar{p}(t,y),   
\end{equation}
for any $(t,y)$ in $(0,T_0]\times\R^N$, where $C>0$ is a constant depending only on $N$. We also remark that, at the cost of further reducing the time interval, it is possible to show that the density $p^{\text{tr}}(t,\cdot)$ is even more regular in space and the associated density $\bar{p}(t,\cdot)$ presents stronger smoothing properties but that the choice in \eqref{INTRO:smoothing_effect_Levy} is the one adapted to our purposes. In particular, the exponent $N+3$ in the density $\bar{p}(t,\cdot)$ is the minimal necessary to be able to integrate the contributions coming from the functions $\phi$ in $C^{\beta}_{b,d}(\R^N)$ with index $\beta< 1+\alpha<3$ (cf.\ Equation \eqref{INTRO:Control_Cgamma_density} below). Thanks to Equation \eqref{INTRO:convolution_measures}, we can now rewrite the representation for the Markov semigroup $P_t$ in \eqref{INTRO:eq:Representation_semi-group1} in the following way:
\begin{multline}\label{INTRO:eq:Representation_semi-group2}
P_t\phi(x) \, = \, 
\int_{\R^N}\int_{\R^N}\phi\bigl(e^{tA}x+\mathbb{M}_t(y_1+y_2)\bigr)p^{\text{tr}}(t,y_1) \, dy_1 \pi_t(dy_2) \\
= \, \int_{\R^N}\int_{\R^N}\phi\bigl(y_1+\mathbb{M}_ty_2\bigr)\frac{p^{\text{tr}}(t,\mathbb{M}^{-1}_t(y_1-e^{tA}x))}{\det \mathbb{M}_t} \, dy_1 \pi_t(dy_2).
\end{multline}
Furthermore, the explicit form of the density $\bar{p}(t,\cdot)$ and the corresponding controls in \eqref{INTRO:smoothing_effect_Levy} allow one to easily show the smoothing effects associated with the ``partial'' density $p^{\text{tr}}(t,\cdot)$, even if only along the component $y_1$. In particular, given $\gamma$ in $(0,3)$, $k$ in $\llbracket 0,2\rrbracket$, $i$ in $\llbracket 1,n \rrbracket$ and $t$ small enough, it holds that
\begin{equation}
\label{INTRO:eq:Smoothing_effects_Levy}
\int_{\R^{nd}} \frac{|D^k_{x_i}p^{\text{tr}}(t,\mathbb{M}^{-1}_t(y_1-e^{tA}x))|}{\det \mathbb{M}_t} \mathbf{d}^\gamma(y_1,e^{tA}x)\, dy_1 \, \le \, C(s-t)^{\frac{\gamma}{\alpha}-k\frac{1+\alpha(i-1)}{\alpha}}. 
\end{equation}

Finally, it is clear that the main difference with the techniques presented in Section \ref{Sec:INTRO:Stime_Schauder_stabili} (for the Schauder estimates of the proxy operator) consists in the fact that here we can only exploit the regularising effects on the variable $y_1$, associated with the small jumps contributions of the truncated $\alpha$-stable process, given that we do not know any particular properties of the remainder term (associated with the integration variable $y_2$). In fact, we only know that the probability $\pi_t$ has a finite total measure, uniformly in $t$. While this peculiarity certainly introduces further additional difficulties in some steps of the proof, for example on the cancellation techniques for the a priori Hölder controls, we highlight that exploiting only the regularising effects associated with small jumps allows us to drop here some of the assumptions (cf.\ hypotheses [\textbf{P}]) on the parameters $\alpha$, $\beta$ needed in Chapter \ref{Chap:Schauder_estimates_Stable} precisely because we considered the density associated with the full $\alpha$-stable process.

\subsubsection{Controls on the associated Markov semigroup}
The regularising effects of the density $p^{\text{tr}}(t,\cdot)$ shown in \eqref{INTRO:eq:Smoothing_effects_Levy} naturally imply the first controls on $P_t$ and its derivatives when the semigroup acts on the space $C_b(\R^N)$. In more detail, we have that
\begin{equation}
\label{INTRO:eq_control_infty}
\Vert D^k_{x_i} P_t \phi \Vert_\infty \, \le \, C\Vert \phi \Vert_\infty \bigl(1+t^{-k\frac{1+\alpha(i-1)}{\alpha}}\bigr), \quad t>0,
\end{equation}
where $k$ is in $\llbracket 0,3\rrbracket$, $i$ in $\llbracket 1,n\rrbracket$ and $\phi$ is a function in $C_b(\R^N)$. In particular, the above estimates implies the continuity of the semigroup $P_t$, at any fixed time $t>0$, as an operator on the space $C_b(\R^N)$ with values in the same or in the anisotropic Hölder-Zygmund space $C^3_{b,d}(\R^N)$, a natural generalisation to integer indexes of the anisotropic Hölder spaces we defined. This result will then be extended to general anisotropic Hölder spaces by noticing that any space $C^\gamma_{b,d}(\R^N)$ (with $\gamma$ in $(0,3)$) can be seen as an interpolation space between $C_b(\R^N)$ and $C^3_{b,d}(\R^N)$. In fact, we will exploit the following identity:
\begin{equation}
\label{INTRO:interpolazione}
\bigl(C_b(\R^N),C^3_{b,d}(\R^N)\bigr)_{\gamma,\infty} \, = \, C^{\gamma}_{b,d}(\R^N), 
\end{equation}
with the equivalence between the respective norms, where for any two Banach spaces $X,Y$, the symbol $(X,Y)_{\gamma,\infty}$ represents the \emph{real interpolation} space with infinite norm between $X$ and $Y$. For more details on the subject, we refer to Triebel's books \cite{book:Triebel83} and Lunardi's one \cite{book:lunardi18} or to \cite[Theorem $2.2$]{Lunardi97} in an anisotropic diffusive context. Thanks to interpolation techniques similar to those in \eqref{INTRO:interpolazione}, we will then show that the semigroup $P_t$ is continuous as an operator on $C_b(\R^N)$ with values in the anisotropic Hölder space $C^\gamma_{b,d}(\R^N)$ for each $\gamma$ in $[0,1+\alpha)$:
\begin{equation}\label{INTRO:eq:coroll:C0_Cgamma_control}
\Vert P_t \Vert_{\mathcal{L}(C_b,C^\gamma_{b,d})} \, \le \, C\bigl(1+t^{-\frac{\gamma}{\alpha}}\bigr), \quad t>0,
\end{equation}
where $\Vert \cdot \Vert_{\mathcal{L}(X,Y)}$ denotes the usual operator norm between two Banach spaces $X$ and $Y$.

For now, we have only considered the semigroup behaviour on the space $C_b(\R^N)$. To extend the analysis to $C^{\beta}_{b,d}(\R^N)$ as well, the first step is to obtain controls similar to those in \eqref{INTRO:eq_control_infty} when $P_t$ acts instead on Hölder continuous functions $\phi$. Similarly to the previous chapter, we will now use cancellation techniques to exploit the Hölder regularity of the functions $\phi$. However, we highlight that the main difficulty in this case will be linked to the fact that the smoothing effects associated with the density $p^\text{tr}$ are only available with respect to the ``good'' variable $y_1$. We will then have to introduce \emph{partial} cancellation techniques which in fact allow us to isolate only the components of $\phi$ along $y_1$. In more detail and focusing only on the case $\beta<1$, $k=1$ for simplicity, the fundamental idea will be to exploit, similarly to \eqref{INTRO:eq:cancellazione1}, that
\[\int_{\R^N}\phi(\mathbb{M}_ty_2+e^{tAx})D_{x_i}\int_{\R^N}\frac{p^{\text{tr}}(t,\mathbb{M}^{-1}_t(y_1-e^{tA}x))}{\det \mathbb{M}_t} \, dy_1 \pi_t(dy_2) \, = \, 0\]
to add the term $\phi(\mathbb{M}_ty_2+e^{tAx})$ in the control of $|D_{x_i}P_t|$ and conclude as in \eqref{INTRO:eq:cancellazione1}, thanks to the smoothing effect associated with the density $p^{\text{tr}}(t,\cdot)$ given in \eqref{INTRO:eq:Smoothing_effects_Levy}. Arguments similar to the one mentioned now will allow in particular to show that for any $\beta$ in $[0,3)$, $i$ in $\llbracket 1, n \rrbracket$ and $k$ in $\llbracket0,3 \rrbracket$, it holds that
\begin{equation}
\label{INTRO:Control_Cgamma_density}
   \Vert D^k_{x_i} P_t \phi \Vert_\infty \, \le \, C\Vert \phi \Vert_{C^\beta_{b,d}} \bigl(1+t^{\frac{\beta}{\alpha}-k\frac{1+\alpha(i-1)}{\alpha}}\bigr), \quad t>0.
\end{equation}
Exploiting again interpolation techniques as in \eqref{INTRO:interpolazione} starting from the controls on the derivatives of $P_t\phi$ in \eqref{INTRO:Control_Cgamma_density}, we will finally be able to show the continuity of the semigroup $P_t$ as an operator between anisotropic Hölder spaces. More precisely, we will obtain that for any $\beta<\gamma$ in $[0,1+\alpha)$, the following estimate hold:
\begin{equation}\label{INTRO:eq:Intro:Gradient_Estimates}
\Vert P_t \Vert_{\mathcal{L}(C^\beta_{b,d},C^\gamma_{b,d})} \, \le \,C\Bigl(1+t^{\frac{\beta-\gamma}{\alpha}}\
\Bigr), \quad t>0.
\end{equation}
To the best of our knowledge, such estimates appear to be new in the degenerate Lévy context and of independent interest. We also point out that the possibility here of a smoothing effect that is independent from the index $\alpha$ (unlike in the previous section) essentially reflects the fact that the density is in this case associated only with the small jumps contribution, while we know that they are the tails, corresponding to to the process big jumps, that impose the integrability conditions seen previously.

\subsubsection{Schauder estimates in the elliptic case}
Once the necessary a priori controls on the semigroup $P_t$ have been shown, the Schauder estimates for the solution $u$ in the elliptic context \eqref{INTRO:eq:Schauder_ellittico} and in the parabolic one \eqref{INTRO:eq:Schauder_parabolico} can be obtained from the relative representations of $u$ in terms of the Markov semigroup $P_t$ (cf.\ Equations \eqref{INTRO:eq:representation_elliptic_sol} and \eqref{INTRO:eq:representation_parabolic_sol}). To give the reader an idea of the method we have followed, we will now briefly present the elliptic case.

Given a solution $u$ to the elliptic problem in \eqref{INTRO:eq:Elliptic_IPDE}, we know from the Laplace formula in \eqref{INTRO:eq:representation_elliptic_sol} that is then necessary to control the following term:
\[u(t,x) \, = \,\int_{0}^{\infty}e^{-\lambda t}\bigl(P_tg\bigr)(z) \, dt, \]
in the anisotropic Hölder space of order $\alpha+\beta$ with a source $g$ in $C^{\beta}_{b,d}(\R^N)$. For this kind of problem, it is actually convenient to exploit an equivalent norm to the Hölder one introduced in \eqref{INTRO:def_holder_seminorm} which does not require considering the derivatives along each direction but only the finite differences of order $3$ for the function $u$. More precisely, we introduce for an initial point $x_0$ in $\R^N$ and $z$ in $E_i(\R^N)$,
\[
\Delta^3_{x_0}\phi(z) \, := \, \phi(x_0+3z)-3\phi(x_0+2z)+ 3\phi(x_0+z)-\phi(x_0).
\]
It has been shown in \cite{Lunardi97} that a function $\phi$ is in $C^\gamma_{b}(E_i(\R^N))$ if and only if
\[\sup_{x_0 \in \R^N}\sup_{z \in E_i(\R^N);z\neq 0}\frac{\bigl{\vert} \Delta^3_{x_0}\phi(z)\bigr{\vert}}{\vert z \vert^\gamma} \, < \,
\infty.\]
Finally, it is then necessary to show that for each $i$ fixed in $\llbracket 1,n\rrbracket$, it holds that
\begin{equation}
\label{INTRO:eq:finale_Schauder_Levy}
\vert \Delta^3_{x_0}u(z)\vert \, = \, \Bigl{\vert}\int_{0}^{\infty}e^{-\lambda t}\Delta^3_{x_0}\bigl(P_tg\bigr)(z) \, dt\Bigr{\vert} \,
\le \, C\Vert g \Vert_{C^\beta_{b,d}}\vert z\vert^{\frac{\alpha+\beta}{1+\alpha(i-1)}},
\end{equation}
for some constant $C>0$ independent from $x_0$ in $\R^N$ and from $z$ in $E_i(\R^N)$.

As often happens in the proof of Hölder norm controls, we will firstly divide the analysis according to three possible regimes, dictated by the relationship between the spacial point $z$ in $E_i(\R^N)$ and the time $t$ at the intrinsic scale of the system along the $i$-th direction considered. On the one hand, the \emph{macroscopic} regime will appear when $|z|\ge 1$ and will be the easiest to handle. Indeed, it will only require showing the boundedness of the solution $u$ starting from that of $g$, since in this case, $\Vert g \Vert_\infty \le \Vert g \Vert_\infty|z|^{ \gamma}$ for each $\gamma>0$. On the other hand, we will say that we are in an \emph{off-diagonal} regime if $t^{\frac{1+\alpha(i-1)}{\alpha}} \le \vert z \vert \le $1. In this case, the spacial distance along the $i$-th component will be bigger than the associated characteristic time. Finally, a \emph{diagonal} regime will appear when $t^{\frac{1+\alpha(h-1)}{\alpha}} \ge \vert z \vert$ and the space point will instead be smaller than the typical intensity of the characteristic time. In particular, we will denote by $t_0$ the transition time between the diagonal and off-diagonal regimes, with respect to the dilation operators $\delta_\lambda$ (defined in \eqref{INTRO:def_operatore_dilatazione_generale}) along the $i$-th direction of the system, that is:
\[
t_0 \, := \, |z|^{\frac{\alpha}{1+\alpha(i-1)}}.
\]

As already mentioned before, the control in the macroscopic regime $(\vert z \vert\ge 1)$ immediately follows from the contraction property of the semigroup $P_t$ on $C_b(\R^N)$:
\begin{equation}
\label{INTRO:control_macro_schauder}
\Bigl{\vert}\int_{0}^{\infty}e^{-\lambda t}\Delta^3_{x_0}\bigl(P_tg\bigr)(z) \, dt\Bigr{\vert} \, \le \, 3\int_{0}^{\infty}e^{-\lambda t}\Vert P_tg\Vert_\infty\, dt \, \le \, C\Vert g \Vert_\infty\vert z \vert^{\frac{\alpha+\beta}{1+\alpha(i-1)}}.
\end{equation}
To study separately the diagonal and off-diagonal regimes, we will split, as in \cite{Priola09}, the term $\Delta^3_{x_0}u(z)$ into two components $R_1(z)+R_2(z) $, where
\begin{align*}
  R_1(z) \, &:= \, \int_{0}^{t_0}e^{-\lambda t}\Delta^3_{x_0}\bigl(P_tg\bigr)(z) \, dt; \\
  R_2(z) \, &:= \, \int_{t_0}^{\infty}e^{-\lambda t}\Delta^3_{x_0}\bigl(P_tg\bigr)(z) \, dt.
\end{align*}
In the off-diagonal case associated with the first component $R_1$, we will eventually have to integrate over $[0,t_0]$ and it will therefore be important to pay attention not to introduce singularities in time that are not integrable. In practice, we will exploit the continuity of the semigroup $P_t$ on the space $C^\beta_{b,d}(\R^N)$ to control the term $R_1$ as follows:
\begin{equation}
\label{INTRO:control_off-diag_schauder}
    \begin{split}
\vert R_1(z)\vert \, &\le \, \int_{0}^{t_0}\vert\Delta^3_{x_0}\bigl(P_tg\bigr)(z)\vert \, dt \\
&\le \, \vert z \vert^{\frac{\beta}{1+\alpha(h-1)}}\int_{0}^{t_0}\Vert P_tg\Vert_{C^{\beta}_{b,d}}\, dt \\
&\le \, C\Vert g \Vert_{C^\beta_{b,d}}\vert z \vert^{\frac{\alpha+\beta}{1+\alpha(h-1)}}.
    \end{split}
\end{equation}
On the other hand, the diagonal regime associated with the contribution $R_2$ will not present problems of integrability in time. Indeed, since the point $z$ is, in this case, closer to zero than the characteristic time scales, it makes sense to iteratively apply a Taylor expansion on $\Delta^3_{x_0}\bigl(P_t\phi\bigr)$ so that a third order derivative appears along the $i$-th component. More precisely,
\[\begin{split}
\bigl{\vert}\Delta^3_{x_0}\bigl(&P_tg\bigr)(z)\bigr{\vert}\\
&= \, \Bigl{\vert}\int_{0}^{1} \langle
D_{x_i}P_t g(x_0+\lambda z)-2D_{x_i}P_t g(x_0+z+\lambda z)+D_{x_i}P_t g(x_0+2z+\lambda z), z \rangle \,d\lambda\Bigr{\vert} \\
&\le \, \Bigl{\vert}\int_{0}^{1}\int_{0}^{1} \langle \bigl[D^2_{x_i}P_tg(x_0+(\lambda+\mu)z)-D^2_{x_i}P_tg(x_0+z+(\lambda+\mu)z)\bigr]z,
z \rangle \,d\lambda d\mu\Bigr{\vert} \\
&\le \, \Bigl{\vert}\int_{0}^{1}\int_{0}^{1}\int_{0}^{1} \langle \bigl[D^3_{x_i}P_tg(x_0+(\lambda+\mu+\nu)z)\bigr](z,z),z \rangle
\,d\lambda d\mu d\nu\Bigr{\vert},
\end{split}
\]
where, for simplicity, we have identified the differential operator $D^3_{x_i}P_t\phi$ with the associated $3$-tensor. The a priori controls in \eqref{INTRO:Control_Cgamma_density} on the derivatives of the semigroup $P_t$ now imply that
\[
\bigl{\vert}\Delta^3_{x_0}\bigl(P_tg\bigr)(z)\bigr{\vert}\, \le \, C\Vert D^3_{x_i}P_tg\Vert_\infty \vert z \vert^3 \, \le \, C\Vert g \Vert_{C^{\gamma}_{b,d}}\bigl(1+t^{\frac{\gamma-3(1+\alpha(i-1))}{\alpha}}\bigr) \vert z \vert^3.
\]
We can thus conclude that in the diagonal case the following estimates hold:
\begin{equation}
\label{INTRO:control_diag_schauder}
\begin{split}
\vert R_2(z) \vert \, &\le \, \int_{t_0}^{\infty}e^{-\lambda t}\vert\Delta^3_{x_0}\bigl(P_tg\bigr)(z)\vert \, dt \\
&\le \, C\Vert g \Vert_{C^\beta_{b,d}}\vert z \vert^3 \int_{t_0}^{\infty}e^{-\lambda t}\bigl(1+t^{\frac{\beta-3(1+\alpha(i-1))}{\alpha}}\bigr)\, dt \\
&\le C\Vert g \Vert_{C^\beta_{b,d}}\vert z \vert^3\bigl(\lambda^{-1}+\vert z\vert^{\frac{\alpha+\beta-3(1+\alpha(i-1))}{1+\alpha(i-1)}}\bigr) \\
&\le C\Vert g \Vert_{C^\beta_{b,d}}\vert z \vert^{\frac{\alpha+\beta}{1+\alpha(i-1)}},
\end{split}
\end{equation}
where, in the last step, we exploited again that $\vert z \vert \le 1$. Finally, by bringing together the contributions \eqref{INTRO:control_macro_schauder}, \eqref{INTRO:control_off-diag_schauder} and \eqref{INTRO:control_diag_schauder} associated with the three possible regimes, we can conclude that Equation \eqref{INTRO:eq:finale_Schauder_Levy} holds and that in particular, Schauder estimates \eqref{INTRO:eq:Schauder_ellittico} in the elliptic case hold true.

Schauder estimates in the parabolic case will be obtained following a similar procedure. We also remark that this type of decompositions into diagonal and off-diagonal regimes also appears, even if in a less explicit form, in the Hölder norm controls for the Schauder estimates in Chapter \ref{Chap:Schauder_estimates_Stable}. Finally, we mention that we could have used the method exposed here also in Chapter \ref{Chap:Schauder_estimates_Stable} in order to show the Schauder estimates \eqref{INTRO:eq:Schauder_Estimates_for_proxy} for the proxy solution $\tilde{u}^{\tau, \xi}$, under the more general stable domination condition [\textbf{SD}]. However, to obtain the estimates in \eqref{INTRO:eq:Schauder_estimates} for the solution $u$ to the original non-linear system \eqref{INTRO:Degenerate_Stable_PDE} using the perturbation method described above, the tricky part would have been to prove, under the more general hypotheses, the actual independence of the constant from the frozen parameters $(\tau,\xi)$. This difficulty will appear even more clearly in the model assumptions in the next section, when we will additionally consider a multiplicative noise.

\setcounter{equation}{0}

\section{Weak well-posedness for degenerate SDEs driven by Lévy processes}
\fancyhead[RO]{Section \thesection. Weak well-posedness for Lévy degenerate SDEs}
\label{Sec:INTRO:Buona_posizione_Levy}

We briefly summarise here the main elements of Chapter \ref{Chap:Weak_Well-Posedness} of the current thesis. Written in collaboration with my PhD supervisor, Prof.\ Stéphane Menozzi, this work has appeared as a pre-print \cite{Marino:Menozzi21}. We want to study here the effect of the propagation of a non-degenerate Lévy noise through a chain of interconnected oscillators, where the aforementioned noise only acts on the first of them, as shown in Figure \ref{fig:my_label}. More precisely, we are interested in a stochastic dynamics on $\R^N$ of the form:
\begin{equation}\label{INTRO:eq:SDE}
\begin{cases}
  dX_s \, = \, G(s,X_s)ds + B\sigma(s,X_{s-})dZ_s, \quad s\ge t,\\
  X_t\, = \, x,
\end{cases}
\end{equation}
for an initial point $(t,x)$ in $[0,\infty)\times \R^N$, where the deterministic drift $G\colon [0,\infty)\times \R^N\to\R^N$ and the diffusion matrix $\sigma \colon [0,\infty)\times \R^N\to\R^d\otimes \R^d$ are given. From now on, we will say that Equation \eqref{INTRO:eq:SDE} is \emph{degenerate} in the sense that the Lévy noise $\{Z_s\}_{s\ge 0 }$ initially acts only on the ``small'' space $\R^d$. The non-degenerate nature of the matrix $\sigma$ then ensures that the noise associated with the stochastic integral $\int_0^t\sigma(s,X_{s-})dZ_s$ is transmitted through the whole space $\R^d$. Finally, the latter acts on the ``large'' space $\R^N$ through the embedding matrix $B$ in $\R^N\otimes \R^d$, given by:
\[B:=\begin{pmatrix}
            I_{d\times d}, & 0_{d\times (N-d)}
\end{pmatrix}^t.\]
In particular, we will assume that $\{Z_s\}_{s\ge 0}$ is a pure jump process, i.e.\ that $\Sigma=0$ in the associated Lévy triplet $(b,\Sigma,\nu)$. To highlight the effective dependence on the given initial point $(t,x)$, from now on we will denote by $\{X^{t,x}_s\}_{s\ge0}$ a generic solution process of the stochastic dynamics in \eqref{INTRO:eq:SDE} extended up to time zero, i.e.\ imposing $X^{t,x}_s=x$ if $s$ in $[0,t]$.

In \cite{Chaudru:Menozzi17}, the authors considered a degenerate chain as in \eqref{INTRO:eq:SDE} perturbed by a Brownian motion $\{Z_t\}_{t\ge0}$ and they then characterised the minimal Hölder regularity on the drift $G$ ensuring its well-posedness in a weak sense. The initial goal of this work was indeed to extend this result to SDEs like \eqref{INTRO:eq:SDE} whose noise was only Lévy, under the same stable domination assumption [\textbf{SD}] presented in the previous section (cf.\ Equation \eqref{INTRO:def_non-deg1}). We have actually managed to obtain only a partial generalisation of the results in \cite{Chaudru:Menozzi17}. Indeed, we will show later the existence of an optimal threshold for the Hölder regularity of the drift $G$, but only for a particular type of diagonal deterministic structure. Moreover, the method of proof we followed, through a perturbative approach by backward expansions, required to strengthen the stable domination condition [\textbf{SD}] and to add some further assumptions. In this section, we briefly summarise the reasons that led us to these considerations.

The deterministic structure of  the stochastic dynamics (cf.\ \eqref{INTRO:eq:SDE} for $\sigma=0$) will be similar to the one presented in the first two chapters. In particular, we will suppose again that we can decompose the ``large'' space $\R^N$ into $n$ subspaces $\R^{d_i}$, $i\in \llbracket 1, n\rrbracket$ such that $d_1=d$ and $d_1+\dots+d_n=N$, as shown in Section \ref{Sec:INTRO:Stime_Schauder_Levy}. Moreover, the deterministic drift $G$ will present, with respect to this decomposition, a particular ``upper triangular'' structure and its elements on the subdiagonal will be considered non-degenerate and linear. In practice, we will impose that $G$ has the following form:
\begin{equation}
\label{INTRO:decomposizione_deriva_G}
G(s,x)\,:= \, A_sx+F(s,x), 
\end{equation}
where $A\colon[0,\infty)\to \R^N\otimes\R^N$ and $F\colon [0,\infty)\times \R^N\to \R^N$ are two functions such that
\begin{description}
  \item[{[H]}]
  \begin{itemize}
      \item for any level $i$ in $\llbracket 1,n\rrbracket$, $F_i$ only depends on time and on the last $n-(i-1)$ spacial variables, i.e.\ $F_i(s,x_i,\dots,x_n)$;
      \item $A\colon [0,\infty)\to \R^N\otimes \R^N$ is bounded and 
      \[[A_s]^{i,j} \, = \, 
      \begin{cases}
      \text{non-degenerate, uniformly in } s, &\mbox{ if } j= i-1;\\
      0, &\mbox{ if } j< i-1.
      \end{cases}\]
  \end{itemize}
\end{description}
As already mentioned above, this assumption in the linear case with additive noise (i.e.\ $F=0$ and $\sigma=1$) can be understood as an Hörmander type condition, or equivalently a Kalman rank one (cf.\ condition [\textbf{K}] in Section \ref{Sec:INTRO:Stime_Schauder_Levy}), ensuring the hypoellipticity of the infinitesimal generator for the process $\{X^{t,x}_s\}_{s\ge0}$, solution to Equation \eqref{INTRO:eq:SDE}.

Under Condition [\textbf{H}], we then notice that the matrix $A$ can be rewritten as a ``time dependent'' version of the one that appeared in Section \ref{Sec:INTRO:Stime_Schauder_Levy}, Equation \eqref{INTRO:eq:Lancon_Pol}. Furthermore, this condition allows to represent the stochastic dynamics in \eqref{INTRO:eq:SDE} in the following, more explicit, form:
\[
  \begin{cases}\vspace{2pt}
    dX^1_t \, = \,\left[A^{1,1}_t X_t^1+\cdots+A^{1,n}_tX_t^n+ F_1(t,X^1_t,\dots,X^n_t)\right]dt + \sigma(t,X^1_{t-},\dots,X^n_{t-})dZ_t,\\ \vspace{2pt}
    dX^2_t \, = \, \left[ A^{2,1}_tX^1_t+\cdots+A^{2,n}_tX^n_t +F_2(t,X^2_t,\dots,X^n_t)\right] dt,\\
    dX^3_t \, = \, \left[ A^{3,2}_tX^2_t+\cdots+A^{3,n}_tX^n_t+F_3(t,X^3_t,\dots,X^n_t)\right]dt,\\
    \vdots\\
    dX^n_t \, = \, \left[ A^{n,n-1}_t X^{n-1}_t+A^{n,n}_t X^{n}_t+F_n(t,X^n_t)\right]dt,
  \end{cases}
\]
where we have decomposed $X_t=(X^1_t,\dots,X^n_t)$ such that each $X^i_t$ belongs to $\R^{d_i}$.

Since the main aim of this chapter is to determine the minimal Hölder regularity on the drift $F$ which guarantees the weak well-posedness of the stochastic dynamics and such optimal thresholds will involve only the degenerate components ($i>1$) of $F$, we can now state separately the conditions on the other coefficients. In particular, we will assume that:
\begin{description}
\item[{[R]}] there exist an index $\beta^1$ in $(0,1)$ and a constant $K>0$ such that
\begin{itemize}
    \item $\sigma(t,\cdot)$ is $\beta^1$-Hölder continuous, uniformly in $t$;
    \item $F_1(t,x)$ is $\beta^1$-Hölder continuous, uniformly in $t$;
    \item for any $i$ in $\llbracket 1,n\rrbracket$, it holds that
    \[|F_i(t,0)|\, \le \, K, \quad t \in [0,T].\]
\end{itemize}
\end{description}

To expect a minimal regularising effect (of order $\alpha$) to be carried by the driving noise, we will assume that $\{Z_t\}_{t\ge0}$ can be associated with an $\alpha$-stable process of the ``tempered'' type. However, we remark that this class does not include only the classical tempered $\alpha$-stable processes. More in detail, we will supposed that the Lévy measure $\nu$ of the driving process $\{Z_s\}_{s\ge 0}$ is \emph{symmetric} and
\begin{description}
  \item[{[ND']}] there exists a Borel measurable function $Q\colon\R^d\to \R$ such that
  \begin{itemize}
      \item $Q$ is positive and bounded, i.e.\ $Q\ge 0$ and $\sup_{z\in \R^d}Q(z)<\infty$;
      \item $Q$ is bounded away from zero and Lipschitz regular in a neighbourhood of the origin, i.e.\ there exist $r_0>0$ and $c>0$ such that $Q(z)\ge c$ and $Q$ is Lipschitz continuous in $B(0,r_0)$;
      \item there exist $\alpha\in (1,2) $ and a finite, non-degenerate (cf.\ \eqref{INTRO:def_non-deg}) measure $\mu$ on $\mathbb{S}^{d-1}$ such that
       \[
\nu(\mathcal{A}) \,=\,\int_{0}^{\infty}\int_{\mathbb{S}^{d-1}}\mathds{1}_{\mathcal{A}}(\rho s)Q(\rho s)\,\mu(ds) \frac{d\rho}{\rho^{1+\alpha}},\quad \mathcal{A} \in\mathcal{B}(\R^d_0).
\]
\end{itemize}
\end{description}
Bearing in mind the classical decomposition of an $\alpha$-stable Lévy measure given, for example, in \eqref{INTRO:eq_decomposition_misura_stabile}, the \emph{non-degeneration} condition [\textbf{ND'}] intuitively imposes that the Lévy measure of $\{Z_t\}_{t\ge0}$ is absolutely continuous with respect to that of a non-degenerate $\alpha$-stable process and that their Radon-Nikodym derivative is given by a ``tempering'' function $Q$ with suitable properties. Clearly, a natural possibility is given by the constant function $Q=c$, that is, when the process $\{Z_t\}_{t\ge0}$ is a symmetric $\alpha$-stable process. However, we point out that the class of noises considered here also covers several symmetric \emph{quasi-stable} processes including, for example, the relativistic $\alpha$-stable process or the Lamperti one (see Chapter \ref{Chap:Weak_Well-Posedness}, Section \ref{Sec:Well-Pos:Complete_model_assumptions} for more details).

When considering SDEs driven by multiplicative noise, a natural assumption is the \emph{uniform ellipticity} of the non-degenerate component associated with the diffusion matrix $\sigma(t,x)$, at each fixed point of the space-time. Indeed, we will assume that:
\begin{description}
  \item[{[UE]}] there exists a constant $\eta>1$ such that for any $s\ge 0$ and any $x$ in $\R^N$, it holds that
  \[\eta^{-1}\vert \xi \vert^2 \, \le \, \sigma(s,x)\xi\cdot \xi  \, \le\,\eta
\vert \xi \vert^2, \quad \xi \in \R^d,
\]
\end{description}
where we recall that ``$\cdot$'' denotes the inner product on the ``small'' space $\R^d$. Intuitively, the uniform ellipticity assumption ensures that the diffusion matrix $\sigma$ effectively preserves the noise $Z_t$ over the whole space $\R^d$, uniformly in time and space.

When considering SDE \eqref{INTRO:eq:SDE} driven by a truly multiplicative noise, i.e. in the presence of a diffusion matrix $\sigma$ which is also dependent on $x$, we will also have to assume that the measure $\nu$ is absolutely continuous with respect to the Lebesgue measure on $\R^d$ and that their Radon-Nykodim derivative is Lipschitz continuous. More precisely, we will impose the following condition:
\begin{description}
\item[{[AC]}] if $\sigma(t,\cdot)$ is non-constant for some $t\ge0$, then there exists a Lipschitz continuous function $g\colon \mathbb{S}^{d-1}\to \R$ such that
\[\nu(dz) \, = \, Q(z)\frac{g\left(\frac{z}{|z|}\right)}{|z|^{d+\alpha}}dz.\]
\end{description}
Even if Condition [\textbf{AC}] significantly reduces the class of Lévy measures, and therefore processes, which we can consider as a driving noise in Equation \eqref{INTRO:eq:SDE}, we highlight that this assumption seems to be necessary in our framework, at least with respect to the approach considered (cf.\ Equation \eqref{INTRO:eq_per_AC} below) and also natural, since it appeared in other past works on similar topics (cf.\ \cite{Huang:Menozzi16, Frikha:Konakov:Menozzi21}). However, we remark that at least in the additive case, or more generally, when the diffusion matrix $\sigma$ does not depend on space, the non-degeneracy condition [\textbf{ND}] considered in Chapter \ref{Chap:Schauder_estimates_Stable} for the Schauder estimates in an $\alpha$-stable framework, can be understood as a special case of Condition [\textbf{ND}'] assumed here. In particular, we have not imposed any regularity on the Lévy measure $\nu$ of the process $\{Z_t\}_{t\ge0}$ which, a priori, could have a strongly singular support on the space $\R^d$. For example, our model allows us to choose as driving noise $\{Z_t\}_{t\ge 0}$ the cylindrical $\alpha$-stable process, corresponding to the infinitesimal generator presented in \eqref{INTRO:eq:def_stabile_cilindrico}, whose spectral measure has a support concentrated only on the axes of $\R^d$.

We can now summarise the main results that we will prove in Chapter \ref{Chap:Weak_Well-Posedness}. We begin by showing that under minimal Hölder regularity assumptions on the deterministic drift $F$, it is possible to show the well-posedness in a weak sense of the stochastic dynamics in \eqref{INTRO:eq:SDE}.

\begin{teorema}
\label{INTRO:thm:main_result}
For any $j$ in $\llbracket 2, n\rrbracket$, let $\beta^j$ be an index in $(0,1)$ such that
\begin{itemize}
    \item $x_j \to F_i(t,x_i,\dots,x_j,\dots,x_n)$ is $\beta^j$-Hölder continuous, uniformly in time and in the other spacial variables, for any $i$ in $\llbracket 1, j\rrbracket$.
\end{itemize}
Then, SDE \eqref{INTRO:eq:SDE} is weakly well-posed if
\begin{equation}
\label{INTRO:eq:thresholds_beta}
\beta^j\, > \, \frac{1+\alpha(j-2)}{1+\alpha(j-1)},\quad j \ge 2.
\end{equation}
\end{teorema}

To show the above result, we will exploit the equivalence, explained in Section \ref{Sec:INTRO:Unicita_in_legge}, between the well-posedness in a weak sense for SDEs and the corresponding well-posedness for the martingale problem associated with the operator $L_s$, where $L_s$ is (formally) the infinitesimal generator of the solution process $\{X^{t,x}_s\}_{s\ge 0}$ to \eqref{INTRO:eq:SDE}. More precisely, remembering that we consider a pure jump driving process $\{Z_s\}_{s\ge0}$, the operator $L_s$ can be represented for any regular enough function $\phi\colon \R^N\to \R$ as
\begin{align}
   \label{INTRO:eq:def_generator}
L_s\phi(s,x) \, &:= \, \langle G(s,x),D_x\phi(x)\rangle +\mathcal{L}_s\phi(s,x) \\\notag
&:= \,\langle A_sx+F(s,x),D_x\phi(x)\rangle + \text{p.v.}\int_{\R^d}\bigl[\phi(x+B\sigma(s,x)z)-\phi(x) \bigr] \,\nu(dz),
\end{align}
where, for simplicity, we absorbed the term involving $b$ inside $F$ and, similarly to Chapter \ref{Chap:Schauder_estimates_Stable}, we exploited the symmetry of $\nu$ to cancel the first order term $\langle D_x\phi(x),B\sigma(t,x)z\rangle$ in the integral. 

Our method of proof will also require, as an intermediate step, to prove a particular type of controls, called Krylov estimates, for the solution to SDE \eqref{INTRO:eq:SDE}. Such controls owe their name to N.V.\ Krylov who was the first to prove them in \cite{Krylov71} for It\^o diffusions. Since then, they have become a versatile tool in many different frameworks, from showing the well-posedness in a weak or strong sense of SDEs to applications to control theory and non-linear filtering problems. In a non-degenerate, $\alpha$-stable ($\alpha>1$) context where such estimates has been exploited, we cite, for example, \cite{Kurenok08} where the existence of weak solutions is showed for a one-dimensional stochastics dynamic driven by a multiplicative noise and with a bounded measurable drift, or \cite{Zhang13}, where the strong well-posedness for a stochastic equation driven by an additive noise and a singular drift in appropriate Sobolev spaces (under conditions similar to \eqref{INTRO:homogeneous_condizione_C} for $n=1$) is proved. For other similar results in the non-degenerate multiplicative jump noise case, see also the following works: \cite{Anulova:Pragarauskas77,Melkinov83,Lepeltier:Marchal76}. To precisely state the Krylov-type estimates in our context, we will have to impose some conditions on the integrability indexes of the space $L^p(0,T;L^q(\R^N))$. Intuitively, this threshold will guarantee the necessary integrability with respect to the intrinsic scales given by the degenerate nature of the considered dynamics. For simplicity, we will say that two real numbers $p$, $q$ in $(1,+\infty)$ satisfy the integrability condition $(\mathscr{C})$ if:
\[ \tag{$\mathscr{C}$}
\bigl(\frac{1-\alpha}{\alpha}N+\sum_{i=1}^n id_i\bigr)\frac{1}{q}+\frac{1}{p} \, < \, 1.
\]
Such a threshold actually becomes clearer if we consider the homogeneous case, i.e. when all the components $\R^{d_i}$ of the system have the same size ($d_i=d$ and $N=nd$). Indeed in this case, the condition $(\mathscr{C})$ can then be rewritten as
\begin{equation}
\label{INTRO:homogeneous_condizione_C}
\left(\frac{2+\alpha(n-1)}{\alpha}\right)\frac{nd}{q}+\frac{2}{p} \, < \, 2.
\end{equation}
In this form, the above threshold can naturally be understood as a generalisation of the condition imposed in \cite{Chaudru:Menozzi17}, to obtain the same type of estimates in the degenerate diffusive case ($\alpha=2$). We also mention that Condition $(\mathscr{C})$ also appears in \cite{Krylov:Rockner05} in a non-degenerate diffusive context ($\alpha=2$ and $n=1$) as the assumption (of integrability on $f$) necessary to integrate a function $f$ against the Gaussian density (cf.\ Equation $(3.2)$ in the proof of Lemma $3.2$ in \cite{Krylov:Rockner05}).
\begin{teorema}
\label{INTRO:thm:Krylov}
Under the hypothesis of Theorem \ref{INTRO:thm:main_result}, let $T>0$ and $p,q$ be in $(1,+\infty)$ such that Condition $(\mathscr{C})$ holds. Then, there exists a constant $C:=C(T,p,q)$ such that for any $f$ in $L^p\bigl(0,T;L^q(\R^N)\bigr)$,
\begin{equation}
\label{INTRO:Stime_Krylov}
\Bigl{|}\mathbb{E}\bigl[\int_t^Tf(s,X^{t,x}_s) \, ds\bigr] \Bigr{|} \, \le \, C\Vert f \Vert_{L^p_tL^q_x}, \quad (t,x) \in [0,T]\times \R^N,    
\end{equation}
where $\{X^{t,x}_s\}_{s\ge0}$ is the unique weak solution to SDE \eqref{INTRO:eq:SDE} with initial condition $(t,x)$.
\end{teorema}
This type of estimates also emphasises that the solution process $\{X_s\}_{s\ge 0}$ actually has a density with suitable integrability properties up to a certain threshold determined by Condition $(\mathscr{C})$.

In \cite{Chaudru18}, Chaudru de Reynal characterised, through suitable counter-examples, the minimal Hölder regularity for the weak well-posedness of a kinetic stochastic dynamics with degenerate diffusive noise (cf.\ Equation \eqref{INTRO:eq:SDE} with $\alpha=2$ and $n=2$). Through an extension of such Peano-type counterexamples, we also have been able to show a non-uniqueness result for the degenerate chain in \eqref{INTRO:eq:SDE}. We recall that $\{e_i\colon i \in \llbracket1,n\rrbracket\}$ is the canonical basis on the space $\R^N$.
\begin{teorema}
\label{INTRO:thm:counterexample}
Fixed $j$ in $\llbracket 2,n\rrbracket$ and $i$ in $\llbracket 2,j\rrbracket$, there exists a drift $F(t,x)=e_i\text{sgn}(x_j)|x_j|^{\beta^j_i}$ with
\[
\beta^j_i \, < \, \frac{1+\alpha(i-2)}{1+\alpha(j-1)},
\]
for which the weak uniqueness of solutions fails in SDE \eqref{INTRO:eq:SDE}.
\end{teorema}

Unlike the degenerate Gaussian case analysed in \cite{Chaudru:Menozzi17}, it has not been possible to show here that the Hölder regularity on the coefficients $F_i$ is actually the minimal one necessary to guarantee the well-posedness of the stochastic dynamics. As better explained further on, this problem is intrinsically connected to the $\alpha$-stable nature of the process $\{Z_s\}_{s\ge0}$, to the degenerate nature of the dynamics in \eqref{INTRO:eq:SDE} and in particular, to the possibly very singular geometry of the spectral measure associated with the underlying proxy process, understood as the linearized model (i.e.\ Ornstein-Uhlenbeck type process) around which we will apply the ``perturbative'' method. However, we remark that at least for a degenerate chain perturbed by a nonlinear drift $F$ only on the diagonal, our results show the wanted optimal thresholds. Indeed, let us consider a system of the form:
\begin{equation}
\label{INTRO:eq:DIAG}
  \begin{cases}
    dX^1_t \, = \, F_1(t,X^1_t,\dots,X^n_t)dt + \sigma(t,X^1_{t-},\dots,X^n_{t-})dZ_t,\\
    dX^2_t \, = \, \left[A^2_tX^1_t +F_2(t,X^2_t)\right] dt,\\
    dX^3_t \, = \, \left[A^3_tX^2_t+F_3(t,X^3_t)\right]dt,\\
    \vdots\\
    dX^n_t \, = \, \left[A^n_t X^{n-1}_t+F_n(t,X^n_t)\right]dt,
  \end{cases}
\end{equation}
that is, such that the function $F$ only depends on the current level of the chain. We then notice that Theorems \ref{INTRO:thm:main_result} and \ref{INTRO:thm:counterexample} together present an (almost) complete characterisation of the weak well-posedness for degenerate stochastic dynamics such as \eqref{INTRO:eq:DIAG}, with respect to the Hölder regularity of their coefficients. Indeed, we know that:
\begin{itemize}
    \item if $\beta^j >\frac{1+\alpha(j-2)}{1+\alpha(j-1)}$ for any $j\ge 0$, the well-posedness in a weak sense of SDE \eqref{INTRO:eq:SDE} can be derived from Theorem \ref{INTRO:thm:main_result};
    \item if there exists $j\ge 2$ such that $\beta^j=\beta_j^j<\frac{1+\alpha(j-2)}{1+\alpha(j-1)}$, then there exist some counter-examples (cf.\ Theorem \ref{INTRO:thm:counterexample}) for which the weak uniqueness of solutions fails for a dynamics as  \eqref{INTRO:eq:SDE}.
\end{itemize}
We finally mention that the \emph{critical} case, corresponding to the following exponents
\[\overline{\beta}^j_j \, = \, \frac{1+\alpha(j-2)}{1+\alpha(j-1)}, \quad
 j \in \llbracket 2,n\rrbracket, \]
remains to be investigated and seems to be a delicate problem already in a kinetic diffusive framework (cf.\ \cite{Zhang18}).

\subsection{Sketch of the proof}

We have decided to follow here a perturbative approach by \emph{backward} parametrix expansions as originally introduced by McKean-Singer in \cite{Mckean:Singer67} in a non-degenerate diffusive context and then extended to the degenerate case with unbounded drifts in \cite{Delarue:Menozzi10, Menozzi18}. This terminology comes from the fact that the underlying proxy process will be associated with a backward in time flow, that is, we will fix $(\tau,\xi)=(s,y)$ in \eqref{INTRO:eq:flusso_associato}. This method appears to be particularly useful to show the weak uniqueness in a degenerate context on $L^p_t-L^q_x$ spaces (cf.\ \cite{Chaudru:Menozzi17} in the Brownian case). In fact, it intuitively requires only to establish bounds on the gradients (in a weak sense) of the solutions to the associated Cauchy problem, in order to to apply the operator inversion technique, as originally exploited in \cite{book:Stroock:Varadhan79} in a diffusive case.

Despite what was initially suggested in Section \ref{Sec:INTRO:Unicita_in_legge}, we finally decided not to exploit the Schauder estimates, shown in the previous two sections, to prove the uniqueness in law of SDE \eqref{INTRO:eq:SDE}. We are however convinced that the perturbative approach by forward parametrix expansions presented in Section \ref{Sec:INTRO:Stime_Schauder_stabili} could have been extended here to establish the Schauder estimates for the class of processes considered now. From such estimates, we could have then proved the well-posedness of the stochastic dynamics in a weak sense, through a Zvonkin type argument. Such a method usually appears in the proof for the strong well-posedness of SDEs but we remark that it has also been used to establish the weak well-posedness for degenerate chains, as done, for example, in \cite {Chaudru18}. It is the latter type of reasoning that we will refer to as a Zvonkin-type argument below. However, this method immediately appears very long and complicated since it would have required extending the Schauder estimates, presented in \eqref{INTRO:eq:Schauder_estimates}, under the expected optimal regularity in $C^{\beta}_{b,d} (\R^N)$ on the source $f$, to a more generic Hölder space with multiple unrelated regularity indexes, like $C^{\mathbf{\gamma}}_{b,d}( \R^N)$ with $\mathbf{\gamma}=(\gamma_1,\dots,\gamma_d)$ in $\N^d$, precisely because of the Zvonkin method. In particular, it would have required establishing point-wise estimates on the first order derivatives of the solution to the associated Cauchy problem also with respect to the degenerate components of the system. To do this, we would have needed to extend our duality reasoning on Besov spaces also with respect to the degenerate derivatives and to consider a source $f$ with the same regularity as the drift $F$ in a space $C^\gamma_{ b,d}(\R^N)$ with multi-indexes of regularity. Another possible advantage of the backward perturbative method is that it allows to easily obtain, as an intermediate step in our proof, the Krylov estimates in \eqref{INTRO:Stime_Krylov} for the process $\{X^{t,x}_s\}_ {s\ge 0}$, solution to SDE \eqref{eq:SDE}. This type of controls seems to be new for degenerate stochastic chains driven by quasi-stable noises and it would have been very difficult to establish from the Schauder estimates. Finally, we remark that our method, compared with the approach through Zvonkin transform as in \cite{Chaudru18}, allows to obtain a more precise analysis of the degenerate chain at least along the first component, in the sense that we can highlight here that no minimal thresholds are required on the spacial regularity of $F_1$. Intuitively, the method through Zvonkin transform requires to put each component $F_i$ of the drift as source. This leads in particular to the same global optimal thresholds on the Hölder regularity for $F$ at each level of the chain (cf.\ Equation \eqref{INTRO:eq:thresholds_beta} with $j\ge 1$). As shown for example in \cite{Chaudru17,Fedrizzi:Flandoli:Priola:Vovelle17, Chaudru:Menozzi17} in the diffusive case, this approach seems to be more adequate when one wants to obtain the strong well-posedness of the stochastic dynamics.

For the purposes of this presentation, we will limit ourselves to consider a matrix $A_t$ in \eqref{INTRO:decomposizione_deriva_G} independent from time and such that only its elements on the sub-diagonal are non-zero, i.e., we will assume that $A$ is given in \eqref{INTRO:eq:matrix_A_sub-diag}. In fact, as explained in the previous section, the additional difficulties given by the presence of non-zero elements above the sub-diagonal can be easily solved through an additional argument in small time. Furthermore, the possible dependence on time would essentially require introducing, instead of the exponential matrix $e^{A(s-t)}$, the resolvent $\mathcal{R}_{s,t}$ associated with the matrix $A_t$. In practice, $\mathcal{R}_{s,t}$ is the time dependent matrix in $\R^N\otimes \R^N$ solution to the following matrix differential equation:
\[
\begin{cases}
  \partial_s \mathcal{R}_{s,t} \, = \, A_s \mathcal{R}_{s,t} ,\quad s \in [t,T]; \\
  \mathcal{R}_{t,t} \, = \, \Id_{N\times N}.
\end{cases}
\]
Finally, we mention that decompositions similar to those presented in \eqref{INTRO:decomposition_exponential} have also been shown to be valid for the resolvent $\mathcal{R}_{s,t}$. See for example \cite{Huang:Menozzi16}, Lemmas $5.1$ and $5.2$ or \cite{Delarue:Menozzi10}, Proposition $3.7$.

\subsubsection{Perturbative approach by backward parametrix expansions}

As already explained in Section \ref{Sec:INTRO:Stime_Schauder_stabili}, the crucial element of the perturbative method consists in carefully choosing a suitable proxy operator with known properties and bounds, around which to expand the infinitesimal generator $L_s$, at the price of an additional expansion error to control. When the drift $F$ is sufficiently regular, for example globally Lipschitz continuous, it has been shown in \cite{Delarue:Menozzi10,Menozzi11,Menozzi18} that a suitable proxy is given by the linearization of SDE \eqref{INTRO:eq:SDE} around the deterministic flow associated with the dynamics (i.e.\ when $\sigma=0$ in \eqref{INTRO:eq:SDE}), which leds, in the works above in a degenerate diffusive context, to consider a multi-scale Gaussian process as a proxy. The natural generalisation in our context will instead lead to choose a multi-scale Lévy  process with time dependent symbol. More precisely, given the freezing parameters $(s,y)$ in $[0,T]\times \R^N$, let $\theta_{t,s}(y)$ be one of the possible solutions to the following equation:
\begin{equation}
\label{INTRO:eq:flusso_retrogrado}
\theta_{t,s}(y) \, = \,  y - \int_{t}^{s}\bigl[ A \theta_{u,s}( y)+F(u, \theta_{u,s}(y))\bigr] \,  du.
\end{equation}
Since the freezing point $y$ will be also the integration variable (see for example the definition of the Green kernel in \eqref{INTRO:eq:def_Green_Kernel}), it is important to remark right now that it is always possible, among the possible solutions $\theta_{t,s}(y)$ to \eqref{INTRO:eq:flusso_retrogrado}, to choose one that is measurable with respect to $(s,y)$ in $[0,T]\times \R^N$ (cf.\ Lemma \ref{lemma:measurability_flow} in Chapter \ref{Chap:Weak_Well-Posedness}). From now on we will assume that we have fixed such a version of $\theta_{t,s}(y)$.

The next step will be to introduce the \emph{linearized} stochastic dynamics around the backward flow $\theta_{t,s}(y)$. More precisely, we will consider for any starting point $(t,x)$ in $[0,s]\times \R^N$, the proxy process $\{\tilde{X}^{t,x,s,y }_u\}_{u\ge 0}$ solution to the following stochastic dynamics:
\begin{equation}
\label{INTRO:eq:SDE_frozen_SDE}
\begin{cases}
d\tilde{X}^{t,x,s,y}_u\, = \, \bigl[A\tilde{X}^{t,x,s,y}_u+\tilde{F}^{s, y}_u\bigr]\, du +B\tilde{\sigma}^{s,y}_u\,  dZ_u, \quad u\in  [t,T],\\
 \tilde{X}^{t,x,s,y}_t\, = \, x,
\end{cases}
\end{equation}
where, for simplicity, we denoted $\tilde{F}^{s,y}_u:=F(u,\theta_{u,s}(y))$ and $\tilde{\sigma}^{s,y}_u:=\sigma(u,\theta_{u,s}(y))$. 
A direct integration then allow to obtain a more explicit representation of the proxy process:
\begin{equation}
\label{INTRO:eq:processo_proxy}
\tilde{X}^{t,x,s,y}_s \, = \, \tilde{m}^{s,y}_{s,t}(x) +\int_t^se^{A(s-u)}B\tilde{\sigma}^{s,y}_u dZ_u,
\end{equation}
where, analogously to \eqref{INTRO:eq:rappresentazione_X_tilda}, the frozen transport term $\tilde{m}^{s,y}_{s,t}(x)$ is given by
\[
\tilde{{m}}^{s,y}_{s,t}(x) \, = \, e^{A(s-t)} x + \int_{t}^{s}e^{A(s-u)} \tilde{F}^{s,y}_u \, du.
\]

By reasoning in Fourier spaces, similarly to those in \eqref{INTRO:prova1}-\eqref{INTRO:prova2}, we can then show that the following crucial identity in law holds even in this case:
\begin{equation}
\label{INTRO:identita_in_legge2}
\tilde{X}^{t,x,s,y}_s \, \overset{(\text{law})}{=} \, \tilde{{m}}^{s,y}_{s,t}(x) +\mathbb{M}_{s-t} \tilde{S}^{s,y}_{s-t},
\end{equation}
where for any fixed freezing parameter $(s,y)$, $\{\tilde{S}^{s,y}_u\}_{u\ge 0}$ can be understood as a non-degenerate (in the sense indicated in [\textbf{ND}']), $\alpha$-stable process on $\R^N$. We also point out that the dependence on the freezing parameters, a crucial difficulty in this part of the proof, is essentially linked to the presence of the frozen diffusion matrix $\tilde{\sigma}^{s,y}$ in the stochastic convolution in \eqref{INTRO:eq:processo_proxy}. In fact, this dependence would disappear in the case of additive noise (i.e.\ $\sigma(t,x)=1$) or, more generally, when considering a diffusion coefficient homogeneous in space. As already explained in Section \ref{Sec:INTRO:Stime_Schauder_stabili}, the non-degeneracy of the spectral measure of the process $\{\tilde{S}^{s,y}_{u}\}_{u\ge 0 }$ ensures in particular the existence of a regular (in space) density $p_{\tilde{S}^{s,y}}(u,\cdot)$ for this process. The identity in \eqref{INTRO:identita_in_legge2} then implies that the following function
\[\tilde{p}^{s,y}(t,s,x,y) \, := \,\frac{ p_{\tilde{S}^{s,y}}(s-t,\mathbb{M}^{-1}_{s-t}(y-\tilde{m}^{s,y}_{t,s}(x)))}{\det \mathbb{M}_{s-t}},\]
is the ``frozen'' density for the proxy process $\{\tilde{X}^{t,x,s,y}_s\}_{s\ge 0}$, frozen at the final space-time point $(s,y)$.

Similar to the previous sections, we will next focus on determining the smoothing properties associated with $\tilde{p}^{s,y}(t,s,x,y)$. In particular, we will show that the derivatives of the frozen ``density'' are bounded from above by another density at the price of additional singularities in time and importantly, that such control holds \emph{uniformly} in the freezing parameters $(s,y)$. More precisely, it holds for any $k$ in $\llbracket 0,2 \rrbracket$ and any $i$ in $\llbracket 1,n \rrbracket$, that
  \begin{equation}
  \label{INTRO:eq:derivative_prop_of_q}
      \vert D^k_{x_i}\tilde{p}^{s,y}(t,s,x,y) \vert \le C\frac{\overline{p}\left(1,\mathbb{T}^{-1}_{s-t}(y-\tilde{m}^{s,y}_{s,t}(x))\right)}{\det \mathbb{T}_{s-t}}(s-t)^{-k\frac{1+\alpha(i-1)}{\alpha}},
  \end{equation}
where, for simplicity, we denoted $\mathbb{T}_u=u^{\frac{1}{\alpha}}\mathbb{M}_u$. This estimate is essentially analogous to the one obtained in \eqref{INTRO:eq:Smoothing_effects}, where for convenience we have already rescaled the system to the unit time $t=1$, using the self-similarity of the $\alpha$-stable density $\overline{p}(u,\cdot)$, i.e.:
\[\overline{p}(u,z) \, = \, u^{-N/\alpha}\overline{p}(1,u^{-1/\alpha}z), \quad (u,z) \in [0,\infty)\times \R^N. \]
Since the freezing point $y$ will also appear as an integration variable (cf.\ definition of $\tilde{G}_\varepsilon f$ in \eqref{INTRO:eq:def_Green_Kernel}), we highlight how important it is to obtain a bound that is independent of this parameter. Estimates on the density as in \eqref{INTRO:eq:derivative_prop_of_q} or in \eqref{INTRO:eq:Smoothing_effects} are often obtained through the It\^o-Lévy decomposition of the random variable $\tilde{S}^{s,y}_u$ at the corresponding stable characteristic time. More precisely, let $\tilde{M}^{s,y}_u$ and $\tilde{N}^{s,y}_u$ be the two independent random variables associated with the small and big jumps of the process $\{\tilde{S}^{s,y}_u\}$, truncated at time $u^{1/\alpha}$. This truncation will allow in particular to rewrite the density $p_{\tilde{S}^{ s,y}}(u,z)$ of $\tilde{S}^{s,y}_u$ as follows:
\begin{equation}
\label{INTRO:decomposizione_Ito_Levy}
p_{\tilde{S}^{s,y}}(u,z) \, = \, \int_{\R^N}p_{\tilde{M}^{s,y}}(u,z-w)\, P_{\tilde{N}^{s,y}_u}(dw)
\end{equation}
where $p_{\tilde{M}^{s,y}}(u,\cdot)$ is the density generated by $\tilde{M}^{s,y}_u$ and $P_{\tilde{N}^{s,y}_u}$ is the law of $\tilde{N}^{s,y}_u$. While arguments such as those carried out in Section \ref{Sec:INTRO:Stime_Schauder_Levy} (cf.\ Equation \eqref{INTRO:smoothing_effect_Levy}) on the truncated density $p^{\text{tr}}$ can also be applied in this case leading to the following estimates:
\begin{equation}
\label{INTRO:proof:control_pM_segnato}
\left| D^k_{z} p_{\tilde{M}^{s,y}}(u,z) \right| \, \le \, Cu^{-(N+k)/\alpha}\left(\frac{u^{1/\alpha}}{u^{1/\alpha}+|z|}\right)^{N+3} \, =: \, Cu^{-\frac{k}{\alpha}}p_{\bar{M}}(u,z),
\end{equation}
where $C$ is independent from the freezing parameters $(s,y)$; it will be much more delicate to establish a uniform control on the probability measures $P_{\tilde{N}^{s,y}_u}(dy)$ of the type:
\begin{equation}
\label{INTRO:proof_control_N_segnato}
P_{\tilde{N}^{s,y}_u}(\mathcal{A}) \, \le \, C\overline{P}_u(\mathcal{A}), \quad \mathcal{A} \in \mathcal{B}(\R^N),
\end{equation}
where $\{\overline{P}_u\}_{u\ge 0}$ is a family of probability measures which preserve the same integrability properties as the tails of and $\alpha$-stable process. This difficulty is the main reason why we could not consider, as in \cite{Chaudru:Menozzi17}, a completely non-linear drift, i.e.\ $G_i(t,x)=G_i(t,x_{i-1},\cdots,x_n) $ with a non-linear dependence on the variable $x_{i-1}$ that transmits the noise, but only a semi-linear version of $G$, given in \eqref{INTRO:decomposizione_deriva_G}. In fact, the more general model analyzed in \cite{Chaudru:Menozzi17} would have required in particular to linearise SDE \eqref{INTRO:eq:SDE} around the sub-diagonal entries of the Jacobian matrix of $G$ frozen at $(s,y)$, or to consider the following proxy process:
\[d\tilde{X}^{t,x,s,y}_u\, = \,  \bigl[\tilde{A}^{s,y}_u\tilde{X}^{t,x,s,y}_u+\tilde{F}^{s, y}_u\bigr]\, du +B\tilde{\sigma}^{s,y}_u\,  dZ_u,\]
where $\tilde{A}^{s,y}_u$ is a matrix such that
\[[\tilde{A}^{s,y}_u]_{i,j} \, = \,
      \begin{cases}
      D_{x_{i-1}}G_i(s,\theta_{t,s}(y)), &\mbox{ if }j=i-1,\\
      0, &\mbox{ otherwise.}
      \end{cases}\]
For such a model, we have not been able to actually show a uniform bound like in \eqref{INTRO:proof_control_N_segnato}. We also remark that it is for this reason that we have not been able to consider a Lévy measure $\nu$ associated with the process $\{Z_t\}_{t\ge0}$ that is asymmetric, as done in Chapter \ref{Chap:Schauder_Estimates_Levy}. In fact, the backward perturbative method requires more delicate regularizing properties associated with the operators involved and more importantly, a certain compatibility between the proxy and the original process. Intuitively, the symmetry of the underlying Lévy measure is not a strong constraint for the contribution in \eqref{INTRO:proof:control_pM_segnato} associated with small jumps and this allowed in fact to consider possibly asymmetric operators in the previous section. In our model, which instead requires a global smoothing effect for the density, it then seems natural to impose the symmetry of the corresponding spectral measure in the control of the process tails, as we will see in the proof of the estimates in \eqref{INTRO:proof_control_N_segnato} (cf.\ Lemma \ref{lemma:Control_P_N} in Chapter \ref{Chap:Weak_Well-Posedness}).

\subsubsection{Analytic properties of the frozen density}

We highlight that it is not immediate to determine which type of Cauchy problem the ``density'' $\tilde{p}^{s,y}(t,s,x,y)$ frozen at the terminal point $(s, y)$ is indeed a fundamental solution for. In fact, the choice of the freezing point $y$ also as an integration variable (for example in \eqref{INTRO:eq:def_Green_Kernel}) makes it much more difficult to determine its analytical properties and above all, prove its convergence to the Dirac mass $\delta_x$ when the time $t$ tends to zero. Let $\tilde{L}^{s,y}_t$ be the infinitesimal generator associated with the frozen process $\{\tilde X_{s}^{s,y,t,x}\}_{s\ge 0}$. More precisely, we write for any regular enough function $\phi\colon \R^N\to\R$, that
\begin{equation}
\label{INTRO:eq:def_frozen_generator}
\begin{split}
 \tilde{L}^{s,y}_t\phi(x)\, &:= \,
\langle A x+\tilde{F}^{s, y}_t, D_x \phi(x)\rangle + \tilde{\mathcal{L}}^{s,y}_t\phi(x) \\
&:= \langle A x+\tilde{F}^{s, y}_t, D_x \phi(x)\rangle + \text{p.v.}\int_{\R^d_0}\bigl[\phi(x+B\tilde{\sigma}^{s, y}_t w)-\phi(x)\bigr] \,\nu(dw),
\end{split}
\end{equation}
where, we recall, we denoted $\tilde{F}^{s,y}_t:=F(t,\theta_{t,s}(y))$ and $\tilde{\sigma}^{s,y}_t:= \sigma(t,\theta_{t,s}(y))$. 
If we now fix the freezing parameter $y$ and vary only the integration variable in $\tilde{p}^{s,y}(t,s,x,\cdot)$, the function becomes a true density and it is then not difficult to show, through a direct calculation, that:
\begin{equation}\label{INTRO:eq:differential-eq}
 \left(\partial_t + \tilde{L}^{s,y}_t\right) \tilde{p}^{s,y} (t,s,x,z) \,  = \, 0, \quad (t,z) \in [0,s)\times \R^N,
\end{equation}
for any fixed $(s,x,y)$ in $[0,T]\times \R^{2N}$.

To determine which kind of parabolic system is solved by the frozen density $(s,y)\mapsto\tilde{p}^{s,y}(t,s,x,y)$, we now consider an operator $\tilde{G}_\varepsilon$ which can be understood as the Green kernel associated with the ``density'' $\tilde{p}^{s,y}(t,s,x,y)$ and concentrated away from the initial time $t$. In particular, given $\varepsilon>0$ small enough for our purposes, we have that
\begin{equation}\label{INTRO:eq:def_Green_Kernel}
\tilde{G}_\varepsilon f(t,x) \, := \, \int_{t+\varepsilon}^T \int_{\R^N}  \tilde{p}^{s,y}(t,s,x,y) f(s,y) \, dy ds, \quad (t,x) \in [0,T)\times \R^N,
\end{equation}
for any smooth enough function $f\colon [0,T)\times \R^N\to \R$ with compact support. The expression in \eqref{INTRO:eq:def_Green_Kernel} is well-defined, since we know that the frozen density $\tilde{p}^{s,y}(t,s,x,y)$ is measurable in $(s,y)$.
The localisation in $\varepsilon$ away from $t$ in the integral is crucial since it ensures the wanted regularising effect for the Green kernel $\tilde{G}_\varepsilon$.
Indeed, we show that, in the limiting case (i.e.\ $\varepsilon \to 0$), the regularity of the function $f$ is not a sufficient condition to derive the regularity for the Green kernel $\tilde{G}_0 f$. This further difficulty arises from the dependence of the proxy on the integration variable $y$.

If we now introduce the following quantity:
\[
\tilde{M}_\varepsilon f(t,x)\, := \, \int_{t+\varepsilon}^T \int_{\R^N}  \tilde{L}^{s,y}_t \tilde{p}^{s,y}(t,s,x,y)f(s,y) \, dy ds, \quad (t,x) \in [0,T)\times\R^N,
\]
we can finally derive from Equation
 \eqref{INTRO:eq:differential-eq} that the Green operator $\tilde{G}_{\varepsilon}$ solves the following parabolic equation:
\begin{equation}\label{INTRO:eq:differential_eq2}
\partial_t \tilde{G}_\varepsilon f(t,x)+ \tilde{M}_\varepsilon f(t,x) \, = \, -I_\varepsilon f(t,x),\quad (t,x)\in [0,T)\times \R^N.
\end{equation}
Also in this case, the localisation in $\varepsilon$ has been crucial to obtain \eqref{INTRO:eq:differential_eq2} directly from the Equation \eqref{INTRO:eq:differential-eq}, exploiting classical dominated convergence arguments. Above, the operator $I_\varepsilon$ can be represented for any sufficiently regular function $f\colon [0,T)\times \R^N \to \R$ as
\begin{equation}
\label{INTRO:eq:def_f_epsilon}
I_\varepsilon f(t,x) \, := \, \int_{\R^N} f(t+\varepsilon,y)\mathds{1}_{[0,T-\varepsilon]}(t)
\tilde{p}^{t+\varepsilon,y}(t,t+\varepsilon,x,y) \, dy.
\end{equation}
We now notice that $I_\varepsilon$ can be understood as a version  of the identity operator localised outside $t$. In particular, we will obtain in Chapter \ref{Chap:Weak_Well-Posedness} (cf.\ Lemmas \ref{convergence_dirac} and \ref{prop:convergence_LpLq}), its convergence to the Dirac mass concentrated in $(t,x) $ on the suitable functional norms:
\begin{equation}
\label{INTRO:eq:Control_I_epsilon}
\lim_{\varepsilon\to 0}\Vert I_\varepsilon f -f \Vert_\infty \,= \, 0, \quad \lim_{\varepsilon\to 0}\Vert I_\varepsilon f -f \Vert_{L^p_tL^q_x} \,= \, 0.    
\end{equation}
Although at first sight the convergence properties above seem to be immediate, we remark that the presence of the integration variable $y$ also as a freezing parameter prevents us from obtaining the estimates in \eqref{INTRO:eq:Control_I_epsilon} directly from the convergence in law of the frozen process $\{\tilde{X}_{s}^{s,y,t,x}\}_{s\ge 0}$ to the Dirac mass (cf.\ Equation \eqref{INTRO:eq:differential-eq}).

\subsubsection{Uniqueness of the martingale problem and associated Krylov estimates}

We now present the main steps in the proof for the uniqueness of the martingale problem associated with $\mathcal{L}_t$. As already mentioned before, the analytical theory associated with the proxy process will be the fundamental tool in our method of proof.

We will initially exhibit the Krylov-type estimates in \eqref{INTRO:Stime_Krylov} using the perturbative method with backward proxy, assuming however that the indices $p,q$ are sufficiently large but finite. More precisely, given a sufficiently regular function $f$ and a solution $\{X^{t,x}_s\}_{s\ge 0}$ of the stochastic dynamics in \eqref{INTRO:eq:SDE}, the first step of our method consists in applying the It\^o formula to the frozen Green operator $\tilde{G}_\varepsilon f$ with respect to the process $\{X^{t,x}_s \}_{s\ge 0}$:
\begin{equation}
\label{INTRO:eq1}
\mathbb{E}\left[\tilde{G}_\varepsilon f(t,x)+ \int_t^{T} (\partial_s+L_s)\tilde{G}_\varepsilon f(s,X^{t,x}_s)   ds\right] \, =\, 0.
\end{equation}
Recalling that $\tilde{G}_\varepsilon f$ solves Equation, we can rewrite the above equation \eqref{INTRO:eq1} as follows:
\begin{equation}
\begin{split}
\label{INTRO:eq:NUMBER_PREAL_REG_DENS}
\mathbb{E}\Bigl[\int_t^{T} \int_{\R^N}I_\varepsilon f(s,X^{t,x}_s) \, ds \Bigr] \,
&= \, \tilde{G}_\varepsilon f
(t,x)+\mathbb{E}\Bigl[\int_t^T\left[L_s\tilde{G}_\varepsilon f- \tilde{M}_{\varepsilon}f\right] (s,X^{t,x}_s)\, ds\Bigr]\\
&=: \, \tilde{G}_\varepsilon f
(t,x)+\mathbb{E}\Bigl[\int_t^T\tilde{R}_\varepsilon f(s,X^{t,x}_s)\, ds\Bigr].
\end{split}
\end{equation}
To infer the Krylov estimates in \eqref{INTRO:Stime_Krylov} from the above equation, it will be necessary to show the following controls:
\begin{equation}
\label{INTRO:control_pointwise}
\Vert \tilde{G}_\varepsilon f \Vert_{\infty} \, \le \, C \Vert f \Vert_{L^p_tL^q_x}, \quad \Vert \tilde{R}_\varepsilon f \Vert_{\infty} \, \le \, C \Vert f \Vert_{L^p_tL^q_x}.
\end{equation}
In particular, we need to initially impose a minimal threshold on the indexes $p$, $q$ precisely due to the point-wise bound on the remainder term $\tilde{R}_\varepsilon f$, valid only for $p$ and $q$ large enough. After establishing the estimates in \eqref{INTRO:control_pointwise}, we can then show from Equation \eqref{INTRO:eq:NUMBER_PREAL_REG_DENS} that
\[
\left|\mathbb{E}\left[\int_t^T I_\varepsilon f(s,X_s^{t,x})\,  ds \right]\right|\, \le \,  C\Vert f\Vert_{L^p_tL^q_x}.
\]
Now exploiting the point-wise convergence of $I_\varepsilon$ in \eqref{INTRO:eq:Control_I_epsilon} to make $\varepsilon$ tend to zero, and a regular approximation reasoning on the function space $L^p(0, T;L^q(\R^N))$, it will finally be possible to conclude that Krylov estimates of the form:
\begin{equation}
\label{INTRO:eq_parzial_Krylov}
\left|\mathbb{E}\left[\int_t^T f(s,X_s^{t,x})\,  ds \right]\right|\, \le \,  C\Vert f\Vert_{L^p_tL^q_x}, 
\end{equation}
hold for any $f$ in $L^p(0,T;L^q(\R^N))$, even if under the additional assumption that $p$ and $q$ are large enough.

The ``partial'' result obtained above is sufficient to show the uniqueness of the martingale problem associated with the operator $L_t$. Through a duality reasoning on the $L^p_t-L^q_x$ spaces, the estimates in \eqref{INTRO:eq_parzial_Krylov} imply in particular the existence of a ``density'' for the process $\{X^{t ,x}_s\}_{s\ge 0}$, solution of the stochastic dynamics \eqref{INTRO:eq:SDE} starting at $(t,x)$, but defined only for almost every $(s,y)$. Called $p(t,s,x,y)$ this density, it easily follows from Equation \eqref{INTRO:eq:NUMBER_PREAL_REG_DENS} that
\[\tilde{G}_\varepsilon f(t,x) \, = \, \int_t^{T} \int_{\R^N}(I_\varepsilon-\tilde{R}_\varepsilon)f(s,y) p(t,s,x,y) \, dy  ds.\]
At this point, we want to invert the operator $I_\varepsilon - \tilde{R}_\varepsilon$ on the functional space $L^p(0,T;L^q(\R^N))$ . To do this, we will establish estimates on the remainder term $\tilde{R}_\varepsilon$ of the following form:
\begin{equation}
\label{INTRO:controllo_LpLq_resto}
 \Vert \tilde{R}_\varepsilon \Vert_{L^p_tL^q_x} \, \le \, C_T\Vert f \Vert_{L^p_tL^q_x},   
\end{equation}
for some constant $C_T>0$ independent from $\varepsilon$ and such that $C_T\to 0$, when $T$ tends to zero. Choosing a sufficiently small time interval, we can then assume that the norm of $\tilde{R}_\varepsilon$ as an operator on $L^p(0,T;L^q(\R^N))$ is smaller than $1$ and therefore, that the operator $I_\varepsilon - \tilde{R}_\varepsilon$ is indeed invertible:
\[
\mathbb{E}\left[\int_t^{T} f(s,X_s^{t,x}) \, ds\right]\, = \, \tilde{G}_\varepsilon\circ (I_\varepsilon-\tilde{R}_\varepsilon)^{-1}f(t,x).\]
Exploiting again the convergence properties of the operator $I_\varepsilon$ in \eqref{INTRO:eq:Control_I_epsilon}, 
this time on the space $L^p(0,T;L^q(\R^N))$, we will conclude that
\[
\mathbb{E}\left[\int_t^{T} f(s,X_s^{t,x}) \, ds\right]\, = \, \tilde{G}\circ (I-\tilde{R})^{-1}f(t,x).
\]
The above identity immediately implies the uniqueness of the martingale problem over a sufficiently small time interval. Through a concatenation argument in time, we will then extend it to the global uniqueness for the martingale problem, thus concluding the proof of Theorem \ref{INTRO:thm:main_result}. For more details on this type of localisation, see for example \cite[Section $4.6$]{book:Ethier:Kurtz86}.

We will only show later that the Krylov estimates in \eqref{INTRO:Stime_Krylov} actually hold for any pair $(p,q)$ satisfying the integrability condition ($\mathscr{C}$), through a stochastic mollification argument. More in detail, we will fix an isotropic $\alpha$-stable process $\{\overline{Z}_s\}_{s\ge0}$ on $\R^N$ and a small regularisation parameter $\delta$ and we will consider a regularised version of the solution process $\{X^{t,x}_s\}_{s\ge 0}$ to SDE \eqref{INTRO:eq:SDE}, given by
\begin{equation}
\label{INTRO:eq_mollificazione}
\overline{X}^{t,x,\delta}_s\, := \, X^{t,x}_s+\delta\mathbb{M}_{s-t} \overline{Z}_{s-t}.
\end{equation}
Intuitively, the idea of regularising the solution process $\{X^{t,x}_s\}_{s\ge0}$ serves to obtain, at the level of the associated densities, a controllability in the functional space dual to $L^ {p}_t-L^{q}_x$ up to the sought threshold (cf.\ integrability condition $(\mathscr{C})$). In fact, if we denote respectively with $p'$, $q'$ the conjugate exponents of $p$ and $q$, we know from the ``partial'' Krylov estimates in \eqref{INTRO:eq_parzial_Krylov} that the density $p(t,s,x,y)$ has finite norm in the space $L^{p'}(0,T,L^{q'}(\R^N))$ for large enough values of $p,q$. Furthermore, if we denote by $p^\delta(t,s,x,\cdot)$ the density associated with the random variable $\overline{X}^{t,x,\delta}_{s}$, the identity in \eqref{INTRO:eq_mollificazione}  immediately implies that
\begin{equation}
\label{INTRO:eq:convolution}
p^{\delta}(t,s,x,y) \, = \, \left[p(t,s,x,\cdot)\ast q^\delta(s-t,\cdot)\right](y), 
\end{equation}
where $q^\delta(t,\cdot)$ denotes the density for the process $\{\delta \mathbb{M}_s\overline{Z}_s\}_{s\ge0}$. Exploiting now convolution inequalities arguments in \eqref{INTRO:eq:convolution}, we will show that for any pair $(p,q)$ satisfying condition $(\mathscr{C})$, the quantity $\Vert p^\delta\Vert_{L^{p'}_tL^{q'}_x}$ is finite, even if possibly explosive, when $\delta$ tends to zero. Such a bound will then show that the mollified process $\overline{X}^{t,x,\delta}_s$ satisfies the Krylov estimates in \eqref{INTRO:Stime_Krylov} for all pairs $(p,q)$ in the wanted interval but for a constant $C$ possibly dependent on the mollification parameter $\delta$ (and explosive with respect to this parameter). By reproducing the perturbative analysis carried out in the first part of the proof, we will finally show that the controls on the regularised density $p^\delta(t,s,x,y)$, and therefore the constant in the Krylov estimates, do not actually depend on $\delta$. Finally, letting $\delta$ tend to zero, we will be able to establish the Krylov estimates for the solution $X^{t,x}_s$ to the original SDE \eqref{INTRO:eq:SDE}, under the natural condition $(\mathscr{C})$ on $(p,q)$.

\subsubsection{On the a priori estimates for the Green operator and the remainder term}

As mentioned before, our method of proof relies on some crucial estimates on the Green kernel $\tilde{G}_\varepsilon f$ and the remainder term $\tilde{R}_\varepsilon f$ appeared in \eqref{INTRO:control_pointwise} and in \eqref{INTRO:controllo_LpLq_resto}. The proof of these estimates will be quite long and complex and, together with the proof of the convergence of $I_\varepsilon$ in \eqref{INTRO:eq:Control_I_epsilon}, will occupy a substantial part of the technical section in Chapter \ref{Chap:Weak_Well-Posedness}. Since we do not think it is possible to summarise them here in a coherent way for the reader, we have decided instead to focus only on some of the most salient steps in the reasoning that allow us to show the necessity of the assumptions we have made.

For example, the integrability thresholds in Condition $(\mathscr{C})$ clearly appears in the proof for the point-wise controls (cf.\ \eqref{INTRO:control_pointwise}) of the Green operator $\tilde{G}_\varepsilon$. Indeed, recalling the definition of $\tilde{G}_\varepsilon f$ in \eqref{INTRO:eq:def_Green_Kernel}, one can start by splitting the integral through an Hölder inequality:
\begin{equation}
\label{INTRO:COntrol_G} 
\begin{split}
|\tilde{G}_\varepsilon f(t,x)| \,
&\le \, C \Vert f \Vert_{L^p_tL^q_x} \left(\int_{t+\varepsilon}^T \left (\int_{\R^N} \left|\tilde{p}^{s,y}(t,s,x,y)\right|^{q'} \,dy \right)^{\frac{p'}{q'}} ds \right)^\frac{1}{p'} \\
&=: \, C \Vert f \Vert_{L^p_tL^q_x}\left(\int_{t+\varepsilon}^T \left(|\mathcal{J}(t,s,x)| \right)^{\frac{p'}{q'}} ds \right)^\frac{1}{p'},
\end{split}
\end{equation}
where, we recall, we denoted by $p'$, $q'$ the conjugate exponents of $p$ and $q$, respectively. Exploiting the controls in \eqref{INTRO:eq:derivative_prop_of_q} on the density $\tilde{p}^{s,y}(t,s,x,\cdot)$, it is possible to bound the term $\mathcal{J}(t,s,x)$ in the following way:
\[
|\mathcal{J}(t,s,x)| \, \le \, C \left(\det \mathbb{T}_{s-t}\right)^{-q'}\int_{\R^N}  \left[\bar{p}\left(1, \mathbb{T}^{-1}_{s-t}(y-\tilde{m}^{s,y}_{s,t}(x))\right)\right]^{q'}\, dy \, \le\,  C\left(\det \mathbb{T}_{s-t}\right)^{1-q'}.
\] 
It then follows that
\begin{equation}
\label{INTRO:eq_control_G}
    |\tilde{G}_\varepsilon f(t,x)| \,
\le \, C \Vert f \Vert_{L^p_tL^q_x} \left(\int_{t+\varepsilon}^T \left (\det \mathbb{T}_{s-t}\right)^{(1-q')\frac{p'}{q'}}\, ds \right)^\frac{1}{p'}.
\end{equation}
To conclude, one needs to show that the above integral in time is finite. From the definition $\mathbb{T}_t:=t^{1/\alpha}\mathbb{M}_t$, we now notice that
\[
\det \mathbb{T}_{s-t}\, = \, (s-t)^{\frac{1}{\alpha}+\sum_{i=1}^nd_i(i-1)} \, = \, (s-t)^{\sum_{i=1}^nd_i\frac{1+\alpha(i-1)}{\alpha}}.
\]
Since it holds that $(1-q')\frac{p'}{q'} =- \frac{p'}{q}$, we can finally conclude that the integral in time in Equation \eqref{INTRO:eq_control_G} is finite if
\[
\bigl(\sum_{i=1}^nd_i\frac{1+\alpha(i-1)}{\alpha}\bigr)\frac{p'}{q} \, < \, 1  \, \Leftrightarrow \, \bigl(\sum_{i=1}^nd_i\frac{1+\alpha(i-1)}{\alpha}\bigr)\frac{1}{q}+\frac{1}{p} \, <\,  1.
\]
The right-hand side inequality is precisely the integrability condition $(\mathscr{C})$ we assumed in Theorem \ref{INTRO:thm:Krylov}.
Intuitively, we can then explain the threshold in $(\mathscr{C})$ as the condition of integrability in time necessary to control the different intrinsic scales associated with the degenerate nature of the dynamics, when considering estimates in $L^p_t-L^q_x$ norm.

We focus now on the point-wise control on the remainder term $\tilde{R}_\varepsilon$ given in \eqref{INTRO:control_pointwise}. From its definition in \eqref{INTRO:eq:NUMBER_PREAL_REG_DENS}, we can split $\tilde{R}_\varepsilon$ as follows:
\begin{align}
\label{INTRO:decomposizione_errore_R}
  \tilde{R}_\varepsilon f(t,x) \, &= \, \int_{t+\varepsilon}^T  \int_{\R^N}(\mathcal{L}_t -\tilde{\mathcal{L}}^{s,y}_t)\tilde{p}^{s,y}(t,s,x,y) f(s,y) \, dy ds \\\notag
   &\qquad \qquad\qquad + \int_{t+\varepsilon}^T  \int_{\R^N}\langle F(t,x)-\tilde{F}^{s,y}_t,D_x\tilde{p}^{s,y}(t,s,x,y)\rangle f(s,y) \, dy ds \\\notag
   &=: \,  \tilde{R}^0_\varepsilon f(t,x)+ \tilde{R}^1_\varepsilon f(t,x), 
\end{align}
where, we recall, the operators $\mathcal{L}_s$ and $\tilde{\mathcal{L}}^{s,y}_s$ are the non-local component of the original infinitesimal generator (cf.\ \eqref{INTRO:eq:def_generator}) and of the frozen proxy one (cf.\ \eqref{INTRO:eq:def_frozen_generator}), respectively. In particular, we highlight that the first term $\tilde{R}^0_\varepsilon f$ in the decomposition  is non-zero only if the diffusion matrix $\sigma(t,\cdot)$ is not constant in space. A bound for the difference between the two generators in $\tilde{R}^0_\varepsilon f$, given by
\begin{multline}
\label{INTRO:eq_per_AC}
(\mathcal{L}_t-\tilde{\mathcal{L}}^{s,y}_t)\tilde{p}^{s,y}(t,s,\cdot,y)(x) \\
= \, \int_{\R^d_0}\left[\tilde{p}^{s,y}(t,s,x+B\sigma(t,x)z,y)-\tilde{p}^{s,y}(t,s,x+B\tilde{\sigma}^{s,y}_tz,y)\right]\nu(dz)
\end{multline}
immediately appears delicate, especially because we cannot directly exploit the regularity of the frozen density $\tilde{p}^{s,y}(t,s,x,y)$ in $x$, since its smoothing effects in the variable $y$ 
will be necessary subsequently in the bound of the external integral (with respect to $y$).
Instead, when the absolutely continuity condition in [\textbf{AC}] holds, we know that
\[\nu(dz) \, = \, Q(z)\frac{g\left(\frac{z}{|z|}\right)}{|z|^{d+\alpha}}dz,\]
for some Lipschitz continuous function $g$ on $\mathbb{S}^{d-1}$. Moreover, the uniform ellipticity condition in [\textbf{UE}] implies that $\det \sigma(t,x) \neq 0$. Assuming without loss of generality that $\det \sigma(t,x)>0$, the changes of variables $\tilde{z}=\sigma(t,x)z$ in the integral within $\mathcal{L}_s$ and $\tilde{z}= \tilde{\sigma}^{s,y}_tz$ in the one for $\tilde{\mathcal{L}}^{s,y}_s $ allow to rewrite the difference between the two infinitesimal generators in the following way:
\[\left(\mathcal{L}_t-\tilde{\mathcal{L}}^{s,y}_t \right)\tilde{p}^{s,y}(t,s,\cdot,y)(x) \, = \, \int_{\R^d_0}\left[\tilde{p}^{s,y}(t,s,x+Bz,y)-\tilde{p}^{s,y}(t,s,x,y) \right]\tilde{H}^{s,y}_{t,x}(z)\, \frac{dz}{|z|^{d+\alpha}},\]
where, for simplicity, we denoted
\[
   \tilde{H}^{s,y}_{t,x}(z) :=  Q(\sigma^{-1}(t,x)z)\frac{g\left(\frac{\sigma^{-1}(t,x)z}{|\sigma^{-1}(t,x)z|}\right)}{\det \sigma(t,x)|\sigma^{-1}(t,x)\frac{z}{|z|}|^{d+\alpha}}
   - Q((\tilde{\sigma}^{s,y}_t)^{-1}z)\frac{g\left(\frac{(\tilde{\sigma}^{s,y}_t)^{-1}z}{|(\tilde{\sigma}^{s,y}_t)^{-1}z|}\right)}{\det \tilde{\sigma}^{s,y}_t|(\tilde{\sigma}^{s,y}_t)^{-1}\frac{z}{|z|}|^{d+\alpha}}
\]
Intuitively, Condition [\textbf{AC}] allows to transfer the error we want to bound on the tempering functions $g,Q$, and then, to exploit their properties. We also mention that in the control of the difference between the infinitesimal generators $\mathcal{L}_t-\tilde{\mathcal{L}}^{s,y}_t$, we will exploit that $\tilde{H}^{ s,y}_{t,x}$ is even, thanks to the symmetry of the Lévy measure $\nu$, as it will allow to add the first order terms necessary for the Taylor expansions, at any cut level.

The bound in \eqref{INTRO:decomposizione_errore_R} on the error term $\tilde{R}^1_\varepsilon f$ shows, at least heuristically, why the Hölder regularity thresholds for $F$, given in Theorem \ref{INTRO:thm:counterexample}, appear natural
Indeed, let us suppose, as in the above theorem, that $x_j\to F_i(t,x)$ is $\beta^j_i$-Hölder continuous, uniformly in time and in the other spacial variables. A typical argument we will have to follow in the estimates of $\tilde{R}^1_\varepsilon f$ focuses on showing that quantities of the following form:
\[\tilde{R}^i_\varepsilon f(t,x) \, := \, \int_{\R^N}\left[F_i(s,x)-F_i(s,\theta_{t,s}(y))\right]D_{x_i}\tilde{p}^{s,y}(t,s,x,y) \, dy\]
generates an integrable singularity in time. Exploiting the properties of the frozen density $\tilde{p}^{s,y}(t,s,x,y)$ in \eqref{INTRO:eq:derivative_prop_of_q}, we can then write that
\[
\begin{split}
|\tilde{R}^i_\varepsilon f (t,x)| \, &= \, \sum_{j=i}^n\int_{\R^N}\frac{\left|\left(x-\theta_{t,s}(y)\right)_j\right|^{\beta^j_i}}{(s-t)^{\frac{1+\alpha(i-1)}{\alpha}}}\frac{\bar{p}(1,\mathbb{T}^{-1}_{s-t}(y-\tilde{m}^{s,y}_{s,t}(x))}{\det \mathbb{T}_{s-t}} \, dy\\
&= \, \sum_{j=i}^n(s-t)^{-\zeta^j_i}\int_{\R^N}\left(\frac{\left|\left(x-\theta_{t,s}(y)\right)_j\right|}{(s-t)^{\frac{1+\alpha(i-1)}{\alpha}}}\right)^{\beta^j_i}\frac{\bar{p}(1,\mathbb T^{-1}_{s-t}(y-\tilde{m}^{s,y}_{s,t}(x))}{\det \mathbb T_{s-t}} \, dy \\
&\le \, C\sum_{j=i}^n(s-t)^{-\zeta^j_i}\int_{\R^N}\left|\mathbb{T}^{-1}_{s-t}\left(x-\theta_{t,s}(y)\right)\right|^{\beta^j_i}\frac{\bar{p}(1,\mathbb T^{-1}_{s-t}(y-\tilde{m}^{s,y}_{s,t}(x))}{\det \mathbb T_{s-t}} \, dy
\end{split}\]
where, for simplicity, we have denoted
\[\zeta^j_i \,:= \, \frac{1+\alpha(i-1)}{\alpha}-\beta^j\frac{1+\alpha(j-1)}{\alpha}.\]
At this point, we firstly notice that the choice of the freezing parameters $(\tau,\xi)=(s,y)$ is necessary to ensure an homogeneity between the difference in the drifts at $x-\theta_{t,s}(y)$ and the argument inside the density at $y-\tilde{m}^{s,y}_{s,t}(x)$. In particular, it holds that
\[y-\tilde{m}^{s,y}_{t,s}(x) \, = \, \tilde{m}^{s,y}_{t,s}(y)-x \, = \, \theta_{t,s}(y)-x.\]
Assuming for the moment that from the smoothing effect of the frozen density  $\bar{p}(t,\cdot)$, it is possible to infer that
\begin{equation}
\label{INTRO:controllo_dubbio}
\int_{\R^N} \left|\mathbb{T}^{-1}(\theta_{t,s}(y)-x)\right|^{\beta^j_i}\frac{\bar{p}(1,\mathbb{T}^{-1}(\theta_{t,s}(y)-x))}{\det \mathbb{T}_{s-t}}\, dy \, < \, +\infty,    
\end{equation}
it immediately follows that
\[
\begin{split}
|\tilde{R}^i_\varepsilon f (t,x)| \, &\le \, C\sum_{j=i}^n(s-t)^{\zeta^j_i}\int_{\R^N}\left|\mathbb{T}^{-1}_{s-t}\left(x-\theta_{t,s}(y)\right)\right|^{\beta^j_i}\frac{\bar{p}(1,\mathbb T^{-1}_{s-t}(x-\theta_{t,s}(y))}{\det \mathbb T_{s-t}} \, dy \\
&\le \, C\sum_{j=i}^n(s-t)^{-\zeta^j_i}.
\end{split}\]
In order to obtain an integrable singularity in time for $\tilde{R}^i_\varepsilon f$, the \emph{natural} thresholds on the Hölder regularity on $F_i$ should be given by:
\begin{equation}
\label{INTRO:soglie_naturali}
    \zeta^j_i <1 \Leftrightarrow \beta^j_i \, > \, \frac{1+\alpha(i-2)}{1+\alpha(j-1)}.
\end{equation}
When $\alpha=2$, i.e.\ in the diffusive context, these thresholds are actually found in \cite{Chaudru:Menozzi17}.

The type of controls obtained above also shows why we will have to assume initially that the indexes $p$, $q$ are large enough in the point-wise bound on the remainder term $\tilde{R}_\varepsilon f$. Indeed, a reasoning similar to the one in \eqref{INTRO:COntrol_G} allows us to control the term $\tilde{R}^1_\varepsilon f$ associated with the drift $F$ in \eqref{INTRO:decomposizione_errore_R} in the following way:
\[|\tilde{R}^1_\varepsilon f(t,x)| \, \le \, C\Vert f \Vert_{L^p_tL^q_x} \sum_{j=i}^n\left(\int_t^T (s-t)^{-\zeta^j_i p'}\left(\det \mathbb{T}_{s-t}\right)^{-\frac{p'}{q}} ds\right)^{\frac{1}{p'}},\]
where, we recall, $p'$ and $q'$ are the conjugate exponents of $p$ and $q$, respectively. Then, the minimal threshold on $p$ and $q$ will be necessary to impose that $p'$ and $q'$ are small enough to ensure that
\[\zeta^j_ip'+\bigl(p'-\frac{p'}{q'}\bigr)\bigl(\sum_{i=1}^nd_i\frac{1+\alpha(i-1)}{\alpha}\bigr)\, <\, 1,\]
and obtain an integrable quantity within the integral in time.

In conclusion, we now briefly explain why we have not managed to obtain the originally wanted result, i.e.\ the weak well-posedness of SDE \eqref{INTRO:eq:SDE} under the natural conditions in \eqref{INTRO:soglie_naturali} on the regularity of drift $F$, but instead we had to assume the same regularity $\beta^j$ for each component $F_i$ along the variable $x_j$. The crucial element precisely lies in being able to prove Equation \eqref{INTRO:controllo_dubbio} for the point-wise estimates of the remainder term $\tilde{R}^1_\varepsilon f$. We immediately point out that for our model, the low regularity of the flow $y\to \theta_{t,s}(y)$ prevents us to directly derive \eqref{INTRO:controllo_dubbio} through the change of variables $\tilde{y} =\theta_{t,s}(y)-x$. Indeed, the most natural approach, exploited for example in \cite{Chaudru:Menozzi17}, is to shift the flow on the variable $x$ through an ``approximate'' Lipschitz property of the form:
\begin{equation}
\label{INTRO:eq:Lipschitz_property}
|\mathbb{T}^{-1}_{s-t}(\theta_{t,s}(y)-x)|\, \asymp (1+ |\mathbb{T}^{-1}_{s-t}(y-\theta_{s,t}(x))|).
\end{equation}
The main problem in our case is that we have not been able to establish, in complete generality, that:
\begin{equation}\label{INTRO:equiv_density_flows}
\overline{p}(1,\mathbb{T}^{-1}_{s-t}(\theta_{t,s}(y)-x))\le C \check{p}(1,\mathbb{T}^{-1}_{s-t}(y-\theta_{s,t}( x)),
\end{equation}
for a density $\check p $ which possesses the same regularising properties of $\bar{p}$. We mention that this difficulty is intrinsically linked to the degenerate and $\alpha$-stable nature of our dynamics. Indeed, a control like in \eqref{INTRO:equiv_density_flows} can be easily obtained in the Gaussian case starting from the explicit expression of the density $\bar{p}$ and the bound in \eqref{INTRO:eq:Lipschitz_property}. To establish a point-wise estimate, as in \eqref{INTRO:equiv_density_flows}, in the stable case, the difficulty precisely lies in being able to obtain a sufficiently precise description of the behaviour of the tails (associated with large jumps) which, as it is well-known, are associated with the geometry of the corresponding spectral measure. In this regard, we cite  Watanabe's work \cite{Watanabe07} in the stable case and Sztonyk's one \cite{Sztonyk10} for an extension to the tempered stable case. In particular, the tricky part appears when considering the behaviour of the Poisson measure (large jumps) associated with the density $\bar{p}$ in an off-diagonal regime, i.e. when $|\mathbb{T}^{-1 }_{s-t}(\theta_{t,s}(y)-x)|> K $. In this case, one could infer from \eqref{INTRO:decomposizione_Ito_Levy}, \eqref{INTRO:proof:control_pM_segnato} and \eqref{INTRO:proof_control_N_segnato}, that:
\begin{align}
\overline{p}(1,\mathbb{T}^{-1}_{s-t}(\theta_{t,s}(y)-&x))
\, \le \, C\int_{\R ^N} \frac{1}{(1+|\mathbb{T}_{s-t}^{-1}(\theta_{t,s}(y)-x)-w|)^{N+3}} P_{\bar N_1}(dw)\notag\\
&\le \, C\int_{0}^1P_{\bar N_{1}}(\{w\in \R^N:(1+|\mathbb{T}_{s-t}^{-1}(\theta_{t,s}(y)-x)-w|)^{-(N+3)} |>u  \})du\notag\\
&\le \, C\int_{0}^1P_{\bar N_1}(B(\mathbb{T}_{s-t}^{-1}(\theta_{t,s}(y)-x), u^{-1/(N+3)})du. \label{INTRO:eq:spiegazione1}
 \end{align}
Recalling from the arguments below Equation \eqref{INTRO:decomposition} that the support of the spectral measure of the process $\{\tilde{S}^{s,y}_{u}\}_{u\ge 0}$ on $\mathbb S^{N-1}$  effectively has dimension $d$, Watanabe in \cite[Lemma $3.1$]{Watanabe07} proved that there exists a constant $C>0$ such that, for every $z$ in $\R^{N}$ and any $r>0$:
\begin{equation}
\label{INTRO:EST_POISSON}
P_{\bar N_1}(B(z,r))\le Cr^{d+1} (1+r^\alpha)|z|^{-(d+1+\alpha)}.
\end{equation}
Said differently, the worst decay in the global bound is given by the size of the support of the associated spectral measure. We also show that these estimates are, in a certain sense, optimal, at least along some directions of the system (cf.\ \cite[Lemma $3.1$]{Watanabe07}). In this regards, see also \cite{Pruitt:Taylor69}. Exploiting Control \eqref{INTRO:EST_POISSON} in the estimate in \eqref{INTRO:eq:spiegazione1}, we would now obtain that
\[
\begin{split}
\overline{p}(1,\mathbb{T}^{-1}_{s-t}(\theta_{t,s}(y)-x))& \, \le \, C|\mathbb{T}_{s-t}^{-1}(\theta_{t,s}(y)-x)|^{-(d+1+\alpha)}\int_0^1u^{-\frac{d+1}{N+3}}
 (1+
 u^{-\frac{\alpha}{N+3}}  ) du
\\
& \le \,  C (1+|\mathbb{T}_{s-t}^{-1}(\theta_{t,s}(y)-x)|)^{-(d+1+\alpha)}.
\end{split}
\]
The approximated Lipschitz property (cf.\  Equation \eqref{INTRO:eq:Lipschitz_property}) would then imply that
\begin{equation}\label{INTRO:BAD_BOUND}
\begin{split}
\overline{p}(1,\mathbb{T}^{-1}_{s-t}(\theta_{t,s}(y)-x)) &\le  C (1+|\mathbb{T}_{s-t}^{-1}(y-\theta_{s,t}(x))|)^{-(d+1+\alpha)} \\
&=: \, C \check{p}(1,\mathbb{T}_{s-t}^{-1}(y-\theta_{s,t}(x))).
\end{split}
\end{equation}
A bound of this type would actually be sufficient for our purposes but it would impose very strong conditions on the dimensions $d$, $n$ of the space to ensure the integrability of the density $\check{p}(t,\cdot)$. This phenomenon also appeared in \cite{Huang:Menozzi16} where the authors limited the analysis to the case $d=1$, $n=3$ to prove the well-posedness of the martingale problem associated with a stochastic dynamics with linear drift and driven by a degenerate multiplicative (isotropic) $\alpha$-stable noise. Finally, let us mention that this difficulty would also arise in the ``classical'' tempered case, i.e.\ if we imposed further conditions on the function $Q$. The advantage in this context would have been to keep the tempering function $Q$ inside the density $\bar{p}$, so as to exploit the advantages of the tempering at infinity, and thus find the same concentration problems as in \eqref {INTRO:BAD_BOUND}. As discussed in \cite[Corollary $6$]{Sztonyk10}, we would then have obtained bounds of the form:
\[\tilde p^{\tau,\xi} (t,s,x,y)\le  C (1+|\mathbb{T}_{s-t}^{-1}(y-\theta_{s,t}(x))|)^{-(d+1+\alpha)}Q\left( |\mathbb M_{s-t}^{-1}(y- \theta_{s,t}(x))|\right)
\]
which clearly improve the integrability in space but at the same time deteriorate the one in time, since it is no longer possible to exploit the $\alpha$-self-similarity of the underlying stable density. This type of difficulty would have appeared even if we had considered only the truncated case, i.e. when $Q(z)=\mathds{1}_{B(0,r_0)}$ for some $r_0>0$. For more details, see for example, \cite{Chen:Kim:Kumagai08} in a non-degenerate context.

To solve this crucial difficulty, we will then follow in Chapter \ref{Chap:Weak_Well-Posedness} an alternative approach. Recalling that the natural change of variables $\tilde{y}:=\mathbb{T}^{-1}_{s-t}(\theta_{t,s}(y)-x) $ in \eqref{INTRO:controllo_dubbio} is not directly possible in our context because the coefficients are not regular enough, we will introduce a regularised flow $\theta_{t,s}^\delta(y)$ associated with a smoothed version of the coefficients. By applying the desired change of variables (with respect to the regularised flow), it will be possible to obtain the estimates in \eqref{INTRO:controllo_dubbio} with respect to $\theta^\delta_{t,s}(y)$, and control the difference between the flows similarly to how we will have already done to establish the approximate Lipschitz condition (cf.\ Equation \eqref{INTRO:eq:Lipschitz_property}). Finally, it will only remain to bound $\det(\nabla \theta_{t,s}^\delta(y))$, uniformly with respect to the mollification parameter. Indeed, it will be this last control to push us to reinforce the expected Hölder regularity thresholds for $F$ in \eqref{INTRO:soglie_naturali} and assume that each component $F_i$, ($i\in \llbracket 2,n\rrbracket$ ) have the same regularity along the variable $x_j$ ($j\in \llbracket 2,n\rrbracket$), uniformly in time and in the other spacial variables (cf.\ Equation \eqref{INTRO:eq:thresholds_beta}).

\subsubsection{Peano counter-examples for uniqueness in law}

We briefly explain now the heuristic reasoning behind the proof for the non-uniqueness result (cf.\ Theorem \eqref{INTRO:thm:counterexample}).
The idea is to adapt the Peano counter-examples shown in \eqref{INTRO:eq:solution_Peano} and then exploited as well in \cite{Chaudru:Menozzi17}, to our degenerate Lévy framework. If for example, one wants to test the threshold $\beta^j_i$ associated with the critical Hölder exponent for the $i$-th component of the drift $F$ with respect to the variable $x_j$ (with $j\ge i>1$),  we are going to consider the following model (with $d_1=\dots=d_n=1$ and $N=n$):
\begin{equation}
\label{INTRO:eq:SDE_for_counter}
\begin{cases}
dX^1_t \, =\, dZ_t, &\mbox{ if } k=1;\\
dX^k_t \, =\, X^{k-1}_tdt, &\mbox{ if } k\in\llbracket 2,i-1\rrbracket;\\
dX^i_t \, =\, X^{i-1}_tdt+\text{sgn}(X^j_t)|X^j_t|^{\beta^j_i} dt, &\mbox{ if } k=i;\\
dX^k_t \, =\, X^{k-1}_tdt, &\mbox{ if } k\in\llbracket i+1,n\rrbracket,
\end{cases}
\end{equation}
where $\{Z_t\}_{t\ge 0}$ is a real valued, symmetric, $\alpha$-stable process. It is not difficult to notice that Equation \eqref{INTRO:eq:SDE_for_counter} can be rewritten in the form of \eqref{INTRO:eq:SDE} by choosing $\sigma=1$ and $G(t,x)=Ax+e_i\text{sgn}(x_j)|x_j|^{\beta^j_i}$, where 
$e_i$ is the $i$-th element of the canonical basis on $\R^N$ and $A$ is the matrix in $\R^N\otimes \R^N$ given by:
\[
A \, := \, \begin{pmatrix}
               0& \dots         & \dots         & \dots     & 0 \\
              1       & 0 & \dots         & \dots     & 0\\
               0 & 1       & \ddots & \ddots     & \vdots \\
               \vdots        & \ddots        & \ddots        & \ddots    & \vdots        \\
               0 & \dots         & 0 & 1 & 0
             \end{pmatrix}.
\]
If we now focus on the $i$-th component of Equation \eqref{INTRO:eq:SDE_for_counter}, 
it can be rewritten in an integral form as:
\begin{equation}
\label{INTRO:eq:SDE_for_counter_integral}
X^j_t \, = \, \int_0^t\text{sgn}\bigl(I^{j-i}_t(X^j)\bigr)\bigl{\vert}I^{j-i}_t(X^j)\bigr{\vert}^{\beta^j_i} dt + I^{i-1}_t(Z), \quad t\ge 0,
\end{equation}
where, for any càdlàg path $y\colon [0,\infty)\to \R$, the notation $I^k_t(y)$ denotes $k$-th times iterated integral at time $t$. 
As already explained above, for a regularisation by noise phenomenon to occur it is necessary that, at least in small time, the mean fluctuations of the random perturbation dominate the irregularity of the deterministic drift. More precisely for our model, we can compare the fluctuations of the noise $I^{i-1}_t(Z)$, of order $i-1+\frac{1}{\alpha}$, with the extremal deterministic solutions obtained from the dynamics without noise (i.e.\ $B=0$ in the Equation above). Thus, we will have that
\[t^{i-1+\frac{1}{\alpha}} \, > \, t^{\frac{(j-1){\beta^j_i}-1}{1-\beta^j_i}}.\]
Since such condition has to hold for a small time $t$, we can then conclude that the condition
\[i-1+\frac{1}{\alpha} \, < \, \frac{(j-1)\beta^j_i-1}{1-\beta^j_i} \,\, \Leftrightarrow \,\, \beta^j_i \, > \, \frac{1+\alpha(i-1)}{1+\alpha(j-2)}\]
is the heuristic threshold ensuring that the noise dominates the deterministic drift so that a regularisation by noise phenomenon occurs.

\setcounter{equation}{0}

\section{On the constants in Schauder and Sobolev estimates for degenerate Kolmogorov operators}
\fancyhead[RO]{Section \thesection. On the constants in Schauder and Sobolev estimates}
\label{Sec:INTRO:Sulle_Costanti_Ottimali}

We now briefly present the main results of Chapter \ref{Chap:About_Sharp_Constant}. This work, written in collaboration with my PhD supervisors, Prof.\ Stéphane Menozzi and Prof.\ Enrico Priola, has recently appeared in pre-publication (cf.\ \cite{Marino:Menozzi:Priola21}). We are interested here in studying the effects of a second order perturbation on a degenerate diffusive Ornstein-Uhlenbeck operator. In particular, we want to determine how the constants of some ``known'' estimates for this class of operators, such as the Schauder estimates presented in Sections \ref{Sec:INTRO:Stime_Schauder_stabili} and \ref{Sec:INTRO:Stime_Schauder_Levy}, actually depend on the perturbation. Our method of proof will be based on a suitable space transformation, which allows to cancel the first order transport term, and on the perturbation method through Poisson processes, introduced in \cite{Krylov:Priola17} and adapted to our context. More precisely, fixed a positive integer $N$, we consider the following family of diffusive Ornstein-Uhlenbeck operators: 
\begin{equation}
\label{INTRO:DEF_OU_OP_PROXY}
L^\text{ou} \, := \,  
{\rm Tr}(B D^2_z)
+ \langle A z , D_z \rangle,\quad \text{ on }\R^N,
\end{equation}
where $\langle \cdot,\cdot\rangle$ denotes the inner product on $\R^N$ and $A,B$ are two matrices in $\R^N\otimes \R^N$ such that $B$ is symmetric.

We will assume, as in Chapter \ref{Chap:Schauder_Estimates_Levy}, that the two matrices $A$, $B$ satisfy the Kalman rank condition ensuring the hypoellipticity of the system:
\begin{description}
\item[\textbf{[K]}] there exists a non-negative integer $n$ such that
\[
{\rm rank} [ B,AB,\cdots,
  A^{n-1}B
] \, = \, N,
\]
where, we recall, $[B,AB,\dots,A^{n-1}B]$ is the matrix in $\R^N\otimes \R^{Nn}$ whose blocks are $B,AB,\cdots A^{n-1}B$.
\end{description}
As already mentioned in Section \ref{Sec:INTRO:Stime_Schauder_Levy}, Condition [\textbf{K}] allows to decompose the space $\R^N$ according to the image space obtained by the successive iterations of the commutator between $ A$ and $B$ (cf.\ Equation \eqref{INTRO:eq:iterazioni_commutatori}). Assuming that $\text{rank}(B)=d$ for some $d>0$, we know in particular that there exist non-negative integers $\{d_1,\dots,d_n\}$ such that $d_1=d$, $\sum_{i=1}^n d_i=N$ and the two matrices $A$, $B$ can be rewritten, after a possible change of variables, in the following, more explicit, form (cf.\ Equation \eqref{INTRO:eq:Lancon_Pol}):
\[B \, = \,
    \begin{pmatrix}
        B_0 & 0 & \dots & 0   \\
        0  &  0 & \ddots & \vdots    \\
        \vdots &  \ddots & \ddots & \vdots  \\
        0 & \dots & \dots & 0
    \end{pmatrix}
\,\, \text{ and } \,\, A \, = \,
    \begin{pmatrix}
        \ast   & \ast  & \dots  & \dots  & \ast   \\
         A_2   & \ast  & \ddots & \ddots  & \vdots   \\
        0      & A_3   & \ast  & \ddots & \vdots \\
        \vdots &\ddots & \ddots& \ddots & \ast \\
        0      & \dots & 0     & A_n    & \ast
    \end{pmatrix}
\]
where $B_0$ is a non-degenerate matrix in $\R^{d}\otimes \R^{d}$, $A_i$ is a matrix in $\R^{d_i}\otimes \R^{d_{i-1}}$ such that
$\text{rank}(A_i)=d_i$ for any $i$ in $\llbracket 2,n\rrbracket$ and the elements $\ast$ may be non-zero.

We are then interested in the solutions to the following Cauchy problem:
\begin{equation}
\label{INTRO:eq:OU_initial:intro}
\begin{cases}
 \partial_tu(t,z) \, = \, L^{\text{ou}}u(t,z) +f(t,z) &\mbox{ on } (0,T)\times \R^N;\\
u(0,z) \, = \, 0, &\mbox{ on } \R^N.
\end{cases}
\end{equation}
We will show the existence and uniqueness of bounded continuous solutions to Equation \eqref{INTRO:eq:OU_initial:intro} assuming, as in \cite{Krylov:Priola17}, that the source $f$ belongs to the space $ B_b\left(0,T;C^\infty_0(\R^N)\right)$. This space can be understood as the family of functions that are measurable, bounded in time and smooth, compactly supported in space, uniformly in time. For a precise definition of such a space, we refer the reader to Section \ref{SEC_NOT} of Chapter \ref{Chap:About_Sharp_Constant}. We highlight that we could not consider the source $f$ in a more ``usual'' class of functions, such as $C^\infty_c([0,T]\times \R^N)$, since the perturbative method through Poisson processes in \cite{Krylov:Priola17}  which we will also exploit, will require to consider sources $f$ in \eqref{INTRO:eq:OU_initial:intro} that are possibly discontinuous in time (cf.\ Section $2 $ in \cite{Krylov:Priola17}). Due to the low regularity in time for the source $f$, the system \eqref{INTRO:eq:OU_initial:intro} will be understood only in an integral sense, i.e. a bounded, continuous function $u\colon [0,T]\times \ R^N\to \R$ will be a solution to Equation \eqref{INTRO:eq:OU_initial:intro} if $u(t,\cdot)$ belongs to $C^2(\R^N)$ for each fixed $t$ and 
\begin{equation}
\label{INTRO:eq:sol_integrale}
u(t,x) \, = \, \int_0^t \left[\text{Tr}(BD^2_zu(s,z))+\langle Az,D_zu(s,z)\rangle +f(s,z)\right] \, ds.    
\end{equation}

Given a continuous function $t\mapsto S(t)$ such that $S(t)$ is a  non-negative definite matrix in $\R^N\otimes\R^N$, we will be interested in the second order perturbation given by $S(t)$ to the Ornstein-Uhlenbeck operator, i.e. in the following operator:
\[
L_t^{\text{ou},S} \, :=\, L^\text{ou} +
{\rm Tr}(S(t) D^2_z ) \, = \, {\rm Tr}\left(\left[B+S(t)\right] D^2_z\right)+ \langle A z , D_z \rangle, \quad \text{ on } \R^N,
\]
and in the associated \emph{perturbed} Cauchy problem:
\begin{equation}
\label{INTRO:eq:OU_PERT0}
\begin{cases}
    \partial_t u_S(t,z) \, = \,  {\rm Tr}\left(\left[B+S(t)\right] D^2_z u_S(t,z)\right)+ \langle A z , D_z  u_S(t,z)\rangle +f(t,z);\\
   u_S(0,z)\, = \, 0.
\end{cases}
\end{equation}
In this work, we will mainly focus on two families of estimates for the solutions $u$ to Cauchy problem \eqref{INTRO:eq:OU_initial:intro}: Sobolev-type estimates, i.e. controls in the $L^p$ norm for the first non-degenerate component of the solution gradient and Schauder estimates with respect to Hölder spaces with multi-index of regularity, already met in Sections \ref{Sec:INTRO:Stime_Schauder_stabili} and \ref{Sec:INTRO:Stime_Schauder_Levy}. In particular, Bramanti \emph{et al.} showed in \cite[Theorem $3$]{Bramanti:Cupini:Lanconelli:Priola10} that for any $p$ in $(1,\infty)$, there exists a constant $C_p>0$, independent of $f$, such that
\begin{equation}\label{INTRO:de3}
\| B^{1/2} D^2 u \, B^{1/2} \|_{L^p ((0,T) \times \R^N)} \, \le \, C_p  \| f    \|_{L^p ((0,T) \times \R^N)}.
\end{equation}
However, we highlight that the estimates in \eqref{INTRO:de3} are actually obtained in \cite{Bramanti:Cupini:Lanconelli:Priola10} assuming that the source $f$ is smooth in space and time. Through some explicit properties on the  underlying Gaussian heat kernel, we will show in Section \ref{KH_OU} of Chapter \ref{Chap:About_Sharp_Constant} that indeed such estimates can be extended to the more general source $f$ we consider here.

Under the Kalman condition in \textbf{[K]}, Lunardi (cf.\ \cite[Theorem $1.2$]{Lunardi97}) showed instead that for any $\beta$ in $(0,1)$, there exists a constant $C_\beta$, independent from $f$, such that
\begin{equation}\label{INTRO:SCHAU_OU}
\Vert u \Vert_{L^\infty((0,T),C^{2+\beta}_{b,d})} \, \le \, C_\beta \Vert f \Vert_{L^\infty((0,T),C^{\beta}_{b,d})},
\end{equation}
where, we recall, the anisotropic Hölder spaces $C^{\gamma}_{b,d}(\R^N)$ are defined exactly as in Section \ref{Sec:INTRO:Stime_Schauder_Levy}.

We can now resume the main results of Chapter \ref{Chap:About_Sharp_Constant} in the following theorem:
\begin{teorema}
\label{thm:costanti}
Let $f$ be in $B_b\left(0,T;C^\infty_0(\R^N)\right) $. Then, there exists a unique integral solution $u_S$ to Cauchy problem \eqref{INTRO:eq:OU_PERT0} such that, for any $p$ in $(1,+\infty)$ and $\beta$ in $(0,1)$ it holds that
\begin{align}\label{INTRO:de4}
\| B^{1/2} D^2 u_S \, B^{1/2} \|_{L^p ((0,T) \times \R^N)} \, &\le \,C_p  \| f
\|_{L^p ((0,T) \times \R^N)}\\
\Vert u_S \Vert_{L^\infty(C^{2+\beta}_{b,d})} \, &\le \, C_\beta \Vert f \Vert_{L^\infty(C^{\beta}_{b,d})},
\label{INTRO:de41}
\end{align}
with the same constants $C_p$, $C_\beta$ appeared in in \eqref{INTRO:de3} and \eqref{INTRO:SCHAU_OU}, respectively. In particular, the constants $C_p$, $C_\beta$ does not depend on the matrix $S(t)$.
\end{teorema}

Beyond the preservation of constants property shown in Theorem \ref{thm:costanti} above, the Sobolev estimates in \eqref{INTRO:de4} seem, to the best of our knowledge, to be new for an operator like $L_t^{\text{ou},S}$ and of independent interest. We also mention the recent work by Fornaro \textit{et al.} \cite{Fornaro:Metafune:Pallara:Schnaubelt21} in which they outline a complete description of the spectrum of hypoelliptic Ornstein-Uhlenbeck operators in $L^p$ spaces. We also point out that if we consider a time independent matrix $S$, analogous results can also be obtained for the corresponding elliptic estimates, following the method shown in \cite[Corollary $3.5$]{Krylov:Priola17}. Finally, we remark that in Chapter \ref{Chap:About_Sharp_Constant} we will also show that more general $L^p$ estimates, which also consider degenerate directions, are independent of second order perturbations, even if only for matrices $A$ that are invariant under dilations (cf.\ Equation \eqref{INTRO:eq:matrix_A_sub-diag} in Section \ref{Sec:INTRO:Stime_Schauder_stabili}). For more details on this, see Section \ref{ESTENSIONI} in Chapter \ref{Chap:About_Sharp_Constant}.

\subsection{Sketch of the proof}
To give the reader an idea of the method we used, we now briefly illustrate the main steps in the proof for the Sobolev estimates in Equation \eqref{INTRO:de4}. As already mentioned at the beginning of this section, a fundamental tool will be a known (cf.\ \cite{Daprato:Lunardi95}) transformation of the space $\R^N$ (at any fixed time) which will allow to get rid of the drift term $\langle Az,D_zu\rangle$ in \eqref{eq:OU_initial:intro}. More in detail, given a bounded solution $u$ to Cauchy problem \eqref{eq:OU_initial:intro}, we will introduce the function $v\colon[0,T]\times\R^N$ given by
\[v(t,z) \,:=  \, u(t, e^{-tA} z).\]  
Indeed, recalling that $u$ is a solution to \eqref{eq:OU_initial:intro} and noting that $u(t,z) = v(t, e^{tA} z)$, it is not difficult to check that
\[
 \begin{split}
 f (t,z)\, &= \,  \partial_t u(t,z)  - L^{\text{ou}} u(t,z)\\
 &= \, v_t(t,e^{tA} z)+\langle D v(t, e^{tA} z) , A e^{tA} z\rangle  - {\rm Tr} \big( e^{tA}Be^{tA^*} D^2  v(t, e^{tA} z) \big) \\
 &\quad- \langle D v (t, e^{tA} z) , A e^{tA} z\rangle\\
&= \, v_t(t,e^{tA} z)    - {\rm Tr} \big( e^{tA} Be^{tA^*} D^2 v(t, e^{tA} z) \big).
 \end{split}\]
for any $(t,z)$ in $(0,T)\times \R^N$. Denoting for simplicity $\tilde{f}(t,z):= f(t,e^{-tA}z)$, it immediately follows from the above calculations that $v$ is then a solution to the following Cauchy problem:
 \[
 \begin{cases}
 \partial_tv(t, z)   \, = \,  {\rm Tr} \left(e^{tA}Be^{tA^*} D^2 v(t, z) \right) + \tilde{f}(t,z) &\mbox{ on }(0,T)\times \R^N;\\
 v(0,z) \, = \, 0 &\mbox{ on } \R^N.
 \end{cases}
\]
Furthermore, the estimates in \eqref{INTRO:de3} can be rewritten with respect to $v$ as
\begin{equation} \label{INTRO:1}
\| B^{1/2} e^{- t A^*} D^2 v (t, e^{tA} \cdot ) \, e^{tA} B^{1/2} \|_{L^p((0,T)\times \R^N)} \, \le \,  C_p  \|  \tilde f(t, e^{tA}  \cdot )  \|_{L^p((0,T)\times \R^N)}.
\end{equation}
Through a change of variables in the integral and denoting by $L^p((0,T)\times \R^N, m)$ the usual $L^p$ space but with respect to the measure $m$ given by
\[m(dt,dx) \, :=\, {\rm det}(e^{-At}) dt dx,\]
it is immediate to notice that the estimates in \eqref{INTRO:1} are equivalent to:
\begin{equation}\label{INTRO:2}
\| B^{1/2} e^{tA^*} D^2 v (t,  \cdot ) \, e^{tA} B^{1/2} \|_{L^p((0,T)\times \R^N,m)} \le C_p  \|  \tilde f  \|_{L^p((0,T)\times \R^N, m)}.
\end{equation}

We can now focus on the Cauchy problem perturbed by the matrix $S(t)$:
\begin{equation} \label{INTRO:d2}
\begin{cases}
 \partial_t w(t, z)    + {\rm Tr} \big( e^{tA} B
 e^{tA^*} D^2 w(t, z) \big) + {\rm Tr} \big( e^{tA} S(t)
 e^{tA^*} D^2 w(t, z) \big) \, = \, \tilde f(t,z);\\
 w(0,z)\, = \, 0.
 \end{cases}
\end{equation}
Through an additional probabilistic reasoning, we will show in particular that there exists a unique solution $w\colon[0,T]\times \R^N\to \R$ to the above problem.

Adapting now some of the arguments presented in \cite{Krylov:Priola17}, we will be able to deduce that the $L^p$-norm estimates in \eqref{INTRO:2} also hold for the solution $w$ to the perturbed Cauchy problem \eqref{INTRO:d2}, independently from the matrix $S(t)$. More precisely, we will get that the estimate
\begin{equation} \label{INTRO:11}
\| B^{1/2} e^{tA^*} D^2 w (t, \cdot ) \, e^{tA} B^{1/2} \|_{L^p((0,T)\times \R^N, {m})} \le C_p \|  \tilde f(t,  \cdot )  \|_{L^p((0,T)\times \R^N, {m})}
\end{equation}
holds again with the \emph{same} constant $C_p$ appearing in \eqref{INTRO:2}. The crucial element in the method of proof of \cite{Krylov:Priola17} consists in introducing a small random perturbation on the source $f$ through a suitable Poisson type process and study the properties associated with the corresponding equation. Taking the expected value in the integral formulation of the equation, the contributions associated with the process jumps generate, for an appropriate intensity of the underlying Poisson process, a finite difference operator. Importantly, the initial estimates remain preserved for the system solved by the expected value of the solution and which also involves the finite difference operator. Arguments of compactness finally allow us to conclude that the initial estimates also hold at the limit, when we exchange the finite difference operator with the corresponding second order differential operator.

Finally, we will then have to go back to our original Ornstein-Uhlenbeck model, applying the inverse transformation with respect to the spacial variable. In particular, we will introduce $\tilde u(t,z) := w(t, e^{tA} z)$ which solves, by definition, the following equation: 
\[
\begin{cases}
 \partial_t \tilde u(t,z) + L_t^{\text{ou}, S} \tilde u(t,z)\, = \,  f  (t,z), &\mbox{ on } (0,T)\times \R^N,\\
 \tilde u(0,z)=0, &\mbox{ on } \R^N.
\end{cases}
\]
Exploiting the following identity:
\[D^2 w(t, \cdot) = D^2 [\tilde u(t, e^{-tA} \, \cdot )] = e^{-tA^*}D^2 \tilde u(t, e^{-tA} \cdot )e^{-tA},\]
we will finally conclude from \eqref{INTRO:11} that the following Sobolev estimates hold:
\[
\| B^{1/2} D^2 \tilde u \, B^{1/2} \|_{L^p ((0,T) \times \R^N)} \, \le \,
C_p  \|  f \|_{L^p ((0,T) \times \R^N)} .
\]
Summarising, the arguments shown above actually allow one to construct a solution $\tilde u$ to Cauchy problem \eqref{eq:OU_PERT0} which satisfies the Sobolev estimates in \eqref{de4} with the \emph{same} constant $C_p$ appearing in the analogous estimates for the \emph{proxy} operator  $L^{\text{ou}}$. Then, noting that the maximum principle also holds for the perturbed Ornstein-Uhlenbeck operator $L_t^{\text{ou},S}$, we can finally show the uniqueness of such solution $\tilde u $.

In conclusion, we highlight that to apply the perturbation method briefly summarised above, only a few specific properties are in fact required on the underlying semi-norms. Intuitively, the reasoning presented in \cite{Krylov:Priola17} only exploits the invariance under translations of the semi-norms and a kind of commutativity property between the norms (or a function of the norm as in the case $L^p$) and the expected value operator. Indeed, it seems natural that this approach could be extended to a much more general class of estimates on other functional spaces, such as, for example, Besov spaces (cf.\ Equation \eqref{INTRO:alpha-thermic_Characterization}). Furthermore, this type of controls seems to be promising for a more detailed analysis of the weak well-posedness for some related SDEs. These aspects will be researched in the near future.

\setcounter{equation}{0}
\selectlanguage{english}
\chapter{Schauder estimates for degenerate stable Kolmogorov equations}

\fancyhead[LE]{Chapter \thechapter. Schauder estimates for non-linear stable equations}

\label{Chap:Schauder_estimates_Stable}

\paragraph{Abstract:} We provide here global Schauder-type estimates for a chain of integro-partial differential equations (IPDE) driven by a degenerate
stable Ornstein-Uhlenbeck operator possibly perturbed by a deterministic drift, when the coefficients lie in some suitable anisotropic
H\"older spaces. Our approach mainly relies on a perturbative method based on forward parametrix expansions and, due to the low
regularizing properties on the degenerate variables and to some integrability constraints linked to the stability index, it also exploits
duality results between appropriate Besov Spaces. In particular, our method also applies in some super-critical cases. Thanks to these
estimates, we show in addition the well-posedness of the considered IPDE in a suitable functional space.

\section{Introduction}
\fancyhead[RO]{Section \thesection. Introduction}
For a fixed time horizon $T>0$ and two integers $n,d$ in $\N$, we are interested in proving global Schauder estimates for the following
parabolic integro-partial differential equation (IPDE):
\begin{equation}
\label{Degenerate_Stable_PDE}
\begin{cases}
   \partial_t u(t,x) + \langle A x + F(t,x), D_{x}u(t,x)\rangle +
\mathcal{L}_\alpha u(t,x) \, = \, -f(t,x) & \mbox{on } [0,T]\times \R^{nd}; \\
    u(T,x) \, = \, u_T(x) & \mbox{on }\R^{nd}.
  \end{cases}
\end{equation}
where $x:=(x_1,\dots,x_n)$ is in $\R^{nd}$ with each $x_i$ in $\R^d$ and $\langle \cdot,\cdot \rangle$ represents the
inner product on  $\R^{nd}$. We consider a symmetric, $\alpha$-stable operator $\mathcal{L}_\alpha$ acting non-degenerately only on the
first $d$ variables and a matrix $A$ in $\R^{{nd}} \otimes\R^{{nd}}$ with the following sub-diagonal structure:
\begin{equation}\label{eq:def_matrix_A_stable}
A \, := \, \begin{pmatrix}
               0_{d\times d} & \dots         & \dots         & \dots     & 0_{d\times d} \\
               A_{2,1}       & 0_{d\times d} & \dots         & \dots     & 0_{d\times d} \\
               0_{d\times d} & A_{3,2}       & 0_{d\times d} & \dots     & 0_{d\times d} \\
               \vdots        & \ddots        & \ddots        & \vdots    & \vdots        \\
               0_{d\times d} & \dots         & 0_{d\times d} & A_{n,n-1} & 0_{d\times d}
             \end{pmatrix}.
\end{equation}
We will assume moreover that it satisfies a H\"ormander-like condition, allowing the smoothing effect of $\mathcal{L}_\alpha$ to propagate into the system. \newline
Above, the source $f\colon [0,T]\times \R^{nd} \to \R$ and the terminal condition $u_T\colon \R^{nd}\to \R$ are assumed to be bounded and to
belong to some suitable anisotropic H\"older space. \newline
The additional drift term $F(t,x)=\bigl(F_1(t,x),\dots,F_n(t,x)\bigr)$ can be seen as a perturbation of the
Ornstein-Ulhenbeck operator $\mathcal{L}_\alpha+\langle Ax, D_{x}\rangle$ and it has  structure "compatible" with $A$, i.e.\ at level $i$,
it depends only on the super diagonal entries:
\[F_i(t,x) \, := \, F_i(t,x_i,\dots,x_n).\]
It may be unbounded but we assume it to be H\"older continuous with an index depending on the level of the chain.

\paragraph{Related results.}
A large literature on the topic of Schauder estimates in the $\alpha$-stable non-local framework has been developed in the recent years (see e.g.\ Lunardi and R\"ockner \cite{Lunardi:Rockner19} for an overview of the field), mainly
in the non-degenerate setting and assuming that $\alpha\ge 1$, the so called sub-critical case. We mention for instance the stable-like setting, corresponding to time-inhomogeneous operators of the form
\begin{multline} 
\label{eq_Mik_Prag}
\bar{L}_t\phi(x) \, = \, \int_{\R^{nd}}\bigl[\phi(x+y)-\phi(x) - \mathds{1}_{1\le \alpha <2}\langle y, D_{x} \rangle
\bigr] m(t,x,y) \frac{dy}{\vert y \vert^{d+\alpha}} \\
+ \mathds{1}_{1\le \alpha <2}\langle F(t,x), D_{x}u(t,
x)\rangle
\end{multline}
where the diffusion coefficient $m$ is bounded from above and below, H\"older continuous in the spatial variable $x$ and even in $y$ if $\alpha=1$.
Under these conditions and assuming the drift $F$ to be bounded and H\"older continuous in space, Mikulevicius and Pragarauskas in \cite{Mikulevicius:Pragarauskas14} obtained parabolic Schauder type bounds on the whole space and derived
from those estimates the well-posedness of the corresponding martingale problem. We notice however that for the super-critical case
(when $\alpha<1$), the drift term in \eqref{eq_Mik_Prag} is set to zero. This is mainly due to the fact that in the super-critical case,
$\mathcal{L}_\alpha$ is of order $\alpha$ (in the Fourier space) and does not dominate the drift term $F$ which is roughly speaking of order one.
\newline
In the non-degenerate, driftless framework (i.e.\ when $Ax +F=0$ and $n=1$ in \eqref{Degenerate_Stable_PDE}), Bass \cite{Bass09} was the first to derive
elliptic Schauder estimates for stable like operators. We can refer as well to the recent work of Imbert and collaborators
\cite{Imbert:Jin:Shvydkoy18} concerning Schauder estimates for stable-like operator \eqref{eq_Mik_Prag} with $\alpha=1$ and some related applications to non-local
Burgers equations. Eventually, still in the driftless case, Ros-Oton and Serra worked in \cite{Ros-Oton:Serra16} for interior and boundary
elliptic-regularity in a general, symmetric $\alpha$-stable setting, assuming that the L\'evy measure $\nu_\alpha$ associated with $\mathcal{L}_\alpha$
writes in polar coordinates $y=\rho s$, $(\rho,s)\in [0,\infty)\times \mathbb{S}^{d-1}$ as
\[\nu_\alpha(dy) \, = \, \tilde{\mu}(ds)\frac{d\rho}{\rho^{1+\alpha}}\]
where $\tilde{\mu}$ is a non-degenerate, symmetric measure on the sphere $\mathbb{S}^{d-1}$. Related to the above, we can mention also the
associated work of Fernandez-Real and Ros-Oton \cite{Fernanadez:Ros-Oton17} for parabolic equations.

In the elliptic setting, when $\alpha\in [1,2)$ and $\mathcal{L}_\alpha$ is a non-degenerate, symmetric $\alpha$-stable operator and for bounded H\"older
drifts, global Schauder estimates were obtained by Priola in \cite{Priola12} or in \cite{priola18} for respective applications to
the strong well-posedness and Davie's uniqueness for the corresponding SDE. We notice furthermore that in the
sub-critical case, elliptic Schauder estimates can be proven for more general, translation invariant, L\'evy-type generators for
 following \cite{priola18} (see Section $6$, and Remark $5$ therein).

In the super-critical case, parabolic Schauder estimates were established by Chaudru de Raynal, Menozzi and Priola in
\cite{Chaudru:Menozzi:Priola19} under similar assumptions to \cite{Ros-Oton:Serra16}. An existence result is also provided therein. \newline
We  mention as well the work of Zhang and Zhao \cite{Zhang:Zhao18} who address through probabilistic arguments the parabolic
Dirichlet problem for stable-like operators of the form \eqref{eq_Mik_Prag} with a non-trivial bounded drift, i.e.\ getting rid of the
indicator function for the drift. They also obtain interior Schauder estimates and some boundary decay estimates (see e.g.\ Theorem $1.5$
therein).

As we have seen, most of the literature is focused on the non-degenerate case. In the degenerate diffusive setting,  Lunardi
\cite{Lunardi97} was the first one to prove Schauder estimates for linear Kolmogorov equations under weak H\"ormander assumptions,
exploiting anisotropic H\"older spaces (where the H\"older index depends on the variable considered), in order exactly to control the
multiple scales appearing in the different directions, due to the degeneracy of the system.\newline
After, in \cite{Lorenzi05} and \cite{Priola09}, the authors established Schauder-like estimates for hypoelliptic  Kolmogorov equations
driven by partially nonlinear smooth drifts.
On the other hand, let us also mention \cite{Chaudru:Honore:Menozzi18_Sharp} where the authors first establish Schauder estimates
for nonlinear Kolmogorov equations under some weak H\"ormander-type assumption. Their method is based on a perturbative approach through
proxies that we here adapt and exploit.
In the degenerate, stable setting, we have to refer also to a recent work of Zhang and collaborators \cite{Hao:Wu:Zhang19} who show Schauder
estimates for the degenerate kinetic  dynamics ($n=2$ above) extending a method based on Littlewood-Paley decompositions already used in
other works by Zhang (see e.g. \cite{Zhang:Zhao18}), to the degenerate, multi-scaled framework.  Even with different approaches and frameworks, we consider here a
generic $d$-level chain and we exploit thermic characterizations of Besov norms, our and their works bring to the same results in the
intersecting cases, at least to the best of our knowledge.
About a different but correlated argument, we mention that the L$^p$-maximal regularity for degenerate non-local Kolmogorov equations with
constant coefficients was also obtained in \cite{Chen:Zhang19} for the kinetic dynamics ($n=2$ above) and in \cite{Huang:Menozzi:Priola19}
for the general $n$-levels chain.

In the diffusive setting, Equation \eqref{Degenerate_Stable_PDE} appears naturally as a microscopic model for heat diffusion phenomena (see
\cite{Rey-Bellet:Thomas00}) or, in the kinetic case ($n=2$), it can be naturally associated with speed/position (or Hamiltonian)
dynamics where the speed component is noisy. It can be found in many fields of application from physics to finance, see for example
\cite{Herau:Nier04} or \cite{Barucci:Polidoro:Vespri01}. When noised by stable processes, it can be used to model the appearance
 of turbulence (cf. \cite{Cushman:Park:Kleinfelter:Moroni05}) or some abnormal diffusion phenomena. \newline
Moreover, the Schauder estimates will be a fundamental first step in order to study the weak and strong well-posedness for the following
stochastic differential equation (SDE):
\begin{equation}\label{SDE_Associated}
\begin{cases}
  dX^1_t \, = \,  F_1(t,X^1_t,\dots,X^n_t)dt +dZ_t \\
  dX^2_t \, = \, A_{2,1}X^1_t+F_2(t,X^2_t,\dots,X^n_t)dt \\
 \vdots \\
  dX^n_t \, = \, A_{n,n-1}X^{n-1}_t+F_n(t,X^n_t)dt
\end{cases}
\end{equation}
where $Z_t$ is a symmetric, $\R^d$-valued $\alpha$-stable process with non-degenerate L\'evy measure $\nu_\alpha$ on some filtered
probability space $(\Omega,(\mathcal{F}_t)_{t\ge0},\mathbb{P})$. The complete operator $\mathcal{L}_\alpha + \langle A x + F(t,x),
D_{x}\rangle$ then corresponds to the infinitesimal generator of the process $\{ X_t\}_{t\ge 0}$, solution of Equation
\eqref{SDE_Associated}.

\paragraph{Mathematical outline.} In this work, we will establish global Schauder estimates for the solution of the IPDE
\eqref{Degenerate_Stable_PDE} exploiting the perturbative approach firstly introduced in \cite{Chaudru:Honore:Menozzi18_Sharp} to derive
such estimates for degenerate Kolmogorov equations. Roughly speaking, the idea is to perform a first order parametrix expansion, such as a
Duhamel-type representation, to a solution of the IPDE \eqref{Degenerate_Stable_PDE} around a suitable proxy. The main idea behind consists
in exploiting this easier framework in order to
subsequently obtain a tractable control on the error expansion. When applying such a strategy, we basically have two ways to proceed.\newline
On the one hand, one can adopt a backward parametrix approach, as introduced by McKean and Singer \cite{Mckean:Singer67} in the
non-degenerate, diffusive setting. This technique has been extended to the degenerate Brownian case involving unbounded
perturbation, and successfully exploited for handling the corresponding martingale problem in \cite{Chaudru:Menozzi17}.
Anyway, this approach does not seem very adapted to our framework especially because it does not allow to deal easily with point-wise
gradient estimates which will, at least along the non-degenerate variable $x_1$,  be fundamental to establish our result. \newline
On the other hand, the so-called forward parametrix approach has been successfully used by Friedman \cite{book:Friedman08} or Il'in et al.\ \cite{ilcprimein:Kalavsnikov:Oleuinik62} in the non-degenerate, diffusive
setting to obtain point-wise bounds on the fundamental solution and its derivatives for the corresponding heat-type equation or in
\cite{Chaudru17} to derive strong uniqueness for the associated SDE \eqref{SDE_Associated} (i.e. $n=2$ with the previous
notations). Especially, this approach is better tailored to exploit cancellation techniques that are crucial when derivatives come in,
 as opposed to the backward one.

The main difficulties to overcome in order to prove Schauder estimates in our framework will be linked to the degeneracy of the operator
$\mathcal{L}_\alpha$ that acts only on the first $d$ variables, as well as the unboundedness of the perturbation $F$. Concerning this second
issue, let us also mention that Schauder estimates for unbounded non-linear drift coefficients in the non-degenerate diffusive setting were
obtained under mild smoothness assumptions by Krylov and Priola \cite{Krylov:Priola10} who heavily used an auxiliary, deterministic flow associated with the
transport term in
\eqref{Degenerate_Stable_PDE}, i.e.\ for a fixed couple $(t,x)$,
\begin{equation}\label{eq:_INtro_def_flusso}
\begin{cases}
\partial_s\theta_s(x) =A\theta_s(x)+F(s,\theta_s(x)); & \mbox{if } s>t \\
  \theta_t(x) \, = \, x,
\end{cases}
 \quad
\end{equation}
to precisely get rid of the unbounded terms.

%
The drawback of this approach is that we will need at first to establish Schauder estimates in a small time interval.
This seems quite intuitive since the expansion along the chosen proxy on which the method relies is precisely designed for small times
because it requires that the original operator and the proxy are "close" enough in a suitable sense. To obtain the result for an arbitrary
but finite time, we will then iterate the reasoning, which is quite natural since Schauder estimates provide a sort of stability in the
considered functional space. We are therefore far from the optimal constants for the Schauder estimates established in the non-degenerate,
diffusive setting for time dependent coefficients by Krylov and Priola \cite{Krylov:Priola17}.

On the other hand, we want to establish the Schauder estimates in the sharpest possible H\"older setting for the coefficients of the IPDE
\eqref{Degenerate_Stable_PDE}. To do so, we will need to establish some subtle controls, in particular we have no true derivatives of the
coefficients. This is the reason why we will heavily rely on duality results on Besov spaces (see Section $4.1$ below, Chapter $3$ in
\cite{book:Lemarie-Rieusset02} or \cite{book:Triebel83} for a more complete survey of the argument). However, in contrast with the
non-degenerate case (cf. \cite{Chaudru:Menozzi:Priola19}), we will need to ask
for the perturbation $F$ some additional regularity, represented by parameter $\gamma_i$ in assumption [\textbf{R}] below, on the
degenerate entries $F_i$ $(i>1)$. This assumption seems quite natural if we think that, due to the  degenerate structure of the system
(cf. Section $2.2$ below), the more we descend on the chain, the lower the smoothing effect of $\mathcal{L}_\alpha$ will be.
The additional smoothness on $F$ can be then seen as the "price" to pay to re-equilibrate the increasing time
singularities appearing along the chain.

\paragraph{Organization of the paper.} The article is organized as follows. We state our precise framework and give our main results in the
following Section $2$. Section $3$ is then dedicated to the perturbative approach which is the central argument to derive our estimates. In
particular, we obtain therein some Schauder estimates for drifted operators along the inhomogeneous flow $\theta_{s,t}$ defined above
in \eqref{eq:_INtro_def_flusso}, as well as the key Duhamel representation for solutions. Since the arguments to show the Schauder estimates
will be quite long and involved, we postpone the proofs of these results in the next Sections $4$ and $5$. The existence results are then
established in Section $6$. In the last Section $7$, we are going to explain briefly how the perturbative approach presented before could be
applied with slight modifications to prove Schauder-type estimates for a class of completely non-linear, locally H\"older continuous drifts
with an additional "diffusion" coefficient. \newline
Finally, the proof of some technical results concerning the stability properties of H\"older flows are postponed to the Appendix.

\section{Setting and main results}
\fancyhead[RO]{Section \thesection. Setting and main results}
\subsection{Main operators considered}

The operator $\mathcal{L}_\alpha$ we consider is the generator of a non-degenerate, symmetric, stable process and it acts only on the first $d$
coordinates of the system. More precisely, $\mathcal{L}_\alpha$ can be represented for any sufficiently regular $\phi\colon [0,T]\times\R^{nd}\to \R$
as
\[\mathcal{L}_\alpha\phi(t,x)\, := \, \text{p.v.}\int_{\R^d}\bigl[\phi(t,x+By)-\phi(t,x) \bigr] \,\nu_\alpha(dy), \,\, \text{ where }
\,\, B \, := \,
    \begin{bmatrix}
           I_{d\times d} \\
            0_{d\times d}\\
            \vdots\\
           0_{d\times d}
    \end{bmatrix}\]
and $\nu_\alpha$ is a symmetric, stable L\'evy measure on $\R^d$ of order $\alpha$ that we assume to be non-degenerate in a sense that we
are going to specify below. \newline
Passing to polar coordinates $y=\rho s$ where $(\rho,s) \in [0,\infty)\times \mathbb{S}^{d-1}$, it is well-known (see for example Chapter $3$ in \cite{book:Sato99}) that the stable L\'evy measure $\nu_\alpha$ can be decomposed as
\begin{equation}\label{eq:decomposition_measure}
\nu_\alpha(dy) \, :=\, \tilde{\mu}(ds)\frac{d\rho}{\rho^{1+\alpha}},
\end{equation}
where $\tilde{\mu}$ is a symmetric measure on $\mathbb{S}^{d-1}$ which represents the spherical part of $\nu_\alpha$.\newline
We remember now that the L\'evy symbol associated with $\mathcal{L}_\alpha$ is defined through the Levy-Khitchine formula (see, for instance
\cite{book:Jacob05}) as:
\[\Phi(p) \, := \, \int_{\R^d}\bigl[e^{i p\cdot y}-1\bigr]\, \nu_\alpha(dy), \quad \text{ for any $p$ in $\R^d$},\]
where ``$\cdot$'' represents the inner product on the smaller space $\R^d$. In the current symmetric setting, it can be rewritten (cf.
Theorem $14.10$ in \cite{book:Sato99}) as
\begin{equation}\label{def:Levy_Symbol_stable}
\Phi(p) \, = \, -\int_{\mathbb{S}^{d-1}}\vert p\cdot s \vert^\alpha \, \mu(ds),
\end{equation}
where $\mu=C_{\alpha,d}\tilde{\mu}$ is usually called the spherical measure associated with $\nu_\alpha$ . Following \cite{Kolokoltsov00}, we then say that
$\nu_\alpha$ is non-degenerate if the associated L\'evy symbol $\Phi$ is equivalent, up to some multiplicative constant, to $\vert p \vert^\alpha$.
More precisely, we suppose that $\mu$ is non-degenerate if
\begin{description}
  \item[{[ND]}] there exists a constant $\eta\ge 1$ such that for any $p$ in $\R^d$.
\begin{equation}\label{non-degeneracy_of_measure}
\eta^{-1}\vert p \vert^\alpha \, \le \, \int_{\mathbb{S}^{d-1}}\vert p\cdot s \vert^\alpha \, \mu(ds) \, \le\,\eta
\vert p \vert^\alpha.
\end{equation}
\end{description}
It is important to remark that such a condition does not restrict our model too much. Indeed, there are many
different kind of spherical measures $\mu$ that are non-degenerate in the above sense, from the stable-like case, i.e.\ measures that are
absolutely continuous with respect to the Lebesgue measure on $\mathbb{S}^{d-1}$, to very singular ones such that the spherical measure induced by the sum of Dirac masses along the canonical directions:
\[\sum_{i=1}^{d}(\partial^2_{x_k})^{\alpha/2}.\]
We can introduce now the complete Ornstein-Uhlenbeck operator $L^{\text{ou}}$, defined for any sufficiently regular $\phi\colon \R^{nd}\to \R$ as
\begin{equation}\label{eq:def_of_OU_operator}
L^{\text{ou}}\phi(x) \, := \, \langle Ax, D_{x}\phi(x)\rangle +\mathcal{L}_\alpha\phi(x),
\end{equation}
where $A$ is the matrix in $\R^{nd}\times\R^{nd}$ defined in Equation \eqref{eq:def_matrix_A_stable}. We assume
that $A$ satisfies the following H\"ormander-like condition of non-degeneracy:
\begin{description}
  \item[{[H]}] $A_{i,i-1}$ is non-degenerate (i.e.\ it has full rank $d$) for any $i$ in $\llbracket2,n \rrbracket$.
\end{description}
Above, $\llbracket2,n \rrbracket$ denotes the set of all the integers in the interval. It is well known (see for example \cite{book:Sato99}) that under these assumptions, the operator $L^{\text{ou}}$ generates a convolution Markov
semigroup $\{P^{\text{ou}}_t\}_{t\ge 0}$ on $B_b(\R^{nd})$, the family of all the bounded and Borel measurable functions on $\R^{nd}$,
defined by
\[
\begin{cases}
  P^{\text{ou}}_t\phi(x) \, = \, \int_{\R^{nd}}\phi(x+y) \, \mu_t(dy); \\
  P^{\text{ou}}_0\phi(x) \, = \, \phi(x),
\end{cases}
\]
where $\{\mu_t\}_{t>0}$ is a family of Borel probability measures on $\R^{nd}$. In particular, the function
$P^{\text{ou}}_t\phi$ provides the classical solution to the Cauchy problem

\begin{equation}\label{PDE_of_OU}
\begin{cases}
  \partial_tu(t,x)+\mathcal{L}_\alpha u(t,x)+\langle Ax, D_{x}u(t,x)\rangle \, = \, 0 \,\, &\mbox{ on } (0,\infty)\times \R^{nd}; \\
 u(0,x) \, = \, \phi(x) \,\,  &\mbox{ on } \R^{nd}.
\end{cases}
\end{equation}

Moving to the stochastic counterpart if necessary, it is readily derived from \cite{Priola:Zabczyk09} that the semigroup $(P^{\text{ou}}_t)_{t\ge0}$
admits a smooth density $p^{\text{ou}}(t,\cdot)$ with respect to the Lebesgue measure on $\R^{nd}$. Moreover, such a density $p^{\text{ou}}$ has the
following useful representation:
\begin{equation}\label{eq:Representation_of_p_ou}
  p^{\text{ou}}(t,x,y) \, = \, \frac{1}{\det \mathbb{M}_t}p_S(t,\mathbb{M}^{-1}_t\bigl(e^{At}x-y)\bigr),
\end{equation}
where $p_S$ is the density of $\{S_t\}_{t\ge 0}$, a stable process in $\R^{nd}$ whose L\'evy measure satisfies the assumption
[\textbf{ND}] above on $\R^{nd}$ and $\mathbb{M}_t$ is a diagonal matrix on $\R^{nd}\times\R^{nd}$ given by
\begin{equation}\label{eq:def_of_Mt}
\bigl[\mathbb{M}_t\bigr]_{i,j} \, := \,
\begin{cases}
  t^{i-1}I_{d\times d}, & \mbox{if } i=j; \\
  0_{d\times d}, & \mbox{otherwise}.
\end{cases}
\end{equation}
We remark already that the appearance of the matrix $\mathbb{M}_t$ in Equation \eqref{eq:Representation_of_p_ou} and its particular
structure reflect the multi-scaled structure of the dynamics considered (cf. Paragraph ($2.2$) below for a more precise explanation).
\newline
Moreover, the density $p_S$ shows a useful property we will call the \emph{smoothing effect} since it will be fundamental to reduce the
singularities appearing when working with time integrals. Fixed $\gamma$ in $[0,\alpha)$, there exists a constant
$C:=C(\gamma)$ such that for any $l$ in $\llbracket0,3\rrbracket$,
\begin{equation}\label{Smoothing_effect_of_S}
\int_{\R^{nd}}\vert y \vert^\gamma \vert D^l_{y}p_S(t,y) \vert \, dy \, \le \, Ct^{\frac{\gamma-l}{\alpha}} \, \, \text{
for any }t>0.
\end{equation}
These results can be proven following the arguments of Proposition $2.3$ and Lemma $4.3$ in \cite{Huang:Menozzi:Priola19}. We will provide
however a complete proof in the Appendix for the sake of completeness.

\subsection{Intrinsic time scale and associated H\"older spaces}
In this section, we are going to choose which is the most suitable functional space in which to state our Schauder
estimates.\newline
To answer this question, we need firstly to understand how the system typically behaves. We focus for the moment on the Ornstein-Uhlenbeck
case:
\[\bigl(\partial_t  + L^{\text{ou}}\bigr)u(t,x) \, = \, -f(t,x) \quad \text{ on } (0,\infty)\times \R^{nd}, \]
and search for a dilation operator $\delta_\lambda\colon (0,\infty)\times \R^{nd} \to (0,\infty)\times \R^{nd}$ that is invariant for the
considered dynamics, i.e.\ a dilation that transforms solutions of the above equation into other solutions of the same equation.\newline
Due to the structure of $A$ and the $\alpha$-stability of $\nu$, we can consider for any fixed $\lambda>0$, the following
\[ \delta_\lambda(t,x) := (\lambda^\alpha t,\lambda x_1,\lambda^{1+\alpha} x_2,\dots,\lambda^{1+\alpha(n-1)}x_n),\]
i.e.\ with a possible slight abuse of notation, $\bigl(\delta_\lambda(t,x)\bigr)_0:=\lambda^\alpha t$ and for any $i$ in $\llbracket
1,n\rrbracket$, $\bigl(\delta_\lambda(t,x)\bigr)_i := \lambda^{1+\alpha(i-1)}x_i$. It then holds that
\[\bigl(\partial_t +L^{\text{ou}}\bigr) u = 0 \, \Longrightarrow \bigl(\partial_t +L^{\text{ou}} \bigr)(u \circ \delta_\lambda) = 0.\]
The previous reasoning suggests us to introduce a parabolic distance $\mathbf{d}_P$ that is homogenous with respect to the dilation $\delta_\lambda$,
so that $\mathbf{d}_P\bigl(\delta_\lambda(t,x);\delta_\lambda(s,x')\bigr) = \lambda \mathbf{d}_P\bigl((t,x);(s,x')\bigr)$. Precisely,
following the notations in \cite{Huang:Menozzi:Priola19}, we set for any $s,t$ in $[0,T]$ and any $x,x'$ in $\R^{nd}$,
\begin{equation}\label{Definition_distance_d_P_stable}
\mathbf{d}_P\bigl((t,x),(s,x')\bigr)  \, := \, \vert s-t\vert^\frac{1}{\alpha}+\sum_{j=1}^{n} \vert(x-x')_j\vert^{
\frac{1}{1+\alpha(j-1)}}.
\end{equation}
The idea of a dilation $\delta_\lambda$ that summarizes the multi-scaled behaviour of the dynamics was firstly introduced by
Lanconelli and Polidoro in \cite{Lanconelli:Polidoro94} for degenerate Kolmogorov equations in the diffusive setting. Since then, it has
become a ``standard'' tool in the analysis of degenerate equations (see for example \cite{Lunardi97}, \cite{Huang:Menozzi:Priola19}
or \cite{Hao:Wu:Zhang19}). \newline
Since we will quite always use only the spatial part of the distance $\mathbf{d}_P$, we denote for simplicity
\begin{equation}\label{Definition_distance}
\mathbf{d}(x,y) \, = \, \sum_{j=1}^{n} \vert(x-x')_j\vert^{\frac{1}{1+\alpha(j-1)}}.
\end{equation}
Technically speaking, $\mathbf{d}_P$ (and thus, $\mathbf{d}$) does not however induce a norm on $[0,T]\times \R^{nd}$ in the usual sense since it lacks of
linear homogeneity. We remark anyhow again that for any $\lambda>0$, it precisely holds that
$\mathbf{d}\bigl(\delta_\lambda(t,x);\delta_\lambda(s,x')\bigr) = \lambda \mathbf{d}\bigl((t,x);(s,x')\bigr)$.
As it can be seen, $\mathbf{d}_P$ is an extension of the standard parabolic distance in the stable case, adapted to respect the multi-scaled nature
of our dynamics. Indeed, the exponents appearing in \eqref{Definition_distance_d_P_stable} are those which make each space component homogeneous to
the characteristic time scale $t^{1/\alpha}$.\newline
 The appearance of this kind of phenomena is due essentially by the particular structure of the
matrix $A$ (cf. Equation \eqref{Degenerate_Stable_PDE}) that allows the smoothing effect of $\mathcal{L}_\alpha$, acting only on the first variable,
to propagate in the system, as it can be seen in the following lemma:

\begin{lemma}[Scaling Lemma]\label{lemma:Scaling_Lemma}
Let $i$ be in $\llbracket 1,n \rrbracket$. Then, there exist $\{C_j\}_{j \in \llbracket 1,n\rrbracket}$ positive constants, depending only
from $A$ and $i$, such that
\[D_{x_i}p^{\text{ou}}(t,x,y) \, = \, -\sum_{j=i}^{n}C_jt^{j-i}D_{ y_j}p^{\text{ou}}(t,x,y)\]
for any $t>0$ and any $x$, $y$ in $\R^{nd}$.
\end{lemma}
\begin{proof}
Recalling the representation of $p^{\text{ou}}$ in Equation \eqref{eq:Representation_of_p_ou}, it is easy to see that
\[D_{x_i}p^{\text{ou}}(t,x,y) \,=\, \frac{1}{\det \mathbb{M}_t}D_{z}p_S(t,\cdot)\bigl(\mathbb{M}^{-1}_t(e^{At}x- y)
\bigr)\mathbb{M}^{-1}_tD_{x_i}\bigl[e^{At}x-y\bigr].\]
Hence, in order to conclude, we need to show that
\begin{equation}\label{Proof:Scaling_Lemma}
D_{x_i}\bigl[e^{At}x-y\bigr] \, = \, -\sum_{j=i}^{n}C_jt^{j-i}D_{ y_j}\bigl[e^{At}x-y\bigr].
\end{equation}
To prove the above equality, we need to analyze more in depth the structure of the resolvent $e^{At}$. Recalling from Equation
\eqref{eq:def_matrix_A_stable} that $A$ has a sub-diagonal structure, we notice that for any $i,j$ in $\llbracket 1,n\rrbracket$,
\begin{equation}\label{Proof:Scaling_Lemma1}
\Bigl[e^{At}\Bigr]_{i,j} \, = \,
    \begin{cases}
     C_{i,j} t^{j-i}, & \mbox{if } j\ge i;\\
     0, & \mbox{otherwise},
    \end{cases}
\end{equation}
for a family of constants $\{C_{i,j}\}_{i,j \in \llbracket 1,n\rrbracket}$ depending only from $A$. It then follows that for any
$x,y$ in $\R^{nd}$, it holds that
\begin{equation}\label{Proof:Scaling_Lemma2}
\Bigl[e^{At}x-y\Bigr]_{i} \, = \, \sum_{k=1}^{i}C_{i,k}t^{i-k} x_k- y_i.
\end{equation}
Equation \eqref{Proof:Scaling_Lemma} then follows immediately. For a more detailed proof of this result, see also
\cite{Huang:Menozzi16} or \cite{Huang:Menozzi:Priola19}.
\end{proof}

We finally remark the link with the stochastic counterpart of equation \eqref{Degenerate_Stable_PDE}. From a more
probabilistic point of view, the exponents in Equation \eqref{Definition_distance_d_P_stable} can be related to the characteristic
time scales of the iterated integrals of an $\alpha$-stable process.

We are now ready to define the suitable H\"older spaces for our estimates. We start recalling some useful notations we will need below.
Fixed $k$ in $\N\cup \{0\}$ and $\beta$ in $(0,1)$, we follow Krylov \cite{book:Krylov96}, denoting the usual \emph{homogeneous}
H\"older space $C^{k+\beta}(\R^d)$ as the family of functions $\phi\colon \R^d \to \R$ such that
\[\Vert \phi \Vert_{C^{k+\beta}} \,:=\, \sum_{i=1}^{k}\sup_{\vert\vartheta\vert = i}\Vert D^\vartheta\phi
\Vert_{L^\infty}+\sup_{\vert\vartheta\vert=k}[\mathbf{d}^{\vartheta}\phi]_\beta \, < \, \infty,\]
where
\[[\mathbf{d}^{\vartheta}\phi]_\beta \, := \, \sup_{x\neq y}\frac{\vert \mathbf{d}^{\vartheta}\phi(x) -
\mathbf{d}^{\vartheta}\phi(y) \vert}{\vert x - y\vert^{\beta}}.\]
Additionally, we are going to need the associated subspace $C^{k+\beta}_b(\R^d)$ of bounded functions in $C^{k+\beta}(\R^d)$, equipped with
the norm
\[\Vert \cdot \Vert_{C^{k+\beta}_b} \, = \, \Vert \cdot \Vert_{L^\infty} + \Vert \cdot \Vert_{C^{k+\beta}}.\]
We can now define the anisotropic H\"older space with multi-index of regularity associated with the distance $d$.
For sake of brevity and readability, we firstly define for a function $\phi \colon \R^{nd} \to \R$, a point $z$ in $\R^{d(n-1)}$ and $i$ in
$\llbracket 1,n\rrbracket$, the function
\[\Pi_z^i \phi \colon x \in \R^d \, \mapsto\,  \phi(z_1,\dots,z_{i-1},x,z_{i+1},z_n) \in \R,\]
with the obvious modifications if $i=1$ or $i=n$.
Intuitively speaking, the function $\Pi_z^i\phi$ is the restriction of $\phi$ on its $i$-th $d$-dimensional variable while fixing all the other coordinates
in $z$. The space $C^{k+\beta}_d(\R^{nd})$ is then defined as the family of all the function $\phi \colon \R^{nd} \to \R$ such that
\begin{equation}\label{eq:def_anisotropic_norm}
\Vert \phi \Vert_{C^{k+\beta}_d} \, := \, \sum_{i=1}^{n}\sup_{z\in \R^{d(n-1)}}\Vert \Pi_z^i\phi
\Vert_{C^{\frac{k+\beta}{1+\alpha(i-1)}}} \, < \, \infty.
\end{equation}
The modification to the bounded subspace $C^{k+\beta}_{b,d}(\R^{nd})$ is straightforward.\newline
Roughly speaking, the anisotropic norm works component-wise, i.e.\ we firstly fix a coordinate and then calculate the standard H\"older norm along that particular direction, but with
index scaled according to the dilation of the system in that direction, uniformly over time and the other space components. We conclude
summing the contributions associated with each component.
\newline We highlight however that it is possible to recover the expected joint regularity for the partial derivatives, when
they exist. In such a case, they actually turn out to be H\"older continuous in the pseudo-metric $\mathbf{d}$ with order one less than the function.
(cf. Lemma \ref{lemma:Reverse_Taylor_Expansion} in the Appendix for the case $i=1$).

Since we are working with evolution equations, the functions we consider will quite often depend on time, too. For this reason, we denote by
$L^\infty\bigl(0,T;C^{k+\beta}_d(\R^{nd})\bigr)$ (respectively, $L^\infty\bigl(0,T;C^{k+\beta}_{b,d}(\R^{nd})\bigr)$) the family of
functions $\psi \colon [0,T]\times \R^{nd} \to \R$ with finite $C^{k+\beta}_d$-norm (respectively, $C^{k+\beta}_{b,d}$-norm), uniformly in time. It is endowed with the following norm:
\begin{equation}\label{eq:def_anisotropic_norm_in_time}
\Vert \phi \Vert_{L^\infty(C^{k+\beta}_d)} \, := \, \sup_{t\in [0,T]}\Vert \phi(t,\cdot)\Vert_{C^{k+\beta}_d},
\end{equation}
with a straightforward modification for the bounded subspace $L^\infty\bigl(0,T;C^{k+\beta}_{b,d}(\R^{nd})\bigr)$.

\subsection{Assumptions and main results}

From this point further, we consider two fixed numbers $\alpha$ in $(0,2)$ and $\beta$ in $(0,1)$ such that $\alpha$ will represent the
index of stability of the operator $\mathcal{L}_\alpha$ while $\beta$ will stand for the index of H\"older regularity of the coefficients.

We will assume the following:
\begin{description}
  \item[{[S]}] assumptions [\textbf{ND}] and  [\textbf{H}] are satisfied and the drift $F=(F_1,\dots,F_n)$ is such that for
      any $i$ in $\llbracket 1,n\rrbracket$, $F_i$ depends only on time and on the last $n-(i-1)$ components, i.e.\
      $ F_i(t,x_i,\dots,x_n)$;
  \item[{[P]}] $\alpha$ is a number in $(0,2)$, $\beta$ is in $(0,1)$ such that $\alpha+\beta\in (1,2)$ and if $\alpha<1$ (super-critical
      case),
  \[\beta<\alpha, \quad  1-\alpha <\frac{\alpha-\beta}{1+\alpha(n-1)};\]
  \item[{[R]}] Recalling the notations in \eqref{eq:def_anisotropic_norm}-\eqref{eq:def_anisotropic_norm_in_time}, the source $f$ is in $L^\infty(0,T;C^{\beta}_{b,d}(\R^{nd}))$, the terminal
      condition $u_T$ is in $C^{\alpha+\beta}_{b,d}(\R^{nd})$ and for any $i$ in $\llbracket 1,n\rrbracket$, the drift $F_i$ belongs to
      $L^\infty(0,T;C^{\gamma_i+\beta}_{d}(\R^{nd}))$, where
    \begin{equation}\label{Drift_assumptions}
    \gamma_i \,:= \,
    \begin{cases}
        1+ \alpha(i-2), & \mbox{if } i>1; \\
        0, & \mbox{if } i=1.
   \end{cases}
   \end{equation}
\end{description}
From now on, we will say that assumption [\textbf{A}] holds when the above conditions [\textbf{S}], [\textbf{P}] and
[\textbf{R}] are in force.

\begin{remark}[About the assumptions]
We remark that the constraints [\textbf{P}] we are imposing in the super-critical case ($\alpha<1$) seem quite natural for our system. The
condition $\beta<\alpha$ reflects essentially the low integrability properties of the stable density $p_S$ (cf. Equation
\eqref{Smoothing_effect_of_S}). Even if one is interested only on the fractional Laplacian case, i.e. $\mathcal{L}_\alpha=\Delta^{\alpha/2}$, such a
condition cannot be dropped in general, since it does not refer to the integrability property of $p_\alpha$ and its derivatives but instead
to those of its ``projection'' $p_S$ on the bigger space $\R^{nd}$ (cf. Equation
\eqref{eq:Representation_of_p_ou}).\newline
The second condition $\alpha+\beta>1$ is necessary to give a point-wise definition of the gradient of a solution $u$ with respect
to the non-degenerate variable $x_1$. Moreover, there is a famous counterexample of Tanaka and his collaborators
\cite{Tanaka:Tsuchiya:Watanabe74} that shows that even in the scalar case, weak uniqueness (a direct consequence of Schauder
estimates) may fail for the associated SDE if $\alpha+\beta$ is smaller than one. \newline
The last assumption is indeed a technical constraint and it is necessary to work properly with the perturbation $F$ at any level
$i=1,\dots,n$. In particular, it seems the minimal threshold that allows us to exploit the smoothing effect of the density (see for example
Equation \eqref{Proof:Second_Besov_Control_I_2} in the proof of Lemma \ref{lemma:Second_Besov_COntrols} for more details). We conclude
highlighting that these assumptions are always fulfilled if $\alpha\ge 1$ (sub-critical case).
\end{remark}

At this stage, it should be clear that under our assumptions [\textbf{A}], IPDE \eqref{Degenerate_Stable_PDE} will be understood in a
\emph{distributional} sense. Indeed, we cannot hope to find a ``classical'' solution for Equation \eqref{Degenerate_Stable_PDE}, since for such a
function $u$ in $L^\infty(0,T;C^{\alpha+\beta}_{b,d}(\R^{nd}))$, the total gradient $D_{x}u$ is not defined point-wise.\newline
Let us denote for any function $\phi\colon [0,T]\to\R^{nd}$ regular enough, the complete operator $L_t$ as
\begin{equation}\label{eq:Complete_Operator}
L_t\phi(t,x) \, := \, \langle A x + F(t,x), D_{x}u(t,x)\rangle + \mathcal{L}_\alpha u(t,x).
\end{equation}
We will say that a function $u$ in $L^\infty(0,T;C^{\alpha+\beta}_{b,d}(\R^{nd}))$ is a distributional (or weak) solution of IPDE
\eqref{Degenerate_Stable_PDE} if for any $\phi$ in $C^\infty_0((0,T]\times\R^{nd})$, it holds that
\begin{equation}\label{eq:Def_weak_sol}
\int_{0}^{T}\int_{\R^{nd}}\Bigl(-\partial_t+L_t^*\Bigr)\phi(t,y)u(t,y) \, dy +
\int_{\R^{nd}}u_T(y)\phi(T,y) \, dy \, = \, - \int_{0}^{T}\int_{\R^{nd}}\phi(t,y)f(t,y) \, dy,
\end{equation}
where $\mathcal{L}^*_\alpha$ denotes the formal adjoint of $L_t$.
On the other hand, denoting from now on,
\begin{equation}\label{eq:norm_H_for_F}
\Vert F\Vert_H \, := \sup_{i\in\llbracket 1,n \rrbracket}\Vert  F_i\Vert_{L^\infty(C^{\gamma_i+\beta}_{d})},
\end{equation}
we will quite often use the following other notion of solution:
\begin{definition}
\label{definition:mild_sol}
A function $u$ is a mild solution in $L^\infty\bigl(0,T;C^{\alpha+\beta}_{b,d}(\R^{nd})\bigr)$ of Equation \eqref{Degenerate_Stable_PDE} if
for any triple of sequences $\{f_m\}_{m\in \N}$, $\{u_{T,m}\}_{m\in \N}$ and $\{F_m\}_{m\in \N}$ such that
\begin{itemize}
  \item $\{f_m\}_{m\in \N}$ is in $C^\infty_b((0,T)\times\R^{nd})$ and $f_m$ converges to $f$ in $L^\infty\bigl(0,T;C^\beta_{b,d}(\R^{nd})\bigr)$;
  \item $\{u_{T,m}\}_{m\in \N}$ is in $C^\infty_b(\R^{nd})$ and $u_{T,m}$ converges to $u_T$ in $C^{\alpha+\beta}_{b,d}(\R^{nd})$;
   \item $\{F_m\}_{m\in \N}$ is in $C^\infty_b((0,T)\times\R^{nd};\R^{nd})$ and $\Vert F_m-F\Vert_H$ converges to $0$,
\end{itemize}
there exists a sub-sequence $\{u_{m}\}_{m\in \N}$ in $C^\infty_b\bigl((0,T)\times\R^{nd}\bigr)$ such that
\begin{itemize}
  \item $u_{m}$ converges to $u$ in $L^\infty\bigl(0,T;C^{\alpha+\beta}_{b,d}(\R^{nd})\bigr)$;
  \item for any fixed $m$ in $\N$, $u_{m}$ is a classical solution of the following ``regularized'' IPDE:
  \begin{equation}\label{Regularizied_PDE}
  \begin{cases}
    \begin{aligned}
       &  \partial_t u_m(t,x)+\mathcal{L}_\alpha u_m(t,x) \\
       &\qquad\qquad\qquad+ \langle A x + F_m(t,x), D_{x}
   u_m(t,x)\rangle\, = \, -f_m(t,x),
    \end{aligned}           &\mbox{on }   (0,T)\times \R^{nd}; \\
 u_m(T,x) \, = \, u_{T,m}(x) & \mbox{on }\R^{nd}.
  \end{cases}
  \end{equation}
\end{itemize}
\end{definition}

We can now state our main result:
\begin{theorem}(Schauder Estimates)
\label{theorem:Schauder_Estimates}
Let $u$ be a mild solution in $L^{\infty}\bigl(0,T;C^{\alpha+\beta}_{b,d}(\R^{nd})\bigr)$ of IPDE \eqref{Degenerate_Stable_PDE}. Under
[\textbf{A}], there exists a constant $C:=C(T)$ such that
\begin{equation}\label{equation:Schauder_Estimates}
\Vert u \Vert_{L^\infty(C^{\alpha+\beta}_{b,d})} \, \le \, C \bigl[\Vert f \Vert_{L^\infty(C^{\beta}_{b,d})} + \Vert u_T
\Vert_{C^{\alpha+\beta}_{b,d}}\bigr].
\end{equation}
\end{theorem}

Associated with an existence result we will exhibit in Section $6$, we will eventually derive the well-posedness for Equation
\eqref{Degenerate_Stable_PDE}.

\begin{theorem}
Under [\textbf{A}], there exists a unique mild solution $u$ of IPDE
\eqref{Degenerate_Stable_PDE} belonging to $L^{\infty}\bigl(0,T;C^{\alpha+\beta}_{b,d}(\R^{nd})\bigr)$. Moreover, such a function $u$ is a weak solution, too.
\end{theorem}

In the following, we will denote for sake of brevity
\begin{equation}\label{eq:def_alpha_i_and_beta_i}
\alpha_i \, := \, \frac{\alpha}{1+\alpha(i-1)} \, \, \text{ and } \, \, \beta_i \, := \, \frac{\beta}{1+\alpha(i-1)}\,\,  \text{ for any
$i$ in }\llbracket 1,n\rrbracket.
\end{equation}
Clearly, these quantities were introduced to reflect exactly the relative scale of the system at every considered level $i$ (cf. Section
$2.2$ above).\newline
In the following, as well as in Theorem \ref{theorem:Schauder_Estimates} above, $C$ denotes a generic constant that may
change from line to line but depending only on the parameters in assumption [\textbf{A}]. Other dependencies that may occur are explicitly
specified.

\setcounter{equation}{0}
\mysection{Proof through  perturbative approach}
As already said in the introductive section, our method of proof relies on a perturbative approach introduced in
\cite{Chaudru:Honore:Menozzi18_Sharp} for the degenerate, Kolmogorov, diffusive setting.\newline
Roughly speaking, we will firstly choose a suitable proxy for the equation of interest, i.e.\ an operator whose associated semigroup and
density are known and that is close enough to the original one:
\[\mathcal{L}_\alpha+\langle A x +F(t,x), D_{x}\rangle.\]
Furthermore, we will exhibit suitable regularization properties for the proxy and in particular, we will show that it satisfies the
Schauder estimates \eqref{equation:Schauder_Estimates}. This will be the purpose of Sub-section $3.1$.\newline
In Sub-section $3.2$ below, we will then expand a solution $u$ of IPDE \eqref{Degenerate_Stable_PDE} along the chosen proxy through a
Duhamel-type formula and eventually show that the expansion error only brings a negligible contribution, so that the Schauder estimates still
hold for $u$. Due to our choice of method, this will be possible only adding some more assumptions on the system. Namely, we will assume in
addition to be in a small time interval, so that the proxy and the original operator do not differ too much. \newline
The last Sub-section $3.3$ will finally show how to remove the additional assumption in order to prove the Schauder
estimates (Theorem \ref{theorem:Schauder_Estimates}) through a scaling argument.

\subsection{Frozen semigroup}
The crucial element in our approach consists in choosing wisely a suitable proxy operator along which to expand a solution $u$ in
$L^\infty(0,T;C^{\alpha+\beta}_{b,d}(\R^{nd}))$ of IPDE \eqref{Degenerate_Stable_PDE}. In order to deal with potentially unbounded
perturbations $F$, it is natural to use a proxy involving a non-zero first order term  associated with a flow
representing the dynamics driven by $Ax+F$, the transport part of Equation \eqref{Degenerate_Stable_PDE} (see e.g.
\cite{Krylov:Priola10} or \cite{Chaudru:Menozzi:Priola19}).

Remembering that we assume $F$ to be H\"older continuous, we know that there exists a solution of
\[
\begin{cases}
  d \theta_{s,\tau}(\xi) \, = \, \bigl[A\theta_{s,\tau}(\xi)+F(s,\theta_{s,\tau}(\xi))\bigr]ds \,\,\, \mbox{ on
  } [\tau,T];\\
  \theta_{\tau,\tau}(\xi) \, = \, \xi,
\end{cases}
\]
even if it may be not unique. For this reason, we are going to choose one particular flow, denoted by
$\theta_{s,\tau}(\xi)$, and consider it fixed throughout the work.\newline
More precisely, given a freezing couple $(\tau,\xi)$ in $[0,T]\times\R^{nd}$, the flow will be defined on $[\tau,T]$ as
\begin{equation}\label{Flow}
  \theta_{s,\tau}(\xi) \, = \, \xi + \int_{\tau}^{s}\bigl[A\theta_{v,\tau}(\xi)+F(v,\theta_{v,\tau}(\xi))\bigr] \, dv.
\end{equation}
We can now introduce the ``frozen'' IPDE on $(0,T)\times \R^{nd}$, associated with the chosen proxy:
\begin{equation}\label{Frozen_PDE}
\begin{cases}
   \partial_s \tilde{u}^{\tau,\xi}(s,x)+\mathcal{L}_\alpha \tilde{u}^{\tau,\xi}(s,x) + \langle A x +
   F(s,\theta_{s,\tau}(\xi)), D_{x} \tilde{u}^{\tau,\xi}(s,x)\rangle  \, = \, -f(s,x); \\
    \tilde{u}^{\tau,\xi}(s,x) \, = \, u_T(x).
  \end{cases}
\end{equation}
Noticing that the proxy operator $\mathcal{L}_\alpha+\langle A x + F(s,\theta_{s,\tau}(\xi)), D_{x}\rangle$ can be seen as
an Ornstein-Ulhenbeck operator with an additional time-dependent component $F(s,\theta_{s,\tau}(\xi))$,  it is clear that
under assumption [\textbf{A}], it generates a two parameters semigroups we will denote by $\{\tilde{P}^{\tau,\xi}_{s,t}\}_{t\le
s}$. Moreover, it admits a density given by
\begin{equation}\label{eq:definition_tilde_p}
\tilde{p}^{\tau,\xi}(t,s,x,y)\,  = \, \frac{1}{\det \mathbb{M}_{s-t}} p_S\bigl( s-t,\mathbb{M}^{-1}_{s-t} (y -
\tilde{m}^{ \tau,\xi}_{s,t}(x))\bigr),
\end{equation}
remembering Equation \eqref{eq:Representation_of_p_ou} for the definition of $p_S$ and with the following notation for the ``frozen shift''
$\tilde{m}^{\tau,\xi}_{s,t}(x)$:
\begin{equation}\label{eq:def_tilde_m_stable}
\tilde{m}^{\tau,\xi}_{s,t}(x) \, = \, e^{A(s-t)}x + \int_{t}^{s}e^{A(s-v)} F(v, \theta_{v,\tau}
(\xi)) \, dv.
\end{equation}

We point out already the following important property of the shift $\tilde{m}^{\tau,\xi}_{s,t}(x)$:

\begin{lemma}
\label{lemma:identification_theta_m_stable}
Let $t<s$ in $[0,T]$ and $x$ a point in $\R^{nd}$. Then,
\begin{equation}\label{eq:identification_theta_m_stable}
\tilde{m}^{\tau,\xi}_{s,t}(x) \, = \, \theta_{s,\tau}(\xi),
\end{equation}
taking $\tau=t$ and $\xi=x$.
\end{lemma}
\begin{proof}
We start noticing that by construction, $\tilde{m}^{\tau,\xi}_{s,t}(x)$ satisfies
\[\tilde{m}^{\tau,\xi}_{s,t}(x) \, = \, x +\int_{t}^{s} \bigl[A\tilde{m}^{t,x}_{v,t}(x)+F(v, \theta_{v,\tau} (\xi))\bigr] \, dv.\]
It then holds that
\[\vert \tilde{m}^{t,x}_{s,t}(x) - \theta_{s,t}(x) \vert \, \le \, \int_{t}^{s}A\vert
\tilde{m}^{t,x}_{v,t}(x)-\theta_{v,t}(x) \vert \, dv.\]
The above Equation \eqref{eq:identification_theta_m_stable} then follows immediately applying the Gr\"onwall lemma.
\end{proof}

Moreover, we can extend the smoothing effect \eqref{Smoothing_effect_of_S} of $p_S$ to the frozen density $\tilde{p}^{\tau,\xi}$
through the representation \eqref{eq:definition_tilde_p}:

\begin{lemma}[Smoothing effects of the frozen density]
\label{lemma:Smoothing_effect_frozen}
Under [\textbf{A}], let $\vartheta,\varrho$ be two multi-indexes in $\N^n$ such that $\vert \varrho+\vartheta\vert \le 3$ and $\gamma$ in
$[0,\alpha)$. Then, there exists a constant $C:=C(\vartheta,\varrho,\gamma)$ such that
\begin{equation}\label{eq:Smoothing_effects_of_tilde_p}
\int_{\R^{nd}} \vert D^\varrho_{y}D^\vartheta_{x}
\tilde{p}^{ \tau,\xi} (t,s,x,y) \vert \mathbf{d}^\gamma\bigl(y,\tilde{m}^{\tau,\xi}_{s,t}(x)\bigr) \,
dy \, \le \, C (s-t)^{\frac{\gamma}{\alpha}-\sum_{k=1}^{n} \frac{\vartheta_k+\varrho_k}{\alpha_k}}
\end{equation}
for any $t<s$ in $[0,T]$, any $x$ in $\R^{nd}$ and any frozen couple $(\tau,\xi)$ in $[0,T]\times\R^{nd}$.
In particular, if $\vert \vartheta \vert \neq 0$, it holds for any $\phi$ in $C^\gamma_d(\R^{nd})$ that
\begin{equation}\label{eq:Control_of_semigroup}
\bigl{\vert} D^\vartheta_{x}\tilde{P}^{\tau,\xi}_{s,t}\phi(x) \bigr{\vert} \, \le \, C \Vert\phi\Vert_{C^\gamma_d}
(s-t)^{\frac{\gamma}{\alpha}-\sum_{k=1}^{n}\frac{\vartheta_k}{\alpha_k}}.
\end{equation}
\end{lemma}
\begin{proof}
Since $p_S$ is the density of an $\alpha$-stable process, we remember that the following $\alpha$-scaling property
\begin{equation}\label{eq:alpha_scaling_p_S}
p_S(t,y) \, = \, t^{-\frac{nd}{\alpha}}p_S(1,t^{-\frac{1}{\alpha}}y)
\end{equation}
holds for any $t>0$ and any $y$ in $\R^{nd}$. Fixed $i$ in $\llbracket1,n\rrbracket$, we then denote for simplicity
\[\mathbb{T}_{s-t}\, := \, (s-t)^{\frac{1}{\alpha}}\mathbb{M}_{s-t}\]
and we calculate the derivative of $\tilde{p}^{\tau,\xi}$ with respect to $x_i$ through
\[\begin{split}
\vert D_{x_i}\tilde{p}^{\tau,\xi}(t,s,x,y) \vert \, &= \, \Bigl{\vert} \frac{1}{\det \mathbb{M}_{s-t}} D_{x_i} \bigl[p_S\bigl(s-t,\mathbb{M}^{-1}_{s-t} (\tilde{m}^{\tau,\xi}_{s,t}(x)-y)\bigr)] \Bigr{\vert} \\
&= \, \Bigl{\vert} \frac{1}{\det\mathbb{T}_{s-t}}D_{x_i}\bigl[p_S\bigl(1,\mathbb{T}^{-1}_{s-t}(\tilde{m}^{\tau,\xi}_{s,t}
(x)-y)\bigr)] \Bigr{\vert} \\
&= \,\Bigl{\vert} \frac{1}{\det (\mathbb{T}_{s-t})} \langle D_{z}p_S\bigl(1,\cdot\bigr)(\mathbb{T}^{-1}_{s-t}
(\tilde{m}^{\tau,\xi}_{s,t}(x)-y));\mathbb{T}^{-1}_{s-t}D_{x_i}(\tilde{m}^{\tau,\xi}_{s,t}
(x))\rangle \Bigr{\vert},
\end{split}
\]
where in the second equality we exploited the $\alpha$-scaling property in \eqref{eq:alpha_scaling_p_S}. From Equation
\eqref{Proof:Scaling_Lemma1} in the Scaling Lemma \ref{lemma:Scaling_Lemma}, we now notice that
\[
\begin{split}
\bigl{\vert} \mathbb{T}^{-1}_{s-t}D_{x_i}(\tilde{m}^{\tau,\xi}_{s,t} (x)) \bigr{\vert} \, &= \, \bigl{\vert}
\mathbb{T}^{-1}_{s-t}D_{x_i}\bigl(e^{A(t-s)}(x)\bigr) \bigr{\vert} \\
&= \,
(s-t)^{-\frac{1}{\alpha}}\sum_{k=i}^{n}C_k(s-t)^{-(k-1)}(s-t)^{k-i} \\
&\le \,  C(s-t)^{-\frac{1+\alpha(i-1)}{\alpha}},  
\end{split}
\]
and we use it to show that
\[\vert D_{x_i}\tilde{p}^{\tau,\xi}(t,s,x,y) \vert \, \le \, C(s-t)^{-\frac{1+\alpha(i-1)}{\alpha}}\frac{1}{\det
(\mathbb{T}_{s-t})} \bigl{\vert}D_{z}p_S\bigl(1,\cdot\bigr)(\mathbb{T}^{-1}_{s-t}
(\tilde{m}^{\tau,\xi}_{s,t}(x)-y))\bigr{\vert}.\]
Similarly, if we fix $j$ in $\llbracket1,n\rrbracket$, it holds that
\[\vert D_{y_j}D_{x_i}\tilde{p}^{\tau,\xi}(t,s,x,y) \vert \, \le \, C(s-t)^{-\frac{1}{\alpha_i}-\frac{1}{\alpha_j}}
\frac{1}{\det (\mathbb{T}_{s-t})} \bigl{\vert}D^2_{z}p_S\bigl(1,\cdot\bigr)(\mathbb{T}^{-1}_{s-t}(\tilde{m}^{\tau,\xi}_{s,t}
(x)-y))\bigr{\vert}.\]
It is then easy to show by iteration of the same argument that
\begin{multline}\label{eq:derivative_frozen_density}
\vert D^\varrho_{y}D^\vartheta_{x}\tilde{p}^{\tau,\xi}(t,s,x,y) \vert \\
\le \,
C(s-t)^{-\sum_{k=1}^{n}\frac{\varrho_k+\vartheta_k}{\alpha_k}}\frac{1}{\det (\mathbb{T}_{s-t})}
\bigl{\vert}D^{\vert \varrho+\vartheta\vert}_zp_S\bigl(1,\cdot\bigr)(\mathbb{T}^{-1}_{s-t} (\tilde{m}^{\tau,\xi}_{s,t}
(x)-y))\bigr{\vert}.
\end{multline}
Control \eqref{eq:Smoothing_effects_of_tilde_p} immediately follows from the analogous smoothing effect for $p_S$ (cf.\ Equation
\eqref{Smoothing_effect_of_S}) and the change of variables $z=\mathbb{T}^{-1}_{s-t} (\tilde{m}^{\tau,
\xi}_{s,t} (x)-y)$. Indeed,
\[
\begin{split}
\int_{\R^{nd}} \vert D^\varrho_{y}D^\vartheta_{x}&\tilde{p}^{ \tau,\xi} (t,s,x,y) \vert
\mathbf{d}^\gamma\bigl(y,\tilde{m}^{\tau,\xi}_{s,t}(x)\bigr) \,
dy \, \le \, C(s-t)^{-\sum_{k=1}^{n}\frac{\varrho_k+\vartheta_k}{\alpha_k}}\\
&\qquad  \quad \times 
\int_{\R^{nd}}\frac{1}{\det(\mathbb{T}_{s-t})}\bigl{\vert}D^{\vert
\varrho+\vartheta\vert}_{z}p_S\bigl(1,\cdot\bigr)(\mathbb{T}^{-1}_{s-t}(\tilde{m}^{\tau,\xi}_{s,t} (x)-y))
\bigr{\vert} \mathbf{d}^\gamma\bigl(y,\tilde{m}^{\tau,\xi}_{s,t}(x)\bigr) \, dy \\
&= \,(s-t)^{-\sum_{k=1}^{n}\frac{\varrho_k+\vartheta_k}{\alpha_k}}\int_{\R^{nd}}\bigl{\vert}D^{\vert \varrho+\vartheta\vert}_{z}p_S(1,
z)\bigr{\vert}\mathbf{d}^\gamma\bigl(\mathbb{T}_{s-t}z+\tilde{m}^{\tau,\xi}_{s,t}(x),\tilde{m}^{\tau,\xi}_{s,t}(x)\bigr) \,dz.
\end{split}\]
To conclude, we notice that
\[\mathbf{d}^\gamma\bigl(\mathbb{T}_{s-t}z+\tilde{m}^{\tau,\xi}_{s,t}(x),\tilde{m}^{\tau,\xi}_{s,t}(x)\bigr) \,
\le \, C\sum_{i=1}^{n}\vert (s-t)^{\frac{1+\alpha(i-1)}{\alpha}} z_i\vert^{\frac{\gamma}{1+\alpha(i-1)}} \, = \,
(s-t)^{\frac{\gamma}{\alpha}}\sum_{i=1}^{n}\vert  z_i\vert^{\frac{\gamma}{1+\alpha(i-1)}}\]
and use it to write that
\begin{multline*}
\int_{\R^{nd}} \bigl{\vert} D^\varrho_{y}D^\vartheta_{x}\tilde{p}^{\tau,\xi} (t,s,x,y) \bigr{\vert}
\mathbf{d}^\gamma\bigl(y,\tilde{m}^{\tau,\xi}_{s,t}(x)\bigr) \, dy \\
\le \,C(s-t)^{\frac{\gamma}{\alpha}-\sum_{k=1}^{n}\frac{\varrho_k+\vartheta_k}{\alpha_k}}\sum_{i=1}^{n}\int_{\R^{nd}}\bigl{\vert}D^{\vert
\varrho+\vartheta\vert}_{z}p_S(1,z)\bigr{\vert}\vert  z_i\vert^{\frac{\gamma}{1+\alpha(i-1)}} \, dz\, \le \, C(s-t)^{\frac{
\gamma}{\alpha}- \sum_{k=1}^{n}\frac{\varrho_k+\vartheta_k}{\alpha_k}},
\end{multline*}
where in the last passage we used the smoothing effect for $p_S$ (Equation \eqref{Smoothing_effect_of_S}), recalling that for any $i$ in
$\llbracket 1,n\rrbracket$, it holds that
\[\frac{\gamma}{1+\alpha(i-1)} \, \le \, \gamma \, <\, \alpha\]
and we have thus the required integrability.\newline
To prove instead Inequality \eqref{eq:Control_of_semigroup}, we use a cancellation argument to write
\[
\begin{split}
\bigl{\vert} D^\vartheta_{x}\tilde{P}^{\tau,\xi}_{s,t}\phi(x) \bigr{\vert} \, &= \, \Bigl{\vert}\int_{\R^{nd}}
D^\vartheta_{x}\tilde{p}^{\tau,\xi}(t,s,x,y) \bigl[\phi(y)- \phi(\tilde{m}^{\tau,\xi}_{s,t}
(x))\bigr] \, dy\Bigr{\vert} \\
&\le \, \int_{\R^{nd}} \vert D^\vartheta_{x}\tilde{p}^{\tau,\xi}(t,s,x,y)\vert\, \vert
\phi(y)- \phi(\tilde{m}^{\tau,\xi}_{s,t} (x))\vert \, dy.
\end{split}\]
But since we assume $\phi$ to be in $C^\gamma_d(\R^{nd})$, we can control the last expression as
\[
\begin{split}
\bigl{\vert} D^\vartheta_{x}\tilde{P}^{\tau,\xi}_{s,t}\phi(x) \bigr{\vert} \, &\le \, \Vert
\phi\Vert_{C^\gamma_d}\int_{\R^{nd}} \mathbf{d}^\gamma\bigl(y, \tilde{m}^{\tau,\xi}_{s,t} (x)\bigr) \vert
D^\vartheta_{x}\tilde{p}^{\tau,\xi}(t,s,x,y)\vert  \, dy \\ 
&\le \, C \Vert\phi\Vert_{C^\gamma_d}
(s-t)^{\frac{\gamma}{\alpha}-\sum_{k=1}^{n}\frac{\vartheta_k}{\alpha_k}},
\end{split}\]
where in the last passage we used Equation \eqref{eq:Smoothing_effects_of_tilde_p}.
\end{proof}

We can define now our candidate to be the mild solution of the ``frozen'' IPDE. If it exists and it is
smooth enough, such a candidate appears to be a representation of the solution of Equation \eqref{Frozen_PDE} obtained through the Duhamel
principle. For this reason, the following expression
\begin{equation}\label{Duhamel_representation_of_proxy}
\tilde{u}^{\tau,\xi}(t,x) \, := \, \tilde{P}^{\tau,\xi}_{T,t}u_T(x) + \int_{t}^{T}\tilde{P}^{\tau,\xi}_{s,t}f(s,x) \,
ds \quad \text{ for any $(t,x)$ in }[0,T]\times \R^{nd},
\end{equation}
will be called the \emph{Duhamel representation of the proxy}. As it seems, under our assumption [\textbf{A}] such a representation is robust
enough to satisfy Schauder estimates similar to \eqref{equation:Schauder_Estimates}. Since the proof of this result is quite long, we will
postpone it to Section $4.2$ for clarity.

\begin{prop}(Schauder Estimates for Proxy)
\label{prop:Schauder_Estimates_for_proxy}
Under [\textbf{A}], there exists a constant $C:=C(T)$ such that
\begin{equation}\label{equation:Schauder_Estimates_for_proxy}
\Vert \tilde{u}^{\tau,\xi} \Vert_{L^\infty(C^{\alpha+\beta}_{b,d})} \, \le \, C\bigl[\Vert f \Vert_{L^\infty(C^{\beta}_{b,d})}+\Vert u_T \Vert_{C^{\alpha+\beta}_{b,d}} \bigr]
\end{equation}
for any freezing couple $(\tau,\xi)$ in $[0,T]\times\R^{nd}$.
\end{prop}

We conclude this section showing that the function $\tilde{u}^{\tau,\xi}$ is indeed a mild solution in
$L^\infty(0,T;C^{\alpha+\beta}_{b,d}(\R^{nd}))$ of the ``frozen'' IPDE \eqref{Frozen_PDE}.
Moreover, the converse statement is also true. If regular enough, any solution of Equation \eqref{Frozen_PDE} corresponds to the Duhamel
Representation \eqref{Duhamel_representation_of_proxy}.

\begin{prop}
\label{prop:frozen_Duhamel_Formula}
Let us assume to be under assumption [\textbf{A}]. Then,
\begin{itemize}
  \item the function $\tilde{u}^{\tau,\xi}$ defined in \eqref{Duhamel_representation_of_proxy} is a mild solution
      in $L^\infty\bigl(0,T;C^{\alpha+\beta}_{b,d}(\R^{nd})\bigr)$ of the ``frozen'' IPDE \eqref{Frozen_PDE} for any freezing couple $(\tau,\xi)$
      in $[0,T]\times \R^{nd}$;
  \item Fixed $(\tau,\xi)$ in $[0,T]\times \R^{nd}$, let $\tilde{v}^{\tau,\xi}$ be a mild solution in
      $L^\infty\bigl(0,T;C^{\alpha+\beta}_{b,d}(\R^{nd})\bigr)$ of IPDE \eqref{Frozen_PDE}. Then,
    \[\tilde{v}^{\tau,\xi}(t,x) \, = \, \tilde{P}^{\tau,\xi}_{T,t}u_T(x) +  \int_{t}^{T}\tilde{P}^{\tau,\xi}_{s,t}f
    (s,x) \, ds.\]
\end{itemize}
\end{prop}
\begin{proof}
The first assertion is quite straightforward. Let us consider three sequences $\{f_m\}_{m\in \N}$, $\{u_{T,m}\}_{m\in \N}$ and $\{F_m\}_{m\in \N}$  of
smooth and bounded coefficients such that $f_m$ converges to $f$ in $L^\infty\bigl(0,T;C^\beta_{b,d}(\R^{nd})\bigr)$, $u_{T,m}$ to $u_T$ in
$C^{\alpha+\beta}_{b,d}(\R^{nd})$ and $\Vert F_m-F\Vert_H \to 0$. Denoting now by $\{\tilde{P}^{m,\tau, \xi}_{s,t}
\}_{t\le s}$ the semigroup associated with the ``regularized'' operator
\[\mathcal{L}_\alpha +\langle Ax+F_m(t,\theta_{t,\tau}(\xi)),D_{x}\rangle,\]
it is not difficult to show that for any fixed $m$ in $\N$, the following
\[\tilde{u}^{\tau,\xi}_m \, := \, \tilde{P}^{m,\tau,\xi}_{T,t}u_{T,m}(x)+\int_{t}^{T}\tilde{P}^{m,\tau,\xi}_{s,t}f_m
(s,x )\, ds\]
is a classical solution of the ``frozen'' IPDE \eqref{Frozen_PDE} with regularized coefficients $f_m,u_{T,m}$ and $F_m$.
A detailed guide of this result can be found, even if in the diffusive setting, in Lemma $3.3$ in \cite{Krylov:Priola10}. Using now the
Schauder Estimates \eqref{equation:Schauder_Estimates_for_proxy} for the regularized solutions $\tilde{u}^{\tau,\xi}_m$, it follows
immediately that $\tilde{u}^{\tau,\xi}_m \to \tilde{u}^{\tau,\xi}$ in $L^\infty\bigl(0,T;C^{\alpha+\beta}_{b,d}(\R^{nd})\bigr)$
and thus, that
$\tilde{u}^{\tau,\xi}$ is a mild solution of \eqref{Frozen_PDE} in $L^\infty\bigl(0,T;C^{\alpha+\beta}_{b,d}(\R^{nd})\bigr)$.\newline
To prove the second statement, we start fixing a freezing couple $(\tau,\xi)$ in $[0,T]\times \R^{nd}$ and consider three sequences
$\{f_m\}_{m\in \N}$, $\{u_{T,m}\}_{m\in \N}$ and $\{F_m\}_{m\in \N}$ of bounded and smooth coefficients such that $f_m \to f$ in
$L^\infty\bigl(0,T;C^\beta_{b,d}(\R^{nd})\bigr)$, $u_{T,m} \to u_T$ in $C^{\alpha+\beta}_{b,d}(\R^{nd})$ and $\Vert F_m-F\Vert_H
\to 0$. They can be constructed through mollification.\newline
Since $\tilde{v}^{\tau,\xi}$ is a mild solution of the ``frozen'' IPDE \eqref{Frozen_PDE}, we
know that there exists a sequence $\{\tilde{v}^{\tau,\xi}_m\}_{m \in \N}$ of classical solutions of the ``regularized frozen'' IPDE
\eqref{Frozen_PDE} with coefficients $f_m,u_{T,m}$ and $F_m$ such that $\tilde{v}^{\tau,\xi}_m \to \tilde{v}^{\tau,\xi}$ in $L^\infty\bigl(0,T;C^{\alpha+\beta}_{b,d}(\R^{nd})\bigr)$. Fixed
$m$ in $\N$, we then denote
\[h_m(t,x) \, : = \, \tilde{v}^{\tau,\xi}_m\bigl(t,x-\int_{t}^{T}e^{A(t-s)}F_m(s, \theta_{s,\tau}(\xi))\,
ds\bigr)\]
for any $t$ in $[0,T]$ and any $x$ in $\R^{nd}$. Direct calculations imply that
\begin{gather*}
D_{x} h_m(t,x) \, = \, D_{x} \tilde{v}^{\tau,\xi}_m\bigl(t,x-\int_{t}^{T}e^{A(t-s)}F_m(s,
\theta_{s,\tau}(\xi)) \,
ds\bigr); \\
\mathcal{L}_\alpha h_m(t,x) \, = \, \mathcal{L}_\alpha \tilde{v}^{\tau,\xi}_m\bigl(t,x-\int_{t}^{T}e^{A(t-s)}F_m(s,
\theta_{s,\tau}(\xi)) \, ds\bigr)
\end{gather*}
and
\[
\begin{split}
  \partial_t h_m(t,x) \, &= \, \partial_t \tilde{v}^{\tau,\xi}_m\bigl(t,x-\int_{t}^{T}e^{A(t-s)}F_m(s,
  \theta_{s,\tau}(\xi)) \, ds\bigr)\\
  &+ \bigl{\langle} F_m(t,\theta_{t,\tau}(\xi)),
  D_{x}\tilde{v}^{\tau,\xi}_m\bigl(t,x-\int_{t}^{T}e^{A(t-s)}F_m(s,\theta_{s,\tau}(\xi))\,ds\bigr)\bigr{\rangle}\\
  &-\bigl{\langle} A\int_{t}^{T}e^{A(t-s)}F_m(s, \theta_{s,\tau}(\xi)) \, ds,
  D_{x}\tilde{v}^{\tau,\xi}_m\bigl(t,x-\int_{t}^{T}e^{A(t-s)}F_m(s,\theta_{s,\tau}(\xi))\,ds\bigr)\bigr{\rangle}.
\end{split}
\]
Remembering that $\tilde{v}^{\tau,\xi}_m$ is a classical solution of Equation \eqref{Frozen_PDE} replacing therein $f$, $u_T$ and
$F$ with coefficients $f_m,u_{T,m}$ and $F_m$, it follows immediately that the function $h_m$ solves for any $m$ in $\N$ the following:
\begin{equation}\label{Frozen_PDE_No_Perturb}
\begin{cases}
   \partial_t h_m(t,x)+\mathcal{L}_\alpha h_m(t,x) + \langle A x, D_{x} h_m(t,x)\rangle  \, = \, -l_m(t,x);\\
    h_m(T,x) \, = \, u_{T,m}(x)
  \end{cases}
\end{equation}
where $l_m(t,x):=f_m\bigl(t,x-\int_{t}^{T}e^{A(t-s)}F_m(s, \theta_{s,\tau}(\xi)) \, ds\bigr)$.\newline
Since we are going to exploit reasonings in Fourier spaces, we need however to have
integrability properties on the solution $h_m$. For this reason, we introduce now a family $\{\rho_R\}_{R>0}$ of smooth functions such that
any $\rho_R$ is equal to $1$ in
$B(0,R)$ and vanishes outside $B(0,R+1)$. We then denote for any $R>0$,
\[h_{m,R}(t,x) \, := \,h_m(t,x)\rho_R(x).\]
It is then straightforward that $h_{m,R}$ solves
\begin{equation}\label{Frozen_PDE_No_Perturb2}
\begin{cases}
   \partial_t h_{m,R}(t,x)+\mathcal{L}_\alpha h_{m,R}(t,x) + \langle A x, D_{x} h_{m,R}(t,x)\rangle  \, = \, -\tilde{l}_{m,R}(t,x);\\
    h_{m,R}(T,x) \, = \, g_{m,R}(x),
  \end{cases}
\end{equation}
where $g_{m,R}(x)=u_{T,m}(x)\rho_R(x)$ and
\[
\begin{split}
\tilde{l}_{m,R}(t,x) \, &= \, \rho_R(x)l_m(t,x)+h_m(t,x)\mathcal{L}_\alpha\rho_R(x)\\
&\qquad \qquad+ \int_{\R^d}\bigl[h_m(t,x+By)-h_m(t,x)\bigr]
\bigl[\rho_R(x+By)-\rho_R(x) \bigr]\, \nu_\alpha(dy).
\end{split}
\]
Noticing now that $\tilde{l}_{m,R}$ is integrable with integrable Fourier transform, we can  apply the Fourier transform in space to
Equation \eqref{Frozen_PDE_No_Perturb2} in order to write that
\[
\begin{cases}
   \partial_t \widehat{h}_{m,R}(t,p)+\mathcal{F}_x \bigl(\bigl[\mathcal{L}_\alpha + \langle A x, D_{x}\rangle\bigr] h_{m,R}\bigr)(t,p)  \,
   =
   \, -\widehat{\tilde{l}}_{m,R}(t,p); \\
    \widehat{h}_{m,R}(T,p) \, = \, \widehat{u_T}_{m,R}(p).
  \end{cases}
\]
We remember in particular that the above operator $\mathcal{L}_\alpha + \langle A x, D_{x}\rangle$ has an associated L\'evy symbol
$\Phi^{\text{ou}}(p)$ and, following Section $3.3.2$ in \cite{book:Applebaum09}, it holds that
\[\mathcal{F}_x \bigl(\bigl[\mathcal{L}_\alpha + \langle A x, D_{x}\rangle\bigr] h_{m,R}\bigr)(t,p) \, = \,
\Phi^{\text{ou}}(p)\widehat{h}_{m,R}(t,p).\]
We can then use it to show that $\widehat{h}_{m,R}$ is a classical solution of the following equation:
\[
\begin{cases}
   \partial_t \widehat{h}_{m,R}(t,p)+ \Phi^{\text{ou}}(p)\widehat{h}_{m,R}(t,p)  \, = \, -\widehat{\tilde{l}}_{m,R}(t,p); \\
    \widehat{h}_{m,R}(T,p) \, = \, \widehat{u_T}_{m,R}(p).
  \end{cases}
\]
The above equation can be easily solved by integration in time, giving the following representation of $\widehat{h}_{m,R}(t,p)$:
\[\widehat{h}_{m,R}(t,p) \, = \, e^{(T-t)\Phi^{\text{ou}}(p)}\widehat{u_T}_{m,R}(p) + \int_{t}^{T} e^{(s-t)\Phi^{\text{ou}}(p)} \widehat{\tilde{l}}_{m,R} (s,p)\, ds.\]
In order to go back to $\tilde{v}^{\tau,\xi}_m$, we apply now the inverse Fourier transform to write that
\[h_{m,R}(t,x) \, = \, P^{\text{ou}}_{T-t}g_{m,R}(x) + \int_{t}^{T}P^{\text{ou}}_{s-t}\tilde{l}_{m,R}(s,x) \, ds,\]
remembering that $\{P^{\text{ou}}_t\}_{t\ge 0}$ is the convolution Markov semigroup associated with the Ornstein-Uhlenbeck operator
$\mathcal{L}_\alpha + \langle A x, D_{x}\rangle $. Letting $m$ go to $\infty$, it then follows immediately that $g_{m,R}\to u_{T,m}$,
$h_{m,R}\to h_{m}$ and $\tilde{l}_{m,R} \to l_m$. A change of variable allows us to show the Duhamel representation, at least in the regularized setting:
\[\begin{split}
\tilde{v}^{\tau,\xi}_m(t,y)  &= \, P^{\text{ou}}_{T-t}u_{T,m}\Bigl(y+\int_{t}^{T}e^{A(t-s)}F_m(s,\theta_{s,\tau}(\xi))\,
ds\Bigr) \\
&\qquad\qquad \qquad+ \int_{t}^{T}P^{\text{ou}}_{s-t}f_m\Bigl(s,y+\int_{t}^{s}e^{A(t-u)}F_m(u, \theta_{\tau,u}(\xi))\, du\Bigr) \, ds.    
\end{split}
\]
Letting $m$ goes to zero and remembering that $\tilde{v}^{\tau,\xi}_m \to \tilde{v}^{\tau,\xi}$, $f_m \to f$ , $u_{T,m} \to u_T$
and $F_m \to F$ in the right functional spaces, we can conclude that
$\tilde{v}^{\tau,\xi}=\tilde{u}^{\tau,\xi}$.
\end{proof}

\subsection{Expansion along the proxy}

We are going to use now the ``frozen'' IPDE \eqref{Frozen_PDE} in order to derive appropriate quantitative
controls of a solution $u$ of Equation \eqref{Degenerate_Stable_PDE}. Up to now, the freezing parameters $(\tau,\xi)$ were set free but
they will be later chosen appropriately depending on the control we aim to establish.\newline
The main idea is to exploit the Duhamel formula (Proposition \ref{prop:frozen_Duhamel_Formula}) for the proxy to expand any solution $u$ of
the original IPDE \eqref{Degenerate_Stable_PDE} along the proxy.
To make things more precise, let $u$ be a mild solution in $L\bigl(0,T;C^{\alpha+\beta}_{b,d}(\R^{nd})\bigr)$ of IPDE
\eqref{Degenerate_Stable_PDE}. Mollifying if necessary, it is possible to construct three sequences $\{f_m\}_{m\in \N}$, $\{u_{T,m}\}_{m\in \N}$
and $\{F_m\}_{m\in \N}$ of bounded and smooth functions with bounded derivatives such that $f_m \to f$ in
$L^\infty\bigl(0,T;C^\beta_{b,d}(\R^{nd})\bigr)$, $u_{T,m} \to u_T$ in $C^{\alpha+\beta}_{b,d}(\R^{nd})$ and
$\Vert F_m-F\Vert_H \to 0$. Since $u$ is a mild solution of \eqref{Degenerate_Stable_PDE}, we know that there exists a smooth
sequence $\{u_m\}_{m\in \N}$ converging to $u$ in $L\bigl(0,T;C^{\alpha+\beta}_{b,d}(\R^{nd})\bigr)$ and such that for any fixed $m$ in $\N$, $u_m$ solves in a classical sense the ``regularized'' IPDE \eqref{Regularizied_PDE}.\newline
Exploiting now that $F_m$ is bounded and smooth, we can define the ``regularized'' flow $\theta^m_{\cdot,\tau}
(\xi)$ as the \emph{unique} flow satisfying
\begin{equation}\label{eq:def_reg_flow}
\theta^m_{t,\tau}(\xi) \, = \, \xi +
\int_{\tau}^{t}\bigl[A\theta^m_{s,\tau}(\xi)+F_m(s,\theta^m_{s,\tau}(\xi))\bigr] \, ds, \quad t \, \in \, [\tau,T].
\end{equation}
It is then easy to notice that $u_m$ is also a classical solution in $L\bigl(0,T;C^{\alpha+\beta}_{b,d}(\R^{nd})\bigr)$ of
\[\partial_t u_m(t,x) +\mathcal{L}_\alpha u_m(t,x)+ \langle A x  +F_m(t,\theta^m_{t,\tau}(\xi)),
D_{x}u_m(t,x)\rangle \, = \, -\bigl[f_m(t,x)+R_m^{\tau,\xi}(s,x)\bigr] \]
on $(0,T)\times \R^{nd}$ with terminal condition $u_{T,m}$. Above, we have denoted
\begin{equation}\label{eq:def_remainder_regul}
R_m^{\tau,\xi}(t,x)\,:=\, \langle F_m(t,x)-F_m(t,\theta^m_{t,\tau}(\xi)),D_{x}u_m(t,x)
\rangle.
\end{equation}

Since clearly, $R^{\tau,\xi}_m$ is in $L^\infty\bigl(0,T;C^{\alpha+\beta}_{b,d}(\R^{nd})\bigr)$, we can use the Duhamel Formula (Proposition \ref{prop:frozen_Duhamel_Formula}) for the proxy to write that
\[
u_m(t,x) \, = \,  \tilde{P}^{m,\tau,\xi}_{T,t}u_{T,m}(x) +  \int_{t}^{T}\tilde{P}^{m,\tau,\xi}_{s,t}\bigl[f_m(s,x)+
R_m^{\tau,\xi}(s,x)\bigr] \, ds, \quad (t,x) \, \in \, (0,T)\times\R^{nd},
\]
where $\{\tilde{P}^{m,\tau,\xi}_{s,t}\}_{t\le s}$ is the semigroup associated with the operator 
\[\mathcal{L}_\alpha+\langle
Ax+F_m(t,\theta^m_{t,\tau}(\xi)),D_{x} \rangle.\]

The reasoning above is summarized in the following Duhamel-type formula that allows to expand any classical solution $u_m$ of the
``regularized'' IPDE \eqref{Regularizied_PDE} along the ``regularized frozen'' proxy.

\begin{prop}[Duhamel-Type Formula]
\label{prop:Expansion_along_proxy}
Let $(\tau,\xi)$ a freezing couple in $[0,T]\times\R^{nd}$. Under [\textbf{A}], any classical solution $u_m$ of the ``regularized'' IPDE \eqref{Regularizied_PDE} can be represented  as
\begin{equation}\label{eq:Expansion_along_proxy}
u_m(t,x) \, = \, \tilde{u}^{\tau,\xi}_m(t,x) + \int_{t}^{T}  \tilde{P}^{m,\tau,\xi}_{s,t}R^{m,\tau,\xi}(s,x)\,
ds, \quad (t,x) \, \in \, (0,T)\times\R^{nd}
\end{equation}
where $R_m^{\tau,\xi}$ is as in \eqref{eq:def_remainder_regul} and $\tilde{u}^{\tau,\xi}_m$ is defined through the Duhamel
Representation \eqref{Duhamel_representation_of_proxy} with the ``regularized'' coefficients $f_m$, $u_{T,m}$.
\end{prop}

Thanks to the above representation (Equation \eqref{eq:Expansion_along_proxy}), we know that, since we have already shown the suitable control for the frozen solution $u^{\tau,\xi}_m$ (namely, Proposition \ref{prop:Schauder_Estimates_for_proxy} with $f_m,u_{T,m}$), the main term which remains to be investigated in order to show the Schauder Estimates (Theorem \ref{theorem:Schauder_Estimates}) is the remainder
\begin{equation}\label{eq:qqq}
  \int_{t}^{T}\tilde{P}^{m,\tau,\xi}_{s,t}R_m^{\tau,\xi}(s,x) \, ds,
\end{equation}
that represents exactly the error in the expansion along the proxy.

To be precise, we could have passed to the limit in Equation \eqref{eq:Expansion_along_proxy} in order to obtain a similar Duhamel-type
formula for a mild solution $u$ in $L^\infty\bigl([0,T];C^{\alpha+\beta}_{b,d}(\R^{nd})\bigr)$. However, a problem appears when trying to
give a precise meaning at the limit for the remainder contribution \eqref{eq:qqq}. We already know that the limit exists point-wise by
difference, but for our approach to work, we need to establish precise quantitative controls on this term. Such estimates could be obtained
through duality techniques in Besov spaces (cf. Section $5.1$) but only at the expense of fixing
already the freezing couple as $(\tau,\xi)=(t,x)$. The drawback of this method is that it does not allow to differentiate Equation
\eqref{eq:Expansion_along_proxy}, which is needed to estimate $D_{ x_1}u$.

In order to show the suitable estimates for Expression \eqref{eq:qqq}, we will need at first an additional constraint on the behaviour of
the system. In particular, we will say to be under assumption [\textbf{A'}] when assumption [\textbf{A}] is considered and if moreover,
\begin{description}
  \item[{[ST]}] we assume to be in a small time interval, i.e.\ $T\le 1$.
\end{description}

Under these stronger assumptions, we will then be  able to show in Section $5$ below that the following control holds:

\begin{prop}[A Priori Estimates]
\label{prop:A_Priori_Estimates}
Let $u$ be a mild solution of IPDE \eqref{Degenerate_Stable_PDE} in $L^\infty\bigl([0,T];C^{\alpha+\beta}_{b,d}(\R^{nd})\bigr)$. Under
[\textbf{A'}], there exists a constant $C\ge 1$ such that
\begin{multline}\label{eq:A_Priori_Estimates}
\Vert u \Vert_{L^\infty(C^{\alpha+\beta}_{b,d})} \, \le \, Cc_0^{\frac{\beta-\gamma_n}{\alpha}}\bigl[\Vert f
\Vert_{L^\infty(C^{\beta}_{b,d})}+\Vert u_T
\Vert_{C^{\alpha+\beta}_{b,d}}\bigr]\\
+C\bigl(c_0^{\frac{\beta-\gamma_n}{\alpha}}\Vert
F\Vert_H + c_0^{\frac{\alpha+\beta-1}{1+\alpha(n-1)}}\bigr)\Vert u \Vert_{L^\infty(C^{\alpha+\beta}_{b,d})},
\end{multline}
where $c_0 \in (0,1)$ is assumed to be fixed but chosen later.
\end{prop}

We remark already that in the above control, the constants multiplying $\Vert u \Vert_{L^\infty(C^{\alpha+\beta}_{b,d})}$ have to be small
if one wants to derive the expected Schauder estimates. If $c_0$ is small enough, then $Cc_0^{\frac{\alpha+\beta-1}{1+\alpha(n-1)}}$ can be
made smaller than $1/4$. Anyhow, for this chosen small $c_0$, the quantity $c_0^{\frac{\beta-\gamma_n}{\alpha}}$ becomes large and
therefore, it needs to be balanced with $C\Vert F\Vert_H$. Namely, we can conclude if for instance,
$Cc_0^{\frac{\beta-\gamma_n}{\alpha}}\Vert F\Vert_H<1/4$ that implies in particular that $\Vert F\Vert_H$ has to be small with
respect to $c_0$.

\subsection{Conclusion of proof}
In the first part of this section, we prove the Schauder estimates (Theorem \ref{theorem:Schauder_Estimates}) from the A Priori estimates
(Proposition \ref{prop:A_Priori_Estimates}) through a suitable scaling procedure. Roughly speaking, the idea is to start from a general
dynamics and then use the scaling procedure to make the H\"older norm $\Vert F \Vert_H$ small enough in order to make a \emph{circular}
argument work. Again, if $c_0$ and $\Vert F\Vert_H$ are small enough in Equation \eqref{eq:A_Priori_Estimates}, the
$L^\infty\bigl(0,T;C^{\alpha+\beta}_{b,d}(\R^{nd})\bigr)$-norm of $u$ on the right-hand side can be absorbed by the left-hand one. Once the
Schauder estimates \eqref{equation:Schauder_Estimates} hold in the scaled dynamics, we will conclude going back to the original IPDE through
the inverse scaling procedure, even if for a small final time horizon $T$. \newline
The second part of the section focuses on showing how to drop the additional assumption [\textbf{A'}]. The key point here is to proceed
through iteration up to an arbitrary, but finite, given time $T$ thanks to the stability of a solution $u$ in the space
$L^{\infty}\bigl(0,T;C^{\alpha+\beta}_{b,d}(\R^{nd})\bigr)$.

\subsubsection{Scaling argument}
Under [\textbf{A}], we start considering a mild solution $u$ of IPDE \eqref{Degenerate_Stable_PDE} on $[0,T]$ for some final time $T\le
1$ to be fixed later. For a scaling parameter $\lambda$ in $(0,1]$ to be chosen later, we would like to analyze IPDE
\eqref{Degenerate_Stable_PDE} under the change of variables
\begin{equation}\label{eq:change_of_variable}
(t,x) \, \mapsto \, (\lambda t,\mathbb{T}_\lambda x),
\end{equation}
where $\mathbb{T}_\lambda:=\lambda^{1/\alpha}\mathbb{M}_\lambda$. Again, the scaling is performed accordingly to the homogeneity induced by
the distance $\mathbf{d}_P$ in \eqref{Definition_distance_d_P_stable}.\newline
To this purpose, we firstly introduce the scaled solution $u_\lambda$ defined by
\[u_\lambda(t,x) \,:= \, u(\lambda t,\mathbb{T}_\lambda x).\]
It then follows immediately that $u_\lambda$ is a mild solution of the following Equation:
\[\begin{split}
    \lambda^{-1}\partial_tu_\lambda(t,x) +\lambda^{-1}\mathcal{L}_\alpha u_\lambda+\bigl{\langle} A \mathbb{T}_\lambda x +F(\lambda
    t,\mathbb{T}_\lambda x), \mathbb{T}^{-1}_\lambda
    &D_{x}u_\lambda(t,x)\bigr{\rangle}  \\
    &= \, -f(\lambda t,\mathbb{T}_\lambda x) 
    \,\,\text{ on } (0,T_\lambda)\times \R^{nd},
  \end{split}\]
with terminal condition  $u_\lambda(T_\lambda,x) = u_T(\mathbb{T}_\lambda x)$, where $T_\lambda:=T/\lambda$. Since we want the scaled dynamics to satisfy assumption $(\textbf{A}')$, we choose now $T$ so that
$T_\lambda \le 1$. It is important to notice that this is possible since we assumed $\lambda$ to be fixed, even if we have not chosen it
yet. Denoting now
\begin{align*}
  f_\lambda(t,x)\, &:= \, \lambda f(\lambda t,\mathbb{T}_\lambda x); \\
  u_{T,\lambda}(x)\, &:= \, u_T(\mathbb{T}_\lambda x);\\
  A_\lambda\, &:= \, \lambda\mathbb{T}^{-1}_\lambda A\mathbb{T}_\lambda; \\
  F_\lambda(t,x)\, &:= \lambda\mathbb{T}^{-1}_\lambda F(\lambda t,\mathbb{T}_\lambda x),
\end{align*}
we can rewrite the scaled dynamics as:
\begin{equation}
\label{Scaled_Degenerate_Stable_PDE}
\begin{cases}
    \partial_tu_\lambda(t,x) +\bigl{\langle} A_\lambda x +F_\lambda(t,x), D_{x}u_\lambda(t,x) \bigr{\rangle} +
    \mathcal{L}_\alpha u_\lambda(t,x) \, =
    \, -f_\lambda(t,x); \\
    u_\lambda(T_\lambda,x) \, = \, u_{T,\lambda}(x).
  \end{cases}
\end{equation}

To continue, we need the following lemma that shows how the scaling procedure reflects on the norms of the
coefficients. Recalling Equation \eqref{eq:norm_H_for_F} for the definition of $\Vert\cdot \Vert_H$, a direct calculation on the norms leads
to the following result:
\begin{lemma}[Scaling Homogeneity of Norms]
\label{lemma:Scaling_Homogeneity_Norms}
Under [\textbf{A}], it holds that
\begin{align} \label{eq:Scaling_Homogeneity_Norms}
  \Vert F_{\lambda}\Vert_H \, &= \, \lambda^{\beta/\alpha}\Vert F\Vert_H; \notag\\
  \lambda^{\frac{\alpha+\beta}{\alpha}}\Vert f\Vert_{L^\infty(C^\beta_{b,d})}\, \le \, \Vert &f_\lambda\Vert_{L^\infty(C^\beta_{b,d})} \, \le
  \, \Vert f\Vert_{L^\infty(C^\beta_{b,d})};\\
  \lambda^{\frac{\alpha+\beta}{\alpha}}\Vert u_T\Vert_{C^{\alpha+\beta}_{b,d}} \, \le \, \Vert &u_{T,\lambda}\Vert_{C^{\alpha+\beta}_{b,d}} \, \le
  \, \Vert u_T\Vert_{C^{\alpha+\beta}_{b,d}};\notag \\
  \lambda^{\frac{\alpha+\beta}{\alpha}}\Vert u\Vert_{L^\infty(C^{\alpha+\beta}_{b,d})}\, \le \, \Vert
  &u_\lambda\Vert_{L^\infty(C^{\alpha+\beta}_{b,d})} \, \le \, \Vert u\Vert_{L^\infty(C^{\alpha+\beta}_{b,d})}. \notag
\end{align}
\end{lemma}
Since the scaled dynamics in \eqref{Scaled_Degenerate_Stable_PDE} satisfies assumption [\textbf{A'}], we know from Proposition
\ref{prop:A_Priori_Estimates} that the scaled solution $u_\lambda$ satisfies the a priori Estimates  \eqref{eq:A_Priori_Estimates}:
\begin{multline}
\label{eq:A_priori_Estimate_Scaled}
\Vert u_\lambda \Vert_{L^\infty(C^{\alpha+\beta}_{b,d})} \, \le \, Cc_0^{\frac{\beta-\gamma_n}{\alpha}}\bigl[\Vert f_\lambda \Vert_{L^\infty(C^{\beta}_{b,d})}+\Vert u_{T,\lambda}
\Vert_{C^{\alpha+\beta}_{b,d}}\bigr]\\
+C\bigl(c_0^{\frac{\beta-\gamma_n}{\alpha}}\Vert
F_\lambda\Vert_H + c_0^{\frac{\alpha+\beta-1}{1+\alpha(n-1)}}\bigr)\Vert u_\lambda \Vert_{L^\infty(C^{\alpha+\beta}_{b,d})}
\end{multline}
for some constant $c_0$ in $(0,1]$ to be chosen later.\newline
We would like now to exploit a circular argument in order to bring to the left-hand side of Equation \eqref{eq:A_priori_Estimate_Scaled} the term involving $u_\lambda$ on the right-hand one.
To do that, we need to choose properly $\lambda$ and $c_0$ in order to have
\[C\bigr(c_0^{\frac{\beta-\gamma_n}{\alpha}}\Vert F_\lambda\Vert_H + c_0^{\frac{\alpha+\beta-1}{1+\alpha(n-1)}}\bigr)\, < \, 1.\]
This is true if, for example, we choose firstly $c_0$ such that
\[Cc_0^{\frac{\alpha+\beta-1}{1+\alpha(n-1)}} \, = \, \frac{1}{4}\]
and fixed $c_0$, we choose $\lambda$ so that
\[Cc_0^{\frac{\beta-\gamma_n}{\alpha}}\lambda^{\beta/\alpha}\Vert F\Vert_H  \, = \,
Cc_0^{\frac{\beta-\gamma_n}{\alpha}}\Vert F_\lambda\Vert_H \, = \, \frac{1}{4}.\]
With this choice, it then follows from Equation \eqref{eq:A_priori_Estimate_Scaled} that
\[\Vert u_\lambda \Vert_{L^\infty(C^{\alpha+\beta}_{b,d})} \, \le \, 2Cc_0^{\frac{\beta-\gamma_n}{\alpha}}\bigl[\Vert f_\lambda \Vert_{L^\infty(C^{\beta}_{b,d})}+\Vert u_{T,\lambda}
\Vert_{C^{\alpha+\beta}_{b,d}}\bigr].\]
We can finally use Lemma \ref{lemma:Scaling_Homogeneity_Norms} to go back to the original dynamics and write that
\[\Vert u \Vert_{L^\infty(C^{\alpha+\beta}_{b,d})} \, \le \,\lambda^{-\frac{\alpha+\beta}{\alpha}}\Vert u_\lambda
\Vert_{L^\infty(C^{\alpha+\beta}_{b,d})} \, \le \, \overline{C}\bigl[\Vert f
\Vert_{L^\infty(C^{\beta}_{b,d})}+\Vert u_T\Vert_{C^{\alpha+\beta}_{b,d}}\bigr]\]
for some constant $\overline{C}>0$ defined by
\[\overline{C} \, := \,  2\lambda^{-\frac{\alpha+\beta}{\alpha}}Cc_0^{\frac{\beta-\gamma_n}{\alpha}}.\]
\subsubsection{Schauder estimates for general time}
Up to this point, we have assumed to be in a small enough final time horizon (i.e. $T\le 1$) to let our procedure work. We are going now to
extend the Schauder estimates (Equation \eqref{equation:Schauder_Estimates}) to an arbitrary but fixed final time $T_0>0$. Our proof
will consist essentially in a backward iterative procedure through a chain of identical differential dynamics on different, small
enough, time intervals. We recall indeed that the Schauder estimates precisely provide a stability result in the chosen functional space.

\begin{prop}
Under [\textbf{A}], let $T_0>T$ and $u$ a mild solution in $L^{\infty}\bigl(0,T_0;C^{\alpha+\beta}_{b,d}(\R^{nd})\bigr)$ of IPDE
\eqref{Degenerate_Stable_PDE} on $[0,T_0]$ that satisfies the Schauder Estimates \eqref{equation:Schauder_Estimates} on $[0,T]$. Then, there
exists a constant $C_0:=C_0(T_0)$ such that
\[\Vert u \Vert_{L^{\infty}(0,T_0;C^{\alpha+\beta}_{b,d})} \, \le \, C_0\Bigl[\Vert f\Vert_{L^{\infty}( 0,T_0;C^{\beta}_{b,d})} +
\Vert u_T \Vert_{C^{\alpha+\beta}_{b,d}}\Bigr].\]
\end{prop}
\begin{proof}
Fixed $N=\lceil\frac{T_0}{T}\rceil$, we are going to consider a system of $N$ Cauchy problems:
\[
\begin{cases}
    \partial_tu_k(t,x) +\bigl{\langle}A x + F(t,x), D_{x}u_k(t,x)\bigr{\rangle} +
\mathcal{L}_\alpha u_k(t,x) \, = \, -f(t,x); \\
    u_k((1-\frac{k-1}{N})T_0,x) \, = \, u_{k-1}((1-\frac{k-1}{N})T_0,x),
  \end{cases}
\]
on $((1-\frac{k}{N})T_0,(1-\frac{k-1}{N})T_0)\times \R^{nd}$, for $k=1,\dots,N$ with the notation that $u_0(T_0,x)=u_T(x)$. Reasoning iteratively, we find that any mild solution of IPDE
\eqref{Degenerate_Stable_PDE} on $[0,T_0]$ is also a mild solution of any of the equations of the system. Moreover, since any solution $u_k$
is defined on $[(1-\frac{k}{N})T_0,(1-\frac{k-1}{N})T_0]$ and
\[(1-\frac{k-1}{N})T_0-(1-\frac{k}N)T_0 \, = \, \frac{k}{N}T_0 - \frac{k-1}{N}T_0 \, = \, \frac{1}{N}T_0 \, \le \, T,\]
the Schauder estimates (Equation \eqref{equation:Schauder_Estimates}) hold for any solution $u_k$ with terminal condition
$u_{k-1}((1-\frac{k-1}{N})T_0,\cdot)$. In particular,
\[
\begin{split}
\Vert u_k &\Vert_{L^{\infty}((1-\frac{k}{N})T_0,(1-\frac{k-1}{N})T_0;C^{\alpha+\beta}_{b,d})} \\
&\qquad\qquad\le \, C\Bigl[\Vert f
\Vert_{L^{\infty}((1-\frac{k}{N})T_0,(1-\frac{k-1}{N})T_0;C^{\beta}_{b,d})} + \Vert u_{k-1}((1-\frac{k-1}{N})T_0,\cdot)
\Vert_{C^{\alpha+\beta}_{b,d}}\Bigr] \\
&\qquad\qquad\le \, C^2\Bigl[\Vert f \Vert_{L^{\infty}((1-\frac{k}{N})T_0,(1-\frac{k-1}{N})T_0;C^{\beta}_{b,d})} + \Vert f \Vert_{ L^{\infty}((1-
\frac{k-1}{N})T_0,(1-\frac{k-2}{N})T_0;C^{\beta}_{b,d})}\\
&\qquad\qquad\qquad\qquad\qquad \qquad \qquad \qquad \qquad \qquad \qquad + \Vert u_{k-2}((1-\frac{k-2}{N})T_0,\cdot)
\Vert_{C^{\alpha+\beta}_{b,d}}\Bigr]\\
&\qquad\qquad\le \, C^2\Bigl[\Vert f \Vert_{L^{\infty}((1-\frac{k}{N})T_0,(1-\frac{k-2}{N})T_0;C^{\beta}_{b,d})}+\Vert
u_{k-2}((1-\frac{k-2}{N})T_0, \cdot) \Vert_{C^{\alpha+\beta}_{b,d}}\Bigr],
\end{split}
\]
exploiting that $u_{k-1}$ satisfies the Schauder estimates with source $f$ and terminal condition $u_{k-2}((1-\frac{k-2}{N})T_0,\cdot)$. Applying the same procedure
recursively, we finally find that
\[\Vert u_k \Vert_{L^{\infty}((1-\frac{k}{N})T_0,(1-\frac{k-1}{N})T_0;C^{\alpha+\beta}_{b,d})} \, \le \, C^k\Bigl[\Vert f
\Vert_{L^{\infty}((1-\frac{k}{N})T_0,T_0;C^{\beta}_{b,d})} + \Vert u_T \Vert_{C^{\alpha+\beta}_{b,d}}\Bigr].\]
Hence,
\[\Vert u \Vert_{L^{\infty}(0,T_0;C^{\alpha+\beta}_{b,d})} \, \le \, C^N\Bigl[\Vert f\Vert_{L^{\infty}(0,T_0;C^{\beta}_{b,d})} +
\Vert u_T\Vert_{C^{\alpha+\beta}_{b,d}}\Bigr]\]
and we have concluded the proof.
\end{proof}

\setcounter{equation}{0}

\section{Schauder estimates for the proxy}
\fancyhead[RO]{Section \thesection. Schauder Estimates for the Proxy}
The aim of this section is to show how to properly control a solution $\tilde{u}^{\tau,\xi}$ of the ``frozen'' IPDE \eqref{Frozen_PDE} in
order to prove the Schauder estimates (Proposition \ref{prop:Schauder_Estimates_for_proxy}) for the proxy.
We recall the definition of $\tilde{u}^{\tau,\xi}$ through the Duhamel Representation \eqref{Duhamel_representation_of_proxy}. Namely,
for any freezing couple $(\tau,\xi)$ in $[0,T]\times\R^{nd}$, it holds that
\begin{equation}\label{align:Representation}
  \tilde{u}^{\tau,\xi}(t,x) \, = \, \tilde{P}^{\tau,\xi}_{T,t}u_T(x) + \tilde{G}^{\tau,\xi}_{T,t}f(t,x)
\end{equation}
where we have denoted for simplicity with $\{\tilde{G}^{\tau,\xi}_{r,v}\}_{t>v\ge0}$ the family of Green kernels associated with
the frozen density $\tilde{p}^{\tau,\xi}$. More in details, we have for any $v<r$ in $[0,T]$ that
\begin{equation}\label{eq:def_Green_Kernel_stable}
\tilde{G}^{\tau,\xi}_{r,v}f(t,x) \, := \, \int_{v}^{r}\int_{\R^{nd}}\tilde{p}^{\tau,\xi}(t,s,x,y) f(s,y) \,
dy\, ds.
\end{equation}
We can then differentiate the above equation with respect to $ x_1$ in order to obtain an analogous Duhamel type representation for the
derivative $D_{ x_1}\tilde{u}^{\tau,\xi}$:
\begin{equation}\label{align:Representation_deriv}
 D_{ x_1}\tilde{u}^{\tau,\xi}(t,x)\, =
  \,D_{ x_1}\tilde{P}^{\tau,\xi}_{T,t}u_T(x)+D_{ x_1}\tilde{G}^{\tau,\xi}_{T,t}f(t,x).
 \end{equation}
It is then clear that in order to control $\tilde{u}^{\tau,\xi}(t,x)$ in the norm $\Vert \cdot
\Vert_{L^\infty(C^{\alpha+\beta}_{b,d})}$, we can analyze separately the contributions appearing from the frozen semigroup
$\tilde{P}^{\tau,\xi}_{T,t}u_T(x)$ and those from the frozen Green kernel $\tilde{G}^{\tau,\xi}_{T,t}f(t,x)$.

\subsection{First Besov control}
We focus for the moment on the contribution in the Duhamel Representation \eqref{align:Representation} associated with the source $u_T$ that
is, as it will be seen, the more delicate to treat. In the non-degenerate setting (i.e. with respect to $ x_1$), it precisely write:
\[D_{ x_1}\tilde{P}^{\tau,\xi}_{T,t}u_T(x) \, = \, \int_{\R^{nd}} D_{ x_1}\tilde{p}^{\tau,\xi}(t,T,x,y)
u_T(y) \, dy.\]
Looking at the particular structure of $\tilde{p}^{\tau,\xi}$ (cf.\ Equation \eqref{eq:definition_tilde_p}), it can be seen from Lemma
\ref{lemma:Scaling_Lemma} that
\begin{lemma}\label{lemma:link_derivative_density}
Let $i$ in $\llbracket 1,n \rrbracket$. Then, there exist constants $\{C_j\}_{j\in \llbracket i,n \rrbracket}$ such that
\begin{equation}\label{eq:link_derivative_density}
D_{x_i}\tilde{p}^{\tau,\xi}(t,s,x,y) \, = \, \sum_{j=i}^{n}C_j(s-t)^{j-i} D_{ y_j}
\tilde{p}^{\tau,\xi}(t,s,x,y)
\end{equation}
for any $t< s$ in $[0,T]$, any $x,y$ in $\R^{nd}$ and any freezing couple $(\tau,\xi)$ in $[0,T]\times \R^{nd}$.
\end{lemma}

We can now use Wquation \eqref{eq:link_derivative_density} to rewrite $D_{ x_1}\tilde{P}^{\tau,\xi}_{T,t}u_T(x)$ as
\begin{align}
\label{Besov:eq_Introduction}
\bigl{\vert}D_{ x_1}\tilde{P}^{\tau,\xi}_{T,t}u_T(x)\bigr{\vert} \, &= \, \Bigl{\vert}\int_{\R^{nd}}
D_{ x_1}\tilde{p}^{\tau,\xi}(t,T,x,y)
u_T(y) \, dy\Bigr{\vert}\\
&\le \, C\sum_{j=1}^{n}(s-t)^{j-1}\Bigl{\vert}
\int_{\R^{nd}}D_{ y_j}\tilde{p}^{\tau,\xi}(t,T,x,y)
u_T(y) \,dy \Bigr{\vert}.\notag
\end{align}
Remembering that $u_T$ is in $C^{\alpha+\beta}_{b,d}(\R^{nd})$ for $\alpha+\beta>1$ by hypothesis,  we know that it is differentiable with
respect to the first (non-degenerate) variable $ x_1$. Then, the above expression can be controlled easily for $j=1$ as
\[
\begin{split}
\Bigl{\vert}\int_{\R^{nd}}D_{ y_1}\tilde{p}^{\tau,\xi}(t,T,x,y)u_T(y) \, dy \Bigr{\vert} \, &= \,
\Bigl{\vert}\int_{\R^{nd}}\tilde{p}^{\tau,\xi}(t,T,x,y)D_{ y_1}u_T(y) \, dy \Bigr{\vert} \\
&\le \,
\Vert D_{ y_1}u_T \Vert_{L^\infty} \\
&\le \, \Vert u_T \Vert_{C^{\alpha+\beta}_{b,d}},
\end{split}\]
using integration by parts formula. We can then focus on the degenerate components in \eqref{Besov:eq_Introduction}, i.e.
\begin{equation}\label{Besov:eq_Introduction2}
\Bigl{\vert}\int_{\R^{nd}}D_{ y_j}\tilde{p}^{\tau,\xi}(t,T,x,y)u_T(y) \, dy \Bigr{\vert}
\end{equation}
for some $j>1$. Since $u_T$ is not differentiable with respect to $ y_j$ if $j>1$, we cannot apply the same reasoning above but we will
need a more subtle control. Our main idea will be to use the duality in Besov spaces to derive bounds for Expression
\eqref{Besov:eq_Introduction2}. Namely, we introduce for a given $y$ in $\R^d$,
\[y_{\smallsetminus j} \, :=  \, ( y_1,\dots, y_{j-1}, y_{j+1},\dots, y_n) \, \in \, \R^{(n-1)d}.\]
With this definition at hand, we then denote for any function $\phi$ on $\R^{nd}$, the function $\phi(y_{\smallsetminus j},\cdot)$ on $\R^d$ with a slight abuse of notation as
\begin{equation}\label{eq:notation_smallsetminus}
\phi(y_{\smallsetminus j},z) \, := \, \phi( y_1,\dots, y_{j-1},z, y_{j+1},\dots, y_n).
\end{equation}
The key point now is to control the H\"older modulus of $u_T(y_{\smallsetminus j},\cdot)$ on $\R^d$, uniformly in $y_{\smallsetminus j}
\in \R^{(n-1)d}$. To do so, we will need the identification $C^{\alpha_j+\beta_j}_b(\R^d) = B^{\alpha_j+\beta_j}_{\infty,\infty}
(\R^d)$ with the usual notations for the Besov spaces.

We recall now some useful definitions/characterizations about Besov spaces $B^{\tilde{\gamma}}_{p,q}(\R^d)$. For a more detailed
analysis of this argument, we suggest the reader to see Section $2.6.4$ of Triebel \cite{book:Triebel83}. For $\tilde{\gamma}$ in $(0,1)$,
$q,p$ in $(0,+\infty]$, we define the Besov space of indexes $(\tilde{\gamma},p,q)$ on $\R^d$ as:
\[B^{\tilde{\gamma}}_{p,q}(\R^d):= \{f \in \mathcal{S}'(\R^d)\colon \Vert f \Vert_{\mathcal{H}^{\tilde{\gamma}}_{p,q}} \, < + \infty\},\]
where $\mathcal{S}(\R^d)$ denotes the Schwartz class on $\R^d$ and
\begin{equation}\label{alpha-thermic_Characterization}
\Vert f \Vert_{\mathcal{H}^{\tilde{\gamma}}_{p,q}} \, := \, \Vert (\phi_0\hat{f})^\vee \Vert_{L^p}+ \Bigl(\int_{0}^{1}
v^{(1-\frac{\tilde{\gamma}}{\alpha})q}\Vert \partial_vp_h(v,\cdot)\ast f \Vert^q_{L^p} \, \frac{dv}{v}\Bigr)^{\frac{1}{q}},
\end{equation}
with $\phi_0$ a function in $C^\infty_0(\R^d)$ such that $\phi_0(0) \neq 0$ and $p_h$ the isotropic $\alpha$-stable heat kernel on
$\R^d$, i.e.\ the stable density on $\R^d$ whose L\'evy symbol is equivalent to $\vert \lambda \vert^\alpha$. \newline
We point out that the quantities in \eqref{alpha-thermic_Characterization} are well-defined for any $q\neq +\infty$. The modifications for
$q=+\infty$ are obvious and can be written passing to the limit. The previous definition of $B^{\tilde{\gamma}}_{p,q}(\R^d)$ is known as the
stable thermic characterization of Besov spaces and it is particularly adapted to our framework. By a little abuse of notation, we will
write $\Vert f \Vert_{B^{\tilde{\gamma}}_{p,q}}:=\Vert f \Vert_{\mathcal{H}^{\tilde{\gamma}}_{p,q}}$ when this quantity is finite.

For the heat-kernel $p_h$, it is possible to show an improvement of the smoothing effect (cf. Equation \eqref{Smoothing_effect_of_S}),
due essentially to its better decay at infinity. Namely, we are no more bounded to the condition $\gamma<\alpha$ but we can integrate up to
an order $\gamma$ strictly smaller than $1+\alpha$.

\begin{lemma}[Smoothing Effect of the Isotropic Stable Heat-Kernel]
\label{lemma:Smoothing_effect_of_Heat_Kern}
Let $l$ be in $\{1,2\}$ and $\gamma$ in $[0,1+\alpha)$. Then, there exists a positive constant $C:=C(\gamma)$ such that
\begin{equation}\label{Smoothing_effect_of_Heat_Kern}
\int_{\R^d}\vert y \vert^\gamma \vert \partial_vD^l_yp_h(v,y) \vert \, dy \, \le \, Ct^{\frac{\gamma-l}{\alpha}-1}.
\end{equation}
\end{lemma}
A proof of the above result can be derived using the estimates of Kolokoltsov \cite{Kolokoltsov00} (see also \cite{Bogdan:Jakubowski07}).

As already indicated before, it can be seen from the Thermic Characterization \eqref{alpha-thermic_Characterization} that
\begin{equation}\label{Besov:ident_Holder_Besov}
C^{\tilde{\gamma}}_b(\R^d) \, = \, B^{\tilde{\gamma}}_{\infty,\infty}(\R^d).
\end{equation}
Moreover, it is well known (see for example Proposition $3.6$ in \cite{book:Lemarie-Rieusset02}) that
$B^{\tilde{\gamma}}_{\infty,\infty}(\R^d)$
and $B^{-\tilde{\gamma}}_{1,1}(\R^d)$ are in duality. Namely, it holds
\begin{equation}\label{Besov:duality_in_Besov}
\bigl{\vert}\int_{\R^d} fg \, dx \bigr{\vert}\, \le \, C\Vert f \Vert_{B^{\tilde{\gamma}}_{\infty,\infty}}\Vert u_T
\Vert_{B^{-\tilde{\gamma}}_{1,1}},
\end{equation}
for any $f$ in $B^{\tilde{\gamma}}_{\infty,\infty}(\R^d)$ and any $u_T$ in $B^{-\tilde{\gamma}}_{1,1}(\R^d)$.

With these definitions and properties at hand, we can now go back at Expression \eqref{Besov:eq_Introduction2} to write that
\[
\begin{split}
\Bigl{\vert}\int_{\R^{nd}}D_{ y_j}\tilde{p}^{\tau,\xi}(t,T,x,y)&u_T(y) \, dy \Bigr{\vert} \, \le \,
\int_{\R^{(n-1)d}}\Bigl{\vert}D_{ y_j}\tilde{p}^{\tau,\xi}(t,T,x,y)u_T(y)d y_j\Bigr{\vert} dy_{\smallsetminus
j} \\
&\le \,\int_{\R^{(n-1)d}}\Bigl{\Vert}D_{ y_j}\tilde{p}^{\tau,\xi} (t,T,x,y_{\smallsetminus j},\cdot) \Bigr{\Vert}_{
B^{-(\alpha_j+\beta_j)}_{1,1}}\Bigl{\Vert}u_T(y_{\smallsetminus j},\cdot)\Bigr{\Vert}_{B^{\alpha_j+\beta_j}_{\infty,\infty}}
dy_{\smallsetminus j} \\
&\le \, \Vert u_T \Vert_{C^{\alpha+\beta}_{b,d}}\int_{\R^{(n-1)d}}\Bigl{\Vert}D_{ y_j}\tilde{p}^{\tau,\xi}(t,T,x,
y_{\smallsetminus j},\cdot) \Bigr{\Vert}_{B^{-(\alpha_j+\beta_j)}_{1,1}} \, dy_{\smallsetminus j}.
\end{split}
\]
In order to control the above quantities, we will then need a control on the integral of the Besov norms of the derivatives of the proxy.
Since however an additional derivative with respect to $ x_1$ will often appear, for example in
Equation \eqref{Proof:H\"older_Frozen_Semigroup_Idj} below, we state the following result in a more general way.

\begin{lemma}[First Besov Control]
\label{lemma:First_Besov_COntrols}
Let $j$ be  in $\llbracket 2,n\rrbracket$ and $l\in \{0,1\}$. Under [\textbf{A}], there exists a constant $C:=C(j,l)$ such that
\[\int_{\R^{(n-1)d}}\Bigl{\Vert}D_{ y_j}D^{l}_{x_1}\tilde{p}^{\tau,\xi}(t,s,x,y_{\smallsetminus
j},\cdot)\Bigr{\Vert}_{B^{-(\alpha_j+\beta_j
)}_{1,1}} \, dy_{\smallsetminus j} \, \le \, C(s-t)^{\frac{\alpha+\beta}{\alpha}-\frac{1}{\alpha_j}-\frac{l}{\alpha}},\]
for any $t<s$ in $[0,T]$, any $x$ in $\R^{nd}$ and any frozen couple $(\tau,\xi)$ in $[0,T]\times \R^{nd}$.
\end{lemma}
\begin{proof}
To control the Besov norm in $B^{-(\alpha_j+\beta_j)}_{1,1}(\R^d)$, we are going to use the Thermic Characterization
\eqref{alpha-thermic_Characterization} with $\tilde{\gamma}=-(\alpha_j+\beta_j)$.
We start considering the second term in the characterization, i.e.
\[\int_{0}^{1}v^{\frac{\alpha_j+\beta_j}{\alpha}}\int_{\R^d}\Bigl{\vert}\int_{\R^d}\partial_vp_h(v,z- y_j)D_{ y_j}
D^{l}_{x_1}\tilde{p}^{\tau,\xi}(t, s,x,y)\,d y_j  \Bigr{\vert} \, dzdv.\]
Fixed a constant $\delta_j\ge1$ to be chosen later, we split the integral with respect to $v$ in two components:
\[
\begin{split}
&\Vert D_{ y_j}D^{l}_{x_1}\tilde{p}^{\tau,\xi}(t,s,x,y_{\smallsetminus j},\cdot)
\Vert_{B^{-(\alpha_j+\beta_j)}_{1,1}} \\
&\qquad \qquad = \,\int_{0}^{(s-t)^{\delta_j}}v^{\frac{\alpha_j+\beta_j}{\alpha}}\int_{\R^d}\Bigl{\vert}\int_{\R^d} \partial_v p_h(v,z- y_j)
D_{ y_j} D^{l}_{x_1}\tilde{p}^{\tau,\xi}(t,s,x,y) \,d y_j  \Bigr{\vert} \, dzdv \\
& \qquad \qquad \qquad \qquad +\int_{(s-t)^{\delta_j}}^{1}v^{\frac{\alpha_j+\beta_j}{\alpha}}\int_{\R^d}\Bigl{\vert}\int_{\R^d}\partial_vp_h(v,z- y_j)
D_{ y_j}D^{l}_{x_1}\tilde{p}^{\tau,\xi}(t,s,x,y)\,d y_j  \Bigr{\vert} \, dzdv \\
&\qquad \qquad =: \,
\bigl(I_1 +I_2\bigr)(y_{\smallsetminus j}).
\end{split}
\]
The second component $I_2$ has no time-singularity and can be easily controlled by
\[I_2(y_{\smallsetminus j}) \, = \,\int_{(s-t)^{\delta_j}}^{1}v^{\frac{\alpha_j+\beta_j}{\alpha}}\int_{\R^d} \Bigl{\vert}
\int_{\R^d}D_z\partial_vp_h(v,z- y_j)\otimes D^{l}_{x_1}\tilde{p}^{\tau,\xi}(t,s,x, y)\,d y_j
\Bigr{\vert} \, dzdv,\]
using integration by parts formula and noticing that $D_{ y_j}p_h(v,z- y_j)=-D_zp_h(v,z- y_j)$. Then,
\[I_2(y_{\smallsetminus j}) \, \le \, \int_{(s-t)^{\delta_j}}^{1}v^{\frac{\alpha_j+\beta_j}{\alpha}}\int_{
\R^d}\int_{\R^d}\vert D_z\partial_vp_h(v,z- y_j)\vert \, \vert D^{l}_{x_1}\tilde{p}^{\tau,\xi}(t,s,x, y) \vert
\,d y_j  \, dzdv.\]
We can then use Fubini theorem to separate the integrals and apply the smoothing effect of the heat-kernel $p_h$ (Lemma
\ref{lemma:Smoothing_effect_of_Heat_Kern}) to show that
\[
\begin{split}
I_2(y_{\smallsetminus j}) \, &\le \, \int_{(s-t)^{\delta_j}}^{1}v^{\frac{\alpha_j+\beta_j}{\alpha}}\int_{\R^d}\Bigl(\int_{
\R^d}\vert D_z\partial_vp_h(v,z- y_j)\vert \, dz\Bigr)  \vert D^{l}_{x_1}\tilde{p}^{\tau,\xi}(t,s,x,
y) \vert \, d y_jdv \\
&\le \, C\Bigl(\int_{(s-t)^{\delta_j}}^{1}v^{\frac{\alpha_j+\beta_j-1}{\alpha}-1} \, dv\Bigr)\Bigl(\int_{\R^d}\vert D^{l}_{x_1}
\tilde{p}^{\tau,\xi}(t,s,x,y)\vert \,d y_j\Bigr) \\
&\le \, C(s-t)^{\frac{\delta_j(\alpha_j+\beta_j-1)}{\alpha}}\int_{\R^d}\vert
D^{l}_{x_1}\tilde{p}^{\tau,\xi}(t,s,x,y)\vert \,d y_j.
\end{split}\]
Using the smoothing effect (Equation \eqref{eq:Smoothing_effects_of_tilde_p}) of the frozen density $\tilde{p}^{\tau,\xi}$, we have
thus found that
\begin{equation}\label{Proof:First_Besov_COntrol_I2}
\begin{split}
\int_{\R^{(n-1)d}}I_2(y_{\smallsetminus j}) \, dy_{\smallsetminus j} \, &\le \,
(s-t)^{\frac{\delta_j(\alpha_j+\beta_j-1)}{\alpha}}
\int_{\R^{nd}}\vert D^{l}_{x_1}\tilde{p}^{\tau,\xi}(t,s,x,y)\vert\,dy \\
&\le \,
C(s-t)^{\frac{\delta_j(\alpha_j+\beta_j-1)-l}{\alpha}}.
\end{split}
\end{equation}
On the other hand, the term $I_1$ needs a more delicate treatment in order to avoid time-integrability problems. We start using a
cancellation argument with respect to the derivative $\partial_vp_h$ of the heat-kernel to rewrite $I_1$ as
\[\begin{split}
I_1(y_{\smallsetminus j}) \, &= \, \int_{0}^{(s-t)^{\delta_j}}v^{\frac{\alpha_j+\beta_j}{\alpha}}\int_{\R^d}\Bigl{\vert}\int_{\R^d}\partial_vp_h(v,z- y_j) \\
&\qquad\qquad\qquad\qquad\times 
\bigl[D_{ y_j}D^l_{ x_1}\tilde{p}^{\tau,\xi}(t,s,x,y)-D_{ y_j}D^l_{ x_1}\tilde{p}^{\tau,\xi}(t,
s,x,y_{\smallsetminus j},z)\bigr]\,d y_j  \Bigr{\vert} \,dzdv \\
&= \, \int_{0}^{(s-t)^{\delta_j}} v^{\frac{\alpha_j+\beta_j}{\alpha}}\int_{\R^d}\Bigl{\vert}\int_{\R^d}D_z\partial_vp_h(v,z- y_j)\\
&\qquad\qquad\qquad\qquad\otimes\bigl[D^l_{ x_1}\tilde{p}^{\tau,\xi}(t,s,x,y)-D^l_{ x_1}\tilde{p}^{\tau,\xi}(t,s,x,
y_{\smallsetminus j},z)\bigr]\,d y_j  \Bigr{\vert} \,dzdv,
\end{split}\]
where in the second passage we used again integration by parts formula to move the derivative to $p_h$ and the equality $D_{ y_j}p_h(v,
z- y_j)=-D_zp_h(v,z- y_j)$. We can then apply a Taylor expansion with respect to variable $ y_j$ in order to write that
\[\begin{split}
I_1(y_{\smallsetminus j}) \, &= \,\int_{0}^{(s-t)^{\delta_i}}  v^{\frac{\alpha_j+\beta_j}{\alpha}}\int_{\R^d}\Bigl{\vert}\int_{\R^d}D_z\partial_vp_h(v,z- y_j)\\
&\qquad\qquad\qquad\times
\int_{0}^{1} D_{ y_j}D^l_{ x_1}\tilde{p}^{\tau,\xi}(t,s,x,y_{\smallsetminus
j}, y_j+\lambda(z- y_j))\cdot(z- y_j) \,  d\mu d y_j \Bigr{\vert} \, dzdv \\
&\le \, \int_{0}^{(s-t)^{\delta_i}}  v^{\frac{\alpha_j+\beta_j}{\alpha}}\int_{\R^d}\int_{\R^d}\int_{0}^{1}\vert
D_z\partial_vp_h(v,z- y_j)\vert\\
&\qquad\qquad\qquad\times\vert D_{ y_j}D^l_{ x_1}\tilde{p}^{\tau,\xi}(t,s,x,
y_{\smallsetminus j}, y_j+\lambda(z- y_j))\vert\, \vert z- y_j\vert \,  d\lambda d y_jdzdv.
\end{split}
\]
We can then use the Fubini theorem and the changes of variables $\tilde{z}=z- y_j$ (fixed $ y_j$) and
$\tilde{y}_j= y_j+\lambda\tilde{z}$ (considering $\tilde{z}$ and $\lambda$ fixed) to separate the integrals so that
\begin{multline*}
I_1(y_{\smallsetminus j}) \, \le \, \int_{0}^{(s-t)^{\delta_i}}v^{\frac{\alpha_j+\beta_j}{\alpha}} \Bigl(\int_{\R^d} \vert D_z\partial_vp_h(v,\tilde{z})\vert \,
\vert \tilde{z} \vert \, dz\Bigr)\\
\times\Bigl(\int_{\R^d} \vert D_{ y_j}D^l_{ x_1}\tilde{p}^{\tau,\xi}(t,s,x, y_{\smallsetminus j },\tilde{y}_j) \vert\,d y_j\Bigr)\, dv.    
\end{multline*}
The smoothing effect of the heat-kernel $p_h$ (Lemma \ref{lemma:Smoothing_effect_of_Heat_Kern}) allows now to control the first term:
\[
\begin{split}
I_1(y_{\smallsetminus j}) \, &\le \, C\Bigl(\int_{0}^{(s-t)^{\delta_i}}v^{\frac{\alpha_j+\beta_j-1}{\alpha}}\, dv\Bigr)
\Bigl(\int_{\R^d} \vert D_{ y_j}D^l_{ x_1}\tilde{p}^{\tau,\xi}(t,s,x, y_{\smallsetminus j },z+\lambda( y_j-z))
\vert\,d y_j\Bigr) \\
&\le \, C(s-t)^{\delta_j\frac{\alpha_j+\beta_j}{\alpha}} \int_{\R^d}\vert D_{ y_j}D^l_{ x_1}\tilde{p}^{\tau,\xi}
(t,s,x, y_{\smallsetminus j },z+\lambda( y_j-z)) \vert \,d y_j.
\end{split}
\]
It then follows using the smoothing effect of the frozen semigroup (Lemma \ref{lemma:Smoothing_effect_frozen}) that
\begin{align}
\int_{\R^{(n-1)d}}I_1(y_{\smallsetminus j}) \, dy_{\smallsetminus j} \, &\le \,
C(s-t)^{\delta_j\frac{\alpha_j+\beta_j}{\alpha}}\int_{\R^{nd}}\vert D_{ y_j}D^l_{ x_1}\tilde{p}^{\tau,\xi}
(t,s,x, y_{\smallsetminus j },z+\lambda( y_j-z)) \vert \,dy \notag \\
&\le \, C(s-t)^{\delta_j\frac{\alpha_j+\beta_j}{\alpha}-\frac{l}{\alpha} -\frac{l}{\alpha_j}}. \label{Proof:First_Besov_COntrol_I1}
\end{align}
Going back to equations \eqref{Proof:First_Besov_COntrol_I2} and \eqref{Proof:First_Besov_COntrol_I1}, we notice that we need $\delta_j$ to
be such that
\[\delta_j\bigl[\frac{\alpha_j +\beta_j}{\alpha}\bigr] \, = \, \frac{\alpha+\beta}{\alpha} \,\, \text{ and } \,\,
\delta_j\bigl[\frac{\alpha_j+\beta_j-1}{\alpha}\bigr] \, = \, \frac{\alpha+\beta}{\alpha}-\frac{1}{\alpha_j}.\]
Recalling Equation \eqref{eq:def_alpha_i_and_beta_i} for the relative definitions, we can thus conclude choosing
\[\delta_j=(\alpha+\beta)/(\alpha_j+\beta_j)=1+\alpha(j-1).\]
Reproducing the previous computations, we can also write for a test function in $\phi_0$ in $C^\infty_0(\R^d)$,
\[
\begin{split}
\int_{\R^{(n-1)d}}\Bigl{\Vert} \bigl(\phi_0\bigl(D_{ y_j}D^{l}_{x_1}\tilde{p}^{\tau,\xi}(t&,s,x,y_{\smallsetminus
j},\cdot)\bigr)\hat{\phantom{A}}\bigr)^\vee \Bigr{\Vert}_{L^1} \, dy_{\smallsetminus j} \\
&= \, \int_{\R^{(n-1)d}}\int_{\R^d}\Bigl{\vert}\int_{\R^d}D_{ y_j}\hat{\phi_0}(z- y_j)\cdot
D^l_{ x_1}\tilde{p}^{\tau,\xi}(t,s,x,y)\, d y_j \Bigr{\vert} \, dzdy_{\smallsetminus j} \\
&\le \, C \int_{\R^{nd}}\vert D^l_{ x_1}\tilde{p}^{\tau,\xi}(t,s,x,y) \vert \, dy \\
&\le \,
C(s-t)^{-\frac{l}{\alpha}}.
\end{split}
\]
The proof is thus concluded.
\end{proof}

\subsection{Proof of Proposition \ref{prop:Schauder_Estimates_for_proxy}}

Thanks to the first Besov control (Lemma \ref{lemma:First_Besov_COntrols}), we are now ready to prove the Schauder estimates for the proxy
(Proposition \ref{prop:Schauder_Estimates_for_proxy}). Such a proof will be divided in three parts: the estimates for the supremum norms of
the solution and its non-degenerate gradient are stated in Lemma \ref{lemma:supremum_norm_proxy} while the controls of the H\"older moduli
of the solution and its gradient with respect to the non-degenerate variable are given in Lemmas \ref{lemma:Holder_modulus_proxy_Non-Deg}
and \ref{lemma:Holder_modulus_proxy_Deg}, respectively.

\begin{lemma}(Controls on Supremum Norm)
\label{lemma:supremum_norm_proxy}
Under [\textbf{A}], there exists a constant $C:=C(T)\ge 1$ such that for any freezing couple $(\tau,\xi)$ in
$[0,T]\times\R^{nd}$, any $t$ in $[0,T]$ and any $x$ in $\R^{nd}$,
\[\vert \tilde{u}^{\tau,\xi}(t,x) \vert + \vert D_{ x_1}\tilde{u}^{\tau,\xi}(t,x)\vert \, \le \,
C\Bigl[\Vert f \Vert_{L^\infty(C^\beta_{b,d})} +\Vert u_T \Vert_{C^{\alpha+\beta}_{b,d}}\Bigr].\]
\end{lemma}
\begin{proof}
We start noticing that $\tilde{P}^{\tau,\xi}_{T,t}u_T(x)$ and $\tilde{G}^{\tau,\xi}_{T,t}f(t,x)$ can be easily
bounded using the supremum norm of $f$ and $u_T$, respectively. \newline
Moreover, we can use the estimates on the frozen semigroup (Equation \eqref{eq:Control_of_semigroup}) to control
$D_{ x_1}\tilde{G}^{\tau,\xi}_{T,t}f(t,x)$. Indeed,
\[\begin{split}
    \bigl{\vert} D_{ x_1}\tilde{G}^{\tau,\xi}_{T,t}f(t,x) \bigr{\vert} \, &\le \, \int_{t}^{T}
\bigl{\vert}D_{ x_1}\tilde{P}^{\tau,\xi}_{s,t}
f(s,x) \bigr{\vert} \, ds \\
&\le \, C (T-t)^{\frac{\alpha+\beta-1}{\alpha}}\Vert f\Vert_{L^\infty(C^\beta_{b,d})} \\ 
&\le \, C
T^{\frac{\alpha+\beta-1}{\alpha}}\Vert f\Vert_{L^\infty (C^\beta_{b,d})},
\end{split}
\]
remembering in the last inequality that $\alpha + \beta -1 >0$ by hypothesis [\textbf{P}]. \newline
It remains to control $D_{ x_1}\tilde{P}^{\tau,\xi}_{T,t}u_T(x)$. As done in the previous Sub-section $4.1$, we start using
the scaling lemma \ref{lemma:link_derivative_density} to write that
\[
\begin{split}
\bigl{\vert} D_{ x_1}\tilde{P}^{\tau,\xi}_{T,t}u_T(x) \bigr{\vert}\, &= \,
\Bigl{\vert}\int_{\R^{nd}}D_{ x_1}\tilde{p}^{\tau,\xi}(t,T,x,y)u_T(y)
\,dy\Bigr{\vert} \\
&\le \, C\sum_{j=1}^{n}(T-t)^{j-1} \Bigl{\vert} \int_{\R^{nd}} D_{ y_j} \tilde{p}^{\tau,\xi} (t,T,x,y)
u_T(y) \, dy \Bigr{\vert}\\
&=: \, C\sum_{j=1}^{n}(T-t)^{j-1}J_j.
\end{split}\]
Since $u_T$ is differentiable in the first, non-degenerate variable $ x_1$, the contribution $J_1$ can be easily bounded using integration
by parts formula:
\begin{equation}\label{Proof:Supremum_for_Frozen_Semigroup_J1}
J_1\, = \,\Bigl{\vert} \int_{\R^{nd}}\tilde{p}^{\tau,\xi}(t,T,x,y)D_{ y_1}u_T(y) \, dy
\Bigr{\vert} \, \le \,\Vert D_{ y_1}u_T\Vert_{L^{\infty}} \, \le \,\Vert u_T\Vert_{C^{\alpha+\beta}_{b,d}}.
\end{equation}
To control the other terms $J_j$ for $j>1$, we use instead the duality in Besov spaces (Equation \eqref{Besov:duality_in_Besov}) and
Identification \eqref{Besov:ident_Holder_Besov}, so that
\begin{equation}
\label{Proof:Supremum_for_Frozen_Semigroup_Jj}
\begin{split}
J_j \, &\le \, C\Vert u_T \Vert_{C^{\alpha+\beta}_{b,d}}\int_{\R^{(n-1)d}}  \Vert D_{ y_j}\tilde{p}^{\tau,\xi}(t,T,x,
y_{\smallsetminus j},\cdot)\Vert_{B^{-(\alpha_j+\beta_j)}_{1,1}} \, dy_{\smallsetminus j}\\
&\le \, C\Vert u_T
\Vert_{C^{\alpha+\beta}_{b,d}}(T-t)^{\frac{\alpha+\beta}{\alpha}-\frac{1}{\alpha_j}},
\end{split}
\end{equation}
where in the last inequality we applied the first Besov control (Lemma \ref{lemma:First_Besov_COntrols}).\newline
Looking back at Equations \eqref{Proof:Supremum_for_Frozen_Semigroup_J1}-\eqref{Proof:Supremum_for_Frozen_Semigroup_Jj}, it finally
holds that
\[
\begin{split}
\bigl{\vert} D_{ x_1}\tilde{P}^{\tau,\xi}_{T,t}u_T(x) \bigr{\vert} \, &\le \, C\Vert u_T\Vert_{C^{\alpha+\beta}_{b,d}}
\bigl(1+\sum_{j=2}^{n}(T-t)^{j-1}(T-t)^{\frac{\alpha+\beta}{\alpha}-\frac{1}{\alpha_j}}\bigr) \\
&\le \,
C\bigl(1+T^{\frac{\alpha+\beta-1}{\alpha}}\bigr)\Vert u_T \Vert_{C^{\alpha+\beta}_{b,d}},
\end{split}\]
where in the last passage we used again that $\alpha+\beta -1 >0$ by hypothesis [\textbf{P}].
\end{proof}

Before starting with the calculations on the H\"older modulus, we will need to
distinguish two cases. Fixed $(t,x,x')$ in $[0,T]\times \R^{2nd}$, we will say that the \emph{off-diagonal regime} holds if $T-t \le c_0\mathbf{d}^\alpha(x,x')$ for a constant $c_0$
to be specified but meant to be smaller than $1$. This means in particular that the spatial distance is larger than the characteristic
time-scale up to the prescribed constant $c_0$ which will be useful further on in the computations for a circular argument. \newline
On the other hand, we will say that a \emph{global diagonal regime} is in force if $T-t \ge c_0\mathbf{d}^\alpha(x,x')$ and the spatial
points are instead closer than the typical time-scale magnitude. In particular, when a time integration is involved (for example, in the
control of the frozen Green kernel), the same two regime appears even in a local base. Considering a
variable $s$ in $[t,T]$, there are again a \emph{local off-diagonal regime} if $s-t \le c_0\mathbf{d}^\alpha(x,x')$ and a \emph{local
diagonal regime} when $s-t \ge c_0\mathbf{d}^\alpha(x,x')$. In particular, we will denote with $t_0$  the critical time at
which a change of regime occurs in the globally diagonal regime. Namely,
\begin{equation}\label{eq:def_t0}
t_0 \, := \, \bigl(t+c_0\mathbf{d}^\alpha(x,x')\bigr)\wedge T.
\end{equation}
We highlight however that this approach was already used in \cite{Chaudru:Honore:Menozzi18_Sharp} to obtain Schauder estimates for
degenerate Kolmogorov equations and can be adapted in the current setting.

Moreover, it is important to notice that the norm $\Vert \cdot\Vert_{C^{\alpha+\beta}_d}$ is essentially defined as the sum of the norms
$\Vert \cdot \Vert_{C^{\frac{\alpha+\beta}{1+\alpha(i-1)}}}$ with respect to the $i$-th variable and uniformly on the other components.
Thus, there is a big difference between the case $i=1$ where
$\alpha+\beta$ is in $(1,2)$ and we have to deal with a proper derivative and the other situations ($i>1$) where instead
$(\alpha+\beta)/(1+\alpha(i-1))<1$ and the norm is calculated directly on the function. For this reason, we are going to analyze the two
cases separately. Lemma \ref{lemma:Holder_modulus_proxy_Non-Deg} will work on the non-degenerate setting ($i=1$) while Lemma
\ref{lemma:Holder_modulus_proxy_Deg} will concern the degenerate one ($i>1$).

\begin{lemma}[Controls on H\"older Moduli: Non-Degenerate]
\label{lemma:Holder_modulus_proxy_Non-Deg}
Let $x,x'$ be in $\R^{nd}$ such that $ x_j= x'_j$ for any $j\neq 1$. Under [\textbf{A}], there exists a constant $C\ge 1$
such that for any $t$ in $[0,T]$ and any freezing couple $(\tau,\xi)$ in $[0,T]\times \R^{nd}$, it holds that
\[\bigl{\vert} D_{ x_1}\tilde{u}^{\tau,\xi}(t,x)- D_{ x_1}\tilde{u}^{\tau,\xi}(t,x')\vert\\
\le \, Cc^{\frac{\alpha+\beta-2}{\alpha}}_0\bigl(\Vert u_T\Vert_{C^{\alpha+\beta}_{b,d}} + \Vert f \Vert_{L^\infty(C^\beta_{b,d})}\bigr)
\mathbf{d}^{\alpha+\beta-1}(x,x').\]
\end{lemma}

Before proving the above result, we point out the control on the H\"older modulus of $\tilde{u}^{\tau,\xi}$ with respect to the
degenerate variables ($i>1$):

\begin{lemma}[Controls on H\"older Moduli: Degenerate]
\label{lemma:Holder_modulus_proxy_Deg}
Let $i$ be in $\llbracket 2,n\rrbracket$ and $x,x'$ in $\R^{nd}$ such that $ x_j= x'_j$ for any $j\neq i$. Under
[\textbf{A}], there exists a constant $C:=C(i)$ such that for any $t$ in $[0,T]$ and any freezing couple $(\tau,\xi)$ in $[0,T]\times
\R^{nd}$, it holds that
\[\bigl{\vert} \tilde{u}^{\tau,\xi}(t,x)- \tilde{u}^{\tau,\xi}(t,x') \vert\, \le \,
Cc^{\frac{\beta-\gamma_i}{\alpha}}_0\bigl(\Vert u_T\Vert_{C^{\alpha+\beta}_{b,d}} + \Vert f \Vert_{L^\infty(C^\beta_{b,d})}\bigr) \mathbf{d}^{\alpha+\beta}(x,x').\]
\end{lemma}
\paragraph{Proof of Lemma \ref{lemma:Holder_modulus_proxy_Non-Deg}}
\emph{Controls on frozen semigroup.} Let us consider firstly the off-diagonal regime, i.e.\ the case $T-t \le c_0\mathbf{d}^\alpha(x,x')$.
Using the scaling lemma
\ref{lemma:link_derivative_density}, it holds that
\[
\begin{split}
D_{ x_1}\tilde{P}^{\tau,\xi}_{T,t}u_T(x) \, &= \,
\int_{\R^{nd}}D_{ x_1}\tilde{p}^{\tau,\xi}(t,T,x,y)u_T(y) \, dy \\ 
&= \,
\sum_{j=1}^{n}C_j(T-t)^{j-1}\int_{\R^{nd}}D_{ y_j}\tilde{p}^{\tau,\xi}(t,T,x,y)u_T(y) \, dy.    
\end{split}
\]
It then follows that
\begin{align}
&\bigl{\vert}D_{ x_1}\tilde{P}^{\tau,\xi}_{T,t}u_T(x) - D_{ x_1}\tilde{P}^{\tau,\xi}_{T,t}u_T(x')
\bigr{\vert} \notag \\
& \qquad \qquad\qquad \,\,\le \,C\sum_{j=1}^n(T-t)^{j-1}\Bigl{\vert}\int_{\R^{nd}}\bigl[D_{ y_j}\tilde{p}^{\tau, \xi}
(t,T,x,y)-D_{ y_j}\tilde{p}^{\tau,\xi}(t,T,x',y) \bigr] u_T(y) \,
dy\Bigr{\vert}\notag \\
& \qquad \qquad\qquad \,\, =:\, C\sum_{j=1}^n (T-t)^{j-1} I^{od}_j. \label{Proof:H\"older_for_Frozen_Semigroup_Decomp}
\end{align}
We are going to treat separately the cases $j=1$ and $j>1$ for the \emph{off-diagonal} contributions $\{I^{od}_j\}_{j\in \llbracket
1,n\rrbracket}$. Indeed, the function $u_T$ is differentiable only with respect to the first
component $ y_1$. In this first case, we can apply integration by parts formula to move the derivative on $u_T$, so that
\[I^{od}_1 \, = \, \Bigl{\vert}\int_{\R^{nd}} \bigl[\tilde{p}^{\tau,\xi}(t,T,x,y) -\tilde{p}^{\tau,\xi}
(t,T,x',y)\bigr] D_{ y_1}u_T(y) \,dy\Bigr{\vert}.\]
Noticing that $D_{ y_1}u_T$ is in $C^{\alpha+\beta-1}_{b,d}(\R^{nd})$ thanks to the reverse Taylor expansion (Lemma
\ref{lemma:Reverse_Taylor_Expansion}), the last expression can be then rewritten as
\begin{align}
\label{zz1} I^{od}_1 \, &\le \, \Bigl{\vert}\int_{\R^{nd}} \tilde{p}^{\tau,\xi}(t,T,x,y)\bigl[D_{ y_1}u_T(y)\pm D_{
 y_1}u_T(\tilde{m}^{\tau,\xi}_{T,t}(x))\bigr] \\
 &\qquad\qquad\qquad\qquad\qquad\qquad-\tilde{p}^{\tau,\xi}(t,T,x',y)\bigl[D_{
 y_1}u_T(y)\pm D_{ y_1}u_T(\tilde{m}^{\tau,\xi}_{T,t}(x')) \bigr] \, dy\Bigr{\vert}\notag \\
&\le \, C\Vert u_T \Vert_{C^{\alpha+\beta}_{b,d}}\Bigl{\{}\int_{\R^{nd}}
\bigl[\tilde{p}^{\tau,\xi}(t,T,x,y)\mathbf{d}^{\alpha+\beta-1}(y, \tilde{m}^{\tau,\xi}_{T,t}(x))\bigr] \, dy \notag\\
&\qquad\qquad\qquad\qquad\qquad\qquad+\int_{\R^{nd}}
\bigl[\tilde{p}^{\tau,\xi}(t,T,x',y) \mathbf{d}^{\alpha+\beta-1}(y,\tilde{m}^{\tau,\xi}_{T,t}(x'))\bigr]
\, dy  \notag\\
& \qquad \qquad \qquad \qquad\qquad \qquad \qquad \qquad\qquad\qquad+ \mathbf{d}^{\alpha+\beta-1}(\tilde{m}^{\tau,\xi}_{T,t}(x), \tilde{m}^{\tau,\xi}_{T,t}(x'))\Bigr{\}}. \notag
\end{align}
Now, we use the smoothing effect of $\tilde{p}^{\tau,\xi}$ (Equation \eqref{eq:Smoothing_effects_of_tilde_p}) to control the two
integrals in the last expression, so that
\[I^{od}_1\, \le \, C\Vert u_T \Vert_{C^{\alpha+\beta}_{b,d}}\bigl[(T-t)^{\frac{\alpha+\beta-1}{\alpha}} +
\mathbf{d}^{\alpha+\beta-1}(\tilde{m}^{\tau,\xi}_{T,t}(x),\tilde{m}^{\tau,\xi}_{T,t}(x'))\bigr].\]
We can then conclude the case $j=1$ recalling that the mapping $x\to \tilde{m}^{\tau,\xi}_{T,t}(x)$ is affine (see
Equation \eqref{eq:def_tilde_m_stable} for definition of $\tilde{m}^{\tau,\xi}_{T,t}(x)$) in
order to show that
\begin{equation}\label{Proof:H\"older_for_Frozen_Semigroup_Iod1}
I^{od}_1 \, \le \, C\Vert u_T \Vert_{C^{\alpha+\beta}_{b,d}}\bigl[(T-t)^{\frac{\alpha + \beta -1}{\alpha}} +
\mathbf{d}^{\alpha+\beta-1}(x,x')\bigr].
\end{equation}
Let us consider now the case $j>1$. Using Duality  \eqref{Besov:duality_in_Besov} in Besov spaces and Identification
\eqref{Besov:ident_Holder_Besov}, we can write from Equation \eqref{Proof:H\"older_for_Frozen_Semigroup_Decomp} that
\begin{align}\notag
I^{od}_j \, &\le \, C\Vert u_T \Vert_{C^{\alpha+\beta}_{b,d}}\int_{\R^{(n-1)d}}  \Vert D_{ y_j}\tilde{p}^{\tau,\xi}(t,T,x,
y_{\smallsetminus j},\cdot)-D_{ y_j}\tilde{p}^{\tau,\xi} (t,T,x',y_{\smallsetminus j},\cdot)
\Vert_{B^{-(\alpha_j+\beta_j)}_{1,1}} \, dy_{\smallsetminus j}\\
\notag &\le \, C\Vert u_T \Vert_{C^{\alpha+\beta}_{b,d}}\int_{\R^{(n-1)d}}  \Vert D_{ y_j}\tilde{p}^{\tau,\xi}(t,T,x,
y_{\smallsetminus j},\cdot)\Vert_{B^{-(\alpha_j+\beta_j)}_{1,1}} \\
\notag
&\qquad\qquad\qquad\qquad\qquad\qquad\qquad\qquad+ \Vert D_{ y_j}\tilde{p}^{\tau,\xi}
(t,T,x',y_{\smallsetminus j},\cdot) \Vert_{B^{-(\alpha_j+\beta_j)}_{1,1}} \, dy_{\smallsetminus j} \\
&\le \, C\Vert u_T \Vert_{C^{\alpha+\beta}_{b,d}}(T-t)^{\frac{\alpha+\beta}{\alpha}-\frac{1}{\alpha_j}},\label{Proof:H\"older_for_Frozen_Semigroup_Iodj}
\end{align}
where in the last inequality we applied the first Besov control (Lemma \ref{lemma:First_Besov_COntrols}). Going back at Equations
\eqref{Proof:H\"older_for_Frozen_Semigroup_Iod1}-\eqref{Proof:H\"older_for_Frozen_Semigroup_Iodj}, we finally conclude
that
\begin{align}
\bigl{\vert}D_{ x_1}\tilde{P}^{\tau,\xi}_{T,t}u_T&(x) - D_{ x_1}\tilde{P}^{\tau,\xi}_{T,t}u_T(x')
\bigr{\vert}  \notag\\
&\le \, C\Vert u_T\Vert_{C^{\alpha+\beta}_{b,d}}\bigl[(T-t)^{\frac{\alpha + \beta-1}{\alpha}}+\mathbf{d}^{\alpha+\beta-1}(x,x')
+\sum_{j=2}^n(T-t)^{j-1}(T-t)^{\frac{\alpha+\beta}{\alpha}-\frac{1}{\alpha_j}}\bigr] \notag\\
&\le \, C\Vert u_T\Vert_{C^{\alpha+\beta}_{b,d}} \bigl[(T-t)^{\frac{\alpha + \beta -1}{\alpha}}+\mathbf{d}^{\alpha+\beta-1}(x,x')
\bigr] \label{a2}\\
&\le \, C\Vert u_T\Vert_{C^{\alpha+\beta}_{b,d}}\mathbf{d}^{\alpha+\beta-1}(x,x'),\notag
\end{align}
where in the last passage we used that $T-t \le c_0 \mathbf{d}^\alpha(x,x')$ for some $c_0\le1$.\newline
We focus now on the diagonal regime, i.e.\ when $T-t > c_0\mathbf{d}^\alpha(x,x')$. Remembering that we assumed that
$ x_j= x'_j$ for any $j$ in $\llbracket 2,n\rrbracket$, we start using a Taylor expansion on
the density $\tilde{p}^{\tau,\xi}$ with respect to the first, non-degenerate variable $ x_1$. Namely,
\begin{multline*}
D_{ x_1}\tilde{P}^{\tau,\xi}_{T,t}u_T(x) - D_{ x_1}\tilde{P}^{\tau,\xi}_{T,t}u_T(x') \, = \,
\int_{\R^{nd}}\bigl[D_{ x_1}\tilde{p}^{\tau,\xi}(t,T,x,y)-D_{ x_1}\tilde{p}^{\tau,\xi}
(t,T,x',y)\bigr]u_T(y) \, dy \\
= \, \int_{\R^{nd}}\int_{0}^{1}D^2_{ x_1}\tilde{p}^{\tau,\xi}\bigl(t,T,x'+\lambda(x-x'),y\bigr)
(x-x')_1 u_T(y) \, d\lambda dy.
\end{multline*}
Moreover, from the Scaling Lemma \ref{lemma:link_derivative_density}, it holds that
\[D^2_{ x_1}\tilde{p}^{\tau,\xi}\bigl(t,T,x'+\lambda(x-x'),y\bigr) \, = \,
\sum_{j=1}^{n}C_j(T-t)^{j-1}D_{ y_j}D_{ x_1}\tilde{p}^{\tau,\xi}\bigl(t,T,x'+\lambda(x-x'),y\bigr)\]
and we can use it to write
\begin{align}
\bigl{\vert}D_{x_1}&\tilde{P}^{\tau,\xi}_{T,t}u_T(x) - D_{ x_1}\tilde{P}^{\tau,\xi}_{T,t}u_T(x')
\bigr{\vert} \notag \\
&\le \,C\vert (x-x')_1\vert\sum_{j=1}^n(T-t)^{j-1}\Bigl{\vert}\int_{0}^{1}\int_{\R^{nd}}D_{ y_j}D_{ x_1} \tilde{p}^{\tau,\xi}\bigl(t,T,x'+\lambda(x-x'),y\bigr)u_T(y)\, dyd\lambda\Bigr{\vert} \notag \\
&=:\, C\vert(x-x')_1 \vert\sum_{j=1}^n (T-t)^{j-1} I^d_j. \label{Proof:H\"older_Frozen_Semigroup_Decomp2}
\end{align}
Similarly to the off-diagonal regime, we are going to treat separately the cases $j=1$ and $j>1$ for the \emph{diagonal} contributions
$\{I^{d}_j\}_{j\in \llbracket 1,n\rrbracket}$. In the first case, we can apply
integration by parts formula to show that
\[I^d_1 \, = \, \Bigl{\vert}\int_{0}^{1} \int_{\R^{nd}}
D_{ x_1}\tilde{p}^{\tau,\xi}\bigl(t,T,x'+\lambda(x-x'),y\bigr)\otimes D_{ y_1}u_T(y) \,
dyd\lambda \Bigr{\vert}.\]
A cancellation argument with respect to $D_{ x_1}\tilde{p}^{\tau,\xi}$ then leads to
\[
\begin{split}
I^d_1 \, &= \, \Bigl{\vert}\int_{0}^{1} \int_{\R^{nd}}D_{ x_1}\tilde{p}^{\tau,\xi}(t,T,x'+\lambda(x-x'),y) \\ &\qquad\qquad\qquad\qquad\qquad\qquad\qquad\otimes\bigl[D_{ y_1}u_T(y)-D_{ y_1} u_T(\tilde{m}^{\tau,\xi}_{T,t}
(x'+\lambda(x-x')))\bigr] \, dyd\lambda\Bigr{\vert}\\
&\le \, C\Vert u_T \Vert_{C^{\alpha+\beta}_{b,d}}\int_{0}^{1} \int_{\R^{nd}} \vert D_{ x_1}\tilde{p}^{\tau,\xi}(t,T,x'+\lambda(x-x'),y) \vert\\
&\qquad\qquad\qquad\qquad\qquad\qquad\qquad\times \mathbf{d}^{\alpha+\beta-1}\bigl(y,\tilde{m}^{\tau,\xi}_{T,t}
(x'+\lambda(x-x'))\bigr) \, dyd\lambda.
\end{split}
\]
Since $\alpha + \beta -1 < \alpha$ by hypothesis [\textbf{P}], we can conclude using the smoothing effect
of $\tilde{p}^{\tau,\xi}$ (Lemma \ref{lemma:Smoothing_effect_frozen}) to show that
\begin{equation}\label{Proof:H\"older_Frozen_Semigroup_Id1}
I^d_1 \, \le \, C\Vert u_T \Vert_{C^{\alpha+\beta}_{b,d}}(T-t)^{\frac{\alpha+\beta-2}{\alpha}}.
\end{equation}
For the case $j>1$, we use instead the duality in Besov spaces (Equation \eqref{Besov:duality_in_Besov}) and Identification
\eqref{Besov:ident_Holder_Besov} to write
\begin{equation}\label{Proof:H\"older_Frozen_Semigroup_Idj}
\begin{split}
I^d_j \, &\le \, \int_{0}^{1} \int_{\R^{(n-1)d}} \Vert
D_{ y_j}D_{ x_1}\tilde{p}^{\xi}(t,T,x'+\lambda(x-x'),y_{\smallsetminus j},\cdot)
\Vert_{B^{-(\alpha_j+\beta_j)}_{1,1}} \, dy_{\smallsetminus j}d\lambda \\ 
&\le \, C\Vert u_T \Vert_{C^{\alpha+\beta}_{b,d}}
(T-t)^{\frac{\alpha+\beta}{\alpha}-\frac{1}{\alpha_j}-\frac{1}{\alpha}},
\end{split}
\end{equation}
where in the last passage we applied the first Besov control (Lemma \ref{lemma:First_Besov_COntrols}). From Equations
\eqref{Proof:H\"older_Frozen_Semigroup_Decomp2}, \eqref{Proof:H\"older_Frozen_Semigroup_Id1} and
\eqref{Proof:H\"older_Frozen_Semigroup_Idj}, it is possible to conclude that
\[
\begin{split}
\bigl{\vert}D_{ x_1}\tilde{P}^{\tau,\xi}_{T,t}u_T(x) - D_{ x_1}\tilde{P}^{\tau,\xi}_{T,t}u_T(x') \bigr{\vert}
\, &\le \, C\Vert u_T\Vert_{ C^{\alpha+\beta}_{b,d}}\vert (x-x')_1 \vert\sum_{j=1}^n (T-t)^{j-1}(T-t)^{\frac{\alpha+\beta-1
}{\alpha}- \frac{1}{\alpha_j}}\\
&\le \, C\Vert u_T \Vert_{C^{\alpha+\beta}_{b,d}}\vert (x-x')_1 \vert(T-t)^{\frac{\alpha+\beta-2}{\alpha}} \\
&\le \,
Cc^{\frac{\alpha+\beta-2}{\alpha}}_0\Vert u_T\Vert_{C^{\alpha+\beta}_{b,d}}\mathbf{d}^{\alpha+\beta-1}(x,x'),
\end{split}
\]
where in the last passage we used that $\vert(x-x')_1\vert=\mathbf{d}(x,x')$ and since
$\frac{\alpha+\beta-2}{\alpha}<0$, that
\[\vert(x-x')_1\vert (T-t)^{\frac{\alpha+\beta-2}{\alpha}} \, \le \, c^{\frac{\alpha+\beta-2}{\alpha}}_0 \mathbf{d}^{\alpha +
\beta-1}(x,x').\]
Remembering that $c_0$ is considered fixed and bigger then zero, the searched control follows immediately.

\emph{Controls on frozen Green kernel.} We recall that, in order to preserve the previous terminology of off-diagonal/diagonal regime for the frozen
semigroup, we have introduced the transition time $t_0$, defined in \eqref{eq:def_t0}. Then, while integrating in $s$ from $t$ to $T$, we
will say that the local off-diagonal regime holds for $\tilde{G}^{\tau,\xi}$ if $s$ is in $[t,t_0]$ and that the local diagonal
regime holds if $s$ is in $[t_0,T]$. With the notations of \eqref{eq:def_Green_Kernel_stable} in mind, it seems quite natural now to decompose the derivative of the frozen Green kernel with respect to $t_0$, i.e.
\[D_{ x_1}\tilde{G}^{\tau,\xi}_{T,t}f(t,x) \, = \, D_{ x_1}\tilde{G}^{\tau,\xi}_{t_0,t}f(t,x)
+D_{ x_1}\tilde{G}^{\tau,\xi}_{T,t_0}f(t,x).\]
We remark however that the globally off-diagonal regime is considered in the above decomposition, too. Indeed, when $T-t\le
c_0\mathbf{d}^\alpha(x,x')$, $t_0$ coincides with $T$ and the second term on the right-hand side vanishes.\newline
We start considering the off-diagonal regime represented by 
\[\bigl{\vert} D_{ x_1} \tilde{G}^{\tau,\xi}_{t_0,t} f(t,x) -
D_{ x_1}\tilde{G}^{\tau,\xi}_{t_0,t}f(t,x') \bigr{\vert}.\]
It holds that
\[\bigl{\vert} D_{ x_1}\tilde{G}^{\tau,\xi}_{t_0,t}f(t,x)-D_{ x_1}\tilde{G}^{\tau,\xi}_{t_0,t}
f(t,x')\bigr{\vert} \, \le \, \int_{t}^{t_0}\Bigl[\bigl{\vert} D_{ x_1}\tilde{P}^{\tau,\xi}_{s,t}f(s,x)\bigr{\vert}
+ \bigl{\vert} D_{ x_1} \tilde{P}^{\tau,\xi}_{s,t}f(s,x') \bigr{\vert}\Bigr] \, ds.\]
We then use the control on the frozen semigroup (Equation \eqref{eq:Control_of_semigroup}) to find that
\[
\begin{split}
   \bigl{\vert} D_{ x_1}\tilde{G}^{\tau,\xi}_{t_0,t}f(t,x)-D_{ x_1}\tilde{G}^{\tau,\xi}_{t_0,t}
f(t,x') \bigr{\vert} \, &\le \, C\Vert f \Vert_{L^{\infty}(C^\beta_{b,d})}\int_{t}^{t_0}(s-t)^{\frac{\beta -1}{\alpha}} \, ds  \\
&\le \, C\Vert f \Vert_{L^{\infty}(C^\beta_{b,d})}(t_0-t)^{\frac{\beta + \alpha -1}{\alpha}}. 
\end{split}
\]
Our choice of $t_0$ (cf. Equation \eqref{eq:def_t0}) allows then to conclude that
\[\bigl{\vert}
D_{ x_1}\tilde{G}^{\tau,\xi}_{t_0,t}f(t,x)-D_{ x_1}\tilde{G}^{\tau,\xi}_{t_0,t}f(t,x')
\bigr{\vert} \, \le \, C\Vert f \Vert_{L^{\infty}(C^\beta_{b,d})}\mathbf{d}^{\alpha +\beta -1}(x,x'),\]
remembering that $c_0\le 1$ by assumption. \newline
We can focus now on the diagonal regime represented by 
\[\bigl{\vert} D_{ x_1} \tilde{G}^{\tau,\xi}_{T,t_0}f(t,x) -
D_{ x_1}\tilde{G}^{\tau,\xi}_{T,t_0}f(t,x')\bigr{\vert}.\]
We start applying a Taylor expansion on the derivative of the
semigroup $\tilde{P}^{\tau,\xi}f(t,x)$ so that
\begin{multline*}
\bigl{\vert}
D_{ x_1}\tilde{G}^{\tau,\xi}_{T,t_0}f(t,x)-D_{ x_1}\tilde{G}^{\tau,\xi}_{T,t_0}f(t,x')
\bigr{\vert}\, = \, \Bigl{\vert}\int_{t_0}^{T} \Bigl[ D_{ x_1}\tilde{P}^{\tau,\xi}_{s,t}f(s,x)
-D_{ x_1}\tilde{P}^{\tau, \xi}_{s,t}f(s,x')\Bigr] \,  ds \Bigr{\vert} \\
= \, \Bigl{\vert} \int_{t_0}^{T} \int_{0}^{1}D^2_{ x_1}\tilde{P}^{\tau,\xi}_{s,t}f(s,x+\lambda(x'-x))
(x'-x)_1 \,d\lambda ds \Bigr{\vert}.
\end{multline*}
Then, Fubini theorem and the control on the frozen semigroup (Equation \eqref{eq:Control_of_semigroup}) allow us to write that
\begin{multline*}
\bigl{\vert} D_{ x_1}\tilde{G}^{\tau,\xi}_{T,t_0}f(t,x)-D_{ x_1}\tilde{G}^{\tau,\xi}_{T,t_0}
f(t,x') \bigr{\vert} \, \le \, C\Vert f \Vert_{L^{\infty}(C^\beta_{b,d})}\vert (x-x')_1 \vert
\int_{t_0}^{T}(s-t)^{\frac{\beta -2}{\alpha}} \, ds  \\
\le \, C\Vert f\Vert_{L^{\infty}(C^\beta_{b,d})}\vert
(x-x')_1\vert\bigl[\frac{\alpha}{\alpha+\beta-2}(s-t)^{\frac{\alpha+\beta-2}{\alpha}}\bigr]_{
t_0}^T.
\end{multline*}
Since by hypothesis [\textbf{P}], it holds that $\alpha/(\alpha+\beta -2)<0$, it follows that
\[\bigl{\vert}D_{ x_1}\tilde{G}^{\tau,\xi}_{T,t_0}f(t,x)-D_{ x_1}\tilde{G}^{\tau,\xi}_{T,t_0}
f(t,x')\bigr{\vert} \, \le \, C\Vert f\Vert_{L^{\infty}(C^\beta_{b,d})}\vert (x-x')_1\vert(t_0-t)^{\frac{
\alpha +\beta -2}{\alpha}}.\]
Using that $\vert (x-x')_1\vert =\mathbf{d}(x,x')$ and remembering our choice of $t_0$ in \eqref{eq:def_t0}, we can then conclude
that
\[\bigl{\vert} D_{ x_1}\tilde{G}^{\tau,\xi}_{T,t_0}f(t,x)-D_{ x_1} \tilde{G}^{\tau,\xi}_{T,t_0}
f(t,x') \bigr{\vert} \, \le \, Cc_0^{\frac{\alpha+\beta-2}{\alpha}}\Vert f\Vert_{L^{\infty}(C^\beta_{b,d})}\mathbf{d}^{\alpha+\beta-1}
(x,x').\]

\paragraph{Proof of Lemma \ref{lemma:Holder_modulus_proxy_Deg}}
\emph{Controls on frozen semigroup.} Using the change of variables $z =
\tilde{m}^{\tau,\xi}_{T,t}(x)-y$, we can rewrite $\tilde{P}^{\tau,\xi}_{T,t} u_T(x)$ as
\[\begin{split}
\tilde{P}^{\tau,\xi}_{T,t}u_T(x) \, &= \, \int_{\R^{nd}}\tilde{p}^{\tau,\xi}(t,T,x,y)u_T(y) \,
dy \\
&= \, \int_{\R^{nd}}\frac{1}{\det\bigl(\mathbb{M}_{T-t}\bigr)}p_S(T-t,\mathbb{M}^{-1}_{T-t}\bigl(
\tilde{m}^{\tau,\xi}_{T,t}(x)-y\bigr)u_T(y)\,dy \\
&=\,\int_{\R^{nd}}\frac{1}{\det\bigl(\mathbb{M}_{T-t}\bigr)}p_S(T-t,\mathbb{M}^{-1}_{T-t}z\bigr)u_T(\tilde{m}^{\tau,\xi}_{T,t}
(x)-z) \, dz.
\end{split}
\]
It then follows that
\begin{multline*}
 \bigl{\vert}\tilde{P}^{\tau,\xi}_{T,t}u_T(x) - \tilde{P}^{\tau,\xi}_{T,t}u_T(x')\bigr{\vert} \\
 =\,
\Bigl{\vert}\int_{\R^{nd}}\frac{1}{\det\bigl(\mathbb{M}_{T-t}\bigr)}p_S\bigl(T-t,\mathbb{M}^{-1}_{T-t}z\bigr)
\bigl[u_T(\tilde{m}^{\tau,\xi}_{T,t}(x)-z)-u_T(\tilde{m}^{\tau,\xi}_{T,t}(x')-z)\bigr]
\,dz\Bigr{\vert}.   
\end{multline*}
We now observe that the function $x\to \tilde{m}^{\tau,\xi}_{T,t}(x)$ is affine (cf.\ Equation \eqref{eq:def_tilde_m_stable})
and thus, that
\[\bigl(\tilde{m}^{\tau,\xi}_{T,t}(x)-z\bigr)_1 \, = \,
\bigl(\tilde{m}^{\tau,\xi}_{T,t}(x')-z\bigr)_1,\]
since $ x_1= x'_1$. It then holds that
\[\begin{split}
\bigl{\vert} u_T(\tilde{m}^{\tau,\xi}_{T,t}(x)-z)-u_T(\tilde{m}^{\tau,\xi}_{T,t}(x')-z)\bigr{\vert} \,
&\le \, C\Vert u_T \Vert_{C^{\alpha+\beta}_{b,d}}\mathbf{d}^{\alpha+\beta}\bigl(\tilde{m}^{\tau,\xi}_{T,t}(x),
\tilde{m}^{\tau,\xi}_{T,t}(x')\bigr) \\
&\le\, C\Vert u_T \Vert_{C^{\alpha+\beta}_{b,d}}\mathbf{d}^{\alpha+\beta}(x, x').
\end{split}\]
Hence, we can conclude using it to write
\[
\begin{split}
\bigl{\vert}\tilde{P}^{\tau,\xi}_{T,t}u_T(x) - \tilde{P}^{\tau,\xi}_{T,t}u_T(x')\bigr{\vert} \, &\le \, C\Vert u_T
\Vert_{C^{\alpha+\beta}_{b,d}}\, \mathbf{d}^{\alpha+\beta}(x, x')\int_{\R^{nd}} \frac{p_S\left(T-t,\mathbb{M}^{-1}_{T-t} z\right)}{\det\mathbb{M}_{T-t}}
 \,dz\\
&\le \, C\Vert u_T \Vert_{C^{\alpha+\beta}_{b,d}}\,\mathbf{d}^{\alpha+\beta}(x, x').
\end{split}
\]

\emph{Controls on frozen Green kernel.} We will assume the same notations appeared in the previous lemma for the frozen Green kernel. In
particular, we decompose it as
\[\tilde{G}^{\tau,\xi}_{T,t}f(t,x) \, = \, \tilde{G}^{\tau,\xi}_{t_0,t}f(t,x)+
\tilde{G}^{\tau,\xi}_{T,t_0}f(t,x)\]
with $t_0$ defined in Equation \eqref{eq:def_t0}. \newline
We start rewriting the off-diagonal regime contribution as
\[\begin{split}
\bigl{\vert} \tilde{G}^{\tau,\xi}_{t_0,t}f(t,x)-\tilde{G}^{\tau,\xi}_{t_0,t}f(t,x') \bigr{\vert} \, &= \, \Bigl{\vert}\int_{t}^{t_0}\int_{\R^{nd}}\tilde{p}^{\tau,\xi}(t,s,x,y)\bigl[f(s,y)\pm
f(s,\tilde{m}^{\tau,\xi}_{s,t}(x))\bigr] \\
&\qquad \qquad\,\, -
\tilde{p}^{\tau,\xi}(t,s,x',y) \bigl[f(s,y)\pm f(s,\tilde{m}^{\tau,\xi}_{s,t}(x'))\bigr]\,
dyds\Bigr{\vert} \\
&\le\,\Bigl{\vert}\int_{t}^{t_0}\int_{\R^{nd}}\tilde{p}^{\tau,\xi}(t,s,x,y)\bigl[f(s,y)-f(s,\tilde{m}^{\tau,\xi}_{s,t}(x))\bigr]\\
& \qquad \qquad\,\, 
 -\tilde{p}^{\tau,\xi}(t,s,x',y)\bigl[f(s,y)- f(s,\tilde{m}^{\tau,\xi}_{s,t}(x'))\bigr]\,
dyds\Bigr{\vert} \\
&\,\,\qquad \qquad\qquad \qquad+ \Bigl{\vert}\int_{t}^{t_0} f(s,\tilde{m}^{\tau,\xi}_{s,t}(x))-f(s,\tilde{m}^{\tau,\xi}_{s,t}(x')) \, ds
\Bigr{\vert}.
\end{split}\]
We can then use the smoothing effect for $\tilde{p}^{\tau,\xi}$ (Equation \eqref{eq:Smoothing_effects_of_tilde_p}) to show that
\begin{multline}\label{zz2}
\bigl{\vert} \tilde{G}^{\tau,\xi}_{t_0,t}f(t,x)-\tilde{G}^{\tau,\xi}_{t_0,t}f(t,x') \bigr{\vert} \\
\le \,
C\Vert f \Vert_{L^{\infty}(C^\beta_{b,d})}\int_{t}^{t_0}\bigl[(s-t)^{\beta/\alpha} +\mathbf{d}^{\beta}(\tilde{m}^{\tau,\xi}_{s,t}
(x),\tilde{m}^{\tau,\xi}_{s,t}(x'))\bigr] \, ds.
\end{multline}
Recalling from Equation \eqref{eq:def_tilde_m_stable} that $x\to \tilde{m}^{\tau,\xi}_{s,t}(x)$ is affine, it follows that
\[
\begin{split}
\bigl{\vert} \tilde{G}^{\tau,\xi}_{t_0,t}f(t,x)-\tilde{G}^{\tau,\xi}_{t_0,t}f(t,x') \bigr{\vert} \, &\le \,
C\Vert f\Vert_{L^{\infty}(C^\beta_{b,d})}\int_{t}^{t_0}\bigl[(s-t)^{\beta/\alpha} +\mathbf{d}^{\beta}(x,x')\bigr] \, ds \\
&\le \, C\Vert f \Vert_{L^{\infty}(C^\beta_{b,d})}\bigl[(t_0-t)\mathbf{d}^{\beta}(x,x') + (t_0-t)^{\frac{\beta +\alpha}{\alpha}}\bigr].
\end{split}\]
Using that $t_0-t\le c_0\mathbf{d}^\alpha(x,x')$ for some $c_0\le 1$, we can finally conclude that
\[\bigl{\vert} \tilde{G}^{\tau,\xi}_{t_0,t}f(t,x)-\tilde{G}^{\tau,\xi}_{t_0,t}f(t,x') \bigr{\vert} \, \le \,
C\Vert f\Vert_{L^{\infty}(C^\beta_{b,d})}\mathbf{d}^{\alpha+\beta}(x,x').\]
Now, we can focus our analysis to the diagonal regime contribution: 
\[\bigl{\vert} \tilde{G}^{\tau,\xi}_{T,t_0} f(t,x) -
\tilde{G}^{\tau,\xi}_{T,t_0}f(t,x')\bigr{\vert}.\]
We start applying a Taylor expansion on the frozen semigroup $\tilde{P}^{\tau,\xi}_{s,t}f$ with respect to the $i$-th variable
$x_i$, which is, by hypothesis, the only one for which the entries of $x$ and $x'$ differ.
Namely,
\[
\begin{split}
\bigl{\vert}\tilde{G}^{\tau,\xi}_{T,t_0}f(t,x)-\tilde{G}^{\tau,\xi}_{T,t_0}f(t,x') \bigr{\vert} \, &= \,
\bigl{\vert} \int_{t_0}^{T} \tilde{P}^{\tau,\xi}_{s,t}f(s,x) - \tilde{P}^{\tau,\xi}_{s,t}f(s,x') \,  ds
\bigr{\vert} \\
&= \, \bigl{\vert} \int_{t_0}^{T} \int_{0}^{1}D_{x_i}\tilde{P}^{\tau,\xi}_{s,t}f(s,x+\lambda(x'-x))
\cdot(x'-x)_i \, d\lambda ds \bigr{\vert}.
\end{split}
\]
The control on the frozen semigroup (Equation \eqref{eq:Control_of_semigroup}) then implies that
\begin{equation}\label{a3}
\bigl{\vert}\tilde{G}^{\tau,\xi}_{T,t_0}f(t,x)-\tilde{G}^{\tau,\xi}_{T,t_0}f(t,x') \bigr{\vert} \, \le \,
C\Vert f\Vert_{L^{\infty}(C^\beta_{b,d})}\vert (x-x')_i\vert\int_{t_0}^{T}(s-t)^{\frac{\beta}{\alpha}-\frac{1}{\alpha_i}}
\, ds.
\end{equation}
Noticing from assumption [\textbf{P}] that $\beta+\alpha-1-\alpha(i-1) <0$ for $i\ge 2$, it holds that
\[
\begin{split}
\int_{t_0}^{T}(s-t)^{\frac{\beta}{\alpha}-\frac{1}{\alpha_i}}\, ds \, &= \, \int_{t_0}^{T}(s-t)^{\frac{\beta-[1+\alpha(i-1)]}{\alpha}}\, ds
\\
&\le \, C\Bigl[-(s-t)^{\frac{\beta+\alpha-1-\alpha(i-1)}{\alpha}}\Bigr]_{t_0}^T\\
&\le \,C(t_0-t)^{\frac{\beta-1-\alpha(i-2)}{\alpha}}.
\end{split}\]
Using that $\vert (x-x')_i\vert=\mathbf{d}^{1+\alpha(i-1)}(x,x')$ and our choice of $t_0$ (cf. Equation \eqref{eq:def_t0}), we
can then conclude from Equation \eqref{a3} that
\begin{multline*}
\bigl{\vert} \tilde{G}^{\tau,\xi}_{T,t_0}f(t,x)-\tilde{G}^{\tau,\xi}_{T,t_0}f(t,x') \bigr{\vert} \\
\le  \,
Cc^{\frac{\beta-1-\alpha(i-2)}{\alpha}}_0\Vert f \Vert_{L^{\infty}(C^\beta_{b,d})}\mathbf{d}^{\alpha+\beta}(x,x') \, \le \,
Cc^{\frac{\beta-\gamma_i}{\alpha}}_0\Vert f \Vert_{L^{\infty}(C^\beta_{b,d})}\mathbf{d}^{\alpha+\beta}(x,x'),
\end{multline*}
remembering the definition of $\gamma_i$  given in \eqref{Drift_assumptions}.

\section{A priori estimates}
\fancyhead[RO]{Section \thesection. A priori estimates}
Since the aim of this section is to prove Proposition \ref{prop:A_Priori_Estimates}, we will assume tacitly from this point further
that assumption [\textbf{A'}] holds. Moreover, we recall here that throughout this section, we are considering the
regularized framework of Section $3.2$.

\textbf{WARNING:} For notational simplicity, we drop here the subscripts and the superscripts in $m$ associated with the
regularization. For any fixed $(\tau,\xi)$ in $[0,T]\times\R^{nd}$, we rewrite, with some abuse in notations, the Duhamel Expansion \eqref{eq:Expansion_along_proxy} as:
\begin{equation}\label{align:Representation2}
u(t,x) \, = \, \tilde{u}^{\tau,\xi}(t,x) + \int_{t}^{T}\tilde{P}^{\tau,\xi}_{s,t}R^{\tau,\xi}(s,x)\, ds,
\end{equation}
where $\tilde{u}^{\tau,\xi}$ is defined through the Duhamel Representation \eqref{Duhamel_representation_of_proxy} and
\[R^{\tau,\xi}(t,x) \, = \, \bigl{\langle}F(t,x)-F(t,\theta_{t,\tau}(\xi)), D_{x}u(t,x)
\bigr{\rangle}, \quad (t,x) \, \in \, (0,T)\times\R^{nd}.\]
It is however important to keep in mind that $f$, $u_T$, $F$ are now smooth and bounded functions so that all the terms above are clearly
defined. We recall however that we aim at obtaining controls in the $L^\infty(C^{\alpha+\beta}_{b,d})$-norm, uniformly with respect to the regularization parameter.\newline
From the expansion above, we know moreover that for any $(t,\xi)$ in $[0,T]\times \R^{nd}$, it holds  that
\begin{equation}\label{align:Representation2_deriv}
D_{ x_1}u(t,x) \, =\,
D_{ x_1}\tilde{u}^{\tau,\xi}(t,x)+\int_{t}^{T}D_{ x_1}\tilde{P}^{\tau,\xi}_{s,t}R^{\tau,\xi}(s,x)\, ds.
\end{equation}
As seen in the previous section, these decompositions will allow us to obtain a control for $u$ in $L^\infty\bigl(0,T;C^{\alpha+\beta}_{b,d}(\R^{nd})\bigr)$
analyzing separately the contributions from the Duhamel representation $\tilde{u}^{\tau,\xi}$ and those from the expansion error $R^{\tau,\xi}(t,x)$, for suitable choices of freezing parameters $(\tau,\xi)$.

\subsection{Second Besov control}
\label{Sec:stabile:Second_Besov}
This sub-section focuses on the contribution associated with the remainder term $R^{m,\tau,\xi}$ appearing in the
Duhamel-type Expansion \eqref{align:Representation2}. We recall that we aim at controlling it with the
$L^\infty(C^{\alpha+\beta}_{b,d})$-norm of the coefficients, uniformly in the regularization parameter.  Let us start decomposing it through
\[\Bigl{\vert}\int_{t}^{T}  \tilde{P}^{\tau,\xi}_{s,t}R^{\tau,\xi}(s,x)\, ds\Bigr{\vert}\, = \, \Bigl{\vert}
\sum_{j=1}^{n}\int_{t}^{T}\int_{\R^{nd}}\tilde{p}^{\tau,\xi}(t,s,x,y)\Delta^{\tau,\xi} F_j(s,y)\cdot D_{ y_j}u(s,y) \, dyds\Bigl{\vert},\]
where we have denoted for simplicity
\begin{equation}
\label{eq:scorciatoia}
\Delta^{\tau,\xi} F_j(s,y) \, := \,F_j(s,y)- F_j(s,\theta_{s,\tau}
(\xi)), \quad j \in \llbracket 1,n\rrbracket.
\end{equation}
We then notice that the non-degenerate contribution in the sum (corresponding to the index $j=1$) can be treated easily, remembering that
$u$ is differentiable with respect to the first component with a bounded derivative.
Indeed, using the smoothing effect for the frozen density $\tilde{p}^{\tau,\xi}$ (Equation \eqref{eq:Smoothing_effects_of_tilde_p}), it
holds that
\[
\begin{split}
&\Bigl{\vert}\int_{t}^{T}\int_{\R^{nd}}\tilde{p}^{\tau,\xi}(t,s,x,y) \Delta^{\tau,\xi} F_1(s,y)\cdot D_{ y_1}u(s,y) \, dyds\Bigl{\vert} \\
&\qquad\qquad\le \, C\Vert D_{ y_1}u(s,y)\Vert_{l^\infty(L^\infty)}\Vert F\Vert_H \int_{t}^{T}\int_{\R^{nd}}\tilde{p}^{\tau,\xi}
(t,s,x,y)\mathbf{d}^{\alpha+\beta}\bigl(y,\theta_{s,\tau}(\xi)\bigr) \, dyds \\
&\qquad\qquad\le \, C\Vert D_{ y_1}u(s,y)\Vert_{l^\infty(L^\infty)}\Vert F\Vert_H \int_{t}^{T} (s-t)^{\frac{\beta}{\alpha}} \, ds \\ &\qquad\qquad\le \,
C\Vert D_{ y_1}u(s,y)\Vert_{l^\infty(L^\infty)}\Vert F\Vert_H(T-t)^{\frac{\alpha+\beta}{\alpha}}.
\end{split}\]
In order to deal with the degenerate indexes, we will use, similarly to the previous subsection, a reasoning in Besov spaces. Since $u$ is
not differentiable with respect to $ y_j$ if $j>1$, we move the
derivative to the other terms using integration by parts formula:
\[\Bigl{\vert}\int_{t}^{T}\int_{\R^{nd}}D_{ y_j}\cdot\Bigl{\{}\tilde{p}^{\tau,\xi}(t,s,x,y)\Delta^{\tau,\xi} F_j(s,y)\Bigr{\}} u(s,y) \, dyds\Bigl{\vert}.\]
In order to rely again on the duality in Besov spaces (Equation \eqref{Besov:duality_in_Besov}), we rewrite the above expression as
\[
\begin{split}
  &\Bigl{\vert}\int_{t}^{T}\int_{\R^{nd}}D_{ y_j}\cdot\Bigl{\{}\tilde{p}^{\tau,\xi}(t,s,x,y)\Delta^{\tau,\xi} F_j(s,y)\Bigr{\}} u(s,y) \, dyds\Bigl{\vert} \\
&\qquad \qquad\le \, \int_{t}^{T}\int_{\R^{(n-1)d}}\Bigl{\Vert}D_{ y_j}\cdot\Bigl{\{}\tilde{p}^{\tau,\xi}(t,s,x,y_{\smallsetminus j},
\cdot)\Delta^{\tau,\xi} F_j(s,y_{\smallsetminus j},\cdot)\Bigr{\}} \Bigr{\Vert
}_{B^{-(\alpha_j+\beta_j)}_{1,1}}\\
&\qquad \qquad \qquad\qquad \qquad\qquad \qquad \qquad\qquad \qquad\qquad \qquad\times \Vert u(s,y_{\smallsetminus j},\cdot) \Vert_{B^{\alpha_j+\beta_j}_{\infty,\infty}} \,
dy_{\smallsetminus j}ds.
\end{split}\]
Remembering the identification in Equation \eqref{Besov:ident_Holder_Besov}, it holds now that
\begin{multline*}
\Bigl{\vert}\int_{t}^{T}\int_{\R^{nd}}D_{ y_j}\cdot\Bigl{\{}\tilde{p}^{\tau,\xi}(t,s,x,y)\Delta^{\tau,\xi} F_j(s,y)\Bigr{\}} u(s,y) \, dyds\Bigl{\vert} \\
\le \, \Vert u\Vert_{L^\infty(C^{\alpha+\beta}_{b,d})} \int_{t}^{T}\int_{\R^{(n-1)d}}\Bigl{\Vert}D_{ y_j} \cdot\Bigl{\{}
\tilde{p}^{\tau,\xi}(t,s,x, y_{\smallsetminus j}, \cdot)\Delta^{\tau,\xi} F_j(s,y_{\smallsetminus
j},\cdot)\Bigr{\}} \Bigr{\Vert}_{B^{-(\alpha_j+\beta_j)}_{1,1}} \,
dy_{\smallsetminus j}ds.
\end{multline*}
It then remains to control the integral of the Besov norm above. To do that, we will need a refinement of the smoothing effect (Equation
\eqref{eq:Smoothing_effects_of_tilde_p}) that involves only partial differences of variables. For a fixed $i$ in $\llbracket2,n\rrbracket$,
we start denoting  by $d_{i:n}(\cdot,\cdot)$ the part of the anisotropic distance considering only the last $n-(i-1)$ variables. Namely,
\[d_{i:n}( x,x') \, := \, \sum_{j=i}^{n} \vert(x-x')_j\vert^{\frac{1}{1+\alpha(j-1)}}.\]

\begin{lemma}[Partial Smoothing Effect]
\label{lemma:Partial_Smoothing_Effect}
Let $i$ be in $\llbracket2,n\rrbracket$, $\gamma$ in $(0,1\wedge \alpha(1+\alpha(i-1)))$ and $\vartheta$, $\varrho$ two $n$-multi-indexes
such that $\vert \vartheta +\varrho \vert\le 3$. Then, there exists a constant $C:=C(\vartheta,\varrho,\gamma)$ such that for any $t<s$ in
$[0,T]$, any $x$ in $\R^{nd}$,
\begin{equation}\label{eq:Partial_Smoothing_Effect}
\int_{\R^{nd}} \vert D^\varrho_{y}D^\vartheta_{x} \tilde{p}^{ \tau,\xi} (t,s,x,y) \vert \mathbf{d}^\gamma_{i:n}
\bigl(y,\theta_{s,\tau}(\xi)\bigr) \, dy \, \le \, C (s-t)^{\frac{\gamma}{\alpha}-\sum_{i=k}^{n}
\frac{\vartheta_k+\varrho_k}{\alpha_k}},
\end{equation}
taking $(\tau,\xi)=(t,x)$.
\end{lemma}

The above assumption on $\gamma$ should not appear too strange. Indeed, in the partial distance $\mathbf{d}^{\gamma}_{i:n}( x,x')$, the
stronger term to be integrated  is at level $i$ with intensity of order $\gamma/(1+\alpha(i-1))$. Since by the smoothing effect (Equation
\eqref{eq:Smoothing_effects_of_tilde_p}) of the frozen density we know we can integrate against $\tilde{p}^{\tau,\xi}$ contributions
of order up to $\alpha$, the condition $\gamma<\alpha(1+\alpha(i-1))$ appears naturally.\newline
A proof of this result can be obtained mimicking with slightly modifications the proof in Lemma \ref{lemma:Smoothing_effect_frozen}.

As done above for the first Besov control (Lemma \ref{lemma:First_Besov_COntrols}), we will however state the result considering a possibly
additional derivative with respect to
$ x_1$. Namely, we would like to control the following:
\[D_{ y_j}\cdot\Bigl{\{}\mathbf{d}^{\vartheta}_{x}\tilde{p}^{\tau,\xi}(t,s,x,
y_{\smallsetminus j},\cdot) \otimes \bigl[ F_j (s,y_{\smallsetminus j},\cdot)- F_j(s,\theta_{s,\tau}(\xi))
\bigr]\Bigr{\}}\]
where we have denoted as in \eqref{eq:notation_smallsetminus}, $ F_j (s,y_{\smallsetminus j},\cdot) :=  F_j
(s,y_{1},\dots, y_{j-1},\cdot, y_{j+1},\dots, y_n)$ and, with a slightly abuse of notation, by $D_{ y_j}\cdot$ an
extended form of the divergence over the $j$-th variable. In other words, this ``enhanced'' divergence form decreases by one the order of the
input tensor.

\begin{lemma}[Second Besov Control]
\label{lemma:Second_Besov_COntrols}
Let $j$ be in $\llbracket 2,n\rrbracket$ and $\vartheta$ a multi-index in $\N^n$ such that $\vert \vartheta \vert\le 2$. Under
(\textbf{A}'), there exists a constant $C:=C(j,\vartheta)$ such that for any $x$ in $\R^{nd}$ and any $t<s$ in $[0,T]$
\begin{multline*}
\int_{\R^{(n-1)d}}\Bigl{\Vert}D_{ y_j}\cdot\Bigl{\{}\mathbf{d}^{\vartheta}_{x}\tilde{p}^{\tau,\xi}(t,s,x,
y_{\smallsetminus j},\cdot) \otimes \bigl[ F_j (s,y_{\smallsetminus j},\cdot)- F_j(s,\theta_{s,\tau}(\xi))
\bigr]\Bigr{\}} \Bigr{\Vert}_{B^{-(\alpha_j+\beta_j)}_{1,1}} \, dy_{\smallsetminus j}\\
\le \, C\Vert F
\Vert_H(s-t)^{\frac{\beta}{\alpha} -\sum_{k=1}^{n}\frac{\vartheta_k}{\alpha_k}},
\end{multline*}
taking $(\tau,\xi)=(t,x)$.
\end{lemma}
\begin{proof}
To control the Besov norm in $B^{-(\alpha_j+\beta_j)}_{1,1}(\R^d)$, we are going to use the Thermic Characterization
\eqref{alpha-thermic_Characterization} with $\tilde{\gamma}=-(\alpha_j+\beta_j)$. Since the first term can be controlled as in the first
Besov control (Lemma \ref{lemma:First_Besov_COntrols}), we will
focus on the second one, i.e.
\[\int_{0}^{1}v^{\frac{\alpha_j+\beta_j}{\alpha}}\int_{\R^d}\Bigl{\vert}\int_{\R^d}\partial_vp_h(v,z- y_j)D_{ y_j}\cdot
\Bigl{\{}D^\vartheta_{x}\tilde{p}^{\tau,\xi}(t,s,x,y)\otimes\Delta^{\tau,\xi} F_j(s,y)\Bigr{\}}
\,d y_j \Bigr{\vert} \, dzdv,\]
where we exploited the same notations for $\Delta^{\tau,\xi} F_j$ given in \eqref{eq:scorciatoia}.\newline 
We start applying integration by parts formula noticing that 
\[D_{ y_j}p_h(v,z- y_j)\, = \, -D_z p_h(v,z- y_j),\]
to write that
\[\int_{0}^{1}v^{\frac{\alpha_j+\beta_j}{\alpha}}\int_{\R^d}\Bigl{\vert}\int_{\R^d}D_z\partial_vp_h(v,z- y_j)\cdot \Bigl{\{}D^\vartheta_{
x} \tilde{p}^{\tau,\xi}(t,s,x,y)\otimes\Delta^{\tau,\xi} F_j(s,y)\Bigr{\}}
 \,d y_j \Bigr{\vert} \,dzdv.\]
Fixed a constant $\delta_j\ge1$ to be chosen later, we then split the above integral with respect to $v$ into two components:
\begin{multline*}
\int_{0}^{(s-t)^{\delta_j}}\!\!v^{\frac{\alpha_j+\beta_j}{\alpha}}\int_{\R^d}\Bigl{\vert}\int_{\R^d}D_{z}\partial_vp_h
(v,z- y_j)\cdot\Bigl{\{}D^\vartheta_{x}\tilde{p}^{\tau,\xi}(t,s,x,y)\Delta^{\tau,\xi} F_j(s,y) \Bigr{\}}\,d y_j  \Bigr{\vert} \, dzdv \\
+\int_{(s-t)^{\delta_j}}^{1}v^{\frac{\alpha_j+\beta_j}{\alpha}}\int_{\R^d}\Bigl{\vert}\int_{\R^d}D_z\partial_vp_h(v,z- y_j)\cdot\Bigl{\{}
D^\vartheta_{x} \tilde{p}^{\tau,\xi}(t,s,x,y)\Delta^{\tau,\xi} F_j(s,y)
\Bigl{\}} \,d y_j \Bigr{\vert} \, dzdv \\
=: \, \bigl(I_1+I_2\bigr)(y_{\smallsetminus j}).
\end{multline*}
The second component $I_2$ has no time-singularity and it can be easily controlled using Fubini theorem in the following way
\begin{multline*}
I_2(y_{\smallsetminus j}) \, \le \,
C\Vert F\Vert_{H}\int_{(s-t)^{\delta_j}}^1 v^{\frac{\alpha_j+\beta_j}{\alpha}}\int_{\R^d}\Bigl(\int_{\R^d}\vert
D_z\partial_vp_h(v,z- y_j) \vert \, dz\Bigr) \vert D^\vartheta_{x} \tilde{p}^{\tau,\xi}(t,s,x,y)\vert\\
\times \mathbf{d}^{1+\alpha(j-2)+\beta}_{j:n}(y, \theta_{s,\tau}(\xi)) \,d y_jdv,
\end{multline*}
remembering that $ F_j(t,\cdot)$ depends only on the last $(n-j)$ variables and it is in $C^{1+\alpha(j-2)+\beta}_{b,d}(\R^{nd})$ by
assumption [\textbf{R}]. We can then use the smoothing effect of the heat-kernel $p_h$ (Equation \eqref{Smoothing_effect_of_Heat_Kern}) and
Fubini theorem again, in order to write that
\begin{multline*}
I_2(y_{\smallsetminus j}) \, \le \, C\Vert F\Vert_{H}\int_{(s-t)^{\delta_j}}^1
\frac{v^{\frac{\alpha_j+\beta_j-1}{\alpha}}}{v}\int_{\R^d} \vert D^\vartheta_{x} \tilde{p}^{\tau,\xi}(t,s,x,y)\vert
\mathbf{d}^{1+\alpha(j-2)+\beta}_{j:n}(y, \theta_{s,\tau}(\xi)) \,d y_jdv \\
\le \, C\Vert F\Vert_{H}\Bigl(\int_{(s-t)^{\delta_j}}^1\frac{v^{\frac{\alpha_j+\beta_j-1}{\alpha}}}{v} \, dv\Bigr)\Bigl(\int_{\R^d}
\vert D^\vartheta_{x} \tilde{p}^{\tau,\xi}(t,s,x,y)\vert \mathbf{d}^{1+\alpha(j-2)+\beta}_{j:n}(y,
\theta_{s,\tau}(\xi)) \,d y_j\Bigr).
\end{multline*}
Noticing from \eqref{eq:def_alpha_i_and_beta_i} that $\alpha_j+\beta_j-1<0$, it holds now that
\[I_2(y_{\smallsetminus j})\,\le\,C \Vert F\Vert_{H}(s-t)^{\delta_j\frac{\alpha_j+\beta_j-1}{\alpha}}\int_{\R^d} \vert
D^\vartheta_{x} \tilde{p}^{\tau,\xi}(t,s,x,y)\vert \mathbf{d}^{1+\alpha(j-2)+\beta}_{j:n}(y, \theta_{s,\tau}(\xi))
\,d y_j.\]
We can finally add the integral with respect to the other components $y_{\smallsetminus j}$. In order to use now the partial smoothing
effect (Equation \eqref{eq:Partial_Smoothing_Effect}), we take $\tau=t$ and $\xi=x$ and notice that
by assumption [\textbf{P}],
\begin{equation}\label{a1}
1+\alpha(j-2)+\beta \, < \, 1+\alpha(j-1)-(1-\alpha)\bigl(1+\alpha(j-1)\bigr) \, =\,
\alpha(1+\alpha(j-1)).
\end{equation}
It then holds that
\begin{equation}\label{Proof:Second_Besov_Control_I_2}
\begin{split}
\int_{\R^{(n-1)d}}&I_2(y_{\smallsetminus j}) \, dy_{\smallsetminus j}\\
&\le\, C \Vert F
\Vert_{H}(s-t)^{\delta_j\frac{\alpha_j+\beta_j-1}{\alpha}} \int_{\R^{nd}} \vert
D^\vartheta_{x} \tilde{p}^{\tau,\xi}(t,s,x,y)\vert \mathbf{d}^{1+\alpha(j-2)+\beta}_{j:n}(y, \theta_{s,\tau}(\xi))
\,dy \\
&\le\, C\Vert F \Vert_H(s-t)^{\delta_j \frac{\alpha_j+\beta_j-1}{\alpha}+\frac{1+\alpha(j-2)+\beta}{\alpha}-
\sum_{k=1}^{n}\frac{\vartheta_k}{\alpha_k}}.
\end{split}
\end{equation}
To control the other term $I_1$, we focus at first on the inner integral with respect to $ y_j$:
\[\int_{\R^d}D_z \partial_vp_h (v,z- y_j)\cdot \Bigl{\{}D^\vartheta_{x}\tilde{p}^{\tau,\xi}(t,s,x,y)\otimes\Delta^{\tau,\xi} F_j(s,y)\Bigr{\}}\,d y_j.\]
We start using a cancellation argument with respect to the density $p_h$ to write that
\begin{multline*}
\int_{\R^d}D_z \partial_vp_h (v,z- y_j)\cdot\Bigl{\{}D^\vartheta_{x}\tilde{p}^{\tau,\xi}(t,s,x,y)\otimes
\bigl[ F_j(s,y) - F_j(s,\theta_{s,\tau}(\xi))\Bigr{\}} \\
- D^\vartheta_{x}\tilde{p}^{\tau,\xi}(t,s,x,y_{\smallsetminus j},z)\otimes\bigl[ F_j(s,y_{\smallsetminus
j},z)- F_j(s,\theta_{s,\tau}(\xi))\bigr]\Bigr{\}}\,d y_j 
\end{multline*}
We can then divide the above integral into two components $J_1+J_2$ given by  in
\[\begin{split}
J_1(v,y_{\smallsetminus j},z) \, &:= \, \int_{\R^d}D_z \partial_vp_h(v,z- y_j)\cdot \Bigl{\{}D^\vartheta_{x}\tilde{p}^{\tau,\xi}(t,s,x,y)\\
&\qquad\qquad\qquad\qquad\qquad\qquad\qquad\qquad\otimes\bigl[
 F_j(s,y)- F_j(s,y_{\smallsetminus j},z)\bigr]\Bigr{\}} \,d y_j; \\
J_2(v,y_{\smallsetminus j},z) \, &:= \, \int_{\R^d}D_z \partial_vp_h(v,z- y_j)\cdot\Bigl{\{}\bigl[D^\vartheta_{x}\tilde{p}^{\tau,\xi}(t,s,x,y)-
D^\vartheta_{x} \tilde{p}^{\tau,\xi}(t,s,x,y_{\smallsetminus j},z)\bigr]\\
&\qquad\qquad\qquad\qquad\qquad\qquad\qquad\qquad\otimes\bigl[ F_j(s,y_{\smallsetminus
j},z)- F_j(s,\theta_{s,\tau}(\xi)) \bigr]\Bigl{\}} \,d y_j.
\end{split}
\]
Remembering the notation for $ F_j(s,y_{\smallsetminus j},z)$ in \eqref{eq:notation_smallsetminus} and that $ F_j$ is
$\frac{1+\alpha(j-2)+\beta}{1+\alpha(j-1)}$-H\"older continuous with respect to its $j$-th variable by assumption [\textbf{R}], the first
component $J_1$ can be easily controlled using Fubini theorem by
\[\begin{split}
\int_{\R^d}\Bigl{\vert}J_1(&v,y_{\smallsetminus j},z)\Bigr{\vert}\,dz\\
&\le \,C\Vert F \Vert_H \int_{\R^d}\Bigl(\int_{\R^d}\vert z- y_j\vert^{\frac{1+\alpha(j-2)+\beta}{1+\alpha(j-1)}}\vert D_{z}
\partial_vp_h(v,z- y_j)\vert \, dz\Bigr)\vert D^\vartheta_{x}\tilde{p}^{\tau,\xi}(t,s,x,y) \vert \,d y_j\\
&\le \,C\Vert F \Vert_H v^{\frac{1}{\alpha}\frac{1+\alpha(j-2)+\beta}{1+\alpha(j-1)}-\frac{1}{\alpha}-1} \int_{\R^d}
\vert D^\vartheta_{x}\tilde{p}^{\tau,\xi}(t,s,x,y) \vert \,d y_j,
\end{split}\]
where in the last passage we used the smoothing effect of the heat-kernel $p_h$ (Equation \eqref{Smoothing_effect_of_Heat_Kern}), noticing
that
\[\frac{1+\alpha(j-2)+\beta}{1+\alpha(j-1)} \, = \, 1+\frac{\beta-\alpha}{1+\alpha(j-1)}\, < \, 1+\alpha,\]
since $\alpha>\beta$ by assumption [\textbf{P}].
Using now the identity
\begin{equation}\label{Proof:Second_Besov_Control}
\frac{\alpha_j+\beta_j}{\alpha}+\frac{1}{\alpha}\Bigl(\frac{1+\alpha(j-2)+\beta}{1+\alpha(j-1)}-1\Bigr) \, = \, \frac{2\beta_j}{\alpha},
\end{equation}
we add the integral with respect to $v$ and write that
\[
\begin{split}
\int_{0}^{(s-t)^{\delta_j}}v^{\frac{\alpha_j+\beta_j}{\alpha}}\int_{\R^d}\Bigl{\vert}J_1(v,y_{\smallsetminus j},z)\Bigr{\vert}\,dzdv\,
&\le \,C\Vert F \Vert_H \int_{0}^{(s-t)^{\delta_j}}\! \frac{v^{\frac{2\beta_j}{\alpha}}}{v}\int_{\R^d}
\vert D^\vartheta_{x}\tilde{p}^{\tau,\xi}(t,s,x,y) \vert \,d y_j dv \\
&\le \, C\Vert F \Vert_H (s-t)^{\delta_j\frac{2\beta_j}{\alpha}} \int_{\R^d} \vert
D^\vartheta_{x}\tilde{p}^{\tau,\xi}(t,s,x,y) \vert \,d y_j.
\end{split}\]
Adding the integral with respect to the other components $y_{\smallsetminus j}$, we can finally conclude that
\begin{equation}\label{Proof:Second_Besov_Control_J1}
\begin{split}
\int_{\R^{(n-1)d}}\int_{0}^{(s-t)^{\delta_j}}v^{\frac{\alpha_j+\beta_j}{\alpha}}\int_{\R^d}&\Bigl{\vert}J_1(v,y_{\smallsetminus j},z)
\Bigr{\vert} \, dzdv \,dy_{\smallsetminus j}\\
&\le \,  C\Vert F\Vert_H(s-t)^{\delta_j\frac{2\beta_j}{\alpha}} \int_{\R^{nd}}
\vert D^\vartheta_{x}\tilde{p}^{\tau,\xi}(t,s,x,y) \vert\,dy\\
&\le \, C\Vert F \Vert_H (s-t)^{\delta_j\frac{2\beta_j}{\alpha}-\sum_{k=1}^{n}\frac{\vartheta_k}{\alpha_k}}.
\end{split}
\end{equation}
To control the second component $J_2$, we start applying a Taylor expansion on $\tilde{p}^{\tau,\xi}$ with respect to $ y_j$:
\begin{multline}\label{eq:J_2_component_for_extension}
J_2(v,y_{\smallsetminus j},z) \,= \, \int_{\R^d}D_z \partial_vp_h(v,z- y_j)\cdot\Bigl{\{}\Delta^{\tau,\xi} F_j(s,y_{\smallsetminus
j},z)\\
\otimes\int_{0}^{1}D_{ y_j}D^\vartheta_{x} \tilde{p}^{\tau,\xi}(t,s,x,y_{\smallsetminus j}, y_j+
\lambda(z- y_j))\cdot(z)\Bigr{\}} \,d\lambda d y_j.
\end{multline}
We then notice that for any fixed $\lambda$ in $[0,1]$, it holds that
\[
\begin{split}
|\Delta^{\tau,\xi} F_j(s,&y_{\smallsetminus j}, z)| \, = \, \vert  F_j(s,y_{\smallsetminus j}, z)- F_j(s,\theta_{s,\tau}(\xi)) \vert \\
&\le \, C\Vert F
\Vert_H\Bigl{\{}\bigl{\vert}\bigl(z-\theta_{s,\tau}(\xi)\bigr)_j\bigr{\vert}^{\frac{1+\alpha(j-2)+\beta}{1
+\alpha(j-1)}}+\sum_{k=j+1}^{n}\bigl{\vert}\bigl(y-\theta_{s,\tau}(\xi)\bigr)_k
\bigr{\vert}^{\frac{1+\alpha(j-2)+\beta}{1+\alpha(k-1)}}\Bigr{\}} \\
&\le \, C\Vert F \Vert_H\Bigl{\{}\bigl{\vert} \lambda( y_j-z) \bigr{\vert}^{\frac{1+\alpha(j-2)+\beta}{1 +\alpha(j-1)}}+\bigl{\vert}
\bigl( y_j+\lambda(z- y_j)-\theta_{s,\tau}(\xi)\bigr)_j\bigr{\vert}^{\frac{1+ \alpha(j-2) +\beta}{1
+\alpha(j-1)}}\\
&\quad\qquad\qquad\qquad\qquad\qquad\qquad\qquad\qquad\qquad+\sum_{k=j+1}^{n}\bigl{\vert}\bigl(y-\theta_{s,\tau}(\xi)\bigr)_k\bigr{\vert}^{\frac{1+\alpha(j-2)+\beta}{1+
\alpha(k-1)}}\Bigr{\}}\\
&\le \, C\Vert F \Vert_H\Bigl{\{}\bigl{\vert} z- y_j \bigr{\vert}^{\frac{1+\alpha(j-2)+\beta}{1 +\alpha(j-1)}} +
\mathbf{d}^{1+\alpha(j-2)+\beta}_{j:n}\bigl((y_{\smallsetminus j}, y_j+\lambda(z- y_j))\bigr),\theta_{s,\tau}(\xi)
\bigr)\Bigr{\}},
\end{split}\]
where as in \eqref{eq:notation_smallsetminus}, we denoted \[(y_{\smallsetminus j}, y_j+\lambda(z- y_j)):=( y_1,\dots,
 y_{j-1}, y_1,\dots, y_{j-1}, y_j+\lambda(z- y_j), y_{j+1}, \dots, y_n).\]
We can thus split $J_2$ as
\begin{equation}\label{Proof:Second_Besov_Control_J'0}
  \vert J_2(v,y_{\smallsetminus j},z) \vert \, \le \,
C\Vert F \Vert_H\int_{0}^{1}\bigl(J_{2,1}+J_{2,2}\bigr)(v,y_{\smallsetminus j},z,\lambda)\, d\lambda,
\end{equation}
where we denoted for simplicity:
\[
\begin{split}
J_{2,1}(v,y_{\smallsetminus j},z,\lambda) \, &:= \, \int_{\R^d}\vert z- y_j \vert^{\frac{1+\alpha(j-2)+\beta}{1
+\alpha(j-1)}+1}\vert D_z \partial_vp_h(v,z- y_j)\vert\\
&\qquad\qquad\qquad\quad\times\vert D_{ y_j}D^\vartheta_{x}
\tilde{p}^{\tau,\xi}(t,s,x,y_{\smallsetminus j}, y_j+\lambda(z- y_j))\vert \,d y_j \\
J_{2,2}(v,y_{\smallsetminus j},z,\lambda) \, &:= \, \int_{\R^d}\vert z- y_j \vert\, \vert D_z \partial_vp_h(v,z- y_j)\vert\\
&\qquad\qquad\qquad\quad\times\vert D_{ y_j}D^\vartheta_{x}
\tilde{p}^{\tau,\xi}(t,s,x,y_{\smallsetminus j}, y_j+\lambda(z- y_j))\vert\\ 
&\qquad\qquad\qquad\qquad\quad\times \mathbf{d}^{1+\alpha(j-2)+\beta}_{j:n} ((y_{\smallsetminus j}, y_j+\lambda(z- y_j)),\theta_{s,\tau}(\xi))
\, d y_j
\end{split}\]
Adding now the integral with respect to $z$, the first term $J_{2,1}$ can be rewritten as
\begin{multline*}
\int_{\R^d} J_{2,1}(v,y_{\smallsetminus j},z,\lambda) \, dz \, \le \,\int_{\R^d}\int_{\R^d}\vert z- y_j \vert^{\frac{1+\alpha(j-2)+\beta}{1+\alpha(j-1)}+1}\vert D_z \partial_vp_h(v,z- y_j)\vert\\
\times \vert
D_{ y_j}D^\vartheta_{x} \tilde{p}^{\tau,\xi}(t,s,x,y_{\smallsetminus j}, y_j+\lambda(z- y_j))\vert \,
d y_j dz.
\end{multline*}
Fubini Theorem and the change of variables $\tilde{z}=z- y_j$ and $\tilde{y}_j= y_j+\lambda\tilde{z}$ allow then to
split the integrals in the following way:
\begin{multline*}
\int_{\R^d} J_{2,1}(v,y_{\smallsetminus j},z,\lambda) \, dz \\
\le \, \Bigl(\int_{\R^d} \vert \tilde{z}
\vert^{\frac{1+\alpha(j-2)+\beta}{1+\alpha(j-1)}+1}\vert D_{\tilde{z}}\partial_vp_h(v,\tilde{z}) \vert \,
d\tilde{z}\Bigr)\Bigl(\int_{\R^d}\vert D_{\tilde{y}_j} D^\vartheta_{x} \tilde{p}^{\tau,\xi}(t,s,x,
y_{\smallsetminus j},\tilde{y}_j)\vert \, d\tilde{y}_j \Bigr).
\end{multline*}
Noticing now from assumption [\textbf{P}] that
\[\frac{1+\alpha(j-2)+\beta}{1+\alpha(j-1)}+1 \, =\, 1-\frac{\beta-\alpha}{1+\alpha(j-1)}+1 \,< \, 2-(1-\alpha) \, = \, 1+\alpha, \]
we can use the smoothing effect of the heat-kernel $p_h$ (Equation \eqref{Smoothing_effect_of_Heat_Kern}) to show that
\[\int_{\R^d} J_{2,1}(v,y_{\smallsetminus j},z,\lambda) \, dz \,\le \,\frac{v^{\frac{1+\alpha(j-2)+\beta}{\alpha(1+\alpha(j-1))}}}{v}
\int_{\R^d}\vert D_{\tilde{y}_j} D^\vartheta_{x}\tilde{p}^{\tau,\xi}(t,s,x,y_{\smallsetminus j},
\tilde{y}_j) \vert \, d\tilde{y}_j.\]
Remembering Equation \eqref{Proof:Second_Besov_Control}, we can add the integral with respect to $v$ and show that
\begin{multline*}
\int_{0}^{(s-t)^{\delta_j}}v^{\frac{\alpha_j+\beta_j}{\alpha}}\int_{\R^d} J_{2,1}(v,y_{\smallsetminus j},z,\lambda) \, dz\, dv \\
\le
\, (s-t)^{\delta_j\frac{2\beta_j+1}{\alpha}}\int_{\R^d}\vert D_{\tilde{y}_j} D^\vartheta_{x}
\tilde{p}^{\tau,\xi}(t,s,x, y_{\smallsetminus j},\tilde{y}_j)\vert \, d\tilde{y}_j.    
\end{multline*}
Adding the integral with respect to $y_{\smallsetminus j}$, we can conclude with $J_{2,1}$ that
\begin{multline}\label{Proof:Second_Besov_Control_J_{2,1}}
\int_{\R^{(n-1)d}}\int_{0}^{(s-t)^{\delta_j}}v^{\frac{\alpha_j+\beta_j}{\alpha}}\int_{\R^d} J_{2,1}(v,y_{\smallsetminus j},z,\lambda) \,
dz\, dv \, dy_{\smallsetminus j} \\
\le \, C (s-t)^{\delta_j\frac{2\beta_j+1}{\alpha}}\int_{\R^{nd}}\vert D_{ y_j}D^\vartheta_{x}\tilde{p}^{\tau,\xi}(t,s,x,y)
\vert \, dy \, \le \, C (s-t)^{\delta_j\frac{2\beta_j+1}{\alpha}-\frac{1}{\alpha_j}-\sum_{k=1}^{n}\frac{\vartheta_k}{\alpha_k}},
\end{multline}
where, for simplicity, we have changed back the variable $\tilde{y}_j$ with $ y_j$. \newline
To control instead the term $J_{2,2}$ (cf. Equation \eqref{Proof:Second_Besov_Control_J'0}), we can use again Fubini theorem and the
changes of variables $\tilde{z}=z- y_j$, $\tilde{y}_j= y_j+\lambda\tilde{z}$ to
split the integrals and show that
\[\begin{split}
\int_{\R^d} J_{2,2}(v,y_{\smallsetminus j},z,\lambda&) \, dz  \, \le \, \Bigl(\int_{\R^d}\vert\tilde{z}\vert \, \vert D_{\tilde{z}} \partial_v
p_h(v,\tilde{z})\vert\, d\tilde{z}\Bigr)\\
&\quad \times\Bigl(\int_{\R^d}\vert D_{\tilde{y}_j}D^\vartheta_{x} \tilde{p}^{\tau,\xi}
(t,s,x,y_{\smallsetminus j},\tilde{y}_j) \vert \mathbf{d}^{1+\alpha(j-2)+\beta}_{j:n} ((y_{\smallsetminus j},\tilde{y}_j),
\theta_{s,\tau}(\xi)) d\tilde{y}_j\Bigr)\\
&\le \, \frac{1}{v}\int_{\R^d}\vert D_{ y_j}D^\vartheta_{x} \tilde{p}^{\tau,\xi} (t,s,x,y) \vert
\mathbf{d}^{1+\alpha(j-2)+\beta}_{j:n} (y, \theta_{s,\tau}(\xi)) d y_j dv,
\end{split}
\]
where in the second inequality we used the smoothing effect of the heat-kernel $p_h$ (Equation \eqref{Smoothing_effect_of_Heat_Kern}) and
changed back the variable $\tilde{y}_j$ with $ y_j$ for simplicity. It then follows that
\begin{multline*}
\int_{0}^{(s-t)^{\delta_j}}v^{\frac{\alpha_j+\beta_j}{\alpha}}\int_{\R^d} J_{2,2}(v,z,y_{\smallsetminus j}) \, dz dv \\
\le \, (s-t)^{\delta_j\frac{\alpha_j+\beta_j}{\alpha}}\int_{\R^d}\vert D_{ y_j}D^\vartheta_{x} \tilde{p}^{\tau,\xi}
(t,s,x,y) \vert \mathbf{d}^{1+\alpha(j-2)+\beta}_{j:n} (y, \theta_{s,\tau}(\xi)) d y_j.
\end{multline*}
Taking now $\tau=t$ and $\xi=x$, we conclude with $J_{2,2}$ applying the partial smoothing effect
(Equation \eqref{eq:Partial_Smoothing_Effect}) of
$\tilde{p}^{\tau,\xi}$ under the hypothesis $1+\alpha(j-2)+\beta \le  \alpha(1+\alpha(j-1))$ (see Equation \eqref{a1}) to write that
\begin{align}\notag
\int_{\R^{(n-1)d}}\int_{0}^{(s-t)^{\delta_j}}&v^{\frac{\alpha_j+\beta_j}{\alpha}}\int_{\R^d} J_{2,2}(v,z,y_{\smallsetminus j}) \, dz \,
dv dy_{\smallsetminus j} \\
\notag
&\qquad\le \, (s-t)^{\delta_j\frac{\alpha_j+\beta_j}{\alpha}}\int_{\R^{nd}}\vert D_{ y_j}D^\vartheta_{x} \tilde{p}^{\tau,\xi}
(t,s,x,y) \vert \mathbf{d}^{1+\alpha(j-2)+\beta}_{j:n} (y, \theta_{s,\tau}(\xi)) dy \\
&\qquad\le \, C(s-t)^{\delta_j\frac{\alpha_j+\beta_j}{\alpha}+\frac{1+\alpha(j-2)+\beta}{\alpha}- \frac{1}{\alpha_j}-\sum_{k=1}^{n}
\frac{\vartheta_k}{\alpha_k}}.
\label{Proof:Second_Besov_Control_J'2}
\end{align}
Looking back to Equations \eqref{Proof:Second_Besov_Control_I_2}, \eqref{Proof:Second_Besov_Control_J1},
\eqref{Proof:Second_Besov_Control_J_{2,1}} and \eqref{Proof:Second_Besov_Control_J'2}, we can finally choose the right $\delta_j$. Since
$s-t\le T-t<1$ by hypothesis [\textbf{ST}], it is enough to take $\delta_j$ such that the quantities
\[\delta_j\frac{\alpha_j+\beta_j-1}{\alpha}+\frac{1+\alpha(j-2)+\beta}{\alpha}, \qquad  \delta_j\frac{2\beta_j}{\alpha
}, \qquad \delta_j\frac{ 2\beta_j+1}{\alpha}-\frac{1}{\alpha_j}
\]
and
\[\delta_j\frac{\alpha_j+\beta_j}{\alpha}+\frac{1+\alpha(j-2)+\beta}{\alpha}-\frac{1}{\alpha_j}\]
are bigger than $\beta/\alpha$. This is true if, for example, we choose
\[\delta_j \, = \, \frac{[1+\alpha(j-2)][1+\alpha(j-1)]}{1+\alpha(j-2)-\beta}.\]
We have thus concluded th proof.
\end{proof}

\subsection{Proof of Proposition \ref{prop:A_Priori_Estimates}}
We have now all the tools necessary to prove the a priori estimates in Proposition \ref{prop:A_Priori_Estimates}. In Lemma
\ref{lemma:supremum_norm} below, we will state the estimates for the supremum norms of the solution and its non-degenerate gradient while
the controls of the H\"older moduli of the solution and its gradient with respect to the non-degenerate variable are given in Lemmas
\ref{lemma:Holder_modulus_Non-Deg} and \ref{lemma:Holder_modulus_Deg}, respectively.\newline
The Schauder estimates (Theorem \ref{theorem:Schauder_Estimates}) for a solution $u$ in
$L^\infty\bigl(0,T;C^{\alpha+\beta}_{b,d}(\R^{nd})\bigr)$ of IPDE \eqref{Degenerate_Stable_PDE} will then follows immediately.

\begin{lemma}[Supremum Estimates]
\label{lemma:supremum_norm}
Let $u$ be as in Equation \eqref{align:Representation2}. Then, there exists a constant $C\ge1$ such that for any $t$ in $[0,T]$ and any $x$ in $\R^{nd}$,
\[\vert u(t,x) \vert + \vert D_{ x_1}u(t,x)\vert \, \le \, C\Bigl[\Vert u_T \Vert_{C^{\alpha+\beta}_{b,d}} + \Vert f
\Vert_{L^\infty(C^\beta_{b,d})} + \Vert F\Vert_H\Vert u \Vert_{L^\infty(C^{\alpha+\beta}_{b,d})}\Bigr].\]
\end{lemma}
\begin{proof}
As indicated above, we can control the supremum norm of $u$ and its gradient with respect to $ x_1$ analyzing separately the
contributions from the proxy $\tilde{u}^{\tau,\xi}$, that have already been handled in Lemma
\ref{lemma:supremum_norm_proxy}, and those from the perturbative term $R^{\tau,\xi}(s,x)$.
To control the contribution $\int_{t}^{T} D_{ x_1}\tilde{P}^{\tau,\xi}_{s,t}R^{\tau,\xi}(s,x)\, ds$, we start splitting
it up in the following way:
\begin{align}\notag
\int_{t}^{T} D_{ x_1}\tilde{P}^{\tau,\xi}_{s,t}&R^{\tau,\xi}(s,x)\, ds \\
\notag &= \, \sum_{j=1}^{n}\int_{t}^{T} \int_{\R^{nd}}
D_{ x_1}\tilde{p}^{\tau,\xi}(t,s,x,y)\bigl[ F_j(s,y) - F_j(s,\theta_{s,\tau}
(\xi)) \bigr]\cdot D_{ y_j}u(s,y) \, dyds \\
&=:\,\sum_{j=1}^{n}I_j(t,x).\label{Proof:Supremum_Schauder0}
\end{align}
Since by hypothesis $u$ has a proper derivative with respect to the first variable $ x_1$, it is possible to bound $I_1$ through
\[\vert I_1(t,x)\vert \, \le \, C\Vert F\Vert_H \Vert u\Vert_{L^{\infty}(C^{\alpha+\beta}_{b,d})}\int_{t}^{T}\int_{\R^{nd}}
\bigl{\vert}D_{ x_1}\tilde{p}^{\tau,\xi}(t,s,x,y)\bigr{\vert} \mathbf{d}^{\beta}(y,\theta_{s,\tau}(\xi)) \, dy
ds.\]
We take now $(\tau,\xi)=(t,x)$ so that $\theta_{s,\tau}(\xi)=\tilde{m}^{\tau,\xi}_{s,t}$ (cf. Equation
\eqref{eq:identification_theta_m_stable} in Lemma \ref{lemma:identification_theta_m_stable}) and we then use the smoothing effect for the frozen density
$\tilde{p}^{\tau,\xi}$ (Equation \eqref{eq:Smoothing_effects_of_tilde_p}) to show that
\begin{equation}\label{Proof:Supremum_Schauder1}
\vert I_1(t,x)\vert \, \le \,C\Vert F\Vert_H \Vert u\Vert_{L^{\infty}( C^{\alpha+\beta}_{b,d})}
(T-t)^{\frac{\beta+\alpha-1}{\alpha}}.
\end{equation}
Hence, it holds that $\vert I_1(t,x)\vert \le  C\Vert F\Vert_H \Vert u \Vert_{L^{\infty}(C^{\alpha+\beta}_{b,d})}$, since
$T\le1$ and $\alpha+\beta>1$ by assumptions [\textbf{ST}] and [\textbf{P}].\newline
The control for the terms $I_j$ with $j>1$ can be obtained easily from the second Besov control (Lemma \ref{lemma:Second_Besov_COntrols}).
For this reason, we start applying integration by parts formula to show that
\[\vert I_j(t,x) \vert \, = \,
\Bigl{\vert}\int_{t}^{T}\int_{\R^{nd}}D_{ y_j}\cdot\Bigl{\{}D_{ x_1}\tilde{p}^{\tau,\xi}(t,s,x,y)
 \Delta^{\tau,\xi} F_j(s,y)\Bigr{\}} u(s,y) \, dyds\Bigr{\vert},\]
where we exploited the same notations for $\Delta^{\tau,\xi} F_j$ given in \eqref{eq:scorciatoia}.
We can then use identification \eqref{Besov:ident_Holder_Besov} and duality in Besov spaces \eqref{Besov:duality_in_Besov} to write
that
\begin{multline*}
\vert I_j(t,x) \vert \, \le  \\
\Vert u \Vert_{L^\infty(C^{\alpha+\beta}_{b,d})}
\int_{\R^{(n-1)d}}\Bigl{\Vert}D_{ y_j}\cdot\Bigl{\{}D_{ x_1}\tilde{p}^{\tau,\xi}(t,s,x,y_{\smallsetminus j},\cdot)
\Delta^{\tau,\xi} F_j(s, y_{\smallsetminus j},\cdot)\Bigl{\}}
\Bigr{\Vert}_{B^{-(\alpha_j+\beta_j)}_{1,1}} \, dy_{\smallsetminus j}.
\end{multline*}
Taking now $(\tau,\xi)=(t,x)$, the second Besov control (Lemma \ref{lemma:Second_Besov_COntrols}) can be applied to show that
\begin{equation}\label{Proof:Supremum_Schauder2}
\begin{split}
\vert I_j(t,x)\vert \, &\le \, C\Vert F\Vert_H\Vert u \Vert_{L^\infty(C^{\alpha+\beta}_{b,d})}\int_{t}^{T}
(s-t)^{\frac{\beta-1}{\alpha}}\,ds  \\
&\le \, C\Vert F\Vert_H\Vert u \Vert_{L^\infty(C^{\alpha+\beta}_{b,d})}(T-t)^{
\frac{\beta+\alpha-1}{\alpha}}.
\end{split}
\end{equation}
Since $T\le 1$ by assumption [\textbf{ST}], we can conclude that $\vert I_j(t,x)\vert \le C\Vert F\Vert_H\Vert u
\Vert_{L^\infty(C^{\alpha+\beta}_{b,d})}$. \newline
The control on the pertubative term
\[\int_{t}^{T}  \tilde{P}^{\tau,\xi}_{s,t}R^{\tau,\xi}(s,x)\, ds\]
can be obtained in a similar way. Namely, Inequalities \eqref{Proof:Supremum_Schauder1} and \eqref{Proof:Supremum_Schauder2} hold again
with $(T-t)^{\frac{\beta+\alpha-1}{\alpha}}$ replaced by $(T-t)^{\frac{\beta+\alpha}{\alpha}}$.
\end{proof}

As already specified in the previous sub-section, there is a big difference between the non-degenerate case $i=1$, where
$\alpha+\beta$ is in $(1,2)$ and we have to deal with a proper derivative, and the other degenerate situations ($i>1$), where instead
$(\alpha+\beta)/(1+\alpha(i-1))<1$ and the norm is calculated directly on the function. Again, we are going to analyze the two
cases separately. Lemma \ref{lemma:Holder_modulus_proxy_Non-Deg} will focus on the non-degenerate setting ($i=1$) while lemma
\ref{lemma:Holder_modulus_proxy_Deg} will concerns the degenerate one ($i>1$).

Moreover, we will need to divide the proofs in two cases, depending on which regime we are considering. Since the global off-diagonal
regime, i.e.\ when $T-t\le c_0\mathbf{d}^\alpha(x,x')$, will work essentially as the already shown Schauder estimates
(Proposition \ref{prop:Schauder_Estimates_for_proxy}) for the proxy, the proof will be quite shorter.\newline
Instead, in the global diagonal case, such that $T-t \ge c_0\mathbf{d}^\alpha(x,x')$, when a time integration is involved (for example in
the control of the frozen Green kernel or the perturbative term), two different situations appear. There are again a local off-diagonal
regime if $s-t \le c_0\mathbf{d}^\alpha(x,x')$ and a local diagonal regime when $s-t \ge c_0\mathbf{d}^\alpha(x,x')$. In order to handle
these terms properly, the key tool is to be able to change the freezing points depending on which regime we are. It seems reasonable that,
when the spatial points are in a local diagonal regime, the auxiliary frozen
densities are considered for the same freezing parameter and conversely that, in the local off-diagonal regime, the densities are frozen
along their own spatial argument. For this reason, we have postponed the relative proofs in two specific sub-sections.

Before presenting the main results of this section, we are going to state three auxiliary estimates we will need
below. We refer to the Section A.$2$ for a precise proof of these results. \newline
The first one concerns the sensitivity of the H\"older flow $\theta_{s,t}$ with respect to the initial point. Indeed,
\begin{lemma}[Controls on the Flows]
\label{lemma:Controls_on_Flow1}
Let $t<s$ be two points in $[0,T]$ and $x,x'$ two points in $\R^{nd}$. Then, there exists a constant $C\ge 1$ such that
\[\mathbf{d}(\theta_{s,t}(x), \theta_{s,t}(x')) \, \le \, C\Vert F\Vert_H \bigl[\mathbf{d}(x,x')+
(s-t)^{1/\alpha}\bigr].\]
\end{lemma}

The second result is the following:

\begin{lemma}
\label{lemma:Controls_on_means1}
Let $t<s$ be two points in $[0,T]$ and $x,x'$ two points in $\R^{nd}$ and $y,y'$ two points in $\R^{nd}$ such that
$ y_1= y'_1$. Then, there
exists a constant $C\ge1$ such that
\[\bigl{\vert}(\tilde{m}^{t,x}_{s,t}(y)-\tilde{m}^{t,x'}_{s,t}(y'))_1 \bigr{\vert} \, \le
\, C\Vert F\Vert_H\Bigl[(s-t)\mathbf{d}^{\beta}(x,x')+(s-t)^{\frac{\alpha+\beta}{\alpha}}\Bigr].\]
\end{lemma}

Finally, the impact of the freezing point in the linearization procedure is the argument of this last Lemma. Namely,

\begin{lemma}
\label{lemma:Controls_on_means}
Let $t$ be in $[0,T]$ and $x,x'$ two points in $\R^{nd}$. Then, there exists a constant $C\ge 1$
such that
\[\mathbf{d}(\tilde{m}^{t,x}_{t_0,t}(x'),\tilde{m}^{t,x'}_{t_0,t} (x')) \, \le \, Cc_0^{\frac{1}{1+\alpha(n-1)}}\Vert
F\Vert_H\mathbf{d}(x,x')\]
where $t_0$ is the change of regime time defined in \eqref{eq:def_t0}.
\end{lemma}

Thanks to the above controls, we will eventually prove the following results:

\begin{lemma}[Controls on H\"older Moduli: Non-Degenerate]
\label{lemma:Holder_modulus_Non-Deg}
Let $x,x'$ be in $\R^{nd}$ such that $ x_j= x'_j$ for any $j\neq 1$ and $u$ as in
 Equation \eqref{align:Representation2}. Then, there exists a constant $C\ge 1$ such that for any $t$ in $[0,T]$,
\begin{multline*}
\bigl{\vert} D_{ x_1}u(t,x)- D_{ x_1}u(t,x')\vert\, \le \, C\Bigl{\{}c_0^{\frac{\alpha+\beta-2}{\alpha}}\bigl(\Vert u_T\Vert_{C^{\alpha+\beta}} + \Vert f \Vert_{L^\infty(C^\beta)}\bigr) \\
+\bigl(c_0^{\frac{\alpha+\beta-1}{1+\alpha(n-1)}} +c_0^{
\frac{\alpha+\beta-2 }{\alpha}}\Vert F\Vert_{H}\bigr)\Vert u \Vert_{L^\infty(C^{\alpha+\beta}_{b,d})}\Bigr{\}} \mathbf{d}^{\alpha+\beta-1}
(x,x').
\end{multline*}
\end{lemma}

We can point out now the analogous result in the degenerate setting.

\begin{lemma}[Controls on H\"older Moduli: Degenerate]
\label{lemma:Holder_modulus_Deg}
Let $i$ be in $\llbracket 1,n\rrbracket$ and $x,x'$ in $\R^{nd}$ such that $x_j=x'_j$ for any $j\neq i$ and  $u$ as in
Equation \eqref{align:Representation2}. Then, there exists a constant $C\ge1$ such that
for any $t$ in $[0,T]$,
\begin{multline*}
\bigl{\vert} u(t,x)- u(t,x') \vert\,\le \, C\Bigl{\{}c_0^{\frac{\beta-\gamma_i}{\alpha}}\bigl(\Vert u_T
\Vert_{C^{\alpha+\beta}} + \Vert f \Vert_{L^\infty(C^\beta)}\bigr) \\
+ \bigl(c_0^{\frac{\alpha+\beta}{1+\alpha(n-1)}}+c_0^{\frac{\beta-\gamma_i }{\alpha}}\Vert F\Vert_{H}\bigr)\Vert u \Vert_{L^\infty(C^{\alpha+\beta}_{b,d})}\Bigr{\}}\mathbf{d}^{\alpha+\beta}
(x,x').
\end{multline*}
\end{lemma}

\subsubsection{Off-diagonal regime}
We focus here on the proof of the controls on the H\"older moduli, either in the non-degenerate setting (Lemma
\ref{lemma:Holder_modulus_Non-Deg}) and in the degenerate one (Lemma \ref{lemma:Holder_modulus_Deg}), when a
off-diagonal regime is assumed. For this reason, all the statements presented in this sub-section will tacitly assume that $T-t \le
c_0\mathbf{d}^\alpha(x,x')$ for some given $(t,x,x')$ in $[0,T]\times \R^{2nd}$. \newline
To show these two controls, we will need to adapt the auxiliary estimates above to the off-diagonal regime case we consider here. Namely,
\begin{align}
\label{eq:Sensitivity_on_flow}
\mathbf{d}(\tilde{m}^{t,x}_{T,t} (x),\tilde{m}^{t,x'}_{T,t}(x')) \, = \, \mathbf{d}(\theta_{T,t}(x),
\theta_{T,t}(x')) \, &\le \, C\Vert F\Vert_H \mathbf{d}(x,x'); \\
\label{eq:Sensitivity_on_flow1}
\text{if } x_1 \, = \,  x_1', \quad \bigl{\vert} \bigl(\tilde{m}^{t,x}_{T,t} (x)-\tilde{m}^{t,x'}_{T,t}(x')\bigr)_1  \bigr{\vert} \, &\le \, C\Vert F\Vert_H\mathbf{d}^{\alpha+\beta}(x,x')
\end{align}
They can be obtained easily from Equation \eqref{eq:identification_theta_m_stable} in Lemma \ref{lemma:identification_theta_m_stable}
and the sensitivity controls (Lemmas \ref{lemma:Controls_on_Flow1} and \ref{lemma:Controls_on_means1}, respectively), taking $s=T$ and
$(y,y')=(x,x')$.

\paragraph{Proof of Lemma \ref{lemma:Holder_modulus_Non-Deg} in the off-diagonal regime.}
From the Duhamel-type Expansion \eqref{align:Representation2_deriv}, we can represent a mild solution $u$ of IPDE
\eqref{Degenerate_Stable_PDE} as
\begin{multline*}
\vert D_{ x_1}u(t,x)-D_{ x_1}u(t,x')\vert \\
\le \, \bigl{\vert}D_{ x_1}\tilde{P}^{\tau,\xi}_{T,t}u_T(x) -
D_{ x_1}\tilde{P}^{\tau' ,\xi'}_{T,t} u_T(x')\bigr{\vert} + \bigl{\vert}D_{ x_1}\tilde{G}^{\tau,\xi}_{T,t}f(t,x) - D_{
x_1}\tilde{G}^{\tau' ,\xi'}_{T,t} f(t,x') \bigr{\vert}\\
+\Bigl{\vert}\int_{t}^{T}D_{ x_1}\tilde{P}^{\tau,\xi}_{s,t}R^{\tau,\xi}(s,x) -
D_{ x_1}\tilde{P}^{\tau' ,\xi'}_{s,t} R^{\tau' ,\xi'} (s,x') \, ds\Bigr{\vert},
\end{multline*}
for any fixed $(\tau,\xi),(\tau',\xi')$ in $[0,T]\times \R^{nd}$
After possible differentiations, we will choose $\tau=\tau'=t$, $\xi = x$ and $\xi' = x'$ in order to exploit the
sensitivity Controls \eqref{eq:Sensitivity_on_flow1} and \eqref{eq:Sensitivity_on_flow}.

\emph{Control on the frozen semigroup.} It can be handled following the analogous part in the proof of the
H\"older control for the proxy (Lemma
\ref{lemma:Holder_modulus_proxy_Non-Deg}).  The only difference is that we cannot control
\[\mathbf{d}(\tilde{m}^{\tau,\xi}_{T,t} (x),\tilde{m}^{\tau' ,\xi'}_{T,t} (x'))\]
in Equation \eqref{zz1} using the affinity of the mapping $x \to \tilde{m}^{\tau,\xi}_{T,t}(x)$, since the two freezing
point are now different. Instead, we can take $\tau = \tau' = t$, $\xi = x$ and $\xi' = x'$ and apply the sensitivity
control \eqref{eq:Sensitivity_on_flow} to write that
\[\mathbf{d}(\tilde{m}^{\tau,\xi}_{T,t} (x),\tilde{m}^{\tau' ,\xi'}_{T,t} (x')) \, = \,
\mathbf{d}(\theta_{T,t}(x), \theta_{T,t}(x')) \, \le \, C\Vert F\Vert_H \mathbf{d}(x,x').\]

\emph{Control on the Green kernel.} It follows immediately from the proof of the H\"older control (Lemma
\ref{lemma:Holder_modulus_proxy_Non-Deg}) for the proxy, noticing that $t_0=T$, since we are in the off-diagonal regime.

\emph{Control on the perturbative error.} Since we do not exploit the difference of the spatial points $(x,x')$ in the
off-diagonal regime but instead we control the two contributions separately, we can rely  on the controls on the supremum norms we have
already shown in Lemma \ref{lemma:supremum_norm}. Namely, we start writing that
\begin{multline}\label{z1}
\Bigl{\vert}\int_{t}^{T}D_{ x_1}\tilde{P}^{\tau,\xi}_{s,t}R^{\tau,\xi}(s,x)- D_{ x_1}
\tilde{P}^{\tau' ,\xi'}_{s,t} R^{\tau' ,\xi'} (s,x')\, ds\Bigr{\vert} \\
\le \, \Bigl{\vert}\int_{t}^{T}D_{ x_1}\tilde{P}^{\tau,\xi}_{s,t}R^{\tau,\xi}(s,x)\, ds\Bigr{\vert}+
\Bigl{\vert}\int_{t}^{T}D_{ x_1} \tilde{P}^{\tau' ,\xi'}_{s,t} R^{\tau' ,\xi'} (s,x')\, ds\Bigr{\vert}.
\end{multline}
Then, we can follow the same reasonings of Lemma \ref{lemma:supremum_norm} concerning the remainder term (cf. Equations
\eqref{Proof:Supremum_Schauder0}, \eqref{Proof:Supremum_Schauder1} and \eqref{Proof:Supremum_Schauder2}) to
show that
\begin{equation}\label{z2}
\Bigl{\vert}\int_{t}^{T}D_{ x_1}\tilde{P}^{\tau,\xi}_{s,t}R^{\tau, \xi}(s,x) \, ds\Bigr{\vert} \, \le \, C\Vert
F\Vert_H\Vert u \Vert_{L^{\infty}(C^{\alpha+\beta}_{b,d})} (T-t)^{\frac{\alpha+\beta-1}{\alpha}}.
\end{equation}
Using it in the above Equation \eqref{z1}, we can finally conclude that
\begin{equation}\label{z3}
\Bigl{\vert}\int_{t}^{T}D_{ x_1}\tilde{P}^{\tau,\xi}_{s,t}R^{\tau,\xi}(s,x)- D_{ x_1}
\tilde{P}^{\tau' ,\xi'}_{s,t} R^{\tau' ,\xi'} (s,x')\, ds\Bigr{\vert} \, \le \,
C\Vert F\Vert_H\Vert u \Vert_{L^{\infty}(C^{\alpha+\beta}_{b,d})} \mathbf{d}^{\alpha+\beta-1}(x,x'),
\end{equation}
remembering that we assumed to be in the off-diagonal regime, i.e.\ $T-t\le c_0\mathbf{d}^\alpha(x,x')$ for some $c_0\le1$.

\paragraph{Proof of Lemma \ref{lemma:Holder_modulus_Deg} in the off-diagonal regime.}
As done before, we are going to analyze separately the single terms appearing from the Duhamel-type Representation
\eqref{align:Representation2} of a solution $u$:
\begin{multline*}
\vert u(t,x)-u(t,x')\vert \, \le \, \bigl{\vert}\tilde{P}^{\tau,\xi}_{T,t}u_T(x) - \tilde{P}^{\tau' ,\xi'}_{T,t} u_T
(x')\bigr{\vert} + \bigl{\vert}\tilde{G}^{\tau,\xi}_{T,t}f(t,x) - \tilde{G}^{\tau' ,\xi'}_{T,t} f(t,x')\bigr{\vert}\\
+\Bigl{\vert}\int_{t}^{T}\tilde{P}^{\tau,\xi}_{s,t}R^{\tau,\xi}(s,x)-\tilde{P}^{\tau' ,\xi'}_{s,t}
 R^{\tau' ,\xi'} (s,x') \, ds\Bigr{\vert},
\end{multline*}
for some $(\tau,\xi),(\tau',\xi')$ in $[0,T]\times \R^{nd}$ fixed but to be chosen later as $\tau = \tau' = t$, $\xi = x$
and $\xi' = x'$.

\emph{Control on the frozen semigroup.} We can essentially follow the proof of the H\"older control (Lemma
\ref{lemma:Holder_modulus_proxy_Deg}) for the proxy. However, this time we cannot exploit the affinity of the mapping $x \to
\tilde{m}^{\tau,\xi}_{T,t}(x)$ to control the difference
\[\bigl{\vert} u_T(\tilde{m}^{\tau,\xi}_{T,t}(x)-z)-u_T(\tilde{m}^{\tau,\xi}_{T,t}(x')-z)\bigr{\vert}.\]
Instead, we notice now that we can bound it as
\begin{multline*}
    \bigl{\vert} u_T(\tilde{m}^{\tau,\xi}_{T,t}(x)-z)-u_T(\tilde{m}^{\tau,\xi}_{T,t}(x')-z)\bigr{\vert} \\
\le \, C\Vert u_T \Vert_{C^{\alpha+\beta}_{b,d}}\Bigl[\mathbf{d}^{\alpha+\beta}\bigl(\tilde{m}^{\tau,\xi}_{T,t}(x),\tilde{m}^{\tau,\xi}_{T,t}(x')\bigr)+\bigl{\vert}\bigl(\tilde{m}^{\tau,\xi}_{T,t}(x)-\tilde{m}^{\tau,\xi}_{T,t}(x')\bigr)_1 \bigr{\vert}\Bigr],
\end{multline*}
since $u_T$ is differentiable and thus Lipschitz continuous, in the first non-degenerate variable.\newline
Taking now  $\tau = \tau' = t$, $\xi = x$ and $\xi' = x'$, we can use the sensitivity Controls
\eqref{eq:Sensitivity_on_flow}-\eqref{eq:Sensitivity_on_flow1} (noticing that by assumption, $ x_1= x'_1$) to write that
\[\bigl{\vert} u_T(\tilde{m}^{\tau,\xi}_{T,t}(x)-z)-u_T(\tilde{m}^{\tau,\xi}_{T,t}(x')-z)\bigr{\vert} \, \le \,
C\Vert F\Vert_H\Vert u_T \Vert_{C^{\alpha+\beta}_{b,d}}\mathbf{d}^{\alpha+\beta}(x,x').\]


\emph{Control on the Green kernel.} It can be obtained following the analogous part in the proof of the H\"older control (Lemma
\ref{lemma:Holder_modulus_proxy_Deg}) for the proxy. Similarly to the paragraph ``Control on the frozen semigroup'' in the previous proof, we
need to take $(\tau,\xi)= (t,x)$, $(\tau,\xi') = (t,x)$ and apply the sensitivity Control
\eqref{eq:Sensitivity_on_flow} to bound the term
\[\mathbf{d}(\tilde{m}^{\tau,\xi}_{T,t} (x),\tilde{m}^{\tau' ,\xi'}_{T,t} (x'))\]
appearing in Equation \eqref{zz2}.

\emph{Control on the perturbative error.} The proof of this estimate essentially matches the previous, analogous one in the
non-degenerate setting. Namely, Equations \eqref{z1}, \eqref{z2} and \eqref{z3} hold again with $(T-t)^{\frac{\beta+\alpha}{\alpha}}$
instead of $(T-t)^{\frac{\beta+\alpha-1}{\alpha}}$.

\subsubsection{Diagonal regime}

Since the aim of this section is to prove Lemmas \ref{lemma:Holder_modulus_Non-Deg} and \ref{lemma:Holder_modulus_Deg} when a
diagonal regime is assumed, we will assume from this point further that $T-t \ge c_0\mathbf{d}^\alpha(x,
x')$ for some given $(t,x,x')$ in $[0,T]\times \R^{2nd}$. \newline
As preannounced in the introduction of this section, we need here a modification of the Duhamel-type Representation
\eqref{eq:Expansion_along_proxy} that allows to change the freezing points along the time integration variable. Remembering the previous
notations for $\tilde{G}^{\tau,\xi}_{r,v}$ and $R^{\tau,\xi}$ in \eqref{eq:def_Green_Kernel_stable} and \eqref{eq:def_remainder_regul}
respectively, it holds that

\begin{lemma}[Change of frozen point]
Let $(\tau,\xi)$ be a freezing couple in $[0,T]\times\R^{nd}$ and $\tilde{\xi}$ another freezing point in $\R^{nd}$. Then, any
classical solution $u$ in $L^\infty(0,T;C^{\alpha+\beta}_{b,d}(\R^{nd}))$ of IPDE \eqref{Degenerate_Stable_PDE} can be represented for
any $(t,x)$ in $[0,T]\times\R^{nd}$ as
\begin{multline}\label{eq:Change_of_Freez_point_HOlder}
u(t,x) \, = \, \tilde{P}^{\tau,\tilde{\xi}}_{T,t}u_T(x) + \tilde{G}^{\tau,\xi}_{t_0,t}f(t,x) +
\tilde{G}^{\tau,\tilde{\xi}}_{T,t_0}f(t,x) \\
+ \int_{t}^{t_0}  \tilde{P}^{\tau,\xi}_{s,t}R^{\tau,\xi} (s,x)\, ds +\int_{t_0}^{T} \tilde{P}^{\tau,\tilde{\xi}}_{s,t}
R^{\tau,\tilde{\xi}}(s,x)\, ds + \tilde{P}^{\tau,\xi}_{t_0,t}u(t_0,x) -\tilde{P}^{\tau,\tilde{\xi}}_{t_0,t}u(t_0,x),
\end{multline}
where $t_0$ is the change of regime time defined in \eqref{eq:def_t0}.
\end{lemma}
\begin{proof}
Fixed $t$ in $(0,T)$, we start considering another point $r$ in $(t,T)$. On $(0,r)$, it is clear that $u$ is again a mild
solution of IPDE \eqref{Degenerate_Stable_PDE} but with terminal condition $u(r,x)$. Then, Duhamel Expansion
\eqref{eq:Expansion_along_proxy} can be applied with respect to the frozen couple $(\tau,\xi)$, allowing us to write that
\[u(t,x) \, = \, \tilde{P}^{\tau,\xi}_{r,t}u_T(x) + \int_{t}^{r}\tilde{P}^{\tau,\xi}_{s,t}f(s,x) \, ds + \int_{t}^{r}
\tilde{P}^{\tau,\xi}_{s,t}R^{\tau,\xi} u(s,x)\, ds.\]
Noticing that $u$ is independent from $r$, it is possible now to differentiate the above equation with respect to $r$ in $(t,T)$ in order to show that
\begin{equation}\label{eq:change_frez_point0}
0 \, = \, \partial_r \bigl[\tilde{P}^{\tau,\xi}_{r,t}u(r,x)\bigr] + \tilde{P}^{\tau,\xi}_{r,t}f(r,x) +
\tilde{P}^{\tau,\xi}_{r,t}R^{\tau,\xi}(r,x).
\end{equation}
We highlight now that the above expression holds for any chosen frozen couple $(\tau,\xi)$ and any fixed time $r$. Thus, it is possible
to integrate it with respect to $r$ for a fixed $\xi$ between $t$ and $t_0$ and for another frozen point
$\tilde{\xi}$ between $t_0$ and $T$, leading to
\begin{multline*}
0 \, = \, \tilde{P}^{\tau,\xi}_{t_0,t}u(t_0,x) - \tilde{P}^{\tau,\xi}_{t,t}u(t,x) +\int_{t}^{t_0}\tilde{P}^{\tau,\xi}_{r,t}f(r,x) \,  dr +  \int_{t}^{t_0} \tilde{P}^{\tau,\xi}_{r,t}R^{\tau,\xi}(r,x) \, dr  \\
+ \tilde{P}^{\tau,\tilde{\xi}}_{T,t}u(T,x) - \tilde{P}^{\tau,\tilde{\xi}}_{t_0,t}u(t_0,x) + \int_{t_0}^{T}
\tilde{P}^{\tau,\tilde{\xi}}_{r,t}f(r,x) \,  dr +  \int_{t_0}^{T} \tilde{P}^{\tau,\tilde{\xi}}_{r,t}
R^{\tau,\tilde{\xi}}(r,x) \,  dr.
\end{multline*}
With our previous notations, the above expression can be finally rewritten as
\begin{multline*}
0 \, = \, \tilde{P}^{\tau,\xi}_{t_0,t}u(t_0,x) - u(t,x) + \tilde{G}^{\tau,\xi}_{t_0,t}f(t,x) +  \int_{t}^{t_0}
\tilde{P}^{\tau,\xi}_{r,t}R^{\tau,\xi}(r,x) \, dr  \\
+ \tilde{P}^{\tau,\tilde{\xi}}_{T,t}u_T(x) - \tilde{P}^{\tau,\tilde{\xi}}_{t_0,t}u(t_0,x) +
\tilde{G}^{\tau,\tilde{\xi}}_{T,t_0}f(t,x) +  \int_{t_0}^{T}
\tilde{P}^{\tau,\tilde{\xi}}_{r,t}R^{\tau,\tilde{\xi}}(r,x) \,  dr
\end{multline*}
and we have concluded.
\end{proof}

Similarly to the off-diagonal case, we are going to apply the auxiliary estimates associated with the proxy (Lemmas
\ref{lemma:Controls_on_means1} and \ref{lemma:Controls_on_means}) in the current diagonal regime. Namely, taking $s=t_0$ and
$(y,y')=(x,x)$ in Lemma \ref{lemma:Controls_on_means1}, we know that there exists a constant
$C\ge1$ such that for any $t$ in $[0,T]$ and any $x,x'$ in $\R^{nd}$,
\begin{equation}\label{eq:Sensitivity_on_mean1} \text{if } x_1 \, = \,  x'_1, \quad
\bigl{\vert}(\tilde{m}^{t,x}_{t_0,t}(x)-
\tilde{m}^{t,x'}_{t_0,t}(x))_1\bigr{\vert} \, \le \, C\Vert F\Vert_H\mathbf{d}^{\alpha+\beta}(x,x').
\end{equation}
Moreover, in order to control the perturbative term when a local diagonal regime appears, i.e.\ when the time integration variable $s$ is in
$[t_0,T]$, we will quite often use a Taylor expansion on the frozen density. To be able to exploit the already proven controls, such that
the smoothing effect for the frozen density (Equation \eqref{eq:Smoothing_effects_of_tilde_p}) or the second Besov control (Lemma
\ref{lemma:Second_Besov_COntrols}), we will need the following:
\begin{equation}
\label{eq:translation_inv_for_density}
\text{if } s-t\ge c_0\mathbf{d}^\alpha(x,x'), \quad \bigl{\vert} D^\vartheta_{x}\tilde{p}^{\tau,\xi'} (t,s,x+ \lambda(x'-
x),y)\bigr{\vert} \, \le \, C \bigl{\vert}D^\vartheta_{x}\tilde{p}^{\tau,\xi'} (t,s,x,y)\bigr{\vert},
\end{equation}
for any multi-index $\vartheta$ in $\N^d$ such that $\vert \vartheta\vert \le 2$ and any $\lambda$ in $[0,1]$.
The proof of these results can be found in Section A.$2$.

We are now ready to prove Lemmas \ref{lemma:Holder_modulus_Non-Deg} and \ref{lemma:Holder_modulus_Deg} when a global
diagonal regime is considered.

\paragraph{Proof of Lemma \ref{lemma:Holder_modulus_Non-Deg} in the diagonal regime.}
We start recalling that in Lemma \ref{lemma:Holder_modulus_Non-Deg} we assumed fixed a time $t$ in $[0,T]$ and two spatial points
$x,x'$ in $\R^{nd}$ such that $ x_j= x'_j$ if $j \neq 1$.\newline
From the above Representation \eqref{eq:Change_of_Freez_point_HOlder} and the Duhamel-type Formula \eqref{eq:Expansion_along_proxy}, we know
that
\begin{multline*}
D_{x_1}u(t,x) - D_{x_1}u(t,x') \, = \, \Bigl(D_{x_1}\tilde{P}^{\tau,\tilde{\xi}}_{T,t}u_T(x) -
D_{x_1}\tilde{P}^{\tau' ,\xi'}_{T,t} u_T(x')\Bigr) \\
+ \, \Bigl(D_{x_1}\tilde{G}^{\tau,\xi}_{t_0,t}f(t,x)+ D_{x_1}\tilde{G}^{\tau,\tilde{\xi}}_{T,t_0}f(t,x)
-D_{x_1}\tilde{G}^{\tau' ,\xi'}_{T,t} f(t,x')\Bigr) \\
+ \Bigl(\int_{t}^{t_0}  D_{x_1}\tilde{P}^{\tau,\xi}_{s,t}R^{\tau,\xi}(s,x)\, ds +\int_{t_0}^{T}D_{x_1}
\tilde{P}^{\tau,\tilde{\xi}}_{s,t}R^{\tau,\tilde{\xi}}(s,x)\, ds -\int_{t}^{T}  D_{x_1}\tilde{P}^{\tau' ,
\xi'}_{s,t} R^{\tau' ,\xi'} (s,x') \, ds\Bigr) \\
+ \Bigl(D_{x_1}\tilde{P}^{\tau,\xi}_{t_0,t}u(t_0,x)
-D_{x_1}\tilde{P}^{\tau,\tilde{\xi}}_{t_0,t}u(t_0,x)\Bigr),
\end{multline*}
for some freezing couples $(\tau,\xi),(\tau,\tilde{\xi}),(\tau',\xi')$ in $[0,T]\times \R^{nd}$ fixed but to be chosen later.
To help the readability of the following, we assume from this point further $\tau=\tau'$ and $\tilde{\xi}=\xi'$.

\emph{Control on frozen semigroup.} We start focusing on the control of the frozen semigroup, i.e.
\[\bigl{\vert}D_{x_1}\tilde{P}^{\tau,\xi'}_{T,t}u_T(x) - D_{x_1}\tilde{P}^{\tau,\xi'}_{T,t}u_T(x')\bigr{\vert}.\]
Since the freezing couples coincide, the control on the frozen semigroup can be obtained following the proof of the H\"older control (Lemma
\ref{lemma:Holder_modulus_proxy_Non-Deg}) for the proxy.

\emph{Control on the Green kernel.} As done before, we split the analysis with respect to the change of regime time $t_0$. Namely, we write
\begin{multline*}
\bigl{\vert} D_{ x_1}\tilde{G}^{\tau,\xi}_{t_0,t}f(t,x)+
D_{ x_1}\tilde{G}^{\tau,\tilde{\xi}}_{T,t_0}f(t,x) -D_{ x_1}\tilde{G}^{\tau,\xi'}_{T,t}f(t,x')\bigl{\vert} \\
\le \, \bigl{\vert} D_{ x_1}\tilde{G}^{\tau,\xi}_{t_0,t}f(t,x)-D_{ x_1}\tilde{G}^{\tau,\xi'}_{t_0,t}
f(t,x') \bigl{\vert}+ \bigl{\vert} D_{ x_1}\tilde{G}^{\tau,\tilde{\xi}}_{T,t_0}f(t,x) -
D_{ x_1}\tilde{G}^{\tau,\xi'}_{T,t_0}f(t,x')\bigl{\vert}.
\end{multline*}
While in the local off-diagonal regime, the first term in the r.h.s.\ of the above expression can be handled as in the
global off-diagonal regime, the local diagonal regime contribution represented by
\[\bigl{\vert} D_{ x_1} \tilde{G}^{\tau,\tilde{\xi}}_{T,t_0}f(t,x) - D_{ x_1}
\tilde{G}^{\tau,\xi'}_{T,t_0}f(t,x')\bigr{\vert} \, = \,\bigl{\vert} D_{ x_1} \tilde{G}^{\tau,\xi'}_{T,t_0}f(t,x) -
D_{ x_1} \tilde{G}^{\tau,\xi'}_{T,t_0}f(t,x')\bigr{\vert}\]
since $\tilde{\xi}=\xi'$, can be controlled following again the proof of the H\"older control (Lemma
\ref{lemma:Holder_modulus_proxy_Non-Deg}) for the proxy.

\emph{Control on the discontinuity term.} We can now focus on the contribution
\[\bigl{\vert}D_{x_1}\tilde{P}^{\tau,\xi}_{t_0,t}u(t_0,x)
-D_{x_1}\tilde{P}^{\tau,\tilde{\xi}}_{t_0,t}u(t_0,x)\bigr{\vert},\]
arising from the change of freezing point in the Representation \eqref{eq:Change_of_Freez_point_HOlder}.\newline
Since at fixed time $t_0$, the function $u$ shows the same spatial regularity of $u_T$, this control can be handled following the paragraph in
the proof of the H\"older control for the proxy (Lemma \ref{lemma:Holder_modulus_proxy_Non-Deg}) concerning the frozen semigroup in the
off-diagonal regime. The only main difference is in Equation \eqref{zz1} where, this time, we need to take $(\tau,\xi,\xi') =
(t,x,x')$ and exploit the sensitivity estimate (Lemma \ref{lemma:Controls_on_means}) to control the quantity
\[\mathbf{d}(\tilde{m}^{\tau,\xi}_{t_0,t}(x),\tilde{m}^{\tau,\xi'}_{t_0,t}(x)).\]
In the end, it is possible to show again (cf. Equation \eqref{a2}) that
\[\bigl{\vert}D_{ x_1}\tilde{P}^{\tau,\xi}_{t_0,t}u(t_0,x)-D_{x_1}\tilde{P}^{\tau,\tilde{\xi}}_{t_0,t}
u(t_0,x)\bigr{\vert} \le \, C\Vert u\Vert_{L^\infty(C^{\alpha+\beta}_{b,d})}c_0^{\frac{\alpha + \beta
-1}{\alpha}}\mathbf{d}^{\alpha+\beta-1}(x,x').\]

\emph{Control on the perturbative term.} We start splitting the analysis into two cases with respect to the
critical time $t_0$ giving the change of regime. Namely, we write
\[
\begin{split}
\Bigl{\vert}\int_{t}^{t_0}  D_{ x_1}\tilde{P}^{\tau,\xi}_{s,t}R^{\tau,\xi}(s,x)\, ds +\int_{t_0}^{T}
&D_{x_1}\tilde{P}^{\tau,\xi'}_{s,t}R^{\tau,\xi'}(s,x)\, ds -\int_{t}^{T}  D_{ x_1}\tilde{P}^{
\tau,\xi'}_{s,t}R^{\tau,\xi'}(s,x') \, ds\Bigl{\vert} \\
&\quad\le \, \Bigl{\vert}\int_{t}^{t_0}D_{ x_1}\tilde{P}^{\tau,\xi}_{s,t}R^{\tau,\xi}(s,x)-D_{ x_1}
\tilde{P}^{\tau,\xi'}_{s,t}R^{\tau,\xi'}(s,x')\, ds\Bigr{\vert} \\
&\qquad \,\,+\Bigl{\vert}\int_{t_0}^{T}D_{ x_1}\tilde{P}^{\tau,\xi'}_{s,t}R^{\tau,\xi'}(s,x) -D_{ x_1}\tilde{P}^{\tau,\xi'}_{s,t}R^{\tau,\xi'}(s,x') \,
ds\Bigl{\vert}.
\end{split}
\]
We then notice that the local off-diagonal regime represented by
\[\Bigl{\vert}\int_{t}^{t_0}D_{ x_1}\tilde{P}^{\tau,\xi}_{s,t}R^{\tau,\xi}(s,x)-D_{ x_1}
\tilde{P}^{\tau,\xi'}_{s,t}R^{\tau,\xi'}(s,x')\, ds\Bigr{\vert}\]
can be handled following the proof in the global off-diagonal regime of Lemma \ref{lemma:Holder_modulus_Non-Deg}.\newline
We can then focus our attention on the local diagonal regime, i.e.
\[\Bigl{\vert}\int_{t_0}^{T}D_{ x_1}\tilde{P}^{\tau,\xi'}_{s,t}R^{\tau,\xi'}(s,x)
-D_{ x_1}\tilde{P}^{\tau,\xi'}_{s,t}R^{\tau,\xi'}(s,x') \, ds\Bigl{\vert}.\]
Since the freezing couples coincide, we can use a Taylor expansion with respect to the first variable $ x_1$ in order to write that
\begin{multline*}
\Bigl{\vert}\int_{t_0}^{T}D_{ x_1}\tilde{P}^{\tau,\xi'}_{s,t}R^{\tau,\xi'}(s,x)-D_{ x_1}\tilde{P}^{
\tau,\xi'}_{s,t}R^{\tau,\xi'}(s,x') \, ds\Bigl{\vert}\\
= \, \Bigl{\vert} \int_{t_0}^{T}\int_{\R^{nd}}\int_{0}^{1}D^2_{ x_1} \tilde{p}^{\tau,\xi'} (t,s,x+
\lambda(x'-x),y)(x'-x)_1R^{\tau,\xi'}(s,y) \, dy ds d\lambda \Bigr{\vert}.
\end{multline*}
Noticing that we are integrating from $t_0$ to $T$, Equation \eqref{eq:translation_inv_for_density} can be rewritten as
\begin{align}
\Bigl{\vert}\int_{t_0}^{T}D_{ x_1}\tilde{P}^{\tau,\xi'}_{s,t}R^{\tau,\xi'}(s,x)&-D_{ x_1}\tilde{P}^{
\tau,\xi'}_{s,t}R^{\tau,\xi'}(s,x') \, ds\Bigl{\vert} \notag\\\notag
&\le \, \vert(x'-x)_1\vert\sum_{j=1}^{n}\int_{0}^{1}\int_{t_0
}^{T}\Bigl{\vert} \int_{\R^{nd}}D^2_{ x_1}\tilde{p}^{\tau,\xi'} (t,s,x,y)\\ \notag
&\qquad\qquad\quad\times\Bigl{\{}\bigl[ F_j(s,y)- F_j(s,\theta_{s,t}(\xi'))\bigr]\cdot D_{ y_j}u(s,y)\Bigr{\}} \, dy
\Bigr{\vert} ds d\lambda \\
\label{Proof:H\"older_Perturbative_term_Id0}
&=: \, \vert(x-x')_1\vert\sum_{j=1}^{n}\int_{t_0}^{T}I^d_j(s) ds.
\end{align}
As done before, we are going to treat separately the cases $j=1$ and $j>1$. In the first case, the term $I^{d}_1$ can be easily controlled by
\begin{multline}\label{Proof:H\"older_Perturbative_term_Id1}
I^{d}_1(s) \, \le \, \Vert D_{ y_1}u\Vert_{L^\infty(L^\infty)} \int_{\R^{nd}} \bigl{\vert}
D^2_{ x_1}\tilde{p}^{\tau,\xi'}(t,s,x,y)\bigr{\vert}\,\bigl{\vert} F_1(s,y)-
 F_1(s,\theta_{s,t}(\xi'))\bigr{\vert}\, dy \\
\le \, C\Vert F \Vert_H\Vert u \Vert_{L^{\infty}(C^{\alpha+\beta}_{b,d})}(s-t)^{\frac{\beta-2}{\alpha}},
\end{multline}
where in the last passage we used the smoothing effect for the frozen density $\tilde{p}^{\tau,\xi}$ (Equation
\eqref{eq:Smoothing_effects_of_tilde_p}).\newline
On the other side, the case $j>1$ can be exploited using the second Besov control (Lemma \ref{lemma:Second_Besov_COntrols}). For this reason,
we start using integration by parts formula to show that
\[I^{d}_j(s) \, = \, \Bigl{\vert}\int_{\R^{nd}}D_{ y_j}\cdot\Bigl{\{}D^2_{ x_1}\tilde{p}^{\tau,\xi'}(t,s,x,y)
\otimes\bigl[ F_j(s,y)- F_j(s,\theta_{s,t}(\xi'))\bigr]\Bigr{\}}u(s,y)\, dy \Bigr{\vert}.\]
Through Duality \eqref{Besov:duality_in_Besov} in Besov spaces and the identification in Equation \eqref{Besov:ident_Holder_Besov}, we
then write that
\begin{multline*}
I^{d}_j(s) \,\le \,C\Vert u \Vert_{L^\infty(C^{\alpha+\beta}_{b,d})}\int_{\R^{(n-1)d}} \Vert D_{ y_j}\cdot\Bigl{\{}D^2_{ x_1} \tilde{p}^{\tau,
\xi'}(t,s,x,y_{\smallsetminus j},\cdot)\\
\otimes\bigl[ F_j(s,y_{\smallsetminus j},\cdot) -  F_j(s,\theta_{s,t}(
\xi')) \bigr]\Bigr{\}}\Vert_{ B^{-(\alpha_j +\beta_j)}_{1,1}}\,dy_{\smallsetminus j}.
\end{multline*}
We can now apply the second Besov control (Lemma \ref{lemma:Second_Besov_COntrols}) to show that
\begin{equation}\label{Proof:H\"older_Perturbative_term_Idj}
I^d_j(s) \, \le \, C\Vert F\Vert_H\Vert u \Vert_{L^\infty(C^{\alpha+\beta}_{b,d})}(s-t)^{\frac{\beta-2}{\alpha}}.
\end{equation}
Going back at Equations \eqref{Proof:H\"older_Perturbative_term_Id0},\eqref{Proof:H\"older_Perturbative_term_Id1} and
\eqref{Proof:H\"older_Perturbative_term_Idj}, we then notice that
\begin{align}\label{Proof:H\"older_Perturbative_term_Id2}
\Bigl{\vert}\int_{t_0}^{T}D_{ x_1}\tilde{P}^{\tau,\xi'}_{s,t} R^{\tau,\xi'}(s,x) -D_{ x_1} \tilde{P}^{
\tau,\xi'}_{s,t}&R^{\tau,\xi'}(s,x') \, ds\Bigl{\vert} \\\notag
&\le \, C\Vert F\Vert_H\Vert u\Vert_{L^\infty(C^{\alpha+\beta}_{b,d})} \vert(x-x')_1\vert\int_{t_0}^{T}
(s-t)^{\frac{\beta-2}{\alpha}} ds \\\notag
&\le \, C\Vert F\Vert_H\Vert u \Vert_{L^\infty(C^{\alpha+\beta}_{b,d})} \vert(x-x')_1\vert
(t_0-t)^{\frac{\alpha+\beta-2}{\alpha}},
\end{align}
where in the last passage we used that $\frac{\alpha +\beta-2}{\alpha}<0$ to pick the starting point $t_0$ in the integral.\newline
Using that $t_0-t=c_0\mathbf{d}^{\alpha}(x,x')$, we can finally write that
\begin{multline*}
\Bigl{\vert}\int_{t_0}^{T}D_{ x_1}\tilde{P}^{\tau,\xi'}_{s,t}R^{\tau,\xi'}(s,x)
-D_{ x_1}\tilde{P}^{\tau,\xi'}_{s,t}R^{\tau,\xi'}(s,x') \, ds\Bigl{\vert}
\\
\le \, C c^{\frac{\alpha+\beta-2}{\alpha}}_0\Vert F\Vert_H\Vert u \Vert_{L^{\infty}(C^{\alpha+\beta}_{b,d})}
\mathbf{d}^{\alpha+\beta-1}(x,x').
\end{multline*}

\paragraph{Proof of Lemma \ref{lemma:Holder_modulus_Deg} in diagonal regime.}
We conclude this section showing the H\"older control in the degenerate setting when a diagonal regime is assumed. We start recalling that in Lemma \ref{lemma:Holder_modulus_Deg}, we assumed fixed a time $t$ in $[0,T]$ and two spatial points
$x,x'$ in $\R^{nd}$ such that $ x_j= x'_j$ if $j \neq i$ for some $i$ in $\llbracket 2 ,n \rrbracket$.

Representation \eqref{eq:Change_of_Freez_point_HOlder} and Duhamel-type Expansion \eqref{eq:Expansion_along_proxy} allows to control the
Holder modulus of a solution $u$ analyzing separately the
different terms:
\[
\begin{split}
u(t,x) - &u(t,x') \\
&= \, \Bigl(\tilde{P}^{\tau,\tilde{\xi}}_{T,t}u_T(x) - \tilde{P}^{\tau' ,\xi'}_{T,t} u_T(x')
\Bigr) + \Bigl(\tilde{G}^{\tau,\xi}_{t_0,t}f(t,x)+ \tilde{G}^{\tau,\tilde{\xi}}_{T,t_0}f(t,x) -
\tilde{G}^{\tau' ,\xi'}_{T,t} f(t,x')\Bigr) \\
&\qquad \,\, + \Bigl(\int_{t}^{t_0}  \tilde{P}^{\tau,\xi}_{s,t}R^{\tau,\xi}(s,x)\, ds +\int_{t_0}^{T}
\tilde{P}^{\tau,\tilde{\xi}}_{s,t}R^{\tau,\tilde{\xi}}(s,x)\, ds -
\int_{t}^{T}  \tilde{P}^{\tau',\xi'}_{s,t} R^{\tau' ,\xi'} (s,x') \, ds\Bigr) \\
&\qquad \qquad\quad+ \Bigl(\tilde{P}^{\tau,\xi}_{t_0,t}u(t_0,x)
-\tilde{P}^{\tau,\tilde{\xi}}_{t_0,t}u(t_0,x)\Bigr),
\end{split}\]
for some freezing couples $(\tau,\xi),(\tau,\tilde{\xi}),(\tau,\xi')$ fixed but to be chosen later. As done before, we assume
however from this point further that $\tau=\tau'$ and $\tilde{\xi}=\xi'$.

\emph{Control on the frozen semigroup.} Noticing that we have taken the same freezing couples since $\tilde{\xi}=\xi'$, the
control on the frozen semigroup $\bigl{\vert}\tilde{P}^{\tau,\xi'}_{T,t}u_T(x) -
\tilde{P}^{\tau,\xi'}_{T,t}u_T(x')\bigr{\vert}$ can be obtained exploiting the same argument used in the proof of the H\"older
control (Lemma \ref{lemma:Holder_modulus_proxy_Deg}) for the proxy.

\emph{Control on the Green kernel.} The proof of this estimate essentially matches the previous, analogous one in the
non-degenerate setting. Namely, we follow the proof in the global off-diagonal regime of Lemma \ref{lemma:Holder_modulus_Deg} to
control the local off-diagonal regime contribution $\bigl{\vert} \tilde{G}^{\tau,\xi}_{t_0,t}f(t,x)-\tilde{G}^{\tau,\xi'
}_{t_0,t} f(t,x') \bigl{\vert}$ while in the locally diagonal regime term
\[\bigl{\vert} \tilde{G}^{\tau,\xi'}_{T,t_0} f(t,x) - \tilde{G}^{\tau,\xi'}_{T,t_0}f(t,x')\bigr{\vert},\]
the freezing couples coincide and we can thus exploit the same argument used in the proof of the H\"older control (Lemma
\ref{lemma:Holder_modulus_proxy_Deg}) for the proxy.

\emph{Control on the discontinuity term.}
The proof of this result will follow essentially the one about the off-diagonal regime of the frozen semigroup with respect to the
degenerate variables. It holds that
\[
\begin{split}
\tilde{P}^{\tau,\xi}_{t_0,t}u(t_0,x) \, &= \, \int_{\R^{nd}}\tilde{p}^{\tau,\xi}(t,t_0,x,y)u(t_0,y)
\, dy \\
&= \, \int_{\R^{nd}}\frac{1}{\det \bigl(\mathbb{M}_{t_0-t}\bigr)}p_S\bigl(t_0-t,\mathbb{M}^{-1}_{t_0-t}(\tilde{m}^{\tau,\xi}_{t_0,t}(x)-y)  \bigr) u(t_0,y) \, dy \\
&= \, \int_{\R^{nd}}\frac{1}{\det \mathbb{M}_{t_0-t}}p_S(t_0-t,\mathbb{M}^{-1}_{t_0-t}z)u(t_0, \tilde{m}^{\tau,\xi}_{t_0,t}(x)-z) \,dz,
\end{split}\]
where in the last passage we used the change of variable $z = \tilde{m}^{\tau,\xi}_{t_0,t}(x)-y$.
Since a similar argument works also for $\tilde{P}^{\xi'}_{t_0,t}u(t_0,x)$, it then follows that
\begin{multline*}
\bigl{\vert}\tilde{P}^{\tau,\xi}_{t_0,t}u(t_0,x) - \tilde{P}^{\tau,\xi'}_{t_0,t}u(t_0,x)\bigr{\vert} \\
=\, \Bigl{\vert}\int_{\R^{nd}}\frac{1}{\det \mathbb{M}_{t_0-t}}p_S\bigl(t_0-t,\mathbb{M}^{-1}_{t_0-t}z\bigr)
\bigl[u(t_0,\tilde{m}^{\tau,\xi}_{t_0,t}(x)-z)-u(t_0,\tilde{m}^{\tau,\xi'}_{t_0,t}(x)-z)
\bigr]\,dz\Bigr{\vert}.
\end{multline*}
Remembering that $u(t_0,\cdot)$ is Lipschitz with respect to the first non-degenerate variable, we can write now that
\[\begin{split}
\bigl{\vert}\tilde{P}^{\tau,\xi}_{t_0,t}u(t_0,x) - &\tilde{P}^{\tau,\xi'}_{t_0,t}u(t_0,x)\bigr{\vert} \, \le \, C\Vert u \Vert_{L^\infty(C^{\alpha+\beta}_{b,d})}
\left(\int_{\R^{nd}}p_S\bigl(t_0-t,
\mathbb{M}^{-1}_{t_0-t} z\bigr)\,\frac{dz}{\det \mathbb{M}_{t_0-t}}\right)\\
&\qquad \qquad\qquad\qquad \times \bigl[\mathbf{d}^{\alpha+\beta}(\tilde{m}^{\tau,\xi}_{t_0,t}(x),\tilde{
m}^{\tau,\xi'}_{t_0,t}(x))+\bigl{\vert} (\tilde{m}^{\tau,\xi}_{t_0,t}(x)-\tilde{m}^{
\tau,\xi'}_{t_0,t}(x))_1\bigr{\vert}\bigr]\\
&\quad\le \,C\Vert u \Vert_{L^\infty(C^{\alpha+\beta}_{b,d})}\bigl[\mathbf{d}^{\alpha+\beta}(\tilde{m}^{\tau,\xi}_{t_0,t}(x),\tilde{
m}^{\tau,\xi'}_{t_0,t}(x))+\bigl{\vert} (\tilde{m}^{\tau,\xi}_{t_0,t}(x)-\tilde{m}^{
\tau,\xi'}_{t_0,t}(x))_1\bigr{\vert}\bigr].
\end{split}
\]
Taking $\xi=\xi'=x$, we can then use the sensitivity controls (Lemma \ref{lemma:Controls_on_means} and Equation
\eqref{eq:Sensitivity_on_mean1}) to show that
\[\bigl{\vert}\tilde{P}^{\tau,\xi}_{t_0,t}u(t_0,x) - \tilde{P}^{\tau,\xi'}_{t_0,t}u(t_0,x)\bigr{\vert} \, \le \,
C\Vert u\Vert_{L^\infty(C^{\alpha+\beta}_{b,d})}\Vert F \Vert_Hc_0^{\frac{\alpha+\beta}{1+\alpha(n-1)}}\mathbf{d}^{\alpha+\beta}
(x,x).\]

\emph{Control on the perturbative term.}
The proof of this estimate essentially matches the previous, analogous one in the non-degenerate setting. Namely, Inequalities
\eqref{Proof:H\"older_Perturbative_term_Id1}, \eqref{Proof:H\"older_Perturbative_term_Idj} and \eqref{Proof:H\"older_Perturbative_term_Id2}
hold again with $(s-t)^{\frac{\beta-2}{\alpha}}$ replaced by $(s-t)^{\frac{\beta}{\alpha}-\frac{1}{\alpha_i}}$.
\setcounter{equation}{0}

\subsubsection{Mollifying procedure}
We now make the mollifying parameter $m$ appear again using the notations introduced in Section $3.2$ (see
Equation \eqref{eq:Expansion_along_proxy}). Then, Lemmas \ref{lemma:supremum_norm}, \ref{lemma:Holder_modulus_Non-Deg} and
\ref{lemma:Holder_modulus_Deg} rewrite together in the following way. There exists a constant $C>0$ such that for any $m$ in $\N$,
\begin{multline}\label{eq:A_Priori_Estimates_Reg}
 \Vert u_m \Vert_{L^\infty(C^{\alpha+\beta}_{b,d})} \,
 \le \, Cc_0^{\frac{\beta-\gamma_n}{\alpha}}\bigl[\Vert u_{T,m}
\Vert_{C^{\alpha+\beta}_{b,d}}+\Vert f_m \Vert_{L^\infty(C^{\beta}_{b,d})}\bigr]\\
+C\bigl(c_0^{\frac{\beta-\gamma_n}{\alpha}}\Vert
F_m\Vert_H + c_0^{\frac{\alpha+\beta-1}{1+\alpha(n-1)}}\bigr)\Vert u_m \Vert_{L^\infty(C^{\alpha+\beta}_{b,d})},
\end{multline}
where $c_0$ is assumed to be fixed but chosen later. Importantly, $c_0$ and $C$ does not depends on the regularizing parameter $m$. Thus,
letting $m$ go to $\infty$ and remembering Definition \ref{definition:mild_sol} of a mild solution $u$, the above expression immediately
implies the a priori estimates (Proposition \ref{prop:A_Priori_Estimates}).

\section{Existence result}
\fancyhead[RO]{Section \thesection. Existence result}
The aim of this section is to show the well-posedness in a mild sense of the original IPDE \eqref{Degenerate_Stable_PDE}.
Recalling Definition \ref{definition:mild_sol} for a mild solution of IPDE \eqref{Degenerate_Stable_PDE}, let us consider three
sequences $\{f_m\}_{m\in \N}$, $\{u_{T,m}\}_{m\in \N}$
and $\{F_m\}_{m\in \N}$ of ``regularized'' coefficients such that
\begin{itemize}
  \item $\{f_m\}_{m\in \N}$ is in $C^\infty_b((0,T)\times\R^{nd})$ and $f_m$ converges to $f$ in
      $L^\infty\bigl(0,T;C^\beta_{b,d}(\R^{nd})\bigr)$;
  \item $\{u_{T,m}\}_{m\in \N}$ is in $C^\infty_b(\R^{nd})$ and $u_{T,m}$ converges to $u_T$ in $C^{\alpha+\beta}_{b,d}(\R^{nd})$;
   \item $\{F_m\}_{m\in \N}$ is in $C^\infty_b((0,T)\times\R^{nd};\R^{nd})$ and $\Vert F_m-F\Vert_H$ converges to $0$.
\end{itemize}
It can be derived through stochastic flows techniques (see e.g. \cite{book:Kunita04}) that there exists a solution $u_m$ in $C^\infty_b\bigl((0,T)\times\R^{nd}\bigr)$ of the ``regularized'' IPDE:
\[
\begin{cases}
   \partial_t u_m(t,x)+\mathcal{L}_\alpha u_m(t,x) + \langle A x + F_m(t,x), D_{x}
   u_m(t,x)\rangle  \, = \, -f_m(t,x) &\mbox{on}  \, (0,T)\times \R^{nd}; \\
    u_m(T,x) \, = \, u_{T,m}(x) & \mbox{on}\, \R^{nd}.
  \end{cases}
\]
In order to pass to the limit in $m$, we notice now that the arguments used above for the proof of the Schauder estimates (Equation
\eqref{equation:Schauder_Estimates}) can be applied to the above dynamics, too. Namely, there exists a constant $C>0$ such that
\[\Vert u_m \Vert_{L^\infty(C^{\alpha+\beta}_{b,d})} \, \le \, C\bigl[\Vert f_m\Vert_{L^\infty(C^{\beta}_{b,d})}+\Vert u_{T,m}
\Vert_{C^{\alpha+\beta}_{b,d}}\bigr] \, \le \,C\bigl[\Vert f\Vert_{L^\infty(C^{\beta}_{b,d})}+\Vert u_T
\Vert_{C^{\alpha+\beta}_{b,d}}\bigr].\]
Importantly, the above estimates is uniformly in $m$ and thus, the sequence $\{u_m\}_{m\in \N}$ is bounded in the space
$L^\infty\bigl(0,T;C^{\alpha+\beta}_{b,d}(\R^{nd})\bigr)$.
From Arzel\`a-Ascoli Theorem, we deduce now that there exists $u$ in $L^\infty\bigl(0,T;C^{\alpha+\beta}_{b,d}(\R^{nd})\bigr)$ and a sequence
$\{u_{m_k}\}_{k \in \N}$ of smooth and bounded functions converging to $u$ in $L^\infty\bigl(0,T;C^{\alpha+\beta}_{b,d}(\R^{nd})\bigr)$ and
such that $u_{m_k}$ is solution of the ``regularized'' IPDE \eqref{Regularizied_PDE}. It is then clear that $u$ is a mild solution of the
original IPDE \eqref{Degenerate_Stable_PDE}.

\paragraph{From mild to weak solutions}
We conclude showing that any mild solution $u$ of the IPDE \eqref{Degenerate_Stable_PDE} is indeed a weak solution. The proof
of this result will be essentially an application of the arguments presented before, especially the second Besov control (Lemma
\ref{lemma:Second_Besov_COntrols}).
Let $u$ be a mild solution of the IPDE \eqref{Degenerate_Stable_PDE} in $L^\infty\bigl(0,T;C^{\alpha+\beta}_{b,d}(\R^{nd})\bigr)$. Recalling
the definition of weak solution in \eqref{eq:Def_weak_sol}, we start fixing a test function $\phi$ in $C^\infty_0\bigl((0,T]\times
\R^{nd}\bigr)$ and passing to the ``regularized'' setting (see Definition \ref{definition:mild_sol}), we then notice that it holds that
\[\int_{0}^{T}\int_{\R^{nd}}\phi(t,y)\Bigl(\partial_t+L_t^m\Bigr)u_m(t,y) \, dy \, = \,
-\int_{0}^{T}\int_{\R^{nd}}\phi(t,y)f_m(t,y) \, dy,\]
where $L_t^m$ is the ``complete'' operator defined in \eqref{eq:Complete_Operator} but with respect to the regularized
coefficients. Integration by parts formula allows now to move the operators to the test function. Indeed, remembering that
$u_m(T,\cdot)=u_{T,m}(\cdot)$, it holds that
\begin{multline}
\label{TOweak:1}
\int_{0}^{T}\int_{\R^{nd}}\Bigl(-\partial_t+(L_t^m)^*\Bigr)\phi(t,y)u_m(t,y) \, dydt +
\int_{\R^{nd}}\phi(T,y)u_{T,m}(y) \, dy \\
= \, -\int_{0}^{T}\int_{\R^{nd}}\phi(t,y)f_m(t,y) \, dydt,
\end{multline}
where $\mathcal{L}^*_{m,\alpha}$ denotes the formal adjoint of $L_t^m$. We would like now to go back to the solution $u$,
letting $m$ go to $\infty$. We start rewriting the right-hand side term  in the following way:
\begin{multline*}
\int_{0}^{T}\int_{\R^{nd}}\phi(t,y)f_m(t,y) \, dydt \\
= \, \int_{0}^{T}\int_{\R^{nd}}\phi(t,y)f(t,y) \,
dydt + \int_{0}^{T}\int_{\R^{nd}}\phi(t,y)\bigl[f_m-f\bigr](t,y) \, dydt.
\end{multline*}
Exploiting that  $f_m$ converges to $f$ in $L^\infty\bigl(0,T;C^{\beta}_{b,d}(\R^{nd})\bigr)$ by assumption, it is easy to see that the second
contribution above goes to $0$ if we let $m$ go to $\infty$. A similar argument can be used to show that
\[\int_{\R^{nd}}\phi(T,y)u_{T,m}(y) \, dy \, \overset{m}{\to} \, \int_{\R^{nd}}\phi(T,y)u_T(y) \, dy.\]
On the other hand, we can decompose the first term in the left-hand side of Equation \eqref{TOweak:1} as
\begin{multline}\label{weak:decom}
\int_{0}^{T}\int_{\R^{nd}}\Bigl(-\partial_t+(L_t^m)^*\Bigr)\phi(t,y)u_m(t,y) \, dydt \\
= \,
\int_{0}^{T}\int_{\R^{nd}}\Bigl(-\partial_t+L_t^*\Bigr)\phi(t,y)u(t,y) \, dydt + R^1_m + R^2_m,
\end{multline}
where we have denoted
\begin{align*}
  R^1_m \, & = \, \int_{0}^{T}\int_{\R^{nd}}\Bigl[\mathcal{L}_{\alpha}^*-(L_t^m)^*\Bigr]\phi(t,y)u_m(t,y) \,
  dydt; \\
  R^2_m \, & = \, \int_{0}^{T}\int_{\R^{nd}}\Bigl(-\partial_t+\mathcal{L}_{\alpha}^*\Bigr)\phi(t,y)\bigl[u_m(t,y)-u(t,y)
  \bigr]\, dydt,
\end{align*}
with $\mathcal{L}_{\alpha}^*$ as the formal adjoint of the complete operator $L_t$. Noticing that
\[\Bigl[\mathcal{L}_{\alpha}^*-(L_t^m)^*\Bigr]\phi(t,y) \, =
\,D_{y}\cdot\bigl{\{}\phi(t,y)[F(t,y)-F_m(t,y)] \bigr{\}},\]
it is clear that the remainder contribution $R^1_m$ can be essentially handled as in the introduction of Section $5.1$, exploiting that
$\Vert F-F_m\Vert_H\to 0$.\newline
To control instead the second contribution $R^2_m$, we start decomposing it as
\[\begin{split}
R^2_m \, &= \, -\int_{0}^{T}\int_{\R^{nd}}\partial_t\phi(t,y)\bigl[u_m(t,y)-u(t,y)\bigr] \, dydt \\
&\qquad \qquad \qquad\qquad + \sum_{j=1}^{n}
\int_{0}^{T}\int_{\R^{nd}}D_{ y_j}\bigl[\phi  F_j\bigr](t,y)\bigl[u_m(t,y)-u(t,y)\bigr] \, dydt \\
&=: \, R^2_{0,m}+\sum_{j=1}^{n}R^2_{j,m}.
\end{split}\]
We firstly observe that $\vert R^2_{0,m}\vert$ goes to $0$ if we let $m$ go to $\infty$, since $\Vert
u-u_m\Vert_{L^\infty(C^{\alpha+\beta}_{b,d})
}\overset{m}{\to} 0$. On the other hand, integration by parts formula allows to show that
\[\vert R^2_{1,m} \vert \, = \, \Bigl{\vert}\int_{0}^{T}\int_{\R^{nd}}\bigl[\phi F\bigr](t,y)D_{ y_j}\bigl[u_m-u\bigr](t,y) \,
dydt\Bigr{\vert}\]
which again tends to $0$ when $m$ goes to $\infty$.
To control instead the contributions $R^m_{j,m}$ for $j>1$, the point is to use the Besov duality argument again. Namely, from Equations
\eqref{Besov:duality_in_Besov}, \eqref{Besov:ident_Holder_Besov} and with the notations in \eqref{eq:notation_smallsetminus}, it holds that
\[\begin{split}
\vert R^2_{j,m} \vert \, &\le \, \int_{0}^{T}\int_{\R^{d(n-1)}}\bigl{\Vert}D_{ y_j}\bigl[\phi F\bigr](t,y_{\smallsetminus
j},\cdot)\bigr{\Vert}_{B^{-(\alpha_j+\beta_j)}_{1,1}}\bigl{\Vert}\bigl[u_m-u\bigr](t,y_{\smallsetminus
j},\cdot)\bigr{\Vert}_{B^{\alpha_j+\beta_j}_{\infty,\infty}} \, dy_{\smallsetminus j}dt \\
&\le \, \int_{0}^{T}\int_{\R^{d(n-1)}}\bigl{\Vert}D_{ y_j}\bigl[\phi F\bigr](t,y_{\smallsetminus
j},\cdot)\bigr{\Vert}_{B^{-(\alpha_j+\beta_j)}_{1,1}}\bigl{\Vert}\bigl[u_m-u\bigr](t,y_{\smallsetminus
j},\cdot)\bigr{\Vert}_{C^{\alpha_j+\beta_j}_b} \, dy_{\smallsetminus j}dt.
\end{split}\]
Following the same arguments used in the proof of the second Besov control (Lemma \ref{lemma:Second_Besov_COntrols}), we know that there
exists a constant $C$ such that $\bigl{\Vert}D_{ y_j}\bigl[\phi F\bigr](t,y_{\smallsetminus
j},\cdot)\bigr{\Vert}_{B^{\alpha_j+\beta_j}_{1,1}} \le C\psi_j(t,y_{\smallsetminus j})$, where $\psi_j$ has compact support on
$\R^{d(n-1)}$. \newline
Since moreover $\Vert u_m -u \Vert$ goes to zero with $m$, we easily deduce that $R^2_{m,j}\overset{m}{\to}0$ for any $j$ in $\llbracket
2,n\rrbracket$. From the above controls, we can deduce now that $R^1_m + R^2_m\overset{m}{\to} 0$. From Equation \eqref{weak:decom}, it then
follows that
\[\int_{0}^{T}\int_{\R^{nd}}\Bigl(-\partial_t+(L_t^m)^*\Bigr)\phi(t,y)u_m(t,y) \, dydt \, \overset{m}{\to} \,
\int_{0}^{T}\int_{\R^{nd}}\Bigl(-\partial_t+L_t^*\Bigr)\phi(t,y)u(t,y) \, dydt\]
and the proof is concluded.

\setcounter{equation}{0}
\section{Extensions}
\fancyhead[RO]{Section \thesection. Extensions}
As already said in the introduction, our assumption of (global) H\"older regularity on the drift $\bar{F}$, as well as the choice of
considering a perturbed Ornstein-Uhlenbeck operator instead of a more general non-linear dynamics, was done to
preserve, as possible, the clarity and understandability of the article. In this conclusive section, we would like to explain briefly how it
possible to naturally extend it.
\subsection{General drift}
Here, we illustrate how the perturbative method explained above can be easily adapted to work in a
more general setting. In particular, the same results (well-posedness of the IPDE \eqref{Degenerate_Stable_PDE} and associated Schauder estimates)
can be proven to hold also for an equation of the form:
\begin{equation}
\label{Ext:Degenerate_Stable_PDE}
\begin{cases}
   \partial_t u(t,x)+ \mathcal{L}_\alpha u(t,x) + \langle\bar{F}(t,x), D_{x}u(t,x)\rangle \, =
   \, -f(t,x), & \mbox{on } (0,T)\times \R^{nd}; \\
    u(T,x) \, = \, u_T(x) & \mbox{on }\R^{nd},
  \end{cases}
\end{equation}
where $\bar{F}(t,x)=\bigl(\bar{F}_1(t,x),\dots,\bar{F}_n(t,x)\bigr)$ has the following structure
\[\bar{F}_i(t,x_{(i-1)\vee1},\dots,x_n).\]
We remark in particular that if for any $i$ in
$\llbracket 2,n\rrbracket$, $\bar{F}_i$ is linear with respect to $ x_{i-1}$ and independent from time, the previous analysis
works since we can rewrite $\bar{F}(t,x) =Ax+F(t,x)$.

In order to deal with this more general dynamics addressed in the diffusive setting in \cite{Chaudru:Honore:Menozzi18_Sharp}, we will need however to add some additional constraints and to modify slightly the ones
presented in assumption $[\textbf{A}]$. First of all, the non-degeneracy assumption $ [\textbf{H}]$ does not make sense in this new framework and
it will be replaced by the following one:
\begin{description}
  \item[{[H']}] the matrix $D_{ x_{i-1}}\bar{F}_i(t,x)$ has full rank $d$ for any $i$ in $\llbracket2,n
      \rrbracket$ and any $(t,x)$ in $[0,T]\times \R^{nd}$.
\end{description}
In particular, we will say that assumption $[\bar{\textbf{A}}]$ is in force when
\begin{description}
   \item[{[S']}] assumption [\textbf{ND}] and  [\textbf{H'}] are satisfied and the drift $\bar{F}=(\bar{F}_1,\dots,\bar{F}_n)$ is such that for any
       $i$ in $\llbracket 2,n\rrbracket$, $\bar{F}_i$ depends only on time and on the last $n-(i-2)\vee 0$ components, i.e.\
       $\bar{F}_i(t, x_{i-1},\dots,x_n)$;
  \item[{[P']}] $\alpha$ is a number in $(0,2)$, $\beta$ is in $(0,1)$ and it holds that
  \[\beta<\alpha, \quad \alpha+\beta\in (1,2) \, \, \text{ and }\,\, \beta < (\alpha-1)(1+\alpha(n-1));\]
\item[{[R']}] Recalling the notations in \eqref{eq:def_anisotropic_norm}-\eqref{eq:def_anisotropic_norm_in_time}, the source $f$ is in $L^\infty(0,T;C^{\beta}_{b,d}(\R^{nd}))$, the terminal condition $u_T$ is in
      $C^{\alpha+\beta}_{b,d}(\R^{nd})$ and for any $i$ in $\llbracket 1,n\rrbracket$, the drift $\bar{F}_i$ belongs to
      $L^\infty(0,T;C^{\gamma_i+\beta}_{d}(\R^{nd}))$ where $\gamma_i$ was defined in \eqref{Drift_assumptions}.
    \end{description}

To prove Schauder-type estimates for a solution of IPDE \eqref{Ext:Degenerate_Stable_PDE}, our idea is to adapt the perturbative
approach to this new dynamics. In particular, we can exploit the differentiability of $\bar{F}_i$ with respect to
$ x_{i-1}$ to ``linearize'' it along a flow that takes into account the perturbation (cf. Section $3.1$). Namely, we
are interested in the following equation:
\begin{multline}\label{Extension:Frozen_PDE}
    \partial_t \bar{u}^{\tau,\xi}(t,x)+\mathcal{L}_\alpha \bar{u}^{\tau,\xi}(t,x) + \bigl{\langle}\bar{A}^{\tau,\xi}_t
   \bigl(x-\bar{\theta}_{t,\tau}(\xi)\bigr) + \bar{F}(t,\bar{\theta}_{t,\tau}(\xi)), D_{x}
   \bar{u}^{\tau,\xi}(t,x)\bigr{\rangle} \\
   = \, -f(t,x);
\end{multline}
with initial condition $\bar{u}^{\tau,\xi}(T,x) = u_T(x)$, where the time-dependent matrix $\bar{A}^{\tau,\xi}_t$ is defined through
\[\bigl[\bar{A}^{\tau,\xi}_t\bigr]_{i,j} \, = \,
\begin{cases}
  D_{ x_{i-1}}\bar{F}_i(t,\theta_{t,\tau}(\xi)), & \mbox{if } j=i-1; \\
  0_{d\times d}, & \mbox{otherwise}
\end{cases}
\]
and $\bar{\theta}_{t,\tau}(\xi))$ is a fixed flow satisfying the dynamics
\begin{equation}
\bar{\theta}_{t,\tau}(\xi) \, = \, \xi +\int_{\tau}^{t} \bar{F}(v,\bar{\theta}_{v,\tau}(\xi)) \, dv.
\end{equation}
A first significant difference with respect to the previous approach consists in handling a time-dependent matrix
$\bar{A}^{\tau,\xi}_t$.
Indeed, it is possible to modify slightly the presentation in \cite{Priola:Zabczyk09} (allowing time-dependency on $A$) in order to
show that under assumption [\textbf{S'}], the two parameters semigroup $\{\bar{P}^{\tau,\xi}_{s,t}\}_{t\le s}$  associated with
the proxy operator
\[\mathcal{L}_\alpha + \langle\bar{A}^{\tau,\xi}_t \bigl(x-\bar{\theta}_{t,\tau}(\xi)\bigr) +
\bar{F}(t,\bar{\theta}_{t,\tau}(\xi)), D_{x} \rangle\]
admits a density $\bar{p}^{\tau,\xi}$ and that it can be written as
\[ \bar{p}^{\tau,\xi}(t,s,x,y)\,  = \, \frac{1}{\det \mathbb{M}_{s-t}} p_S\bigl( s-t,\mathbb{M}^{-1}_{s-t} (y - \bar{m}^{\tau,\xi}_{s,t}(x))\bigr).\]
Here, the notations for $p_S$ and $\mathbb{M}_{t}$ remain the same of above while this time the shift $\bar{m}^{\tau,\xi}_{s,t}$
is defined through
\[\bar{m}^{\tau,\xi}_{s,t}(x) \, = \, \mathcal{R}^{\tau,\xi}_{s,t}x + \int_{t}^{s}\mathcal{R}^{\tau,\xi}_{s,v}
\bigl[\bar{F}(v, \bar{\theta}_{v,\tau}(\xi))-\bar{A}^{\tau,\xi}_v\bar{\theta}_{v,\tau}(\xi)\bigr] \, dv,\]
where $\mathcal{R}^{\tau,\xi}_{s,t}$ is the time-ordered resolvent of $\bar{A}^{\tau,\xi}_s$ starting at time $t$, i.e.
\[\begin{cases}
    d\mathcal{R}^{\tau,\xi}_{s,t} \, = \,  \bar{A}^{\tau,\xi}_s\mathcal{R}^{\tau,\xi}_{s,t}ds, & \mbox{on } [t,T];\\
    \mathcal{R}^{\tau,\xi}_{t,t}\, = \, I.
  \end{cases}\]
We can as well refer to \cite{Huang:Menozzi16} for related issues (see Proposition $3.2$ and Section C about the
linearization, therein).

Following the same reasonings of Propositions \ref{prop:frozen_Duhamel_Formula} and \ref{prop:Expansion_along_proxy}, it is then possible to
state a Duhamel type formula suitable for IPDE \eqref{Ext:Degenerate_Stable_PDE}:
\begin{equation}\label{Ext:Expansion_along_proxy}
u(t,x) \, = \, \bar{P}^{\tau,\xi}_{T,t}u_T(x) + \int_{t}^{T}\bar{P}^{\tau,\xi}_{s,t} \bigl[f(s,\cdot) +\bar{R}^{\tau,
\xi}(s,\cdot)\bigr](x) \, ds,
\end{equation}
where the remainder term is given now by
\[\bar{R}^{\tau,\xi}(t,x) \, = \, \langle F(t,x)-F(t,\bar{\theta}_{t,\tau}(\xi))-
\bar{A}_t^{\tau,\xi}\bigl(x-\bar{\theta}_{t,\tau}(\xi)\bigr),D_{x}u(t,x)\rangle.\]

Looking back at the first part of the article, it is important to notice that the main steps of proof (cf. Equation
\eqref{eq:Smoothing_effects_of_tilde_p}, Propositions \ref{prop:Schauder_Estimates_for_proxy} and \ref{prop:A_Priori_Estimates} and Section
$3.3$) does not rely on the explicit formulas for $\bar{m}^{\tau,\xi}_{s,t}(x)$ and
$\bar{R}^{\tau,\xi}$ but instead, they
exploit only the Besov controls for the remainder $\bar{R}^{\tau,\xi}$ (cf. Section $5.1$) and the controls on the shift
$\bar{m}^{\tau,\xi}_{s,t}(x)$ (Section A.$2$).
Hence, once we have proven the suitable controls, the proofs of the analogous results for the new IPDE \eqref{Ext:Degenerate_Stable_PDE}
can be obtained easily modifying slightly the notations and following the same reasonings above.\newline
For example, exploiting that
\[\bar{m}^{\tau,\xi}_{s,t}(x) \, = \, x + \int_{t}^{s}\mathcal{R}^{\tau,\xi}_{v,t}\Bigl(\bar{m}^{\tau,\xi}_{v,t}(x) - \theta_{v,\tau}(\xi)\Bigr) + F(v, \theta_{v,\tau}(\xi)) \, dv,\]
we can follow the same method of proof in the above lemma \ref{lemma:identification_theta_m_stable} to show again that
\[\bar{m}^{\tau,\xi}_{s,t}(x) \, = \, \bar{\theta}_{s,\tau}(\xi),\]
taking $\tau=t$ and $\xi=x$.

Letting the interested reader look in the appendix for the suggestions on how to extend the controls on the shift
$\bar{m}^{\tau,\xi}_{s,t}(x)$ in this more general setting, we will focus now on proving the Besov controls.
First of all, we notice immediately that the proof of the first Besov control (Lemma \ref{lemma:First_Besov_COntrols}) relies essentially
only on the smoothing effect (Equation \eqref{eq:Smoothing_effects_of_tilde_p}) and thus, it can be obtained following the same reasoning above. The proof of
the second Besov control (Lemma \ref{lemma:Second_Besov_COntrols}) in this framework is a bit more involved and we are going to
explain it below more in details.\newline
We start noticing that Lemma \ref{lemma:Second_Besov_COntrols} can be reformulated for the new dynamics in the
following way:
\begin{multline*}
\int_{\R^{(n-1)d}}\Bigl{\Vert}D_{ y_j}\cdot\Bigl{\{}\mathbf{d}^{\vartheta}_{x}\bar{p}^{\tau,\xi}(t,s,x,
y_{\smallsetminus j},\cdot) \otimes \bar{\Delta}^{\tau,\xi}_j(s,y_{\smallsetminus j},\cdot) \Bigr{\}}
\Bigr{\Vert}_{B^{-(\alpha_j+\beta_j)}_{1,1}} \, dy_{\smallsetminus j}\\
\le \, C\Vert \bar{F}\Vert_H
(s-t)^{\frac{\beta}{\alpha} -\sum_{k=1}^{n}\frac{\vartheta_k}{\alpha_k}},
\end{multline*}
taking $(\tau,\xi)=(t,x)$, where we have denoted for simplicity
\[\bar{\Delta}^{\tau,\xi}_j(s,y) \, := \bar{F}_j(s,y)-\bar{F}_j(s,\theta_{s,\tau}(\xi))
-D_{x_{j-1}}\bar{F}_j(s,\theta_{s,\tau}(\xi))\bigl(y-\theta_{s,\tau}(\xi)\bigr)_{j-1},\]
for any $j$ in $\llbracket 2,n\rrbracket$.
The above control can be obtained mimicking the proof in the second Besov control (Lemma \ref{lemma:Second_Besov_COntrols}), exploiting this
time that
\[\vert \bar{\Delta}^{\tau,\xi}_j(s,y) \vert \le
C\Vert\bar{F}\Vert_H\mathbf{d}^{1+\alpha(j-2)+\beta}_{j-1:n}\bigl(y,\bar{\theta}_{s,\tau}(\xi)\bigr)\]
and the additional assumption [\textbf{P'}] in order to make the partial smoothing effect (Equation \eqref{eq:Partial_Smoothing_Effect})
work in this framework, too.\newline
The main difference in the proof is related to the control of the component $J_2(v,y_{\smallsetminus j},z)$ appearing in
Equation \eqref{eq:J_2_component_for_extension}. Namely,
\begin{multline*}
\int_{\R^d}D_z \partial_vp_h(v,z- y_j)\cdot\Bigl{\{}\bar{\Delta}^{\tau,\xi}_j(s,y_{\smallsetminus j},z)\\
\otimes
\int_{0}^{1}D_{ y_j}D^\vartheta_{x} \bar{p}^{\tau,\xi}(t,s,x,y_{\smallsetminus j},z+\lambda( y_j-z)) \cdot
( y_j-z)\Bigr{\}} \,d\lambda d y_j
\end{multline*}
with our new notations. Indeed, the dependence of $\bar{F}$ on $x_{j-1}$ pushes us to add a new term in the difference $\vert
\bar{F}_j(s,y_{\smallsetminus j}, z)-\bar{F}_j(s,\theta_{s,\tau}(\xi)) \vert$ (now,
$\vert\bar{\Delta}^{\tau,\xi}_j (s,y_{\smallsetminus j},z)\vert$) before splitting it up. In particular, it holds that
\[\begin{split}
\vert \bar{\Delta}^{\tau,\xi}_j (s,y_{\smallsetminus j},z) &\vert \, = \, \bigl{\vert} \bar{F}_j(s,y_{\smallsetminus j},z)-\bar{F}_j(s,\theta_{s,\tau}(\xi))-D_{x_{j-1}}
\bar{F}_j(s,\bar{\theta}_{s,\tau}(\xi))\bigl(y-\bar{\theta}_{s,\tau}(\xi)\bigr)_{j-1} \\
&  \qquad \qquad\quad \pm
\bar{F}_j(s,y_{1:j-1},\bigl(\bar{\theta}_{s,\tau} (\xi)\bigr)_{j:n}) \bigr{\vert} \\
&\le \, C\Vert \bar{F}\Vert_H\bigl(\vert z-\bigl(\bar{\theta}_{s,\tau}(\xi)\bigr)_j\vert^{\frac{1+\alpha(j-2)+\beta}{1
+\alpha(j-1)}}+\sum_{k=j+1}^{n}\vert \bigl(y-\bar{\theta}_{s,\tau}(\xi)\bigr)_k
\vert^{\frac{1+\alpha(j-2)+\beta}{1+\alpha(k-1)}}\\
& \qquad \qquad\quad +\vert \bigl(y-\bar{\theta}_{s,\tau}(\xi)\bigr)_{j-1}\vert^{\frac{1+\alpha (j-2)+\beta}{1+\alpha (j-2)}}\bigr) \\
&\le \, C\Vert \bar{F} \Vert_H\bigl(\vert \lambda(z- y_j) \vert^{\frac{1+\alpha(j-2)+\beta}{1 +\alpha(j-1)}}+\vert
z+\lambda( y_j-z)-\theta_{s,\tau}(\xi)_j\vert^{\frac{1+ \alpha(j-2) +\beta}{1 +\alpha(j-1)}}\\
&  \qquad \qquad\quad+\sum_{k=j+1}^{n}\vert
y-\bar{\theta}_{s,\tau}(\xi)_k\vert^{\frac{1+\alpha(j-2)+\beta}{1+\alpha(k-1)}}
+\vert \bigl(y - \bar{\theta}_{s,\tau}(\xi)
\bigr)_{j-1}\vert^{\frac{1+\alpha (j-2)+\beta}{1+\alpha (j-2)}}\bigr)\\
&\le \, C\Vert \bar{F} \Vert_H\bigl(\vert z- y_j
\vert^{\frac{1+\alpha(j-2)+\beta}{1 +\alpha(j-1)}} +
\mathbf{d}^{1+\alpha(j-2)+\beta}_{j+1:n}((y_{\smallsetminus j},z+\lambda( y_j-z)),\bar{\theta}_{s,\tau}(\xi))\bigr).
\end{split}\]
The remaining part of the proof exactly matches the original method in Lemma \ref{lemma:Second_Besov_COntrols}.

Even in this more general framework, it is thus possible to obtain the following:

\begin{theorem}[Well-posedness]
Under $[\bar{\rm\textbf{A}}]$, there exists a unique mild solution $u$ of IPDE \eqref{Ext:Degenerate_Stable_PDE}. Moreover, there exists a
constant $C:=C(T)$ such that
\[\Vert u \Vert_{L^\infty(C^{\alpha+\beta}_d)} \, \le \, C \bigl[\Vert f \Vert_{L^\infty(C^{\beta}_{b,d})} + \Vert u_T
\Vert_{C^{\alpha+\beta}_{b,d}}\bigr].\]
\end{theorem}

\subsection{Locally H\"older drift}
This part is designed to give a brief explanation on how it is possible to deal with the general IPDE \eqref{Ext:Degenerate_Stable_PDE} when
the drift $\bar{F}$ is only locally H\"older continuous in space. Namely, we assume with the notations in \eqref{Drift_assumptions} that
\begin{description}
  \item[{[LR']}] there exists a constant $K_0>0$ such that for any $i$ in $\llbracket1,n\rrbracket$
\[\mathbf{d}(\bar{F}(t,x),\bar{F}(t,x')) \,\le \, K_0\mathbf{d}^{\beta+\gamma_i}(x,x'), \quad t\, \in [0,T], \,
x,x'\in \R^{nd} \, \text{ s.t. } \mathbf{d}(x,x')<1.\]
\end{description}
In other words, it is required that $\bar{F}_i$ is in $L^\infty(0,T;C^{\beta+\gamma_i}(B(x_0,1/2)))$, uniformly in $x_0\in \R^{nd}$.

Under assumption $[\bar{\textbf{A}}]$ (with condition [\textbf{R'}] replaced by [\textbf{LR'}]), it is possible to recover the Schauder-type
estimates (Theorem \ref{theorem:Schauder_Estimates}), following the approach developed successfully in \cite{Chaudru:Menozzi:Priola19} for
the non-degenerate, super-critical stable setting.  Roughly speaking, in order to handle the local assumption, as well as the potentially unboundedness of the drift $\bar{F}$, we need to introduce a ``localized'' version of the Duhamel formulation (cf. Equation \eqref{eq:Expansion_along_proxy}). The key point here is
to multiply a solution $u$ by a suitable bump function $\bar{\eta}^{\tau,\xi}$ that ``localizes'' in space along the deterministic flow
$\bar{\theta}_{t,\tau}(\xi)$ that characterizes the proxy. Namely, we fix a smooth function $\rho$ that is equal to $1$ on $B(0,1/2)$
and vanishes outside $B(0,1)$) and then define for any $(\tau,\xi)$ in $[0,T]\times\R^{nd}$,
\[\bar{\eta}^{\tau,\xi}(t,x) \, := \, \rho(x-\bar{\theta}_{t,\tau}(\xi)).\]

We mention however that in the setting of \cite{Chaudru:Menozzi:Priola19}, the ``localization'' with the cut-off function
$\bar{\eta}^{\tau,\xi}$ is not simply motivated by the local H\"older continuity condition but it is also needed to give a proper
meaning to the Duhamel formulation for a solution (cf. Proposition \ref{prop:Expansion_along_proxy}) when $\alpha<1/2$, because of the low
integrability properties of the underlying stable density. Such a problem does not however appear here since condition [\textbf{P}] forces
us to consider only the case $\alpha>1/2$.

Given a mild solution $u$ of IPDE \eqref{Ext:Degenerate_Stable_PDE} and assuming $\bar{F}$ to be only locally H\"older
continuous as in [\textbf{LR'}], it is possible to show, at least formally, that the function
$\bar{v}^{\tau,\xi}:=u\bar{\eta}^{\tau,\xi}$ solves the following equation on $(0,T)\times \R^{nd}$:
\begin{equation}
\label{Localized_Degenerate_Stable_PDE}
\begin{cases}
 \partial_t \bar{v}^{\tau,\xi}(t,x) + \langle \bar{F}(t,x),D_{x}\bar{v}^{\tau,\xi}(t,x)\rangle +
 \mathcal{L}_\alpha \bar{v}^{\tau,\xi}(t,x) \, = \,
 -\bigl[\bar{\eta}^{\tau,\xi}f+\bar{\mathcal{S}}^{\tau,\xi}\bigr](t,x); \\
 \bar{v}^{\tau,\xi}(T,x) \, = \, \bar{\eta}^{\tau,\xi}(T,x)u_T(x),
  \end{cases}
\end{equation}
where we have denoted
\begin{multline*}
\bar{\mathcal{S}}^{\tau,\xi}(t,x) \, := \,\int_{\R^{d}}\bigl[u(t,x+By)-u(t,x) \bigr] \bigl[ \bar{\eta}^{
\tau,\xi}(t,t,x+By)- \bar{\eta}^{\tau,\xi}(t,x)\bigr] \, \nu_\alpha(dy)\\
-u(t,x)\langle\bar{F}(t,x)-\bar{F}(t,\bar{\theta}_{t,\tau}(\xi)), D\rho(x-\bar{\theta}_{t,\tau}(\xi))\rangle.
\end{multline*}
Equations as \eqref{Localized_Degenerate_Stable_PDE} can be essentially seen as a ``local'' version of the original one
\eqref{Ext:Degenerate_Stable_PDE}, depending on the freezing parameter $(\tau,\xi)$. In particular, it is important to notice that the
difference
\[\bar{F}(t,x)-\bar{F}(t,\bar{\theta}_{t,\tau}(\xi))\]
appearing in the ``localizing'' error $\bar{\mathcal{S}}^{\tau,\xi}$ can be controlled exactly because it is multiplied by the
derivative of the bump function $\rho$ in the right point $x-\bar{\theta}_{t,\tau}(\xi)$, allowing us to exploit the
\emph{local} H\"older regularity. On the other hand, the first integral term in the r.h.s. can be seen as a commutator which involves only
the non-degenerate variables and thus, that can be handled with interpolation techniques as in \cite{Chaudru:Menozzi:Priola19}.

Even with the additional difficulty in controlling the remainder term, the perturbative approach explained
in Section $3$ can be applied, leading to show Schauder-type estimates as in Theorem \ref{theorem:Schauder_Estimates} and the
well-posedness of the IPDE \eqref{Ext:Degenerate_Stable_PDE} when assuming $\bar{F}$ to be only locally H\"older continuous.

Our procedure could be also used in order to establish Schauder-type estimates for the full Ornstein-Uhlenbeck operator as done, for
example, in \cite{Lunardi97} for the diffusive case. Indeed, a general operator of the form $\langle
\bar{A}x,D_{x}\rangle+\mathcal{L}_\alpha$ can be treated, decomposing the matrix as $\bar{A}=A+U$ where $A$ is, as before, the sub-diagonal
matrix that makes the Ornstein-Ulhenbeck operator invariant by the dilation operator associated with the distance $d$, while $U$ is an upper
triangular matrix that could be seen as an additional \emph{locally} H\"older term.

\subsection{Diffusion coefficient}
We conclude the article showing briefly how an additional diffusion term $\sigma\colon [0,T]\times\R^{nd}\to
\R^d\otimes \R^d$ can be handled in the IPDE \eqref{Ext:Degenerate_Stable_PDE} with an operator $\mathcal{L}_{\alpha,t}$ of the form:
\[\mathcal{L}_{\alpha,t}\phi(t,x)\, := \, \text{p.v.}\int_{\R^d}\bigl[\phi(t,x+B\sigma
(t,x)y) -\phi(t,x)\bigr]\nu_\alpha(dy).\]
In this framework, it is quite standard (cf. \cite{Hao:Wu:Zhang19} and \cite{Zhang:Zhao18}) to assume the L\'evy measure
$\nu_\alpha$ to be absolutely continuous with respect to the Lebesgue measure on $\R^d$ i.e.\ $\nu_{\alpha}(dy)=f(y)dy$, for some Lipschitz
function $f\colon \R^d \to \R$. In particular, since $\nu_\alpha$ is a symmetric, $\alpha$-stable, L\'evy measure, it holds passing to polar
coordinates $y=\rho s$ where $(\rho,s)\in [0,\infty)\times \mathbb{S}^{d-1}$ that
\[f(y) \, = \, \frac{g(s)}{\rho^{d+\alpha}}\]
for an even, Lipschitz function $g$ on $\mathbb{S}^{d-1}$ (see also Equation \eqref{eq:decomposition_measure}). Moreover, $\sigma$ is
considered uniformly elliptic and in
$L^\infty(0,T;C^{\beta}(\R^n,\R))$. \newline
Introducing now the ``frozen'' operator \[\bar{\mathcal{L}}^{\tau,\xi}_{\alpha,t}\phi(t,x) \, = \, \text{p.v.}\int_{\R^d}\bigl[\phi(t,x+B\sigma
(t,\bar{\theta}_{t,\tau}(\xi))y) -\phi(t,x)\bigr]\nu_\alpha(dy),\]
this would lead to consider for the IPDE an
additional term in the Duhamel formula (cf. Equation \eqref{Ext:Expansion_along_proxy}) that would write:
\begin{equation}\label{Duhamel_formula_for_additional}
u(t,x) \, = \, \breve{P}^{\tau,\xi}_{T,t}u_T(x) +
\int_{t}^{T}\breve{P}^{\tau,\xi}_{s,t}f(s,x)+\breve{P}^{\tau,\xi}_{s,t}\bar{R}^{\tau,\xi}(s,x)+
\breve{P}^{\tau,\xi}_{s,t}\bigl[\bigl(\mathcal{L}_{\alpha,t}-\bar{\mathcal{L}}^{\tau,\xi}_{\alpha,t}\bigr)u(s,\cdot)\bigr](x) \, ds.
\end{equation}
Here, $\{\breve{P}^{\tau,\xi}_{s,t}\}_{t\le s}$ denotes the two parameter semigroup associated with the proxy operator
\[\bar{L}^{\tau,\xi}_\alpha+\langle\bar{A}^{\tau,\xi}_t \bigl(x-\bar{\theta}_{t,\tau}(\xi)\bigr) +
\bar{F}(t,\bar{\theta}_{t,\tau}(\xi)), D_{x} \rangle.\]
Let us focus on the last term in the integral of Equation \eqref{Duhamel_formula_for_additional}. Looking back at the proof of the a priori
estimates (Proposition \ref{prop:A_Priori_Estimates}), we notice in
particular that we aim to establish the following control:
\begin{equation}\label{zzz1}
\bigl{\vert} \bigl(\mathcal{L}_{\alpha,t}-\bar{\mathcal{L}}^{\tau,\xi}_{\alpha,t}\bigr)u(t,x)\bigr{\vert} \, \le \, C\Vert \sigma
\Vert_{L^\infty(C^{\beta}_{b,d})}\Vert u
\Vert_{L^\infty(C^{\alpha+\beta}_{b,d})}\mathbf{d}^{\beta}(x,\bar{\theta}_{t,\tau}(\xi))
\end{equation}
in order to apply the same reasoning above in this new framework. To this end, we write that
\[\begin{split}
\bigl(\mathcal{L}_{\alpha,t}-\bar{\mathcal{L}}^{\tau,\xi}_{\alpha,t}\bigr)u(t,x) \, &= \, \text{p.v.}\int_{\R^{d}}\bigl{\{}u(t,x+B\sigma(t,x)y)-u(t,x)\bigr{\}} \, \nu_\alpha(dy) \\
&\qquad\qquad\qquad -
\text{p.v.}\int_{\R^{d}}\bigl{\{}u(t,x+B\sigma(t,\bar{\theta}_{t,\tau}(\xi))y)-u(t,x)\bigr{\}} \, \nu_\alpha(dy) \\
&= \, \text{p.v.}\int_{\R^{d}}\bigl{\{}u(t,x+Bz)-u(t,x)\bigr{\}}\frac{f\bigl(\sigma^{-1}(t,x)z\bigr)}{\det \sigma(t,x)} \, dz \\
&\qquad  \qquad\qquad-
\int_{\R^{d}}\bigl{\{}u(t,x+Bz)-u(t,x)\bigr{\}}\frac{f\bigl(\sigma^{-1}(t,\bar{\theta}_{t,\tau}(\xi))z\bigr)}{\det
\sigma(t,\bar{\theta}_{t,\tau}(\xi))} \, dz \\
&= \,\text{p.v.}\int_{0}^\infty \frac{1}{\rho^{1+\alpha}}\int_{\mathbb{S}^{d-1}}\bigl{\{}u(t,x+B\rho
s)-u(t,x)\bigr{\}}\bar{D}^{\tau,\xi}(t,x,s)
\, dsd\rho
\end{split}\]
where we have denoted, for notational convenience
\[\bar{D}^{\tau,\xi}(t,x,s) \, := \, \Bigl{\{}\frac{g\bigl(\frac{\sigma^{-1}(t,x)s}{\vert
\sigma^{-1}(t,x)s\vert}\bigr)}{\vert
\sigma^{-1}(t,x)s\vert^{d+\alpha}\det \sigma(t,x)} -\frac{g\bigl(\frac{\sigma^{-1}(t,\bar{\theta}_{t,\tau}
(\xi))s}{\vert\sigma^{-1}(t,\bar{\theta}_{t,\tau} (\xi))s\vert}\bigr)}{\vert \sigma^{-1}(t,\bar{\theta}_{t,\tau}
(\xi))s\vert^{d+\alpha}\det\sigma(t,\bar{\theta}_{t,\tau}(\xi))}\Bigr{\}}.\]
Using now that $g$ is Lipschitz and the assumptions on $\sigma$, we can show that
\begin{equation}\label{eq:COntrol_on_D}
  \vert \bar{D}^{\tau,\xi}(t,x,s) \vert \, \le \,  C\vert \sigma(t,x) - \sigma(t,\bar{\theta}_{t,\tau} (\xi))\vert
\, \le \, C\Vert \sigma \Vert_{L^\infty(C^{\beta}_{b,d})} \mathbf{d}^{\beta}(x,\bar{\theta}_{t,\tau}(\xi)).
\end{equation}
Finally, Equation \eqref{zzz1} follows from the previous controls using Taylor expansions and the symmetry condition on $\nu_\alpha$.
Namely, considering the case $\alpha\ge 1$, which is the most delicate one for this part and precisely requires the symmetry of $u_T$, we write that
\begin{align}\notag
\Bigl{\vert}\bigl(\mathcal{L}_{\alpha,t}-\bar{\mathcal{L}}^{\tau,\xi}_{\alpha,t}\bigr)u(t,x)\Bigr{\vert} \, &= \, \Bigl{|}\text{p.v.}\int_{0}^\infty
\frac{1}{\rho^{1+\alpha}}\int_{\mathbb{S}^{d-1}}\bigl{\{}u(t,x+B\rho s)-u(t,x)\bigr{\}}\bar{D}^{\tau,\xi}(t,x,s)\,
dsd\rho\Bigr{|} \\\notag
&\le \, \Bigl{|}\text{p.v.}\int_{(0,1)}
\frac{1}{\rho^{1+\alpha}}\int_{\mathbb{S}^{d-1}}\bigl{\{}u(t,x+B\rho s)-u(t,x)\bigr{\}}\bar{D}^{\tau,\xi}(t,x,s)\,
dsd\rho \Bigl{|}\\\notag
&\qquad+ \int_{(1,\infty)} \frac{1}{\rho^{1+\alpha}}\int_{\mathbb{S}^{d-1}}\bigl{|}u(t,x+B\rho
s)-u(t,x)\bigr{|}\,|\bar{D}^{\tau,\xi}(t,x,s)|\, dsd\rho \\
&=: \bigl[\bar{I}^{\tau,\xi}_s + \bar{I}^{\tau,\xi}_l\bigr](t,x).\label{eq:Control_additional}
\end{align}
The \emph{large jump} contribution $\bar{I}^{\tau,\xi}_l$ is easily handled from Equation \eqref{eq:COntrol_on_D}. We get that
\begin{equation}\label{eq:Control_additional1}
\begin{split}
\bar{I}^{\tau,\xi}_l(t,x) \,&\le \, 2C\Vert \sigma \Vert_{L^\infty(C^{\beta}_{b,d})}\Vert u \Vert_{L^\infty(L^\infty)}
\mathbf{d}^{\beta}(x,\bar{\theta}_{t,\tau}(\xi)) \\
&\le \, 2C\Vert \sigma \Vert_{L^\infty(C^{\beta}_{b,d})}\Vert u
\Vert_{L^\infty(C^{\alpha+\beta}_{b,d})} \mathbf{d}^{\beta}(x,\bar{\theta}_{t,\tau}(\xi)).
\end{split}
\end{equation}
On the other hand, from the symmetry assumption on $\nu_\alpha$, which transfers to $u_T$, we can control the \emph{small jump} contribution
$\bar{I}^{\tau,\xi}_s$ through Taylor expansion and a centering
argument. Indeed,
\begin{align}\notag
&\bar{I}^{\tau,\xi}_s(t,x) \\
\notag
&= \, \Bigl{|}\text{p.v.}\int_{(0,1)}
\frac{1}{\rho^{1+\alpha}}\int_{\mathbb{S}^{d-1}}\int_{0}^{1}\bigl[D_{ x_1}u(t,x+\lambda B\rho s)-D_{ x_1}u(t,x)\bigr]\rho
s\bar{D}^{\tau,\xi}(t,x,s)\, d\lambda dsd\rho  \Bigl{|}\\\notag
&\le \, C\Vert \sigma \Vert_{L^\infty(C^{\beta}_{b,d})} \mathbf{d}^{\beta}(x,\bar{\theta}_{t,\tau}(\xi))\int_{(0,1)}
\frac{1}{\rho^\alpha}\int_{\mathbb{S}^{d-1}}\int_{0}^{1}\bigl{\vert}D_{ x_1}u(t,x+\lambda B\rho
s)-D_{ x_1}u(t,x)\bigr{\vert} \, d\lambda dsd\rho \\\notag
&\le \, C\Vert \sigma \Vert_{L^\infty(C^{\beta}_{b,d})} \Vert D_{ x_1}u
\Vert_{L^\infty(C^{\alpha+\beta-1}_{b,d})}\mathbf{d}^{\beta}(x,\bar{\theta}_{t,\tau}(\xi))\int_{(0,1)}
\frac{1}{\rho^{\alpha}} \rho^{\alpha+\beta -1} \, d\rho \\
&\le \, C\Vert \sigma \Vert_{L^\infty(C^{\beta}_{b,d})} \Vert u
\Vert_{L^\infty(C^{\alpha+\beta}_{b,d})}\mathbf{d}^{\beta}(x,\bar{\theta}_{t,\tau}(\xi)).
\label{eq:Control_additional2}
\end{align}
Using Controls \eqref{eq:Control_additional1} and \eqref{eq:Control_additional2} in Equation \eqref{eq:Control_additional}, we
obtain the expected bound (Equation \eqref{zzz1}). We remark that the case $\alpha<1$ could be handled similarly for
the contribution $\bar{I}^{\tau,\xi}_l$ and even more directly for $\bar{I}^{\tau,\xi}_s$. Indeed, in that case, the centering argument is not needed since the Taylor
expansion already yields an integrable singularity.

\setcounter{equation}{0}
\section{Appendix: proofs of complementary results}
\fancyhead[RO]{Section \thesection. Appendix}

\subsection{Smoothing effects for Ornstein-Ulhenbeck operator}

We state and prove here some of the key properties of the Ornstein-Uhlenbeck operator. Namely, we will prove the representation
\eqref{eq:Representation_of_p_ou} and the associated $\alpha$-smoothing effect \eqref{Smoothing_effect_of_S}. We highlight however that
these results are only a slight modification to our purpose of those in  \cite{Huang:Menozzi:Priola19}.

The two lemma below presents a deep connection with stochastic analysis and their proofs relies on tools that are more familiar in the
probabilistic realm. For this reason, we are going to consider the stochastic counterpart of the Ornstein-Ulhenbeck operator $L^{\text{ou}}$.
Namely, for a given starting point $x$ in $\R^{nd}$, we are interested in the following dynamics
\begin{equation}\label{eq:Stochastic_OU_equation}
\begin{cases}
dX_t \, = \, AX_tdt+BdZ_t, & \mbox{on } [0,T] \\
X_0 \, = \, x
\end{cases}
\end{equation}
where $(Z_t)_{t\ge 0}$ is an $\alpha$-stable, $\R^{nd}$-dimensional process with L\'evy measure $\nu_\alpha$, defined on some complete
probability space $(\Omega,\mathcal{F},\mathbb{P})$.

\begin{lemma}[Representation]
Under [\textbf{A}], the semigroup $\{P^{\text{ou}}_t\}_{t>0}$ generated by the Ornstein-Ulehnbeck operator $L^{\text{ou}}$ (defined in \eqref{eq:def_of_OU_operator})  admits for any fixed $t>0$, a density $p^{\text{ou}}(t,\cdot)$ which writes for any $t>0$ and any $x,y$
in $\R^{nd}$
\[p^{\text{ou}}(t,x,y) \, = \, \frac{1}{\det \mathbb{M}_t}p_S(t,\mathbb{M}^{-1}_t\bigl(e^{At}x-y)\bigr)\]
where $\mathbb{M}_t$ is the matrix defined in \eqref{eq:def_of_Mt} and $p_S$ is the smooth density of an $\R^{nd}$-valued, symmetric and
$\alpha$-stable process $S$ whose L\'evy measure $\mu_S$ satisfies the non-degeneracy assumption [\textbf{ND}] on $\R^{nd}$.
\end{lemma}
\begin{proof}
We start noticing that the above dynamics \eqref{eq:Stochastic_OU_equation} can be explicitly integrated and gives
\[X_t \, = \, e^{tA}x+\int_{0}^{t}e^{(t-s)A}B \, dZ_s.\]
It is then readily derived from \cite{Priola:Zabczyk09} that, for any $t>0$, the random variable $X_t$ has a density $p_X(t,x,\cdot)$
with respect to the Lebesgue measure on $\R^{nd}$ and it is moreover well known (see for example \cite{book:Dynkin65}) that $p_X$ coincides
with the density $p^{\text{ou}}$ of the Ornstein-Ulhenbeck operator $L^{\text{ou}}$ .\newline
For this reason, we fix $t\ge 0$ and consider, for a given $N$ in $\N$, a uniform partition $\{t_i\}_{i\in\llbracket 0,N\rrbracket}$ of
$[0,t]$. Then, it holds for any $p$ in $\R^{nd}$,
\begin{multline*}
\mathbb{E}\Bigl[\text{exp}\Bigl(i\langle p, \sum_{i=1}^{N}e^{(t-t_{i-1})A}B\bigl(Z_{t_i}-Z_{t_{i-1}}\bigr)\Bigr)\rangle\Bigr] \\
= \,
\text{exp}\Bigl(-\frac{1}{N}\sum_{i=1}^{N}\int_{\mathbb{S}^{d-1}}\vert \langle B^\ast e^{(t-t_{i-1})A^\ast}p,s\rangle\vert^\alpha \,
\mu(ds)\Bigr)
\end{multline*}
where $\mu$ is the spherical measure associated with $\nu_\alpha$ (see Equation \eqref{def:Levy_Symbol_stable}). By dominated convergence
theorem, we let $m$ goes to infinity and show that
\[\mathbb{E}\Bigl[\text{exp}\Bigl(i\langle p,\int_{0}^{t}e^{(t-s)A}B \, dZ_s\Bigr)\Bigr] \, = \,
\text{exp}\Bigl(-\int_{0}^{t}\int_{\mathbb{S}^{d-1}}\vert \langle e^{uA^\ast}p,Bs\rangle \, \mu(ds)du\Bigr).\]
Thanks to the above equation, we can rewrite the characteristic function of $X_t$ as:
\[\begin{split}
\psi_{X_t}(p) \, &= \,\mathbb{E}\Bigl[\text{exp}\Bigl(i\langle p,e^{tA}x+\int_{0}^{t}e^{(t-s)A}B \, dZ_s\Bigr)\Bigr] \\
&= \, \text{exp}\Bigl(i\langle p,e^{tA}x\rangle -\int_{0}^{t}\int_{\mathbb{S}^{d-1}}\vert \langle e^{uA^\ast} p, Bs \rangle
\vert^\alpha \, \mu(ds)du \Bigr)\\
&= \, \text{exp}\Bigl(i\langle p,e^{tA}x\rangle -t\int_{0}^{1}\int_{\mathbb{S}^{d-1}}\vert \langle e^{vtA^\ast} p, Bs \rangle
\vert^\alpha \, \mu(ds)dv \Bigr)
\end{split}\]
where in the last passage we used the change of variables $u=vt$. For the next step, we firstly notice that it holds
\[e^{tA} \, = \, \mathbb{M}_te^A\mathbb{M}^{-1}_t,\]
shown using the definition of matrix exponential and the trivial relation $\mathbb{M}_tA\mathbb{M}^{-1}_t = tA$. Exploiting the above
identity, we then find that
\begin{multline*}
  \psi_{X_t}(p) \, = \, \text{exp}\Bigl(i\langle p,e^{tA}x\rangle -t\int_{0}^{1}\int_{\mathbb{S}^{d-1}}\vert \langle
  \mathbb{M}_tp,e^{vA}\mathbb{M}_t^{-1} Bs\rangle\vert^\alpha\, \mu(ds)dv \Bigr)\\
  = \, \text{exp}\Bigl(i\langle p,e^{tA}x\rangle -t\int_{0}^{1}\int_{\mathbb{S}^{d-1}}\vert \langle
  \mathbb{M}_tp,e^{vA}Bs\rangle\vert^\alpha\, \mu(ds)dv \Bigr)
\end{multline*}
where in the last passage we used the straightforward identity $\mathbb{M}_t^{1}By=By$.
We focus now only on the double integral
\[\int_{0}^{1}\int_{\mathbb{S}^{d-1}}\vert \langle \mathbb{M}_tp,e^{vA}Bs\rangle\vert^\alpha\, \mu(ds)dv.\]
If we consider the measure $m_\alpha(dv,ds) := \vert e^{vA}Bs \vert^\alpha \mu(ds)dv$ on $[0,1]\times \mathbb{S}^{d-1}$ and the
normalized lift function $l\colon [0,1]\times \mathbb{S}^{d-1}\to \mathbb{S}^{nd-1}$ given by
\[l(v,s) \, := \, \frac{e^{vA}Bs}{\vert e^{vA}Bs \vert},\]
it then follows that
\[
\begin{split}
\int_{0}^{1}\int_{\mathbb{S}^{d-1}}\vert \langle \mathbb{M}_tp,e^{vA}Bs\rangle\vert^\alpha\, \mu(ds)dv \, &= \,
\int_{0}^{1}\int_{\mathbb{S}^{d-1}}\vert \langle \mathbb{M}_tp,\frac{e^{vA}Bs}{\vert e^{vA}Bs\vert}\rangle\vert^\alpha\,
m_\alpha(ds,dv)\\
&= \, \int_{\mathbb{S}^{nd-1}}\vert \langle \mathbb{M}_tp,\xi\rangle\vert^\alpha\, \mu_S(d\xi),
\end{split}\]
where $\mu_S\, := \, \text{Sym}(l_\ast(m_\alpha))$ is the symmetrized version of the measure $m_\alpha$ push-forwarded through $l$.\newline
Noticing that $\mu_S$ is the L\'evy measure of a symmetric $\alpha$-stable process $\{S_t\}_{t\ge 0}$ satisfying assumption [\textbf{ND}] on $\R^{nd}$,
we can finally write that
\[\psi_{X_t}(p) \, = \, \text{exp}\Bigl(i\langle p,e^{tA}x\rangle -t\Phi_S(\mathbb{M}_tp) \Bigr)\]
where $\Phi_S$ is the L\'evy symbol associated with $S_t$ (cf. Equation \eqref{def:Levy_Symbol_stable}).\newline
From Lemma A.$1$ in \cite{Huang:Menozzi:Priola19}, we know that under assumption [\textbf{ND}], the above calculations implies that
\[\int_{0}^{1}\int_{\mathbb{S}^{d-1}}\bigl{\vert}(\mathbb{M}_tp) \cdot (e^{Av}Bs)\bigr{\vert}^\alpha \, \mu_S(ds)dv \, \ge \, C\vert
\mathbb{M}_tp\vert^\alpha\]
for some constant $C> 0$. It follows in particular that the function $p\to \psi_{X_t}(p)$ is in $L^1(\R^{nd})$. Thus, by inverse
fourier transform and a change of variables, we can prove that
\[ \begin{split}
\mathcal{F}^{-1}\bigl[\psi_{X_t}\bigr](y) \, &= \, \frac{1}{(2\pi)^{nd}}\int_{\R^{nd}} e^{-i\langle p,
y\rangle}\text{exp}\Bigl(i\langle p,e^{tA}x\rangle -t\Phi_S(\mathbb{M}_tp) \Bigr)\, dp \\
&=\, \frac{\det(\mathbb{M}^{-1}_t)}{(2\pi)^{nd}}\int_{\R^{nd}} \text{exp}\Bigl(-i\bigl{\langle}\mathbb{M}^{-1}_tp, y-e^{tA}x
\bigr{\rangle}\Bigr)e^{-t\Phi(p)}\, dp \\
&= \, \frac{\det (\mathbb{M}^{-1}_t)}{(2\pi)^{nd}}\int_{\R^{nd}} \text{exp}\Bigl(-i\bigl{\langle}p,
\mathbb{M}^{-1}_t\bigl(y-e^{tA}x\bigr)\bigr{\rangle}\Bigr) e^{-t\Phi(p)}\, dp \\
&= \,
\frac{1}{\det\mathbb{M}_t}p_S(t,\mathbb{M}^{-1}(y-e^{At}x))
\end{split}\]
and we have concluded since $p_S$ is symmetric.
\end{proof}

We can now point out the smoothing effect (Equation \eqref{Smoothing_effect_of_S}) associated with the Ornstein-Uhlenbeck density $p^{\text{ou}}$.

\begin{lemma}[Smoothing effect]
\label{Appendix:Smoothing_effect}
Under $[\textbf{A}]$, there exists a family $\{q(t,\cdot)\colon t \in [0,T]\}$ of densities on $\R^{nd}$ such that
\begin{itemize}
  \item  for any $l$ in $\llbracket 0,3 \rrbracket$, there exists a constant $C:=C(l,nd)$ such that $\vert D^l_{y}p_S(t,y)\vert \, \le \,
      Cq(t,y)t^{-l/\alpha}$ for any $t$ in $[0,T]$ and any $y$ in $\R^{nd}$;
  \item (stable scaling property) $q(t,y)=t^{-nd/\alpha}q(1,t^{-1/\alpha}y)$ for any $t$ in $[0,T]$ and any $y$ in
      $\R^{nd}$;
  \item (stable smoothing effect) for any $\gamma$ in $[0,\alpha)$, there exists a constant $c:=c(\gamma,nd)$ such that
  \begin{equation}\label{equation:integration_prop_of_q_stable}
    \int_{\R^{nd}} q(t,y) \, \vert y \vert^\gamma \, dy \, \le \, ct^{\gamma/\alpha} \,\, \text{for any }t>0.
  \end{equation}
\end{itemize}
\end{lemma}
\begin{proof}
Fixed a time $t>0$, we start applying the Ito-L\'evy decomposition to $S$ at the associated characteristic stable time
scale, i.e. we choose to truncate at threshold $t^{1/\alpha}$, so that we can write $S_t = M_t+N_t$ for some $M_t,N_t$ independent random
variables corresponding to the small jumps part and the large jumps part, respectively. Namely, we denote for any $s>0$
\[N_s \, := \, \int_{0}^{s}\int_{\vert x \vert >t^{1/\alpha}}xP(du,dx) \,\, \text{ and } \,\, M_s: \,= \, S_s-N_s
\]
where $P$ is the Poisson random measure associated with the process $S$. We can thus rewrite the density $p_S$ in the following way
\[p_S(t,x) \, = \, \int_{\R^{nd}}p_M(t,x-y)P_{N_t}(dy)\]
where $p_M(t,\cdot)$ corresponds to the density of $M_t$ and $P_{N_t}$ is the law of $N_t$.\newline
It is important now to notice that it is precisely our choice of the cutting threshold $t^{1/\alpha}$ that gives $M$ and $N$ the
$\alpha$-similarity property (for any fixed $t$)
\[N_t \, \overset{law}{=} \, t^{1/\alpha}N_1 \, \, \text{ and } M_t \, \overset{law}{=} \, t^{1/\alpha}M_1\]
we will need below. Indeed, to show the assertion for $N$, we can start from the L\'evy-Khintchine formula for the characteristic function of
$N$:
\[\mathbb{E}\bigl[e^{i \langle p,N_t\rangle}\bigr] \, = \, \exp\Bigl[t\int_{\mathbb{S}^{nd-1}}\int_{t^{1/\alpha}}^{\infty}\bigl(\cos(\langle
p,r\xi \rangle)-1\bigr)\frac{dr}{r^{1+\alpha}}\overline{\mu}_S(d\xi)\Bigr]\]
for any $p$ in $\R^{nd}$. We then use the change of variable $rt^{-1/\alpha}=s$ to get that
\[\mathbb{E}\bigl[e^{i \langle p,N_t\rangle}\bigr] \, = \, \mathbb{E}\bigl[e^{i \langle p,t^{1/\alpha}N_1\rangle}].\]
This implies in particular our assertion on $N$. In a similar way, it is possible to get the analogous assertion on $M$.\newline
From lemma A.$2$ in \cite{Huang:Menozzi:Priola19}  with $m =3$, we know that there exist a family $\{p_{\overline{M}}(t,\cdot)\}_{t>0}$ of
densities on $\R^{nd}$ and a constant $C:=C(m,\alpha)$ such that
\[\vert D^l_{y}p_M(t,y) \vert \, \le \, Cp_{\overline{M}}(t,y)t^{-l/\alpha}\]
for any $t>0$, any $x$ in $\R^{nd}$ and any $l\in\{0,1,2\}$.\newline
Moreover, denoting $\overline{M}_t$ the random variable with density $p_{\overline{M}}(t,\cdot)$ and independent from $N_t$, we can easily
check from $p_{\overline{M}}(t,y)=t^{-nd/\alpha}p_{\overline{M}}(1,t^{-1/\alpha}x)$
that $\overline{M}$ is $\alpha$-selfsimilar
\[\overline{M}_t \, \overset{law}{=} \, t^{1/\alpha}\overline{M}_1.\]
We can finally define the family $\{q(t,\cdot)\}_{t>0}$ of densities as
\[q(t,x) \,:=\, \int_{\R^{nd}}p_{\overline{M}}(t,x-y)P_{N_t}(dy)\]
corresponding to the density of the random variable
\[\overline{S}_t \, ; = \, \overline{M}_t +N_t\]
for any fixed $t>0$.
Using Fourier transform and the already proven $\alpha$-selfsimilarity of $\overline{M}$ and $N$, we can show now that
\[\overline{S}_t \, \overset{law}{=} \, t^{1/\alpha}\overline{S}_1\]
or equivalently, that
\[q(t,y)=t^{-nd/\alpha}q(1,t^{-1/\alpha}y)\]
for any $t$ in $[0,T]$ and any $y$ in $\R^{nd}$. Moreover,
\[\mathbb{E}[\vert \overline{S}_t\vert^\gamma] \, = \, \mathbb{E}[\vert \overline{M}_t+N_t\vert^\gamma] \, = \,
Ct^{\gamma/\alpha}\bigl(\mathbb{E}[\vert \overline{M}_1\vert^\gamma]+\mathbb{E}[\vert N_t\vert^\gamma]\bigr) \, \le \,  Ct^{\gamma/\alpha}.\]
This shows in particular that equation \eqref{equation:integration_prop_of_q_stable} holds.
\end{proof}

We conclude this sub-section showing Control \eqref{eq:translation_inv_for_density} appearing in the proof of Proposition
\ref{prop:A_Priori_Estimates} for the diagonal regime. First of all, we will need the following lemma:

\begin{lemma}
Let $t$ in $[0,T]$, $x,b$ in $\R^{nd}$ such that $\vert b\vert \le ct^{1/\alpha}$ for some constant $c>0$.
Under [\textbf{A}], there exists a constant $C:=C(c)$ such that
\[\vert D^l_{x}p_S\bigl(t,x+b)\vert \, \le \, \tilde{C} \vert D^l_{x}p_S\bigl(t,x)\vert\]
\end{lemma}
\begin{proof}
Looking back at the proof of the previous lemma \ref{Appendix:Smoothing_effect}, we know that
\[D^l_{x}p_S(t,x+b) \, = \, \int_{\R^{nd}}D^l_{x}p_M(t,x+b-y)P_{N_t}(dy)\]
where $p_M(t,\cdot)$ is the density of $M_t$ and $P_{N_t}$ is the law of $N_t$, corresponding to the small and big jumps in the
Ito-L\'evy decomposition.\newline
From lemma $A.2$ in \cite{Huang:Menozzi:Priola19} we know moreover that
\[\vert D^l_{x}p_M(t,x+b-y)\vert \, \le \, \frac{C}{t^{\frac{l}{\alpha}}}p_{\overline{M}}(t,x+b-y) \,\,
\text{ where } \,\, p_{\overline{M}}(t,z) \, = \, \frac{C}{t^{\frac{nd}{\alpha}}}\frac{1}{\Bigl(1+\frac{\vert z\vert}{t^
\frac{1}{\alpha}}\Bigr)^3}.\]
It is then enough to show that
\[\begin{split}
p_{\overline{M}}(t,z+b) \, &= \, \frac{C}{t^{\frac{nd}{\alpha}}}\frac{1}{\Bigl(1+\frac{\vert z+b\vert}{t^
\frac{1}{\alpha}}\Bigr)^3} \, \le \, \frac{\tilde{C}}{t^{\frac{nd}{\alpha}}}\frac{1}{\Bigl(1+c+\frac{\vert z+b\vert}{t^
\frac{1}{\alpha}}\Bigr)^3} \\
&\le \, \frac{C}{t^{\frac{nd}{\alpha}}}\frac{1}{\Bigl(1+c\frac{\vert z\vert}{t^\frac{1}{\alpha}}-\frac{\vert b\vert}{t^\frac{1}{
\alpha}}\Bigr)^3} \,
\le \, \frac{C}{t^{\frac{nd}{\alpha}}}\frac{1}{\Bigl(1+\frac{\vert z\vert}{t^\frac{1}{\alpha}}\Bigr)^3} \\
&\le \,
Cp_{\overline{M}}(t,z).
\end{split}\]
to conclude the proof.
\end{proof}

\begin{proof}[Proof of Equation \eqref{eq:translation_inv_for_density}]
We start looking back to the proof of Lemma \ref{lemma:Smoothing_effect_frozen} to find that
\begin{multline*}
\bigl{\vert} D^\vartheta_{x}\tilde{p}^{\tau,\xi'} (t,s,x+ \lambda(x'- x),y)\bigr{\vert} \\
= \,
C(s-t)^{-\sum_{k=1}^{n}\frac{\vartheta_k}{\alpha_k}}\frac{1}{\det \mathbb{M}_{s-t}}
\bigl{\vert}\mathbf{d}^{\vert \vartheta\vert}_{z}p_S\bigl(s-t,\cdot)(\mathbb{M}^{-1}_{s-t} (\tilde{m}^{\tau,\xi}_{s,t}
(x)-y)\bigr)\bigr{\vert}
\end{multline*}
Moreover, we notice that
\[\mathbb{M}^{-1}_{s-t} \bigl(\tilde{m}^{\tau,\xi}_{s,t} (x+\lambda(x-x'))-y\bigr) \, = \,\mathbb{M}^{-1
}_{s-t} \bigl(\tilde{m}^{\tau,\xi}_{s,t} (x)-y\bigr)+\lambda\mathbb{M}^{-1}_{s-t} e^{A(s-t)}(x-x').\]
Then, Control \eqref{eq:translation_inv_for_density} follows immediately from the previous lemma once we have shown that
\[\bigl{\vert}\lambda\mathbb{M}^{-1}_{s-t} e^{A(s-t)}(x-x')\bigr{\vert} \, \le \, C(s-t)^{1/\alpha}\]
for some constant $C:=C(A)$. Indeed, fixed $i$ in $\llbracket 1,n\rrbracket$, we can exploit the structure of $A$ and $\mathbb{M}_{s-t}$
(cf. Equation \eqref{Proof:Scaling_Lemma2} in Scaling Lemma \ref{lemma:Scaling_Lemma}) to write that
\[
\begin{split}
\bigl[\mathbb{M}^{-1}_{s-t} e^{A(s-t)}(x-x') \bigr]_i \, &= \, \sum_{j=1}^{n}\sum_{k=1}^{n}\bigl[\mathbb{M}^{-1}_{s-t}\bigr]_{i,k}
\bigl[e^{A(s-t)}\bigr]_{k,j}(x-x')_j \\
&= \, \sum_{j=i}^{n}(s-t)^{-(i-1)}C_j(s-t)^{i-j}(x-x')_j.
\end{split}\]
Since moreover we assumed to be in a local diagonal regime, i.e.\ $\mathbf{d}^\alpha(x,x')\le (s-t)^{1/\alpha}$, we have that
\[
\begin{split}
\bigl{\vert}\bigl[\mathbb{M}^{-1}_{s-t} e^{A(s-t)}(x-x') \bigr]_i\bigr{\vert} \, &\le \, C\sum_{j=i}^{n}(s-t)^{-(j-1)}
\vert(x-x')_j\vert \\ 
&\le \, C\sum_{j=i}^{n}(s-t)^{-(j-1)}(s-t)^{\frac{1+\alpha(j-1)}{\alpha}} \\
&= \,C(s-t)^{1/\alpha}.
\end{split}\]
The proof is thus concluded.
\end{proof}

\subsection{Technical tools}

In this section, we present the proof of some technical results already used in the article, for the sake of completeness. \newline
We recall moreover that the results below can be proven also for the flow $\bar{\theta}_{s,\tau}(\xi)$ driven by a more general
perturbation $F$ under assumption $[\bar{\textbf{A}}]$ (cf. Section $7.1$), exploiting that
$\bar{F}_i$ is Lipschitz continuous in the $ x_{i-1}$ variable for any $i$ in $\llbracket 2, n \rrbracket$.

We begin proving Lemma \ref{lemma:Controls_on_Flow1} about the sensitivity of the H\"older flows, appearing in the proof of the
a priori estimates \eqref{eq:A_Priori_Estimates} of Proposition \ref{prop:A_Priori_Estimates}. For this reason, we will assume from this
point further to be under assumption [\textbf{A'}].

\paragraph{Proof of Lemma \ref{lemma:Controls_on_Flow1}.}
We start noticing that our result follows immediately using Young inequality, once we have shown that it holds
\begin{equation}\label{equation:Controls_on_Flow}
\bigl{\vert} (\theta_{s,t}(x)-\theta_{s,t}(x'))_i\bigr{\vert} \, \le \, C\Bigl[(s-t)^{\frac{1+\alpha(i-1)}{\alpha}}+\mathbf{d}^{1+\alpha(i-1)}(
x,x')\Bigr] \quad \text{for any $i$ in }\llbracket 1,n \rrbracket.
\end{equation}
Our proof will rely essentially in iterative applications of the Gr\"onwall lemma. We notice however that under [\textbf{A}], the
perturbation $ F_i$ is only H\"older continuous with respect to its $i$-th variable. To overcome this problem, we are going to mollify
(but only with respect to the variable of interest) the function $F$ in the following way: fixed a mollifier $\rho$ on $\R^d$, i.e. a
compactly supported, non-negative, smooth function such that $\Vert \rho \Vert_{L^1}=1$ and a family $\delta_i$ of positive constants to be
chosen later, the mollified version of the perturbation is given by $F^\delta=( F_1,F^{\delta_2}_2,\dots,F^{\delta_n}_n)$
where
\[F^{\delta_i}_i(t,z_{i:n}) \, := \,  F_i \ast_i \rho_{\delta_i}(t,z_{i:n}) \, = \, \int_{\R^d}
 F_i(t, z_i-\omega,z_{i+1},\dots,z_n)\frac{1}{\delta_i^d} \rho(\frac{\omega}{\delta_i}) \, d\omega.\]
We remark in particular that we do not need to mollify the first component $ F_1$ since it is regular enough, say $\beta$-H\"older
continuous in the first $d$-dimensional variable $ x_1$, by assumption [\textbf{R}].\newline
Then, standard results on mollifier theory and our current assumptions on $F$ show us that the following controls hold
\begin{align}
\label{Proof:Controls_on_flows_mollifier_stable} \vert  F_i(u,z) - F^\delta_i(u,z)\vert \, &\le \,\Vert  F_i \Vert_{L^\infty(C^{\gamma_+\beta}_d)}\delta_i^{\frac{\gamma_i+\beta}{1+\alpha(i-1)}}, \\
\label{Proof:Controls_on_flows_mollifier_stable1} \vert F^\delta_i(u,z) - F^\delta_i(u,z')\vert \, &\le \,C\Vert  F_i \Vert_{L^\infty(C^{\gamma_+\beta}_d)} \bigl[\delta_i^{
\frac{\gamma_i+\beta}{1+\alpha(i-1)}-1}\vert(z-z')_i\vert +\sum_{j=i+1}^{n}\vert(z-z')_{j}\vert^{\frac{\gamma_i+\beta}{1+
\alpha(j-1)}}\bigr].
\end{align}
We choose now $\delta_i$ for any $i$ in $\llbracket 2,n\rrbracket$ in order to have any contribution associated with the mollification
appearing in \eqref{Proof:Controls_on_flows_mollifier_stable} at a good current scale time. Namely, we would like $\delta_i$ to satisfy
\[\bigl{\vert} \bigl((s-t)^{\frac{1}{\alpha}}\mathbb{M}_{s-t}\bigr)^{-1}\bigl(F(u,z)-F^\delta(u,z)\bigr) \bigr{\vert} \,
\le \, C(s-t)^{-1}\]
for any $u$ in $[t,s]$ and any $z$ in $\R^{nd}$.
Using the mollifier controls \eqref{Proof:Controls_on_flows_mollifier_stable}, it is enough to ask for
\[\sum_{i=2}^n(s-t)^{-\frac{1}{\alpha_i}}\delta_i^{\frac{\gamma_i+\beta}{1+\alpha(i-1)}} \, \le \, C(s-t)^{-1}.\]
Recalling that $\gamma_i:=1+\alpha(i-2)$ by assumption [\textbf{R}], this is true if we fix for example,
\begin{equation}\label{Proof:Controls_on_Flows_Choice_delta_stable}
\delta_i \, = \, (s-t)^{\frac{\gamma_i}{\alpha}\frac{1+\alpha(i-1)}{\gamma_i+\beta}} \quad \text{for $i$ in $\llbracket 2,n\rrbracket$.}
\end{equation}
After this introductive part, we start controlling the last component of the flow. By construction of $\theta_{s,t}$, we can write that
\begin{align}
\label{Proof:Control_on_Flow_eq1}
\bigl{\vert} (&\theta_{s,t}(x)-\theta_{s,t}(x'))_n\bigr{\vert} \\ \notag
&= \, \Bigl{\vert} (x-x')_n +
\int_{t}^{s}\bigl{\{}\bigl[A(\theta_{v,t}(x)-\theta_{v,t}(x'))\bigr]_n + F_n(v,\theta_{v,t}(x))
-F_n(v,\theta_{v,t}(x'))\bigr{\}} \,  dv \Bigr{\vert} \\ \notag
&\le \, \vert (x-x')_n\vert \\\notag
&\quad \qquad+ \int_{t}^{s}\bigl{\{}A_{n,n-1}\vert(\theta_{v,t}(x)- \theta_{v,t}(x'))_{n-1}
\vert + \bigl{\vert}F_n(v,\theta_{v,t}(x)) - F_n(v,\theta_{v,t}(x'))\bigr{\vert}\bigr{\}} \,  dv
\end{align}
where in the last passage we have exploited the sub-diagonal structure of $A$ (cf. Equation \eqref{eq:def_matrix_A_stable}). If we focus only on
the last term involving the difference of the drifts, It holds now that
\begin{multline*}
\bigl{\vert} F_n(v,\theta_{v,t}(x)) - F_n(v,\theta_{v,t}(x'))\bigr{\vert} \,
\le \, \bigl{\vert} F_n(v,\theta_{v,t}(x)) - F^\delta_n(v,\theta_{v,t}(x))\bigr{\vert} \\
+
\bigl{\vert}F_n(v,\theta_{v,t}(x')) -
F^\delta_n(v,\theta_{v,t}(x') )\bigr{\vert} 
+ \bigl{\vert} F^\delta_n( v,\theta_{v,t}(x))-F^\delta_n(
v,\theta_{v,t}(x'))\bigr{\vert}.
\end{multline*}
Using the controls \eqref{Proof:Controls_on_flows_mollifier_stable}, \eqref{Proof:Controls_on_flows_mollifier_stable1} on the mollified drifts, we then
write from \eqref{Proof:Control_on_Flow_eq1} and the previous equation that
\begin{multline*}
\bigl{\vert} (\theta_{s,t}(x)-\theta_{s,t}(x'))_n\bigr{\vert} \,\le \,
\vert (x-x')_n\vert + 2(s-t)\delta_{n}^{\frac{\gamma_n+\beta}{1+\alpha(n-1)}}\\
+C\int_{t}^{s}
\bigl{\{}\bigl{\vert}(\theta_{v,t}(x)- \theta_{v,t}(x'))_{n-1} \bigr{\vert} +
\delta_n^{\frac{\gamma_n+\beta}{1+\alpha(n-1)}-1} \bigl{\vert}(\theta_{v,t}(x)-
\theta_{v,t}(x'))_n \bigr{\vert}\bigr{\}} \,  dv.
\end{multline*}
We now apply the Gr\"onwall lemma to show that
\[\bigl{\vert} (\theta_{s,t}(x)-\theta_{s,t}(x'))_n\bigr{\vert} \, \le \, C\Bigl[\vert(x-x')_n\vert+(s-t)
\delta_{n}^{\frac{\gamma_n+ \beta}{1+\alpha(n-1)}} + \int_{t}^{s} \bigl{\vert}(\theta_{v,t}(x)- \theta_{v,t}(x'))_{n-1}
\bigr{\vert} \,  dv\Bigr].\]
From our previous choice for $\delta_n$ (cf. Equation \eqref{Proof:Controls_on_Flows_Choice_delta_stable}), we know that
\[
(s-t)^{-\frac{1}{\alpha_n}}\delta_n^{\frac{\gamma_n+\beta}{1+\alpha(n-1)}}\, \le \,  C(s-t)^{-1}\] and thus, we can rewrite the last inequality as
\begin{equation}\label{Proof:Controls_on_flows1}
\bigl{\vert}(\theta_{s,t}(x) - \theta_{s,t}(x'))_n\bigr{\vert} \, \le \,C\Bigl[\bigl{\vert}(x -
x')_n\bigr{\vert} +(s-t)^{\frac{1+\alpha(n-1)}{\alpha}} + \int_{t}^{s}\bigl{\vert}(\theta_{v,t}(x) -
\theta_{v,t}(x'))_{n-1}\bigr{\vert} \, dv\Bigr].
\end{equation}
We would like now to obtain a similar control on the $(n-1)$-th term. As already done at the beginning of the proof, we can write that
\begin{multline*}
\bigl{\vert}(\theta_{s,t}(x) - \theta_{s,t}(x'))_{n-1}\bigr{\vert} \, \le \, \bigl{\vert}(x - x')_{n-1}\bigr{\vert} + C\delta_{n-1}^{\frac{\gamma_{n-1} + \beta}{1+\alpha(n-2)}}(s-t)+\int_{t}^{s}
\bigl{\vert}(\theta_{v,t}(x) - \theta_{v,t}(x'))_{n-2}\bigr{\vert} \\
+ \delta_{n-1}^{\frac{\gamma_{n-1} +\beta}{1+\alpha(n-2)}-1} \bigl{\vert}(\theta_{v,t}(x)
-\theta_{v,t}(x'))_{n-1}\bigr{\vert} +
\bigl{\vert}(\theta_{v,t}(x)-\theta_{v,t}(x'))_{n}\bigr{\vert}^{\frac{\gamma_{n-1} +\beta}{1+\alpha(n-1)}} \, dv
\end{multline*}
We then apply the Gr\"onwall lemma to find that
\begin{multline*}
\bigl{\vert}(\theta_{s,t}(x) - \theta_{s,t}(x'))_{n-1}\bigr{\vert} \, \le \, C\Bigl[\bigl{\vert}(x -
x')_{n-1}\bigr{\vert} + \delta_{n-1}^{\frac{\gamma_{n-1} + \beta}{1+\alpha(n-2)}}(s-t)\\
+\int_{t}^{s}\bigl{\{}\bigl{\vert}(\theta_{v,t}(x) - \theta_{v,t}(x'))_{n-2}\bigr{\vert} + \bigl{\vert}
(\theta_{v,t}(x) - \theta_{v,t}(x'))_{n} \bigr{\vert}^{\frac{\gamma_{n-1} +\beta}{1+\alpha(n-1)}}\bigr{\}} \, dv\Bigr].
\end{multline*}
Remembering our previous choice of $\delta_{n-1}$, it holds now that
\begin{multline}
\label{Proof:Controls_on_flows1/2}
\bigl{\vert}(\theta_{s,t}(x) - \theta_{s,t}(x'))_{n-1}\bigr{\vert} \, \le \, C
\Bigl[\vert(x-x')_{n-1}\vert+(s-t)^{\frac{1+\alpha(n-2)}{\alpha}}+\int_{t}^{s}\bigl{\vert}(\theta_{v,t}(x) -
\theta_{v,t}(x'))_{n-2}\bigr{\vert} \\
+ \bigl{\vert}(\theta_{v,t}(x) - \theta_{v,t}(x'))_n\bigr{\vert}^{\frac{\gamma_{n-1}+\beta}{1+\alpha(n-1)}} \, dv\Bigr].
\end{multline}
We then use equation \eqref{Proof:Controls_on_flows1} and the Jensen inequality to write
\begin{multline}
\label{Proof:Controls_on_flows1/3}
\bigl{\vert}(\theta_{s,t}(x) - \theta_{s,t}(x'))_{n-1}\bigr{\vert} \\
\le \, C \Bigl[\vert(x-x')_{n-1}\vert + (s-t)^{\frac{1+\alpha(n-2)}{\alpha}} +\int_{t}^{s}\bigl{\{}\bigl{\vert}(\theta_{v,t}(x) -
\theta_{v,t}(x'))_{n-2} \bigr{\vert} + \bigl{\vert}(x - x')_n\bigr{\vert}^{\frac{\gamma_{n-1}+\beta}{1+\alpha(n-1)}} \\
+(v-t)^{\frac{\gamma_{n-1}+\beta}{\alpha}}+ \Bigl(\int_{t}^{v}\bigl{\vert}(\theta_{\omega,t}(x) - \theta_{\omega,t}(x'))_{n-1}
\bigr{\vert}\,d\omega\Bigr)^{\frac{\gamma_{n-1}+\beta}{1+\alpha(n-1)}} \bigr{\}}\, dv\Bigr].
\end{multline}
The idea now is to use Gr\"onwall lemma again. To do so, we firstly move the exponent from the last integral term involving the $(n-1)$-th
term using the Young inequality:
\begin{multline*}
\Bigl(\int_{t}^{v}\bigl{\vert}(\theta_{\omega,t}(x) - \theta_{\omega,t}(x'))_{n-1}\bigr{\vert} \, d\omega
\Bigr)^{\frac{\gamma_{n-1}+\beta}{1+\alpha(n-1)}} \\
\le \, B^{-\frac{1+\alpha(n-1)}{\gamma_{n-1}+\beta}} \int_{t}^{v} \bigl{\vert}
(\theta_{\omega,t}(x) - \theta_{\omega,t}(x'))_{n-1}\bigr{\vert} \, d\omega + B^{\frac{1+\alpha(n-1)}{2\alpha - \beta}}
\end{multline*}
for a quantity $B$ to be fixed later.\newline
Since we need homogeneity with respect to time in equation \eqref{Proof:Controls_on_flows1/2}, we choose $B$ such that
\[B^{\frac{1+\alpha(n-1)}{2\alpha - \beta}} \, = \, (v-t)^{\frac{\gamma_{n-1}+\beta}{\alpha}} \, \Leftrightarrow \, B=(v
-t)^{\frac{\gamma_{n-1}+\beta}{\alpha}\frac{2\alpha - \beta}{1+\alpha(n-1)}}.\]
Plugging it into the general expression in \eqref{Proof:Controls_on_flows1/3}, we find that
\[\begin{split}
\bigl{\vert}(\theta_{s,t}&(x) - \theta_{s,t}(x'))_{n-1}\bigr{\vert} \, \le \, C \Bigl[\vert(x-x')_{n-1}\vert +
(s-t)^{\frac{1+\alpha(n-2)}{\alpha}} \\
&\qquad\qquad\qquad +\int_{t}^{s}\Bigl{\{}\bigl{\vert}(\theta_{v,t}(x) - \theta_{v,t}(x'))_{n-2} \bigr{\vert} + \bigl{\vert}(x
-x')_n\bigr{\vert}^{\frac{\gamma_{n-1}+\beta}{1+\alpha(n-1)}} +(v - t)^{\frac{\gamma_{n-1}+\beta}{\alpha}} \\
&\qquad\qquad\qquad\qquad\qquad +(v-t)^{\frac{\beta}{\alpha}-2}\int_{t}^{v}\bigl{\vert}(\theta_{\omega,t}(x) - \theta_{\omega,t}(x'))_{n-1}
\bigr{\vert} \, d\omega \Bigr{\}}\, dv\Bigr] \\
&\le \, C \Bigl[\vert(x-x')_{n-1}\vert +(s-t)^{\frac{1+\alpha(n-1)}{\alpha}} + (s-t)\bigl{\vert}(x-x')_n \bigr{\vert}^{ \frac{\gamma_{n-1}+
\beta}{1+\alpha(n-1)}} + (s - t)^{\frac{\gamma_{n-1}+\beta+\alpha}{\alpha}}\\
&\qquad + \int_{t}^{s}\Bigl{\{}\bigl{\vert} (\theta_{v,t}(x) - \theta_{v,t}(x'))_{n-2} \bigr{\vert} + (v-t)^{\frac{\beta}{\alpha}-1} \sup_{\omega\in[t,v]}
\bigl{\vert}(\theta_{\omega,t}(x) -\theta_{\omega,t}(x'))_{n-1} \bigr{\vert}\Bigr{\}} \, dv\Bigr].
\end{split}\]
Since the previous inequality is also true for any $\overline{s}$ in $[t,s]$, it follows that
\[\begin{split}
\sup_{\overline{s}\in[0,s]}&\bigl{\vert}(\theta_{\overline{s},t}(x) - \theta_{\overline{s},t}(x'))_{n-1}\bigr{\vert} \\
&\le \, C \Bigl[\vert(x-x')_{n-1}\vert + (s-t)^{\frac{1+\alpha(n-2)}{\alpha}} + (s-t)\bigl{\vert}(x - x')_n \bigr{\vert}^{
\frac{\gamma_{n-1}+\beta}{1+\alpha(n-1)}}+(s - t)^{\frac{\gamma_{n-1}+\beta+\alpha}{\alpha}} \\
&\quad \quad+\int_{t}^{s} \Bigl{\{}\bigl{\vert} (\theta_{v,t}(x) -
\theta_{v,t}(x'))_{n-2}\bigr{\vert} + (v-t)^{\frac{\beta}{\alpha}-1} \sup_{\omega \in[t,v]}\bigl{\vert}(\theta_{\omega,t}
(x) - \theta_{\omega,t}(x'))_{n-1} \bigr{\vert}\Bigr{\}} \, dv\Bigr].
\end{split}
\]
We can finally apply the Gr\"onwall lemma to show that for any $s$ in $[t,T]$, there exists a constant $C$ such that
\begin{multline*}
\bigl{\vert}(\theta_{s,t}(x) - \theta_{s,t}(x'))_{n-1}\bigr{\vert} \,\le \, C \Bigl[\vert(x-x')_{n-1}\vert + (s-t)^{\frac{1+\alpha(n-2)}{\alpha}} + (s-t)\vert(x - x')_n\vert^{\frac{\gamma_{n-1}+\beta}{
1+\alpha(n-1)}} \\+ \int_{t}^{s} \bigl{\vert} (\theta_{v,t}(x) - \theta_{v,t}(x'))_{n-2}
\bigr{\vert} \, dv\Bigr].
\end{multline*}
Moreover, thanks to the Young inequality we know that
\[(s-t)\bigl{\vert}(x - x')_n \bigr{\vert}^{\frac{\gamma_{n-1}+\beta}{1+\alpha(n-1)}} \, \le \,
C\bigl{\{}(s-t)^{\frac{1+\alpha(n-2)}{\alpha}} +\vert(x - x')_n \vert^{\frac{\gamma_{n-1}+\beta}{1+\alpha(n-1)}
\frac{1+\alpha(n-2)}{1+\alpha(n-3)}}\bigr{\}}\]
and remembering that $\mathbf{d}(x,x')\le 1$ by hypothesis,
\[\vert(x - x')_n \vert^{\frac{\gamma_{n-1}+\beta}{1+\alpha(n-1)}\frac{1+\alpha(n-2)}{1+\alpha(n-3)}} \, \le \, \vert(x -
x')_n \vert^{\frac{\gamma_{n-1}+\beta}{\gamma_{n-1}}\frac{1+\alpha(n-2)}{1+\alpha(n-1)}} \, \le \, \vert(x - x')_n
\vert^{\frac{1+\alpha(n-2)}{1+\alpha(n-1)}}.\]
We then use it to write for any $v$ in $[t,T]$,
\begin{multline*}
\bigl{\vert}(\theta_{v,t}(x) - \theta_{v,t}(x'))_{n-1}\bigr{\vert} \, \le \, C \Bigl[\vert(x-x')_{n-1}\vert + (v-t)^{\frac{1+\alpha(n-2)}{\alpha}} + \vert(x - x')_n \vert^{\frac{
1+\alpha(n-2)}{1+\alpha(n-1)}} \\
+ \int_{t}^{v} \bigl{\vert} (\theta_{\omega,t}(x) - \theta_{\omega,t}(x'))_{n-2}
\bigr{\vert} \, d\omega\Bigr].
\end{multline*}
Going back to equation \eqref{Proof:Controls_on_flows1}, we plug in the last bound to find that
\[\begin{split}
\bigl{\vert}(\theta_{s,t}(x) - \theta_{s,t}(x'))_n\bigr{\vert} \, &\le \, C\Bigl[\vert(x - x')_n \vert
+(s-t)^{\frac{1+\alpha(n-1)}{\alpha}} + (s-t)\vert(x - x')_{n-1}\vert \\
&\,\,+ (s-t)\vert(x - x')_n \vert^{\frac{1+\alpha(n-2)}{1+\alpha(n-1)}} + \int_{t}^{s}\int_{t}^{v} \bigl{\vert} (\theta_{t,
\omega}(x) - \theta_{\omega,t}(x'))_{n-2}\bigr{\vert}\, d\omega dv \Bigr] \\
&\le \, C\Bigl[\vert(x - x')_n\vert + (s-t)^{\frac{1+\alpha(n-1)}{\alpha}} + \vert(x - x')_{n-1} \vert^{\frac{1+
\alpha(n-1)}{1+\alpha(n-2)}} \\
&\,\,+ \int_{t}^{s}\int_{t}^{v} \bigl{\vert}(\theta_{\omega,t}(x) - \theta_{\omega,t}(x'))_{n-2}
\bigr{\vert}\, d\omega dv \Bigr]
\end{split}\]
where in the last passage we used again the Young inequality to show that
\[(s-t)\vert(x - x')_{n-1}\vert \, \le \, C(s-t)^{\frac{1+\alpha(n-1)}{\alpha}}+ \vert(x - x')_{n-1}\vert^{\frac{
1+\alpha(n-1)}{1+\alpha( n-2)}}\]
and
\[(s-t)\vert(x - x')_n\vert^{\frac{1+\alpha(n-2)}{1+\alpha(n-1)}} \, \le \, C(s-t)^{\frac{1+\alpha(n-1)}{\alpha}}+ \vert(x -
x')_n \vert.\]
This approach may be naturally iterated up to the first term of the chain, so that
\begin{multline*}
\bigl{\vert}(\theta_{s,t}(x) - \theta_{s,t}(x'))_n\bigr{\vert} \,
\le \, C \Bigl[\sum_{j=2}^{n}\vert(x-x')_j\vert^{\frac{1+\alpha(n-1)}{1+\alpha(j-1)}} + (s-t)^{\frac{1+\alpha(n-1)}{\alpha}} \\ 
+\int_{t}^{v_n=s} dv_{n-1}\dots\int_{t}^{v=2}dv_1\bigl{\vert} (\theta_{v_1,t}(x) - \theta_{v_1,t}(x'))_1 \bigr{\vert}\Bigr].
\end{multline*}
In a similar manner, we can show for any $i$ in $\llbracket 2,n\rrbracket$,
\begin{multline}\label{Proof:Controls_on_flows2}
\bigl{\vert}(\theta_{s,t}(x) - \theta_{s,t}(x'))_i\bigr{\vert} \,
\le \, C \Bigl[\sum_{j=2}^{n}\vert(x-x')_j \vert^{ \frac{1+\alpha(i-1)}{1+\alpha(j-1)}} + (s-t)^{\frac{1+\alpha(i-1)}{\alpha}} \\
+ \int_{t}^{v_i=s}dv_{i-1}\dots\int_{t}^{v=2}dv_1\bigl{\vert} (\theta_{v_1,t}(x) - \theta_{v_1,t}(x'))_1 \bigr{\vert}\Bigr].
\end{multline}
Since all the non-integral terms in \eqref{Proof:Controls_on_flows2} are compatible with the statement of the lemma, it remains to find the proper bound for the first component of the flow. As before, let us consider $\overline{s}$ in $[t,s]$. We can write
\[\vert (\theta_{\overline{s},t}(x)-\theta_{\overline{s},t}(x'))_1\vert \, \le \, \vert (x-x')_1\vert
+C\sum_{j=1}^{n}\int_{t}^{\overline{s}}
\vert (\theta_{v,t}(x)-\theta_{v,t}(x'))_j\vert^{\frac{\beta}{1+\alpha(j-1)}} \, dv\]
or, passing to the supremum on both sides,
\begin{multline*}
\sup_{\overline{s}\in[t,s]}\vert( \theta_{\overline{s},t}(x)-\theta_{\overline{s},t}(x'))_1\vert \,
\le \, \vert (x-x')_1\vert +C\Bigl{\{}(s-t)\bigl(\sup_{v\in[t,s]}\vert (\theta_{v,t}(x)-\theta_{v,t}(x'))_1\vert\bigr)^\beta \\
+ \sum_{j=2}^{n}\int_{t}^{s} \vert (\theta_{v,t}(x)-\theta_{v,t}(x'))_j\vert^{\frac{\beta}{1+\alpha(j-1)}} \, dv\Bigr{\}}.
\end{multline*}
Using equation \eqref{Proof:Controls_on_flows2}, it holds now that
\begin{align}\notag
\sup_{\overline{s}\in[t,s]}(\vert \theta_{\overline{s},t}(x)-\theta_{\overline{s},t}(x'))_1\vert \, \le \, \vert
(x-&x')_1\vert +C\Bigl{\{}(s-t)\bigl(\sup_{v\in[t,s]}\vert
(\theta_{v,t}(x)-\theta_{v,t}(x'))_1\vert\bigr)^\beta\\ \notag
&+ \sum_{j=2}^{n}\Bigl[(s-t)\bigl((s-t)^{\frac{1+\alpha(j-1)}{\alpha}} +\sum_{k=2}^{n} \vert(x-x')_k \vert^{
\frac{1+\alpha(j-1)}{1+\alpha(k-1)}} \\
&+ (s-t)^{j-1}\sup_{v\in[t,s]}\vert (\theta_{v,t}(x)-\theta_{v,t}(x'))_1
\vert\bigr)^{\frac{\beta}{1+\alpha(j-1)}}\Bigr]\Bigr{\}}.
\label{Proof:Controls_on_flows3}
\end{align}
We then apply the Jensen inequality to show that
\begin{align}\notag
\sup_{\overline{s}\in[t,s]}(\vert \theta_{\overline{s},t}(x)-\theta_{\overline{s},t}(x'))_1\vert \, &\le \, \vert
(x-x')_1\vert + C\Bigl{\{}(s-t)\bigl[\sup_{v\in[t,s]}\vert(\theta_{v,t}(x)-\theta_{v,t}(x'))_1\vert
\bigr]^\beta\\\notag
&\qquad+\sum_{j=2}^{n}C(s-t)\bigl[(s-t)^{\frac{\beta}{\alpha}}
+ \sum_{k=2}^{n} \vert(x-x')_k \vert^{ \frac{\beta}{1+\alpha(k-1)}} \\\notag
&\qquad + (s-t)^{\frac{(j-1)\beta}{1+\alpha(j-1)}}\sup_{v\in[t,s]}\vert
(\theta_{v,t}(x)-\theta_{v,t}(x'))_1 \vert^{\frac{\beta}{1+\alpha(j-1)}} \bigr]\Bigr{\}} \\\notag
&\le \,C\Bigl{\{}\vert (x-x')_1\vert + (s-t)^{\frac{\alpha+\beta}{\alpha}} + (s-t)\sum_{k=2}^{n}
\vert(x-x')_k\vert^{\frac{\beta}{1+\alpha(k-1)}}\\
& \qquad+\sum_{j=1}^{n}(s-t)^{1+\frac{(j-1)\beta}{1+\alpha(j-1)}}\sup_{v\in[t,s]}\vert (\theta_{v,t}(x)-\theta_{v,t}(x'))_1
\vert^{\frac{\beta}{1+\alpha(j-1)}}\Bigr{\}}.\label{Proof:Controls_on_flows4}
\end{align}
From Young inequality, we can deduce now that
\[(s-t)\vert(x-x')_k\vert^{\frac{\beta}{1+\alpha(k-1)}} \, \le \,C\bigl((s-t)^{\frac{1}{1-\beta}}+ \vert(x-x')_k
\vert^{\frac{1}{1+\alpha(k-1)}}\bigr)\]
and
\begin{multline*}
(s-t)^{1+\frac{(j-1)\beta}{1+\alpha(j-1)}}\sup_{v\in[t,s]}\vert
(\theta_{v,t}(x)-\theta_{v,t}(x'))_1\vert^{\frac{\beta}{1+\alpha(j-1)}}\\
\le \, C\Bigl{\{}(s-t)^{\frac{1+(\alpha+\beta)(j-1)}{1+\alpha(j-1)-\beta}}+\sup_{v\in[t,s]}\vert (\theta_{v,t}(x)-\theta_{v,t}(x'))_1\vert
\Bigr{\}}
\end{multline*}
Plugging these inequalities into the main one \eqref{Proof:Controls_on_flows4}, we find that
\[\begin{split}
\sup_{\overline{s}\in[t,s]}(\vert \theta_{\overline{s},t}(x)-\theta_{\overline{s},t}(x'))_1\vert \, &\le \,
C\Bigl{\{}\vert (x-x')_1\vert +
(s-t)^{\frac{\alpha+\beta}{\alpha}}+\sum_{k=2}^{n} \vert(x-x')_k\vert^{\frac{1}{1+\alpha(k-1)}}\\
&\qquad\qquad+\sum_{j=1}^{n}(s-t)^{\frac{1+(\alpha+\beta)(j-1)}{1+\alpha(j-1)-\beta}}+\sup_{v\in[t,s]}\vert
(\theta_{v,t}(x)-\theta_{v,t}(x'))_1\vert\Bigr{\}} \\
&\le \,C\Bigl{\{}(s-t)^{\frac{\alpha+\beta}{\alpha}} +(s-t)^{\frac{1}{1-\beta}}+d(x,x')\\
&\qquad\qquad+\sum_{j=1}^n(s-t)^{\frac{1+(\alpha+\beta)(j-1)}{1+
\alpha(j-1)-\beta}}+\sup_{v\in[t,s]}\vert (\theta_{v,t}(x)-\theta_{v,t}(x'))_1\vert\Bigr{\}}
\end{split}\]
Remembering that $s-t\le T-t\le 1$, it finally holds that
\[\vert \theta_{s,t}(x)-\theta_{s,t}(x'))_1\vert \, \le \, C\bigl((s-t)^{1/\alpha} + \mathbf{d}(x,x')\bigr)\]
since by assumption [\textbf{P}],
\[\frac{\alpha+\beta}{\alpha}\, > \, \frac{1}{1-\beta} \, > \, \frac{1}{\alpha}\]
and
\[\frac{1+(\alpha+\beta)(j-1)}{1+\alpha(j-1)-\beta} \, = \, 1+\frac{\beta j}{1+\alpha j-(\alpha+\beta)} \, > \, 1+\frac{\beta j}{\alpha j}
\, > \, 1+\Bigl(\frac{1-\alpha}{\alpha}\Bigr)\, = \, \frac{1}{\alpha}.\]
Plugging this control in equation \eqref{Proof:Controls_on_flows2}, we then conclude since
\[\begin{split}
\bigl{\vert}(\theta_{s,t}(&x) - \theta_{s,t}(x'))_i\bigr{\vert} \\
&\le \,C\Bigl(\mathbf{d}^{1+\alpha(i-1)}(x,x')+ (s-t)^{\frac{1+\alpha(i-1)}{\alpha}} +(s-t)^{i-1}\sup_{\overline{s}\in[t,s]}(\vert
\theta_{\overline{s},t}(x)-\theta_{\overline{s},t}(x'))_1\vert\Bigr) \\
&\le \, C\Bigl(\mathbf{d}^{1+\alpha(i-1)}(x,x')+ (s-t)^{\frac{1+\alpha(i-1)}{\alpha}}+(s-t)^{i-1}\bigl((s-t)^{1/\alpha} +
\mathbf{d}(x,x')\bigr)\Bigr) \\
&\le \, C\Bigl((s-t)^{\frac{1+\alpha(i-1)}{\alpha}} + \mathbf{d}^{1+\alpha(i-1)}(x,x')\Bigr),
\end{split}
\]
using again the Young inequality in the last passage. The proof is complete.

We can now prove the two results (Lemmas \ref{lemma:Controls_on_means1} and Lemma \ref{lemma:Controls_on_means}) concerning the sensitivity
of the frozen shift $\tilde{m}^{\tau,\xi}_{s,t}$.

\paragraph{Proof of Lemma \ref{lemma:Controls_on_means1}.}
From the integral representation of $\tilde{m}^{t,x}_{s,t}(y)$ (cf. Equation \eqref{eq:def_tilde_m_stable}), we can write that
\begin{multline*}
\bigl{\vert}\bigl(\tilde{m}^{t,x}_{s,t}(y)-\tilde{m}^{t,x'}_{s,t} (y')\bigr)_1\bigr{\vert} \, \le \, \int_{t}^{s}\bigl{\vert}  F_1(v,\theta_{v,t}(x)) -  F_1(v,\theta_{v,t}(x')) \bigr{\vert} \, dv \\
 \le \, C\Vert F\Vert_H\int_{t}^{s}\mathbf{d}^{\beta}\bigl(\theta_{v,t}(x),\theta_{v,t}(x')\bigr) \, dv
\end{multline*}
where in the second passage we used that $ F_1$ is in $C^{\beta}_{b,d}(\R^{nd})$.
Thanks to the Control on the flows (Lemma \ref{lemma:Controls_on_Flow1}), it then holds that
\[\bigl{\vert} \bigl(\tilde{m}^{t,x}_{s,t}(y)-\tilde{m}^{t,x'}_{s,t} (y')\bigr)_1\bigr{\vert} \, \le \,C\Vert
F\Vert_H(s-t)\bigl[\mathbf{d}^{\beta}(x,x')+(s-t)^{\frac{\beta}{\alpha}}\bigr]\]
and we have concluded.

\paragraph{Proof of Lemma \ref{lemma:Controls_on_means}.}
We know from Lemma \ref{lemma:identification_theta_m_stable} that $\tilde{m}^{t,x'}_{t_0,t} (x') =
\theta_{t_0,t}(x')$. Fixed $i$ in $\llbracket 1, n\rrbracket$, we can then write that
\[
\begin{split}
\bigl( \tilde{m}^{t,x}_{t_0,t}(x') -\tilde{m}^{t,x'}_{t_0,t} (x')\bigr)_i \, &= \,
\bigl(\tilde{m}^{t,x}_{t_0,t}(x') -\theta_{t_0,t}(x')\bigr)_i \\
&= \, \bigl(\tilde{m}^{t,x}_{t_0,t}(x')-\theta_{t_0,t}(x)\bigr)_i +  \bigl(\theta_{t_0,t}(x)
-\theta_{t_0,t}(x')\bigr)_i.
\end{split}\]
We start focusing on the first term of the above expression. From the integral representation of $\tilde{m}^{t,x}_{t_0,t}(
x')$ and $\theta_{t_0,t}(x)$, it holds that
\begin{equation}\label{Proof:Controls_on_means_1}
\tilde{m}^{t,x}_{t_0,t}(x') -\theta_{t_0,t}(x) \, = \, x'-x+\int_{t}^{t_0}A\bigl[
\tilde{m}^{t,x}_{v,t}(x') -\theta_{v,t}(x)\bigr] \, dv.
\end{equation}
Remembering from \eqref{eq:def_matrix_A_stable} that $A$ is sub-diagonal, it follows that
\begin{equation}\label{Proof:Controls_on_means_2}
\bigl(\tilde{m}^{t,x}_{t_0,t}(x') -\theta_{t_0,t}(x)\bigr)_i \, = \, (x'-x)_i+A_{i,i-1}\int_{t}^{t_0}
\bigl( \tilde{m}^{t,x}_{v,t}(x')-\theta_{v,t}(x)\bigr)_{i-1} \, dv
\end{equation}
for any $i$ in $\llbracket 2,n\rrbracket$ and
\[\bigl(\tilde{m}^{t,x}_{t_0,t}(x') -\theta_{t_0,t}(x)\bigr)_1 \, = \, (x'-x)_1.\]
Iterating the process, we  can find that
\[\bigl{\vert}\bigl(\tilde{m}^{t,x}_{t_0,t}(x') -\theta_{t_0,t}(x)\bigr)_i \bigr{\vert} \, \le \,
C\sum_{k=1}^{i}\bigl{\vert}(x'-x)_k\bigr{\vert}(t_0-t)^{i-k}.\]
On the other side, the integral representation of $\theta_{s,\tau}(\xi)$ (Equation \eqref{Flow}) allows us to write that
\begin{multline}\label{Proof:Controls_on_means_3}
\bigl(\theta_{t_0,t}(x) -\theta_{t_0,t}(x')\bigr)_i \,
= \, (x-x')_i+A_{i,i-1}\int_{t}^{t_0} \bigl{\{}\bigl(
\theta_{t_0,t}(x) -
\theta_{t_0,t}(x')\bigr)_{i-1} \\
+  F_i(v,\theta_{v,t}(x))- F_i(v,\theta_{v,t}(x'))\bigr{\}} \, dv
\end{multline}
for any $i$ in $\llbracket 2,n\rrbracket$ and
\begin{equation}\label{Proof:Controls_on_means_4}
\bigl(\theta_{t_0,t}(x) -\theta_{t_0,t}(x')\bigr)_1 \, = \, (x-x')_1 +\int_{t}^{t_0} \bigl{\{}
 F_1(v,\theta_{v,t}(x))- F_1(v,\theta_{v,t}(x')) \bigr{\}} \,
dv.
\end{equation}
Fixed $i$ in $\llbracket 2,n\rrbracket$, it then follows from \eqref{Proof:Controls_on_means_1} and \eqref{Proof:Controls_on_means_3} that
\begin{multline*}
\bigl{\vert}\bigl( \tilde{m}^{t,x}_{t_0,t}(x') -\tilde{m}^{t,x'}_{t_0,t} (x')\bigr)_i\bigr{\vert} \, \le \,
C\Vert F \Vert_H \Bigl(\sum_{k=1}^{i-1}\vert(x'-x)_k\vert(t_0-t)^{i-k}\\
+\int_{t}^{t_0}\Bigl{\{}\bigl{\vert}\bigl(\theta_{v,t}(x) -\theta_{v,t}(x')\bigr)_{i-1}\bigr{\vert}+
\sum_{j=i}^{n}\bigl{\vert} \bigl(\theta_{v,t}(x) - \theta_{v,t}(x') \bigr)_{j}
\bigr{\vert}^{\frac{\gamma_i+\beta}{1+\alpha(j-1)}}\Bigr{\}} \, dv \Bigr).
\end{multline*}
Also, from \eqref{Proof:Controls_on_means_2} and \eqref{Proof:Controls_on_means_4}, it holds that
\[\bigl{\vert}\bigl( \tilde{m}^{t,x}_{t_0,t}(x') -\tilde{m}^{t,x'}_{t_0,t} (x')\bigr)_1\bigr{\vert} \, \le
\,C\Vert F \Vert_H\int_{t}^{t_0}\sum_{j=1}^{n}\bigl{\vert}\bigl(\theta_{v,t}(x) - \theta_{v,t}(x') \bigr)_{j}
\bigr{\vert}^{\frac{\beta}{1+\alpha(j-1)}} \, dv.\]
Using now Lemma \ref{lemma:Controls_on_Flow1}, we can show that
\begin{multline*}
\bigl{\vert}\bigl( \tilde{m}^{t,x}_{t_0,t}(x') -\tilde{m}^{t,x'}_{t_0,t} (x')\bigr)_i\bigr{\vert}\, \le \,
C\Vert F \Vert_H\Bigl(\sum_{k=1}^{i-1}\vert(x'-x)_k\vert(t_0-t)^{i-k} +(t_0-t)^{\frac{1+\alpha(i-2)}{\alpha}+1}\\
+(t_0-t)\mathbf{d}^{1+\alpha(i-2)}(x,x')+(t_0-t)^{\frac{1+\alpha(i-2)+\beta}{\alpha}+1}+
(t_0-t)\mathbf{d}^{1+\alpha(i-2)+\beta}(x,x')\Bigr)
\end{multline*}
for any $i$ in $\llbracket 2,n\rrbracket$ and
\[\bigl{\vert}\bigl( \tilde{m}^{t,x}_{t_0,t}(x') -\tilde{m}^{t,x'}_{t_0,t} (x')\bigr)_1\bigr{\vert} \, \le
\,C\Vert F \Vert_H(t_0-t)^{\frac{\beta+\alpha}{\alpha}}+(t_0-t)\mathbf{d}^{\beta}(x,x').\]
Since $t_0-t=c_0\mathbf{d}^{\alpha}(x,x')$ by Equation \eqref{eq:def_t0}, we can conclude that
\[\begin{split}
\bigl{\vert}\bigl( \tilde{m}^{t,x}_{t_0,t}(x') -\tilde{m}^{t,x'}_{t_0,t} (x')\bigr)_i\bigr{\vert}\, &\le \,
C\Vert F \Vert_H \Bigl{\{}\sum_{k=1}^{i-1}\mathbf{d}^{1+\alpha(k-1)}(x',x)c^{i-k}_0\mathbf{d}^{\alpha(i-k)}(x,x')\\ 
&\qquad\quad +
c^{\frac{1+\alpha(i-1)}{\alpha}}_0\mathbf{d}^{1+\alpha(i-1)}(x,x') 
+
c_0\mathbf{d}^{1+\alpha(i-1)}(x,x')\\
&\qquad\quad+ c_0^{\frac{1+\alpha(i-2)+\beta}{\alpha}+1}\mathbf{d}^{1+\alpha(i-1)+\beta}(x,x')+c_0
\mathbf{d}^{1+\alpha(i-1)+\beta}(x,x')\Bigr{\}} \\
&\le \, C\Vert F \Vert_H \Bigl{\{}\bigl(c_0+c^{\frac{1+\alpha(i-1)}{\alpha}}_0\bigr)\mathbf{d}^{1+\alpha(i-1)}(x,x')
\\
&\qquad\quad+\bigl(c_0+ c_0^{\frac{1+\alpha(i-1)+\beta}{\alpha}}\bigr)\mathbf{d}^{1+\alpha(i-1)+\beta}(x,x')\Bigr{\}} \\
&\le \, Cc_0\Vert F \Vert_H\mathbf{d}^{1+\alpha(i-1)}(x,x')
\end{split}\]
for any $i$ in $\llbracket 2,n\rrbracket$ and
\[\bigl{\vert}\bigl( \tilde{m}^{t,x}_{t_0,t}(x') -\tilde{m}^{t,x'}_{t_0,t} (x')\bigr)_1\bigr{\vert} \, \le
\,C\Vert F\Vert_H\bigl(c^{\frac{\beta+\alpha}{\alpha}}_0+c_0\bigr)\mathbf{d}^{\alpha+\beta}(x,x')\, \le \,Cc_0\Vert
F \Vert_H\mathbf{d}^{\alpha+\beta}(x,x'),\]
where in the last passage we used that $c_0\le 1$ and $\mathbf{d}(x,x')\le 1$.
After summing all the terms together at the right scale, we finally show that
\[\mathbf{d}(\tilde{m}^{t,x}_{t_0,t}(x'),\tilde{m}^{t,x'}_{t_0,t} (x')) \, \le \, Cc_0^{\frac{1}{1+\alpha(n-1)}}
\Vert F\Vert_H\mathbf{d}(x,x')\]
thanks to convexity inequalities and $c_0\le 1$.

We conclude this section showing the reverse Taylor formula which was used in the proof of Lemma \ref{lemma:Holder_modulus_Non-Deg} in
the diagonal regime to handle the discontinuity term:

\begin{lemma}[Reverse Taylor expansion]
\label{lemma:Reverse_Taylor_Expansion}
Let $\gamma$ be in $(1,2)$, $\phi$ a function in $C^\gamma_{b,d}(\R^{nd})$ and $x,x'$ two points in $\R^{nd}$. Then, there
exists a constant $C:=C(\gamma)$ such that
\[\vert D_{ x_1}\phi(x)-D_{ x_1}\phi(x')\vert \, \le \, C \Vert \phi \Vert_{C^\gamma_{b,d}}
\mathbf{d}^{\gamma-1}(x,x').\]
\end{lemma}
\begin{proof}
We start decomposing the left-hand side $D_{ x_1}\phi(x)-D_{ x_1}\phi(x')$ into $I_1+I_2+I_3$ where we denoted
\[\begin{split}
I_1 \, &:= \,  \Bigl(\int_{0}^{1} D_{ x_1}\phi(x)-D_{ x_1}\phi( x_1+\lambda \mathbf{d}(x,x'),(x)_{2:n})
\, d\lambda\Bigr)\\
I_2 \, &:= \,- \Bigl(\int_{0}^{1} D_{ x_1}\phi(x')-D_{ x_1}\phi( x_1+\lambda
\mathbf{d}(x,x'),(x')_{2:n}) \, d\lambda\Bigr) \\
I_3 \, &:= \,- \Bigl(\int_{0}^{1} D_{ x_1}\phi( x_1+\lambda \mathbf{d}(x,x'),(x')_{2:n}) -
D_{ x_1}\phi( x_1+\lambda \mathbf{d}(x,x'),(x)_{2:n}) \, d\lambda\Bigr).
\end{split}
\]
The first two components can be treated directly using that $D_{ x_1}\phi$ is in $C^{\gamma-1}(\R^d)$ with respect to the first
non-degenerate variable. Indeed,
\begin{multline*}
\vert I_1\vert \, \le \, \int_{0}^{1} \vert D_{ x_1}\phi(x)-D_{ x_1}\phi( x_1+\lambda \mathbf{d}(x,x'),(x)_{2:n})
\vert\, d\lambda \\
\le \, C\Vert \phi \Vert_{C^\gamma}\int_{0}^{1} \vert\lambda \mathbf{d}(x,x') \vert^{\gamma-1} \, d\lambda \, \le \, C\Vert \phi
\Vert_{C^\gamma}\mathbf{d}^{
\gamma-1}(x,x')
\end{multline*}
and
\[\begin{split}
\vert I_2\vert \, &\le \, \int_{0}^{1} \vert D_{ x_1}\phi(x')-D_{ x_1}\phi( x_1+\lambda \mathbf{d}(x,x'),(x')_{2:n})
\vert\, d\lambda \\
&\le\,C\Vert\phi\Vert_{C^\gamma}\int_{0}^{1} \vert (x'-x)_1 + \lambda \mathbf{d}(x,x')
\vert^{\gamma-1}\,d\lambda\\
&\le\,C\Vert\phi\Vert_{C^\gamma}\mathbf{d}^{
\gamma-1}(x,x')
\end{split}\]
where in the last expression we used Young inequality.\newline
To control the last term, we assume for the sake of brevity to be in the scalar case, i.e.\ $d=1$. In the general setting, the proof below
can be reproduced component-wise. The idea is to use a reverse Taylor expansion to pass from the derivative to the function itself. Namely,
\[
\begin{split}
\vert I_3 \vert \, &= \, \frac{1}{\mathbf{d}(x,x')}\Bigl{\vert}\int_{0}^{1}\bigl[\partial_{\lambda}\phi( x_1+\lambda
\mathbf{d}(x,x'),(x')_{2:n}) - \partial_\lambda\phi( x_1+\lambda \mathbf{d}(x,x'),(x)_{2:n})\bigr] \, d\lambda\Bigr{\vert} \\
&\le \, \frac{1}{\mathbf{d}(x,x')} \bigl{\vert} \phi( x_1+ \mathbf{d}(x,x'),(x')_{2:n}) - \phi( x_1,(x')_{2:n}) +
\phi( x_1+ \mathbf{d}(x,x'),(x)_{2:n}) - \phi(x)\bigr{\vert} \\
&\le \, C \Vert \phi \Vert_{C^\gamma}\mathbf{d}^{\gamma-1}(x,x').
\end{split}\]
We have thus concluded the proof.
\end{proof}

\setcounter{equation}{0}
\chapter{Schauder estimates for degenerate L\'evy Ornstein-Uhlenbeck operators}
\fancyhead[LE]{Chapter \thechapter. Schauder estimates for linear L\'evy operators}

\label{Chap:Schauder_Estimates_Levy}

\paragraph{Abstract:}
We establish global Schauder estimates for integro-partial differential equations (IPDE) driven by a possibly degenerate
L\'evy Ornstein-Uhlenbeck operator, both in the elliptic and parabolic setting, using some suitable anisotropic H\"older spaces. The class
of operators we consider is composed by a linear drift plus a L\'evy operator that is comparable, in a suitable sense, with a
possibly truncated stable operator. It includes for example, the relativistic, the tempered, the layered or the Lamperti
stable operators. Our method does not assume neither the symmetry of the L\'evy operator nor the invariance for dilations
of the linear part of the operator. Thanks to our estimates, we prove in addition the well-posedness of the considered IPDE in  suitable
functional spaces.
In the final section, we extend some of these results to more general operators involving non-linear, space-time dependent
drifts.

\section{Introduction}
\fancyhead[RO]{Section \thesection. Introduction}
Fixed an integer $N$ in $\N$, we consider the following integro-partial differential operator of Ornstein-Uhlenbeck type:
\begin{equation}
\label{Degenerate_Stable_PDE1}
  L^{\text{ou}} \, := \, \mathcal{L}+ \langle A x , D_x\rangle \quad \text{ on } \R^N,
\end{equation}
where  $\langle \cdot,\cdot \rangle$ denotes the Euclidean inner product on $\R^N$, $A$ is a matrix in $\R^N\otimes \R^N$ and $\mathcal{L}$ is a
possibly degenerate, L\'evy operator acting non-degenerately only on a subspace of $\R^N$. We are interested in showing the
\emph{well-posedness} and the associated \emph{Schauder estimates} for elliptic and parabolic equations involving the operator
$L^{\text{ou}}$ and with coefficients in a generalized family of H\"older spaces. \newline
We only assume that $A$ satisfies a natural controllability assumption, the so-called Kalman rank condition (condition
[\textbf{K}] below), and that the operator $\mathcal{L}$ is comparable, in a suitable sense, to a non-degenerate,
truncated $\alpha$-stable operator on the same subspace of $\R^N$, for some $\alpha<2$ (condition [\textbf{SD}]
below).\newline
The topic of Schauder estimates for Ornstein-Uhlenbeck operators has been widely studied in the last decades, especially in the diffusive,
local setting, i.e.\ when $\mathcal{L}=\frac{1}{2}\text{Tr}\bigl(Q D^2_x\bigr)$ for some suitable matrix $Q$, and it is now
quite well-understood. See e.g.\ \cite{book:Gilbarg:Trudinger01}. \newline
On the other hand, a literature on the topic for the pure jump, non-local framework has been developed only in the recent years
(\cite{Bass09}, \cite{Dong:Kim13}, \cite{Bae:Kassmann15}, \cite{Ros-Oton:Serra16}), \cite{Fernanadez:Ros-Oton17},
\cite{Imbert:Jin:Shvydkoy18}, \cite{Chaudru:Menozzi:Priola19}, \cite{Kuhn19},  but mainly
in the non-degenerate, $\alpha$-stable setting, i.e.\ when $\mathcal{L}=\Delta^{\alpha/2}_x$ is the fractional Laplacian on $\R^N$ or
similar. To the best of our knowledge, the only two articles dealing with the degenerate, non-local framework (if
$\mathcal{L}=\Delta^{\alpha/2}_x$ acts non-degenerately only on a sub-space of $\R^N$) are \cite{Hao:Peng:Zhang19}, that takes into account the kinetics dynamics ($N=2d$), and \cite{Marino20}, for the general chain. In
order to use \cite{Hao:Peng:Zhang19} or \cite{Marino20} for our operator \eqref{Degenerate_Stable_PDE1}, we would need to impose the
additional strong assumption of invariance for dilations of the matrix $A$.\newline
The analysis of Ornstein-Uhlenbeck operators has been mainly developed following two different approaches. On the one hand, Da Prato and Lunardi in 
\cite{Daprato:Lunardi95} have been the first to use a priori estimates for the corresponding semi-group between suitable function spaces (See also \cite{Lunardi97,Lorenzi05,Chaudru:Honore:Menozzi18_Sharp, priola18}). Such a semi-group approach only adresses the regularity in space and indeed, the associated anisotropic 
H\"older spaces and Schauder estimates reflect this fact. In particular, the parabolic Schauder estimates do not present a bootstrap effect with respect to the 
initial condition.\newline
The second approach, introduced by Manfredini in \cite{Manfredini97}, exploits instead the general analysis on Lie groups to construct intrinsic H\"older spaces 
(see \cite{Pagliarani:Pascucci:Pignotti16} for a definition) that takes into account the joint space-time regularity of the involved functions. For a more thorough explanation along this direction, we suggest the interested reader to see, for example, \cite{Pascucci03}, \cite{DiFrancesco:Polidoro06} or the recent paper \cite{Imbert:Mouhot20}.

Even if the Ornstein-Uhlenbeck operator is usually exploited as a "toy model" for more general operators with space-time dependent, non-linear coefficients, we highlight that they appear naturally in various scientific contexts: for example in physics, for the
analysis of anomalous diffusions phenomena or for Hamiltonian models in a turbulent regime (see e.g.
\cite{Baeumer:Benson:Meerschaert01}, \cite{Cushman:Park:Kleinfelter:Moroni05} and the references therein) or in mathematical finance and
econometrics (see e.g.\ \cite{Brockwell01}, \cite{Barndorff-Nielsen:Shephard01}). The interest in Schauder estimates involving this type of operator also follows from the natural application which consists in establishing
the well-posedness of stochastic differential equations (SDE) driven by L\'evy processes and the associated stochastic control theory. See
e.g.\ \cite{Fleming:Mitter82}, \cite{Chaudru:Menozzi17}, \cite{Hao:Wu:Zhang19}.

Under our assumptions, we have been able to consider more general L\'evy operators not usually included in the literature, such as the
relativistic stable process, the layered stable process or the Lamperti one (see Paragraph "Main Operators Considered" below for details).
Moreover, we do not require the operator $\mathcal{L}$ to be symmetric. Here, we only mention one important example that satisfies our hypothesis, the Ornstein-Uhlenbeck operator on $\R^2$ driven by the
relativistic fractional Laplacian $\Delta^{\alpha/2}_{\text{rel}}$ and acting only on the first component:
\begin{multline}\label{eq:example_new}
x_1(D_{x_1}\phi(x)+D_{x_2}\phi(x)) + \text{p.v.}\int_{\R}\bigl[\phi(\begin{pmatrix} x_1+z \\ x_2 \end{pmatrix})-\phi(\begin{pmatrix} x_1
\\ x_2 \end{pmatrix})\bigr]\frac{1+\vert z \vert^{\frac{d+\alpha -1}{2}}}{\vert z \vert^{d+\alpha}}e^{-\vert z \vert}\, dz \\
= \, \langle A
x , D_x\phi(x) \rangle + \mathcal{L}\phi(x)
\end{multline}
where $x=(x_1,x_2)$ in $\R^2$. Such an example is included in the framework of Equation \eqref{Degenerate_Stable_PDE1}
considering $A=\begin{pmatrix}1 & 0 \\ 1 & 0 \end{pmatrix}$. This operator appears naturally as a fractional
generalization of the relativistic Schr\"odinger operator (See \cite{Ryznar02} for more details).\newline
We remark that example \eqref{eq:example_new} cannot be considered in \cite{Hao:Peng:Zhang19} or in our previous work
\cite{Marino20}. Indeed, the matrix
$A_0$ is not "dilation-invariant" (see example \ref{example_basic} below) and thus, it cannot be rewritten in the form used in
\cite{Marino20} (see also \cite{Lanconelli:Polidoro94} Proposition $2.2$ for a more thorough explanation). Furthermore, operators like the
relativistic fractional Laplacian cannot be treated  in \cite{Hao:Peng:Zhang19} or \cite{Marino20} that indeed have taken into account only
stable-like operators on $\R^N$. Another useful advantage of our technique is that we do not need anymore the symmetry of
the L\'evy measure $\nu$ which was, again, a key assumption in \cite{Marino20}.

More in details, given an integer $d\le N$ and a matrix $B$ in $\R^N\otimes \R^d$ such that $\text{rank}(B)=d$, we consider a
family of  operators $\mathcal{L}$ that can be represented for any sufficiently regular function $\phi\colon \R^N\to \R$ as
\begin{multline}\label{eq:def_operator_L}
\mathcal{L}\phi(x)\, := \, \frac{1}{2}\text{Tr}\bigl(BQB^\ast D^2_x\phi(x)\bigr) +\langle Bb,D_x \phi(x)\rangle \\
+ \int_{\R^d_0}\bigl[\phi(x+Bz)-\phi(x)-\langle D_x\phi(x), Bz\rangle\mathds{1}_{B(0,1)}(z) \bigr] \,\nu(dz),
\end{multline}
where $b$ is a vector in $\R^d$, $Q$ is a symmetric, non-negative definite matrix in $\R^d\otimes \R^d$ and $\nu$ is a L\'evy measure on
$\R^d_0:=\R^d\smallsetminus \{0\}$, i.e.\ a $\sigma$-finite measure on $\mathcal{B}(\R^d_0)$, the Borel $\sigma$-algebra on $\R^d_0$, such
that $\int(1\wedge \vert z\vert^2) \, \nu(dz)$ is finite. We then suppose $\nu$ to satisfy the following \emph{stable domination} condition:
\begin{description}
  \item[{[SD]}] there exists $r_0>0$, $\alpha$ in $(0,2)$ and a finite, non-degenerate measure $\mu$ on the unit sphere $\mathbb{S}^{d-1}$
      such that
      \[\nu(\mathcal{A}) \, \ge \, \int_{0}^{r_0}\int_{\mathbb{S}^{d-1}} \mathds{1}_{\mathcal{A}}(r\theta)\, \mu(d\theta)\frac{dr}{r^{1+\alpha}}, \quad
      \mathcal{A}\in \mathcal{B}(\R^d_0).\]
\end{description}
We recall that a measure $\mu$ on $\R^d$ is non-degenerate if there exists a constant $\eta \ge 1$ such that
\begin{equation}\label{eq:non_deg_measure_Levy}
\eta^{-1}\vert p \vert^\alpha \, \le \, \int_{\mathbb{S}^{d-1}}\vert p\cdot s \vert^\alpha \, \mu(ds) \, \le\,\eta
\vert p \vert^\alpha, \quad p \in \R^d,
\end{equation}
where $"\cdot"$ stands for the inner product on the smaller space $\R^d$. Since any $\alpha$-stable L\'evy measure $\nu_\alpha$  can be
decomposed into a spherical part $\mu$ on
$\mathbb{S}^{d-1}$ and a radial part $r^{-(1+\alpha)}dr$ (see e.g.\ Theorem $14.3$ in \cite{book:Sato99}), assumption
[\textbf{SD}] roughly states that the L\'evy measure of the integro-differential part of $\mathcal{L}$ is bounded from below by the L\'evy
measure of a possibly truncated, $\alpha$-stable operator on $\R^d$. \newline
It is assumed moreover that the matrices $A$, $B$ satisfy the following \emph{Kalman condition}:
\begin{description}
  \item[{[K]}] It holds that $N \, = \, \text{rank}\bigl[B,AB,\dots,A^{N-1}B\bigr]$,
\end{description}
where $\bigr[B,AB,\dots,A^{N-1}B\bigr]$ is the matrix in $\R^N\otimes\R^{dN}$ whose columns are given by $B,AB,\dots,A^{N-1}B$. \newline
Such an assumption is equivalent, in the linear framework, to the  H\"ormander condition (see \cite{Hormander67}) on the
commutators, ensuring the hypoellipticity of the operator $\partial_t-L^{\text{ou}}$. Moreover, condition [\textbf{K}] is
well-known in control theory (see e.g.\ \cite{book:Zabczyk95}, \cite{Priola:Zabczyk09}).

\paragraph{Mathematical outline.} In the present paper, we aim at establishing global Schauder estimates for equations involving the operator
$L^{\text{ou}}$ on $\R^N$, both in the elliptic and parabolic settings. Namely, we consider for a fixed $\lambda>0$ the following
elliptic equation:
\begin{equation}\label{eq:Elliptic_IPDE}
\lambda u(x) - L^{\text{ou}}u(x) \, = \, g(x), \quad x \in \R^N,
\end{equation}
and, for a fixed time horizon $T>0$, the following parabolic Cauchy problem:
\begin{equation}\label{eq:Parabolic_IPDE}
\begin{cases}
  \partial_tu(t,x) \, = \, L^{\text{ou}}u(t,x)+f(t,x), \quad (t,x) \in (0,T)\times \R^N; \\
  u(0,x) \, = \, u_0(x), \quad x \in \R^N,
\end{cases}
\end{equation}
where $f,g,u_0$ are given functions. Since our aim is to show optimal regularity results in H\"older spaces, we will assume for the elliptic
case (Equation \eqref{eq:Elliptic_IPDE}) that the source $g$ belongs to a suitable \emph{anisotropic}
H\"older space $C^\beta_{b,d}(\R^N)$ for some $\beta$ in $(0,1)$, where the H\"older exponent depends on the "direction" considered. The space
$C^\beta_{b,d}(\R^N)$ can be understood as composed by the bounded functions on $\R^N$
that are H\"older continuous with respect to a distance $\mathbf{d}$ somehow induced by the operator $L^{\text{ou}}$. We refer to
Section $2$ for a detailed exposition of such an argument but we
highlight already that the above mentioned distance $\mathbf{d}$ can be seen as a generalization of the classical parabolic distance, adapted to our
degenerate, non-local framework. It is precisely assumption [\textbf{K}], or equivalently the hypoellipticity of
$\partial_t+L^{\text{ou}}$, that ensures the existence of such a distance $\mathbf{d}$ and gives it its anisotropic nature. Roughly speaking, it
allows the smoothing effect of the L\'evy operator $\mathcal{L}$ acting non-degenerately only on some components, say $B\R^N$, to spread in
the whole space $\R^N$, even if with lower regularizing properties. \newline
Concerning the parabolic problem \eqref{eq:Parabolic_IPDE}, we assume similarly that $u_0$ is in $C^{\alpha+\beta}_{b,d}(\R^N)$ and that
$f(t,\cdot)$ is in $C^{\beta}_{b,d}(\R^N)$, uniformly in $t\in(0,T)$. The typical estimates we want to prove can be stated in the parabolic
setting in the following way: there exists a constant $C$, depending only on the parameters of the model, such that any
distributional solution $u$ of the Cauchy problem \eqref{eq:Parabolic_IPDE} satisfies
\begin{equation}
\label{eq:SchauderEstimatesintro} \tag{$\mathscr{S}$}
\Vert u \Vert_{L^\infty(C^{\alpha+\beta}_{b,d})} \, \le \, C\bigl[\Vert u_0 \Vert_{C^{\alpha+\beta}_{b,d}}+\Vert f
\Vert_{L^\infty(C^{\beta}_{b,d})} \bigr].
\end{equation}
As a by-product of the Schauder Estimates $(\mathscr{S})$, we will obtain the well-posedness of the
Cauchy problem \eqref{eq:Parabolic_IPDE} in the space $L^\infty\bigl(0,T;C^{\alpha+\beta}_{b,d}(\R^N)\bigr)$ , once the existence of a
solution is established. The additional regularity for the solution $u$ with respect to the source $f$ reflects the appearance of a
smoothing effect associated with $L^{\text{ou}}$ of order $\alpha$, as it is expected by condition [\textbf{SD}]. It can be seen
as a generalization of the "standard" parabolic bootstrap to our degenerate, non-local setting. We highlight that the parabolic
bootstrap in ($\mathscr{S}$) is precisely derived from the non-degenerate stable-like part in $\mathcal{L}$ (lowest regularizing effect in
the operator).  \newline
To show our result, we will follow the semi-group approach as firstly introduced in \cite{Daprato:Lunardi95}, which became
afterwards a very robust tool to study Schauder estimates in a wide variety of frameworks (\cite{Lunardi97}, \cite{Lorenzi05}, \cite{Saintier07}, \cite{Priola09}, \cite{Priola12}, \cite{Dong:Kim13}, \cite{Kim:Kim15}, \cite{Chaudru:Honore:Menozzi18_Sharp}, \cite{Kuhn19}). The main idea is to consider the Markov transition semi-group $P_t$ associated with
$L^{\text{ou}}$ and then, in the elliptic case, to use the Laplace transform
formula in order to represent the unique distributional solution $u$ of Equation \eqref{eq:Elliptic_IPDE} as:
\[u(x) \, = \, \int_{0}^{\infty}e^{-\lambda t}\bigl[P_tg\bigr](x) \, dt \, =: \, \int_{0}^{\infty}e^{-\lambda t}P_tg(x) \, dt.\]
In the parabolic setting, we exploit instead the variation of constants (or Duhamel) formula in order to show a similar representation for
the weak solution of the Cauchy problem \eqref{eq:Parabolic_IPDE}:
\[u(t,x) \, = \, P_tu_0(x) +\int_{0}^{t}\bigl[P_{t-s}f(s,\cdot)\bigr](x)\, ds \, =: \, P_tu_0(x) +\int_{0}^{t}P_{t-s}f(s,x)\, ds.\]
In order to prove global regularity estimates for the solutions, the crucial point is to understand the
action of the operator $P_t$ on the anisotropic H\"older spaces. In particular, we will show in Corollary
\ref{corollary:continuity_between_holder} the continuity of $P_t$
as an operator from $C^\beta_{b,d}(\R^N)$ to $C^\gamma_{b,d}(\R^N)$ for $\beta<\gamma$ and, more precisely, that it holds:
\begin{equation}\label{eq:Intro:Gradient_Estimates}
\Vert P_t \phi\Vert_{C^\gamma_{b,d}} \, \le \,C\Vert \phi\Vert_{C^{\beta}_{b,d}}\Bigl(1+t^{-\frac{\gamma-\beta}{\alpha}}\Bigr), \quad t>0.
\end{equation}
The above estimate can be obtained through interpolation techniques (see Equation \eqref{eq:Interpolation_ineq}), once
sharp controls in supremum norm (Theorem \ref{prop:Control_in_Holder_norm} below) are established for the spatial derivatives of $P_t\phi$
when $\phi \in C^\beta_{b,d}(\R^N)$. We think that such an estimate \eqref{eq:Intro:Gradient_Estimates} and the controls in
Theorem \ref{prop:Control_in_Holder_norm} can be of independent interest and used also beyond our scope in other contexts. \newline
We face here two main difficulties to overcome. While in the gaussian setting, $L^\infty$-estimates of this type have been
established exploiting, for example, explicit formulas for the density of the semi-group $P_t$ (\cite{Lunardi97}), a priori controls of
Bernstein type combined with interpolation methods (\cite{Lorenzi05} and \cite{Saintier07}, when $n=2$ in \eqref{eq:def_of_n_Levy} below) or
probabilistic representations of the semi-group $P_t$, allowing Malliavin calculus (\cite{Priola09}), we cannot rely on
these techniques in our non-local framework, mainly due to the lower integrability properties for $P_t$. Instead,
we are going to use a \emph{perturbative approach} which consists in considering the L\'evy operator $\mathcal{L}$ as a perturbation, in a
suitable sense, of an $\alpha$-stable operator, at least for the associated small jumps. Indeed, we can "decompose" the operator
$\mathcal{L}$ in a smoother part, $\mathcal{L}^\alpha$, whose L\'evy measure is given by
\[\mu(d\theta)\frac{\mathds{1}_{(0,r_0]}(r)}{r^{1+\alpha}}dr\]
and a remainder part. It is precisely condition [\textbf{SD}] that allows such a decomposition, since it ensures the positivity of the
L\'evy measure
\[d\nu-d\mu \frac{\mathds{1}_{[0,r_0]}}{r^{1+\alpha}}dr\]
associated with the remainder term. The main difference with the previous techniques in the diffusive setting is that we will work mainly on
the truncated $\alpha$-stable contribution $\mathcal{L}^\alpha$, being the remainder term only bounded. \newline
Following \cite{Schilling:Sztonyk:Wang12}, we will establish that the Hartman-Winter condition holds, ensuring the
existence of a smooth density for the semi-group associated with $\mathcal{L}^\alpha$ and then, the required gradient estimates.
Indeed, assumption [\textbf{SD}] roughly states that the small jump
contributions of $\nu$, the ones responsible for the creation of a density, are controlled from below by an $\alpha$-stable measure, whose
absolute continuity is well-known in our framework. \newline
On the other hand, we will have to deal with the degeneracy of the operator $\mathcal{L}$, that acts non-degenerately, through the embedding
matrix $B$, only on a subspace of dimension $d$. It will be managed by adapting the reasonings firstly appeared in \cite{Huang:Menozzi16}.
Namely, we will show that the semi-group associated with the Ornstein-Uhlenbeck operator $L^{\text{ou}}$ coincides with a
non-degenerate one but "multiplied" by a time-dependent matrix that precisely takes into account the original degeneracy of the operator
(see definition of matrix $\mathbb{M}_t$ in Section $2.1$).

\paragraph{Main operators considered.} We conclude this introduction showing that assumption [\textbf{SD}] applies to a large class of L\'evy
operators on $\R^d$. As already pointed out in \cite{Schilling:Sztonyk:Wang12}, it is satisfied by any L\'evy measure $\nu$ that can be
decomposed in polar coordinates as
\[\nu(\mathcal{A}) \, = \, \int_{0}^{\infty}\int_{\mathbb{S}^{d-1}}\mathds{1}_{\mathcal{A}}(r\theta)Q(r,\theta)\,\mu(d\theta) \frac{dr}{r^{1+\alpha}}, \quad \mathcal{A}
\in \mathcal{B}(R^d_0),\]
for a finite, non-degenerate (in the sense of Equation \eqref{eq:non_deg_measure_Levy}), measure $\mu$ on $\mathbb{S}^{d-1}$ and a Borel
function $Q\colon(0,\infty)\times \mathbb{S}^{d-1}\to \R$ such
that there exists $r_0>0$ so that
\[ Q(r,\theta) \, \ge  \, c>0, \quad \text{a.e. in }[0,r_0]\times \mathbb{S}^{d-1}.\]
In particular, assumption [\textbf{SD}] holds for the following families of "stable-like" examples with $\alpha \in (0,2)$:
\begin{enumerate}
  \item Stable operator \cite{book:Sato99}:
  \[Q(r,\theta) \, = \, 1;\]
  \item Truncated stable operator with $r_0>0$ \cite{Kim:Song08}:
  \[Q(r,\theta) \, = \, \mathds{1}_{(0,r_0]}(r);\]
  \item Layered stable operator with $\beta$ in $(0,2)$ and $r_0>0$ \cite{Houdre:Kawai07}:
  \[Q(r,\theta) \, = \, \mathds{1}_{(0,r_0]}(r)+\mathds{1}_{(r_0,\infty)}(r)r^{\alpha-\beta};\]
  \item Tempered stable operator \cite{Rosinski09}:
  \[Q(\cdot,\theta) \text{ completely monotone, $Q(0,\theta)>0$ and $Q(\infty,\theta)=0$ a.e.\ in }S^{d-1};\]
  \item Relativistic stable operator \cite{Carmona:Masters:Simon90}, \cite{Byczkowski:Malecki:Ryznar09}:
  \[Q(r,\theta) \, = \, (1+r)^{(d+\alpha-1)/2}e^{-r};\]
  \item Lamperti stable operator with $f\colon S^{d-1}\to \R$ such that $\sup f(\theta)<1+\alpha$ \cite{Caballero:Pardo:Perez10}:
  \[Q(r,\theta) \, = \, e^{rf(\theta)}\Bigl(\frac{r}{e^r-1}\Bigr)^{1+\alpha}.\]
\end{enumerate}

\setcounter{equation}{0}
\paragraph{Organization of the paper.} The article is organized as follows. Section $2$ introduces some useful notations and then, the
anisotropic distance $\mathbf{d}$ induced by the dynamics as well as Zygmund-H\"older spaces associated with such a distance. In Section $3$, we are
going to show some analytical properties of the semi-group $P_t$ generated by $L^{\text{ou}}$, such as the existence of a
smooth density and, at least for small times, some controls for its derivatives. Section $4$ is then
dedicated to different estimates in the $L^\infty$-norm for $P_tf$ and its spatial derivatives, involving the supremum or the H\"older norm
of the function $f$. In particular, we show here the
continuity of $P_t$ as an operator between anisotropic Zygmund-H\"older spaces. In Section $5$, we use the controls established in the
previous parts in order to prove the elliptic Schauder estimates and show that Equation \eqref{eq:Elliptic_IPDE} has a unique solution.
Similarly, we establish the well-posedness of the Cauchy problem \eqref{eq:Parabolic_IPDE} as well as the associated parabolic
Schauder estimates. In the final section of the article, we briefly explain some possible extensions of the previous results to non-linear,
space-time dependent operators.

\section{Geometry of the dynamics}
\fancyhead[RO]{Section \thesection. Geometry of the dynamics}
In this section, we are going to choose the right functional space "in which" to state our Schauder
estimates. The idea is to construct an H\"older space $C^\beta_{b,d}(\R^N)$ with respect to a distance $\mathbf{d}$ that it is
homogeneous to the dynamics, i.e. such that for any $f$ in $C^\beta_{b,d}(\R^N)$, any
distributional solution $u$ of
\begin{equation}\label{eq:distrib.sol}
L^{\text{ou}}u(x) \, = \, \mathcal{L}u(x)+\langle Ax,Du(x)\rangle \, = \, f(x), \quad x \in \R^N
\end{equation}
is in $C^{\alpha+\beta}_{b,d}(\R^N)$, the expected parabolic bootstrap associated to this kind of operator. We recall in particular that the Kalman rank condition [\textbf{K}] is equivalent to the hypoellipticity (in the sense of H\"ormander \cite{Hormander67}) of the operator $L^{\text{ou}}$ that ensures the existence and smoothness of a distributional solution of Equation \eqref{eq:distrib.sol} for sufficiently regular $f$. See e.g.\ \cite{book:Ishikawa16} or \cite{Hao:Peng:Zhang19} for more details.

\subsection{The distance associated with the dynamics}

To construct the suitable distance $\mathbf{d}$, we start noticing that the Kalman rank condition [\textbf{K}] allows us to denote
\begin{equation}\label{eq:def_of_n_Levy}
n \,:= \, \min\{r \in \N \colon N \, = \, \rank \bigl[B,AB,\dots,A^{r-1}B\bigr]\}.
\end{equation}
Clearly, $n$ is in $\llbracket 1,N\rrbracket$, where $\llbracket \cdot, \cdot \rrbracket$ denotes the set of all the integers in the
interval, and $n=1$ if and only if $d=N$, i.e. if the dynamics is non-degenerate.\newline
As done in \cite{Lunardi97}, the space $\R^N$ will be decomposed with respect to the family of linear operators $B, AB,\dots,
A^{n-1}B$. We start defining the family $\{V_h\colon h\in \llbracket 1,n \rrbracket\}$ of subspaces of $\R^N$ through
\[
V_h \,  := \, \begin{cases}
            \text{Im} (B), & \mbox{if } h=1, \\
            \bigoplus_{k=1}^{h}\text{Im}(A^{k-1}B), & \mbox{otherwise}.
        \end{cases}\]
It is easy to notice that $V_h\neq V_k$ if $k\neq h$ and $V_1\subset V_2\subset\dots V_n=\R^N$. We can then construct iteratively the family
$\{E_h \colon h\in \llbracket 1,n \rrbracket\}$ of orthogonal projections from $\R^N$ as
\[E_h \,  := \,
        \begin{cases}
            \text{projection on } V_1, & \mbox{if } h=1; \\
            \text{projection on }(V_{h-1})^\perp \cap V_h, & \mbox{otherwise}.
        \end{cases}
\]
With a small abuse of notation, we will identify the projection operators $E_h$ with the corresponding matrices in
$\R^N\otimes \R^N$. It is clear that $\dim E_1(\R^N)=d$. Let us then denote $d_1:=d$ and
\[d_h \,:= \, \dim E_h(\R^N), \quad \text{ for $h$}>1.\]

We can define now the distance $\mathbf{d}$ through the decomposition $\R^N=\bigoplus_{h=1}^nE_h(\R^N)$ as
\[\mathbf{d}(x,x') \,:= \, \sum_{h=1}^{n}\vert E_h(x-x')\vert^{\frac{1}{1+\alpha(h-1)}}.\]
The above distance can be seen as a generalization of the usual Euclidean distance when $n=1$ (non-degenerate dynamics) as well as an extension of the standard parabolic distance for $\alpha=2$. It is important to highlight that it does
not induce a norm since it lacks of linear homogeneity. \newline
The anisotropic distance $\mathbf{d}$ can be understood direction-wise: we firstly fix a "direction" $h$ in $\llbracket
1,n\rrbracket$ and then calculate the standard Euclidean distance on the associated subspace $E_h(\R^N)$, but scaled according to the
dilation of the system in that direction. We conclude summing the contributions associated with each component. The choice of such a
dilation will be discussed thoroughly in the example at the end of this section.

As emphasized by the result from Lanconelli and Polidoro recalled below (cf. \cite{Lanconelli:Polidoro94}, Proposition $2.1$), the decomposition of
$\R^N$ with respect to the projections $\{E_h \colon h \in \llbracket 1,n\rrbracket\}$ determines a particular structure of the matrices
$A$ and $B$. It will be often exploited in the following.

\begin{theorem}[\cite{Lanconelli:Polidoro94}]
\label{Thm:Lancon_Pol_Levy}
Let $\{e_i \colon i \in \llbracket 1,N\rrbracket\}$ be an orthonormal basis consisting of generators of $\{E_h(\R^N) \colon h\in \llbracket
1,n\rrbracket\}$. Then, the matrices $A$ and $B$ have the following form:
\begin{equation}\label{eq:Lancon_Pol_Levy}
B \, = \,
    \begin{pmatrix}
        B_0    \\
        0      \\
        \vdots  \\
        0
    \end{pmatrix}
\,\, \text{ and } \,\, A \, = \,
    \begin{pmatrix}
        \ast   & \ast  & \dots  & \dots  & \ast   \\
         A_2   & \ast  & \ddots & \ddots  & \vdots   \\
        0      & A_3   & \ast  & \ddots & \vdots \\
        \vdots &\ddots & \ddots& \ddots & \ast \\
        0      & \dots & 0     & A_n    & \ast
    \end{pmatrix}
\end{equation}
where $B_0$ is a non-degenerate matrix in $\R^{d_1}\otimes \R^{d_1}$ and $A_h$ are matrices in $\R^{d_h}\otimes \R^{d_{h-1}}$ with
$\text{rank}(A_h)=d_h$ for any $h$ in $\llbracket 2,n\rrbracket$. Moreover, $d_1\ge d_2\ge \dots\ge d_n\ge 1$.
\end{theorem}
Applying a change of variables if necessary, we will assume from this point further to have fixed such a canonical basis $\{e_i \colon i \in
\llbracket1,N\rrbracket\}$. For notational simplicity, we denote by $I_h$, $h \in \llbracket 1,n\rrbracket$, the family of indexes $i$ in
$\llbracket 1,N\rrbracket$ such that $\{e_i\colon i\in I_h\}$ spans $E_h(\R^N)$.

The particular structure of $A$ and $B$ given by Theorem \ref{Thm:Lancon_Pol_Levy} allows us to decompose accurately the
exponential $e^{tA}$ of the matrix $A$ in order to make the intrinsic scale of the system appear. Further on, we will consider fixed a time-dependent matrix $\mathbb{M}_t$ on $\R^N\otimes \R^N$ given by
\[\mathbb{M}_t := \text{diag}(I_{d_1\times d_1},tI_{d_2\times d_2},\dots,t^{n-1}I_{d_{n}\times d_{n}}), \quad t\ge0.\]

\begin{lemma}
\label{lemma:Decomposition_exp_A}
There exists a time-dependent matrix $\{R_t\colon t\in[0,1]\}$ in $\R^N\otimes \R^N$ such that
 \begin{equation}\label{eq:Decomposition_exp_A}
  e^{tA}\mathbb{M}_t \,   = \mathbb{M}_tR_t, \quad t \in [0,1].
 \end{equation}
Moreover, there exists a constant $C>0$ such that for any $t$ in $[0,1]$,
\begin{itemize}
  \item any $l,h$ in $\llbracket 1,n \rrbracket$ and any $\theta$ in $\mathbb{S}^{N-1}$, it holds that
  \[\bigl{\vert} E_le^{tA}E_h \theta\bigr{\vert} \, \le \, \begin{cases}
                                                      Ct^{l-h}, & \mbox{if } l\ge h \\
                                                      Ct, & \mbox{if } l< h.
                                                    \end{cases}\]
  \item any $\theta$ in $\mathbb{S}^{d-1}$, it holds that
  \[\bigl{\vert} R_tB\theta \bigr{\vert} \, \ge \, C^{-1}.\]
\end{itemize}
\end{lemma}
\begin{proof}
By definition of the matrix exponential, we know that
\begin{equation}\label{Proof:eq:def_Exponential}
E_le^{tA}E_h  \, = \, \sum_{k=0}^{\infty}\frac{t^k}{k!}E_lA^kE_h.
\end{equation}
Using now the representation of $A$ given by Theorem \ref{Thm:Lancon_Pol_Levy}, it is easy to check that $E_lA^kE_h=0$ for $k<l-h$  (when $l-h$
is non-negative). Thus, for $l\ge h$, it holds that
\[\bigl{\vert}E_le^{tA}E_h\theta\bigr{\vert}  \, = \, \bigl{\vert}\sum_{k=l-h}^{\infty}\frac{t^k}{k!}E_lA^kE_h\theta\bigr{\vert} \,
\le \, Ct^{l-h},\]
where we exploited that $t$ is in $ [0,1]$ and $\vert \theta \vert =1$. Assuming instead that $l<h$, it is clear that $E_lI_{N\times N}E_h$
vanishes. We can then write that
\[\bigl{\vert}E_le^{tA}E_h\theta\bigr{\vert}  \, = \, \bigl{\vert}\sum_{k=l}^{\infty}\frac{t^k}{k!}E_lA^kE_h\theta\bigr{\vert} \, \le \,
Ct,\]
using again that $t$ is in $ [0,1]$ and $\vert \theta \vert =1$.\newline
To show the other control, we highlight that the matrix $\mathbb{M}_t$ is not invertible in $t=0$ and for this reason, we
define the time-dependent matrix $R_t$ as
\[R_t \, := \, \begin{cases}
                 I_{N\times N}, & \mbox{if } t=0; \\
                 \mathbb{M}^{-1}_te^{tA}\mathbb{M}_t, & \mbox{if } t \in (0,1].
               \end{cases}\]
We could have also defined $R_t:=\bigr(\tilde{R}^t_s\bigr)_{|s=1}$ where $\tilde{R}^t_s$ solves the following ODE:
\[\begin{cases}
    \partial_s\tilde{R}^t_s \, = \, \mathbb{M}^{-1}_ttA\mathbb{M}_t\tilde{R}^t_s, \quad \text{on }(0,1], \\
    \tilde{R}^t_0 \, = \, I_{N\times N}.
  \end{cases}\]
Equivalently, $\tilde{R}^t_s$ is the resolvent matrix associated with $\mathbb{M}^{-1}_ttA\mathbb{M}_t$, whose sub-diagonal entries are
"macroscopic" from the structure of $A$ and $\mathbb{M}_t$. \newline
It follows immediately that Equation \eqref{eq:Decomposition_exp_A} holds.  Moreover, we notice that
\[\bigl{\vert} R_tB\theta\bigr{\vert} \, \ge \,\bigl{\vert} E_1R_tB\theta\bigr{\vert} \, = \, \bigl{\vert} E_1e^{tA}E_1B\theta\bigr{\vert}.\]
Remembering the definition of matrix exponential (Equation \eqref{Proof:eq:def_Exponential} with $l=h=1$), we use now that
\[E_1A^kE_1 \, = \, (E_1AE_1)^k \, = \, (A_{1,1})^kE_1,\]
where in the last expression the multiplication is meant block-wise, in order to conclude that
\[\bigl{\vert} R_t B\theta \bigr{\vert} \, \ge \,\bigl{\vert} e^{tA_{1,1}}B_0\theta\bigr{\vert}.\]
Using that $e^{tA_{1,1}}B_0$ is non-degenerate and continuous in time and that $\theta$ is in $\mathbb{S}^{d-1}$, it is easy to conclude.
\end{proof}

We conclude this sub-section with a simpler example taken from \cite{Huang:Menozzi:Priola19}. We hope that it will help the reader to
understand the introduction of the anisotropic distance $\mathbf{d}$.
\begin{example}\label{example_basic}
Fixed $N=2d$, $n=2$ and $d=d_1=d_2$, we consider the following operator:
\[L^{\text{ou}}_\alpha \, = \, \Delta^{\frac{\alpha}{2}}_{x_1}+x_1\cdot\nabla_{x_2} \quad \text{ on }\R^{2d},\]
where $(x_1,x_2)\in\R^{2d}$ and $\Delta^{\frac{\alpha}{2}}_{x_1}$ is the fractional Laplacian with respect to $x_1$. In our framework, it is associated with the matrices
\[A \, := \,
    \begin{pmatrix}
               0 & 0 \\
               I_{d\times d} & 0
    \end{pmatrix} \,\, \text{ and } \,\, B \, := \,
    \begin{pmatrix}
                I_{d\times d} \\
                 0
    \end{pmatrix}.
             \]
The operator $L^{\text{ou}}_\alpha$ can be seen as a generalization of the classical Kolmogorov example (see e.g. \cite{Kolmogorov34}) to
our non-local setting. \newline
In order to understand how the system typically behaves, we search for a dilation
\[\delta_\lambda\colon [0,\infty)\times \R^{2d} \to [0,\infty)\times \R^{2d}\]
which is invariant for the considered dynamics, i.e.\ a dilation that transforms solutions of the equation
\[\partial_tu(t,x)-L^{\text{ou}}_\alpha u(t,x) \, = \, 0 \quad \text{ on }(0,\infty)\times\R^{2d}\]
into other solutions of the same equation.\newline
Due to the structure of $A$ and the $\alpha$-stability of $\Delta^{\frac{\alpha}{2}}$, we can consider for any fixed $\lambda>0$, the
following
\[ \delta_\lambda(t,x_1,x_2) := (\lambda^\alpha t,\lambda x_1,\lambda^{1+\alpha}x_2).\]
It then holds that
\[\bigl(\partial_t -L^{\text{ou}}_\alpha\bigr) u = 0 \, \Longrightarrow \bigl(\partial_t -L^{\text{ou}}_\alpha \bigr)(u \circ
\delta_\lambda) = 0.\]
Introducing now the complete time-space distance $\mathbf{d}_P$ on $[0,\infty)\times \R^{2d}$ given by
\begin{equation}\label{Definition_distance_d_P1}
\mathbf{d}_P\bigl((t,x),(s,x')\bigr)  \, := \, \vert s-t\vert^\frac{1}{\alpha}+ \mathbf{d}(x,x') \, = \, \vert s-t\vert^\frac{1}{\alpha}+ \vert x_1-x'_1\vert +\vert x_2-x'_2\vert^{\frac{1}{1+\alpha}},
\end{equation}
we notice that it is homogeneous with respect to the dilation $\delta_\lambda$, so that
\[\mathbf{d}_P\bigl(\delta_\lambda(t,x);\delta_\lambda(s,x')\bigr) = \lambda \mathbf{d}_P\bigl((t,x);(s,x')\bigr).\]
Precisely, the exponents appearing in Equation \eqref{Definition_distance_d_P1} are those which make each space-component homogeneous to
the characteristic time scale $t^{1/\alpha}$. From a more probabilistic point of view, the exponents in Equation
\eqref{Definition_distance_d_P1}, can be related to the characteristic time scales of the iterated integrals of an $\alpha$-stable process.
It can be easily seen from the example, noticing that the operator $L^{\text{ou}}_\alpha$ corresponds to the generator of an isotropic $\alpha$-stable
process and its time integral. \newline
Going back to the general setting, the appearance of this kind of phenomena is due essentially to the particular structure of the
matrix $A$ (cf. Theorem \ref{Thm:Lancon_Pol_Levy}) that allows the smoothing effect of the operator $\mathcal{L}$, acting only on the first "component" given by $B_0$, to propagate into the system.
\end{example}

\subsection{Anisotropic Zygmund-H\"older spaces}
We are now ready to define the Zygmund-H\"older spaces $C^{\gamma}_{b,d}(\R^N)$ with respect to the distance $\mathbf{d}$. We start recalling some useful
notations we will need below.\newline
Given a function $f\colon \R^N\to \R$, we denote by $Df(x)$, $D^2f(x)$ and $D^3f(x)$ the first, second and third Fr\'echet derivative of $f$
at a point $x$ in $\R^N$ respectively, when they exist. For simplicity, we will identify $D^3f(x)$ as a $3$-tensor so that $[D^3f(x)](u,v)$
is a vector in $R^N$ for any $u,v$ in $\R^N$. Moreover, fixed $h$ in $\llbracket 1, n\rrbracket$, we will denote by $D_{E_h}f(x)$ the
gradient of $f$ at $x$ along the direction $E_h(\R^N)$. Namely,
\[D_{E_h}f(x) \, \, := \, \, E_h Df(x).\]
A similar notation will be used for the higher derivatives, too.\newline
Given $X,Y$ two real Banach spaces, $\mathcal{L}(X,Y)$ will represent the family of linear continuous operators between $X$ and $Y$. \newline
In the following, $c$ or $C$ denote generic \emph{positive} constants whose precise value is unimportant. They may change from line
to line and they will depend only on the parameters given by the model and assumptions [\textbf{SD}], [\textbf{K}]. Namely,
$d,N,A,B,\alpha,\nu,r_0$ and $\mu$. Other dependencies that may occur will be explicitly specified.

Let us introduce now some function spaces we are going to use. We denote by $B_b(\R^N)$ the family of Borel measurable and bounded
functions $f\colon \R^N\to \R$. It is a Banach space endowed with the supremum norm $\Vert \cdot \Vert_\infty$. We will consider also its
closed subspace $C_b(\R^N)$ consisting of all the uniformly continuous functions. \newline
Fixed some $k$ in $\N_0:=\N\cup\{0\}$ and $\beta$ in $(0,1]$, we follow Lunardi \cite{Lunardi97} denoting the Zygmund-H\"older semi-norm for a function $\phi\colon \R^N\to \R$ as
\[[\phi]_{C^{k+\beta}} \, := \,
\begin{cases}
    \sup_{\vert \vartheta \vert= k}\sup_{x\neq y}\frac{\vert D^\vartheta\phi(x)-D^\vartheta\phi(y)\vert}{\vert x-y\vert^\beta} , & \mbox{if
    }\beta \neq 1; \\
    \sup_{\vert \vartheta \vert= k}\sup_{x\neq y}\frac{\bigl{\vert}D^\vartheta\phi(x)+D^\vartheta\phi(y)-2D^\vartheta\phi(\frac{x+y}{2})
    \bigr{\vert}}{\vert x-y \vert}, & \mbox{if } \beta =1.
                             \end{cases}\]
Consequently, The Zygmund-H\"older space $C^{k+\beta}_b(\R^N)$ is the family of functions $\phi\colon \R^N
\to\R$ such that $\phi$ and its derivatives up to order $k$ are continuous and the norm
\[\Vert \phi \Vert_{C^{k+\beta}_b} \,:=\, \sum_{i=1}^{k}\sup_{\vert\vartheta\vert = i}\Vert D^\vartheta\phi
\Vert_{L^\infty}+[\phi]_{C^{k+\beta}_b} \,\text{ is finite.}\]

We can define now the anisotropic Zygmund-H\"older spaces associated with the distance $\mathbf{d}$. Fixed $\gamma>0$, the space $C^{\gamma}_{b,d}(\R^N)$ is
the family of functions $\phi\colon \R^N\to \R$ such that for any $h$ in $\llbracket 1,n\rrbracket$ and any $x_0$ in $\R^N$, the function
\[z\in  E_h(\R^N)\,  \to \, \phi(x_0+z) \in \R \,\text{ belongs to }C^{\gamma/(1+\alpha(h-1))}_b\bigl(E_h(\R^N)\bigr),\]
with a norm bounded by a constant independent from $x_0$. It is endowed with the norm
\begin{equation}\label{eq:def_anistotropic_norm_Levy}
\Vert\phi\Vert_{C^{\gamma}_{b,d}} \,:=\,\sum_{h=1}^{n}\sup_{x_0\in \R^N}\Vert\phi(x_0+\cdot)\Vert_{C^{\gamma/(1+\alpha(h-1))}_b}.
\end{equation}

We highlight that it is possible to recover the expected joint regularity for the partial derivatives, when
they exist, as in the standard H\"older spaces. In such a case, they actually turn out to be H\"older continuous with respect to the
distance $\mathbf{d}$ with order one less than the function (See Lemma $2.1$ in \cite{Lunardi97} for more details).

It will be convenient in the following to consider an equivalent norm in the "standard" H\"older-Zygmund spaces
$C^\gamma_b(E_h(\R^N))$ that does not take into account the derivatives with respect to the different directions.
We suggest the interested reader to see \cite{Lunardi97}, Equation $(2.2)$ or \cite{Priola09} Lemma $2.1$ for further details.

\begin{lemma}
\label{lemma:def_equivalent_norm}
Fixed $\gamma$ in $(0,3)$ and $h$ in $\llbracket 1,n\rrbracket$ and $\phi$ in $C_b(E_h(\R^N))$, let us introduce
\begin{equation}\label{eq:def_equivalent_norm}
\Delta^3_{x_0}\phi(z) \, := \, \phi(x_0+3z)-3\phi(x_0+2z)+ 3\phi(x_0+z)-\phi(x_0), \,\quad x_0 \in \R^N; \, z \in E_h(\R^N).
\end{equation}
 Then, $\phi$ is in $C^\gamma_{b}(E_h(\R^N))$ if and only if
\[\sup_{x_0 \in \R^N}\sup_{z \in E_h(\R^N);z\neq 0}\frac{\bigl{\vert} \Delta^3_{x_0}\phi(z)\bigr{\vert}}{\vert z \vert^\gamma} \, < \,
\infty.\]
\end{lemma}

We conclude this subsection with a result concerning the interpolation between the anisotropic
Zygmund-H\"older spaces $C^{\gamma}_{b,d}(\R^N)$. We refer to  Theorem $2.2$ and Corollary $2.3$ in \cite{Lunardi97} for details.

\begin{theorem}\label{Thm:Interpolation}
Let $r$ be in $(0,1)$ and $\beta$, $\gamma$ in $[0,\infty)$ such that $\beta\le \gamma$. Then, it holds that
\[\bigl(C^\beta_{b,d}(\R^N),C^\gamma_{b,d}(\R^N)\bigr)_{r,\infty} \, = \, C^{r\gamma+(1-r)\beta}_{b,d}(\R^N)\]
with equivalent norms, where we have denoted for simplicity: $C^0_{b,d}(\R^N):=C_b(\R^N)$.
\end{theorem}

\setcounter{equation}{0}
\mysection{Smoothing effect for truncated density}

We present here some analytical properties of the semi-group generated by the operator $L^{\text{ou}}$. Following
 \cite{Schilling:Sztonyk:Wang12} and \cite{Schilling:Wang12}, we will show the existence of a smooth density for such a
semi-group and its anisotropic smoothing effect, at least for small times.

Throughout this section, we consider fixed a stochastic base $\bigl(\Omega,\mathcal{F},(\mathcal{F}_t)_{t\ge 0},\mathbb{P}\bigr)$ satisfying
the usual assumptions (see
\cite{book:Applebaum09}, page $72$). Let us then consider the (unique in law) L\'evy process $\{Z_t\}_{t\ge 0}$ on $\R^d$ characterized by
the L\'evy symbol
\[\Phi(p) \, = \, -ib\cdot p +\frac{1}{2}p\cdot Qp + \int_{\R^d_0}\bigl(1-e^{ip \cdot z}+ip\cdot z\mathds{1}_{B(0,1)}(z)\bigr) \, \nu(dz),
\quad p \in \R^d.\]
It is well-known by the L\'evy-Kitchine formula (see
\cite{book:Jacob05}), that the infinitesimal generator of the process $\{BZ_t\}_{t\ge 0}$ is then given by $\mathcal{L}$ on $\R^N$.\newline
Fixed $x$ in $\R^N$, we denote by $\{X_t\}_{t\ge 0}$ the $N$-dimensional Ornstein-Uhlenbeck process driven by $BZ_t$, i.e.\ the unique
(strong) solution of the following stochastic differential equation:
\[X_t \, = \, x + \int_{0}^{t}AX_s \, ds +BZ_t, \quad t\ge0, \,\, \mathbb{P}\text{-almost surely.}\]
By the variation of constants method, it is easy to check that
\begin{equation}\label{eq:Def_OU_Process}
X_t \, = \, e^{tA}x + \int_{0}^{t}e^{(t-s)A}B\, dZ_s, \quad \quad t\ge0, \,\, \mathbb{P}\text{-almost surely.}
\end{equation}
The \emph{transition semi-group} associated with $L^{\text{ou}}$ is then defined as the family
$\{P_t\colon t \ge 0\}$ of linear contractions on $B_b(\R^N)$ given by
\begin{equation}\label{eq:Def_transition_semi-group}
P_t\phi(x) \, = \, \mathbb{E}\bigl[\phi(X_t)\bigr], \quad x \in \R^N,\, \phi \in B_b(\R^N).
\end{equation}
We recall that $P_t$ is generated by $L^{\text{ou}}$ in the sense that its infinitesimal generator
$\mathcal{A}$ coincides with $L^{\text{ou}}$ on $C^\infty_c(\R^N)$, the family of smooth functions with compact support.

The next result shows that the random part of $X_t$ (see Equation \eqref{eq:def_Lambda_t}) satisfies again the non-degeneracy assumption [\textbf{SD}], even if re-scaled
with respect to the anisotropic structure of the dynamics.

\begin{prop}[Decomposition]
\label{prop:Decomposition_Process_X_Levy}
For any $t$ in $(0,1]$, there exists a L\'evy process $\{S^t_u\}_{u\ge 0}$ such that
\[X_t \, \overset{law}{=} \, e^{tA}x +\mathbb{M}_t S^t_t.\]
Moreover, $\{S^t_u\}_{u\ge 0}$ satisfies assumption [\textbf{SD}] with same $\alpha$ as before.
\end{prop}
\begin{proof}
For simplicity, we start denoting
\begin{equation}\label{eq:def_Lambda_t}
\Lambda_t \,:= \, \int_{0}^{t}e^{(t-s)A}B\, dZ_s,\quad t> 0,
\end{equation}
so that $X_t=e^{tA}x+\Lambda_t$. To conclude, we need to construct a L\'evy process $\{S^t_u\}_{u\ge 0}$ on $\R^N$ satisfying assumption [\textbf{SD}] and
\begin{equation}\label{eq:identity_in_law_Levy}
\Lambda_t \, \overset{law}{=} \, \mathbb{M}_tS_t^t.
\end{equation}
To show the identity in law, we are going to reason in terms of the characteristic functions. By Lemma $2.2$ in \cite{Schilling:Wang12}, we
know that $\Lambda_t$ is an infinitely divisible random variable with associated L\'evy symbol
\[\Phi_{\Lambda_t}(\xi) \, := \, \int_{0}^{t}\Phi\bigl((e^{sA}B)^*\xi\bigr)\,ds, \quad \xi \in \R^N.\]
Remembering the decomposition $e^{sA}B=e^{sA}\mathbb{M}_sB=\mathbb{M}_{s}R_{s}B$ from Lemma \ref{lemma:Decomposition_exp_A}, we can now
rewrite $\Phi_{\Lambda_t}$ as
\[\Phi_{\Lambda_t}(\xi) \, = \, t\int_{0}^{1}\Phi\bigl((e^{stA}B)^*\xi\bigr)\,ds\, = \,
t\int_{0}^{1}\Phi\bigl((R_{st}B)^*\mathbb{M}_s\mathbb{M}_t\xi\bigr)\,ds.\]
The above equality suggests us to define, for any fixed $t$ in $(0,1]$, the (unique in law) L\'evy process  $\{S^t_u\}_{u\ge0}$ associated with the L\'evy symbol
\[\tilde{\Phi}^t(\xi) \, := \, \int_{0}^{1}\Phi\bigl((R_{st}B)^*\mathbb{M}_s\xi\bigr)\,ds, \quad \xi \in \R^N.\]
It is not difficult to check that $\tilde{\Phi}^t$ is indeed a L\'evy symbol associated with the L\'evy triplet
$(\tilde{Q}^t,\tilde{b}^t,\tilde{\nu}^t)$ given by
\begin{align}
  \tilde{Q}^t \, &= \, \int_{0}^{1}\mathbb{M}_sR_{st}BQ(\mathbb{M}_sR_{st}B)^\ast \, ds; \\
  \tilde{b}^t \, = \, \int_{0}^{1}\mathbb{M}_sR_{st}Bb \, ds &+ \int_{0}^{1}\int_{\R^d}
  \mathbb{M}_sR_{st}Bz\bigl[\mathds{1}_{B(0,1)}(\mathbb{M}_sR_{st}Bz)-\mathds{1}_{B(0,1)}(z)\bigr]
  \, \nu(dz)ds;\\
  \tilde{\nu}^t(\mathcal{A}) \, &= \, \int_{0}^{1} \nu\bigl((\mathbb{M}_sR_{st}B)^{-1}\mathcal{A}\bigr)\, ds, \quad \mathcal{A} \in \mathcal{B}(\R^d_0).
\end{align}
 Since we have that
\[\mathbb{E}\bigl[e^{i\langle \xi,\Lambda_t\rangle}\bigr] \, = \, e^{-\Phi_{\Lambda_t}(\xi)}\, = \, e^{-t\tilde{\Phi}^t(\mathbb{M}_t\xi)} \,
= \,\mathbb{E}\bigl[e^{i\langle \xi,\mathbb{M}_tS^t_t\rangle}\bigr],\]
it follows immediately that the identity \eqref{eq:identity_in_law_Levy} holds.\newline
It remains to show that the family of L\'evy measure $\{\tilde{\nu}^t\colon t \in (0,1]\}$ satisfies the assumption [\textbf{SD}].
Recalling that condition [\textbf{SD}] is assumed to hold for $\nu$, we know that
\begin{equation}\label{Proof:eq:Decomposition}
\tilde{\nu}^t(\mathcal{A}) \, = \, \int_{0}^{1} \nu\bigl((\mathbb{M}_sR_{st}B)^{-1}\mathcal{A}\bigr)\, ds \, \ge \,
\int_{0}^{1}\int_{0}^{r_0}\int_{\mathbb{S}^{d-1}} \mathds{1}_{\mathcal{A}}(r\mathbb{M}_sR_{st}B\theta) \mu(d\theta)\frac{dr}{r^{1+\alpha}}ds,
\end{equation}
for any $\mathcal{A}$ in $\mathcal{B}(\R^d_0)$. Furthermore, it holds from Lemma \ref{lemma:Decomposition_exp_A} that
\begin{equation}\label{Proof:eq:Decomposition2}
\inf_{s\in (0,1),\, t \in (0,1],\, \theta \in S^{d-1}}\bigl{\vert} \mathbb{M}_sR_{st}B\theta \bigr{\vert}\, =:\, R_0>0.
\end{equation}
It allows us to define two functions $l^t\colon [0,1]\times
S^{d-1}\to S^{N-1}$, $m^t\colon [0,1]\times S^{d-1}\to \R$, given by
\[l^t(s,\theta) \,:= \, \frac{\mathbb{M}_sR_{st}B\theta}{\vert \mathbb{M}_sR_{st}B\theta\vert} \,\, \text{ and } \,\,  m^t(s,\theta) \,:= \,
\vert \mathbb{M}_sR_{st}B\theta\vert.\]
Using the Fubini theorem, we can now rewrite Equation \eqref{Proof:eq:Decomposition} as
\[
\begin{split}
\tilde{\nu}^t(\mathcal{A}) \, &\ge \,\int_{0}^{1}\int_{\mathbb{S}^{d-1}}\int_{0}^{r_0}\mathds{1}_{\mathcal{A}}(l^t(s,\theta)m^t(s,\theta)r)
\frac{dr}{r^{1+\alpha}}\mu(d\theta)ds \\
&= \, \int_{0}^{1}\int_{\mathbb{S}^{d-1}} \int_{0}^{m^t(s,\theta)r_0}\mathds{1}_{\mathcal{A}}(l^t(s,\theta)r)
\frac{dr}{r^{1+\alpha}}[m^t(s,\theta)]^\alpha\mu(d\theta)ds.
\end{split}
\]
Exploiting again Control \eqref{Proof:eq:Decomposition2}, we can conclude that
\begin{equation}\label{eq:ND_assumption_for_St}
\tilde{\nu}^t(C) \, \ge \, \int_{0}^{1}\int_{\mathbb{S}^{d-1}} \int_{0}^{R_0}\mathds{1}_{C}(l^t(s,\theta)r)
\frac{dr}{r^{1+\alpha}}\tilde{m}^t(ds,d\theta) \, = \, \int_{0}^{R_0}\int_{\mathbb{S}^{N-1}}
\mathds{1}_{C}(\tilde{\theta} r) \tilde{\mu}^t(d\tilde{\theta})\frac{dr}{r^{1+\alpha}},
\end{equation}
where $\tilde{m}^t(ds,d\theta)$ is a measure on $[0,1]\times S^{d-1}$ given by
\[\tilde{m}^t(ds,d\theta) \, := \, [m^t(s,\theta)]^{\alpha}\mu(d\theta)ds\]
and $\tilde{\mu}^t:=(l^t)_{\ast}\tilde{m}^t$ is the measure $\tilde{m}^t$ push-forwarded through $l^t$ on $S^{N-1}$.
It is easy to check that the measure $\tilde{\mu}^t$  is finite and non-degenerate in the sense of \eqref{eq:non_deg_measure_Levy},
replacing therein $d$ by $N$.
\end{proof}

An immediate application of the above result is a first representation formula for the transition semi-group $\{P_t\colon t\ge0\}$ associated
with the Ornstein-Uhlenbeck process $\{X_t\}_{t\ge 0}$, at least for small times. Indeed, denoting by $\mathbb{P}_{X}$ the law of a
random variable $X$, Equation \eqref{eq:identity_in_law_Levy} implies that for any $\phi$ in $B_b(\R^N)$, it holds that
\begin{equation}\label{eq:Representation_semi-group1}
P_t\phi(x) \, = \, \int_{\R^N}\phi(e^{tA}x+y) \, \mathbb{P}_{\Lambda_t}(dy) \, = \, \int_{\R^N}\phi(e^{tA}x+\mathbb{M}_ty) \,
\mathbb{P}_{S^t_t}(dy), \quad x \in \R^N, \, t \in (0,1].
\end{equation}
Moreover, condition [\textbf{SD}] for $\{S^t_u\}_{u\ge 0}$ allows us to decompose it into two components: a truncated, $\alpha$-stable part
and a remainder one. Indeed, if we denote by $\nu^t_{\alpha}$ the measure serving as lower bound to the L\'evy measure $\tilde{\nu}^t$ in \eqref{eq:ND_assumption_for_St}, i.e.
\begin{equation}\label{eq:def_levy_measure_for_stable}
\nu^t_{\alpha}(\mathcal{A}) \, := \, \int_{0}^{R_0}\int_{\mathbb{S}^{N-1}} \mathds{1}_{\mathcal{A}}(\theta r) \tilde{\mu}^t(d\theta)\frac{dr}{r^{1+\alpha}},
\quad C \in \mathcal{B}(\R^N_0),
\end{equation}
we can consider $\{Y^t_u\}_{u\ge 0}$, the L\'evy process on $\R^N$ associated with the L\'evy triplet $(0,0, \nu^t_{\alpha})$. We recall now a useful fact involving the L\'evy symbol $\Phi^t_\alpha$ of the process $Y^t$. The non-degeneracy of the
measure $\tilde{\mu}^t$ is equivalent to the existence of a constant $C>0$ such that
\begin{equation}\label{eq:equivalence_non-deg_measure}
\Phi^t_\alpha(\xi) \, \ge \, C \vert \xi \vert^\alpha, \quad \xi \in \R^N.
\end{equation}
A proof of this result can be found, for example, in \cite{Priola12} p.$424$.\newline
In order to apply the results in \cite{Schilling:Sztonyk:Wang12}, we are going to truncate the above process at the typical
time scale for an $\alpha$-stable process. This is $t^{1/\alpha}$ when considering the process at time $t$ (cf. Example \ref{example_basic}). Namely, we consider the family $\{\mathbb{P}^{\text{tr}}_t\}_{t\ge 0}$ of infinitely
divisible probabilities whose characteristic function has the form
$\widehat{\mathbb{P}^{\text{tr}}_t}(\xi):=\text{exp}[-\Phi^{\text{tr}}_t(\xi)]$, where
\[\Phi^{\text{tr}}_t(\xi) \, := \, \int_{\vert z\vert \le t^{\frac{1}{\alpha}}} \bigl[1-e^{i\langle \xi,z\rangle}+i\langle
\xi,z\rangle\bigr]\, \nu^t_\alpha(dz).\]
On the other hand, since the measure $\tilde{\nu}^t$ satisfies assumption [\textbf{SD}], we know that the remainder
$\tilde{\nu}^t-\mathds{1}_{B(0,t^{1/\alpha})}\nu^t_\alpha$ is again a L\'evy measure on $\R^N$. Let $\{\pi_t\}_{t\ge0}$ be the family of
infinitely divisible probability associated with the following L\'evy triplet:
\[(\tilde{Q}^t,\tilde{b}^t,\tilde{\nu}^t-\mathds{1}_{B(0,t^{1/\alpha})}\nu^t_\alpha).\]
It follows immediately that $\mathbb{P}_{S^t_t} \, = \, \mathbb{P}^{\text{tr}}_t \ast \pi_t$ for any $t>0$. We can now disintegrate the
measure $\mathbb{P}_{S^t_t}$ in Equation \eqref{eq:Representation_semi-group1} in order to obtain
\begin{equation}\label{eq:Representation_semi-group2}
P_t\phi(x) \, = \, \int_{\R^N}\int_{\R^N}\phi\bigl(e^{tA}x+\mathbb{M}_t(y_1+y_2)\bigr) \, \mathbb{P}^{\text{tr}}_t(dy_1) \pi_t(dy_2).
\end{equation}

The next step is to use Proposition $2.3$ in \cite{Schilling:Sztonyk:Wang12} to show a smoothing effect for the family of truncated
stable measures $\{\mathbb{P}^{\text{tr}}_t\colon t\ge 0\}$, at least for small times. Namely,

\begin{prop}
\label{prop:control_on_truncated_density}
Fixed $m$ in $\N_0$, there exists $T_0:=T_0(m)>0$ such that for any $t$ in $(0,T_0]$, the
probability $\mathbb{P}^{\text{tr}}_{t}$ has a density $p^{\text{tr}}(t,\cdot)$ that is $m$-times continuously differentiable on $\R^N$. \newline
Moreover, for any $\vartheta$ in $\mathbb{N}^N$ such that $\vert \vartheta \vert \le m$, there exists a constant $C:=C(m,\vert
\vartheta\vert)$ such that
\[\vert D^\vartheta p^{\text{tr}}(t,y) \vert \, \le \, Ct^{-\frac{N+\vert \vartheta \vert}{\alpha}} \Bigl(1+\frac{\vert y
\vert}{t^{1/\alpha}}\Bigr)^{\vert \vartheta \vert-m}, \quad t \in (0,T_0], \, y \in \R^N.\]
\end{prop}
\begin{proof}
The result follows immediately applying Proposition $2.3$ in \cite{Schilling:Sztonyk:Wang12}. To do so, we need to show that the L\'evy symbol
$\Phi^t_\alpha$ of the process $\{Y^t_u\}_{u\ge 0}$ satisfies the following assumptions:
\begin{itemize}
  \item \emph{Hartman-Wintner condition.} There exists $T>0$ such that
  \[\liminf_{\vert \xi \vert \to \infty}\frac{\text{Re}\Phi^t_{\alpha}(\xi)}{\ln(1+\vert \xi \vert)}\, = \, \infty, \quad t \in (0,T];\]
  \item \emph{Controllability condition.} There exist $T>0$ and $c>0$ such that
  \[\int_{\R^N}e^{-t\text{Re}\Phi^t_{\alpha}(\xi)}\vert \xi \vert^m \, \le \, ct^{-\frac{m+N}{\alpha}}, \quad t \in (0,T].\]
\end{itemize}
In order to show that the above conditions hold, we fix $T\le 1$ and we recall that the L\'evy symbol $\Phi^t_\alpha$ of $Y_t$, the
truncated $\alpha$-stable process with L\'evy measure introduced in \eqref{eq:def_levy_measure_for_stable}, can be written through the
L\'evy-Kitchine formula as
\[\Phi^t_\alpha(\xi) \, = \, \int_{\R^N_0}\bigl(1-e^{i\langle \xi,z\rangle}+i\langle \xi,z\rangle\bigr)\nu_\alpha^t(dz) \, = \,
\int_{0}^{R_0}\int_{\mathbb{S}^{N-1}}\bigl(1-\cos(\langle \xi, r\theta\rangle)\bigr)\, \tilde{\mu}^t(d\theta)\frac{dr}{r^{1+\alpha}}.\]
We have seen in Equation \eqref{eq:equivalence_non-deg_measure} that the non-degeneracy of $\tilde{\mu}^t$ implies that
$\Phi^t_\alpha(\xi)  \ge C \vert \xi \vert^\alpha$. The Hartman-Wintner condition then follows immediately since
\[\liminf_{\vert \xi \vert \to \infty}\frac{\text{Re}\Phi^t_{\alpha}(\xi)}{\ln(1+\vert \xi \vert)}\, \ge \, \liminf_{\vert \xi \vert \to
\infty}\frac{c\vert\xi\vert^\alpha}{\ln(1+\vert \xi \vert)}\, = \,\infty.\]
To show instead the controllability assumption, let us firstly notice that
\[e^{-t\text{Re}\Phi^t_{\alpha}(\xi)} \, \le \,
            \begin{cases}
                1, & \mbox{if } \vert \xi \vert \le R; \\
                e^{-ct\vert \xi \vert^\alpha}, & \mbox{if } \vert \xi \vert > R,
            \end{cases}
\]
for some $R>0$. It then follows that
\begin{align*}
 \int_{\R^N} e^{-t\text{Re}\Phi^t_{\alpha}(\xi)}\vert \xi \vert^m \, d\xi \, &= \, \int_{\vert \xi \vert\le R}\vert \xi \vert^m \, d\xi+
 \int_{\vert \xi \vert>R} e^{-ct\vert\xi\vert^\alpha}\vert \xi \vert^m \, d\xi \\
   &\le \, C+t^{-\frac{m+N}{\alpha}}\int_{\vert \xi \vert>t^{1/\alpha}R}e^{-c\vert\xi\vert^\alpha}\vert \xi \vert^m \, d\xi\\
   &\le \, C+t^{-\frac{m+N}{\alpha}}\int_{\R^N}e^{-c\vert\xi\vert^\alpha}\vert \xi \vert^m \, d\xi \\
   &\le \, Ct^{-\frac{m+N}{\alpha}},
\end{align*}
where in the last step we used that $1\le t^{-\frac{m+N}{\alpha}}$.
\end{proof}

\setcounter{equation}{0}
\mysection{Estimates for transition semi-group}

The results in the previous section (Proposition \ref{prop:control_on_truncated_density} and Equation \eqref{eq:Representation_semi-group2})
allow us to represent the semi-group $P_t$ of the Ornstein-Uhlenbeck process $\{X_t\}_{t\ge 0}$ as
\begin{equation}\label{eq:Representation_semi-group}
P_t\phi(x) \, = \, \int_{\R^N}\int_{\R^N}\phi(\mathbb{M}_t(y_1+y_2)+e^{tA}x)p^{\text{tr}}(t,y_1)\, dy_1\pi_t(dy_2), \quad x \in \R^N,
\end{equation}
at least for small time intervals. \newline
Here, we will focus on estimates in $\Vert \cdot\Vert_\infty$-norm of the transition semi-group $\{P_t\colon t\ge 0\}$  given in Equation
\eqref{eq:Def_transition_semi-group} and its derivatives. The main result in this section is Corollary
\ref{corollary:continuity_between_holder} that shows the continuity of $P_t$ between anisotropic Zygmund-H\"older spaces. These controls
will be fundamental in the next section to prove Schauder Estimates in the elliptic and parabolic settings.\newline
As we will see in the following result, the derivatives of the semi-group $P_t$ with respect to a component $i$ in $I_h$ induces an additional time singularity of order $\frac{1+\alpha(h-1)}{\alpha}$, corresponding to the intrinsic time scale of the considered component.

\begin{prop}
\label{prop:control_in_inf_norm}
Let $h,h',{h''}$ be in $\llbracket 1,n\rrbracket$ and $\phi$ in $B_b(\R^N)$. Then, there exists a constant $C>0$ such that for any $i$ in
$I_h$, any $j$ in $I_{h'}$ and any $k$ in $I_{h''}$, it holds that
\begin{align}
\label{eq:Control_in_inf_norm_deriv1}
    \Vert D_i P_t \phi \Vert_\infty \, &\le \, C\Vert \phi \Vert_\infty \bigl(1+t^{-\frac{1+\alpha(h-1)}{\alpha}}\bigr), \quad t>0;\\
\label{eq:Control_in_inf_norm_deriv2}
    \Vert D^2_{i,j} P_t \phi \Vert_\infty \, &\le \, C\Vert \phi \Vert_\infty \bigl(1+t^{-\frac{2+\alpha(h+h'-2)}{\alpha}}\bigr),\quad t>0;\\
\label{eq:Control_in_inf_norm_deriv3}
  \Vert D^3_{i,j,k} P_t \phi \Vert_\infty \, &\le \, C\Vert \phi \Vert_\infty \bigl(1+t^{-\frac{3+\alpha(h+h'+h''-3)}{\alpha}}\bigr),\quad
  t>0.
\end{align}
\end{prop}
\begin{proof}
We start fixing a time horizon  $T := 1\wedge  T_0(N+4) > 0$, where $T_0(m)$ was defined in Proposition
\ref{prop:control_on_truncated_density}. Our choice of $N+4$ is motivated by the fact that we consider derivatives up to order $3$. \newline
On the interval $(0,T]$, the representation formula \eqref{eq:Representation_semi-group} holds and $P_t\phi$ is three times differentiable for any $\phi$ in $B_b(\R^N)$. We are going to show only Estimate \eqref{eq:Control_in_inf_norm_deriv1} since the controls for the higher derivatives can be obtained similarly.\newline
Fixed $t\le T$, let us consider $i$ in $I_h$ for some $h$ in $\llbracket 1,n\rrbracket$. When $t\le T$, we recall from Equation
\eqref{eq:Representation_semi-group} that, up to a change of variables, it holds that
\[\bigl{\vert}D_i P_t\phi(x)\bigr{\vert} \, = \, \Bigl{\vert} D_i\int_{\R^N}\int_{\R^N} \phi(\mathbb{M}_t(y_1+y_2)) p^{\text{tr}}(t,
y_1-\mathbb{M}^{-1}_te^{tA}x)\, dy_1\pi_t(dy_2) \Bigr{\vert}.\]
We can then move the derivative inside the integral and write that
\begin{align} \notag
\bigl{\vert}D_i P_t\phi(x)\bigr{\vert} \, &= \, \Bigl{\vert} \int_{\R^N}\int_{\R^N} \phi(\mathbb{M}_t(y_1+y_2)) \langle \nabla p^{\text{tr}} (t,
y_1-\mathbb{M}^{-1}_te^{tA}x),\mathbb{M}^{-1}_te^{tA}e_i\rangle\, dy_1\pi_t(dy_2)\Bigr{\vert} \\\notag
&\le \, \vert \mathbb{M}^{-1}_te^{tA} e_i \vert \int_{\R^N}\int_{\R^N} \vert \phi(\mathbb{M}_t(y_1+y_2))\vert \, \vert \nabla p^{\text{tr}}(t,y_1-\mathbb{M}^{-1}_te^{tA}x))\vert\, dy_1\pi_t(dy_2)\\ \label{eq:Control_deriv1}
&\le \, Ct^{-(h-1)}\Vert \phi\Vert_\infty \int_{\R^N}\int_{\R^N} \vert \nabla p^{\text{tr}} (t, y_1)\vert \, dy_1\pi_t(dy_2),
\end{align}
where in the last step we exploited Lemma \ref{lemma:Decomposition_exp_A} to control
\begin{equation}
    \label{eq:control_usual}
\vert \mathbb{M}^{-1}_te^{tA}e_i\vert \, \le \, \sum_{k=1}^{n}\vert  \mathbb{M}^{-1}_tE_ke^{tA}E_he_i\vert \, \le \,
C\bigl[\sum_{k=1}^{h-1}t^{k-h}t+\sum_{k=h}^{n}t^{-(k-1)}t^{-(h-1)}\bigr] \, \le \, Ct^{-(h-1)},
\end{equation}
remembering that $t\le 1$.
We conclude the case $t\le T$ using the control on $p^{\text{tr}}$ (Proposition \ref{prop:control_on_truncated_density} with $m=N+2$) to
write that
\begin{equation}\label{eq:Control_deriv3}
\begin{split}
\bigl{\vert}D_i P_t\phi(x)\bigr{\vert}\, &\le \,  C\Vert \phi \Vert_\infty\pi_t(\R^N)t^{-(h-1)}  \int_{\R^N}
t^{-\frac{N+1}{\alpha}}\bigl(1+\frac{\vert y_1 \vert}{t^{1/\alpha}}\bigr)^{-(N+1)} \, dy_1 \\
&\le \, C\Vert \phi \Vert_\infty t^{-\frac{1+\alpha(h-1)}{\alpha}}  \int_{\R^N} (1+\vert z \vert)^{-(N+1)}\, dz\\
&\le \,C\Vert \phi \Vert_\infty t^{-\frac{1+\alpha(h-1)}{\alpha}}.
\end{split}
\end{equation}
Above, we used the change of variables $z=t^{-1/\alpha}y_1$. When $t>T$, we can exploit the already proven controls for small times, the
semi-group and the contraction properties of $\{P_t\colon t\ge0\}$
on $B_b(\R^N)$ to write that
\begin{equation}\label{eq:COntrol_for_t>1}
\Vert D_i P_t \phi \Vert_\infty \, = \, \Vert D_i P_{T}\bigl(P_{t-T} \phi\bigr) \Vert_\infty \, \le \, C_T\Vert P_{t-T} \phi
\Vert_\infty \, \le \, C\Vert \phi \Vert_\infty.
\end{equation}
We have thus shown Control \eqref{eq:Control_in_inf_norm_deriv1} for any $t>0$.
\end{proof}

The following interpolation inequality (see e.g.\ \cite{book:Triebel83})
\begin{equation}\label{eq:Interpolation_ineq}
\Vert \phi \Vert_{C^{r\delta_1+(1-r)\delta_2}_b} \, \le \, C\Vert \phi\Vert^r_{C^{\delta_1}_b}\Vert \phi \Vert^{1-r}_{C^{\delta_2}_b}
\end{equation}
valid for $0\le \delta_1 <\delta_2$, $r$ in $(0,1)$ and $\phi$ in $C^{\delta_2}(\R^N)$, allows us to extend easily the above
result.

\begin{corollary}
\label{coroll:C0_Cgamma_control}
Let $\gamma$ be in $[0,1+\alpha)$. Then, there exists a constant $C>0$ such that
\begin{equation}\label{eq:coroll:C0_Cgamma_control}
\Vert P_t \Vert_{\mathcal{L}(C_b,C^\gamma_{b,d})} \, \le \, C\bigl(1+t^{-\frac{\gamma}{\alpha}}\bigr), \quad t>0.
\end{equation}
\end{corollary}
\begin{proof}
Let us firstly assume that $\gamma$ is in $(0,1]$. Remembering the definition of $C^\gamma_{b,d}$-norm in \eqref{eq:def_anistotropic_norm_Levy}, we start
fixing a point $x_0$ in $\R^N$ and $h$ in $\llbracket 2,n\rrbracket$. Then, the contraction property of the semi-group implies that
\[\Vert P_t\phi(x_0+\cdot)_{|E_h(\R^N)}\Vert_\infty \, \le \, C\Vert \phi \Vert_\infty.\]
Moreover, Control \eqref{eq:Control_in_inf_norm_deriv1} in Proposition \ref{prop:control_in_inf_norm} ensures
that
\[ \quad \Vert D_iP_t \phi(x_0+ \cdot)_{|E_h(\R^N)} \Vert_\infty \, \le \, C\Vert\phi \Vert_\infty \bigl(1+
t^{-\frac{1+\alpha(h-1)}{\alpha}}\bigr).\]
It follows immediately that
\[\Vert P_t\phi(x_0+\cdot)_{|E_h(\R^N)}\Vert_{C^1_b} \, \le \, C\Vert \phi \Vert_\infty(1+t^{-\frac{1+\alpha(h-1)}{\alpha}}).\]
We can now apply the interpolation inequality \eqref{eq:Interpolation_ineq} with $\delta_1=0$, $\delta_2=1$ and $r=\gamma/(1+\alpha(h-1))$
in order to obtain that
\[\begin{split}
\Vert P_t\phi(x_0+\cdot)_{|E_h(\R^N)}\Vert_{C^r_b} \, &\le \, C\Vert P_t\phi(x_0+\cdot)_{|E_h(\R^N)}\Vert^{r}_{C^1_b} \Vert
P_t\phi(x_0+\cdot)_{|E_h(\R^N)}\Vert^{1-r}_{\infty} \\ 
&\le \, C\Vert \phi\Vert_\infty\bigl(1+t^{-\frac{\gamma}{\alpha}}\bigr).
\end{split}\]
The argument is analogous for $\gamma$ in $(1,3) $, considering only the case $h=0$.
\end{proof}

The next result allows us to extend the controls in Proposition \ref{prop:control_in_inf_norm} to functions in the anisotropic
Zygmund-H\"older spaces. Roughly speaking, it states that the anisotropic $\gamma$-H\"older regularity induces a
"homogeneous" gain in time of order $\gamma/\alpha$ that can be used to weaken, at least partially, the time singularities associated
with the derivatives.
The general argument of proof will mimic the one of Proposition \ref{prop:control_in_inf_norm} even if, this time, we will need to make the
H\"older modulus of $\phi$ appear. It will be managed exploiting some ``partial'' cancellation arguments (cf.\ \eqref{eq:Control_function_K}).

\begin{theorem}\label{prop:Control_in_Holder_norm}
Let $h,h',{h''}$ be in $\llbracket 1,n\rrbracket$ and $\phi$ in $C^\gamma_{b,d}(\R^N)$ for some $\gamma$ in $[0,1+\alpha)$. Then, there
exists a constant $C>0$ such that for any $i$ in $I_h$, any $j$ in $I_{h'}$ and any $k$ in $I_{h''}$, it holds that
\begin{align}
\label{eq:Control_in_Holder_norm_deriv1}
 \Vert D_i P_t \phi \Vert_\infty \, &\le \, C\Vert \phi \Vert_{C^\gamma_{b,d}}\bigl(1+t^{\frac{\gamma-(1+\alpha(h-1))}{\alpha}}\bigr),
  \quad t>0;\\
\label{eq:Control_in_Holder_norm_deriv2}
 \Vert D^2_{i,j}P_t \phi \Vert_\infty \, &\le \,C\Vert \phi\Vert_{C^\gamma_{b,d}}\bigl(1+t^{\frac{\gamma-(2+\alpha(h+h'-2))}{\alpha}}\bigr),
  \quad t>0;\\
\label{eq:Control_in_Holder_norm_deriv3}
 \Vert D^3_{i,j,k}P_t\phi\Vert_\infty \,&\le\,C\Vert\phi\Vert_{C^\gamma_{b,d}}\bigl(1+t^{\frac{\gamma-(3+\alpha(h+h'+h''-3))}{\alpha}}\bigr),
  \quad t>0.
\end{align}
\end{theorem}
\begin{proof}
Similarly to Proposition \ref{prop:control_in_inf_norm}, we start fixing a time horizon
\begin{equation}\label{def:eq:time_horizon_T}
T\,:= \,1\wedge T_0(N+6)>0.
\end{equation}
Then, Corollary \ref{coroll:C0_Cgamma_control} implies the continuity of $P_t$ on
$C^\gamma_{b,d}(\R^N)$,  for any $t\ge T/2$. Indeed,
\begin{equation}\label{Proof:eq:COntinuity_Pt_t>1}
\Vert P_t\phi \Vert_{C^\gamma_{b,d}} \, \le \, C\Vert \phi \Vert_\infty\bigl(1+t^{-\frac{\gamma}{\alpha}}\bigr)\, \le \, C_T\Vert \phi
\Vert_{C^\gamma_{b,d}}.
\end{equation}
The same argument shown in Equation \eqref{eq:COntrol_for_t>1} can now be applied to prove Control \eqref{eq:Control_in_Holder_norm_deriv1} for $t>T$. Namely,
\[\Vert D_i P_t \phi \Vert_\infty \, = \, \Vert D_i P_{T/2}\bigl(P_{t-T/2} \phi\bigr) \Vert_\infty \, \le \, C_T\Vert P_{t-T/2} \phi
\Vert_\infty \, \le \, C\Vert \phi \Vert_\infty.\]
The same reasoning can be used for the higher derivatives, too.\newline
When $t\le T$, let us assume $\alpha>1$, so that $1+\alpha>2$. The case $\alpha\le1$ can be handled similarly taking into
account one less derivative. Moreover, we notice that we need to prove Controls
\eqref{eq:Control_in_Holder_norm_deriv1}-\eqref{eq:Control_in_Holder_norm_deriv3} only for $\gamma$ in $(2,1+\alpha)$ thanks to
interpolation techniques. Indeed, if we want, for example, to prove Estimates \eqref{eq:Control_in_Holder_norm_deriv1} for some $\gamma'$ in
$(0,2]$, we can use Theorem \ref{Thm:Interpolation} to show that
\[\Vert D_iP_t \Vert_{\mathcal{L}(C^{\gamma'}_{b,d};B_b)} \, \le \, \bigl(\Vert D_iP_t
\Vert_{\mathcal{L}(B_b)}\bigr)^{1-\gamma'/\gamma}\bigl(\Vert D_iP_t \Vert_{\mathcal{L}(C^{\gamma}_{b,d},B_b)}\bigr)^{\gamma'/\gamma}\,
\le \, C\bigl(1+t^{\frac{\gamma'-(1+\alpha(h-1))}{\alpha}}\bigr),\]
once we have proven Estimate \eqref{eq:Control_in_Holder_norm_deriv1} for $\gamma>2$.\newline
We are only going to show Control \eqref{eq:Control_in_Holder_norm_deriv1} for $t \le T$ and $\gamma$ in $(2,1+\alpha)$. The estimates
\eqref{eq:Control_in_Holder_norm_deriv2} and \eqref{eq:Control_in_Holder_norm_deriv3} involving the higher derivatives can be obtained in an analogous way.\newline
Fixed $i$ in $I_h$ for some $h$ in $\llbracket 1,n\rrbracket$, we start noticing from Equation \eqref{eq:Representation_semi-group} that, up to the change of variables $\tilde{y}_1=y_1+\mathbb{M}_t^{-1}e^{tA}x$, it holds that
\begin{equation}
\label{eq:inizio}
\begin{split}
D_iP_t\phi(x) \, = \, D_i\int_{\R^N}\int_{\R^N}\phi(\mathbb{M}_t(\tilde{y}_1+y_2))p^{\text{tr}}(t,\tilde{y}_1-\mathbb{M}^{-1}_te^{tA}x)\,
d\tilde{y}_1\pi_t(dy_2) \\
=\, \int_{\R^N}\int_{\R^N}\phi(\mathbb{M}_t(\tilde{y}_1+y_2))D_i\left[p^{\text{tr}}(t,\tilde{y}_1-\mathbb{M}^{-1}_te^{tA}x)\right]\,
d\tilde{y}_1\pi_t(dy_2).
\end{split}
\end{equation}
Recalling that here, $D_i$ stands for the derivative with respect to the variable $x_i$, we then notice that
\[\int_{\R^N}D_i\left[p^{\text{tr}}(t,\tilde{y}_1-\mathbb{M}^{-1}_te^{tA}x)\right]\,
d\tilde{y}_1 \, = \,D_i\int_{\R^N}\left[p^{\text{tr}}(t,\tilde{y}_1-\mathbb{M}^{-1}_te^{tA}x)\right]\,
d\tilde{y}_1 \, = \, 0.\]
In particular, it immediately follows that
\[\int_{\R^N}\int_{\R^N}\phi(\mathbb{M}_ty_2+e^{tA}x) D_i\left[p^{\text{tr}}(t,\tilde{y}_1-\mathbb{M}^{-1}_te^{tA}x)\right]\, dy_1\pi_t(dy_2) \, = \, 0.\]
This property will allow to use a cancellation argument in Equation \eqref{eq:Control_function_K} below, once we split the small jumps in the non-degenerate contributions and the other ones. We thus get from \eqref{eq:inizio} that
\begin{align}
&D_iP_t\phi(x)\notag \\\notag
&= \, 
\int_{\R^N}\int_{\R^N}\left((\phi(\mathbb{M}_t(\tilde{y}_1+y_2)))-\phi(\mathbb{M}_ty_2+e^{tA}x)\right) D_i\left[p^{\text{tr}}(t,\tilde{y}_1-\mathbb{M}^{-1}_te^{tA}x)\right]\, dy_1\pi_t(dy_2)\\
&= \,\int_{\R^N}\int_{\R^N}\Delta \phi(t,y_1,y_2,x) \langle \nabla p^{\text{tr}}(t,y_1),\mathbb{M}^{-1}_te^{tA}e_i\rangle\, dy_1\pi_t(dy_2), \label{eq:Control_function_K}
\end{align}
where, after the backward change of variables, we have also denoted:
\[\Delta \phi(t,y_1,y_2,x) \, :=\,\phi(\mathbb{M}_t(y_1+y_2)+e^{tA}x)-\phi(\mathbb{M}_ty_2+e^{tA}x).\]
We can then decompose the difference $\Delta \phi$ in the following way:
\begin{equation}
    \label{Proof:eq:Control_diff0}
\Delta\phi(t,y_1,y_2,x) \, = \, \Lambda_0(t,y_1,y_2,x)+\Lambda_1(t,y_1,y_2,x)    
\end{equation}
where we denoted
\[
\begin{split}
\Lambda_0(t,y_1,y_2,x)\, &:= \, \phi(E_1y_1+\mathbb{M}_ty_2+e^{tA}x)
-\phi(\mathbb{M}_ty_2+e^{tA}x);\\
\Lambda_1(t,y_1,y_2,x)\, &:= \, \phi(\mathbb{M}_ty_1+\mathbb{M}_ty_2+e^{tA}x)- \phi(E_1y_1+\mathbb{M}_ty_2+e^{tA}x).
\end{split}
\]
Noticing now that the first contribution can be easily controlled
\begin{equation}\label{Proof:eq:Control_diff1}
|\Lambda_1(t,y_1,y_2,x)| \, \le \, \Vert \phi \Vert_{C^\gamma_{b,d}}\sum_{k=2}^{n}\vert E_k\mathbb{M}_ty_1\vert^{\frac{\gamma}{1+\alpha(k-1)}},
\end{equation}
we can then write that
\[
\begin{split}
\Bigl|\int_{\R^N}\int_{\R^N}\Lambda_1(t,y_1,&y_2,x) \langle \nabla p^{\text{tr}}(t,y_1),\mathbb{M}^{-1}_te^{tA}e_i\rangle\, dy_1\pi_t(dy_2) \Bigr|\\
&\le \, \Vert \phi \Vert_{C^\gamma_{b,d}}t^{-(h-1)}\sum_{k=2}^{n}
\int_{\R^N}\int_{\R^N} |\nabla p^{\text{tr}}(t,y_1)|
\vert E_k\mathbb{M}_ty_1\vert^{\frac{\gamma}{1+\alpha(k-1)}}\, dy_1\pi(dy_2),
\end{split}
\]
where we also exploited the control in \eqref{eq:control_usual}.\newline
The above expression allows us to conclude the control for $\Lambda_1$ as in \eqref{eq:Control_deriv3}, using Proposition \ref{prop:control_on_truncated_density} with
$m=N+4$ and $\vert\vartheta\vert=1$. Namely,
\begin{align}\notag
\Bigl|\int_{\R^N}\int_{\R^N}\Lambda_1(t,y_1,&y_2,x) \langle \nabla p^{\text{tr}}(t,y_1),\mathbb{M}^{-1}_te^{tA}e_i\rangle\, dy_1\pi_t(dy_2) \Bigr|\\\notag
&\le \,  C\Vert \phi \Vert_{C^\gamma_{b,d}}t^{-(h-1)}\sum_{k=2}^{n}\int_{\R^N}
t^{-\frac{N+1}{\alpha}}\bigl(1+ \frac{\vert y_1\vert}{t^{\frac{1}{\alpha}}} \bigr)^{-(N+3)} \vert
E_k\mathbb{M}_ty_1\vert^{\frac{\gamma}{1+\alpha(k-1)}} \, dy_1 \\\notag
&\le \, C \Vert \phi \Vert_{C^\gamma_{b,d}} t^{\frac{\gamma-(1+\alpha(h-1))}{\alpha}}\sum_{k=2}^{n}\int_{\R^N}\bigl(1+\vert z
\vert\bigr)^{-(N+3)} \vert z\vert^{\frac{\gamma}{1+\alpha(k-1)}} \, dz \\
&\le \,  C \Vert \phi \Vert_{C^\gamma_{b,d}}t^{\frac{\gamma -(1+\alpha(h-1))}{\alpha}},\label{eq:Control_new1}
\end{align}
where in the second step we used again the change of variable $z=y_1t^{-1/\alpha}$.\newline
For the term $\Lambda_0$, we start instead applying a Taylor expansion of second order along $E_1y_1$:
\begin{align}\notag
\Lambda_0(t,x,&y_1,y_2) \, = \, \phi(E_1y_1+\mathbb{M}_ty_2+e^{tA}x)
-\phi(\mathbb{M}_ty_2+e^{tA}x)\\\notag
&= \,\langle \nabla \phi(\mathbb{M}_ty_2+e^{tA}x),E_1y_1\rangle +\int_0^1\langle D^2\phi(\mathbb{M}_ty_2+e^{tA}x+\lambda E_1y_1)R_1y_1,E_1y_1\rangle \, d\lambda\\
&=: \, \Lambda_{2}(t,x,y_1,y_2)+\Lambda_{3}(t,x,y_1,y_2).
\label{Proof:eq:Control_diff2}
\end{align}
Now, we want to control the component involving the second term $\Lambda_{2}$:
\[\int_{\R^N}\int_{\R^N}\Lambda_{2}(t,y_1,y_2,x) \langle \nabla p^{\text{tr}}(t,y_1),\mathbb{M}^{-1}_te^{tA}e_i\rangle\, dy_1\pi_t(dy_2).\]
To do it, we exploit an integration by part formula with respect to the derivative of $p^{\text{tr}}$. Indeed, reasoning component-wise if necessary, it is not difficult to check that
\begin{align}
\label{eq:Control_strano}
\int_{\R^N}\int_{\R^N}\Lambda_{2}(t,y_1,&y_2,x) \langle \nabla p^{\text{tr}}(t,y_1),\mathbb{M}^{-1}_te^{tA}e_i\rangle\, dy_1\pi_t(dy_2) \\\notag
&= \, \int_{\R^N}\int_{\R^N}\langle \nabla \phi(\mathbb{M}_ty_2+e^{tA}x),E_1y_1\rangle \langle \nabla p^{\text{tr}}(t,y_1),\mathbb{M}^{-1}_te^{tA}e_i\rangle\, dy_1\pi_t(dy_2)\\\notag
&= \, \int_{\R^N}\int_{\R^N}\langle \nabla \phi(\mathbb{M}_ty_2+e^{tA}x),E_1\mathbb{M}^{-1}_te^{tA}e_i\rangle  p^{\text{tr}}(t,y_1)\, dy_1\pi_t(dy_2).
\end{align}
It then follows immediately that
\begin{equation}
\label{eq:Control_new2}
\begin{split}
\Bigl{\vert} \int_{\R^N}\int_{\R^N}\Lambda_{2}(t,y_1,y_2,x) \langle \nabla p^{\text{tr}}(t,y_1),\mathbb{M}^{-1}_t&e^{tA}e_i\rangle\, dy_1\pi_t(dy_2)\Bigr{\vert} \\
&\le \, C\Vert \phi\Vert_{C^\gamma_{b,d}}|E_1e^{tA}e_i|\int_{\R^N} |p^{\text{tr}}(t,y_1)| \, dy_1\\
&\le \, C\Vert \phi\Vert_{C^\gamma_{b,d}},
\end{split}
\end{equation}
where in the last passage we also used Lemma \ref{lemma:Decomposition_exp_A}. \newline
To control the contribution involving the third term $\Lambda_{3}$, we will need  an additional cancellation argument. Let us assume for the moment that the family of truncated probabilities $\{\mathbb{P}^{\text{tr}}_t\}_{t\ge 0}$ has zero mean value, so that it holds that:
\[\int_{\R^N}E_1y_1p^{\text{tr}}(t,y_1) \, dy_1 \, = \, 0_{\R^d}.\]
Under this additional hypothesis, it is possible to show that
\begin{equation}
\label{eq:cancellation2}
    \int_{\R^N}\int_{\R^N}\langle D^2\phi(\mathbb{M}_ty_2+e^{tA}x)E_1y_1,E_1y_1\rangle \langle \nabla p^{\text{tr}}(t,y_1),\mathbb{M}^{-1}_te^{tA}e_i\rangle\, dy_1\pi_t(dy_2) \, = \, 0.
\end{equation}
Indeed, applying again an integration by parts formula with respect to the derivative on $p^\text{tr}$, we notice, reasoning as well component-wise as in \eqref{eq:Control_strano}, that:
\[\begin{split}
\int_{\R^N}\langle D^2\phi(\mathbb{M}_ty_2+e^{tA}x)E_1y_1,&E_1y_1\rangle \langle \nabla p^{\text{tr}}(t,y_1),\mathbb{M}^{-1}_te^{tA}e_i\rangle\, dy_1 \\
&= \, \int_{\R^N}\langle D^2\phi(\mathbb{M}_ty_2+e^{tA}x)E_1y_1,E_1\mathbb{M}^{-1}_te^{tA}e_i\rangle p^{\text{tr}}(t,y_1)\, dy_1\\
&= \, \left\langle D^2\phi(\mathbb{M}_ty_2+e^{tA}x)\left[\int_{\R^N}E_1y_1 p^{\text{tr}}(t,y_1)\, dy_1\right],E_1e^{tA}e_i \right\rangle
\\
&= \, 0.
\end{split}\]
The cancellation argument in \eqref{eq:cancellation2} allows now to write that
\[\begin{split}
\int_{\R^N}\int_{\R^N}&\Lambda_3(t,y_1,y_2,x) \langle \nabla p^{\text{tr}}(t,y_1),\mathbb{M}^{-1}_te^{tA}e_i\rangle\, dy_1\pi_t(dy_2) \\
&= \, \int_{\R^N}\int_{\R^N}\int_0^1\left\langle \left[D^2\phi(\mathbb{M}_ty_2+e^{tA}x+\lambda E_1y_1)-D^2\phi(\mathbb{M}_ty_2+e^{tA}x)\right]E_1y_1,E_1y_1\right\rangle \\
&\qquad \qquad\qquad \qquad\times \langle \nabla p^{\text{tr}}(t,y_1),\mathbb{M}^{-1}_te^{tA}e_i\rangle\, dy_1\pi_t(dy_2).
\end{split}\]
The same arguments presented in Control \eqref{eq:Control_new1} can be also applied here to show that
\begin{align}\notag
\Bigl{\vert} \int_{\R^N}\int_{\R^N}\Lambda_3(t,&y_1,y_2,x) \langle \nabla p^{\text{tr}}(t,y_1),\mathbb{M}^{-1}_te^{tA}e_i\rangle\, dy_1\pi_t(dy_2) \Bigr{\vert} \\ \notag
&\le \,C\Vert \phi\Vert_{C^\gamma_{b,d}}t^{-(h-1)} \int_{\R^N}\int_{\R^N}|E_1y_1|^{\gamma-2}|E_1y_1|^2 |\nabla p^{\text{tr}}(t,y_1)|\, dy_1\pi_t(dy_2)\\
&\le \, C\Vert \phi\Vert_{C^\gamma_{b,d}}t^{\frac{\gamma}{\alpha}-\frac{1+\alpha(h-1)}{\alpha}}.
\label{eq:Control_new3}
\end{align}
Going back to Expression \eqref{eq:Control_function_K}, we notice that the decompositions in \eqref{Proof:eq:Control_diff0} and \eqref{Proof:eq:Control_diff2} implies immediately that
\[|D_iP_t\phi(x) |\,\le \,\sum_{i=1}^3\Bigl{\vert}\int_{\R^N}\int_{\R^N}\Lambda_i(t,y_1,y_2,x) \langle \nabla p^{\text{tr}}(t,y_1),\mathbb{M}^{-1}_te^{tA}e_i\rangle\, dy_1\pi_t(dy_2)\Bigr{\vert}.\]
Then, we can finally use the estimates in \eqref{eq:Control_new1},  \eqref{eq:Control_new2} and
\eqref{eq:Control_new3}
to conclude that
\[\bigl{\vert}D_iP_t\phi(x)\bigr{\vert} \, \le \, C\Vert \phi \Vert_{C^\gamma_{b,d}}\left(1+t^{\frac{\gamma}{\alpha}-\frac{1+\alpha(i-1)}{\alpha}}\right).\]
In the general case, when the family of probabilities $\{\mathbb{P}^{\text{tr}}_t\}_{t\ge0}$ possibly has non-zero mean value, it is then enough to follow the same reasoning above, taking into account, in the first cancellation argument, the additional error given by the mean.\newline
Namely, we will need to consider $\phi(\mathbb{M}_ty_2+e^{tA}x-m_t)$, where $m_t$ is the mean value associated with $\mathbb{P}_t$, in Equation \eqref{eq:Control_function_K}.
\end{proof}

Next, we are going to use the controls in Theorem \ref{prop:Control_in_Holder_norm} to show the main result of this section. It states
the continuity of the semi-group $P_t$ between anisotropic Zygmund-H\"older spaces at a cost of additional time singularities.

\begin{corollary} \label{corollary:continuity_between_holder}
Let $\beta,\gamma$ be in $[0,1+\alpha)$ such that $\beta\le \gamma$. Then, there exists a constant $C>0$ such that
\begin{equation}\label{eq:coroll:Cbeta_Cgamma_control}
\Vert P_t \Vert_{\mathcal{L}(C^\beta_{b,d},C^\gamma_{b,d})} \, \le \, C\bigl(1+t^{\frac{\beta-\gamma}{\alpha}}\bigr), \quad
t>0.
\end{equation}
\end{corollary}
\begin{proof}
It is enough to show the result only for $\gamma=\beta$ non-integer, thanks to interpolation techniques. Indeed,
fixed $\beta< \gamma$, we can use Theorem \ref{Thm:Interpolation} to show that
\[\Vert P_t\Vert_{\mathcal{L}(C^\beta_{b,d}(\R^N),C^{\gamma}_{b,d}(\R^N))} \,\le \, \Bigl(\Vert
P_t\Vert_{\mathcal{L}(C_b(\R^N),C^{\gamma}_{b,d}(\R^N))}\Bigr)^{1-\frac{\beta}{{\gamma}}}\Bigl(\Vert
P_t\Vert_{\mathcal{L}(C^{\gamma}_{b,d}(\R^N))}\Bigr)^{\frac{\beta}{{\gamma}}}.\]
On the other hand, if we fix $\gamma$ integer, we can take $\gamma'$ in $(\gamma,1+\alpha)$ non-integer such that Theorem \ref{Thm:Interpolation} implies:
\[\Vert P_t\Vert_{\mathcal{L}(C^{\gamma}_{b,d}(\R^N))} \, \le  \Bigl(\Vert
P_t\Vert_{\mathcal{L}(C_b(\R^N)}\Bigr)^{1-\frac{\gamma}{\gamma'}}\Bigl(\Vert
P_t\Vert_{\mathcal{L}(C^{\gamma'}_{b,d}(\R^N))}\Bigr)^{\frac{\gamma}{\gamma'}}.\]
The general result will then follows from the two above controls and Equation \eqref{eq:coroll:C0_Cgamma_control}, once we
have shown  Estimate \eqref{eq:coroll:Cbeta_Cgamma_control}  for $\gamma=\beta$ non-integer.\newline
Fixed again the time horizon $T$ given in \eqref{def:eq:time_horizon_T}, we start noticing that Control \eqref{eq:coroll:Cbeta_Cgamma_control} for $t\ge T$ has already been shown in Equation \eqref{Proof:eq:COntinuity_Pt_t>1} .\newline
To prove it when $t\le T$, we are going to exploit the equivalent norm defined in \eqref{eq:def_equivalent_norm} of Lemma \ref{lemma:def_equivalent_norm}. For this reason, we fix $h$ in $\llbracket 1,n \rrbracket$, a point $x_0$ in $\R^N$ and
$z\neq 0$ in $E_h(\R^N)$ and we would like to show that
\begin{equation}\label{proof_eq:Corollary1}
\bigl{\vert} \Delta^3_{x_0}\bigl(P_t\phi\bigr)(z) \bigr{\vert}  \, \le \, C\Vert \phi\Vert_{C^\gamma_{b,d}}\vert
z\vert^{\frac{\gamma}{1+\alpha(h-1)}},
\end{equation}
for some constant $C>0$ independent from $x_0$. Before starting with the calculations, we highlight the presence of three
different "regimes" appearing below. On the one hand, we will firstly consider a \emph{macroscopic regime} appearing for $\vert z \vert\ge
1$. On the other hand, we will say that the \emph{off-diagonal regime} holds if
$t^{\frac{1+\alpha(h-1)}{\alpha}} \le \vert z \vert \le 1$. It will mean in particular that the spatial distance is larger than the
characteristic time-scale. Finally, a \emph{diagonal regime} will be in force when $t^{\frac{1+\alpha(h-1)}{\alpha}} \ge \vert z \vert$ and
the spatial point will be instead smaller than the typical time-scale magnitude. While for the two first regimes, we are going to use the
contraction property of the semi-group, the third regime will require to exploit the controls in H\"older norms given by Theorem
\ref{prop:Control_in_Holder_norm}. \newline
As said above, Estimate \eqref{proof_eq:Corollary1} in the macroscopic regime (i.e.\ $\vert z \vert \ge 1$) follows immediately from
the contraction property of $P_t$ on $B_b(\R^N)$. Indeed,
\begin{equation}\label{proof_eq:Corollary1-1}
\bigl{\vert}\Delta^3_{x_0}\bigl(P_t\phi\bigr)(z) \bigr{\vert} \, \le \, C\Vert P_t\phi \Vert_\infty \, \le \, C\Vert
\phi\Vert_{C^\gamma_{b,d}}\vert z\vert^{\frac{\gamma}{1+\alpha(h-1)}}.
\end{equation}
For $t^{\frac{1+\alpha(h-1)}{\alpha}} \le \vert z \vert \le 1$ and $l$ in $\llbracket 0,3\rrbracket$, we start noticing from Equation
\eqref{eq:Representation_semi-group} that
\[
\begin{split}
P_t\phi(x_0+lz) \, &=\, \int_{\R^N}\int_{\R^N}\phi\bigl(\mathbb{M}_t(y_1+y_2)+e^{tA}(x_0+lz)\bigr) p^{\text{tr}}(t,y_1)\, dy_1\pi_t(dy_2)\\
&=\, \int_{\R^N}\int_{\R^N}\phi\bigl(\xi_0+le^{tA}z)\bigr) p^{\text{tr}}(t,y_1)\, dy_1\pi_t(dy_2),
\end{split}
\]
where we have denoted for simplicity $\xi_0=\mathbb{M}_t(y_1+y_2)+e^{tA}x_0$. We can then exploit Lemma \ref{lemma:Decomposition_exp_A} to
write that
\begin{equation}\label{proof_eq:Corollary1-2}
\begin{split}
\bigl{\vert}\Delta^3_{x_0}\bigl(P_t\phi\bigr)(z) \bigr{\vert} \, &\le \,
\int_{\R^N}\int_{\R^N}\bigl{\vert}\Delta^3_{\xi_0}\phi(e^{tA}z) \bigr)\bigr{\vert}
p^{\text{tr}}(t,y_1)\,dy_1\pi_t(dy_2) \\
&\le \,  \pi_t(\R^N)\Vert \phi \Vert_{C^\gamma_{b,d}}\sum_{k=1}^{n}\vert E_ke^{tA}z\vert^{\frac{\gamma}{1+\alpha(k-1)}} \\
&\le \, C\Vert \phi \Vert_{C^\gamma_{b,d}}\Bigl[\sum_{k=1}^{h-1}\bigl(t\vert
z\vert\bigr)^{\frac{\gamma}{1+\alpha(k-1)}}+\sum_{k=h}^{n}\bigl(t^{k-h}\vert z\vert\bigr)^{\frac{\gamma}{1+\alpha(k-1)}}\Bigr] \\
&\le C\Vert \phi \Vert_{C^\gamma_{b,d}}\vert z \vert^{\frac{\gamma}{1+\alpha(h-1)}}.
\end{split}
\end{equation}
For $\vert z \vert\le t^{\frac{1+\alpha(h-1)}{\alpha}}$, we are going to apply Taylor expansion three times in order to make $D^3_{I_h}$
appear. Namely,
\begin{align}\notag
\bigl{\vert}\Delta^3_{x_0}\bigl(&P_t\phi\bigr)(z)\bigr{\vert}\\\notag
&= \, \Bigl{\vert}\int_{0}^{1} \langle
D_{I_h}P_t\phi(x_0+\lambda z)-2D_{I_h}P_t\phi(x_0+z+\lambda z)+D_{I_h}P_t\phi(x_0+2z+\lambda z), z \rangle \,d\lambda\Bigr{\vert} \\\notag
&\le \, \Bigl{\vert}\int_{0}^{1}\int_{0}^{1} \langle \bigl[D^2_{I_h}P_t\phi(x_0+(\lambda+\mu)z)-D^2_{I_h}P_t\phi(x_0+z+(\lambda+\mu)z)\bigr]z,
z \rangle \,d\lambda d\mu\Bigr{\vert} \\\notag
&\le \, \Bigl{\vert}\int_{0}^{1}\int_{0}^{1}\int_{0}^{1} \langle \bigl[D^3_{I_h}P_t\phi(x_0+(\lambda+\mu+\nu)z)\bigr](z,z),z \rangle
\,d\lambda d\mu d\nu\Bigr{\vert} \\ \label{Proof:eq:Taylor_Expansion}
&\le \, C\Vert D^3_{I_h}P_t\phi\Vert_\infty \vert z \vert^3 \\
&\le \, C\Vert \phi \Vert_{C^{\gamma}_{b,d}}\bigl(1+t^{\frac{\gamma-3(1+\alpha(h-1))}{\alpha}}\bigr) \vert z \vert^3, \notag
\end{align}
where in the last step we used Control \eqref{eq:Control_in_Holder_norm_deriv3} with $h=h'={h''}$. Since $\vert z \vert \le
t^{\frac{1+\alpha(h-1)}{\alpha}}$ and noticing that $\gamma-3(1+\alpha(h-1))<0$, it holds that
\[\bigl(1+t^{\frac{\gamma-3(1+\alpha(h-1))}{\alpha}}\bigr)\vert z \vert^3 \, \le \, \vert z
\vert^{\frac{\gamma-3(1+\alpha(h-1))}{1+\alpha(h-1)}}\vert z \vert^3 \, = \, \vert z \vert^{\frac{\gamma}{1+\alpha(h-1)}}.\]
We can then conclude that
\begin{equation}\label{proof_eq:Corollary1-3}
\bigl{\vert}\Delta^3_{x_0}\bigl(P_t\phi\bigr)(z)\bigr{\vert}\, \le \, C\Vert \phi \Vert_{C^\gamma_{b,d}}\vert
z\vert^{\frac{\gamma}{1+\alpha(h-1)}}.
\end{equation}
Going back to Controls \eqref{proof_eq:Corollary1-1}, \eqref{proof_eq:Corollary1-2} and \eqref{proof_eq:Corollary1-3}, we have thus proven
Estimate \eqref{proof_eq:Corollary1} for any non-integer $\gamma=\beta$.
\end{proof}

\setcounter{equation}{0}
\section{Elliptic and parabolic Schauder estimates}
\fancyhead[RO]{Section \thesection. Elliptic and parabolic Schauder estimates}
In this section, we use the controls shown before to prove Schauder Estimates both for the elliptic and the parabolic equation driven by the Ornstein-Ulhenbeck operator $L^{\text{ou}}$.

Fixed $\lambda>0$ and $g$ in $C_b(\R^N)$, we say that a function $u\colon \R^N\to \R^N$ is a \emph{distributional solution} of Elliptic
Equation \eqref{eq:Elliptic_IPDE} if $u$ is in $C_b(\R^N)$ and for any $\phi$ in $C^\infty_c(\R^N)$ (i.e.\ smooth functions with compact
support), it holds that
\begin{equation}\label{eq:weak_solution_elliptic}
\int_{\R^N}u(x)\bigl[\lambda\phi(x)-\bigl(L^{\text{ou}}\bigr)^*\phi(x)\bigr] \, dx \, = \,
\int_{\R^N}\phi(x)g(x) \, dx,
\end{equation}
where $\bigl(L^{\text{ou}}\bigr)^*$ denotes the formal adjoint of $L^{\text{ou}}$ on $L^2(\R^N)$, i.e.
\begin{equation}\label{eq:def_adjoint_op}
\bigl(L^{\text{ou}}\bigr)^*\phi(x) \, = \, \mathcal{L}^\ast \phi(x) -\langle Ax,D_x \phi(x)\rangle - \text{Tr}(A)\phi(x), \quad
(t,x) \in [0,T]\times \R^N,
\end{equation}
and $\mathcal{L}^\ast$ is the adjoint of the operator $\mathcal{L}$ on $L^2(\R^N)$. It is well-known (see e.g. Section $4.2$ in
\cite{book:Applebaum19}) that it can be represented for any $\phi$ in $C^\infty_c(\R^N)$ as
\begin{multline*}
\mathcal{L}^\ast\phi(x) \, = \,\frac{1}{2}\text{Tr}\bigl(BQB^\ast D^2\phi(x)\bigr) -\langle Bb,D \phi(x)\rangle \\
+ \int_{\R^d_0}\bigl[\phi(x-Bz)-\phi(x)+\langle D\phi(x),
Bz\rangle\mathds{1}_{B(0,1)}(z) \bigr] \,\nu(dz).
\end{multline*}
We state now the main result for the elliptic case, ensuring the well-posedness (in a distributional sense) for Equation \eqref{eq:Elliptic_IPDE}.

\begin{theorem}
\label{thm:well_posedness_elliptic}
Fixed $\lambda>0$, let $g$ be in $C_b(\R^N)$. Then, the function $u\colon \R^N\to \R$ given by
\begin{equation}\label{eq:representation_elliptic_sol}
u(x) \, := \, \int_{0}^{\infty}e^{-\lambda t}P_tg(x) \, dt, \quad x \in \R^N,
\end{equation}
is the unique distributional solution of Equation \eqref{eq:Elliptic_IPDE}.
\end{theorem}
\begin{proof}
\emph{Existence.} We are going to show that the function $u$ given in Equation \eqref{eq:representation_elliptic_sol} is indeed a
distributional solution of the elliptic problem \eqref{eq:Elliptic_IPDE}. It is straightforward to notice that $u$ is in $C_b(\R^N)$, thanks to the contraction property of $P_t$ on $C_b(\R^N)$. Fixed $\phi$ in $C^\infty_c(\R^N)$, we then use Fubini Theorem to write that
\[
\begin{split}
\int_{\R^N}u(x)\bigl(L^{\text{ou}}\bigr)^\ast\phi(x) \, dx \, &= \,
\lim_{\epsilon\to0^+}\int_{\epsilon}^{\infty}\int_{\R^N}e^{-\lambda t}
P_tg(x)\bigl(L^{\text{ou}}\bigr)^\ast\phi(x) \, dx dt \\
&= \,\lim_{\epsilon\to0^+}\int_{\epsilon}^{\infty}\int_{\R^N}e^{-\lambda t}
L^{\text{ou}}P_tg(x)\phi(x) \, dx dt,
\end{split}
\]
where, in the last step, we exploited that $P_tg$ is differentiable and bounded for $t>0$ (Proposition
\ref{prop:control_in_inf_norm}). Since $L^{\text{ou}}$ is the infinitesimal generator of the semi-group $\{P_t\colon
t\ge0\}$, we know that $\partial_t(P_tg)$ exists for any $t>0$ and $\partial_t(P_tg)(x)=L^{\text{ou}}P_tg(x)$ for any $x$ in
$\R^N$. Integration by parts formula allows then to conclude that
\[
\begin{split}
\int_{\R^N}u(x)\bigl(L^{\text{ou}}\bigr)^\ast\phi(x) \, dx \, &=
\,\lim_{\epsilon\to0^+}\int_{\R^N}\phi(x)\int_{\epsilon}^{\infty}e^{-\lambda
t}\partial_tP_tg(x) \, dt dx \\
&= \, \lim_{\epsilon\to0^+}\int_{\R^N}\Bigl(-e^{-\lambda \epsilon}P_\epsilon g(x)+\lambda\int_{\epsilon}^{\infty}e^{-\lambda t}P_t g(x) \,
dt\Bigr) \, dx \\
&= \, \int_{\R^N}-g(x)\phi(x) \, dx + \int_{\R^N}\lambda u(x)\phi(x) \, dx.
\end{split}
\]

\emph{Uniqueness.} It is enough to show that any distributional solution $u$ of Equation \eqref{eq:Elliptic_IPDE} for $g=0$ coincides with
the zero function, i.e.\ $u=0$. To do so, we fix a function $\rho$ in $C^{\infty}_c(\R^N)$ such that $\Vert
\rho \Vert_{L^1}=1$, $0\le \vert \rho \vert \le 1$ and we then define the \emph{mollifier} $\rho_m:=m^N\rho(mx)$ for any $m$ in $\N$. Denoting now, for simplicity, $u_m:=u\ast \rho_m$, we define the function
\begin{equation}\label{Proof:eq:def_gm}
g_m(x) \, := \, \lambda u_m(x)-L^{\text{ou}}u_m(x).
\end{equation}
Using that $u$ is in $C_b(\R^N)$, it is easy to notice that $g_m$ is also in $C_b(\R^N)$ for any fixed $m$ in $\N$.
Truncating the functions if necessary, we can assume that $u_m$ and $g_m$ are integrable with integrable Fourier transform so that we can
apply the Fourier transform in Equation \eqref{Proof:eq:def_gm}:
\begin{equation}\label{proof:Schauder_elliptic1}
\lambda\widehat{u}_m(\xi)- \mathcal{F}_x \bigl(L^{\text{ou}} u_m\bigr)(\xi) \, = \, \widehat{g}_m(\xi).
\end{equation}
We remember in particular that the above operator $L^{\text{ou}}$ has an associated L\'evy symbol
$\Phi^{\text{ou}}(\xi)$ and, following Section $3.3.2$ in \cite{book:Applebaum09}, it holds that
\begin{equation}\label{proof:Schauder_elliptic2}
\mathcal{F}_x \bigl(L^{\text{ou}} u_m\bigr)(\xi) \, = \,
\Phi^{\text{ou}}(\xi)\widehat{u}_m(\xi).
\end{equation}
We can then use it to show that $\widehat{u}_m$ is a classical solution of the following equation:
\[\bigl[\lambda- \Phi^{\text{ou}}(\xi)\bigr]\widehat{u}_m(\xi) \, = \, \widehat{g}_m(\xi).\]
The above equation can be easily solved by direct calculation as
\[\widehat{u}_m(\xi) \, = \, \int_{0}^{\infty} e^{-\lambda t}e^{t\Phi^{\text{ou}}(\xi)} \widehat{g}_m(\xi)\, ds.\]
In order to go back to $u_m$, we apply now the inverse Fourier transform to write that
\[u_m(x) \, = \, \int_{0}^{\infty}e^{-\lambda t}P_{t}g_m(x) \, dt.\]
The contraction property of the semi-group $P_t$  then implies that $\Vert u_m \Vert_\infty \le C \Vert g_m \Vert_\infty$.\newline
In order to conclude, we need to show that
\begin{equation}\label{Proof:convergence_gm}
\lim_{m \to \infty}\Vert g_m \Vert_{\infty} \, = \,  0.
\end{equation}
We start noticing that, since $u$ is a solution of Equation \eqref{eq:Elliptic_IPDE} with $g=0$, it holds that
\[
\begin{split}
g_m(x) \, &= \, \int_{\R^N}u(y)\bigl{\{}\lambda \rho_m(x-y) - \mathcal{L}[\rho_m(\cdot-y)](x)-\langle Ax,D_x\rho_m(x-y)\rangle\bigr{\}} \,
dy\\
&= \, \int_{\R^N}u(y)\bigl{\{}\mathcal{L}^\ast[\rho_m(x-\cdot)](y) - \mathcal{L}[\rho_m(\cdot-y)](x)+\langle
A(x-y),D_x\rho_m(x-y)\rangle\\
&\qquad\quad\qquad \qquad\qquad \qquad\qquad \qquad\qquad \qquad\qquad \qquad\quad+\text{Tr}(A)\rho_m(x-y)\bigr{\}} \, dy\\
&=: \, R^1_m(x)+R^2_m(x)+R^3_m(x),
\end{split}\]
where we have denoted
\begin{align*}
  R^1_m(x) \, &:= \, \int_{\R^N}u(y)\bigl[\mathcal{L}^\ast[\rho_m(x-\cdot)](y) - \mathcal{L}[\rho_m(\cdot-y)](x)\bigr] \, dy;  \\
  R^2_m(x) \, &:= \,  \int_{\R^N}u(y)\langle A(x-y),D_x\rho_m(x-y)\rangle \, dy;\\
  R^3_m(x) \, &:= \,  \int_{\R^N}u(y)\text{Tr}(A)\rho_m(x-y) \, dy.
\end{align*}
On the one hand, it is easy to notice that $R^1_m=0$, since $\mathcal{L}^\ast[\rho_m(x-\cdot)](y)=\mathcal{L}[\rho_m(\cdot-y)](x)$ for any
$m$ in $\N$ and any $y$ in $\R^N$. Indeed, it holds that
\begin{multline*}
\frac{1}{2}\text{Tr}\bigl(BQB^\ast D^2_y[\rho_m(x-\cdot)](y)\bigr)-\langle Bb,D_y[\rho_m(x-\cdot)](y)\rangle \\
= \,
\frac{1}{2}\text{Tr}\bigl(BQB^\ast D^2_x\rho_m(x-y)\bigr)+\langle Bb,D_x\rho_m(x-y)\rangle
\end{multline*}
and
\begin{multline*}
  \int_{\R^d_0}\bigl[\rho_m(x-y+Bz)-\rho_m(x-y)+\langle
D_y[\rho_m(x-\cdot)](y),Bz\rangle\mathds{1}_{B(0,1)}(z)\bigr]\, \nu(dz) \\
= \, \int_{\R^d_0}\bigl[\rho_m((x+Bz)-y)-\rho_m(x-y)-\langle
D_x\rho_m(x-y),Bz\rangle\mathds{1}_{B(0,1)}(z)\bigr]\, \nu(dz).
\end{multline*}
On the other hand, it can be checked (see e.g. \cite{Priola09}) that $\Vert R^2_m+R^3_m\Vert_\infty\to 0$ if $m$ goes to infinity.
Indeed, we firstly notice that $R^3_m$ converges, when $m$ goes to infinity, to the function
$u\text{Tr}(A)$, uniformly in $x$. On the other hand, applying the change of variables $y=x-z/m$ in $R^2_m$, we can obtain that
\[R^2_m(x) \, = \,  m\int_{\R^N}u(x-z/m)\langle A(z/m),D_x\rho(z)\rangle \, dy.\]
Letting $m$ goes to infinity above, we can then conclude that $R^2_m$ converges to the function $-u\text{Tr}(A)$, uniformly in $x$.
\end{proof}

Let us deal now with the parabolic setting. Since we are working with evolution equations, the functions we consider will often depend
on time, too. We denote for any $\gamma>0$ the space $L^\infty(0,T,C^{\gamma}_{b,d}(\R^{N}))$ as the family of functions
$\phi$ in $B_b\bigl([0,T]\times \R^{N})$ such that $\phi(t,\cdot)$ is in $C^{\gamma}_{b,d}(\R^N)$ at any fixed $t$ and the norm
\[\Vert \phi \Vert_{L^\infty(C^\gamma_{b,d})} \, := \, \sup_{t\in [0,T]}\Vert \phi(t,\cdot)\Vert_{C^\gamma_{b,d}} \text{ is finite.}\]

We define now the notion of solution we are going to consider. Fixed $T>0$, $u_0$ in $C_b(\R^N)$ and $f$ in
$L^\infty\bigl(0,T;C_b(\R^N)\bigr)$, we
say that a function $u\colon [0,T]\times\R^N\to \R^N$ is a \emph{weak solution} of the Cauchy problem \eqref{eq:Parabolic_IPDE} if $u$ is in
$L^\infty\bigl(0,T;C_b(\R^N)\bigr)$ and for any $\phi$ in $C^\infty_c([0,T)\times\R^N)$, it holds that
\begin{multline}\label{eq:weak_solution_Parabolic}
 \int_{0}^{T}\int_{\R^N}u(t,x)\Bigl[\partial_t\phi(t,x)+\bigl(L^{\text{ou}}\bigr)^*\phi(t,x)\Bigr]
+f(t,x)\phi(t,x)\, dxdt\\
+\int_{\R^N}u_0(x)\phi(0,x)\, dx\, = \, 0,
\end{multline}
where $\bigl(L^{\text{ou}}\bigr)^*$ denotes the formal adjoint of $L^{\text{ou}}$ on $L^2(\R^N)$ given in Equation
\eqref{eq:def_adjoint_op}.

Similarly to the elliptic setting, we show firstly the weak well-posedness of the Cauchy problem \eqref{eq:Parabolic_IPDE}.

\begin{theorem}
\label{thm:well_posedness_parabolic}
Fixed $T>0$, let $u_0$ be a function in $C_b(\R^N)$ and $f$ in $L^\infty\bigl(0,T;C_b(\R^N)\bigr)$. Then, the function ${u\colon [0,T]\times\R^N\to
\R}$ given by
\begin{equation}\label{eq:representation_parabolic_sol}
u(t,x) \, := \, P_tu_0(x) + \int_{0}^{t}P_{t-s}f(s,x) \, ds, \quad (t,x) \in [0,T]\times\R^N,
\end{equation}
is the unique weak solution of the Cauchy problem \eqref{eq:Parabolic_IPDE}.
\end{theorem}
\begin{proof}
\emph{Existence.} We start considering a "regularized" version of the coefficients appearing in Equation \eqref{eq:Parabolic_IPDE}.
Namely, we consider a family $\{u_{0,m}\}_{m\in \N}$ in $C^\infty_b(\R^N)$ such that $u_{0,m}\to u_0$ uniformly in
$x$ and a family $\{f_m\}_{m\in \N}$ in $L^\infty\bigl(0,T;C^\infty_b(\R^N)\bigr)$ such that $f_m\to f$ uniformly in $t$ and $x$. They can
be obtained through standard mollification methods in space.\newline
Fixed $m$ in $\N$, we denote now by $u_m\colon[0,T]\times\R^n\to\R$ the function given by
\[u_m(t,x) \, := \, P_tu_{0,m}(x) + \int_{0}^{t}P_{t-s}f_m(s,x) \, ds, \quad t \in [0,T], x \in \R^N.\]
On the one hand, we use again that $\partial_t(P_tu_m)(t,x)=L^{\text{ou}}P_tu_m(t,x)$ for any $(t,x)$ in
$[0,T]\times\R^N$ to check that $u_m$ is indeed a \emph{classical} solution of the "regularized" Cauchy Problem:
\[
\begin{cases}
  \partial_tu_m(t,x) \, = \, L^{\text{ou}}u_m(t,x)+f_m(t,x), \quad (t,x) \in (0,T)\times \R^N; \\
  u_m(0,x) \, = \, u_{0,m}(x), \quad x \in \R^N.
\end{cases}
\]
On the other hand, we exploit the linearity and the continuity of the semi-group $P_t$ on $C_b(\R^N)$ to show that
\[u_m\, = \, P_tu_{0,m}(x) + \int_{0}^{t}P_{t-s}f_m(s,x) \, ds \,  \overset{m}{\to} \, P_tu_0(x) + \int_{0}^{t}P_{t-s}f(s,x) \, ds \, = \, u,\]
uniformly in $t$ and $x$, where $u$ is the function given in \eqref{eq:representation_parabolic_sol}. \newline
We fix now a test function $\phi$ in $C^\infty_0\bigl([0,T)\times\R^N\bigr)$ and we then notice that
\[\int_{0}^{T}\int_{\R^N}\phi(t,y)\Bigl(\partial_t-L^{\text{ou}}\Bigr)u_m(t,y) \, dydt \, = \,
\int_{0}^{T}\int_{\R^N}\phi(t,y)f_m(t,y) \, dydt.\]
An integration by parts allows now to move the operator to the test function, being careful to remember that $u_m(0,\cdot)=u_{0,m}(\cdot)$.
Indeed, it holds that
\begin{multline}\label{TOweak}
-\int_{0}^{T}\int_{\R^N}\Bigl(\partial_t+\bigl(L^{\text{ou}}\bigr)^*\Bigr)\phi(t,y)u_m(t,y) \, dydt \\
= \,
\int_{\R^N}\phi(0,y)u_{0,m}(y) \, dy + \int_{0}^{T}\int_{\R^N}\phi(t,y)f_m(t,y) \, dydt,
\end{multline}
where $\bigl(L^{\text{ou}}\bigr)^*$ denotes the formal adjoint of $L^{\text{ou}}$ on $L^2(\R^N)$.\newline
We would like now to go back to the solution $u$, letting $m$ go to infinity. We start rewriting the right-hand side term of \eqref{TOweak} as $R^1_m+R^2_m$,
where
\begin{align*}
R^1_m \, &:= \,\int_{\R^N}\phi(0,y)u_{0,m}(y) \, dy;   \\
R^2_m \, &:= \, \int_{0}^{T}\int_{\R^N}\phi(t,y)f_m(t,y) \, dydt.
\end{align*}
We can rewrite $R^2_m$ as
\[R^2_m \, = \, \int_{0}^{T}\int_{\R^N}\phi(t,y)f(t,y) \, dydt + \int_{0}^{T}\int_{\R^N}\phi(t,y)\bigl[f_m-f\bigr](t,y) \, dydt.\]
Exploiting that, by assumption, $f_m$ converges to $f$ uniformly in $t$ and $x$, it is easy to see that the second
contribution above converges to $0$. A similar argument can be used to show that
\[\int_{\R^N}\phi(0,y)u_{0,m}(y) \, dy \, \overset{m}{\to} \, \int_{\R^N}\phi(0,y)u_{0}(y) \, dy.\]
On the other hand, we can rewrite the left-hand side of Equation \eqref{TOweak} as
\begin{multline*}
-\int_{0}^{T}\int_{\R^N}\Bigl(\partial_t+\bigl(L^{\text{ou}}\bigr)^*\Bigr)\phi(t,y)u_m(t,y) \, dydt \\
= \,
-\int_{0}^{T}\int_{\R^N}\Bigl(\partial_t+\bigl(L^{\text{ou}}\bigr)^*\Bigr)\phi(t,y)u(t,y) \, dydt +L^1_m+L^2_m+L^3_m,
\end{multline*}
where we have denoted
\begin{align} \notag
L^1_m \, &:= \, \int_{0}^{T}\int_{\R^N} \bigl[\frac{1}{2}\text{Tr}\bigl(BQB^\ast D^2_y\phi(t,y)\bigr)+\langle Ay+Bb,D_y\phi(t,y)\rangle+\text{Tr}(A)\phi(t,y)\bigr]\\\notag
&\qquad\qquad \qquad\qquad \qquad\qquad \qquad\qquad \qquad\qquad\qquad \qquad \qquad\times[u_m-u](t,y) \, dydt;   \\
L^2_m \, &:= \, \int_{0}^{T}\int_{\R^N}\partial_t\phi(t,y)[u-u_m](t,y) \, dydt; \label{eq:remainder_in_existence}\\
L^3_m \, &:= \, \int_{0}^{T}\int_{\R^N}\Bigl[\int_{\R^d_0}\phi(t,y-Bz)-\phi(t,y)+\langle
D_y\phi(t,y),Bz\rangle\mathds{1}_{B(0,1)}(z)\, \nu(dz)\Bigr]\notag \\
&\qquad\qquad \qquad\qquad \qquad\qquad \qquad\qquad \qquad\qquad \qquad\qquad \qquad\times[u-u_m](t,y)dydt. \notag
\end{align}
To conclude, we need to show that the remainder $L^1_m+L^2_m+L^3_m$ is negligible, if $m$ goes to infinity. Exploiting that $\phi$ has a
compact support and that $\Vert u_m - u\Vert_\infty\overset{m}{\to} 0$, it is easy to show that
$\vert L^1_m+L^2_m\vert \overset{m}{\to} 0$.\newline
In order to control $L^3_m$, we need firstly to decompose it as $L^{3,1}_m+L^{3,2}_m$, where
\begin{align*}
L^{3,1}_m \, &:= \, \int_{0}^{T}\int_{\R^N}\bigl[u(t,y)-u_m(t,y)\bigr]\\
&\qquad \qquad \qquad \times\Bigl[\int_{0<\vert z\vert<1}\phi(t,y-Bz)-\phi(t,y)+\langle
D_y\phi(t,y),Bz\rangle\, \nu(dz)\Bigr]dydt;   \\
L^{3,2}_m \, &:= \, \int_{0}^{T}\int_{\R^N}\bigl[u(t,y)-u_m(t,y)\bigr]\Bigl[\int_{\vert z\vert>1}\phi(t,y-Bz)-\phi(t,y)\, \nu(dz)\Bigr]dydt.
\end{align*}
The second term $L^{3,2}_m$ can be controlled easily using the Fubini Theorem. Indeed, denoting by $K$ the support of $\phi$ and by $\lambda$ the Lebesgue measure on $\R^N$, we notice that
\[
\begin{split}
\vert L^{3,2}_m\vert \, &\le \, \Vert u-u_m\Vert_\infty\int_{0}^{T}\int_{\vert z\vert>1}\int_{\R^N}\vert\phi(t,y-Bz)-\phi(t,y)\vert\,
dy\nu(dz)dt\\
&\le \,  CT2\lambda(K)\nu\bigl(B^c(0,1)\bigr)\Vert u-u_m\Vert_\infty.
\end{split}
\]
Exploiting that $\nu\bigl(B^c(0,1)\bigr)$ is finite since $\nu$ is a L\'evy measure, we can then conclude that $\vert L^{3,2}_m\vert$ tends to zero if $m$ goes to infinity.\newline
The argument for $L^{3,1}_m$ is similar but we need firstly to apply a Taylor expansion twice to make a term $\vert z\vert^2$ appear in the
integral and exploit that $\vert z\vert^2\nu(dz)$ is finite on $B(0,1)$.

\emph{Uniqueness.} This proof will follow essentially the same arguments as for Theorem \ref{thm:well_posedness_elliptic}.\newline
Let $u$ be any weak solution of Cauchy problem \eqref{eq:Parabolic_IPDE} with $u_0=f=0$. We are going to show
that $u=0$. \newline
We start considering a mollyfing sequence $\{\rho_m\}_{m \in \N}$ in $C^\infty_c((0,T)\times\R^N)$. Denoting for simplicity $u_m(t,x)=u\ast \rho_m(t,x)$, we then notice that $u_m$ is
continuously differentiable in time and that $u_m(0,x)=0$. It makes sense to define now the function
\begin{equation}\label{Proof:eq:def_fm}
f_m(t,x) \, := \, \partial_t u_m(t,x)-L^{\text{ou}}u_m(t,x).
\end{equation}
Moreover, we can truncate $f_m$ and $u_m$ if necessary, so that they are integrable with integrable Fourier transform. Then, the same
reasoning in Equations \eqref{proof:Schauder_elliptic1}, \eqref{proof:Schauder_elliptic2} allows us to write that
\[
\begin{cases}
   \partial_t \widehat{u}_{m}(t,\xi)- \Phi^{\text{ou}}(\xi)\widehat{u}_m(t,\xi) \, = \, \widehat{f}_{m}(t,\xi), \\
    \widehat{u}_{m}(0,\xi) \, = \, 0.
  \end{cases}
\]
The above equation can be easily solved integrating in time, giving the following representation:
\[\widehat{u}_m(t,\xi) \, = \, \int_{0}^{t} e^{(t-s)\Phi^{\text{ou}}(\xi)}
\widehat{f}_m (s,\xi)\, ds.\]
In order to go back to $u_m$, we apply now the inverse Fourier transform to write that
\[u_m(t,x) \, = \, \int_{0}^{t}P_{t-s}f_m(s,x) \, ds.\]
The contraction property of $P_t$  allows us to conclude that $\Vert u_m \Vert_\infty \, \le \, C \Vert f_m\Vert_\infty$.
Letting $m$ goes to zero, we obtain the desired result. Indeed, it is possible to show that
\[\lim_{m \to \infty}\Vert f_m \Vert_{\infty} \, = \,  0,\]
relying on the same reasonings used in the analogous elliptic case
(Theorem \ref{thm:well_posedness_elliptic}).
\end{proof}

The next two conclusive theorems provide the Schauder estimates both in the elliptic and in the parabolic setting.

\begin{theorem}[Elliptic Schauder estimates]
Fixed $\lambda>0$ and $\beta$ in $(0,1)$, let $g$ be in $C^{\alpha+\beta}_{b,d}(\R^N)$. Then, the distributional solution $u$ of Equation
\eqref{eq:Elliptic_IPDE} is in $C^{\beta}_{b,d}(\R^N)$ and there exists a positive constant $C$ such that
\begin{equation}\label{eq:Schauder_Estim_elliptic}
\Vert u \Vert_{C^{\alpha+\beta}_{b,d}} \, \le \, C\bigl(1+\frac{1}{\lambda}\bigr)\Vert g \Vert_{C^{\beta}_{b,d}}.
\end{equation}
\end{theorem}
\begin{proof}
Thanks to Theorem \ref{thm:well_posedness_elliptic}, we know that the unique solution $u$ of the elliptic equation \eqref{eq:Elliptic_IPDE} is given in \eqref{eq:representation_elliptic_sol}. In order to show that such a function $u$ satisfies Schauder estimates
\eqref{eq:Schauder_Estim_elliptic}, we exploit again the equivalent norm defined in \eqref{eq:def_equivalent_norm} of Lemma \ref{lemma:def_equivalent_norm}. Namely, we fix $h$ in $\llbracket 1,n\rrbracket$ and $x_0$ in $\R^N$ and we show that
\[\vert \Delta^3_{x_0}u(z)\vert \, = \, \Bigl{\vert}\int_{0}^{\infty}e^{-\lambda t}\Delta^3_{x_0}\bigl(P_tg\bigr)(z) \, dt\Bigr{\vert} \,
\le \, C\Vert g \Vert_{C^\beta_{b,d}}\vert z\vert^{\frac{\alpha+\beta}{1+\alpha(h-1)}},\quad z\in E_h(\R^N),\]
for some constant $C>0$ independent from $x_0$. For $\vert z \vert\ge 1$, it can be easily obtained from the contraction property of $P_t$
on $B_b(\R^N)$:
\begin{equation}\label{Proof:Schauder:z>1}
\Bigl{\vert}\int_{0}^{\infty}e^{-\lambda t}\Delta^3_{x_0}\bigl(P_tg\bigr)(z) \, dt\Bigr{\vert} \, \le \, 3\int_{0}^{\infty}e^{-\lambda t\Vert P_tg\Vert_\infty\, dt} \, \le \,  \frac{3}{\lambda}\Vert g \Vert_\infty\vert z \vert^{\frac{\alpha+\beta}{1+\alpha(h-1)}}.
\end{equation}
When $\vert z \vert\le 1$, we start fixing a \emph{transition time} $t_0$ given by
\begin{equation}\label{eq:def:time_t0}
t_0 \:=\, \vert z \vert^{\frac{\alpha}{1+\alpha(h-1)}}.
\end{equation}
Notably, $t_0$ represents the transition time between the diagonal and the off-diagonal regime, accordingly to the intrinsic time scales of the system.
We then decompose $\Delta^3_{x_0}u(z)$ as $R_1(z)+R_2(z)$, where
\begin{align*}
  R_1(z) \, &:= \, \int_{0}^{t_0}e^{-\lambda t}\Delta^3_{x_0}\bigl(P_tg\bigr)(z) \, dt; \\
  R_2(z) \, &:= \, \int_{t_0}^{\infty}e^{-\lambda t}\Delta^3_{x_0}\bigl(P_tg\bigr)(z) \, dt.
\end{align*}
The first component $R_1$ is controlled easily using Corollary \ref{corollary:continuity_between_holder} for $\beta=\gamma$.
Indeed,
\begin{equation}\label{Proof:COntrol_tildeLambda1}
\vert R_1(z)\vert \, \le \, \int_{0}^{t_0}\vert\Delta^3_{x_0}\bigl(P_tg\bigr)(z)\vert \, dt \, \le \, \vert z \vert^{\frac{\beta}{1+\alpha(h-1)}}\int_{0}^{t_0}\Vert P_tg\Vert_{C^{\beta}_{b,d}}\, dt \, \le \, C\Vert g \Vert_{C^\beta_{b,d}}\vert z
\vert^{\frac{\alpha+\beta}{1+\alpha(h-1)}}.
\end{equation}
On the other hand, the control for $R_2$ can be obtained following Equation \eqref{Proof:eq:Taylor_Expansion} in order to write that
\begin{equation}\label{Proof:COntrol_tildeLambda2}
\begin{split}
\vert R_2(z) \vert \, &\le \,C\Vert g \Vert_{C^\beta_{b,d}}\vert z \vert^3 \int_{t_0}^{\infty}e^{-\lambda
t}\bigl(1+t^{\frac{\beta-3(1+\alpha(h-1))}{\alpha}}\bigr)\, dt \\
&\le C\Vert g \Vert_{C^\beta_{b,d}}\vert z \vert^3\bigl(\lambda^{-1}+\vert
z\vert^{\frac{\alpha+\beta-3(1+\alpha(h-1))}{1+\alpha(h-1)}}\bigr) \\
&\le C{ \bigl(1+\frac{1}{\lambda}\bigr)}\Vert g \Vert_{C^\beta_{b,d}}\vert z \vert^{\frac{\alpha+\beta}{1+\alpha(h-1)}},
\end{split}
\end{equation}
where, in the last step, we exploited that $\vert z \vert \le 1$.
\end{proof}

\begin{theorem}[Parabolic Schauder estimates]
\label{thm:Parabolic_Schauder_Estimates}
Fixed $T>0$ and $\beta$ in $(0,1)$, let $u_0$ be in $C^{\alpha+\beta}_{b,d}(\R^N)$ and $f$ in
$L^\infty\bigl(0,T;C^{\beta}_{b,d}(\R^N)\bigr)$.
Then, the weak solution $u$ of Cauchy Problem \eqref{eq:Parabolic_IPDE} is in $L^\infty\bigl(0,T;C^{\alpha+\beta}_{b,d}(\R^N)\bigr)$ and
there exists a constant $C:=C(T)>0$ such that
\begin{equation}\label{eq:Schauder_Estimate_parabolic}
\Vert u \Vert_{L^\infty(C^{\alpha+\beta}_{b,d})} \, \le \, C\bigl[\Vert u_0 \Vert_{C^{\alpha+\beta}_{b,d}}+\Vert f
\Vert_{L^\infty(C^{\beta}_{b,d})} \bigr].
\end{equation}
\end{theorem}
\begin{proof}
We are going to show that any function $u$ given by Equation \eqref{eq:representation_parabolic_sol}
satisfies the Schauder Estimates \eqref{eq:Schauder_Estimate_parabolic}. We start splitting the function $u$ in $u_1+u_2$, where
\begin{align}\label{Proof:Decomposition_Parabolic}
  u_1(t,x) &:= P_tu_0(x); \\
  u_2(t,x) &:= \int_{0}^{t}P_{s}f(t-s,x) \, ds.
\end{align}
Corollary \ref{corollary:continuity_between_holder} allows then to control $u_1$ in the
following way:
\[\Vert u_1 \Vert_{L^{\infty}(C^{\alpha+\beta}_{b,d})} \, = \, \sup_{t \in [0,T]}\Vert P_tu_0\Vert_{C^{\alpha+\beta}_{b,d}} \, \le \, C\Vert
u_0 \Vert_{C^{\alpha+\beta}_{b,d}}.\]
In order to deal with the contribution $u_2$, we will follow essentially the same reasoning for the Schauder Estimates in the elliptic
setting. Namely, we use again the equivalent norm defined in \eqref{eq:def_equivalent_norm} of Lemma \ref{lemma:def_equivalent_norm} in order to estimate
\[\Vert u_2\Vert_{L^\infty(C^{\alpha+\beta}_{b,d})} \, \le \, C\Vert f\Vert_{L^\infty(C^{\beta}_{b,d})}.\]
Fixed $h$ in $\llbracket 1,n\rrbracket$ and $x_0$ in $\R^N$, our aim is to show that
\[\vert \Delta^3_{x_0}u_2(z)\vert \, = \, \Bigl{\vert}\int_{0}^{t}\Delta^3_{x_0}\bigl(P_{t-s}f\bigr)(s,z) \, ds\Bigr{\vert} \,
\le \, C\Vert f \Vert_{L^\infty(C^\beta_{b,d}\beta)}\vert z\vert^{\frac{\alpha+\beta}{1+\alpha(h-1)}},\quad z\in E_h(\R^N),\]
for some constant $C>0$ independent from $x_0$. When $\vert z \vert\ge 1$, it can be obtained easily from the contraction property of $P_t$
on $C_b(\R^N)$ as in \eqref{Proof:Schauder:z>1}. For $\vert z \vert\le 1$, we fix again the transition time $t_0$ given in \eqref{eq:def:time_t0}
and we then decompose $\Delta^3_{x_0}u_2(t,z)$ as $\tilde{R}_1(t,z)+\tilde{R}_2(t,z)$, where
\begin{align*}
  \tilde{R}_1(t,z) \, &:= \, \int_{0}^{t\wedge t_0}\Delta^3_{x_0}\bigl(P_{s}f\bigr)(t-s,z) \, ds: \\
  \tilde{R}_2(t,z) \, &:= \, \int_{t\wedge t_0}^{t}\Delta^3_{x_0}\bigl(P_{s}f\bigr)(t-s,z) \, ds.
\end{align*}
The first component $R_1$ can be controlled easily as in \eqref{Proof:COntrol_tildeLambda1}:
\[
\begin{split}
\vert \tilde{R}_1(t,z)\vert \, &\le \, \int_{0}^{t\wedge t_0}\vert\Delta^3_{x_0}\bigl(P_{s}f\bigr)(t-s,z)\vert \, ds \\
&\le \, \vert z \vert^{\frac{\beta}{1+\alpha(h-1)}}\int_{0}^{t\wedge t_0} \Vert P_{s}f(t-s,\cdot)\Vert_{C^{\beta}_{b,d}}\, ds \\
&\le \, C\Vert f \Vert_{L^\infty(C^\beta_{b,d})}\vert z \vert^{\frac{\alpha+\beta}{1+\alpha(h-1)}}.
\end{split}
\]
On the other hand, the control for $R_2$ is obtained following the same steps used in Equation \eqref{Proof:COntrol_tildeLambda2}. Namely,
\[
\begin{split}
\vert \tilde{R}_2(t,z) \vert \, &\le \,C\Vert f \Vert_{L^\infty(C^\beta_{b,d})}\vert z \vert^3 \int_{t\wedge
t_0}^{\infty}\bigl(1+s^{\frac{\beta-3(1+\alpha(h-1))}{\alpha}}\bigr)\, ds \\
&\le C\Vert f \Vert_{L^\infty(C^\beta_{b,d})}\vert z \vert^{\frac{\alpha+\beta}{1+\alpha(h-1)}}.
\end{split}
\]
\end{proof}

\setcounter{equation}{0}
\section{Extensions to time-dependent operators}
\fancyhead[RO]{Section \thesection. Extensions to time dependent operators}
\label{Sec:Levy:Extensions_time_dependent}
In this final section, we would like to show some possible extensions of our method in order to include more general
operators with non-linear, space-time dependent coefficients. Even in this framework, we will prove the well-posedness of the parabolic
Cauchy problem and show the associated Schauder estimates. \newline
Following \cite{Krylov:Priola10}, our first step is to consider a time-dependent Ornstein-Uhlenbeck operator of the following form:
\begin{multline*}
L^{\text{ou}}_t\phi(t,x) \, := \,
\frac{1}{2}\text{Tr}\bigl(B_tQB_t^\ast D^2\phi(x)\bigr) +\langle A_tx,D \phi(x)\rangle \\
+ \int_{\R^d_0}\bigl[\phi(x+B_tz)-\phi(x)-\langle D_x\phi(x), B_tz\rangle\mathds{1}_{B(0,1)}(z) \bigr] \,\nu(dz),
\end{multline*}
where $B_t:=B\sigma_0(t)$ and $A_t$, $\sigma_0(t)$ are two time-dependent matrices in $\R^N\otimes \R^N$ and $\R^d\otimes\R^d$, respectively.
From this point further, we assume that the matrices $A_t$, $\sigma_0(t)$ are measurable in time and that they satisfy the following conditions:
\begin{description}
  \item[{[tK]}] for any fixed $t$ in $[0,T]$, it holds that  $N \, = \, \text{rank}\bigl[B,A_tB,\dots,A^{N-1}_tB\bigr]$;
  \item[{[B]}] the matrix $A_t$ is bounded in time, i.e.\ there exists a constant $\eta>0$ such that
  \[\vert A_t \xi\vert \, \le \,  \eta \vert \xi \vert, \quad \xi \in \R^N;\]
  \item[{[UE]}] the matrix $\sigma_0$ is uniformly elliptic, i.e.\ it holds that
  \[\eta^{-1}\vert \xi \vert^2 \, \le \, \langle \sigma_0(t)\xi,\xi \rangle \, \le \, \eta \vert \xi \vert^2 , \quad (t,\xi) \in [0,T]\times \R^d. \]
\end{description}
It is important to highlight already that this new "time-dependent" version [\textbf{tK}] of the Kalman rank condition [\textbf{K}] allows us to reproduce the same reasonings of Section $2$. In particular, the anisotropic distance $\mathbf{d}$ and the Zygmund-H\"older spaces $C^\beta_{b,d}(\R^N)$ can be constructed under these assumptions, even if only at any \emph{fixed} time $t$. A priori, the number of sub-divisions of the space $\R^N$ may change for different times, leading to consider a time-dependent $n(t)$ in Equation \eqref{eq:def_of_n_Levy} and, consequently, time-dependent anisotropic distances and H\"older spaces. We will however drop the subscript in $t$ below since it does not add any difficulty in the arguments but it may damage the readability of the article.

\begin{prop}
\label{thm:well_posedness_parabolic_time}
Let $u_0$ be in $C_b(\R^N)$ and $f$ in $L^\infty\bigl(0,T;C_b(\R^N)\bigr)$. Then, there exists a unique solution $u\colon [0,T]\times\R^N\to
\R$ of the following Cauchy problem:
\begin{equation}\label{eq:Parabolic_IPDE:time}
\begin{cases}
  \partial_tu(t,x) \, = \, L^{\text{ou}}_tu(t,x)+f(t,x), \quad (t,x) \in (0,T)\times \R^N; \\
  u(0,x) \, = \, u_0(x), \quad x \in \R^N.
\end{cases}
\end{equation}
Furthermore, if $u_0$ is in $C^{\alpha+\beta}_{b,d}(\R^N)$ and $f$ in
$L^\infty\bigl(0,T;C^{\beta}_{b,d}(\R^N)\bigr)$, then the solution $u$ is in $L^\infty\bigl(0,T;C^{\alpha+\beta}_{b,d}(\R^N)\bigr)$ and
there exists a constant $C:=C(T,\eta)>0$ such that
\begin{equation}\label{eq:Schauder_Estimate_parabolic:time}
\Vert u \Vert_{L^\infty(C^{\alpha+\beta}_{b,d})} \, \le \, C\bigl[\Vert u_0 \Vert_{C^{\alpha+\beta}_{b,d}}+\Vert f
\Vert_{L^\infty(C^{\beta}_{b,d})} \bigr].
\end{equation}
\end{prop}
\begin{proof}
The proof of this result can be obtained mimicking the arguments already presented in the first part of the article with some slight modifications. The main difference is the introduction of the resolvent $\mathcal{R}_{s,t}$ associated with the matrix $A_t$ in place of the matrix exponential $e^{tA}$. Namely, $\mathcal{R}_{t,u}$ is a time-dependent matrix in $\R^N\otimes \R^N$ that is solution of the following ODE:
\begin{equation}\label{eq:def_resolvent}
\begin{cases}
  \partial_t \mathcal{R}_{t,u} \, = \, A_t \mathcal{R}_{t,u} ,\quad t \in [u,T]; \\
  \mathcal{R}_{u,u} \, = \, \Id_{N\times N}.
\end{cases}
  \end{equation}
As said before, Section $2$ follows exactly in the same manner as above except for Lemma \ref{lemma:Decomposition_exp_A} (structure of the resolvent), whose proof can be found in \cite{Huang:Menozzi16}, Lemmas $5.1$ and $5.2$.
The arguments in Section $3$ and $4$ can be applied again, even if the formulation of some objects presented there changes slightly. For example in Equation \eqref{eq:Def_OU_Process}, the N-dimensional Ornstein-Uhlenbeck process $\{X_t\}_{t\ge 0}$ driven by $B_tZ_t$ should be now represented by
\[X_t \, = \, \mathcal{R}_{t,0}x+\int_{0}^{t}\mathcal{R}_{t,u}B_u \, dZ_u, \quad t\ge 0, x \in \R^N.\]
Finally in Section $5$, the uniform ellipticity [\textbf{UE}] of $\sigma_0(t)$ and the boundedness [\textbf{B}] of $A_t$ allow us to control the remainder terms appearing in Equation \eqref{eq:remainder_in_existence} as done above and thus, to conclude as in Theorems \ref{thm:well_posedness_parabolic} and \ref{thm:Parabolic_Schauder_Estimates}.
\end{proof}

Once we have shown our results for the time-dependent Ornstein-Uhlenbeck operator $L^{\text{ou}}_t$, we add now a non-linearity to the problem, even if only dependent in time. Namely, we are interested in operators of the following form:
\begin{equation}\label{eq:def_operator_Lt}
L_t\phi(t,x) \, := \, L^{\text{ou}}_t\phi(t,x)+\langle F_0(t),D_x\phi(x)\rangle -c_0(t)\phi(x), \quad (t,x) \in [0,T]\times
\R^N,
\end{equation}
where $c_0\colon [0,T]\to \R$ and $F_0\colon [0,T]\to \R^N$ are two functions. For any sufficiently regular function $\phi\colon[0,T]\to \R$, we are going to denote
\begin{equation}\label{eq:relation_weak_solutions}
\mathcal{T}\phi(t,x) \, := \, e^{-\int_{0}^{t}c_0(s) \,ds}\phi\Bigl(t,x+\int_{0}^{t}F_0(s) \, ds\Bigr), \quad (t,x) \in [0,T]\times \R^N.
\end{equation}
We will see in the next result that the "operator" $\mathcal{T}$ transforms solutions of the Cauchy problem associated with
$\mathcal{L}^\text{ou}_t$ to solutions of the Cauchy problem driven by $L_t$, even if for a modified drift $\mathcal{T}f$.

\begin{lemma}
\label{lemma:relation_weak_solutions}
Fixed $T>0$, let $u_0$ be in $C_b(\R^N)$, $f$ in $L^\infty\bigl(0,T;C_b(\R^N)\bigr)$ and $c_0$, $F_0$ in $C_b([0,T])$. Then, a function
$u\colon [0,T]\times \R^N\to \R$ is a weak solution of Cauchy Problem \eqref{eq:Parabolic_IPDE:time} if and only if the function $v\colon
[0,T]\times \R^N\to \R$ given by $v(t,x)=\mathcal{T}u(t,x)$ is a weak solution of the following Cauchy problem:
\begin{equation}\label{eq:Parabolic_IPDE:ext0}
\begin{cases}
  \partial_tu(t,x) \, = \, L_tu(t,x)+\mathcal{T}f(t,x), \quad (t,x) \in (0,T)\times \R^N; \\
  u(0,x) \, = \, u_0(x), \quad x \in \R^N.
\end{cases}
\end{equation}
In particular, there exists a unique weak solution of Cauchy Problem \eqref{eq:Parabolic_IPDE:ext0}.
\end{lemma}
\begin{proof}
Given a weak solution $u$ of Cauchy problem \eqref{eq:Parabolic_IPDE:time}, we are going to show that the function $v$ given in
\eqref{eq:relation_weak_solutions} is indeed a weak solution of Cauchy Problem \eqref{eq:Parabolic_IPDE:ext0}. The inverse implication  can
be obtained in a similar manner and we will not prove it here. \newline
By mollification if necessary, we can take two sequences $\{c_m\}_{m\in \N}$ and $\{F_m\}_{m\in \N}$ in $C^\infty_b([0,T])$ such that $c_m\to c_0$
and $F_m \to F_0$ uniformly in $t$. Furthermore, we denote for simplicity
\[\tilde{c}_m(t) \, := \, \int_{0}^{t}c_m(s) \, ds; \quad \tilde{F}_m(t) \, := \, \int_{0}^{t}F_m(s) \, ds. \]
Given a test function $\phi$ in $C^\infty_c\bigl([0,T)\times \R^N\bigr)$, let us consider for any $m$ in $\N$, the following function
\[\psi_m(t,x) \, := \, e^{-\tilde{c}_{m}(t)}\phi(t,x-\tilde{F}_m(t)) \quad (t,x) \in [0,T]\times\R^N.\]
Since $\tilde{c}_{m}$ and $\tilde{F}_{m}$ are smooth and bounded, it is easy to check that $\psi_m$ is in $C^\infty_c\bigl([0,T)\times
\R^N\bigr)$. We can then use $\psi_m$ in Equation \eqref{eq:weak_solution_Parabolic} (with time-dependent $A_t$ and $B_t$) to show that
\[\int_{0}^{T}\int_{\R^N}\Bigl[\partial_t+\bigl(L^{\text{ou}}_t\bigr)^*\Bigr]\psi_m(t,y)u(t,y)+f(t,y) \, dydt +
\int_{\R^N}\psi_m(0,y)u_0(y) \, dy \, = \, 0.\]
A direct calculation then show that $\psi_m(0,y)=\phi(0,y)$ and
\begin{align*}
   \bigl(L^{\text{ou}}_t\bigr)^*\psi_m(t,y) \, &= \,
   e^{-\tilde{c}_m(t)}\bigl(L^{\text{ou}}_t\bigr)^*\phi(t,y-\tilde{F}_m(t));\\
  \partial_t\psi_m(t,y) \, &= \, e^{-\tilde{c}_m(t)}\Bigl[\partial_t\phi(t,y-\tilde{F}_m(t))-\langle F_m(t),D_y\phi(t,y-\tilde{F}_m(t))
  \rangle\\
  &\qquad\qquad\qquad\qquad\qquad\qquad\qquad\qquad\qquad-c_m(t)\phi(t,y-\tilde{F}_m(t))\Bigr].
\end{align*}
The above calculations and a change of variable then imply that
\begin{multline*}
\int_{0}^{T}\int_{\R^N}\Bigl[\Bigl(\partial_t+\bigl(L^{\text{ou}}_t\bigr)^*\Bigr)\phi(t,y)-\langle F_m(t),D_y\phi(t,y) \rangle -
c_m(t)\phi(t,y)\Bigr]\mathcal{T}_mu(t,y)\\
+\phi(t,y)\mathcal{T}_mf(t,y) \, dydt +\int_{\R^N}u_0(y)\phi(0,y) \, dy \, = \, 0,
\end{multline*}
where, analogously to \eqref{eq:relation_weak_solutions}, we have denoted for any function $\varphi\colon [0,T]\times \R^N\to \R$,
\[\mathcal{T}_m\varphi(t,y) \, := \, e^{-\tilde{c}_m(t)}\varphi(t,y+\tilde{F}_m(t)).\]
Following similar arguments exploited in the "existence" part in the proof of Theorem \ref{thm:well_posedness_parabolic}, i.e.\ exploiting
the compact support of $\phi$ and the uniform convergence of the coefficients, it is possible to show that the above expression converges,
when $m$ goes to infinity, to
\[\int_{0}^{T}\int_{\R^N}\Bigl[\partial_t+\bigl(L_t\bigr)^*\Bigr]\phi(t,y)v(t,y)+\mathcal{T}f(t,y) \, dxdt + \int_{\R^N}\phi(0,x)u_0(x) \,
dx]\, = \, 0\]
and thus, that $v$ is a weak solution of Cauchy problem \eqref{eq:Parabolic_IPDE:ext0}.
\end{proof}

Thanks to the previous lemma, we are now able to show the Schauder estimates for the solution $v$ of the Cauchy problem
\eqref{eq:Parabolic_IPDE:ext0} and, more importantly, without changing the constant $C$ appearing in Equation
\eqref{eq:Schauder_Estimate_parabolic:time}.

\begin{prop}
\label{thm:well_posedness_parabolic_Lt}
Fixed $T>0$, $\beta\in (0,1)$, let $u_0$ in $C^{\alpha+\beta}_{b,d}(\R^N)$, $f$ in $L^\infty\bigl(0,T;C^{\beta}_b(\R^N)\bigr)$ and $c_0$,
$F_0$ in $B_b([0,T])$. Then, the unique solution $v$ of Cauchy Problem \eqref{eq:Parabolic_IPDE:ext0} is in $L^\infty\bigl(0,T;
C^{\alpha+\beta}_b(\R^N)\bigr)$ and it holds that
\begin{equation}\label{eq:Schauder_estimates_Lt}
\Vert v \Vert_{L^\infty(C^{\alpha+\beta}_{b,d})} \, \le \, C\bigl[\Vert u_0 \Vert_{C^{\alpha+\beta}_{b,d}}+\Vert f
\Vert_{L^\infty(C^{\beta}_{b,d})}\bigr],
\end{equation}
where $C:=C(T,\eta)>0$ is the same constant appearing in Theorem \ref{thm:Parabolic_Schauder_Estimates}.
\end{prop}
\begin{proof}
We start denoting for simplicity
\[\tilde{c}_0(t) \, := \, \int_{0}^{t}c_0(s) \, ds \,\, \text{ and }\,\, \tilde{F}_0(t) \, := \, \int_{0}^{t}F_0(s) \, ds.\]
By Lemma \ref{lemma:relation_weak_solutions}, we know that if $v$ is a weak solution of Cauchy problem \eqref{eq:Parabolic_IPDE:ext0}, then
the function
\[u(t,x)\,:=\,e^{\tilde{c}_0(t)}v(t,x-\tilde{F}_0(t))\]
is the weak solution of Cauchy problem \eqref{eq:Parabolic_IPDE:time} with $\tilde{f}$ instead of
$f$, where
\[\tilde{f}(t,x) \, := \, e^{\tilde{c}_0(t)}f(t,x-\tilde{F}_0(t)), \quad (t,x) \in (0,T)\times\R^N.\]
Moreover, we have that $\tilde{f}$ is in $L^\infty\bigl(0,T;C^{\beta}_b(\R^N)\bigr)$. Considering, if necessary, a smaller time interval $[0,t]$ for some $t\le T$, it is not difficult to check from Proposition
\ref{thm:well_posedness_parabolic_time} that
\[\Vert e^{\tilde{c}_0(t)}v(t,\cdot-\tilde{F}_0(t))\Vert_{C^{\alpha+\beta}_{b,d}} \, \le \, C \bigl[\Vert u_0
\Vert_{C^{\alpha+\beta}_{b,d}}+\sup_{s \in [0,t]}\Vert e^{\tilde{c}_0(s)}f(s,\cdot-\tilde{F}_0(t))\bigr].\]
Using now the invariance of the H\"older norm under translations, we can show that
\[
\begin{split}
\Vert v(t,\cdot)\Vert_{C^{\alpha+\beta}_{b,d}} \, &\le \, C \bigl[e^{-\tilde{c}_0(t)}\Vert u_0
\Vert_{C^{\alpha+\beta}_{b,d}}+e^{-\tilde{c}_0(t)}\sup_{s \in [0,t]}\Vert e^{\tilde{c}_0(s)}f(s,\cdot)\bigr]\\
&\le C \bigl[\Vert u_0
\Vert_{C^{\alpha+\beta}_{b,d}}+\sup_{s \in [0,t]}\Vert f(s,\cdot)\bigr],
  \end{split}
\]
where in the last step we exploited that $\tilde{c}_0(t)$ is non-decreasing. Taking the supremum with respect to $t$ on both sides of the above
inequality, we obtain our result.
\end{proof}

\begin{remark}[About space-time dependent coefficients]
We briefly explain here how to extend the Schauder estimates \eqref{eq:Schauder_Estimate_parabolic} to a class of non-linear, space-time
dependent operators, whose coefficients are only locally H\"older continuous in space and may be unbounded. Namely, we are interested in
operators of the following form:
\begin{multline*}
L_t\phi(t,x) \, := \langle
F(t,x),D_x\phi(x)\rangle \\ +\int_{\R^d_0}\bigl[\phi(x+B\sigma(t,x)z)-\phi(x)-\langle D_x\phi(x), B\sigma(t,x)z\rangle\mathds{1}_{B(0,1)}(z) \bigr] \,\nu(dz),
\end{multline*}
where $B$ is as in \eqref{eq:Lancon_Pol_Levy} and $\sigma \colon [0,T]\times\R^N\to \R^d\otimes\R^d$, $F\colon [0,T]\times\R^N\to \R^N$ are two
measurable functions such that $F(t,0)$ is locally bounded in time and $\sigma$ satisfies assumption [\textbf{UE}] at any fixed $(t,x)$ in
$[0,T]\times \R^N$.\newline
We would like now the operator $L_t$ to present a similar "dynamical" behaviour as above, i.e.\ the transmission of
the smoothing effect of the L\'evy operator to the degenerate components of the system (cf.\ Example \ref{example_basic}). For this reason, we
suppose the following:
\begin{itemize}
  \item the drift $F=(F_1,\dots,F_n)$ is such that for any $i$ in $\llbracket 1,n\rrbracket$, $F_i$ depends only on time and on the last
      $n-(i-2)$ components, i.e.\ $F_i(t, x_{i-1},\dots,x_n)$;
  \item the matrices $D_{ x_{i-1}}F_i(t,x)$ have full rank $d_i$ at any fixed $(t,x)$ in $[0,T]\times \R^N$.
\end{itemize}
As said before, the functions $F$ and $\sigma$ are assumed to be only locally H\"older in space, uniformly in time. Namely, there exists a
positive constant $K_0$ such that
\begin{equation}\label{eq:Local_Holder}
\mathbf{d}\bigl(\sigma(t,x),\sigma(t,y)\bigr) \, \le \, K_0\mathbf{d}^\beta(x,y); \qquad \mathbf{d}\bigl(F_i(t,x),F_i(t,y)\bigr) \, \le \, K_0\mathbf{d}^{\beta+\gamma_i}(x,y)
\end{equation}
for any $i$ in $\llbracket 1,n \rrbracket$, any $t$ in $[0,T]$ and any $x,y$ in $\R^N$ such that $\mathbf{d}(x,y)\le 1$, where
\begin{equation}\label{Drift_assumptions1}
    \gamma_i \,:= \,
    \begin{cases}
        1+ \alpha(i-2), & \mbox{if } i>1; \\
        0, & \mbox{if } i=1.
   \end{cases}
   \end{equation}
We remark in particular that the function $F$ may be unbounded in space.\newline
In order to recover Schauder-type estimates even in this framework, we can follow a perturbative method firstly introduced in
\cite{Krylov:Priola10} that allows to exploit the already proven results for time-dependent operators.
Let us assume for the moment that $\sigma$ and $F$ are \emph{globally} H\"older continuous in space, i.e. they satisfy
\eqref{eq:Local_Holder} for any $x,y$ in $\R^N$. Informally speaking, the method links the
operator $L_t$ with the space independent operator $L_t$ defined in \eqref{eq:def_operator_Lt}, by "freezing" the coefficients of
$L_t$ along a \emph{reference path} $\theta\colon [0,T]\to \R^N$ given by
\[\theta_t \, := \, x_0 +\int_{t_0}^{t}F(s,\theta_s) \, ds,\]
for some $(t_0,x_0)$ in $[0,T]\times \R^N$. It is important to highlight that, since $F$ is only H\"older continuous, we need to fix one of
the possible paths satisfying the above dynamics. We point out that the deterministic flow $\theta_t$ associated with the drift $F$ is
introduced precisely to handle the possible unboundedness of $F$. We could then consider a proxy operator $L_t$ whose coefficients are given
by $\sigma_0(t):=\sigma(t,\theta_t)$, $F_0(t):=F(t,\theta_t)$ and
\[\bigl[A_t\bigr]_{i,j} \, = \,
\begin{cases}
  D_{ x_{i-1}}F_i(t,\theta_t), & \mbox{if } j=i-1; \\
  0, & \mbox{otherwise}
\end{cases}
\]
In particular, Theorem \ref{thm:well_posedness_parabolic_Lt} assures the well-posedness and the Schauder estimates for the Cauchy problem
associated with $L_t$.\newline
The final step of the proof would be to expand a solution $u$ of the Cauchy problem associated with $L_t$ around the proxy $L_t$ through
a Duhamel-like formula and finally show that the expansion error only brings a negligible contribution so that the Schauder estimates still
hold for the original problem.\newline
The a priori estimates for the expansion error are however quite involved (and they are the main reason why we have decided to not show here
the complete proof), since they rely on some non-trivial controls in appropriate Besov norms.\newline
In order to deal with coefficients that are only locally H\"older in space, we need in addition to introduce a "localized" version of the
above reasoning. It would be necessary to multiply a solution $u$ by a suitable bump function $\delta$ that localizes in space along the
deterministic flow $\theta_t$ that characterizes the proxy. Namely, to fix a smooth function $\rho$ that is equal to $1$ on $B(0,1/2)$
and vanishes outside $B(0,1)$ and define $\delta(t,x) := \rho(x-\theta_t)$. We would then follow the above method but with respect to the
"localized" solution
\[v(t,x) \, := \, \delta(t,x)u(t,x), \quad (t,x) \in [0,T]\times \R^N.\]
We suggest the interested reader to see \cite{Chaudru:Honore:Menozzi18_Sharp} for a detailed treatise of the argument in the degenerate
diffusive setting, \cite{Chaudru:Menozzi:Priola19} in the non-degenerate stable framework or \cite{Marino20} for the precise assumptions on
the coefficients.
\end{remark}

\setcounter{equation}{0}

\chapter{Weak well-posedness for degenerate SDEs driven by L\'evy processes}
\fancyhead[LE]{Chapter \thechapter. Weak well-posedness for degenerate L\'evy SDEs}
\label{Chap:Weak_Well-Posedness}
\paragraph{Abstract:}
We study the effects of the propagation of a non-degenerate L\'evy noise through a chain of deterministic differential equations whose coefficients are H\"older continuous and satisfy a weak H\"ormander-like condition. In particular, we assume some non-degeneracy with respect to the components which transmit the noise.
Moreover,  we characterize, for some specific dynamics, through suitable counter-examples, the almost sharp regularity exponents that ensure the weak well-posedness for the associated SDE. As a by-product of our approach, we also derive some Krylov-ype estimates for the density of the weak solutions of the considered SDE.

\section{Introduction}
\fancyhead[RO]{Section \thesection. Introduction}
We investigate the effects of the propagation of a $d$-dimensional L\'evy noise through a chain of $n\ge 2$ differential equations.
Namely, we are interested in a degenerate,  L\'evy-driven stochastic differential equation (SDE in short) of the following form:
\begin{equation}
\label{eq:Intro}
  \begin{cases}
    dX^1_t \, = \,\left[[A_t]_{1,1}X_t^1+\cdots+[A_t]_{1,n}X_t^n+ F_1(t,X^1_t,\dots,X^n_t)\right]dt + \sigma(t,X^1_{t-},\dots,X^n_{t-})dZ_t,\\
    dX^2_t \, = \, \left[ [A_t]_{2,1}X^1_t+\cdots+[A_t]_{2,n}X^n_t +F_2(t,X^2_t,\dots,X^n_t)\right] dt,\\
    dX^3_t \, = \, \left[ [A_t]_{3,2}X^2_t+\cdots+[A_t]_{2,n}X^n_t+F_3(t,X^3_t,\dots,X^n_t)\right]dt,\\
    \vdots\\
    dX^n_t \, = \, \left[ [A_t]_{n-1,n} X^{n-1}_t+[A_t]_{n,n} X^{n}_t+F_n(t,X^n_t)\right]dt,
  \end{cases}
\end{equation}
where for $i\in \llbracket 1,n\rrbracket$ ($\llbracket \cdot, \cdot \rrbracket$ denotes the set of all the integers in the interval),  $X_t^i$ is $\R^{d_i} $ valued, with $d_1=d $ and $ d_i\ge 1,\ i\in \llbracket 2,n\rrbracket$. Set $N=\sum_{i=1}^n d_i $. We suppose that the
$F_i\colon [0,+\infty)\times\R^{{\scriptscriptstyle{\sum_{j=i}^n}} d_j}\to \R^d$,  $\sigma\colon [0,+\infty)\times\R^N\to \R^d\otimes \R^d$ are  Borel and respectively locally bounded
and uniformly elliptic and bounded.

We also assume the entries $([A_t]_{ij})_{1\le i\le n,\ i-1\le j\le n} $ are Borel bounded and such that
 the blocks $[A_t]_{i,i-1}$ in $\R^{d_{i}}\otimes \R^{d_{i-1}}$, $2\le i \le n $ have rank $d_i$, uniformly in time.
 This is a kind of non-degeneracy assumption which can be viewed as weak H\"ormander-like condition. It actually precisely allows the noise to propagate into the system.

Eventually, the noise $\{Z_t\}_{t\ge 0}$ belongs to a class of $d$-dimensional, symmetric, L\'evy processes with suitable properties. In particular, to handle non trivial diffusion coefficients, we will assume that the L\'evy measure of $\{Z_t\}_{t\ge 0}$ is absolutely continuous with respect to the L\'evy measure of a \textcolor{black}{rotationally invariant} $\alpha$-stable process (with $\alpha$ in $(1,2]$) and its Radon-Nikodym derivative enjoys some natural properties. The class of processes $\{Z_t\}_{t\ge0}$ we can consider, includes for example, the tempered, the layered or the relativistic $\alpha$-stable processes. In the case of an additive noise, cylindrical stable processes could be handled as well.

Here, the major issue is linked with the specific degenerate framework we consider. Indeed, the noise only acts on the first component of the dynamics \eqref{eq:Intro} and it implies, in particular, that the random perturbation on the $i$-th
line of SDE \eqref{eq:Intro} only comes from the previous $(i-1)$-th one, through the non-degeneracy of the matrixes $[A_t]_{i,i-1}$. Hence, the smoothing effect associated with the L\'evy noise decreases along the chain,  making thus more and more difficult to regularize by noise the furthest lines of Equation \eqref{eq:Intro}. \newline
We nevertheless prove the weak well-posedness, i.e.\ the existence and uniqueness in law, for the above SDE \eqref{eq:Intro} when the drift $F=(F_1,\dots,F_n)$ and $ \sigma$ lie  in a suitable anisotropic H\"older space with multi-indices of regularity. We assume that $F_1$ and $\sigma $ have spatial H\"older regularity $\beta^1>0 $ with respect to the $j$-th variable. We highlight already that
we could have considered different regularity indexes $\beta_j^1 $ for the regularity of $F_1$ with respect to the $j$-th variable. We keep only one common index for notational simplicity. We also suppose that for fixed $j\in \llbracket2,n \rrbracket  $, $(F_2,\cdots,F_j) $ has H\"older regularity $\beta^j$  with respect to the $j$-th variable, where:
\[ \beta^j\,\, \in \,\, \Bigl(\frac{1+\alpha(j-2)}{1+\alpha(j-1)};1\Bigr].\]
We indeed recall that from the dynamics \eqref{eq:Intro} the variable $x_j$ does not appear in the chain after level $j$.

Furthermore, we will show through suitable counter-examples that the above threshold for the regularity exponents $\beta^j$ is ``almost'' sharp for a perturbation of the $j^{{\rm th}} $ level of the chain.
Such counter-examples are based on Peano-type dynamics adapted to our degenerate, fractional framework.

Models of the form \eqref{eq:Intro} naturally appear in various scientific contexts: in physics, for the
analysis of anomalous diffusions phenomena or for Hamiltonian models in turbulent regimes (see e.g.\ \cite{Baeumer:Benson:Meerschaert01}, \cite{Cushman:Park:Kleinfelter:Moroni05}, \cite{Eckmann:Pillet:Rey-Bellet99}); in mathematical
finance and econometrics, for example in the pricing of Asian options (see e.g.\ \cite{Jeamblanc:Yor:Chesney09}, \cite{Brockwell01}, \cite{Barndorff-Nielsen:Shephard01}). In particular, models that consider L\'evy noises, such as SDE \eqref{eq:Intro}, seem more natural and realistic for many
applications since they allow the presence of jumps.

\paragraph{Weak well-posedness for non-degenerate stable SDEs.}
The topic of weak well-posedness for non-degenerate (i.e.\ $d=N$) SDEs of the form:
\begin{equation}
\label{eq:non-deg_SDE}
X_t \, = \, x +\int_0^tF(X_s) ds +Z_t, \quad t\ge 0,
\end{equation}
where $\{Z_t\}_{t\ge0}$ is a symmetric $\alpha$-stable process on $\R^N$, has been widely studied in the last decades, especially in the diffusive, local setting, i.e.\ when $\alpha=2$ and $\{Z_t\}_{t\ge0}$ is a Brownian Motion, and it is now
quite well-understood. We can first refer to the seminal work \cite{book:Stroock:Varadhan79} where the Authors considered additionally a multiplicative noise with bounded drift and  non-degenerate, continuous in space diffusion coefficient. We recall moreover that in the framework of \eqref{eq:non-deg_SDE} with bounded drift, strong uniqueness also holds (cf. \cite{Veretennikov80}).\newline
SDEs like \eqref{eq:non-deg_SDE} with a proper $\alpha$-stable process ($\alpha<2$) were firstly investigated in \cite{Tanaka:Tsuchiya:Watanabe74} where the weak well-posedness was obtained for the one-dimensional case when the drift $F$ is bounded, continuous and the L\'evy exponent $\Phi$ of $\{Z_t\}_{t\ge 0}$ satisfies $\Re \Phi(\xi)^{-1}= 0(1/\vert \xi \vert)$ if $\vert \xi \vert \to \infty$. The multidimensional case ($d>1$) can be similarly obtained  following \cite{Komatsu83} if the drift is bounded, continuous and the law of $\{Z_t\}_{t\ge 0}$ admits a density with respect to the Lebesgue measure on $\R^d$. Equations as \eqref{eq:non-deg_SDE} with drift in some suitable $L^p$-spaces and non-degenerate noise were also considered in \cite{Jin18} (see also the references therein). We can eventually quote the recent work by Krylov who obtained even strong uniqueness for Brownian SDEs with drifts in critical $L^p$-spaces, see \cite{Krylov21}.

In recent years, SDEs driven by singular (distributional) drift have gained a lot of interest, especially in the Brownian setting, where they arise as a model for diffusions in
random media (see e.g. \cite{Flandoli:Russo:Wolf03},\cite{Flandoli:Russo:Wolf04},\cite{Flandoli:Issoglio:Russo17}, \cite{Delarue:Diel16},
\cite{Cannizzaro:Chouk18}).\newline
In the non-local $\alpha$-stable framework, a first work to appear was \cite{Athreya:Butkovsky:Mytnik18} where the authors considered the one-dimensional case with a time-homogeneous drift of (negative) H\"older regularity strictly greater than $(1-\alpha)/2$. We remark that in the one-dimensional framework, the regularity thresholds on the drift are the same for the strong and the weak well-posedness, since it is possible to exploit local time arguments (see also \cite{Bass:Chen01} in the diffusive setting). On the same side, the almost simultaneous works \cite{Ling:Zhao19} and \cite{Chaudru:Menozzi20} take into account time-homogeneous and time-inhomogeneous, respectively, distributional drift in general Besov spaces with suitable
conditions on the parameters. These results rely on Young integrals in order to give a meaningful sense to the dynamics. Beyond the Young regime, we instead
refer to \cite{kremp:Perkowski20} \textcolor{black}{where techniques such as paracontrolled products (which have also been popular in the recent developments in the SPDE theory)} are exploited to analyze the martingale problem associated \textcolor{black}{with} a
time-inhomogeneous drift of regularity index  strictly greater than $(2-2\alpha)/3$.

Moreover, we would like to remark that the above works concerned the so-called \emph{sub-critical} case, i.e.\ when $\alpha>1$. Indeed, SDEs like
\eqref{eq:non-deg_SDE} are much more difficult to handle if $\alpha\le 1$ since in this case, the noise does not dominate the system \textcolor{black}{for small time scales}. Two recent
works along this path are \cite{Zhao19} and \cite{Chaudru:Menozzi:Priola20} where \textcolor{black}{the authors} consider $\alpha<1$, $(1-\alpha)$-H\"older drift $F$ and $\alpha=1$,
continuous drift, respectively. We also mention that for H\"older drifts, the well-posedness of the associated martingale problem can be obtained following
\cite{Mikulevicius:Pragarauskas14} if $F$ is bounded or through the Schauder estimates given in \cite{Chaudru:Menozzi:Priola19} when $F$ is unbounded.

\paragraph{Regularization by noise in a degenerate setting.}
All the above results present a common phenomenon that, following the terminology in \cite{book:Flandoli11}, is usually called \emph{regularization by
noise}. This occurs when a deterministic ODE is ill-posed (for example if the drift is less than Lipschitz) but its stochastic counterpart
(SDE) is well posed in a strong or a weak sense. \newline
To obtain such phenomenon, the noise plays a fundamental role. A usual assumption is that the noise should act on every line of the dynamics, regularizing the
coefficients. It is then clear that in our degenerate framework, when the noise acts only on the first component of the chain \eqref{eq:Intro}, the situation is
even more delicate. In order to obtain some kind of regularization effect in this case, we need that the noise propagates through the system, reaching all the
lines of Equation \eqref{eq:Intro}. A typical assumption ensuring such type of behaviour is the so-called H\"ormander condition for hypoellipticity (cf.
\cite{Hormander67}).\newline
From the structure of the equation \eqref{eq:Intro} at hand, we will consider a \textit{weak} type H\"ormander condition, i.e. up to some regularization of the diffusion coefficient, the drift is needed to span the space through
Lie bracketing.


In the Hamiltonian setting $n=2$, when $\{Z_t\}_{t\ge0}$ is a Brownian Motion and for a more general, non-linear, drift than in \eqref{eq:Intro} which still satisfies a weak H\"ormander type condition,  Chaudru de Raynal showed in \cite{Chaudru18} that the associated SDE  is weakly well-posed as soon as the drift is H\"older continuous in the degenerate variable with regularity index strictly greater than $1/3$. It was also established through an appropriate counter-example,  that the $1/3$-threshold is (almost) sharp for the second component of the drift. Such a result has been then extended in \cite{Chaudru:Menozzi17} in order to consider the more general case of $n$ oscillators. Therein, the  regularity thresholds that ensure weak uniqueness also depend on the variable and the level of the chain. This seems intuitively clear, the further the variable in the oscillator chain, the larger its typical time scale, the weaker the regularity needed to regularize components which are above that variable in the chain. Also, some corresponding Krylov type estimates, giving existence and integrability properties of the density of the SDE are derived.
We can mention as well the recent work by Gerencs\'er \cite{gerencser20} who obtain similar regularization properties for the iterated time integrals of a fractional Brownian motion.

In the jump case, the situation is much more delicate. Within the proper regularization by noise framework (when the coefficients are less than Lipschitz continuous),  we cite \cite{Huang:Menozzi16} where the
Authors showed the weak well-posedness for \eqref{eq:Intro} with $F=0$ and a H\"older continuous diffusion coefficient, under  some constraints on the dimensions
$d$,$n$. In that framework, the Authors obtained as well same point-wise density estimates. The driving noises considered were stable and tempered stable processes.
\newline
Finally, we mention that it is possible to derive the weak well-posedness of dynamics \eqref{eq:Intro} via the martingale formulation, exploiting the Schauder
estimates given in \cite{Hao:Wu:Zhang19} for the kinetic model ($n=2$). In that work, the Authors actually characterized conditions for strong uniqueness, using Littlewood-Paley decomposition techniques.

We will here proceed through a perturbative approach. Namely, we will expand the formal generator associated with \eqref{eq:Intro} around the one of a well understood process, with possibly time inhomogeneous coefficients which are anyhow frozen in space. We will call such a process a \textit{proxy}. The most natural candidate to be a \textit{proxy} for \eqref{eq:Intro} is a degenerate Ornstein-Uhlenbeck process. In the case of time homogeneous coefficients, Priola and Zabczyk  established in \cite{Priola:Zabczyk09} existence of the density for such processes under the same previously indicated non-degeneracy conditions on the matrix $A$ (which turn out to be equivalent in that setting to the well known Kalman condition).
\paragraph{Intrinsic difficulties associated with large jumps.}
When $Z$ is a strict stable process, the  density of the corresponding degenerate Ornstein-Uhlenbeck process can somehow be related to the one of a multi-scale stable process which has however a very singular associated \textit{spectral measure} (spherical part of the $\alpha $-stable L\'evy measure) on $\mathbb S^{N-1}$, see e.g. \cite{Huang:Menozzi16}, \cite{Huang:Menozzi:Priola19} and Proposition \ref{prop:Decomposition_Process_X} below.
From Watanabe \cite{Watanabe07}, it is known that the tails of stable densities are highly related to the nature of this spectral measure. Specifically, the concentration properties worsen when the measure becomes singular. This renders delicate the characterization of the smoothing properties for the proxy, especially when it depends on parameters and that one would like to obtain  estimates which are uniform w.r.t. those parameters (see Proposition \ref{prop:Smoothing_effect} and Section \ref{SEZIONE_SPIEGAZIONE_CONGELAMENTO_E_ALTRI} below).

Even for smooth coefficients, the stable like jump setting is much more delicate to establish the existence of the density for (degenerate) SDEs. For multiplicative noises, we cannot indeed rely on the flow techniques considered in \cite{book:bichteler:Gravereaux:Jacod87} or \cite{book:kunita19} in the non-degenerate case, and L\'eandre in the degenerate one, see \cite{Leandre85},\cite{Leandre88}.
Still for smooth coefficients, we can refer to the work of Zhang \cite{Zhang14} who obtained existence and smoothness results for the density of equations of type \eqref{eq:Intro} in arbitrary dimension, for a possibly more general non linear drift, still satisfying a weak H\"ormander type condition when the driving process is a rotationally invariant stable process. The strategy therein is based on the  \textit{subordinated} Malliavin calculus, which consists in applying the \textit{usual} Malliavin calculus techniques on a Brownian motion observed along the path of an independent $\alpha $-stable subordinator. In whole generality a \textit{complete} version of the H\"ormander theorem in the jump case seems to lack.
We can refer to the work by Cass \cite{Cass09} who gets smoothness of the density in the weak H\"ormander framework under technical restrictions.


\subsection{Complete model and assumptions}
\label{Sec:Well-Pos:Complete_model_assumptions}
Let us now specify the assumptions on equation \eqref{eq:Intro} that we rewrite in the shortened form:
\begin{equation}\label{eq:SDE}
dX_t \, = \, G(t,X_t)dt + B\sigma(t,X_{t-})dZ_t, \quad t\ge 0,
\end{equation}
where  $B$ is the embedding from $\R^d$ to $\R^N$ given in matricial form as
\[B:=\begin{pmatrix}
            I_{d\times d}, & 0_{d\times (N-d)},
\end{pmatrix}^t\]
and $G(t,x)=A_tx+F(t,x) $ with:
      \begin{equation}\label{eq:def_matrix_A}
A_t \, := \, \begin{pmatrix}
               [A_t]_{1,1}& \dots         & \dots         & \dots     & [A_t]_{1,n} \\
               [A_t]_{2,1}       & [A_t]_{2,2} & \dots         & \dots     & [A_t]_{2,n}\\
               0 & [A_t]_{3,2}       & \ddots & \ddots     & \vdots \\
               \vdots        & \ddots        & \ddots        & \ddots    & \vdots        \\
               0 & \dots         & 0 & [A_t]_{n,n-1} & [A_t]_{n,n}
             \end{pmatrix}.
\end{equation}

A classical assumption in this degenerate framework (cf.\ \cite{book:Stroock:Varadhan79}, \cite{Krylov04}, \cite{Chaudru:Menozzi17}) is the \emph{uniform ellipticity} of the underlying non-degenerate component of the diffusion matrix at any fixed space-time point. Namely,
\begin{description}
  \item[{[UE]}] There exists a constant $\eta>1$ such that for any $t\ge 0$ and any $x$ in $\R^N$, it holds that
  \[\eta^{-1}\vert \xi \vert^2 \, \le \, \sigma(t,x)\xi\cdot \xi  \, \le\,\eta
\vert \xi \vert^2, \quad \xi \in \R^d,
\]
\end{description}
where ``$\cdot$'' stands for the inner product on the smaller space $\R^d$.\newline
We will suppose that the drift $G(t,x)=A_tx+F(t,x)$ has a particular ``upper diagonal'' structure and its sub-diagonal elements are linear and non-degenerate, i.e.\
\begin{description}
  \item[{[H]}]
  \begin{itemize}
      \item $F=(F_1,\dots,F_n)\colon [0,\infty)\times\R^N\to \R^N$ is such that $F_i$ depends only on time and on the last $n-(i-1)$ components, i.e.\ $F_i(t, x_i,\dots,x_n)$,  for any $i$ in $\llbracket 1, n \rrbracket$;
      \item $A\colon [0,\infty)\to \R^N\otimes \R^N$ \textcolor{black}{is bounded} and the blocks $ [A_t]_{i,j}\in \R^{d_i}\otimes \R^{d_j}$ in \eqref{eq:def_matrix_A} are such that \[[A_t]_{i,j} = \begin{cases}
      \text{is non-singular (i.e.\ it has rank  $d_i$) uniformly in } t, \mbox{ if } j= i-1;\\
      0, \mbox{ if } j< i-1.
      \end{cases}\]
  \end{itemize}
\end{description}
Clearly, $n$ is in $\llbracket 1,N\rrbracket$ and $n=1$ if and
only if $d=N$, i.e.\ if the dynamics is non-degenerate.\newline
In the linear framework ($F=0$) and for constant diffusion coefficients ($\sigma(t,x)=\sigma$), this last assumption can be seen as a H\"ormander-type condition, ensuring the hypoellipticity of the infinitesimal generator associated with the process
$\{X_t\}_{t\ge0}$, which is in this setting equivalent to the Kalman condition, see e.g. \cite{Priola:Zabczyk09}. We highlight however that in our framework, the ``classic'' H\"ormander assumption (cf.\ \cite{Hormander67}) cannot be considered, due to the low
regularity of the coefficients we will consider in \eqref{eq:SDE} (see Theorem \ref{thm:main_result}). This prevents us from explicitly calculating the commutators. \newline
In Equation \eqref{eq:SDE} above, $\{Z_t\}_{t\ge 0}$ is a $d$-dimensional, symmetric and adapted L\'evy
process with respect to some stochastic basis $(\Omega, \mathcal{F},\{ \mathcal{F}_t\}_{t\ge 0},\mathbb{P})$.
We recall that a $d$-valued L\'evy process is a stochastically continuous process on $\R^d$ starting from zero and such that its
increments are independent and stationary. Moreover, it is well-known (see e.g.\ \cite{book:Sato99}) that any L\'evy process admits a
c\`adl\`ag modification, i.e.\ a right continuous modification having left limits $\mathbb{P}$-almost surely. We will always assume to
have chosen such a version.\newline
A fundamental tool in the analysis of L\'evy processes is given by the L\'evy-Kitchine formula (see for instance \cite{book:Jacob05}) that allows us to represent
the L\'evy symbol $\Phi(\xi)$ of $\{Z_t\}_{t\ge0}$, given by
\[\mathbb{E}[e^{i\xi\cdot Z_t}] \, = \, e^{t\Phi(\xi)}, \quad \xi \in \R^d\]
in terms of the generating triplet $(b,\Sigma,\nu)$ as:
\[\Phi(\xi) \, = \, i b \cdot \xi -\frac{1}{2} \Sigma\xi\cdot \xi+\int_{\R^d_0}\bigl(e^{i \xi\cdot z}-1-i \xi\cdot z \mathds{1}_{B(0,1)}(z)\bigr)
\,\nu(dy), \quad \xi \in \R^d,\]
where $b$ is a vector in $\R^d$, $\Sigma$ is a symmetric, non-negative definite matrix in $\R^d\otimes \R^d$ and $\nu$ is a L\'evy measure on
$\R^d_0:=\R^d\smallsetminus \{0\}$, i.e.\ a $\sigma$-finite measure on $\mathcal{B}(\R^d_0)$, the Borel $\sigma$-algebra on $\R^d_0$, such
that $\int(1\wedge \vert z\vert^2) \, \nu(dz)$ is finite.
In particular, any L\'evy process is completely determined by its generating triplet $(b,\Sigma,\nu)$. \newline
Importantly, we point out already that a change on the truncation set $B(0,1)$ for the L\'evy-Kitchine formula does not affect the formulation of the L\'evy symbol $\Phi$, since we assumed $\nu$ to be symmetric. Namely, given a threshold $c>0$, the L\'evy symbol $\Phi(\xi)$ of $\{Z_t\}_{t\ge0}$ could be also represented as
\begin{equation}
\label{eq:change_of_truncation}
\Phi(\xi) \, = \, i b \cdot \xi -\frac{1}{2} \Sigma\xi\cdot \xi+\int_{\R^d_0}\bigl(e^{i \xi\cdot z}-1-i \xi\cdot z \mathds{1}_{B(0,c)}(z)\bigr)
\,\nu(dz), \quad \xi \in \R^d,
\end{equation}
where $b$, $\Sigma$ and $\nu$ are as above.
Here, we only consider pure jump processes, i.e.\ $\Sigma=0$. Indeed, the more general case, where a Gaussian component is considered, can be obtained from already existing results (cf. \cite{Chaudru:Menozzi17}).\newline
We will suppose moreover that, additionally to the symmetry, the L\'evy measure $\nu$ of $\{Z_t\}_{t\ge 0}$ satisfies the following
\emph{non-degeneracy condition}:
\begin{description}
  \item[{[ND]}] there exist a Borel function $Q\colon\R^d\to [0,\infty)$ such that
  \begin{itemize}
      \item $\text{ess-sup}\{Q(z)\colon z \in \R^d\}<+\infty$;
      \item there exist $r_0>0$ and $c>0$ such that $Q(z)\ge c$ and Lipschitz continuous in $B(0,r_0)$;
      \item there exists $\alpha\in (1,2) $ and a finite, non-degenerate measure $\mu$ on $\mathbb{S}^{d-1}$ such that
      \[
\nu(\mathcal{A}) \,=\,\int_{0}^{\infty}\int_{\mathbb{S}^{d-1}}\mathds{1}_{\mathcal{A}}(rs)Q(rs)\,\mu(ds) \frac{dr}{r^{1+\alpha}},\quad \mathcal{A} \in\mathcal{B}(\R^d_0),
\]
\end{itemize}
\end{description}
where $\mathcal{B}(\R^d_0)$ stands for the Borelian
$\sigma$-field on $\R^d_0$.
We recall moreover that a spherical measure $\mu$ on $\mathbb{S}^{d-1}$ is non-degenerate (in the sense of Kolokoltsov \cite{Kolokoltsov00}) if there exists a constant $\tilde{\eta} \ge 1$
such that
\begin{equation}\label{eq:non_deg_measure}
\tilde{\eta}^{-1}\vert \xi \vert^\alpha \, \le \, \int_{\mathbb{S}^{d-1}}\vert \xi\cdot s \vert^\alpha \, \mu(ds) \, \le\,\tilde{\eta}\vert \xi \vert^\alpha,
\quad \xi \in \R^d.
\end{equation}
Since any $\alpha$-stable L\'evy measure can be decomposed into a spherical part $\mu$ on $\mathbb{S}^{d-1}$ and a radial part $r^{-(1+\alpha)}dr$ (see e.g.\ Theorem $14.3$ in \cite{book:Sato99}), assumption [\textbf{ND}] roughly states that the L\'evy measure of $\{Z_t\}_{t\ge0}$ is absolutely continuous with respect to the \textcolor{black}{non-degenerate} (in the sense of \eqref{eq:non_deg_measure}), L\'evy measure of a $\alpha$-stable process and that their Radon-Nikodym derivative is given by the function $Q$. From this point further, we will denote such a L\'evy measure by $\nu_\alpha(dz):=\mu(ds)r^{-(1+\alpha)}dr$ with $z=rs$.\newline
In order to deal with a possibly multiplicative noise, i.e.\ in
the presence of a space-inhomogeneous diffusion coefficient $\sigma$ in Equation \eqref{eq:SDE}, we will need the following:
\begin{description}
\item[{[AC]}] If $x\to \sigma(t,x)$ is non-constant for some
$t\ge0$, then the measure $\nu_\alpha$
is absolutely continuous with respect to the Lebesgue measure on $\R^N$ with Lipschitz Radon-Nykodim derivative.
\end{description}
Assumptions [\textbf{ND}] and [\textbf{AC}] together imply in particular that in the multiplicative case, the L\'evy measure $\nu$ of $\{Z_t\}_{t\ge 0}$ can be decomposed as
\begin{equation}
\label{eq:decomposition_nu_AC}
    \nu(dz) \, = \, Q(z)\frac{g(\frac{z}{|z|})}{|z|^{d+\alpha}}dz,
\end{equation}
for some Lipschitz function $g\colon \mathbb{S}^{d-1}\to \R$.\newline
At this point, we would like to remark that no regularity is assumed for the L\'evy measure $\nu$ of $\{Z_t\}_{t\ge0}$ in the additive framework (or more generally, for a space-independent $\sigma$). In particular, the measure $\mu$ in condition [\textbf{ND}] may not be absolutely continuous with respect to the Lebesgue measure on $\mathbb{S}^{d-1}$. Indeed, our model can also include very singular (with respect to the Lebesgue measure) examples such as the cylindrical $\alpha$-stable process associated with $\mu=\sum_{i=1}^d\frac{1}{2}(\delta_{e_i}+\delta_{-e_i})$. See e.g.\ \cite{Bass:Chen06}  for more details.\newline
From this point further, we always assume that the above hypotheses on the coefficients are satisfied.

We would like to conclude the introduction with some comments concerning our assumptions with particular reference with our previous works.\newline
In \cite{Marino21}, the Schauder estimates, an important analytical first step for proving \textcolor{black}{the} well-posedness of SDEs, has been showed for degenerate Ornstein-Uhlenbeck operators driven by a more general class of L\'evy noises. It also includes, for example, the asymmetric version of the stable-like noises we consider in this work. We start highlighting that a similar family of noises could not have been introduced here, as in \cite{Marino20}, due to the non-linear structure of our problem and, especially, our technique of proof through a perturbative approach. Indeed, it requires more delicate regularizing properties for the involved operators and, in particular, a compatibility between some proxy and the original operator, seen as a perturbation of the first one. \newline
Here, we have followed a backward perturbative approach  as firstly introduced by McKean-Singer in \cite{Mckean:Singer67}. This terminology comes from the fact that the underlying proxy process will be associated with a backward in time flow. This method  appears  more natural  for proving weak uniqueness in a degenerate $L^p-L^q$ framework (cf. \cite{Chaudru:Menozzi17} in the diffusive case). Roughly speaking, it only requires controls on the gradients (in a weak sense) for the solutions of the associated PDE in order to apply the inversion technique on the infinitesimal generator.
However, we are confident that the Schauder estimates presented in \cite{Marino20} could be extended to the class of noises we consider here. Relying on them, we could have then proven the uniqueness in law for dynamics such as \eqref{eq:SDE}. This method appears really involved and long since it structurally requires to establish pointwise estimates for the first order derivatives with respect to the degenerate components of the dynamics.
Another useful advantage of the backward perturbative approach  is that it allows us to show Krylov-type estimates on the solution process $X_t$ of SDE \eqref{eq:SDE}. These kind of controls seems of independent interest and new for random dynamics involving degenerate stable-like noises.

The drawback of our approach is that it leads to a specific structure in Equation \eqref{eq:SDE}, \textcolor{black}{given by assumption} \textbf{[H]}. Namely, we cannot consider  drift of the form $F_i(x)=F_i(x_{i-1},\cdots,x_n) $ with non-linear dependence w.r.t. $x_{i-1}$, variable which transmits the noise. This case is often investigated for Brownian noises (see e.g. \cite{Delarue:Menozzi10}, \cite{Chaudru:Menozzi17}).
This feature is specifically linked to the structure of the joint law of a  stable process and its iterated integrals which generate a multi-scale stable process with highly singular associated spectral measure, see e.g. Proposition \ref{prop:Decomposition_Process_X} and Remark \ref{DA_SCRIVERE_VINCOLO_SUL_MODELLO} below or \cite{Huang:Menozzi16}.
Similar issues constrain us to assume in the multiplicative noise case that the driving process has an absolutely continuous spectral measure with respect to the Lebesgue measure on $\mathbb{S}^{d-1}$. This precisely allows us to get estimates which will be uniform with respect to the parameters for the considered class of proxys.

\paragraph{Main driving processes considered.} Here, we highlight that assumption [\textbf{ND}] applies to a large class of L\'evy processes on $\R^d$. As already
pointed out in \cite{Schilling:Sztonyk:Wang12}, it holds for the following families of stable-like examples with $\alpha \in (0,2)$:
\begin{enumerate}
  \item Stable process \cite{book:Sato99}:
  \[Q(z) \, = \, 1;\]
  \item Truncated stable process with $r_0>0$ \cite{Kim:Song08}:
  \[Q(z) \, = \, \mathds{1}_{(0,r_0]}(|z|);\]
  \item Layered stable process with $\beta>\alpha$ and $r_0>0$ \cite{Houdre:Kawai07}:
  \[Q(z) \, = \, \mathds{1}_{(0,r_0]}(|z|)+\mathds{1}_{(r_0,\infty)}(|z|)|z|^{\alpha-\beta};\]
  \item Tempered stable process \cite{Rosinski09} with $Q(z)=Q(rs)$ such that for all $s$ in $\mathbb{S}^{d-1}$,
  \[r\to Q(rs) \text{ is completely monotone, $Q(0)>0$ and $\lim_{r\to +\infty}Q(rs)=0$}.\]
  \item Relativistic stable process \cite{Carmona:Masters:Simon90}, \cite{Byczkowski:Malecki:Ryznar09}:
  \[Q(z) \, = \, (1+|z|)^{(d+\alpha-1)/2}e^{-|z|};\]
  \item Lamperti process with $f\colon \mathbb{S}^{d-1}\to \R$ even such that
  $\sup f(s)<1+\alpha$ \cite{Caballero:Pardo:Perez10}:
  \[Q(z) \, = \,
  \exp\bigl(|z|f(\frac{z}{|z|})\bigr)\Bigl(\frac{|z|}{e^{|z|}-1}\Bigr)^{1+\alpha}, \quad z \in \R^d_0.\]
\end{enumerate}

\paragraph{Organization of the paper.}
The article is organised as follows. In Section $2$, we introduce some useful notations and we present the associated martingale
problem. In particular, we state there our main results.
Section $3$ contains all the associated analytical tools that allow to derive our results. Namely, we follow there a perturbative approach, considering a suitable linearization of our dynamics \eqref{eq:SDE} around a Cauchy-Peano flow which takes into account the deterministic part of our model (corresponding to \eqref{eq:SDE} with $\sigma =0$). Section $4$ is then dedicated to prove the well-posedness of the associated martingale problem, exploiting the analytical results given in Section $3$. In Section $5$, we finally construct an ``ad hoc'' Peano counter-example to the uniqueness in law for SDE \eqref{eq:SDE}.

\setcounter{equation}{0}
\section{Basic notations and main results}
\fancyhead[RO]{Section \thesection. Basic notations and main results}

We start recalling some useful notations we will need below. In the following, $C$ will denote a generic \emph{positive} constant. It may change from line to line and it will depend only on the parameters appearing in the previously stated assumptions, as for instance: $d,N,\alpha,\eta,b,g,r_0,\mu$. We will explicitly specify any other dependence that may occur.\newline
Given a function $f\colon \R^N\to \R$, we denote by $Df(x)$, and $D^2f(x)$ the first and second Fr\'echet derivative of $f$ at a point $x$ in $\R^N$ respectively,
when they exist. We denote by $B_b(\R^N)$ the family of all the Borel and bounded functions $f\colon \R^N\to \R$. It is a Banach space endowed with the supremum
norm $\Vert \cdot \Vert_\infty$. We also consider its closed subspace $C_b(\R^N)$ consisting of all the continuous functions. Moreover, $C^\infty_c(\R^N)\subseteq
C_b(\R^N)$ denotes the space of smooth functions with compact support.

We now recall two correlated definitions of solution associated with SDE \eqref{eq:SDE}. Let us consider fixed $\mu$ in $\mathcal{P}(\R^N)$, the family of
the probability measures on $\R^N$ and an initial time $t\ge0$.

\begin{definition}
A weak solution of SDE \eqref{eq:SDE} with starting condition $(t,\mu)$ is a $N$-dimensional, c\`adl\`ag, adapted process $\{X_s\}_{s\ge 0}$
on some stocastic base $(\Omega, \mathcal{F},\{ \mathcal{F}_s\}_{s \ge 0}, \mathbb{P})$ such that
\begin{itemize}
  \item the law of $X_t$ is $\mu$;
  \item there exists a $d$-dimensional, adapted L\'evy process $\{Z_s\}_{s\ge t}$ satisfying [\textbf{ND}] and [\textbf{AC}] such that
    \begin{equation}
    \label{eq:SDE_Integral}
    X_s \, = \, X_t + \int_t^sG(u,X_u)\,du+\int_t^s \sigma(u,X_{u-})B\, dZ_u, \quad s\ge t,\, \, \mathbb{P}\text{-a.s.}
    \end{equation}
    \end{itemize}
\end{definition}
To state our second definition, we need to consider the infinitesimal generator  $\partial_s+L_s$ (formally) associated with the solutions of SDE \eqref{eq:SDE}. Noticing that the term involving the constant drift $b$ can be absorbed in the expression for $G$ without loss of generality, the operator $L_s$ can be represented for any $\phi$ in $C^\infty_c(\R^N)$ as
\begin{multline}
   \label{eq:def_generator}
L_s\phi(s,x) \, := \, \langle G(s,x),D_x\phi(x)\rangle +\mathcal{L}_s\phi(s,x)
\\
:= \,\langle G(s,x),D_x\phi(x)\rangle
+ \int_{\R^d_0}\bigl[\phi(x+B(s,x)z)-\phi(x) \bigr] \,\nu(dz),
\end{multline}
where $\langle \cdot, \cdot \rangle$ denotes the inner product on the bigger space $\R^N$ and, for brevity, $B(s,x):=B\sigma(s,x)$. As done in \cite{Priola15}, we introduce the following definition:
\begin{definition}
A solution of the martingale problem for $\partial_s+L_s$ with initial condition $(t,\mu)$ is an $N$-dimensional, c\`adl\`ag process $\{X_s\}_{s\ge t}$ on some
probability space $(\Omega, \mathcal{F},\mathbb{P})$ such that
\begin{itemize}
    \item the law of $X_t$ is $\mu$;
    \item for any $\phi$ in $\dom\bigl(\partial_s+L_s\bigr)$, the process
  \[\Bigl{\{}\phi(s,X_s)-\phi(t,X_t)- \int_{t}^{s}\bigl(\partial_u+L_u\bigr)\phi(u,X_u) \, du\Bigr{\}}_{s \ge t}\]
  is a $\mathbb{P}$-martingale with respect to the natural filtration $\{ \mathcal{F}^X_s\}_{s\ge 0}$ of the process $\{X_s\}_{s\ge 0}$.
\end{itemize}
\end{definition}

We can now recall some known results that enlighten the link between the two definitions presented above.
For a more thorough analysis on the topic of martingale problems in a rather abstract and general framework, we refer to Chapter $4$ in \cite{book:Ethier:Kurtz86}.\newline
Given a solution $\{X_s\}_{s\ge0}$ of SDE \eqref{eq:SDE},
an application of the It\^o formula immediately shows that the process $\{X_s\}_{s\ge0}$ is a solution of the martingale problem for
$\partial_s+L_s$ with initial condition $(t,\mu)$, too. \newline
On the other hand, if there exists a solution $\{X_s\}_{s\ge0}$ of the martingale problem for $\partial_t+L_t$ with initial condition $(t,\mu)$, it is possible to
construct an ``enhanced'' filtered probability space $(\tilde{\Omega},\tilde{ \mathcal{F}}, \{\tilde{ \mathcal{F}}\}_{s\ge0},\tilde{\mathbb{P}})$ on which there exists
a solution $\{\tilde{X}_s\}_{s\ge 0}$ of the SDE \eqref{eq:SDE}. Moreover, the two processes $\{X_s\}_{s\ge t}$ and $\{\tilde{X}_s\}_{s\ge t}$ have the same law (See, for more details, \cite{Kurtz11}). Thus, it holds that:

\begin{prop}
\label{thm:equivalence_existence}
Let $\mu$ be in $\mathcal{P}(\R^N)$ and $t\ge0$. The existence of a weak solution for SDE \eqref{eq:SDE} with initial condition $(t,\mu)$ is equivalent to the existence of a
solution to the martingale problem for $\partial_s+L_s$ with initial condition $(t,\mu)$.
\end{prop}

We can now move on the notion of uniqueness associated with our problem.

\begin{definition}
We say that weak uniqueness holds for the SDE \eqref{eq:SDE} with initial condition $(t,\mu)$ if any two solutions $\{X_s\}_{s\ge 0}$, $\{Y_s\}_{s\ge 0}$ of SDE
\eqref{eq:SDE} with initial condition $(t,\mu)$ have same finite dimensional distributions. In particular, we say that  SDE \eqref{eq:SDE} is weakly well-posed if for any
$\mu$ in $\mathcal{P}(\R^N)$ and any $t\ge 0$, there exists a unique weak solution of SDE \eqref{eq:SDE} with initial condition $(t,\mu)$.
\end{definition}

Since the definition above takes into account only the law of the solutions $\{X_s\}_{s\ge t}$, $\{Y_s\}_{s\ge t}$, they may, in general, have been defined on
different stochastic bases or with respect to two different underlying L\'evy processes. The definition of uniqueness for a solution of the martingale problem for
$\partial_s+L_s$ can be stated similarly.

It is not difficult to check that the uniqueness of the martingale problem for $\partial_s+L_s$ with initial condition $(t,\mu)$ implies the weak uniqueness of the SDE
\eqref{eq:SDE}. Furthermore, it has been shown in \cite{Kurtz11}, Corollary $2.5$ that the converse is also true.

\begin{prop}
\label{thm:equivalence_uniqueness}
Let $\mu$ be in $\mathcal{P}(\R^N)$ and $t\ge 0$. Then, weak uniqueness holds for SDE \eqref{eq:SDE} with initial condition $(t,\mu)$ if and only if uniqueness holds for the
martingale problem for $\partial_s+L_s$ with initial condition $(t,\mu)$.
\end{prop}

Thanks to Propositions \ref{thm:equivalence_existence} and \ref{thm:equivalence_uniqueness}, we can conclude that the two approaches, i.e.\ the martingale formulation
and the dynamics given in \eqref{eq:SDE}, are equivalent in specifying a L\'evy diffusion process on $\R^N$. We recall however that a third, yet equivalent, method
is given by the forward Fokker-Plank equation governing the law of the process. We will not  explicitly define it since we will not exploit it afterwards (see,
for more details, e.g.\ \cite{Figalli08}, \cite{LeBris:Lions08}). \newline
From now on, we write $x$ in $\R^N$ as $x=(x_1,\dots,x_n)$ where $x_i=(x^1_i,\dots,x^{d_i}_i)$ is in $\R^{d_i}$ for any $i$ in $\llbracket 1,n\rrbracket$.\newline
We can now state our main theorem.

\begin{theorem}
\label{thm:main_result}
For any $ j$ in $\llbracket 1,n\rrbracket$, let $\beta^j$ be an index in $(0,1]$ such that
\begin{itemize}
    \item $x_j \to \sigma(t,x_1,\dots,x_j,\dots,x_n)$ is $\beta^1$-H\"older continuous, uniformly in $t$ and in $x_i$ for $i\neq j$;
    \item $x_j \to F_1(t,x_1,\dots,x_j,\dots,x_n)$ is $\beta^1$-H\"older continuous, uniformly in $t$ and in $x_i$ for $i\neq j$;

    \item $x_j \to F_i(t,x_i,\dots,x_j,\dots,x_n)$ is $\beta^j$-H\"older continuous, uniformly in $t$ and in $x_k$, for $k\neq j$ and $2\le i\le j$.
\end{itemize}
Additionally, we suppose that there exists $K\ge 1$ such that $|F_i(t,0)|\le K$ for any $i$ in $\llbracket 1,n\rrbracket$ and any $t\ge 0$. Then, the SDE \eqref{eq:SDE} is weakly well-posed if
\begin{equation}
\label{eq:thresholds_beta}
\beta^j\, > \, \frac{1+\alpha(j-2)}{1+\alpha(j-1)},\ j\ge 2.
\end{equation}
\end{theorem}
Theorem \ref{thm:main_result} above will follow from Propositions \ref{thm:equivalence_existence} and \ref{thm:equivalence_uniqueness}, once we have shown
that under the same assumptions, there exists a unique weak solution to the martingale problem for $(\partial_s+L_s,\delta_x)$ at any $x$ in $\R^N$.\newline
As a by-product of our method of proof, we have been able to show a Krylov-type estimates for the solutions of SDE \eqref{eq:SDE}. For notational convenience, we will say that two real numbers $p>1$, $q>1$ satisfy Condition $(\mathscr{C})$ when the following inequality holds:
\[ \tag{$\mathscr{C}$}
\bigl(\frac{1-\alpha}{\alpha}N+\sum_{i=1}^n
id_i\bigr)\frac{1}{q}+\frac{1}{p} \, < \, 1.
\]

Roughly speaking, such a threshold guarantees the necessary integrability in time with respect to the associated intrinsic scale of the system when considering the $L^p_t-L^q_x$ theory (see Equation \eqref{eq:control_det_T2} for more details). Furthermore, when considering the homogeneous case, i.e.\ when all the components of the system has the same dimension ($d_i=d$ and $N=nd$), condition $(\mathscr{C})$ can be rewritten in the following, clearer, way:
\[\left(\frac{2+\alpha(n-1)}{\alpha}\right)\frac{nd}{q}+\frac{2}{p} \, < \, 2.\]
In particular, taking $\alpha=2$ above, we find the same threshold appearing in \cite{Chaudru:Menozzi17} for the diffusive setting. We highlight moreover that our thresholds can be seen as a natural extension of the ones appearing in \cite{Krylov:Rockner05} in the non-degenerate, Brownian setting.

\begin{corollary}\label{coroll:Krylov_Estimates}
Under the same assumptions of Theorem \ref{thm:main_result}, let $T>0$ and $p>1$, $q>1$ such that Condition $(\mathscr{C})$ holds. Then, there exists a constant $C:=C(T,p,q)$ such that for any  $f$ in $L^p\bigl(0,T;L^q(\R^N)\bigr)$, it holds
\begin{equation}
\label{eq:Krylov_Estimates}
\Bigl{|}\mathbb{E}^{\mathbb P_{t,x}}\bigl[\int_t^Tf(s,X_s) \, ds\bigr] \Bigr{|} \, \le \, C\Vert f \Vert_{L^p_tL^q_x}, \quad (t,x) \in [0,T]\times \R^N,
\end{equation}
where $\{X_s\}_{s\ge0}$ is the canonical process associated with $\mathbb P_{t,x}[\cdot]:=\mathbb E[\cdot|X_t=x] $ which is also the unique weak solution of SDE \eqref{eq:SDE} with initial condition $(t,x)$. In particular, the random variable $X_s$ admits a density $p(t,s,x,\cdot)$ for any $t<s$ and any $x$ in $\R^N$.
\end{corollary}

Additionally, we have been able to show the following non uniqueness result.

\begin{theorem}
\label{thm:counterexample}
Let us consider SDE \eqref{eq:SDE} with $\sigma=1$  and assume that
\begin{itemize}
    \item $x_j \to F_i(t,x_i,\dots,x_j,\dots,x_n)$ is $\beta_i^j$-H\"older continuous, uniformly in $t$ and in $x_k$, for $k\neq j$.
\end{itemize}
    Then, for given $i$, $j$ in $\llbracket 2,n\rrbracket$ with $j\ge i$ there exists $F_i(t,x_i,\dots,x_j,\dots,x_n)=F_i(t,x_j) $ with
\[
\beta^j_i \, < \, \frac{1+\alpha(i-2)}{1+\alpha(j-1)},
\]
for which  weak uniqueness fails for the SDE \eqref{eq:SDE}.
\end{theorem}

The above result will be proven in Section $5$, showing a suitable,
explicit Peano-type counter-example.

\begin{remark}
As opposed to the Gaussian driven case, we did not succeed to obtain regularity indexes which are \textit{sharp} at any level of the chain (cf.\ \cite{Chaudru:Menozzi17}). However, we point out that for diagonal systems of the form:
\begin{equation}
\label{eq:DIAG}
  \begin{cases}
    dX^1_t \, = \, F_1(t,X^1_t,\dots,X^n_t)dt + \sigma(t,X^1_{t-},\dots,X^n_{t-})dZ_t,\\
    dX^2_t \, = \, \left[A^2_tX^1_t +F_2(t,X^2_t)\right] dt,\\
    dX^3_t \, = \, \left[A^3_tX^2_t+F_3(t,X^3_t)\right]dt,\\
    \vdots\\
    dX^n_t \, = \, \left[A^n_t X^{n-1}_t+F_n(t,X^n_t)\right]dt,
  \end{cases}
\end{equation}
i.e. the degenerate components are perturbed by a function which only depends of the current level on the chain, we  have that the previous thresholds are \textit{almost} sharp. Indeed, in this case, we are led to consider $\beta^j >\frac{1+\alpha(j-2)}{1+\alpha(j-1)}$ which gives the well-posedness from the conditions in Theorem \ref{thm:main_result}  while Theorem \ref{thm:counterexample} shows that uniqueness fails as soon as $ \beta_j^j<\frac{1+\alpha(j-2)}{1+\alpha(j-1)} $.
For this diagonal system, Theorems \ref{thm:main_result} and
\ref{thm:counterexample} together then provide an ``almost'' complete
understanding of the weak well-posedness for degenerate SDEs of type
\eqref{eq:DIAG} with H\"older coefficients. Indeed,
the problem for the critical exponents
\[\overline{\beta}^j_j \, = \, \frac{1+\alpha(j-2)}{1+\alpha(j-1)}, \quad
 j \in \llbracket 1,n\rrbracket, \]
remains to be investigated and, up to our best knowledge, there are no general
available results even in the diffusive case. We can only mention \cite{Zhang18} in the kinetic case.
\end{remark}

\setcounter{equation}{0}
We present in this section the analytical tools we will need to show the well-posedness of the associated martingale problem. In particular, they will be fundamental in the derivation of our main Theorem \ref{thm:main_result}, thanks to Propositions \ref{thm:equivalence_existence} and \ref{thm:equivalence_uniqueness}. For this reason, we will assume in this section to be under the same conditions of Theorem \ref{thm:main_result}. Moreover, we will suppose that the final time horizon $T$ is small enough for our scopes. Indeed, we could always exploit the Markov property of the involved processes and standard chaining in time arguments to extend the results to arbitrary (but finite) time intervals.

\subsection{The ``frozen'' dynamics}

The crucial element in our approach consists in choosing wisely a suitable proxy operator with well-known properties and controls, along which we can expand the infinitesimal generator $L_s$, with an additional negligible error.\newline
In order to deal with potentially unbounded perturbations $F$, it is natural to use a proxy involving a non-zero first order term  associated with a flow associated with $G(t,x):=Ax+F(t,x)$, the transport part of SDE \eqref{eq:SDE} (see e.g. \cite{Krylov:Priola10} or \cite{Chaudru:Menozzi:Priola19}).\newline
Remembering that we assume $F$ to be H\"older continuous, we know from the classical Peano-Lipschitz Theorem that there exists a solution of
\begin{equation}
\label{eq:def_Cauchy_Peano_flow}
    \begin{cases}
  d \theta_{t,\tau}(\xi) \, = \, \bigl[A_t \theta_{t,\tau}( \xi)+F(t, \theta_{t,\tau}( \xi))\bigr]\, dt \quad \mbox{ on } [0,\tau];\\
   \theta_{\tau,\tau}( \xi) \, = \,  \xi,
\end{cases}
\end{equation}
even if it may be not unique. For this reason, we are going to choose one particular flow, denoted by $\theta_{t,\tau}( \xi)$, and consider it fixed throughout the work. As it will be shown below in Lemma \ref{lemma:measurability_flow}, it is always possible to take a measurable version of such a flow.\newline
More precisely, given a freezing couple $(\tau, \xi)$ in $(0,T]\times\R^N$, the backward flow will be defined on $[0,\tau]$ as
\[\theta_{t,\tau}( \xi) \, = \,  \xi - \int_{t}^{\tau}\bigl[ A_u \theta_{u,\tau}( \xi)+F(u, \theta_{u,\tau}( \xi))\bigr] \, du.\]
Fixed the reference flow, the next step is to consider the stochastic dynamics linearized along the backward flow $\theta_{t,\tau}(\xi)$. Namely, for any fixed starting point $(t,x)$ in $[0,\tau]\times \R^N$, we consider  $\{\tilde{X}^{\tau,\xi,t,x}_s\}_{s\in [t,T]}$ solving the following SDE:
\begin{equation}
\label{eq:SDE_frozen_SDE}
\begin{cases}
d\tilde{X}^{\tau,\xi,t,x}_u\, = \, \bigl[A_u\tilde{X}^{\tau,\xi,t,x}_u+\tilde{F}^{\tau, \xi}_u\bigr]\, du +B\tilde{\sigma}^{\tau, \xi}_u\,  dZ_u, \quad u\in  [t,T],\\
 \tilde{X}^{\tau,\xi,t,x}_t\, = \, x,
\end{cases}
\end{equation}
where $\tilde{\sigma}^{\tau,\xi}_s:=\sigma(s,\theta_{s,\tau}(\xi))$ and $\tilde{F}^{\tau, \xi}_s:=F(s,\theta_{s,\tau}(\xi))$.\newline
In order to obtain an integral representation of the process $\{\tilde{X}^{\tau,\xi,t,x}_s\}_{s\in [t,T]}$, we now introduce the time-ordered resolvent $\mathcal{R}_{s,t}$ of the matrix $A_s$ starting at time $t$. Namely, $\mathcal{R}_{s,t}$ is a time-dependent matrix in $\R^N\otimes \R^N$ that is solution of the following ODE:
\[\begin{cases}
  \partial_s \mathcal{R}_{s,t} \, = \, A_s \mathcal{R}_{s,t} ,\quad s \in [0,T]; \\
  \mathcal{R}_{t,t} \, = \, \Id_{N\times N}.
\end{cases}\]
By the variation of constants method, it is now easy to check that the solution $\tilde{X}^{\tau,\xi,t,x}_s$ of SDE \eqref{eq:SDE_frozen_SDE} satisfies that
\begin{equation}\label{eq:integral_represent_frozen_SDE}
\tilde{X}^{\tau,\xi,t,x}_s \, = \, \tilde{m}^{\tau, \xi}_{s,t}(x) +\int_t^s \mathcal{R}_{s,u}B\tilde{\sigma}^{\tau,\xi}_u dZ_u,
\end{equation}
where the ``frozen shift'' $\tilde{m}^{\tau, \xi}_{s,t}(x)$ is given by:
\begin{equation}\label{eq:def_tilde_m}
\tilde{{m}}^{\tau, \xi}_{s,t}( x) \, = \, \mathcal{R}_{s,t} x + \int_{t}^{s}\mathcal{R}_{s,u} \tilde{F}^{\tau,\xi}_u \, du.
\end{equation}

We point out already two important properties of the shift $\tilde{m}^{\tau,\xi}_{s,t}(x)$.

\begin{lemma}
\label{lemma:identification_theta_m}
Let $s$ in $[0,T]$ and $x,y$ two points in $\R^N$. Then, for any $t<s$, it holds that
\begin{align}
\label{eq:identification_theta_m}
\tilde{m}^{t,x}_{s,t}(x) \, &= \, \theta_{s,t}(x)\\
y-\tilde{m}^{s,y}_{s,t}(x) \, &= \, \theta_{t,s}(y)-x\label{eq:identification_theta_m1}
\end{align}
\end{lemma}
\begin{proof}
We start noticing that by construction in \eqref{eq:def_tilde_m}, $\tilde{m}^{\tau,\xi}_{s,t}(x)$ satisfies
\begin{equation}
\label{eq:dif_differential_m_tilde}
\partial_s\tilde{m}^{\tau,\xi}_{s,t}(x) \, = \,  A_s\tilde{m}^{\tau,\xi}_{s,t}(x)+F(s, \theta_{s,\tau}(\xi)),
\end{equation}
for any freezing parameters $(\tau,\xi)$.
Choosing $\tau=t$, $\xi=x$ above, it then holds that
\[\partial_s\left[\tilde{m}^{s,x}_{s,t}(x) - \theta_{s,t}(x)\right] \, = \, A_s\left[
\tilde{m}^{t,x}_{s,t}(x)-\theta_{s,t}(x) \right].\]
Since, $\tilde{m}^{t,x}_{t,t}(x)=\theta_{t,t}(x)=x$, Equation \eqref{eq:identification_theta_m} then follows immediately applying the Gr\"onwall lemma.\newline
The second identity in \eqref{eq:identification_theta_m1} follows in a similar manner.
\end{proof}

We are now interested in investigating the analytical properties of the ``frozen'' solution process $\tilde{X}^{\tau,\xi,t,x}_s$. In particular, we will show in the next results the existence of a density for such a process and its anisotropic regularizing effect, at least for small times. Further on, we will consider fixed a time-dependent matrix $\mathbb{M}_t$ on $\R^N\otimes \R^N$ given by
\begin{equation}
\label{eq:def_matrix_M}
    \mathbb{M}_t := \text{diag}(I_{d_1\times d_1},tI_{d_2\times d_2},\dots,t^{n-1}I_{d_{n}\times d_{n}}), \quad t\ge0,
\end{equation}
which reflects the multi-scale nature of the underlying dynamics in \eqref{eq:SDE_frozen_SDE}.

\begin{prop}[Decomposition]
\label{prop:Decomposition_Process_X}
Let the freezing couple $(\tau,\xi)$ be in $[0,T]\times \R^N$, $t<s$ in $[0,T]$ and $x$ in $\R^N$. Then, there exists a L\'evy process $\{\tilde{S}^{\tau,\xi,t,s}_{u}\}_{u\ge 0}$ such that
\begin{equation}\label{eq:decomposition_measure_1}
\tilde{X}^{\tau,\xi,t,x}_s \, \overset{(\text{law})}{=} \, \tilde{m}^{\tau, \xi}_{s,t}(x) + \mathbb{M}_{s-t} \tilde{S}^{\tau,\xi,t,s}_{s-t}.
\end{equation}
In particular, the random variable $\tilde{X}^{\tau,\xi,t,x}_s$ admits a continuous density $\tilde{p}^{\tau,\xi}(t,s,x,\cdot)$ given by
\begin{align}
\label{eq:representation_density}
&\tilde{p}^{\tau,\xi}(t,s,x,y) \, = \,  \frac{1}{\det \mathbb{M}_{s-t}}p_{\tilde{S}^{\tau,\xi,t,s}}\left(t-s, \mathbb{M}^{-1}_{s-t}(y-\tilde{m}^{\tau,\xi}_{s,t}(x))\right) \\
&:= \, \frac{\det \mathbb{M}^{-1}_{s-t}}{(2\pi)^N}\int_{\R^N}e^{-i\langle \mathbb{M}^{-1}_{s-t}(y-\tilde{m}^{\tau,\xi}_{s,t}(x)),z\rangle}\exp\left((s-t)\int_{\R^N}\left[\cos(\langle z, p \rangle )-1\right]\nu_{\tilde{S}^{\tau,\xi,t,s}}(dp)\right)\, dz, \notag
\end{align}
where $\nu_{\tilde{S}^{\tau,\xi,t,s}}$ and $p_{\tilde{S}^{\tau,\xi,t,s}}(u,\cdot)$ are  the L\'evy measure and the density associated with the process $\{\tilde{S}^{\tau,\xi,t,s}_{u}\}_{u\ge0}$, respectively.
\end{prop}
\begin{proof}
For simplicity, we start denoting
\[
\tilde{\Lambda}^{\tau,\xi,t,s} \,:= \, \int_t^s \mathcal{R}_{s,u}B\tilde{\sigma}^{\tau,\xi}_u dZ_u,\quad s\ge t,
\]
so that we have from Equation \eqref{eq:integral_represent_frozen_SDE} that $\tilde{X}^{\tau,\xi,t,x}_s \, = \, \tilde{m}^{\tau, \xi}_{s,t}(x) +\tilde{\Lambda}^{\tau,\xi,t,s}$. To conclude, we need to construct a L\'evy process $\{\tilde{S}^{\tau,\xi,t,s}_u\}_{u\ge 0}$ on $\R^N$  such that
\begin{equation}\label{eq:identity_in_law}
\tilde{\Lambda}^{\tau,\xi,t,s} \, \overset{(\text{law})}{=} \, \mathbb{M}_{s-t}\tilde{S}^{\tau,\xi,t,s}_{s-t}.
\end{equation}
To show the identity in law, we are going to reason in terms of the characteristic functions. We start recalling that
the L\'evy process $\{Z_t\}_{t\ge 0}$ on $\R^d$ is characterized by the L\'evy symbol
\[\Phi(p) \, = \,  \int_{\R^d_0}\left[\cos(p\cdot q)-1\right]Q(q) \, \nu_\alpha(dq),
\quad p \in \R^d,\]
where $\nu_\alpha(dq)=\mu(d\theta)\frac{dr}{r^{1+\alpha}}$ is the L\'evy measure of an $\alpha$-stable process.
It is well-known (see e.g.\ Lemma $2.2$ in \cite{Schilling:Wang12}) that at any fixed $t\le s$ in $[0,1]$, $\tilde{\Lambda}^{\tau,\xi,t,s}$ is an infinitely divisible random variable with associated L\'evy symbol
\[\Phi_{\tilde{\Lambda}^{\tau,\xi,t,s}}(z) \, := \, \int_{t}^{s}\Phi\bigl((\mathcal{R}_{s,u}B\tilde{\sigma}^{\tau,\xi}_u)^*z\bigr)\,du, \quad z \in \R^N,\]
where, we recall, we have denoted $\tilde{\sigma}^{\tau,\xi}_u=\sigma(u,\theta_{u,\tau}(\xi))$.\newline
Setting $v:=(u-t)/(s-t)$ and noticing that $u=u(v):=t+v(s-t)$, we can now
rewrite the L\'evy symbol of $\tilde{\Lambda}^{\tau,\xi,t,s}$ as
\begin{equation}
\label{proof:ref_notations}
  \Phi_{\tilde{\Lambda}^{\tau,\xi,t,s}}(z) \, := \, (s-t)\int_{0}^{1}\Phi\bigl((\mathcal{R}_{s,u(v)}B\tilde{\sigma}^{\tau,\xi}_{u(v)})^*z\bigr)\,dv.
\end{equation}
From the analysis performed in \cite{Huang:Menozzi16}, Lemmas $5.1$ and $5.2$ (see also \cite{Delarue:Menozzi10} Proposition $3.7$), we then know that we can decompose the first column of the resolvent $\mathcal{R}_{s,u(v)}$ in the following way:
\[\mathcal{R}_{s,u(v)}B \, = \, \mathbb{M}_{s-t}\widehat{\mathcal{R}}_vB,\]
where $\{\widehat{\mathcal{R}}_v\colon v \in [0,T]\}$ are non-degenerate and bounded matrixes in $\R^N\otimes \R^N$ and the multi-scale matrix $\mathbb{M}_t$ is given in \eqref{eq:def_matrix_M}. We can now rewrite the L\'evy symbol of $\tilde{\Lambda}^{\tau,\xi,t,s}$ as
\[\Phi_{\tilde{\Lambda}^{\tau,\xi,t,s}}(z) \, = \, (s-t)\int_{0}^{1}\Phi\bigl((\widehat{\mathcal{R}}_vB\tilde{\sigma}^{\tau,\xi}_{u(v)})^*\mathbb{M}_{s-t}z\bigr)\,dv, \quad z \in \R^N.\]
The above equality suggests us to define, for any fixed $t\le s$ in $(0,1]$, the (unique in law) L\'evy process  $\{\tilde{S}^{\tau,\xi,t,s}_u\}_{u\ge0}$ associated with the L\'evy symbol
\begin{equation}
\label{proof:eq:def_Levy_symbol_S}
\begin{split}
\Phi_{\tilde{S}^{\tau,\xi,t,s}}(z) \, &:= \, \int_{0}^{1}\Phi\bigl((\widehat{\mathcal{R}}_{v}B\tilde{\sigma}^{\tau,\xi}_{u(v)})^*z\bigr)\,dv \\
&= \, \int_{0}^{1}\int_{\R^d}\left[\cos\left(\langle z,\widehat{\mathcal{R}}_{v}B\tilde{\sigma}^{\tau,\xi}_{u(v)}p\rangle\right)-1\right]\,\nu(dp)dv.
\end{split}
\end{equation}
Since we have that
\begin{equation}
\label{proof:eq:inversion}
\mathbb{E}\bigl[e^{i\langle z,\tilde{\Lambda}^{\tau,\xi,t,s}\rangle}\bigr] \, = \, e^{\Phi_{\tilde{\Lambda}^{\tau,\xi,t,s}}(z)}\, = \, e^{(s-t)\Phi_{\tilde{S}^{\tau,\xi,t,s}}(\mathbb{M}_tz)} \,
= \,\mathbb{E}\bigl[e^{i\langle z,\mathbb{M}_t\tilde{S}^{\tau,\xi,t,s}_{s-t}\rangle}\bigr],
\end{equation}
it follows immediately that Equation \eqref{eq:identity_in_law} holds.\newline
To show the existence of a density for $\tilde{X}^{\tau,\xi,t,x}_s$, we want to exploit the Fourier inversion formula in \eqref{proof:eq:inversion}. To do it, we firstly need to prove that $\text{exp}(\Phi_{\tilde{\Lambda}^{\tau,\xi,t,s}}(z))$ is integrable. From \eqref{proof:eq:def_Levy_symbol_S}, we notice that
\[
\begin{split}
\Phi_{\tilde{S}^{\tau,\xi,t,s}}(z) \,
&= \, \int_{0}^{1}\int_{\R^d}\left[\cos\left(\langle z,\widehat{\mathcal{R}}_{v}B\tilde{\sigma}^{\tau,\xi}_{u(v)}p\rangle\right)-1\right]\,\nu(dp)dv \\
&= \,\int_{0}^{1}\int_{\R^d}\left[\cos\left(\langle z,\widehat{\mathcal{R}}_{v}B\tilde{\sigma}^{\tau,\xi}_{u(v)}p\rangle\right)-1\right]Q(p)\,\nu_\alpha(dp)dv,
\end{split}
\]
where in the last step we used hypothesis [\textbf{ND}]. Exploiting now that the quantities above are non-positive and $Q(p)\ge c>0$ for $p$ in $B(0,r_0)$, we write that
\[
\begin{split}
&\Phi_{\tilde{S}^{\tau,\xi,t,s}}(z)\, \le C \,\int_{0}^{1}\int_{B(0,r_0)}\left[\cos\left(\langle z,\widehat{\mathcal{R}}_{v}B\tilde{\sigma}^{\tau,\xi}_{u(v)}p\rangle\right)-1\right]\,\nu_\alpha(dp)dv \\
&\qquad= \, C\Bigl{\{}-\int_0^1\left|(\widehat{\mathcal{R}}_{v}B\tilde{\sigma}^{\tau,\xi}_{u(v)})^\ast z\right|^\alpha dv+\int_0^1\int_{B^c(0,r_0)}\left[1-\cos\left(\langle z,\widehat{\mathcal{R}}_{v}B\tilde{\sigma}^{\tau,\xi}_{u(v)}p\rangle\right)\right]\nu_\alpha(dp)dv\Bigr{\}} \\
&\qquad\le \, C\Bigl{\{}-\int_0^1\left|(\widehat{\mathcal{R}}_{v}B\tilde{\sigma}^{\tau,\xi}_{u(v)})^\ast z\right|^\alpha\, dv+1\Bigr{\}}.
\end{split}
\]
To conclude, we recall that Lemma $5.4$ in \cite{Huang:Menozzi16}  states that
\[\int_0^1\left|(\widehat{\mathcal{R}}_{v}B\tilde{\sigma}^{\tau,\xi}_{u(v)})^\ast z\right|^\alpha\, dv \, \ge \, C|z|^\alpha,\]
for some positive constant $C$ independent from $t,s,\tau,\xi$. It then follows in particular that
\begin{equation}
\label{eq:control_Levy_symbol_S}
\Phi_{\tilde{S}^{\tau,\xi,t,s}}(z) \, \le \,C\left[1-|z|^\alpha\right], \quad z \in \R^N.
\end{equation}
Since $\text{exp}(\Phi_{\tilde{\Lambda}^{\tau,\xi,t,s}}(z))$ is integrable, it implies that there exists a density $p_{\tilde{\Lambda}^{\tau,\xi,t,s}}(t,s,\cdot)$ of the random variable $\tilde{\Lambda}^{\tau,\xi,t,s}$.  We can now apply the Fourier inversion formula in Equation \eqref{proof:eq:inversion} showing that $p_{\tilde{\Lambda}^{\tau,\xi,t,s}}(t,s,\cdot)$ is given by
\begin{equation}
 \label{eq:definition_density_q}
p_{\tilde{\Lambda}}^{\tau,\xi}(t,s,y) \, := \, \frac{1}{(2\pi)^N}\int_{\R^N}e^{-i\langle y,z\rangle}\text{exp}\left((s-t)\Phi_{\tilde{S}^{\tau,\xi,t,s}}(z)\right)\, dz.
\end{equation}
From Decomposition \eqref{eq:decomposition_measure_1} and Equation \eqref{eq:definition_density_q}, the representation for $\tilde{p}^{\tau,\xi}(t,s,x,\cdot)$ follows immediately.
\end{proof}

Once we have proven the existence of a density $\tilde{p}^{\tau,\xi}(t,s,x,\cdot)$ for the ``frozen'' stochastic dynamics $\tilde{X}^{\tau,\xi,t,x}_s$, we move now on determining its associated smoothing effects. In particular, we show in the following proposition that the derivatives of the ``frozen'' density are controlled by another density at the price of an additional time singularity of the order corresponding to the intrinsic time scale of the considered component in the stable regime. Importantly, such a control holds uniformly in the freezing parameters $(\tau, \xi)$. \newline
Let us introduce for simplicity the following time-dependent scale matrix:
\begin{equation}
\label{eq:def_matrix_T}
    \mathbb{T}_{t} \, := \, t^{\frac{1}{\alpha}}\mathbb{M}_t,\quad t\ge 0.
\end{equation}

\begin{prop}
\label{prop:Smoothing_effect}
There exists a family $\{\overline{p}(u,\cdot)\colon u\ge 0\}$ of densities on $\R^N$ and a positive constant $C:=C(N,\alpha)$ such that
\begin{itemize}
  \item for any $u\ge0$ and any $z$ in $\R^N$,  $\overline{p}(u,z)=u^{-N/\alpha}\overline{p}(1,u^{-1/\alpha}z)$; (stable scaling property)
  \item for any $\gamma$ in $[0,\alpha)$,
  \begin{equation}\label{equation:integration_prop_of_q}
    \int_{\R^N} \overline{p}(u,z) \, \vert z \vert^\gamma \, dz \, \le \, Cu^{\gamma/\alpha}, \quad u> 0;
  \end{equation}
  \item  for any $k$ in $\llbracket 0,2 \rrbracket$, any $i$ in $\llbracket 1,n \rrbracket$, any $t< s$ in $[0,T]$ and any $x,y$ in $\R^N$,
  \begin{equation}
  \label{eq:derivative_prop_of_q}
      \vert D^k_{x_i}\tilde{p}^{\tau,\xi}(t,s,x,y) \vert \le C\frac{(s-t)^{-k\frac{1+\alpha(i-1)}{\alpha}}}{\det \mathbb{T}_{s-t}}\overline{p}\left(1,\mathbb{T}^{-1}_{s-t}(y-\tilde{m}^{\tau, \xi}_{s,t}(x))\right).
  \end{equation}
  where we denoted, coherently with the notations introduced before Theorem \ref{thm:main_result}, $D_{x_i}=\left(D_{x^1_i},\dots, D_{x^{d_i}_i}\right)$.
\end{itemize}
\end{prop}
\begin{remark}[About the freezing parameters]
\label{DA_SCRIVERE_VINCOLO_SUL_MODELLO}
We carefully point out that since we will later on choose as parameters $(\tau,\xi)=(s,y)$, it is particularly important that we manage to obtain an upper bound by a density which is independent from those parameters, since they will be as well the integration variables (see Section \ref{SEZIONE_SPIEGAZIONE_CONGELAMENTO_E_ALTRI} below). This is precisely why we actually impose the specific semi-linear drift structure in SDE \eqref{eq:SDE} (cf.\ assumption \textbf{[H]}), as opposed to the more general one that can be handled  in the Gaussian case \cite{Chaudru:Menozzi17}. This is a framework which naturally gives the independence of the large jumps of the \textit{proxy} process $\tilde{X}^{\tau,\xi,t,x}_s$ as used in \eqref{proof_control_N_segnato} below. The more general case for the first  order dynamics considered in \cite{Chaudru:Menozzi17} would actually lead to linearize around a matrix which would depend on the freezing parameters. For such models, we did not succeed in proving that the corresponding densities can be bounded independently of the parameters (see also the proof of Lemma \ref{lemma:Control_P_N} below for a similar issue regarding the  diffusion coefficient).
\end{remark}

\begin{proof}
Fixed the freezing parameters $(\tau,\xi)$ in $[0,T]\times \R^N$, and the times $t<s$ in $[0,T]$, we start applying the It\^o-L\'evy decomposition to the process $\{\tilde{S}^{\tau,\xi,t,s}_u\}_{u\ge 0}$ introduced in Proposition \ref{prop:Decomposition_Process_X} at the associated characteristic stable time, i.e.\ we choose to truncate at threshold $u^{1/\alpha}$. Thus, we can write
\begin{equation}
\label{eq:decomposition_S}
\tilde{S}^{\tau,\xi,t,s}_u \,= \, \tilde{M}^{\tau,\xi,t,s}_u+\tilde{N}^{\tau,\xi,t,s}_u
\end{equation}
for some $\tilde{M}^{\tau,\xi,t,s}_u,\tilde{N}^{\tau,\xi,t,s}_u$ independent random
variables corresponding to the small jumps part and the large jumps part, respectively. Namely, we denote for any $v>0$,
\begin{equation}\label{DEF_TILDE_N}
\tilde{N}^{\tau,\xi,t,s}_v \, := \, \int_{0}^{v}\int_{\vert z \vert >u^{1/\alpha}}zP_{\tilde{S}^{\tau,\xi,t,s}}(dr,dz) \,\, \text{ and } \,\, \tilde{M}^{\tau,\xi,t,s}_v: \,= \, \tilde{S}^{\tau,\xi,t,s}_v-\tilde{N}^{\tau,\xi,t,s}_v,
\end{equation}
where $P_{\tilde{S}^{\tau,\xi,t,s}}(dr,dz)$ is the Poisson random measure associated with the process $\tilde{S}^{\tau,\xi,t,s}$. It can be shown, similarly to Proposition \ref{prop:Decomposition_Process_X}, that the process $\{\tilde{M}^{\tau,\xi,t,s}_u\}_{u\ge 0}$ admits a density $p_{\tilde{M}^{\tau,\xi,t,s}}(u,\cdot)$. Indeed, it is well-known that the small jump part leads to a density which is in the Schwartz class $\mathcal{S}(\R^N)$ (see Lemma \ref{lemma:Control_p_M} below). We can then rewrite the density $p_{\tilde{S}^{\tau,\xi,t,s}}$ of $\tilde{S}^{\tau,\xi,t,s}$ in the following way:
\begin{equation}
\label{proof:decomposition_density_S}
    p_{\tilde{S}^{\tau,\xi,t,s}}(u,z) \, = \, \int_{\R^N}p_{\tilde{M}^{\tau,\xi,t,s}}(u,z-y)P_{\tilde{N}^{\tau,\xi,t,s}_u}(dy)
\end{equation}
where $P_{\tilde{N}^{\tau,\xi,t,s}_u}$ is the law of $\tilde{N}^{\tau,\xi,t,s}_u$.\newline
We need now to control the modulus of the density $p_{\tilde{S}^{\tau,\xi,t,s}}$ with another density, independently from the parameters $\tau$, $\xi$.
From Lemma \ref{lemma:Control_p_M} in the Appendix (see also Lemma B.$2$ in \cite{Huang:Menozzi16})  with $m =N+1$, we know that there exists a positive constant $C$, independent from $\tau,\xi$ such that
\begin{equation}
\label{proof:control_pM_segnato}
\left| D^k_{z} p_{\tilde{M}^{\tau,\xi,t,s}}(u,z) \right| \, \le \, Cu^{-(N+k)/\alpha}\left(\frac{u^{1/\alpha}}{u^{1/\alpha}+|z|}\right)^{N+2} \, =: \, Cu^{-k/\alpha}p_{\overline{M}}(u,z),
\end{equation}
for any $k$ in $\llbracket 0, 2 \rrbracket$, any $u>0$, and any $z$ in $\R^N$.\newline
Moreover, denoting by $\overline{M}_u$ the random variable with density $p_{\overline{M}}(u,\cdot)$ that is independent from $\tilde{N}^{\tau,\xi,t,s}_u$, we can easily check that $p_{\overline{M}}(u,z)=u^{-N/\alpha}p_{\overline{M}}(1,u^{-1/\alpha}z)$ and thus, that $\overline{M}$ is $\alpha$-selfsimilar:
\[\overline{M}_u \, \overset{{\rm law}}{=} \,
u^{1/\alpha}\overline{M}_1.\]
On the other hand, Lemma \ref{lemma:Control_P_N} in the Appendix (see also Lemma A.$2$ in \cite{Frikha:Konakov:Menozzi21}) ensures the existence of a family $\{\overline{P}_u\}_{u\ge 0}$ of probability measures  such that
\begin{equation}
\label{proof_control_N_segnato}
P_{\tilde{N}^{\tau,\xi,t,s}_u}(\mathcal{A}) \, \le \, C\overline{P}_u(\mathcal{A}), \quad \mathcal{A} \in \mathcal{B}(\R^N),
\end{equation}
for some positive constant $C$ independent from the parameters $\tau,\xi,t,s$.\newline
For any fixed $u\ge0$, let us now denote by $\overline{N}_u$ the random variable with law $\overline{P}_u$ that is independent from $\overline{M}_u$. Thanks to the representation of the measure $\overline{P}_u$ in \eqref{eq:representation_P_N_segnato}, it is then immediate to check that
\[\overline{N}_u \, \overset{(\text{law})}{=} \,
u^{1/\alpha}\overline{N}_1.\]
We can finally define the family $\{\overline{p}(u,\cdot)\}_{u\ge0}$ of densities as
\begin{equation}
\label{proof_decomposition_p_segnato}
\overline{p}(u,z) \,:=\, \int_{\R^N}p_{\overline{M}}(u,z-w)\overline{P}_{u}(dw),
\end{equation}
which corresponds to the density of the following random variable:
\[\overline{S}_u \, := \, \overline{M}_u +\overline{N}_u\]
for any fixed $u\ge 0$.
Using Fourier transform and the already proven $\alpha$-selfsimilarity of $\overline{M}$ and $\overline{N}$, we now show that
\[\overline{S}_u \, \overset{(\text{law})}{=} \, u^{1/\alpha}\overline{S}_1,\]
or equivalently, that
\[\overline{p}(u,z)\, = \, u^{-N/\alpha}\overline{p}(1,u^{-1/\alpha}z)\]
for any $u \ge 0$ and any $z$ in $\R^N$. Moreover,
\[\mathbb{E}[\vert \overline{S}_u\vert^\gamma] \, = \, \mathbb{E}[\vert \overline{M}_u+\overline{N}_u\vert^\gamma] \, = \,
Cu^{\gamma/\alpha}\bigl(\mathbb{E}[\vert \overline{M}_1\vert^\gamma]+\mathbb{E}[\vert \overline{N}_1\vert^\gamma]\bigr) \, \le \,  Cu^{\gamma/\alpha},\]
for any $\gamma <\alpha$. In particular, Equation \eqref{equation:integration_prop_of_q} holds. We emphasize that the integrability constraints precisely come from the Poisson measure $\overline{P}_u$ which behaves as the one associated with the large jumps of an $\alpha$-stable density.\newline
Equation \eqref{eq:derivative_prop_of_q} now follows easily from the previous arguments. From Equation \eqref{proof:decomposition_density_S}, we start noticing that Controls  \eqref{proof:control_pM_segnato}, \eqref{proof_control_N_segnato} and \eqref{proof_decomposition_p_segnato} imply that for any $k$ in $\llbracket 0,2\rrbracket$,
\[\left| D^k_{z} p_{\tilde{S}^{\tau,\xi,t,s}}(u,z) \right| \, \le \, Ct^{-k/\alpha}\overline{p}(u,z), \quad u\ge 0, \, z \in \R^N,\]
for some constant $C>0$, independent from the parameters $\tau,\xi,t,s$. Recalling the decomposition in \eqref{eq:decomposition_measure_1}, Equation \eqref{eq:derivative_prop_of_q} for $k=0$ already follows.\newline
To show instead the case $k=1$, we can write that
\[
\begin{split}
\bigl{\vert} D_{x_i}\tilde{p}^{\tau,\xi}(t,s,\,&x,y) \bigr{\vert} \, = \, \left| \frac{1}{\det (\mathbb{M}_{s-t})} D_{x_i} \left[p_{\tilde{S}^{\tau,\xi,t,s}}\bigl(s-t,\mathbb{M}^{-1}_{s-t} (y-\tilde{m}^{\tau,\xi}_{s,t}(x)))\right] \right| \\
&= \,\left| \frac{1}{\det (\mathbb{M}_{s-t})} \langle D_z p_{\tilde{S}^{\tau,\xi,t,s}}\left(s-t,\cdot\right)(\mathbb{M}^{-1}_{s-t}
(y-\tilde{m}^{\tau,\xi}_{s,t}(x))),D_{x_i}\mathbb{M}^{-1}_{s-t}\tilde{m}^{\tau,\xi}_{s,t}(x)\rangle \right| \\
&= \, \frac{(s-t)^{-1/\alpha}}{\det (\mathbb{T}_{s-t})}  \overline{p}\left(1,\mathbb{T}^{-1}_{s-t}
(y-\tilde{m}^{\tau,\xi}_{s,t}(x))\right)\left|D_{x_i}\mathbb{M}^{-1}_{s-t}\tilde{m}^{\tau,\xi}_{s,t}(x) \right|,
\end{split}
\]
where in the last step we exploited the $\alpha$-scaling property of $\overline{p}$. From Equation \eqref{eq:dif_differential_m_tilde}, we now notice that the function $x\to \tilde{m}^{\tau,\xi}_{s,t}(x)$ is affine, so that
\[\left| D_{x_i}\mathbb{M}^{-1}_{s-t}\tilde{m}^{\tau,\xi}_{s,t} (x) \right|  \, \le \,  C(s-t)^{-(i-1)}.\]
Hence, it follows that
\[\vert D_{x_i}\tilde{p}^{\tau,\xi}(t,s,x,y) \vert \, \le \, C\frac{(s-t)^{-\frac{1+\alpha(i-1)}{\alpha}}}{\det
(\mathbb{T}_{s-t})} \overline{p}\left(1,\mathbb{T}^{-1}_{s-t}
(y-\tilde{m}^{\tau,\xi}_{s,t}(x))\right).\]
The other case ($k=2$) can be derived in an analogous way.
\end{proof}

We conclude this section with a useful control on the powers of the density $\overline{p}(u,z)$.

\begin{corollary}
\label{coroll:control_power_q_density}
Let $q\ge 1$. Then, there exists a positive constant $C:=C(q)$ such that
\begin{equation}
\label{eq:control_power_q_density}
[\overline{p}(u,z)]^q \, \le \, u^{(1-q)\frac{N}{\alpha}}C\overline{p}(u,z), \quad (u,z) \in (0,T]\times \R^N.
\end{equation}
\end{corollary}
\begin{proof}
We start noticing that we can assume without loss of generality that $u=1$, thanks to the scaling property of $\overline{p}(u,z)$ in Proposition \ref{prop:Smoothing_effect}. Moreover, we know that there exists a constant $K$ such that $\overline{p}(1,z)\le 1$ for any $z$ in $B^c(0,K)$, since $\overline{p}(1,\cdot)$ is a density. It then clearly follows that
\[[\overline{p}(1,z)]^q \, \le \, \overline{p}(1,z), \quad z \in B^c(0,K).\]
On the other hand, we recall that $\overline{p}(1,\cdot)$ is continuous. For any $z$ in $B(0,K)$, it then holds that
\[[\overline{p}(1,z)]^q \, = \, \overline{p}(1,z) [\overline{p}(1,z)]^{q-1} \, \le \, C\overline{p}(1,z),\]
where $C$ is the maximum of $[\overline{p}(1,\cdot)]^q$ on $B(0,K)$.
\end{proof}

\subsection{Regularity of the density along the terminal condition} \label{SEZIONE_SPIEGAZIONE_CONGELAMENTO_E_ALTRI}
We briefly explain here how we want to prove the well-posedness of the martingale formulation associated with $\partial_s+L_s$ at some starting point $(t,x)$. We will mainly focus on the problem of uniqueness since the existence of a solution can be easily handled from already known results. Indeed, we recall that under the assumptions we consider, the main part of the operator $L_s$ is of order $\alpha>1$ while the perturbation is sub-linear. Thus, the existence of a solution  can be obtained, for example, from Theorem $2.2$ in \cite{Stroock75}.
\newline
In particular,  uniqueness for the martingale problem will follow  once the Krylov-like estimates \eqref{eq:Krylov_Estimates} have been shown.\newline
Starting from a solution $\{X^{t,x}_s\}_{s\in[0,T]}$ of the martingale problem with starting point $(t,x)$, the idea is to exploit the properties of the frozen dynamics $\{\tilde{X}^{\tau,\xi,t,x}_s\}_{s\in[0,T]}$ in \eqref{eq:integral_represent_frozen_SDE}. For this reason, let us denote by $\tilde{L}^{\tau,\xi}_s$ its infinitesimal generator and define for $f$ in $ C_c^{1,2}([0,T)\times \R^{N})$ the associated Green kernel:
$$\tilde{G}^{\tau,\xi}(t,x)=\int_t^T ds \int_{\R^N}\tilde{p}^{\tau,\xi}(t,s,x,y) f(s,y)dy.$$
Standard results now give that
\begin{equation}
\label{eq:intro_MP}
    \left(\partial_t+\tilde{L}^{\tau,\xi}_t\right)  \tilde{G}^{\tau,\xi}f(t,x) \, = \, -f(t,x), \quad (t,x) \in [0,T)\times \R^N,
\end{equation}
for any (fixed) freezing parameters $(\tau,\xi)$.\newline
The first step of our method then consists in applying the It\^o formula on the function $\tilde{G}^{\tau,\xi} f$, which is indeed smooth enough, and the solution process $\{X^{t,x}_s\}_{s\in[0,T]}$:
\[
\tilde{G}^{\tau,\xi} f(t,x)+ \mathbb{E}\left[\int_t^{T} (\partial_s+L_s)\tilde{G}^{\tau,\xi} f(s,X^{t,x}_s)   \,ds\right] \, =\, 0.
\]
Exploiting \eqref{eq:intro_MP}, we can then write
\[
\tilde{G}^{\tau,\xi} f (t,x)- \mathbb{E}\Bigl[\int_t^{T} f(s,X^{t,x}_s)\,  ds \Bigr]+\mathbb{E}\Bigl[\int_{t}^T\bigl(L_s-\tilde{L}^{\tau,\xi}_s\bigr)                          \tilde{G}^{\tau,\xi} f(s,X^{t,x}_s) \, ds\Bigr] \, = \, 0
\]
or, equivalently,
\[
\mathbb{E}\Bigl[\int_t^{T} f(s,X^{t,x}_s) \, ds \Bigr] \, = \, \tilde{G}^{\tau,\xi} f
(t,x)+\mathbb{E}\Bigl[\int_{t}^T \bigl(L_s-\tilde{L}^{\tau,\xi}_s\bigr) \tilde{G}^{\tau,\xi} f(s,X^{t,x}_s)\, ds\Bigr].\]
While an estimate of the frozen Green kernel $\tilde{G}^{\tau,\xi} f$ can be obtained from Proposition \ref{prop:Smoothing_effect}, the main difficulty of our approach will be to control, uniformly in $(t,x)$, the following quantity:
\[\int_{t}^T\int_{\R^N}\bigl(L_s-\tilde{L}^{\tau,\xi}_s\bigr) \tilde{G}^{\tau,\xi} f(s,x)\, ds.\]
Focusing for example only on the component associated with the deterministic drift $F$, i.e.
\[\int_t^T\int_{\R^N}\langle F(t,x)-F(t,\theta_{t,\tau}(\xi)),D_x \tilde{p}^{\tau,\xi}(t,s,x,y)\rangle f(s,y)\, dyds,\]
it is clear that we need some kind of compatibility between the arguments of the drift $F$ and those of the frozen density $\tilde{p}^{\tau,\xi}(t,s,x,\cdot)$, in order to exploit the associated smoothing effect (Proposition \ref{prop:Smoothing_effect}). Namely, we need to compare the quantities $(x-\theta_{t,\tau}(\xi))$ and $(y-\tilde{m}^{\tau,\xi}_{s,t}(x))$.\newline
Noticing that for $\tau=s$ and $\xi=y$, $(y-\tilde{m}^{\tau,\xi}_{s,t}(x))=\theta_{t,s}(y)-x$, it follows from Proposition \ref{prop:Smoothing_effect} that this choice of freezing parameters gives the natural compatibility between the difference of the generators and the upper-bounds of the derivatives of the corresponding proxy.

The above reasoning requires however a more thorough analysis on the properties of the ``density'' $\tilde{p}^{s,y}(t,s,x,\cdot)$ frozen along the terminal condition $(\tau,\xi)$. Indeed, the freezing parameter $y$ appears also as the integration variable. In other words, with this approach, the freezing parameter cannot be fixed once for all. The present section is precisely dedicated to the handling of such a choice. This will lead us to introduce a \textit{pseudo} Green kernel, see \eqref{eq:def_Green_Kernel} below, from which will then derive uniqueness to the martingale problem following the Stroock and Varadhan approach (see Chapter 7 in \cite{book:Stroock:Varadhan79}), through appropriate inversion in $L_t^q-L_x^p $ spaces,  proving that the remainder has a small corresponding norm.

We start with a lemma showing the existence of at least one version of the flow $\theta_{t,s}(y)$ which is measurable in $s$ and $y$. This result will be fundamental to make licit any integration of this flow along the terminal condition $y$.

\begin{lemma}\label{lemma:measurability_flow}
There exists a measurable mapping $\theta \colon [0,T]^2\times \R^N \to\R^N$ such that
\begin{equation}
\label{eq:measurability_flow}
\theta(t,s,z)\, := \, \theta_{t,s}(z)\, =    z + \int_{t}^{s}\bigl[ A_u \theta_{u,s}( z)+F(u, \theta_{u,s}( z))\bigr] \, du.
\end{equation}
\end{lemma}
\begin{proof}
The result can be obtained from \cite{Zubelevich12} and a standard compactness argument.
\end{proof}
From this point further, we assume without loss of generality to have chosen such a measurable version $\theta_{t,s}(x)$ of the reference flow.

The next Lemma \ref{lemma:bilip_control_flow} (Approximate Lipschitz condition of the flows) will be a key technical tool for our method. Roughly speaking, it says that a kind of equivalence between the rescaled forward and backward flows appears even in our framework (where the drift $F$ is not regular enough), up to an additional constant contribution, for any two measurable flows satisfying Equation \eqref{eq:def_Cauchy_Peano_flow}. We only remark that similar results has been thoroughly exploited in \cite{Delarue:Menozzi10, Menozzi11, Menozzi18} when considering Lipschitz drifts or \cite{Chaudru:Menozzi17} in the degenerate diffusive setting with H\"older coefficients. \newline
This \emph{approximated} Lipschitz property will be fundamental later on in the proof of Lemma \ref{convergence_dirac} (Dirac Convergence of frozen density) below. It will be proved in Appendix $A.1$, adapting the lines of \cite{Chaudru:Menozzi17}.

\begin{lemma}\label{lemma:bilip_control_flow}
Let $\theta \colon [0,T]^2\times \R^N \to\R^N$, $\check{\theta} \colon [0,T]^2\times \R^N \to\R^N$ be two measurable flows satisfying Equation \eqref{eq:measurability_flow}. Then, there exist two positive constants $(C,C'):=(C,C')(T)$ such that for any $t<s$ in $[0,T]$ and any $x,y$ in $\R^N$,
\begin{equation}
\label{eq:bilip_control_flow}
C^{-1} |\mathbb{T}_{s-t}^{-1}(\check{\theta}_{s,t}(x)-y)| - C' \, \le \, |\mathbb{T}_{s-t}^{-1}(x-\theta_{t,s}(y))|\,\le\, C\left[|\mathbb{T}_{s-t}^{-1}(\check{\theta}_{s,t}(x)-y)| +1\right].
\end{equation}
\end{lemma}
From the above lemma, we also derive the following important estimate for the rescaled difference between the forward flow $\theta_{s,t}(x)$ and the linearized forward dynamics $\tilde{m}^{s,y}_{s,t}(x)$ (defined in \eqref{eq:def_tilde_m}) where the linearization is considered along any backward flow.

\begin{corollary}
\label{coroll:control_error_flow}
Let $\theta \colon [0,T]^2\times \R^N \to\R^N$ be a measurable flow satisfying Equation \eqref{eq:measurability_flow}. Then, there exist a positive constant $C:=C(T)$ and $\zeta$ in $(0,1)$ such that for any $t<s$ in $[0,T]$ and any $x,y$ in $\R^N$,
\begin{equation}
\label{eq:control_error_flow}
 |\mathbb{T}_{s-t}^{-1}(\theta_{s,t}(x)-\tilde{m}^{s,y}_{s,t}(x))| \,\le \, C (s-t)^{\frac{1}{\alpha} \wedge \zeta}\left(1+|\mathbb{T}_{s-t}^{-1}(\theta_{s,t}(x)-y)|\right).
\end{equation}
\end{corollary}
\begin{proof}
We start exploiting the differential dynamics given in Equation \eqref{eq:dif_differential_m_tilde} to write that
\begin{equation}
\label{proof:eq:decomposition1}
\begin{split}
   \mathbb{T}_{s-t}^{-1}\left(\theta_{s,t}(x)-\tilde{m}^{s,y}_{s,t}(x)\right)
\, &= \, \mathbb{T}_{s-t}^{-1}
\int_{t}^{s}  \Bigl{\{}\bigl[ F(u,\theta_{u,t}(x))-F(u,\theta_{u,s}(y)) \bigr] \\
&\qquad \quad\qquad \qquad+\left[A_u(\theta_{u,t}(x)-\tilde{m}^{s,y}_{u,t}(x))\right] \Bigr{\}}\, du  \\
&:= \, (\mathcal{I}_{s,t}^{1}+\mathcal{I}_{s,t}^{2})(x,y).
\end{split}
\end{equation}
We start dealing with $\mathcal{I}_{s,t}^{1}(x,y)$. The key idea is to use the sub-linearity of $F$ and the appropriate H\"older exponents.
Namely, using the Young inequality, we derive that
\[
\begin{split}
|&\mathcal{I}_{s,t}^{1}(x,y)| \, \le \, C\sum_{i=1}^n (s-t)^{-\frac{1+\alpha(i-1)}{\alpha}}\sum_{j=i}^n \int_t^{s}  |(\theta_{u,t}(x)-\theta_{u,s}(y) )_{j}|^{\beta^j} \, du\\
&\,\,\le \, C\biggl{\{}(s-t)^{-\frac{1}{\alpha}}\int_t^{s} \left[|\theta_{u,t}(x)-\theta_{u,s}(y) |+1\right] \, du \\
&\quad\quad +\sum_{i=2}^n (s-t)^{-\frac{1+\alpha(i-1)}{\alpha}} \sum_{j=i}^n \int_t^{s} \Bigl[ (s-t)^{-\gamma^j}|((\theta_{u,t}(x)-\theta_{u,s}(y) )_{j})|+(s-t)^{\gamma^j \frac{\beta^j}{1-\beta^j}} \Bigr] \, du\biggr{\}},
\end{split}
\]
for some parameters $\gamma^j>0$ to be specified below. 
Denoting now for simplicity, \[\Gamma_j\,:= \, -\frac{1+\alpha(i-1)}{\alpha}+\gamma^j \frac{\beta^j}{1-\beta^j},\]
we get that
\[
\begin{split}
|\mathcal{I}_{s,t}^{1}(x,
y)|\, &\le \, C\biggl{\{}(s-t)^{\frac{\alpha-1}{\alpha}}+\int_t^{s}|\mathbb{T}_{s-t}^{-1}(\theta_{u,t}(x)-\theta_{u,s}(y) )|\, du \\
& \,\,\quad+\sum_{i=2}^n  \sum_{j=i}^n \int_t^{s} \Bigl[(s-t)^{-i+j -\gamma^j}  \Big(\frac{|((\theta_{u,t}(x)-\theta_{u,s}(y) )_{j})|}{(s-t)^{\frac{1+\alpha(j-1)}{\alpha}}}\Big)+(s-t)^{-\Gamma_j} \Bigr]\, du\biggr{\}}\\
&\le \, C\biggl{\{}(s-t)^{\frac{\alpha-1}{\alpha}}+\int_t^{s}  |\mathbb{T}_(s-t)^{-1}(\theta_{u,t}(x)-\theta_{u,s}(y))|\, du\\
&\quad\,\, +\sum_{i=2}^n  \sum_{j=i}^n \int_t^{s}  \Bigl[(s-t)^{-i+j -\gamma^j}  |\mathbb{T}_{s-t}^{-1}(\theta_{u,t}(x)-\theta_{u,s}(y) )|+(s-t)^{\Gamma_j} \Bigr]\, du\biggr{\}}.
\end{split}
\]
We now use Lemma \ref{lemma:bilip_control_flow} (Approximate Lipschitz condition of the flows) to derive that
\[ |\mathbb{T}_{s-t}^{-1}(\theta_{u,t}(x)-\theta_{u,s}(y) )| \, \le \, C(|\mathbb{T}_{s-t}^{-1}(\theta_{s,t}(x)-y )|+1).\]
We emphasize here that in our current framework we should \textit{a priori} write $\theta_{s,u}(\theta_{u,t}(x)) $ in the above equation since we do not have the flow property. Anyhow, since Lemma \ref{lemma:bilip_control_flow} (Approximate Lipschitz condition of the flows) is valid for any flow starting from $\theta_{u,t}(x) $ at time $u$ associated with the ODE (see Equation \eqref{eq:bilip_control_flow}) we can proceed along the previous one, i.e. $(\theta_{v,t}(x))_{v\in [u,s] }$. The previous reasoning yields that
\begin{equation}
\begin{split}
  \label{PREAL_CTR_I1_EPS}
|\mathcal{I}_{s,t}^{1}(x,y) |\, &\le \,C\biggl{\{}(s-t)^{\frac{\alpha-1}{\alpha}}+(s-t) \bigl[|\mathbb{T}_{s-t}^{-1}(\theta_{s,t}(x)-y)|+1\bigr]\\
& \qquad \qquad \times \Bigl[1+\sum_{i=2}^n  \sum_{j=i}^n  \Bigl((s-t)^{-i+j -\gamma^j}+(s-t)^{-\frac{1+\alpha(i-1)}{\alpha}+\gamma^j \frac{\beta^j}{1-\beta^j}} \Bigr) \Bigr] \biggr{\}}.
\end{split}
\end{equation}
We now choose for $j$ in $\llbracket i,n\rrbracket$,
\[-i+j -\gamma^j \, = \, -\frac{1+\alpha(i-1)}{\alpha}+\gamma^j \frac{\beta^j}{1-\beta^j}\, \Leftrightarrow\, \gamma^j\, = \, \left(j-\frac{\alpha-1}{\alpha}\right)\left(1-\beta^j\right),\]
to balance the two previous contributions associated with the indexes $i,j$. \newline
To obtain a global smoothing effect with respect to ${s-t}$ in \eqref{PREAL_CTR_I1_EPS} we need to impose:
\begin{equation}
\label{eq:natural_threshold}
-i+j-\gamma^j\, >\,-1 \, \Leftrightarrow \, \beta^j\, >\, \frac{1+\alpha(i-2)}{1+\alpha(j-1)}, \quad \forall i \le j.
\end{equation}
Hence, under our assumptions, we have that there exists $\zeta$ in $(0,1)$ depending on $\beta^j$ for any $j\in \llbracket i,n\rrbracket$ such that
\begin{equation}
\label{CTR_I1_EPS}
|\mathcal{I}_{s,t}^{1}(x,y)| \,  \le \, C(s-t)^\zeta\left[1+|\mathbb{T}_{s-t}^{-1}(\theta_{s,t}(x)-y)|\right].
\end{equation}
Recalling from the structure of $A$ that
\[|\mathbb{T}^{-1}_{s-t}A_u\mathbb{T}_{s-t}| \le \, C(s-t)^{-1},\]
Control \eqref{eq:control_error_flow} now follows from
\eqref{proof:eq:decomposition1}, \eqref{CTR_I1_EPS} and the Gr\"onwall lemma.
\end{proof}

Thanks to the Approximate Lipschitz property of the flow presented in Lemma \ref{lemma:bilip_control_flow} above and Corollary \ref{coroll:control_error_flow}, we can now adapt the controls on the derivatives of the frozen density (Proposition \ref{prop:Smoothing_effect}) to the ``density'' $\tilde{p}^{s,y}(t,s,x,y)$. Indeed, we recall again that the function $\tilde{p}^{s,y}(t,s,x,y)$ is not a proper density in $y$ since the integration variable $y$ stands also as freezing parameter. This is one of the main difficulties of the approach.

The following result is the key to our analysis since it precisely quantifies the smoothing effect in time of the proxy we chose.
\begin{corollary}
\label{coroll:Smoothing_effect}
There exists a positive constant $C:=C(N,\alpha)$ such that for any $\gamma$ in $[0,\alpha)$, any $t< s$ in $[0,T]$ and any $x,y$ in $\R^N$,
  \begin{equation}
\label{eq:smoothing_effect_frozen_y}
\int_{\R^N} \frac{|\mathbb{T}^{-1}_{s-t}(\theta_{t,s}(y)-x)|^\gamma}{\det \mathbb{T}_{s-t}}\overline{p}(1,\mathbb{T}^{-1}_{s-t}(\theta_{t,s}(y)-x)) \, dy \, \le \, C.
  \end{equation}
Moreover, if $K>0$ is large enough, it holds that
  \begin{multline}
\label{eq:smoothing_effect_frozen_y_GENERIC_FUNCTION}
\int_{\R^N} \mathds{1}_{|{\mathbb{T}^{-1}_{s-t}(\theta_{t,s}(y)-x)}|\ge K} \frac 1{\det \mathbb{T}_{s-t}}\overline{p}(1,\mathbb{T}^{-1}_{s-t}(\theta_{t,s}(y)-x)) \, dy \\
\le \, C \int_{\R^N} \mathds{1}_{|z|\ge \frac K2}\check p(1,z) dz,
  \end{multline}
  where $\check p $ enjoys the same integrability properties as $\overline{p}$ (stated in Proposition \ref{prop:Smoothing_effect}).
\end{corollary}

The strengthened assumptions concerning the integrability thresholds in Theorem \ref{thm:main_result} with respect to the natural ones appearing in \eqref{eq:natural_threshold} might seem awkward at first sight.
It is actually the specific current framework, which involves as a proxy a stochastic integral with respect to a stable-like jump process and its associated iterated integrals that leads to additional constraints on the regularity indexes needed for our method to work.

The natural approach to get rid of the flow involving the integration variable in \eqref{eq:smoothing_effect_frozen_y} would have been to use the approximate Lipschitz property of the flow established in Lemma \ref{lemma:bilip_control_flow}. This indeed readily yields that:
$$ |\mathbb{T}^{-1}_{s-t}(\theta_{t,s}(y)-x)|^\gamma\le C(1+ |\mathbb{T}^{-1}_{s-t}(y-\theta_{s,t}(x))|^\gamma).$$
The main difficulty is that we  do not actually succeed in establishing in whole generality that:
\begin{equation}\label{equiv_density_flows}
\overline{p}(1,\mathbb{T}^{-1}_{s-t}(\theta_{t,s}(y)-x))\le C \check{p}(1,\mathbb{T}^{-1}_{s-t}(y-\theta_{s,t}( x)),
\end{equation}
for a density $\check p $ which shares the same integrability properties as $\bar p$.

Equation \eqref{equiv_density_flows} is absolutely direct in the diffusive setting from the explicit form of the Gaussian density and it has been thoroughly used in \cite{Chaudru:Menozzi17} to derive sharp thresholds for weak uniqueness. It is clear that the above control has to be considered point-wise and one of the huge difficulties with stable type processes consists in describing precisely their tail behavior which is actually very much related to the geometry of their corresponding spectral measure on the sphere. We refer to the seminal work of Watanabe \cite{Watanabe07} for a precise description of the tails in terms of the dimension of the support of the spectral measure, in the stable case, and to the extension by Sztonyk \cite{Sztonyk10} for the tempered stable case. The delicate point comes of course from the behavior of the Poisson measure (large jumps) as illustrated in  the following computation. From \eqref{proof:control_pM_segnato} and \eqref{proof_decomposition_p_segnato}, we write that
\begin{align*}
\overline{p}(1,\mathbb{T}^{-1}_{s-t}(\theta_{t,s}(y)-x))\, &= \, \int_{\R ^N} p_{\bar M}(1,\mathbb{T}^{-1}_{s-t}(\theta_{t,s}(y)-x)-w) \, \overline{P}_1(dw)\notag\\
&\le \, C\int_{\R^N} \frac{1}{(C+|\mathbb{T}_{s-t}^{-1}(\theta_{t,s}(y)-x)-w|)^M} \, \overline{P}_1(dw).
\end{align*}
\begin{trivlist}
\item[$\bullet$]
Let us first emphasize that, when $|\mathbb{T}^{-1}_{s-t}(\theta_{t,s}(y)-x)|\le K $ (\textit{diagonal type regime}) for some fixed $K$, then Control \eqref{equiv_density_flows} holds. Indeed, since from Corollary \ref{coroll:control_error_flow},
\[|\mathbb{T}_{s-t}^{-1}(\tilde m_{s,t}^{t,y}(x)-\theta_{s,t}(x))|\, \le \, \tilde C(s-t)^{\frac{1}{\alpha} \wedge \zeta}\left(1+|\mathbb{T}_{s-t}^{-1}(\theta_{s,t}(x)-y)|\right),\]
we would get, recalling from Lemma \ref{lemma:identification_theta_m}, Equation \eqref{eq:identification_theta_m1} that $\theta_{t,s}(y)-x =y-\tilde{m}^{s,y}_{s,t}(x)$, that
\begin{align*}
\overline{p}(&1,\mathbb{T}^{-1}_{s-t}(\theta_{t,s}(y)-x))\, \le \, C\int_{\R ^N} \frac{1}{(C+|\mathbb{T}_{s-t}^{-1}(\theta_{t,s}(y)-x)-w|)^M} \,\overline{P}_1(dw)\\
 &\le \,C\int_{\R^N}\frac{1}{([C+ |\mathbb{T}_{s-t}^{-1}(y-\theta_{s,t}(x))-w|-(s-t)^{\frac{1}\alpha\wedge \zeta}|\mathbb{T}_{s-t}^{-1}(\theta_{s,t}(x)-y)|]\vee 1)^M} \,\overline{P}_1(dw)\\
 &\le \,C\int_{\R^N}\frac{1}{([\check C+ |\mathbb{T}_{s-t}^{-1}(y-\theta_{s,t}(x))-w|)^M} \, \overline{P}_1(dw)\\
 &=: \, \check p(1,\mathbb{T}^{-1}_{s-t}(y-\theta_{s,t}(x))),
\end{align*}
and $\check p $  plainly satisfies the required integrability conditions.
These computations actually emphasize that \eqref{equiv_density_flows} holds, up to a modification of $\check C $ above, up to the threshold \[|\mathbb{T}^{-1}_{s-t}(\theta_{t,s}(y)-x)|\, \le\, c_0(s-t)^{-(\frac{1}{\alpha}\wedge \zeta)},\]
for some $c_0>0$ small enough with respect to $C$. It would therefore remain to investigate the complementary \textit{very off-diagonal regime}.
\item[$\bullet $] Let us now concentrate on the \textit{off-diagonal regime} $|\mathbb{T}^{-1}_{s-t}(\theta_{t,s}(y)-x)|> K  $. In that case, we write:
\begin{align}
\overline{p}(1,\mathbb{T}^{-1}_{s-t}(\theta_{t,s}(y)-x))
&\, \le \, C\int_{\R ^N} \frac{1}{(C+|\mathbb{T}_{s-t}^{-1}(\theta_{t,s}(y)-x)-w|)^M} \overline{P}_{1}(dw)\notag\\
&\le \, C\int_{0}^1\overline{P}_{1}(\{w\in \R^N\colon(1+|\mathbb{T}_{s-t}^{-1}(\theta_{t,s}(y)-x)-w|)^{-M} |>u  \})du\notag\\
&\le \, C\int_{0}^1\overline{P}_{1}(B(\mathbb{T}_{s-t}^{-1}(\theta_{t,s}(y)-x), u^{-1/M})du. \label{eq:spiegazione1}
 \end{align}
 It now follows from the proof of Proposition \ref{prop:Decomposition_Process_X} that the support of the spectral measure on $\mathbb S^{N-1} $ associated with  $\{\tilde{S}^{\tau,\xi,t,s}_{u}\}_{u\ge 0}$ has dimension $d$. The related concentration properties also transmit to $\bar N_1 $ (see the proofs of Proposition \ref{prop:Smoothing_effect}
 and Lemma \ref{lemma:Control_P_N}). Thus, we get from \cite{Watanabe07}, \cite{Sztonyk10} (see respectively Lemma 3.1 and Corollary 6 in those references) that there exists a constant $C>0$ such that for all $z$ $\R^{N}$ and $r>0$:
\begin{equation}
\label{EST_POISSON}
\overline{P}_{1}(B(z,r))\le Cr^{d+1} (1+r^\alpha)|z|^{-(d+1+\alpha)}.
\end{equation}
In other words, the global bound is given by the worst decay deriving from the dimension of the support of the spectral measure.
In the current case $|z|\ge K $, this bound is clearly of interest for \textit{large} values of $z$. Hence, from \eqref{eq:spiegazione1} and \eqref{EST_POISSON}, it holds that
 \[
 \begin{split}
\overline{p}(1,\mathbb{T}^{-1}_{s-t}(\theta_{t,s}(y)&-x)) \, \le \, C\int_0^1u^{-(d+1)/M}
 (1+
 u^{-\alpha/M})\, du
|\mathbb{T}_{s-t}^{-1}(\theta_{t,s}(y)-x)|^{-(d+1+\alpha)}
\\
&\le \, C (1+|\mathbb{T}_{s-t}^{-1}(\theta_{t,s}(y)-x)|)^{-(d+1+\alpha)}
\int_0^1 [u^{-(d+1)/M}+u^{-(d+1+\alpha)/M} du],
\end{split}
\]
Choosing  $M>d+1+\alpha $ then gives that there exists $C\ge 1$ such that
$$
\overline{p}(1,\mathbb{T}^{-1}_{s-t}(\theta_{t,s}(y)-x)) \le  C (1+|\mathbb{T}_{s-t}^{-1}(\theta_{t,s}(y)-x)|)^{-(d+1+\alpha)}
.$$
We thus get from Lemma \ref{lemma:bilip_control_flow}, up to a modification of $C$, that:
\begin{align}\label{BAD_BOUND}
\overline{p}(1,\mathbb{T}^{-1}_{s-t}(\theta_{t,s}(y)-x)) \le  C (1+|\mathbb{T}_{s-t}^{-1}(y-\theta_{s,t}(x))|)^{-(d+1+\alpha)}.
\end{align}
This actually leads to strong dimension constraints for this bound to be integrable. This phenomenon already appeared e.g. in  \cite{Huang:Menozzi16}
and induced therein to consider $d=1,n=3$ at most to address the well posedness of the martingale problem associated with a linear drift and a multiplicative isotropic stable noise. Those thresholds and dimension constraints remain with this approach. \newline
Actually, from the threshold appearing in \eqref{eq:thresholds_beta}, we would like to consider the left-hand side of \eqref{eq:smoothing_effect_frozen_y} with $\gamma>\frac{1+\alpha}{1+2\alpha}$ corresponding to $j=3=n$ therein. From Control \eqref{BAD_BOUND}, this would require $-\frac{1+\alpha}{1+2\alpha}+(d+1+\alpha)>3, d=1 \iff \alpha^2-\alpha-1>0$, which in our framework imposes that $\alpha\in (\frac{1+\sqrt 5}2,2)$.

Another possibility would have been, in the tempered case, to keep track of the tempering function, instead of bounding $\tilde p^{\tau,\xi} $ by a self-similar density $\bar p $, in order to benefit from the tempering at infinity to compensate the bad concentration rate in \eqref{BAD_BOUND}. However, see \cite{Huang:Menozzi16} and \cite{Sztonyk10}, we would have obtained bounds of the form
$$\tilde p^{\tau,\xi} (t,s,x,y)\le  C (1+|\mathbb{T}_{s-t}^{-1}(y-\theta_{s,t}(x))|)^{-(d+1+\alpha)}Q\left( |\mathbb M_{s-t}^{-1}(y- \theta_{s,t}(x))|\right).
$$
Such a bound will give space integrability but deteriorates as well the time-integrability.
This difficulty would occur even in the truncated case, thoroughly studied in the non-degenerate case by Chen \textit{et al.} \cite{Chen:Kim:Kumagai08}.
Thus, we will develop here another approach.
%
\end{trivlist}
%
Namely, we would like to change variable to $\bar y:=\mathbb{T}^{-1}_{s-t}(\theta_{t,s}(y)-x) $ in the left-hand side of Equation \eqref{eq:smoothing_effect_frozen_y}. Of course, this is not bluntly possible since the coefficients at hand are not smooth enough. The point is then to introduce a flow $\theta_{t,s}^\delta(y) $ associated with mollified coefficients (for which the difference with respect to the initial flow will be controlled similarly to what is done to establish the approximate Lipschitz property of the flows in Lemma \ref{lemma:bilip_control_flow}) and then, to control  $\det(\nabla \theta_{t,s}^\delta(y)) $ (see Lemma \ref{LEMMA_FOR_DET} below). Since we do not have here a summation with respect to the single rescaled components as in the previous Lemma \ref{lemma:bilip_control_flow} above or as in Corollary \ref{coroll:control_error_flow}, this will conduct to reinforce our assumptions and suppose that $(F_i)_{i\in \llbracket 2,n\rrbracket} $ has the same regularity with respect to the variable $x_j$, $j\in \llbracket 2,n\rrbracket$, whatever the level of the chain.  This is precisely  what leads to consider the condition
\begin{itemize}
\item $x_j \to F_i(t,x_i,\dots,x_j,\dots,x_n)$ is $\beta^j$-H\"older continuous, uniformly in $t$ and in $x_k$ for $k\neq j$, with
 \begin{equation*}
\beta^j \, > \, \frac{1+\alpha(j-2)}{1+\alpha(j-1)}.
\end{equation*}
\end{itemize}


For the sake of clarity the proof of Corollary \ref{coroll:Smoothing_effect} is postponed to the Appendix.

Let us introduce now some useful tools for the study of the martingale problem for $\partial_s+L_s$. The first step is to consider a suitable Green-type kernel associated with the frozen density $\tilde{p}^{s,y}$ and establish which Cauchy-like problem it solves. Namely, we define for any function $f\colon [0,T]\times \R^N \to \R$ regular enough, the \textit{pseudo} Green kernel $\tilde{G}_\epsilon$ given by:
\begin{equation}\label{eq:def_Green_Kernel}
\tilde{G}_\epsilon f(t,x) \, := \, \int_{(t+\epsilon)\wedge T}^T \int_{\R^N}  \tilde{p}^{s,y}(t,s,x,y) f(s,y) \, dy ds, \quad (t,x) \in [0,T)\times \R^N,
\end{equation}
where $\epsilon$ is meant to be small.\newline
We only remark that the above Green kernel $\tilde{G}_\epsilon$ is well-defined, since the frozen density $\tilde{p}^{s,y}(t,s,
x,y)$ is measurable in $(s,y)$ thanks to Lemma \ref{lemma:measurability_flow} (measurability of the flow in these parameters).

\begin{prop}
\label{prop:pointwise_green_kernel}
Let $p$, $q$ in $(1,+\infty)$ such that the integrability Condition $(\mathscr{C})$ holds. Then, there exists a positive constant $C:=C(T,p,q)$  such that for any $f$ in
$L^p\big(0,T;L^q(\R^N)\bigr)$,
\[
\Vert \tilde{G}_{\epsilon}f\Vert_\infty \, \le \, C \Vert f \Vert_{L^p_tL^q_x}.\]
Moreover, it holds that $\lim_{T\to 0}C(T,p,q) = 0$.
\end{prop}
\begin{proof}
We start using the H\"older inequality in order to split the component with $f$ and the part with the density $\tilde{p}(t,s,x,y)$:
\[
\begin{split}
|\tilde{G}_\epsilon f(t,x)| \,
&\le \, C \Vert f \Vert_{L^p_tL^q_x} \left(\int_{(t+\epsilon)\wedge T}^T \left (\int_{\R^N} \left|\tilde{p}^{s,y}(t,s,x,y)\right|^{q'} \,dy \right)^{\frac{p'}{q'}} ds \right)^\frac{1}{p'} \\
&=: \, C \Vert f \Vert_{L^p_tL^q_x}|I_\epsilon(t,x)|,
\end{split}
\]
where we have denoted by $p'$, $q'$ the conjugate of $p$ and $q$, respectively.\newline
In order to control the remainder term $I_\epsilon(t,x)$, we now apply \eqref{eq:derivative_prop_of_q} from Proposition \ref{prop:Smoothing_effect} with $k=0$ and $(\tau,\xi)=(s,y)$ to write that
\[
|I_\epsilon(t,x)|^{p'} \, \le \, C\int_{(t+\epsilon)\wedge T}^T \left(\int_{\R^N} \left(\frac{1}{\det \mathbb{T}_{s-t}} \overline{p}\left(1, \mathbb{T}^{-1}_{s-t}(y-\tilde{m}^{s,y}_{s,t}(x))\right)\right)^{q'} \, dy \right)^{\frac{p'}{q'}} ds,
\]
where we recall that $\mathbb{T}_t=t^{1/\alpha}\mathbb{M}_t$ (see \eqref{eq:def_matrix_T} and \eqref{eq:def_matrix_M}).\newline
From Corollaries \ref{coroll:control_power_q_density} and \ref{coroll:Smoothing_effect}, we then write that
\[\begin{split}
|I_\epsilon(t,x)|^{p'} \, &\le \,
C\int_{(t+\epsilon)\wedge T}^T \left(\int_{\R^N} \frac{1}{(\det \mathbb{T}_{s-t})^{q'}} \overline{p}\left(1, \mathbb{T}^{-1}_{s-t}(y-\tilde{m}^{s,y}_{s,t}(x))\right) \, dy \right)^{\frac{p'}{q'}} ds\\
&\le \, C \int_{(t+\epsilon)\wedge T}^T  (\det \mathbb{T}_{s-t})^{\frac{p'}{q'}-p'} ds\, =  \, C \int_{(t+\epsilon)\wedge T}^T   \frac{1}{(\det \mathbb{T}_{s-t})^{\frac{p'}{q}}} ds.
\end{split}
\]
Since by definition of matrix $\mathbb{T}_t$, it holds that
\begin{equation}
\label{eq:control_det_T1}
\det \mathbb{T}_{s-t} \, = \, (s-t)^{\sum_{i=1}^nd_i\frac{1+\alpha(i-1)}{\alpha}},
\end{equation}
we can conclude that under the integrability assumption $(\mathscr{C})$, we have that
\begin{equation}
\label{eq:control_det_T2}
\bigl(\sum_{i=1}^nd_i\frac{1+\alpha(i-1)}{\alpha}\bigr)\frac{p'}{q} \, < \, 1  \, \Leftrightarrow \, \bigl(\sum_{i=1}^nd_i\frac{1+\alpha(i-1)}{\alpha}\bigr)\frac{1}{q}+\frac{1}{p} \, <\,  1.
\end{equation}
The proof is complete.
\end{proof}

Now, we understand which Cauchy-like problem is solved by the ``density'' $\tilde{p}^{s,y}(t,s,x,y)$ frozen at the terminal point $(s,y)$. We start denoting by $\tilde{L}^{s,y}_t$ the infinitesimal generator of the proxy process $\{\tilde X_{s}^{s,y,t,x}\}_{s\in [t,T]}$. For any smooth function $\phi\colon \R^N\to\R$, it writes:
\begin{equation}
\label{eq:def_frozen_generator}
\begin{split}
 \tilde{L}^{s,y}_t\phi(x)\, &:= \,
\langle A_t x+\tilde{F}^{s, y}_t, D_x \phi(x)\rangle + \tilde{\mathcal{L}}^{s,y}_t \\
&:= \langle A_t x+\tilde{F}^{s, y}_t, D_x \phi(x)\rangle + \int_{\R^d_0}\bigl[\phi(x+B\tilde{\sigma}^{s, y}_t w)-\phi(x) \bigr] \,\nu(dw),
\end{split}
\end{equation}
where, we recall, $\tilde{F}^{s,y}_t:=F(t,\theta_{t,s}(y))$ and $\tilde{\sigma}^{s,y}_t:=\sigma(t,\theta_{t,s}(y))$.\newline
By direct calculation, it is not difficult to check now that for any $(s,x,y)$ in $[0,T]\times \R^{2N}$ it holds that
\begin{equation}\label{eq:differential-eq}
 \left(\partial_t + \tilde{L}^{s,y}_t\right) \tilde{p}^{s,y} (t,s,x,z) \,  = \, 0, \quad (t,z) \in [0,s)\times \R^N.
\end{equation}
However, we carefully point out that some attention is requested to establish the following lemma, which is crucial
to derive which Cauchy-type problem the function $\tilde{G}f := \lim_{\epsilon\to 0}\tilde{G}_\epsilon f$ actually solves. In particular, it is important to highlight that Lemma \ref{convergence_dirac} (Dirac Convergence of frozen density) below cannot be obtained directly from the convergence in law of the frozen process $\tilde{X}_{s}^{s,y,t,x}$ towards the Dirac mass (cf.\ Equation
\eqref{eq:differential-eq}). Indeed, the integration variable $y$ also appears as a freezing parameter which makes the argument more complicated.\newline
The proofs of the following two lemmas is quite involved and technical. For this reason, we decided to postpone them to the Appendix, Section \ref{SEC_TEC_LEMMA_APP}.

\begin{lemma}
\label{convergence_dirac}
Let $(t,x)$ be in $[0,T)\times \R^N$ and $f\colon \R^N\to \R$ a bounded continuous function. Then,
\[
\lim_{\epsilon \to 0}\left| \int_{\R^N} f(y) \tilde{p}^{t+\epsilon,y}(t,t+\epsilon,x,y)\, dy
-f(x) \right| \, = \, 0.
\]
Moreover, the above limit is uniform with respect to $t$ in $[0,T]$.
\end{lemma}

A similar result involving the $L^p_tL^q_x$-norm can also be obtained. For notational simplicity, let us set
\begin{equation}
\label{eq:def_f_epsilon}
I_\epsilon f(t,x) \, := \, \int_{\R^N} f(t+\epsilon,y)\mathds{1}_{[0,T-\epsilon]}(t)
\tilde{p}^{t+\epsilon,y}(t,t+\epsilon,x,y) \, dy
\end{equation}
for any sufficiently regular function $f\colon [0,T]\times \R^N \to \R$.

\begin{lemma}
\label{prop:convergence_LpLq}
Let $p>1$, $q>1$ and $f$ in $C_c^{1,2}([0,T)\times \R^N)$. Then,
\[
\lim_{\epsilon \to 0} \Vert
I_\epsilon f -f \Vert_{L^{p}_tL^q_x} \, = \, 0.
\]
\end{lemma}

We want now to understand which Cauchy-like problem is solved by our frozen Green kernel $\tilde{G}_\epsilon f(t,x)$. For this reason, we introduce for any function $f$ in $C^{1,2}_0([0,T)\times \R^N,\R)$ the following quantity:
\begin{equation}
\label{eq:def_Mepsilon}
\tilde{M}_\epsilon f(t,x)\, := \, \int_{t+\epsilon}^T \int_{\R^N}  \tilde{L}^{s,y}_t \tilde{p}^{s,y}(t,s,x,y)f(s,y) \, dy ds, \quad (t,x) \in [0,T)\times\R^N,
\end{equation}
for some fixed $\epsilon>0$ that is assumed to be small enough.
Then, we can derive from Equation \eqref{eq:differential-eq} and Proposition \ref{prop:Smoothing_effect}
that the following equality holds:
\begin{equation}\label{eq:differential_eq2}
\partial_t \tilde{G}_\epsilon f(t,x)+ \tilde{M}_\epsilon f(t,x) \, = \, -I_\epsilon f(t,x),\quad (t,x)\in [0,T)\times \R^N,
\end{equation}
where we used the same notation in \eqref{eq:def_f_epsilon} for $I_\epsilon f$. We point out that the localization with respect to $\epsilon$ is precisely needed to exploit directly \eqref{eq:differential-eq} and thus, to derive
\eqref{eq:differential_eq2} for any fixed $\epsilon>0$, by usual dominated convergence arguments. In particular, we point out that in the limit case ($\epsilon \to 0$), the smoothness on $f$ is not a sufficient condition to derive the smoothness of $\tilde{G} f$. This is again due to the dependence of the proxy upon the integration variable.
\setcounter{equation}{0}
\section{Well-posedness of the martingale problem}
\fancyhead[RO]{Section \thesection. Well-Posedness of the martingale problem}
This section is devoted to the proof of the well-posedness of the martingale problem for $\partial_s+L_s$ with initial condition $(t,x)$, under the assumptions of Theorem
\ref{thm:main_result}.

Since by definition the paths of any solution $\{X_t\}_{t\ge 0}$ of the martingale problem for $\partial_s+L_s$ are c\`adl\`ag, it will be convenient
afterwards to give an alternative definition. We denote by $\mathcal{D}[0,\infty)$ the family of all the c\`adl\`ag paths from $[0,\infty)$ to $\R^N$, equipped
with the ``standard'' Skorokhod topology. For further details, we suggest the interested reader to see \cite{book:Bass11}, \cite{book:Ethier:Kurtz86} or \cite{book:Jacod:Shiryaev87}.\newline
Fixed a starting point $(t,x)$ in $[0,\infty)\times \R^N$, we will say that a a probability measure $\mathbb{P}$ on $\mathcal{D}[0,\infty)$ is a solution of the martingale problem for
$\partial_t+L_t$ starting at $(t,x)$ if the coordinate process $\{y_t\}_{t\ge 0}$ on $\mathcal{D}[0,\infty)$, defined by
\[y_t(\omega) \, = \, \omega(t), \quad \omega \in  \mathcal{D}[0,\infty)\]
is a solution (in the previous sense) of the martingale problem for $\partial_s+L_s$.\newline
Similarly, we will say that uniqueness holds for the martingale problem for $\partial_s+L_s$ with starting point $(t,x)$ if
\[\mathbb{P}\circ y^{-1} \, = \,\tilde{\mathbb{P}}\circ y^{-1},\]
for any two solutions $\mathbb{P}$, $\tilde{\mathbb{P}}$ of the martingale problem for $\partial_s+L_s$ starting at $(t,x)$.

The existence of a solution $\mathbb{P}$ of the martingale problem for $\partial_s+L_s$ can be obtained adapting the proof of Theorem $2.2$ in \cite{Stroock75} exploiting the sublinear structure of the drift $F$ and localization arguments in order to deal with possibly unbounded coefficients.

\begin{prop}[Existence]
\label{prop:exist_MP_in_x}
Under the assumptions of Theorem \ref{thm:main_result}, let $(t,x)$ be in $[0,\infty)\times\R^N$. Then, there exists a solution $\mathbb{P}$ of the martingale problem for
$\partial_s+L_s$ starting at $(t,x)$.
\end{prop}

We move to the question of uniqueness for the martingale problem associated with $\partial_s+L_s$. As shown already in the introduction of
Section $3$, the analytical properties on the frozen process $(\tilde{X}^{s,y,t,x}_{u})_{u\in [t,T]}$ we presented there will be the crucial tools for the reasoning in the following section.\newline
We will start proving directly that the Krylov-type estimates \eqref{eq:Krylov_Estimates} holds but first for $p$, $q$ big enough (but finite). It will imply in particular the existence of a density for the canonical process associated with any solution of the martingale problem. As a consequence, the weak well-posedness of SDE \eqref{eq:SDE} under our assumptions can be shown to hold.\newline
Only in a second moment, we will then show that the Krylov estimates holds for \emph{any} $p$, $q$ satisfying condition ($\mathscr{C}$) through a regularization technique. Namely, we regularize the driving noise $Z_t$ by introducing an additional isotropic $\alpha$-stable process depending from a regularizing parameter. Following the previous arguments for the regularized dynamics, we will then prove that the solution process satisfies again the Krylov-type estimates for any $p$, $q$ in the considered range, \emph{uniformly} with respect to the regularizing parameter.\newline
Letting the regularizing parameter go to zero, we will then conclude the proof of Corollary \ref{coroll:Krylov_Estimates}.

\subsection{Uniqueness of the martingale problem}
The first step in proving the uniqueness of the Martingale problem for $\partial_s+L_s$ is to show that any solution to the martingale problem satisfies the Krylov-like estimates in Equation \eqref{eq:Krylov_Estimates}. To do so, we prove that the difference operator between the genuine generator $L_t$ and a suitable associated perturbation (associated with the frozen generator $\tilde{L}^{s,y}_t$ given in \eqref{eq:def_frozen_generator}) has small  $L^p_tL^q_x$-norm when considering a sufficiently small final horizon $T$. Namely, we introduce the following remainder:
\begin{equation}\label{eq:def_R}
\tilde{R}_\epsilon f(t,x) \, := \, (L_t\tilde{G}_\epsilon f-\tilde{M}_\epsilon f)(t,x) \, = \, \int_{t+\epsilon}^T  \int_{\R^N}(L_t -\tilde{L}^{s,y}_t)\tilde{p}^{s,y}(t,s,x,y) f(s,y) \, dy ds,
\end{equation}
for some $\epsilon$ to be small enough. We recall that $\tilde{ G}_\epsilon f, \tilde{M}_\epsilon f $ and $\tilde p^{s,y}(t,s,x,y)$ were defined in \eqref{eq:def_Green_Kernel}, \eqref{eq:def_Mepsilon} and \eqref{eq:representation_density}, respectively.
\newline
We firstly present a point-wise control for the remainder term $\tilde{R}_\epsilon f$. Importantly, the constant $C$ below does not depend on $\epsilon$, allowing to pass to the limit in Equation \eqref{eq:def_R}. This will be discussed at the end of the present section.

\begin{prop}
\label{prop:control_tildeR}
There exist  $q_0>1$, $p_0>1$ and $C:=C(T,p_0,q_0)$ such that for any $q\ge q_0$, $p\ge p_0$ and any $f$ in $L^p([0,T];L^q(\R^N))$, it holds that
\begin{equation}
\label{eq:control_infty_R}
\Vert\tilde{R}_\epsilon f\Vert_\infty\, \le \, C \Vert f \Vert_{L^p_tL^q_x}.
\end{equation}
\end{prop}
\begin{proof}
We start recalling from \eqref{eq:def_generator}-\eqref{eq:def_frozen_generator} (exploiting also the change of truncation in \eqref{eq:change_of_truncation})) that we can decompose $\tilde{R}_\epsilon f$ in the following way:
\begin{equation}
\label{proof:remainder_decomposition}
\begin{split}
   \tilde{R}_\epsilon f(t,x) \, &= \, \int_{t+\epsilon}^T  \int_{\R^N}(\mathcal{L}_t -\tilde{\mathcal{L}}^{s,y}_t)\tilde{p}^{s,y}(t,s,x,y) f(s,y) \, dy ds \\
   &\qquad\qquad + \int_{t+\epsilon}^T  \int_{\R^N}\langle F(t,x)-\tilde{F}^{s,y}_t,D_x\tilde{p}^{s,y}(t,s,x,y)\rangle f(s,y) \, dy ds \\
   &=: \,  \tilde{R}^0_\epsilon f(t,x)+ \tilde{R}^1_\epsilon f(t,x)
\end{split}
\end{equation}
where the operators $\mathcal{L}_t$ and $\tilde{\mathcal{L}^{s,y}_t}$ have been defined in \eqref{eq:def_generator} and \eqref{eq:def_frozen_generator}, respectively. 
 \newline
Since by assumptions, $x_j\to F_i(t,x)$ is $\beta^j$-H\"older continuous, we can control the second term $\tilde{R}^1_\epsilon f$, associated with the difference of the drifts, using Proposition \ref{prop:Smoothing_effect} with $(\tau,\xi)=(s,y)$:
\[
\begin{split}
\bigl|\langle F(t,x)-\tilde{F}^{s,y}_t,D_x&\tilde{p}^{s,y}(t,s,x,y)\rangle \bigr| \, \le \, \sum_{i=1}^n\left|F_i(t,x)-F_i(t,\theta_{t,s}(y))\right|\,  |D_{x_i}\tilde{p}^{s,y}(t,s,x,y)| \\
&\le \, C\sum_{i=1}^n(s-t)^{-\frac{1+\alpha(i-1)}{\alpha}}\frac{\overline{p}(1,\mathbb{T}^{-1}_{s-t}(y-\tilde{m}^{s,y}_{s,t}(x)))}{\det \mathbb{T}_{s-t}}\sum_{j=i}^n\left|(x-\theta_{t,s}(y))_j\right|^{\beta^j} \\
&\le \, C\sum_{i=1}^n\sum_{j=i}^n(s-t)^{{\zeta^j_i}}\left|\mathbb{T}^{-1}_{s-t}(x-\theta_{t,s}(y))\right|^{\beta^j}\frac{\overline{p}(1,\mathbb{T}^{-1}_{s-t}(y-\tilde{m}^{s,y}_{s,t}(x)))}{{\det \mathbb{T}_{s-t}}},
\end{split}
\]
with the following notation at hand:
\[\zeta^j_i \,:= \, -\frac{1+\alpha(i-1)}{\alpha}+\beta^j\frac{1+\alpha(j-1)}{\alpha}.\]
Then, we write with the notations of \eqref{proof:remainder_decomposition} that
\begin{align}
\notag
\left|\tilde{R}^1_\epsilon f(t,x)\right| \, &\le \, C\sum_{i=1}^n\sum_{j=i}^n \int_t^T \int_{\R^N} |f(s,y)|\frac{\overline{p}(1,\mathbb{T}^{-1}_{s-t}(y-\tilde{m}^{s,y}_{s,t}(x)))}{\det \mathbb{T}_{s-t}}\frac{\left|\mathbb{T}^{-1}_{s-t}(x-\theta_{t,s}(y))\right|^{\beta^j}}{(s-t)^{-\zeta^j_i}} \, dyds \\
 &=: \, C\sum_{i=1}^n\sum_{j=i}^n \int_t^T \int_{\R^N} |f(s,y)|\mathcal{I}_{ij}(t,s,x,y) \, dyds.\label{proof:remainder:decomposition_R1}
\end{align}
Then, from the H\"older inequality,
\begin{equation}
\label{proof:control_remainder}
\left|\tilde{R}^1_\epsilon f(t,x)\right| \, \le \, C\Vert f \Vert_{L^p_tL^q_x}\sum_{i=1}^n\sum_{j=i}^n \left(\int_t^T \left(\int_{\R^N} \left[\mathcal{I}_{ij}(t,s,x,y) \right]^{q'} dy\right)^{\frac{p'}{q'}}ds\right)^{\frac{1}{p'}},
\end{equation}
where $q'$ and $p'$ are the conjugate exponents of $q$ and $p$, respectively.\newline
Now, the integrals with respect to $y$ can be easily controlled by  Corollary \ref{coroll:control_power_q_density}. Indeed,
\begin{align}
\label{proof:remainder_control_Jij}
\int_{\R^N} \bigl[\mathcal{I}_{ij}(t,s,x,\,&y) \bigr]^{q'} dy  \\\notag
&\le \, C \left(\frac{(s-t)^{\zeta^j_i}}{\det \mathbb{T}_{s-t}}\right)^{q'} \int_{\R^N}
\left|\mathbb{T}^{-1}_{s-t}(x-\theta_{t,s}(y))\right|^{\beta^jq'}\overline{p}(1,\mathbb{T}^{-1}_{s-t}(y-\tilde{m}^{s,y}_{s,t}(x)))
\, dy.
\end{align}
Choosing $q_0>1$ big enough so that $\beta^j q'<\alpha$ for any $ j$ in $\llbracket 1,n\rrbracket$ and any $q\ge q_0$, we can use Corollary \ref{coroll:Smoothing_effect} to show that
\begin{equation}
\label{proof:remainder_control_Jij2}
\int_{\R^N} \left[\mathcal{I}_{ij}(t,s,x,y) \right]^{q'} dy\, \le \, C  (s-t)^{\zeta^j_iq'}\left(\det \mathbb{T}_{s-t}\right)^{1-q'}.
\end{equation}
Going back to Equation \eqref{proof:control_remainder}, we can thus write that
\[
\left|\tilde{R}^1_\epsilon f(t,x)\right| \, \le \, C\Vert f \Vert_{L^p_tL^q_x} \sum_{i=1}^n\sum_{j=i}^n\left(\int_t^T (s-t)^{\zeta^j_i p'}\left(\det \mathbb{T}_{s-t}\right)^{\frac{p'}{q'}-p'} ds\right)^{\frac{1}{p'}}.
\]
Noticing now that for any $i\le j$ in $\llbracket 1, n\rrbracket$
\begin{equation}
\label{proof:eq:control_zeta}
\zeta^j_i >-1 \, \Leftrightarrow \, -\frac{1+\alpha(i-1)}{\alpha}+\beta^j\frac{1+\alpha(j-1)}{\alpha} >-1 \, \Leftrightarrow \, \beta^j >\frac{1+\alpha(i-2)}{1+\alpha(j-1)},
\end{equation}
we can choose $q_0>1$, $p_0>1$ large enough so that $p'$, $q'$ are sufficiently close to $1$ in order to conclude that
\begin{equation}
\label{proof:remainder:final_control_R^1}
\left|\tilde{R}^1_\epsilon f(t,x)\right| \, \le \, C\Vert f \Vert_{L^p_tL^q_x}.
\end{equation}
We can now focus on the control for the first remainder term $\tilde{R}^0_\epsilon f$. Since clearly $\tilde{R}^0_\epsilon f=0$ if $\sigma(t,x)$ is constant in space, we can assume without loss of generality that $\nu$ is absolutely continuous with respect to the Lebesgue measure on $\R^d$ (cf. assumption [\textbf{AC}]). In particular, we know that it can be decomposed as in \eqref{eq:decomposition_nu_AC}:
\[\nu(dz) \, = \ Q(z)\frac{g(\frac{z}{|z|})}{|z|^{d+\alpha}}dz.\]
Given now a smooth enough function $\phi\colon \R^N\to \R$, we start noticing that
\[
\begin{split}
\mathcal{L}_t\phi(x) \, &= \, \int_{\R^d_0}\left[\phi(x+B\sigma(t,x)z)-\phi(x)\right]\, \nu(dz) \\
 &= \, \int_{\R^d_0}\left[\phi(x+B\sigma(t,x)z)-\phi(x) \right]\, Q(z)g\left(\frac{z}{|z|}\right)\frac{dz}{|z|^{d+\alpha}} \\
 &= \, \int_{\R^d_0}\left[\phi(x+B\tilde{z})-\phi(x) \right] Q(\sigma^{-1}(t,x)\tilde{z})\frac{g\left(\frac{\sigma^{-1}(t,x)\tilde{z}}{|\sigma^{-1}(t,x)\tilde{z}|}\right)}{\det \sigma(t,x)}\frac{d\tilde{z}}{|\sigma^{-1}(t,x)\tilde{z}|^{d+\alpha}},
\end{split}
\]
where we assumed, without loss of generality from \textbf{[UE]}, that $ \det \sigma(t,x)>0$. A similar representation holds for $\tilde{\mathcal{L}}^{s,y}_t\phi(x)$, too.
Now, let us introduce for any $z$ in $\R^d$, the following quantity:
\begin{multline*}
\tilde{H}^{s,y}_{t,x}(z) \\
:= \,  Q(\sigma^{-1}(t,x)z)\frac{g\left(\frac{\sigma^{-1}(t,x)z}{|\sigma^{-1}(t,x)z|}\right)}{\det \sigma(t,x)|\sigma^{-1}(t,x)\frac{z}{|z|}|^{d+\alpha}}- Q((\tilde{\sigma}^{s,y}_t)^{-1}z)\frac{g\left(\frac{(\tilde{\sigma}^{s,y}_t)^{-1}z}{|(\tilde{\sigma}^{s,y}_t)^{-1}z|}\right)}{\det \tilde{\sigma}^{s,y}_t|(\tilde{\sigma}^{s,y}_t)^{-1}\frac{z}{|z|}|^{d+\alpha}},
\end{multline*}
where we have normalized $z$ above in order to make the usual isotropic stable L\'evy measure appear. \newline 
Fixed $\eta>0$, local to this section, meant to be small and to be chosen later (and not to be confused with the ellipticity constant in assumption \textbf{[UE]}), we then define
\begin{equation}\label{DEF_ALPHA_ETA}
\alpha_\eta= \alpha/(1-\eta),
\end{equation}
and we decompose the integral in the difference of the generators in the following way:
\[
\begin{split}
 \left(\mathcal{L}_t-\tilde{\mathcal{L}}^{s,y}_t \right)\phi(x) \, &= \, \int_{\R^d_0}\left[\phi(x+Bz)-\phi(x)\right]\tilde{H}^{s,y}_{t,x}(z)\,s \frac{dz}{|z|^{d+\alpha}} \\
 &= \, \int_{\Delta_\eta}\left[\phi(x+Bz)-\phi(x) -\langle D_x\phi(x),Bz\rangle\right]\tilde{H}^{s,y}_{t,x}(z)\, \frac{dz}{|z|^{d+\alpha}} \\
 &\qquad \qquad + \int_{\Delta^c_\eta}\left[\phi(x+Bz)-\phi(x)\right]\tilde{H}^{s,y}_{t,x}(z)\, \frac{dz}{|z|^{d+\alpha}} \\
 &=:\, \sum_{i=1}^2 \left[\Delta_{i}\phi(t,s,\cdot,y)\right](x),
\end{split}
\]
where we have denoted, for simplicity,
\[
\begin{split}
\Delta_\eta \, &:= \, B(0,(s-t)^{\frac{1}{\alpha_\eta}});\\
\Delta^c_\eta \, &:= \, B^c(0,(s-t)^{\frac{1}{\alpha_\eta}}).
\end{split}
\]
We highlight in particular that it is precisely the symmetry of $\nu$ that ensures that the function $\tilde{H}^{s,y}_{t,x}$ is even and that allow us to introduce the odd first order term $\langle D_x\phi(x),Bz\rangle$ in the first integral above on the simmetric space $\Delta_\eta$.\newline
Noticing from Proposition \ref{prop:Decomposition_Process_X} that the frozen ``density'' $\tilde{p}^{s,y}$ is regular enough in $x$, we can now replace $\phi$ in the above decomposition with $\tilde{p}^{s,y}(t,s,\cdot,y)$. Going back to $\tilde{R}^0_\epsilon f$ given in \eqref{proof:remainder_decomposition}, we start rewriting it as
\begin{equation}
\label{proof:remainder:decomposition_R0}
\begin{split}
|\tilde{R}^0_\epsilon f(t,x)| \, &\le \, C\sum_{i=1}^2\int_t^T \int_{\R^N}|f(s,y)|\,\left|\left[\Delta_{i}\tilde{p}^{s,y}(t,s,\cdot,y)\right](x)\right| \, dy ds \\
&=: \, \sum_{i=1}^2\int_t^T \int_{\R^N}|f(s,y)|\mathcal{I}_{0i}(t,s,x,y) \, dy ds.
\end{split}
\end{equation}
As before, we can then apply H\"older inequality to show that
\begin{equation}
\label{proof:control_remainder1}
\left|\tilde{R}^0_\epsilon f(t,x)\right| \, \le \, C\Vert f \Vert_{L^p_tL^q_x}\sum_{i=1}^2 \left(\int_t^T \left(\int_{\R^N} \left[\mathcal{I}_{0i}(t,s,x,y) \right]^{q'} dy\right)^{\frac{p'}{q'}}ds\right)^{\frac{1}{p'}},
\end{equation}
where $q'$ and $p'$ are again the conjugate exponents of $q$ and $p$, respectively.\newline
To control the second term involving $\mathcal{I}_{02}$, we start noticing that
\begin{equation}
\label{proof:control_remainder_H}
    |\tilde{H}^{s,y}_{t,x}(z)| \, \le \, C
\end{equation}
for some constant $C$ independent from the parameters, thanks to assumption [\textbf{UE}] for $\sigma$ and the boundedness of $g$ and $Q$.\newline
Then, we can use Control \eqref{proof:control_remainder_H}, Corollary \ref{coroll:control_power_q_density} and the H\"older inequality to write that
\[
\begin{split}
  |\mathcal{I}_{02}&(t,s,x,y)|^{q'} \, \le \, C\Bigl[\int_{\Delta^c_\eta}|\tilde{p}^{s,y}(t,s,x+Bz,y)-\tilde{p}^{s,y}(t,s,x,y)|\frac{dz}{|z|^{d+\alpha}}\Bigr]^{q'} \\
&\le \, C\left(\int_{\Delta^c_\eta} \frac{dz}{|z|^{d+\alpha}}\right)^{\frac{q'}{q}}\int_{\Delta^c_\eta}|\tilde{p}^{s,y}(t,s,x+Bz,y)-\tilde{p}^{s,y}(t,s,x,y)|^{q'}\frac{dz}{|z|^{d+\alpha}}\\
&\le \, C\frac{(s-t)^{ (\eta-1)\frac{q'}{q}}}{(\det \mathbb{T}_{s-t})^{q'}}\int_{\Delta^c_\eta}\left[\overline{p}(1,\mathbb{T}^{-1}_{s-t}(y-\theta_{s,t}(x+Bz)))+\overline{p}(1,\mathbb{T}^{-1}_{s-t}(y-\theta_{s,t}(x)))\right]\frac{dz}{|z|^{d+\alpha}},
\end{split}
\]
recalling from \eqref{DEF_ALPHA_ETA} that $\alpha_\eta=\alpha/(1-\eta) $ for the last inequality.
The Fubini theorem and the change of variables $\tilde{y}=y-\theta_{s,t}(x+Bz)$ now show that
\begin{equation}
\label{proof:remainder:control_I0j}
\begin{split}
\int_{\R^N} |\mathcal{I}_{02}(t,s,x,y)|^{q'} \, dy
\, &\le \, 2C\frac{(s-t)^{ (\eta-1)\frac{q'}{q}}}{(\det \mathbb{T}_{s-t})^{q'-1}}\int_{B^c(0,(s-t)^{\frac{1}{\alpha_\eta}})}\int_{\R^N}\overline{p}(1,\tilde{y})\, d\tilde{y}\frac{dz}{|z|^{d+\alpha}}
\\
&\le \, C(\det \mathbb{T}_{s-t})^{1-q'}(s-t)^{  (\eta-1)\frac{q'}{q}}\int_{\Delta^c_\eta}\frac{dz}{|z|^{d+\alpha}}
\\
&\le \, C(\det \mathbb{T}_{s-t})^{1-q'}(s-t)^{q' (\eta-1)}.
\end{split}
\end{equation}
Going back to Equation \eqref{proof:control_remainder1}, we can then conclude that
\[
\begin{split}
\int_t^T\left(\int_{\R^N} |\mathcal{I}_{02}(t,s,x,y)|^{q'} \, dy\right)^{\frac{p'}{q'}} ds \,
&\le \, C \int_t^T (\det \mathbb{T}_{s-t})^{-\frac{p'}{q}}(s-t)^{p' (\eta-1)} \, ds  \\
&\le \, C \int_t^T (s-t)^{-p'( 1-\eta+\frac{1}{q}\sum_{i=1}^nd_i\frac{1+\alpha(i-1)}{\alpha})} \, ds,
\end{split}
\]
where in the last step we also exploited \eqref{eq:control_det_T1}.\newline
Assuming now that $\eta<1$ and $p$, $q$ are big enough so that
\[p'( 1-\eta+\frac{1}{q}\sum_{i=1}^nd_i\frac{1+\alpha(i-1)}{\alpha}) \, < \,1,\]
we immediately obtain that
\begin{equation}
\label{proof:control_for_I_02}
\left(\int_t^T\left(\int_{\R^N} |\mathcal{I}_{02}(t,s,x,y)|^{q'} \, dy\right)^{\frac{p'}{q'}} ds\right)^{\frac{1}{p'}} \,
\le \, C_T.
\end{equation}
We can now focus on the integral with respect to $y$ of the first term $\mathcal{I}_{01}$ in Equation \eqref{proof:control_remainder1}. Using the Lipschitz continuity of $Q$ in a neighborhood of zero and the H\"older regularity of the diffusion matrix $\sigma$, it is not difficult to check that
\[ |\tilde{H}^{s,y}_{t,x}(z)| \, \le \, C\sum_{j=1}^n\left|(x-\theta_{t,s}(y))_j\right|^{ \beta^1}.\]
Thanks to the above estimate, we exploit a Taylor expansion on the density $\tilde{p}^{s,y}$ and Proposition \ref{prop:Smoothing_effect} with $k=2$ and $(\tau,\xi)=(s,y)$ to show that
\begin{align}
\Theta(t,s,x,y,z)\, &:= \,  \Bigl|\bigl[\tilde{p}^{s,y}(t,s,x+ Bz,y)-\tilde{p}^{s,y}(t,s,x,y)
-\langle D_x\tilde{p}^{s,y}(t,s,x,y),Bz\rangle\bigr]\tilde{H}^{s
,y}_{t,x}(z) \Bigr|\notag \\
&\le \, \sum_{j=1}^n\int_0^1\left|(x-\theta_{t,s}(y))_j\right|^{ \beta^1}\,  |D^2_{x_1}\tilde{p}^{s,y}(t,s,x+\lambda Bz,y)| |z|^2 \,
d\lambda \notag\\
&\le \, \frac{C}{\det
\mathbb{T}_{s-t}}\int_0^1|z|^2\frac{\overline{p}(1,\mathbb{T}^{-1}_{s-t}(y-\tilde{m}^{s,y}_{s,t}(x+\lambda Bz)))}{(s-t)^{\frac{2}{\alpha}}}\notag\\
&\qquad\qquad\qquad\qquad\times\Bigl[\sum_{j=1}^n\left|(x+\lambda Bz-\theta_{t,s}(y))_j\right|^{\beta^1}+|\lambda B z|^{\beta^1}\Bigr] \, d\lambda \notag \\
\label{proof:remainder:control_H2}
&\le \, \frac{C}{\det
\mathbb{T}_{s-t}}\int_0^1|z|^2\overline{p}(1,\mathbb{T}^{-1}_{s-t}(y-\tilde{m}^{s,y}_{s,t}(x+\lambda Bz)))\\
&\qquad\qquad\qquad\qquad\times\Bigl[\sum_{j=1}^n\frac{\left|\mathbb{T}^{-1
}_{s-t}(x+\lambda Bz-\theta_{t,s}(y))\right|^{\beta^1}}{(s-t)^{{\zeta^j_0}}}+\frac{|z_1|^{\beta^1}}{(s-t)^{\frac{2}{\alpha}}}\Bigr]
 \, d\lambda, \notag
\end{align}
where, similarly to above, we have denoted:
\begin{equation}
\label{proof:def_eta_0}
\zeta^j_0 \,:= \, \frac{2}{\alpha}-\beta^1\frac{1+\alpha(j-1)}{\alpha}.
\end{equation}
It then follows from the H\"older inequality and Corollary \ref{coroll:control_power_q_density} that
\[
\begin{split}
|\mathcal{I}_{01}&(t,s,x,y)|^{q'} \,\le \, \left[\int_{\Delta_\eta}\Theta(t,s, x,y,z)\,\frac{dz}{|z|^{d+\alpha}}\right]^{q'} \\
&\le \, \frac{C}{(\det
\mathbb{T}_{s-t})^{q'}}\left(\int_{\Delta_\eta}1\, dz\right)^{\frac{q'}{q}}\int_0^1\int_{\Delta_\eta}\overline{p}(1,\mathbb{T}^{-1}_{s-t}(y-\tilde{m}^{s,y}_{s,t}(x+\lambda Bz)))\\
&\qquad\qquad\quad\quad \qquad \times \Bigl[\sum_{j=1}^n\frac{\left|\mathbb{T}^{-1
}_{s-t}(x+\lambda Bz-\theta_{t,s}(y))\right|^{\beta^1}}{(s-t)^{{\zeta^j_0}}}+\frac{|z_1|^{\beta^1}}{(s-t)^{\frac{2}{\alpha}}}\Bigr]^{q'}
 \, \frac{dz}{|z|^{q'(d+\alpha-2)}}d\lambda.
\end{split}
\]
If we now add the integral with respect to $y$, Fubini Theorem readily implies that
\[
\begin{split}
\int_{\R^N}|\mathcal{I}_{01}&(t,s,x,y)|^{q'} \, dy \\
&\le \,C\frac{(s-t)^{\frac{d}{\alpha_\eta}(q'-1)}}{(\det
\mathbb{T}_{s-t})^{q'}} \int_0^1\int_{\Delta_\eta}\int_{\R^N}\overline{p}(1,\mathbb{T}^{-1}_{s-t}(y-\tilde{m}^{s,y}_{s,t}(x+\lambda Bz)))\\
&\qquad \qquad \times \Bigl[\sum_{j=1}^n\frac{\left|\mathbb{T}^{-1
}_{s-t}(x+\lambda Bz-\theta_{t,s}(y))\right|^{\beta^1q'}}{(s-t)^{{\zeta^j_0q'}}}+\frac{|z_1|^{\beta^1q'}}{(s-t)^{\frac{2}{\alpha}q'}}\Bigr]
 \, dy\frac{dz}{|z|^{q'(d+\alpha-2)}}d\lambda.
\end{split}
\]
If we assume to have taken $q'$ close enough to $1$ so that $\beta^1q'<\alpha$, we can use Corollary \ref{coroll:Smoothing_effect} to show that
\[
\begin{split}
\int_{\R^N}|\mathcal{I}_{01}(t,s,\,&x,y)|^{q'} \, dy \\
&\le \, C\frac{(s-t)^{\frac{d}{\alpha_\eta}(q'-1)}}{(\det
\mathbb{T}_{s-t})^{q'-1}}\int_{B(0,(s-t)^{\frac{1}{\alpha_\eta}})}\left[\sum_{j=1}^n\frac{1}{(s-t)^{{q'\zeta^j_0}}}+\frac{|z_1|^{q'\beta^1}}{(s-t)^{q'\frac{2}{\alpha}}}\right]
 \,\frac{dz}{|z|^{q'(d+\alpha-2)}}\\
&\le \, C\frac{(s-t)^{\frac{d}{\alpha_\eta}(q'-1)}}{(\det\mathbb{T}_{s-t})^{q'-1}}\int_0^{(s-t)^{\frac{1}{\alpha_\eta}}}\left[\sum_{j=1}^n\frac{r^{d-1-(d+\alpha-2)q'}}{(s-t)^{{q'\zeta^j_0}}}+\frac{r^{d-1-(d+\alpha-2-\beta^1)q'}}{(s-t)^{\frac{2}{\alpha}q'}}\right]
 \,dr.
\end{split}
\]
Similarly, if $q$ is big enough (so that $q'$ is close to $1$), it holds that
\[d-1-q'(d+\alpha-2) > -1 \, \Leftrightarrow \, q' < \frac{d}{d+\alpha-2}\]
 and we can integrate with respect to $r$:
\begin{align}\notag
\int_{\R^N}|\mathcal{I}_{01}(t,\,&s,x,y)|^{q'} \, dy \,
\le \, C\frac{(s-t)^{\frac{d}{\alpha_\eta}(q'-1)}}{(\det\mathbb{T}_{s-t})^{q'-1}}\left[\sum_{j=1}^n\frac{r^{d-q'(d+\alpha-2)}}{(s-t)^{{q'\zeta^j_0}}}+\frac{r^{d-q'(d+\alpha-2-\beta^1)}}{(s-t)^{q'\frac{2}{\alpha}}}\right]\Bigg{|}_0^{(s-t)^{\frac{1}{\alpha_\eta}}} \\
&\le \, C (\det\mathbb{T}_{s-t})^{1-q'}(s-t)^{\frac{q'}{\alpha_\eta}(2-\alpha)}\left[\sum_{j=1}^n(s-t)^{-q'\zeta^j_0}+(s-t)^{q'(\frac{\beta^1}{\alpha_\eta}-\frac{2}{\alpha})}\right].
\label{proof:remainder:control_I0j2}
\end{align}
Hence, it follows from Equation \eqref{eq:control_det_T1} that
\[
\begin{split}
\int_t^T\Bigl(\int_{\R^N}|&\mathcal{I}_{01}(t,s,x,y)|^{q'} \, dy\Bigr)^{\frac{p'}{q'}}ds \\
&\le \, C \int_t^T(\det\mathbb{T}_{s-t})^{-\frac{p'}{q}}(s-t)^{\frac{p'}{\alpha_\eta}(2-\alpha)}\left[\sum_{j=1}^n(s-t)^{-p'\zeta^j_0}+(s-t)^{p'(\frac{\beta^1}{\alpha_\eta}-\frac{2}{\alpha})}\right]\, ds\\
&\le \, C \int_t^T(s-t)^{p'\left( \frac{(2-\alpha)}{\alpha_\eta}-\frac{1}{q}\sum_{i=1}^nd_i\frac{1+\alpha(i-1)}{\alpha}\right)}\left[\sum_{j=1}^n(s-t)^{-p'\zeta^j_0}+(s-t)^{p'(\frac{\beta^1}{\alpha_\eta}-\frac{2}{\alpha})}\right]\, ds
\end{split}
\]
To conclude, we need to show that the two terms above are integrable with respect to $s$. Namely,
\begin{align*}
p'\left(\frac{(2-\alpha)}{\alpha_\eta}-\frac{1}{q}\sum_{i=1}^nd_i\frac{1+\alpha(i-1)}{\alpha}-\zeta^j_0\right) \, &> \, -1, \quad \forall\, j \in \llbracket 1,n \rrbracket;\\
p'\left(\frac{(2-\alpha)}{\alpha_\eta}-\frac{1}{q}\sum_{i=1}^nd_i\frac{1+\alpha(i-1)}{\alpha}+\frac{\beta^1}{\alpha_\eta}-\frac{2}{\alpha}\right) \, &> \, -1.
\end{align*}
Recalling again that we can choose $p$, $q$ big enough as we want, so that Equation \eqref{eq:control_det_T2} holds, it is now sufficient to take $\eta$ in $(0,1)$ in order to have:
\begin{align}
\frac{(2-\alpha)}{\alpha_\eta}-\zeta^j_0 \, = \, \frac{(2-\alpha)}{\alpha_\eta}-\frac{2}{\alpha}+\beta^1\frac{1+\alpha(j-1)}{\alpha} \, &> \, -1,\quad \forall\, j \in \llbracket 1,n \rrbracket; \label{proof:condition1}\\
\frac{(2-\alpha)}{\alpha_\eta}+\frac{\beta^1}{\alpha_\eta}-\frac{2}{\alpha} \, &> \, -1.\label{proof:condition2}
\end{align}
By direct calculations, recalling from \eqref{DEF_ALPHA_ETA} that $\alpha_\eta=\alpha/(1-\eta) $, we now notice that Conditions \eqref{proof:condition1}-\eqref{proof:condition2} can be rewritten as follows
\begin{align*}
\eta \, &< \, \frac{\beta^1(1+\alpha(j-1)}{2-\alpha}, \quad \forall\, j \in \llbracket 1,n \rrbracket;\\
\eta \, &< \, \frac{\beta^1}{2+\beta^1-\alpha}.
\end{align*}
Choosing $\epsilon>0$ so that the above conditions holds, we have that
\begin{equation}
\label{proof:control_for_I_01}
\left(\int_t^T\Bigl(\int_{\R^N}|\mathcal{I}_{01}(t,s,x,y)|^{q'} \, dy\Bigr)^{\frac{p'}{q'}}ds\right)^{\frac{1}{p'}} \, \le \, C_T.
\end{equation}
Going back to Equation \eqref{proof:control_remainder1}, we   use Controls \eqref{proof:control_for_I_02}-\eqref{proof:control_for_I_01} to write that
\begin{equation}
\label{proof:remainder:final_control_R^0}
\left|\tilde{R}^0_\epsilon f(t,x)\right| \, \le \, C\Vert f \Vert_{L^p_tL^q_x}.
\end{equation}
Exploiting Controls \eqref{proof:remainder:final_control_R^1} and \eqref{proof:remainder:final_control_R^0} in Equation \eqref{proof:remainder_decomposition}, we have concluded our proof.
\end{proof}

A similar control in $L^p_tL^q_x$-norms can be obtained.
In particular, we point out that Equation \eqref{eq:control_LpLq_R} below implies  that the operator $I-\tilde{R}_\epsilon$ is invertible in $L^p\left(0,T;L^q(\R^N)\right)$, provided $T$ is small enough. From Lemma \ref{prop:convergence_LpLq}, the same holds for $I_\epsilon-\tilde{R}_\epsilon$.

\begin{prop}
\label{prop:control_tildeR_LpLq}
Let $q>1$, $p>1 $ be such that Condition ($\mathscr{C}$) holds. Then, there exists $C:=C(T,p,q)>0$ such that for any $f$ in $L^p\left(0,T;L^q(\R^N)\right)$,
\begin{equation}
\label{eq:control_LpLq_R}
\Vert \tilde{R}_\epsilon f\Vert_{L^p_tL^q_x}\, \le \, C \Vert f \Vert_{L^p_tL^q_x}.
\end{equation}
In particular, it holds that $\lim_{T\to 0}C(T,p,q)=0$.
\end{prop}
\begin{proof}
We are going to keep the same notations used in the previous proof. In particular, we recall the following decomposition
\[\tilde{R}_\epsilon f(t,x) \, = \, \tilde{R}^0_\epsilon f(t,x)+\tilde{R}^1_\epsilon f(t,x), \]
given in Equation \eqref{proof:remainder_decomposition}. \newline
In order to control the second term $\tilde{R}^1_\epsilon f$ in $L^p_tL^q_x$-norm, we start from Equation \eqref{proof:remainder:decomposition_R1} to write that
\[
    \Vert \tilde{R}^1_\epsilon f(t,\cdot) \Vert_{L^q_x} \, \le \, C\sum_{i=1}^n\sum_{j=i}^n \int_t^T \Bigl{\Vert}\int_{\R^N} |f(s,y)|\mathcal{I}_{ij}(t,s,\cdot,y) \, dy\Bigr{\Vert}_{L^q_x}ds.
\]
The Young inequality now implies that
\begin{align*}
\Bigl{\Vert}\int_{\R^N} |f(s,y)|\mathcal{I}_{ij}(t,s,\cdot,y) \, dy\Bigr{\Vert}_{L^q_x}^q&\, = \, \int_{\R^N} \Big|\int_{\R^N} |f(s,y)| \mathcal{I}_{ij}(t,s,x,y)\, dy\Big|^q \, dx\\
&\le\, \int_{\R^n} 
\!\left(\int_{\R^N} |f(s,y)|^q \mathcal{I}_{ij}(t,s,x,y)\, dy\right) \|\mathcal{I}_{ij}(t,s,x,\cdot)\|_{L^1}^{q/q'}dx\\
&\le\,  C(s-t)^{\zeta_i^j q/q'} \int_{\R^N} |f(s,y)|^q \left(\int_{\R^N} \mathcal{I}_{ij}(t,s,x,y) \,dx\right)dy
\end{align*}
using Control \eqref{proof:remainder_control_Jij2} and the Fubini Theorem for the last inequality. From \eqref{proof:remainder:decomposition_R1}, \eqref{proof:remainder_control_Jij} and the correspondence \eqref{eq:identification_theta_m1} which gives $y-\tilde m_{s,t}^{s,y}(x)=\theta_{t,s}(y)-x $ it is plain to derive that:
$$\int_{\R^N}dx \mathcal{I}_{ij}(t,s,x,y)\le C (s-t)^{\zeta_i^j }.$$
Thus,
\[\Vert \tilde{R}^1_\epsilon f(t,\cdot) \Vert_{L^q_x} \, 
\le \, C\sum_{i=1}^n\sum_{j=i}^n\int_t^T (s-t)^{\zeta^j_i}  \Vert f(t,\cdot) \Vert_{L^q_x} \, ds,\]
where, in the last step, we exploited Equation \eqref{proof:remainder_control_Jij2} with $q'=1$, recalling that $\beta^j<1<\alpha$.\newline
We can then use the above control to write that
\[
\begin{split}
\Vert \tilde{R}^1_\epsilon f\Vert^p_{L^p_tL^q_x} \, &\le \, \sum_{i=1}^n\sum_{j=i}^n \int_0^T \Vert \tilde{R}^1_\epsilon f(t,\cdot) \Vert_{L^q_x}^p \, dt \\
&\le \,C \sum_{i=1}^n\sum_{j=i}^n\int_0^T \Vert f(t,\cdot) \Vert^{p}_{L^q_x} \left( \int_t^T (s-t)^{\zeta^j_i} \, ds  \right)^p dt \\
&\le \, C_T \sum_{i=1}^n\sum_{j=i}^n\int_0^T \Vert f(t,\cdot) \Vert^p_{L^q_x} dt\le C_T\|f\|_{L^p_tL^q_x}^p,
\end{split}
\]
where $C_T:=C(p',q',T)$ denotes a positive constant that tends to zero if $T$ goes to zero (recall indeed from \eqref{proof:eq:control_zeta} that $\zeta^j_i>-1$)
. \newline
The control for $\tilde{R}^0_\epsilon f$ can be obtained following the same arguments above, exploiting Equation \eqref{proof:remainder:decomposition_R0} instead of \eqref{proof:remainder:decomposition_R1} and Equations \eqref{proof:remainder:control_I0j}-\eqref{proof:remainder:control_I0j2} with $q'=1$ for the controls of $\Vert \mathcal{I}_{0j}(t,s,x,\cdot) \Vert_{L^1_x}$.
\end{proof}

Let us fix now a function $f$ in $ C_c^{1,2}([0,T)\times \R^{N})$. The first step of our method consists in applying the It\^o formula on the Green kernel $\tilde{G}_\epsilon f$ and the process $\{X^{t,x}_s\}_{s\in[t,T]}$, solution of the martingale problem with starting point $(t,x)$:
\[
\mathbb{E}\left[\tilde{G}_\epsilon f(t,x)+ \int_t^{T} (\partial_s+L_s)\tilde{G}_\epsilon f(s,X^{t,x}_s)   ds\right] \, =\, 0.
\]
We then exploit Equation \eqref{eq:differential_eq2} to write that
\[
\tilde{G}_\epsilon f (t,x)- \mathbb{E}\left[\int_t^{T} I_\epsilon f(s,X^{t,x}_s)\,  ds \right]+\mathbb{E}\left[\int_{t}^T\left[L_s\tilde{G}_\epsilon f- \tilde{M}_{\epsilon}f\right] (s,X^{t,x}_s) \, ds\right] \, = \, 0.
\]
Thus, it holds that
\begin{equation}
\label{eq:NUMBER_PREAL_REG_DENS}
\mathbb{E}\Bigl[\int_t^{T} I_\epsilon f(s,X^{t,x}_s) \, ds \Bigr] \,
= \, \tilde{G}_\epsilon f
(t,x)+\mathbb{E}\Bigl[\int_t^T\tilde{R}_\epsilon f(s,X^{t,x}_s)\, ds\Bigr].
\end{equation}

Thanks to Proposition \ref{prop:pointwise_green_kernel}, we know that there exists $C(T):=C(T) \underset{T\rightarrow 0}{\longrightarrow} 0$ such that
\begin{equation}
\label{eq:control_infty_Green}
\Vert \tilde{G}_\epsilon f\Vert_\infty \, \le \,  C \Vert f\Vert_{L^p_tL^q_x}.
\end{equation}

Let us assume for now that $p,q$ are large enough so that the control \eqref{eq:control_infty_R} of Lemma \ref{prop:control_tildeR} (pointwise control of the remainder) holds. From \textcolor{black}{Equations} \eqref{eq:NUMBER_PREAL_REG_DENS}, \eqref{eq:control_infty_Green} and \eqref{eq:control_infty_R}, we readily get that
\[
\left|\mathbb{E}\left[\int_t^T I_\epsilon f(X_s^{t,x})\,  ds \right]\right|\, \le \,  C\Vert f\Vert_{L^p_tL^q_x}.
\]

Letting $\epsilon $ go to zero,  we  thus derive that any solution $\{X^{t,x}_s\}_{s\in [t,T]}$ of the martingale problem for $\partial_s+L_s$ with initial condition $(t,x)$ satisfies
\[
\left|\mathbb{E}\left[\int_{t}^T f(s,X_s^{t,x})\, ds\right]\right|\,\le \, C\Vert f\Vert_{L^p_tL^q_x},\]
for any $f$ in  $f$ in $ C_c^{1,2}([0,T)\times \R^{N})$. \textcolor{black}{Above, we have exploited Lemma \ref{convergence_dirac} for the integral in space and the bounded convergence Theorem  for that in time.}\newline
To show the result for a general $f$ in $L^p\left(0,T;L^q(\R^N)\right)$, we now use a density argument and the Fatou Lemma. Indeed, let $\{f_n\}_{n\in \N}$ a sequence of functions in $C_c^{1,2}([0,T)\times \R^{N})$ such that $\Vert f_n -f \Vert_{L^p_tL^q_x}\to 0$.
We then have that:
\begin{equation}\label{eq:Krilov_partial_estimates}
\begin{split}
\left|\mathbb{E}\left[\int_{t}^T f(s,X_s^{t,x})\, ds\right]\right|\,&\le \, \left|\mathbb{E}\left[\int_{t}^T \liminf_n f_n(s,X_s^{t,x})\, ds\right]\right| \\
&\le \, \liminf_n\left|\mathbb{E}\left[\int_{t}^T  f_n(s,X_s^{t,x})\, ds\right]\right| \\
&\le \, C\liminf_n \Vert f_n\Vert_{L^p_tL^q_x} \\
&=  \,C \Vert f\Vert_{L^p_tL^q_x}.
\end{split}
\end{equation}
This is precisely the Estimate \eqref{eq:Krylov_Estimates} in Corollary \ref{coroll:Krylov_Estimates}, provided that $p,q$ are large enough.

Thanks to Estimates \eqref{eq:Krilov_partial_estimates}, we then know that the process $\{X^{t,x}_s\}_{s\in[t,T]}$ has a density we will denote by $p(t,s,x,y)$. From Equation
\eqref{eq:NUMBER_PREAL_REG_DENS} it now follows that
\begin{align}
\notag
\tilde{G}_\epsilon f(t,x) \, &= \, \mathbb{E}\left[\int_t^{T} I_\epsilon f(s,X_s^{t,x}) \, ds\right]-\mathbb{E}\left[\int_{t}^T \tilde{R}_\epsilon f(s,X_s^{t,x}) \, ds\right]\\\notag
&= \,  \int_t^{T} \int_{\R^N}I_\epsilon f (s,y) p(t,s,x,y) \, dy  ds- \int_{t}^T\int_{\R^N}\tilde{R}_\epsilon f(s,y)
p(t,s,x,y) \, dy  ds\\
&= \, \int_t^{T} \int_{\R^N}(I_\epsilon-\tilde{R}_\epsilon)f(s,y) p(t,s,x,y) \, dy  ds.
\label{REP_PRESQUE_FINI}
\end{align}
Then, Proposition \ref{prop:pointwise_green_kernel}, Lemma \ref{prop:convergence_LpLq} (with an additional approximation argument) and Control \eqref{eq:control_LpLq_R} imply that both sides of the above control are bounded in the $L^p_tL^q_x$-norm, uniformly in $\epsilon>0$. Thus, we can conclude that Equation \eqref{REP_PRESQUE_FINI} holds for any $f$ in $L^p\left(0,T;L^q(\R^N)\right)$. We then conclude from Lemma \ref{prop:control_tildeR} (pointwise control of the remainder) that letting $\epsilon$ go to zero, it holds that
\[
\mathbb{E}\left[\int_t^{T} f(s,X_s^{t,x}) \, ds\right]\, = \, \tilde{G}\circ (I-\tilde{R})^{-1}f(t,x),
\]
which gives uniqueness if the final time $T$ is small enough. Global well-posedness is again derived from a chaining argument in time.

To complete the proof of Corollary \ref{coroll:Krylov_Estimates}, it remains to derive the Krylov estimates \eqref{eq:Krylov_Estimates} under Condition ($\mathscr{C}$) and not only for $p,q$ large enough.\newline
Fixed a parameter $\delta>0$ meant to be small, we consider a ``mollified'' version of the solution process $X^{t,x}_s$, given by
\begin{equation}
\label{eq:decomposition_final_Krylov}
\overline{X}^{t,x,\delta}_s\, := \, X^{t,x}_s+\delta\mathbb{M}_{s-t} \overline{Z}_{s-t},
\end{equation}
where $\{\overline{Z}_s\}_{s\ge0}$ is an isotropic $\alpha$-stable process on $\R^N$. \newline
Let us denote now by $p^\delta(t,s,x,\cdot)$ the density associated with the random variable $\overline{X}^{t,s,\delta}_{s}$.
We notice that Equation \eqref{eq:decomposition_final_Krylov} implies in particular that
\[p^{\delta}(t,s,x,y) \, = \, \left[p(t,s,x,\cdot)\ast q^\delta(s-t,\cdot)\right](y),\]
where $q^\delta(t,\cdot)$ is the density of the process $\delta \mathbb{M}_t\overline{Z}_t$ and thus, under the integrability condition $(\mathscr{C})$ and thanks to the Young inequality, the quantity $\Vert p^\delta\Vert_{L^{p'}_tL^{q'}_x}$, where $p',q'$ are the conjugate exponents of $p,q$, respectively, is finite (possibly explosive with $\delta $).
The point is now to reproduce the previous perturbative analysis in order to prove that the controls on $\Vert p^\delta\Vert_{L^{p'}_tL^{q'}_x}$ actually do no depend on $\delta$.

For this reason, we introduce the mollified ``frozen'' process $\tilde{X}^{s,y,t,x,\delta}_{s}$ along the flow $\theta_{t,s}(y)$ as
\begin{equation}
\label{eq:decomposition_delta_frozen}
\tilde{X}^{s,y,t,x,\delta}_{s}\, := \, \tilde{X}^{s,y,t,x}_{s}+\delta\mathbb{M}_{s-t} \overline{Z}_{s-t}.
\end{equation}

Following the same arguments presented in Propositions \ref{prop:Decomposition_Process_X} and \ref{prop:Smoothing_effect}, it is now possible to show that the process $\tilde{X}^{s,y,t,x,\delta}_{s}$ admits a density $\tilde{p}^{s,y,\delta}(t,s,x,y)$ and that it  enjoys a  multi-scale bound similar to \eqref{eq:smoothing_effect_frozen_y}. Namely,
\begin{prop}
\label{coroll:Smoothing_effect_delta}
There exists a positive constant $C:=C(N,\alpha)$ such that for any $k$ in $\llbracket 0,2 \rrbracket$, any $i$ in $\llbracket 1,n \rrbracket$, any $t< s$ in $[0,T]$ and any $x,y$ in $\R^N$,
  \begin{equation}
  \label{eq:smoothing_effect_frozen_delta}
      \vert D^k_{x_i} \tilde{p}^{s,y,\delta}(t,s,x,y) \vert \, \le \,C\frac{\left((s-t)(1+\delta)\right)^{-k\frac{1+\alpha(i-1)}{\alpha}}}{\det \mathbb{T}_{(s-t)(1+\delta)}}\overline{p}\left(1,\mathbb{T}^{-1}_{(s-t)(1+\delta)}(y-\theta_{s,t}(x))\right).
  \end{equation}
\end{prop}
A sketch of proof for the above Proposition has been briefly presented in the Appendix section. \textcolor{black}{Importantly, we highlight that the constant $C$ appearing in \eqref{eq:smoothing_effect_frozen_delta} is independent from the ``smoothing'' parameter $\delta$}.\newline
Then, the same arguments leading to \eqref{eq:NUMBER_PREAL_REG_DENS} can be applied here to show that
\begin{equation}
\label{eq:NUMBER_PREAL_REG_DENS1}
\mathbb E\Bigl[\int_t^{T} I_\epsilon f(s,\overline{X}^{t,x,\delta}_s)\, ds\Bigr]  \, = \, \tilde{G}^\delta_{\epsilon} f
(t,x)+\mathbb{E}\Bigl[\int_{t}^T\int_{\R^N}\bigl(L^\delta_t-\tilde{L}^{s,y,\delta}_t\bigr) \tilde{G}^\delta_{\epsilon} f(s,\overline{X}^{t,x,\delta}_s)\, ds\Bigr],
\end{equation}
where $\tilde{G}^\delta_{\epsilon}$ and $\tilde{\mathcal{L}}^{s,y,\delta}$ are the frozen Green kernel and the frozen infinitesimal generator associated with the process $\tilde{X}^{s,y,t,x,\delta}_{s}$, respectively (cf.\ Equations \eqref{eq:def_Green_Kernel} and  \eqref{eq:def_frozen_generator}).
In particular, we point out that the pointwise bound \eqref{eq:control_infty_Green} on the Green kernel and the controls of Proposition \ref{prop:control_tildeR} (pointwise control of the remainder) are uniform with respect to the additional parameter $\delta$, thanks to Proposition \ref{coroll:Smoothing_effect_delta}.\newline
From Equation \eqref{eq:NUMBER_PREAL_REG_DENS1} and Proposition \ref{prop:control_tildeR_LpLq} ($L_t^pL_x^q $ control of the remainder)
we can then deduce that
\[
\left|\int_{t}^{T}\int_{\R^N} I_\epsilon f(s,y) p^\delta(t,s,x,y)\,  dy ds\right| \, \le \, C_T \left(1+\Vert p^\delta\Vert_{L^{p'}_tL^{q'}_x}\right)\Vert f \Vert_{L^p_tL^q_x}.
\]
\newline
From the Riesz representation theorem and the above inequality, we then deduce that  $\Vert p^\delta\Vert_{L^{p'}_tL^{q'}_x}\le C_T$, for $T$ small enough and uniformly in $\delta$. Hence,
\[\left|\int_t^{T} \int_{\R^N} I_\epsilon f (s,y) p^\delta(t,s,x,y)\,   dy ds\right| \, = \, \left|\int_t^{T}\mathbb{E} [I_\epsilon f (s,X_s^{t,x}+\delta\overline{Z}_s)]\,  ds\right| \le C_T\Vert I_\epsilon f\Vert_{L^p_tL^q_x}.\]
The Krylov-type estimate \eqref{eq:Krylov_Estimates} can be then derived exploiting the dominated convergence theorem and Lemma \ref{convergence_dirac} (Dirac Convergence of frozen density), letting firstly $\epsilon$ and then $\delta $ go to zero. We have thus concluded the proof of Corollary \ref{coroll:Krylov_Estimates}.

\setcounter{equation}{0}
\section{A counter-example to uniqueness}
\fancyhead[RO]{Section \thesection. A counter-example to uniqueness}
In this section, we present a counter-example to the uniqueness in law for the equation \eqref{eq:SDE} when the H\"older regularity in
space of the coefficients is low enough. In particular, we show here the almost sharpness of the thresholds appearing in Theorem
\ref{thm:main_result} for diagonal perturbations, proving also Theorem \ref{thm:counterexample}. In order to test the threshold associated with the critical
H\"older exponent for the $i$-th component of the drift $F$ with respect to the variables $x_j$, we adapt the \emph{ad hoc} Peano
example constructed in \cite{Chaudru:Menozzi17} to our L\'evy framework.\newline
Let us briefly recall it. It is well-known that the following deterministic equation
\begin{equation}
\begin{cases}
\label{eq:deterministic_Peano_example}
dy_t \, = \, \text{sgn}(y_t)\vert y_t\vert^\beta dt, \quad t\ge 0, \\
y_0 \, = \, 0,
\end{cases}
\end{equation}
for some $\beta$ in $(0,1)$, is ill-posed since it admits an infinite number of solutions of the form
\[y_t \, = \, \pm c(t-t_0)^{1/(1-\beta)}\mathds{1}_{[t_0,\infty)}(t), \quad \text{for some $t_0$ in }[0,+\infty).\]
Nevertheless, Bafico and Baldi in \cite{Bafico:Baldi81} proved that the associated SDE, obtained by adding a Brownian Motion
$\{W_t\}_{t\ge 0}$ to the dynamics:
\[\begin{cases}
dX_t \, = \, \text{sgn}(X_t)\vert X_t\vert^\beta dt+\epsilon dW_t, \quad t\ge 0 \\
X_0 \, = \, 0,
\end{cases}\]
is well-posed for any $\epsilon>0$ in a strong (probabilistic) sense. Furthermore, they showed that, letting $\epsilon$ goes to zero, the limit law
concentrates around the two extremal solutions $\pm ct^{1/(1-\beta)}$ of the deterministic equation \eqref{eq:deterministic_Peano_example}, thus providing a \textcolor{black}{selection}
``criterion'' between the infinite deterministic solutions. \newline
In a subsequent article \cite{Delarue:Flandoli14}, Delarue and Flandoli highlighted the hidden dynamical mechanism behind this counter-intuitive behaviour.
Heuristically, this \emph{regularization by noise} happens since, at least in a small time interval, the mean fluctuations of the Brownian noise are stronger than
the irregularity of the deterministic drift. Indeed, they showed that before some transition time $t_\epsilon$, the dominating noise pushes the solution to leave
the drift singularity at $0$, while afterwards, the deterministic part of the system prevails, constraining the (stochastic) solution to fluctuate around one of
the extremal deterministic solutions, given by $\pm ct^{1/(1-\beta)}$. \newline
More quantitatively, we can compare the fluctuations of the noise, say of order $\gamma>0$ with the fluctuations of the deterministic extremal solutions, giving
that
\[t^\gamma \, > \, t^{1/(1-\beta)}.\]
Since it should happen in small times, we then obtain that
\[\beta \, > \, 1-\frac{1}{\gamma},\]
should be the heuristic relation that guarantees the noise dominates in short time.
Clearly, the above inequality holds for any $\beta$ in $(0,1)$ in the Brownian case ($\gamma=1/2$), which would actually give $\beta>-1 $. We can refer to \cite{Delarue:Diel16} which is the closest  work to this threshold since the authors manage to reach $-2/3^+$.

In view the above arguments, we fix $n=N$, $d_i=d=1$ and $i$, $j$ in $\llbracket 1,n\rrbracket$ such that $j\ge i$ and we consider the drift
\[Ax+e_i\text{sgn}(x_j)\vert x_j\vert^\beta\]
where $\{e_i\colon i \in \llbracket 1, N \rrbracket\}$ is the canonical orthonormal basis for $\R^N$,  $A$ is the matrix in $\R^N\otimes \R^N$ given by
\[
A \, := \, \begin{pmatrix}
               0& \dots         & \dots         & \dots     & 0 \\
              1       & 0 & \dots         & \dots     & 0\\
               0 & 1       & \ddots & \ddots     & \vdots \\
               \vdots        & \ddots        & \ddots        & \ddots    & \vdots        \\
               0 & \dots         & 0 & 1 & 0
             \end{pmatrix}.
\]

We will assume moreover that $\beta$ is in $(0,1)$ such that
\[
\beta \, < \, \frac{1+\alpha(i-2)}{1+\alpha(j-1)},
\]
so that we are clearly outside the framework given by condition $(\mathscr{C})$.\newline
Our aim is to prove that uniqueness in law fails for the following equation:
\begin{equation}
\label{eq:SDE_for_counter}
    \begin{cases}
dX_t \, = \, \bigl[Ax+e_i\text{sgn}(X^j_t)\vert X^j_t\vert^\beta\bigr]dt+ BdZ_t, \quad t\ge 0, \\
X_0 \, = \, 0,
\end{cases}
\end{equation}
where $\{Z_t\}_{t\ge 0}$ is a symmetric, $d$-dimensional $\alpha$-stable process such that $\mathbb{E}[\vert Z_1 \vert]$ is finite.\newline
In particular, we are interested on the $i$-th component of the above Equation \eqref{eq:SDE_for_counter} that can be rewritten in integral form as:
\begin{equation}
\label{eq:SDE_for_counter_integral}
X^j_t \, = \, \int_0^t\text{sgn}\bigl(I^{j-i}_t(X^j)\bigr)\bigl{\vert}I^{j-i}_t(X^j)\bigr{\vert}^\beta dt + I^{i-1}_t(Z), \quad t\ge 0,
\end{equation}
where we have denoted  by $I^k_t(y)$ the $k$-th \emph{iterated integral} of a c\`adl\`ag path $y\colon [0,\infty)\to \R$ at a time $t$. Namely,
\begin{equation}
\label{eq:iterate_int_operator}
I^k_t(y) \, := \, \int_{0}^{t_k=t}\dots \int_{0}^{t_2} y_{t_0} \, dt_0\dots dt_{k-1}, \quad t\ge0.
\end{equation}

In order to improve the readability of the next part, we are going to present our reasoning in a slightly more general way. It is not difficult to check that
Equation \eqref{eq:SDE_for_counter_integral} satisfies the assumptions of the following proposition.
\begin{prop}
\label{prop:counter-example-formal}
Let $k$ be in $\N$, $\beta$ in $(0,1)$, $x$ in $\R$ and $\{\mathcal{Z}_t\}_{t\ge0}$ a continuous process on $\R$ such that
\begin{itemize}
\item $\mathbb{E}\bigl[\sup_{s\in [0,1]}|\mathcal{Z}_s|\bigr]<\infty$;
\item it is symmetric and $\gamma$-self-similar in law for some $\gamma>0$. Namely,
\[\bigl(\mathcal{Z}_t\bigr)_{t\ge 0}\, \overset{(\rm{law})}{=} \,\bigl(-\mathcal{Z}_t\bigr)_{t\ge 0}\,  \text{ and } \,\forall \rho>0,\ \bigl(\mathcal{Z}_{\rho t}\bigr)_{t\ge 0}\,
\overset{(\rm{law})}{=} \,\bigl(\mathcal{Z}_t \rho^\gamma\bigr)_{t\ge 0}.\]
\end{itemize}
Then, uniqueness in law fails for the following SDE:
\begin{equation}\label{eq:Peano_SDE}
\begin{cases}
dX_t \, = \, {\rm{sgn}}\bigl(I^k_t(X)\bigr) \bigl{\vert}I^k_t(X)\bigr{\vert}^\beta dt +d\mathcal{Z}_t, \quad t\ge 0 \\
X_0\, = \, x,
\end{cases}
\end{equation}
if $x=0$ and $\beta <\frac{\gamma-1}{\gamma+k}$.
\end{prop}
Since we can clearly apply Proposition \ref{prop:counter-example-formal} to Equation \eqref{eq:SDE_for_counter_integral} taking $\gamma=i-1+\frac{1}{\alpha}$,
$k=j-i$, it implies that SDE \eqref{eq:SDE_for_counter} lacks of uniqueness in law if
\[\beta \, < \,\frac{\gamma-1}{\gamma+k} \, = \, \frac{1+\alpha(i-2)}{1+\alpha(j-1)}.\]
Hence, to complete the proof of Theorem \ref{thm:counterexample}, it suffices to establish Proposition \ref{prop:counter-example-formal}.

Before proving Proposition \ref{prop:counter-example-formal}, we need however an auxiliary result. It roughly states that any solution of SDE \eqref{eq:Peano_SDE}
starting outside zero cannot immediately reach the extremal solutions of the associated deterministic Peano example. Importantly, the constant $\rho$ appearing below does not depend on the starting point $x$.

\begin{lemma}
\label{lemma:Peano_example}
Fixed $x>0$ and $\beta<\frac{\gamma-1}{\gamma+k}$, let $\{X_t\}_{t\ge 0}$ be a solution of Equation \eqref{eq:Peano_SDE} starting from $x$. Then, there exist two
positive constants $\rho:=\rho(k,\beta,\gamma,\mathbb{E}[\sup_{s\in [0,1]}|\mathcal{Z}_s|])$ and $c_0:=c_0(k,\beta)$ such that
\begin{equation}\label{eq:control_on_tau}
\mathbb{P}\bigl(\tau(X)\ge \rho\bigr) \, \ge \, 3/4,
\end{equation}
where $\tau(X)$ is the stopping time on $\Omega$ given by
\begin{equation}\label{eq:Peano_random_time}
\tau(X) \, = \, \inf\{t \ge 0 \colon X_t \, \le \, c_0 t^{\frac{k\beta+1}{1-\beta}}\}.
\end{equation}
\end{lemma}
\begin{proof}
We start noticing that the process $\{X_t\}_{t\ge0}$ is continuous in $0$, since it is c\`adl\`ag. Fixed $c_0>0$ to be chosen later, it implies that $\tau(X)>0$, almost surely. In particular, it makes sense to consider the random interval $(0,\tau(X)]$.\newline
Fixed $t$ in $(0,\tau(X)]$, it holds, by definition of $\tau(X)$, that $X_t>c_0t^{\frac{k\beta+1}{1-\beta}}$. It follows then that
\[\int_{0}^{t}\bigl{\vert} I^k_s(X)\bigr{\vert}^\beta \, ds \, > \, \tilde{C}c^\beta_0t^{\frac{k\beta+1}{1-\beta}} \,\, \text{ where } \,\,
\tilde{C} \,:= \, \bigl(\prod_{i=1}^{k}\frac{k\beta+1}{1-\beta}+(i-1)\bigr)^{-\beta}
.\]
Since $x>0$ by assumption and $X>0$ on $(0,\tau(X)]$, we can now show that
\[X_t \, = \, x + \int_{0}^{t}\text{sgn}\bigl(I^k_s(X)\bigr)\bigl{\vert} I^k_s(X)\bigr{\vert}^\beta \, ds + \mathcal{Z}_t \, > \, \tilde{C}c^\beta_0
t^{\frac{k\beta+1}{1-\beta}}+ \mathcal{Z}_t.\]
The next step is to write $\tilde{C}c^\beta_0=c_0+\hat{C}$ for some constant $\hat{C}>0$. To do so, we need to choose carefully $c_0$. In particular, the
condition above is equivalent to the following
\[\hat C=\tilde{C}c^\beta_0-c_0\, > \, 0 \Leftrightarrow c_0 \, < \, \tilde{C}^{\frac{1}{1-\beta}}.\]
Fixed $c_0=\tilde{C}^{\frac{1}{1-\beta}}/2$, it then holds that
\[X_t \, > \, c_0t^{\frac{k\beta+1}{1-\beta}}+\hat{C}t^{\frac{k\beta+1}{1-\beta}}+ \mathcal{Z}_t\]
for any $t$ in $(0,\tau(X)]$. Fixed $\rho>0$ to be chosen later, we can now define the event $\mathcal A$ in $\Omega$ as
\[\mathcal A \, := \, \{\omega \in \Omega \colon \hat{C}t^{\frac{k\beta+1}{1-\beta}}+ \mathcal{Z}_t>0, \,\, \forall \, t \in (0,\rho]\}.\]
On $\mathcal A$ and for any $t$ in $(0,\tau(X)]$, it then holds that
\[X_t \, > \, c_0t^{\frac{k\beta+1}{1-\beta}}.\]
In particular, we have that $\ \tau(X)\ge \rho$  on $\mathcal A$  and thus, $\mathcal A\subseteq\{\tau(X)\ge \rho\}$ on $\Omega$. It immediately implies that
\[\mathbb{P}\bigl(\tau(X)\ge \rho\bigr) \, \ge \, \mathbb{P}(\mathcal A).\]
It remains to choose $\rho>0$ such that $\mathbb{P}(\mathcal A)\ge 3/4$. Write:
\begin{align*}
\mathbb{P}(\mathcal A)=&\mathbb P[\forall t\in (0,\rho],\  \hat{C}t^{\frac{k\beta+1}{1-\beta}}+\mathcal{Z}_t>0]=\mathbb P[\forall t\in (0,1],\  \hat{C}(\rho t)^{\frac{k\beta+1}{1-\beta}}+\mathcal{Z}_{\rho t}>0]\\
=&\mathbb P[\forall t\in (0,1],\  \hat{C}(\rho t)^{\frac{k\beta+1}{1-\beta}}+\rho^\gamma \mathcal{Z}_{t}>0]=\mathbb P[\forall t\in (0,1],\  \hat{C}\rho^{\frac{k\beta+1}{1-\beta}-\gamma}+t^{-\frac{k\beta+1}{1-\beta}} \mathcal{Z}_{t}>0],
\end{align*}
from the self-similarity assumption on $\mathcal Z$. Since by assumption $\beta<\frac{\gamma-1}{\gamma+k}\iff\frac{k\beta+1}{1-\beta}-\gamma<0 $, the statement will follow taking $\rho $ small enough as soon as we prove the process $\mathcal R_t:=t^{-\frac{k\beta+1}{1-\beta}} \mathcal{Z}_{t},\ t\in (0,1] $, which is continuous on the open set $(0,1]$, can be extended by continuity in $0$ with $\mathcal R_0=0$. Observe that $\mathbb E[|\mathcal R_t|]=t^{\gamma-\frac{k\beta+1}{1-\beta}} \mathbb E[|\mathcal Z_1|]\underset{t\rightarrow 0}{\longrightarrow} 0$. Setting $\delta:= \gamma-\frac{k\beta+1}{1-\beta}>0$ and introducing $ t_n:=n^{-1/\delta(1+\eta)},\eta>0$, we get that for all $\varepsilon>0 $,
$$\mathbb P[|\mathcal R_{t_n}|\ge \varepsilon]\le \varepsilon^{-1}\mathbb E[|\mathcal R_{t_n}|]=\varepsilon^{-1}t_n^\delta \mathbb E[|\mathcal Z_1|]=\varepsilon^{-1}n^{-(1+\eta)}\mathbb E[|\mathcal Z_1|].$$
We thus get from the Borel-Cantelli lemma that $\mathcal R_{t_n}\underset{n,\ a.s.}{\longrightarrow} 0 $. Namely, we have almost sure convergence along the subsequence $t_n$ going to zero with $n$. It now remains to prove that the process $\mathcal R_{t} $ does not fluctuate much between two successive times $t_n $ and $t_{n+1}$. Write for $t\in [t_{n+1},t_n]$:
\begin{align}
|R_t|:=|t^{-\frac{k\beta+1}{1-\beta}} \mathcal{Z}_{t}|\le& t_{n+1}^{-\frac{k\beta+1}{1-\beta}}\Big(|\mathcal{Z}_{t_{n+1}}|+ \sup_{s\in [t_{n+1},t_n]}|\mathcal{Z}_{s}-\mathcal{Z}_{t_{n+1}}|\Big)\notag\\
\le &  t_{n+1}^{-\frac{k\beta+1}{1-\beta}}\Big(2|\mathcal{Z}_{t_{n+1}}|+ \sup_{s\in [0,t_n]}|\mathcal{Z}_{s}|\Big)\label{CTR_REMAIN_BC}.
\end{align}
The first term of the above left hand side tends almost surely to zero with $n$. Observe as well that, from the scaling properties of $\mathcal{Z}$, for any $\varepsilon>0 $:
\begin{align*}
 \mathbb P[ t_{n+1}^{-\frac{k\beta+1}{1-\beta}} \sup_{s\in [0,t_n]}|\mathcal{Z}_{s}|\ge \varepsilon]&=\mathbb P[ t_{n+1}^{-\frac{k\beta+1}{1-\beta}}t_n^\gamma \sup_{s\in [0,1]}|\mathcal{Z}_{s}|\ge \varepsilon]\le \varepsilon^{-1}t_n^{\delta} (\frac{t_n}{t_{n+1}})^{\frac{k\beta+1}{1-\beta}} \mathbb E[\sup_{s\in [0,1]}|\mathcal{Z}_{s}|]\\
 &\le C\varepsilon^{-1}n^{-(1+\eta)}\mathbb E[\sup_{s\in [0,1]}|\mathcal{Z}_{s}|],
\end{align*}
which again gives from the Borel-Cantelli lemma the a.s. convergence with $n$ of the second term in the r.h.s of \eqref{CTR_REMAIN_BC}. We eventually derive
that $\mathcal R_t \underset{t\rightarrow 0, a.s.}{\longrightarrow} 0$.
Again, the key point is that we normalize  the process $\mathcal Z$ at a rate, $t^{\frac{k\beta+1}{1-\beta}} $, which is lower than its own characteristic time scale, $t^{\gamma} $. This is precisely what leaves some margin to establish continuity.

\end{proof}

Exploiting the lower bound for the random time $\tau(X)$ given in Lemma \ref{lemma:Peano_example}, we are now ready to show uniqueness in law fails for SDE \eqref{eq:Peano_SDE} when $x=0$
and $\beta<\frac{\gamma-1}{\gamma+k}$.

\emph{Proof of Proposition \ref{prop:counter-example-formal}.}
By contradiction, we start assuming that uniqueness in law holds for SDE \eqref{eq:Peano_SDE} starting at $x=0$. Fixed any solution $\{X_t\}_{t\ge 0}$ of Equation \eqref{eq:Peano_SDE}
starting at zero, it follows by symmetry that $\{-X_t\}_{t\ge 0}$ is also a solution of the same dynamics.
Since by hypothesis, $-\mathcal{Z}_t\overset{(\text{law})}{=}\mathcal{Z}_t$, uniqueness in law for SDE \eqref{eq:Peano_SDE} implies that the laws of $X$ and $-X$ are identical.\newline
Assuming for the moment that Lemma \ref{lemma:Peano_example} is applicable for $x=0$, we easily find a contradiction. Indeed, it follows from Lemma \ref{lemma:Peano_example}
that
\[\mathbb{P}\bigl(\tau(X)\ge \rho\bigr) \, \ge \, 3/4\]
but on the same time, thanks to the uniqueness in law, we have that
\[\mathbb{P}^0\bigl(\tau(-X)\ge \rho\bigr) \, \ge \, 3/4,\]
which is clearly impossible. To show the validity of Lemma \ref{lemma:Peano_example} in $x=0$, we consider a a sequence $\{\{X^n_t\}_{t\ge 0}\colon n \in \N\}$ of solutions of SDE \eqref{eq:Peano_SDE} starting at $1/n$. It is then easy to check that such a sequence satisfies the Aldous criterion:
\[\mathbb{E}[\vert X^n_t-X^n_0\vert^p] \, \le \, ct^{p\gamma}, \quad t\ge 0\]
for some $p>0$ and $c>0$ independent from $t$ and $n$. It follows (Proposition $34.8$ in \cite{book:Bass11}) that the sequence $\{\mathbb{P}^n\}_{n\in \N}$
of the laws of $\{X^n_t\}_{t\ge0}$ is tight.
Prohorov Theorem (cf. Theorem $30.4$ in \cite{book:Bass11}) ensures now the existence of a converging sub-sequence
$\{\mathbb{P}^{n_k}\}_{k\in \N}$. The uniqueness in law then implies that the sequence $\{\mathbb{P}^{n_k}\}_{k\in \N}$ converges, as
expected, to $\mathbb{P}^0$ the law of the solution starting at $0$. Noticing that inequality \eqref{eq:control_on_tau} holds for any solution $\{X^n_t\}_{t\ge 0}$ and moreover, the constant $\rho$ is independent from the starting points $1/n$, we find that
\[\mathbb{P}\bigl(\tau(X)\ge \rho\bigr) \, \ge \, 3/4.\]
The proof of Proposition \ref{prop:counter-example-formal} is thus concluded.

\setcounter{equation}{0}
\section{Appendix: proofs of complementary results}
\fancyhead[RO]{Section \thesection. Appendix}

\subsection{Controls on the density of the proxy process}
We present here two useful lemmas needed to complete the proof of Proposition \ref{prop:Smoothing_effect}. We will analyze the behavior of the laws of the independent random variables $ \tilde{M}^{\tau,\xi,t,s}$ and $\tilde{N}^{\tau,\xi,t,s}$ obtained in \eqref{eq:decomposition_S} by truncation of the process $\tilde{S}^{\tau,\xi,t,s}$ at the associated stable time scale $u^{1/\alpha}$.

\begin{lemma}
\label{lemma:Control_p_M}
Let $m$ be in $\N$. Then, there exists a positive constant $C:=C(m,T)$ such that for any $k$ in $\llbracket 0, m \rrbracket$,
\[
\left| D^k_{z} p_{\tilde{M}^{\tau,\xi,t,s}}(u,z) \right| \, \le \, Cu^{-(N+k)/\alpha}\left(1+\frac{|z|}{u^{1/\alpha}}\right)^{-m} \, =: \, Cu^{-k/\alpha}p_{\overline{M}}(u,z),
\]
for any $u>0$, any $z$ in $\R^N$, any $t\le s$ in $[0,T]$ and any $(\tau,\xi)$ in $[0,T]\times \R^N$.
\end{lemma}
\begin{proof}
Similarly to the proof of Proposition \ref{prop:Decomposition_Process_X} (see in particular Equation \eqref{eq:definition_density_q}), we start writing
\[
p_{\tilde{M}^{\tau,\xi,t,s}}(u,z) \, = \,\frac{1}{(2\pi)^N}\int_{\R^N}e^{-i\langle z,y\rangle}\text{exp}\left(u\int_{|p|\le u^{1/\alpha}}\left[\cos(\langle y,p \rangle)-1\right]\, \nu_{\tilde{S}^{\tau,\xi,t,s}}(dp)\right)\, dy,
\]
\textcolor{black}{where, we recall, $\nu_{\tilde{S}}^{\tau,\xi,t,s}$ is the L\'evy measure associated with the process $\{\tilde{S}^{\tau,\xi,t,s}_{u}\}_{u\ge0}$ in Proposition \ref{prop:Decomposition_Process_X}}. Setting $u^{1/\alpha}y =\tilde{y}$ then yields
\begin{align}\notag
p_{\tilde{M}^{\tau,\xi,t,s}}(u,z) \, &= \, \frac{u^{-N/\alpha}}{(2\pi)^N}\int_{\R^N}  e^{-i\langle z, \frac{\tilde{y}}{u^{1/\alpha}}\rangle}\text{exp}\left(u\int_{|p|\le u^{1/\alpha}}\left[\cos(\langle \tilde{y},\frac{p}{u^{1/\alpha}} \rangle)-1\right] \nu_{\tilde{S}^{\tau,\xi,t,s}}(dp)\right) d\tilde{y} \\
&=: \, \frac{u^{-N/\alpha}}{(2\pi)^N}\int_{\R^N}  e^{-i\langle \frac{z}{u^{1/\alpha}},\tilde{y}\rangle}\hat{f}^{\tau,\xi,t,s}_u(\tilde{y})\, d\tilde{y}
\label{EXP_DENS_M}
\end{align}
Since the L\'evy measure $\nu_{\tilde{S}}^{\tau,\xi,t,s}$ in the expression above has finite support, Theorem $3.7.13$ in Jacob \cite{book:Jacob05} implies that $\hat{f}^{\tau,\xi,t,s}_u$ is infinitely differentiable in $\tilde{y}$. We can thus calculate
\[
\begin{split}
|\partial_{\tilde{y}} \hat{f}^{\tau,\xi,t,s}_u(\tilde{y})| \, &\le \, u\int_{|p|\le u^{1/\alpha}}\frac{|p|}{u^{1/\alpha}}\left|\sin\left(\bigl{\langle} \tilde{y},\frac{p}{u^{1/\alpha}} \bigr{\rangle}\right)\right|\, \nu_{\tilde{S}^{\tau,\xi,t,s}}(dp) \\
&\qquad\qquad\qquad\qquad \times \text{exp}\left(u\int_{|p|\le u^{1/\alpha}}\left[\cos\left(\bigl{\langle} \frac{\tilde{y}}{u^{1/\alpha}},p \bigr{\rangle}\right)-1\right]\, \nu_{\tilde{S}^{\tau,\xi,t,s}}(dp)\right).
\end{split}
\]
Recalling that $\alpha>1$, we can now write that
\[
\begin{split}
 u\int_{|p|\le u^{1/\alpha}}\frac{|p|}{u^{1/\alpha}}\left|\sin\left(\bigl{\langle} \tilde{y},\frac{p}{u^{1/\alpha}} \bigr{\rangle}\right)\right|\, \nu_{\tilde{S}^{\tau,\xi,t,s}}(dp)\, &\le \, Cu\int_{r\le u^{1/\alpha}}  \frac{r}{u^{1/\alpha}}    \frac{|\tilde{y}|r}{u^{1/\alpha}} \frac{dr}{r^{1+\alpha}} \\
&\le \, Cu\int_{r\le u^{1/\alpha}} |\tilde{y}| \frac{r^{1-\alpha}}{u^{2/\alpha}} \, dr \\
&\le \,  C(1+|\tilde{y}|).
\end{split}
\]
It then follows that
\[\begin{split}
&|\partial_{\tilde{y}} \hat{f}^{\tau,\xi,t,s}_u(\tilde{y})|
\\
&\,\,\qquad\le \, C(1+|\tilde{y}|) \text{exp}\left(u\int_{\R^N}\left[\cos\left(\bigl{\langle} \frac{\tilde{y}}{u^{1/\alpha}},p \bigr{\rangle}\right)-1\right]\, \nu_{\tilde{S}^{\tau,\xi,t,s}}(dp)\right)e^{2u\nu_{\tilde{S}^{\tau,\xi,t,s}}(B^c(0,u^{1/\alpha}))}\\
&\,\,\qquad\le \,
C (1+|\tilde{y}|)\exp(-C^{-1}|\tilde{y}|^\alpha),
\end{split}\]
where in second inequality we exploited Control \eqref{eq:control_Levy_symbol_S} and
\begin{equation}
\label{FINITE_MEAS}
\nu_{\tilde{S}^{\tau,\xi,t,s}}(B^c(0,u^{1/\alpha})) \le C/u.
\end{equation}
Iterating the above reasoning, we can then show that for any $l$ in $\N$,
\[
|\partial^l_{\tilde{y}}\hat{f}^{\tau,\xi,t,s}_u(\tilde{y})| \, \le\,
C_l (1+|\tilde{y}|^l)\exp(-C^{-1}|{\tilde{y}}|^\alpha),
\]
for some positive constant $C:=C(l)$.
It implies in particular that $\hat{f}^{\tau,\xi,t,s}_u(\tilde{y})$ is a Schwartz test function.
Denoting by $f^{\tau,\xi,t,s}_u$ its inverse Fourier transform, we thus have
that for any $m$ in $\N$, there exists a positive constant $C:=C(m)$ such that
\[
 |f^{\tau,\xi,t,s}_u(y)| \le C_m (1+|y|)^{-m}, \quad y \in \R^N.
\]
The result for $k=0$ now follows immediately noticing that 
\[p_{\overline{M}}(t-s,y)\,  = \, (t-s)^{-\frac d\alpha} f_{s,t}(y/(t-s)^{\frac{1}{\alpha}}).\]
The controls on the derivatives can be derived analogously.
\end{proof}

We can now show a similar control on the law of the process $\tilde{N}^{\tau,\xi,t,s}$.

\begin{lemma}
\label{lemma:Control_P_N}
There exists a family $\{\overline{P}_u\}_{u\ge 0}$ of Poisson measures and a positive constant $C:=C(T,N)$ such that for any $\mathcal A$ in $\mathcal{B}(\R^N)$ and $\tilde{N}^{\tau,\xi,t,s}$ as in \eqref{DEF_TILDE_N},
\begin{equation}\label{DEF_POISSON_QUI_DOMINE}
P_{\tilde{N}^{\tau,\xi,t,s}_u}(\mathcal A) \, \le \, C\overline{P}_u(\mathcal A).
\end{equation}
\end{lemma}
\begin{proof}
For notational simplicity, we start introducing \textcolor{black}{the truncated L\'evy measure  associated with the big jumps of the process} $\{\tilde{S}^{\tau,\xi,t,s}_u\}_{u\ge 0}$:
\[\nu_{\text{tr}}^{\tau,\xi,t,s}(dp)\, = \, \mathds{1}_{|p|\ge u^{1/\alpha}}(p)\nu_{\tilde{S}}^{\tau,\xi,t,s}(dp). \]
It follows immediately that $\nu_{\text{tr}}^{\tau,\xi,t,s}$ is a finite measure (see \eqref{FINITE_MEAS} above). With this notation at hand, we can write:
\[
\begin{split}
\widehat {P_{\tilde{N}^{\tau,\xi,t,s}_u}}(y) \, &= \, \exp\left(u \int_{|p|>u^{\frac{1}{\alpha}}}\left[\cos(\langle y,p\rangle)-1\right] \, \nu_{\tilde{S}}^{\tau,\xi,t,s}(dp)\right) \\
&= \, \exp\left( u \widehat{\nu_{\text{tr}}^{\tau,\xi,t,s}}(y)- u \nu_{\text{tr}}^{\tau,\xi,t,s}(\R^N)\right),
\end{split}
\]
where $\widehat{\nu}$ denotes the Fourier-Stieltjes transform of the considered measure $\nu$. Let us introduce then the following measure:
\[\zeta^{\tau,\xi,t,s}\, := \, u\nu_{\text{tr}}^{\tau,\xi,t,s}.\]
Expanding the previous exponential and by termwise Fourier inversion, we now find that
\begin{equation}
\label{proof:appendix_decomposition_P_N}
P_{\tilde{N}^{\tau,\xi,t,s}_u} \, = \, \exp\left(\zeta^{\tau,\xi,t,s}- u \nu_{\text{tr}}^{\tau,\xi,t,s}(\R^N)\right) \, = \, \exp\left(-u\nu_{\text{tr}}^{\tau,\xi,t,s}(\R^N)\right)\sum_{n \in \N}\frac{\left( \zeta^{\tau,\xi,t,s}
\right)^{\star n}}{n!},
\end{equation}
where, for a finite measure $\rho$ on $\R^N$, $(\rho)^{\star n}:= \rho\star\cdots \star \rho$ denotes its $n^{{\rm th}}$ fold convolution.\newline
For now, let us assume that $\sigma(t,x)$ is non-constant in space, so that
\[B\tilde{\sigma}^{\tau,\xi}_{u(v)}\, = \, B\sigma\left(u(v),\theta_{u(v),\tau}(\xi)\right)\]
\textcolor{black}{appearing in the definition of $\nu_{\tilde{S}}^{\tau,\xi,t,s}$}, truly depends on the parameters $\tau,\xi$. Assumption [\textbf{AC}] then ensures the existence of a bounded function $g\colon \mathbb{S}^{d-1}\to \R$ such that
\[\nu(dp) \, = \, Q(p)\frac{g(\frac{p}{|p|})}{|p|^{d+\alpha}}dp.\]
From Equation \eqref{proof:appendix_decomposition_P_N}, it is clear that we need to control the measure $\zeta^{\tau,\xi,t,s}$, uniformly in the parameters $\tau,\xi,t,s$. Namely, for any $\mathcal A$ in $\mathcal{B}(\R^N)$, we write from \eqref{proof:eq:def_Levy_symbol_S} that
\[
\begin{split}
\zeta^{\tau,\xi,t,s}(\mathcal A) \, &= \, u\int_{|p|>u^{\frac{1}{\alpha}}}\mathds{1}_{\mathcal A}(p)\, \nu_{\tilde{S}}^{\tau,\xi,t,s}(dp)\, = \, u\int_0^1\int_{|\widehat{\mathcal{R}}_vB\tilde{\sigma}^{\tau,\xi}_{u(v)}p|>u^{\frac{1}{\alpha}}}\mathds{1}_{\mathcal A}(\widehat{\mathcal{R}}_vB\tilde{\sigma}^{\tau,\xi}_{u(v)}p)\, \nu(dp)dv\\
&= u\int_0^1\int_{|\widehat{\mathcal{R}}_vB\tilde{\sigma}^{\tau,\xi}_{u(v)}p|>u^{\frac{1}{\alpha}}}\mathds{1}_{\mathcal A}(\widehat{\mathcal{R}}_vB\tilde{\sigma}^{\tau,\xi}_{u(v)}p)\frac{g(\frac{p}{|p|})}{|p|^{d+\alpha}}Q(p)\, dpdv\\
&\le \, u\int_0^1 \int_{|\widehat{\mathcal{R}}_vB\tilde{\sigma}^{\tau,\xi}_{u(v)}p|>u^{\frac{1}{\alpha}}}\mathds{1}_{\mathcal A}(\widehat{\mathcal{R}}_vB\tilde{\sigma}^{\tau,\xi}_{u(v)}p)\frac{dp}{|p|^{d+\alpha}}dv
\end{split}
\]
We can then exploit assumption [\textbf{UE}] on $\sigma$ to conclude that
\[
\begin{split}
\zeta^{\tau,\xi,t,s}(\mathcal A) \, &\le \,  u\int_0^1\int_{|\widehat{\mathcal{R}}_vBq|>u^{\frac{1}{\alpha}}}\mathds{1}_{\mathcal A}(\widehat{\mathcal{R}}_vB q)\frac{1}{\det (\tilde{\sigma}^{\tau,\xi}_{u(v)})}\frac{dq}{|(\tilde{\sigma}^{\tau,\xi}_{u(v)})^{-1}q|^{d+\alpha}}dv\\
&\le \, Cu\int_0^1\int_{|\widehat{\mathcal{R}}_vBq|>u^{\frac{1}{\alpha}}}\mathds{1}_{\mathcal A}(\widehat{\mathcal{R}}_vBq)\frac{dq}{|q|^{d+\alpha}}dv.
\end{split}
\]
Denoting now by $\Lambda_{\text{tr}}:=c\mathds{1}_{p>u^{1/\alpha}}\frac{dp}{p^{d+\alpha}}$ the truncated L\'evy measure of the isotropic $\alpha$-stable process and by $\overline{\nu}_{\text{tr}}$ the following push-forward measure
\[\overline{\nu}_{\text{tr}}(\mathcal A) \, := \, \int_0^1\Lambda_{\text{tr}}\left((\widehat{\mathcal{R}}_vB)^{-1}\mathcal A\right)dv, \quad \mathcal A \in \mathcal{B}(\R^N)\]
we derive that there exists a constant $C$ such that for any $(\tau,\xi)$ in $[0,T]\times \R^N$, $t\le s$ in $[0,T]$,
\begin{equation}
\label{proof:appendix_control_N}
\zeta^{\tau,\xi,t,s}(\mathcal A) \, \le \, C u \int_0^1\Lambda_{\text{tr}}\left((\widehat{\mathcal{R}}_vB)^{-1}\mathcal A\right)dv \, = \, u\overline{\nu}_{\text{tr}}(\mathcal A) \, =: \, \overline{\zeta}(\mathcal A).
\end{equation}

Equation \eqref{DEF_POISSON_QUI_DOMINE} now follows from
the above control, \eqref{FINITE_MEAS} and \eqref{proof:appendix_decomposition_P_N}, denoting
\[\overline{P}_u \, := \, \exp\left(-u\overline{\nu}_{\text{tr}}(\R^N)\right)\sum_{n \in \N}\frac{(\overline{\zeta})^{\star n}}{n!},\]
up to a modification of the constant $C$ in \eqref{proof:appendix_control_N}.
Following backwards the same reasoning presented at the beginning of the proof, we then notice that
\[
\begin{split}
\widehat {\overline{P}_u}(y) \, &= \, \exp\left(u \int_0^1\int_{\R^N}\left[\cos(\langle y,p\rangle)-1\right] \, \overline{\nu}_{\text{tr}}(dp)dv\right) \\
&= \, \exp \left(u\int_0^1 \int_{\R^d}\mathds{1}_{\{|\widehat{\mathcal{R}}_vBp|>u^{\frac{1}{\alpha}}\}}\left[\cos(\langle y,\widehat{\mathcal{R}}_vBp\rangle)-1\right] \, \Lambda(dp)dv\right)\\
&= \, \exp \left(u\int_0^1 \int_0^\infty \int_{\mathbb{S}^{d-1}} \mathds{1}_{\{|\widehat{\mathcal{R}}_vB\theta r|>u^{\frac{1}{\alpha}}\}}\left[\cos(\langle y,\widehat{\mathcal{R}}_vB\theta r\rangle)-1\right] \, \mu_{\text{leb}}(d\theta) \frac{dr}{r^{1+\alpha}}dv\right),
\end{split}
\]
where we used the spherical decomposition \textcolor{black}{for the L\'evy measure $\Lambda$ of an isotropic $\alpha$-stable process:}
\begin{equation}
\Lambda(dp) := \frac{dp}{p^{d+\alpha}} \, = \, C\mu_{\text{leb}}(d\theta)\frac{dr}{r^{1+\alpha}},
\end{equation}
with $p=r \theta$ and $\mu_{\text{leb}}$ Lebesgue measure on the sphere $\mathbb{S}^{d-1}$.\newline
We exploit now the non-degeneracy of $\widehat{\mathcal{R}}_v$ to to define two functions $k\colon [0,1]\times \mathbb{S}^{d-1}\to \R$ and $l\colon [0,1]\times
\mathbb{S}^{d-1}\to \mathbb{S}^{N-1}$, given by
\[k(v,\theta) \,:= \,
\vert \widehat{\mathcal{R}}_{v}B\theta\vert \,\, \text{ and } \,\,  l(v,\theta) \,:= \, \frac{\widehat{\mathcal{R}}_{v}B\theta}{\vert\widehat{\mathcal{R}}_{v}B\theta\vert}.\]
Using the Fubini theorem, we can now write that
\[
\begin{split}
&\widehat{\overline{P}_u}(y) \\
&\,\,= \, \exp\left(u\int_0^1\int_{0}^{\infty}\int_{\mathbb{S}^{d-1}} \mathds{1}_{\{|l(v,\theta)k(v,\theta)r|>u^{\frac{1}{\alpha}}\}}
\left[\cos\left(\langle z,l(v,\theta)k(v,\theta)r\rangle\right)-1\right]\, \mu_{\text{leb}}(d\theta)\frac{dr}{r^{1+\alpha}}dv\right) \\
&\,\,= \, \exp\left(u\int_0^1\int_{0}^{\infty}\int_{\mathbb{S}^{d-1}}\mathds{1}_{\{|l(v,\theta)\tilde{r}|>u^{\frac{1}{\alpha}}\}}
\left[\cos\left(\langle z,l(v,\theta)\tilde{r}\rangle\right)-1\right]\, [k(v,\theta)]^\alpha \mu_{\text{leb}}(d\theta)\frac{d\tilde{r}}{\tilde{r}^{1+\alpha}}dv\right).
\end{split}
\]
Denoting now by $\tilde{k}(dv,d\theta)$ the measure on $[0,1]\times \mathbb{S}^{d-1}$ given by
\[\tilde{k}(dv,d\theta) \, := \, [k(v,\theta)]^{\alpha}\mu_{\text{leb}}(d\theta)dv\]
and by $\tilde{\mu}_{\text{sym}}:=\text{Sym}(l)_{\ast}\tilde{k}$ the symmetrization of the measure $\tilde{k}(dv,d\theta)$ push-forwarded through $l$ on $\mathbb{S}^{N-1}$, we can finally conclude that
\begin{align}
\notag \widehat {\overline{P}_u}(y) \, &=\, \exp\left(u\int_0^\infty
\int_{[0,1]\times \mathbb{S}^{d-1}}\mathds{1}_{\{|l(v,\theta)\tilde{r}|>u^{\frac{1}{\alpha}}\}}
\left[\cos\left(\langle z,l(v,\theta)\tilde{r}\rangle\right)-1\right]\, \tilde{k}(dv,d\theta)\frac{d\tilde{r}}{\tilde{r}^{1+\alpha}}\right) \\
&= \, \exp\left(u\int_{|u|^{\frac{1}{\alpha}}}^{\infty}\int_{\mathbb{S}^{N-1}}\left[\cos\left(\langle z, \tilde{\theta} \tilde{r} \rangle\right)-1\right]\tilde{\mu}_{\text{sym}}(d\tilde{\theta})\frac{d\tilde{r}}{\tilde{r}^{1+\alpha}}\right).\label{eq:representation_P_N_segnato}
\end{align}
It is easy to check now that the measure $\tilde{\mu}_{\text{sym}}$ is
finite and non-degenerate in the sense of \eqref{eq:non_deg_measure}. This \textcolor{black}{concludes} the proof of our result under the additional assumption that $\nu$ is absolutely continuous with respect to \textcolor{black}{the} Lebesgue measure. \newline
If this is not the case, assumption [\textbf{AC}] implies immediately that $\sigma(t,x)=:\sigma_t$ does not depends on $x$. Thus,  the ``frozen'' diffusion $\tilde{\sigma}^{\tau,\xi}_t$ does not depends on the parameters $\tau,\xi$ as well. The same arguments above then allow to conclude in a similar manner.
\end{proof}

\paragraph{Sketch of proof for Proposition \ref{coroll:Smoothing_effect_delta}}
We briefly present here the proof of Proposition \ref{coroll:Smoothing_effect_delta} concerning the existence and the associated controls for the density of the mollified frozen process $\tilde{X}^{\tau,\xi,t,x,\delta}_{s}$.\newline
We start noticing that the reasoning in the proof of Proposition \ref{prop:Decomposition_Process_X} can be similarly applied. Indeed, from the definition in \eqref{eq:decomposition_delta_frozen}, it follows immediately that
\[\tilde{X}^{\tau,\xi,t,x,\delta}_{s} \, = \, \tilde{m}^{\tau,\xi}_{s,t}(x)+\mathbb{M}_{s-t}\left(\tilde{S}^{\tau,\xi,t,s}_{s-t}+\delta\overline{Z}_{s-t}\right),\]
and thus, that there exists a density $\tilde{p}^{\tau,\xi,\delta}(t,s,x,y)$ associated with the frozen process $\tilde{X}^{\tau,\xi,t,x,\delta}_{s}$. Moreover, the representation in \eqref{eq:representation_density} holds again if we change there the L\'evy measure $\nu_{\tilde{S}}^{\tau,\xi,t,s}$ with the one associated with the following L\'evy symbol:
\[\Phi_{\tilde{S}^{\tau,\xi,t,s,\delta}}(z) \, := \,
\Phi_{\tilde{S}^{\tau,\xi,t,s}}(z) + c_\alpha\delta |z|^\alpha\, = \,
\int_{0}^{1}\Phi\bigl((\widehat{\mathcal{R}}_{v}B\tilde{\sigma}^{\tau,\xi}_{u(v)})^*z\bigr)\,dv+c_\alpha\delta |z|^\alpha.\]
Namely, it holds that
\begin{multline*}
\tilde{p}^{\tau,\xi,\delta}(t,s,x,y) \, = \, \frac{\det \mathbb{M}^{-1}_{s-t}}{(2\pi)^N}\int_{\R^N}e^{-i\langle \mathbb{M}^{-1}_{s-t}(y-\tilde{m}^{\tau,\xi}_{s,t}(x)),z\rangle}\\
\times\exp\left((s-t)\int_{\R^N}\left[\cos(\langle z, p \rangle )-1\right]\nu_{\tilde{S}^{\tau,\xi,t,s,\delta}}(dp)\right)\, dz,
\end{multline*}
where the L\'evy measure $\nu_{\tilde{S}^{\tau,\xi,t,s,\delta}}$ is given by
\begin{equation}
\label{eq:decomposition_delta}
\nu_{\tilde{S}^{\tau,\xi,t,s,\delta}}(\mathcal{A}) \, = \, \nu_{\tilde{S}^{\tau,\xi,t,s}}(A)+\delta^\alpha \nu_{\overline{Z}}(\mathcal{A}), \quad \mathcal{A} \in \mathcal{B}(\R^N),
\end{equation}
with $\nu_{\overline{Z}}$ L\'evy measure of the isotropic $\alpha$-stable process $\overline{Z}_t$.
In particular, the L\'evy symbol $\Phi_{\tilde{S}^\delta}$ satisfies Control \eqref{eq:control_Levy_symbol_S} for a constant $C$ independent from $\delta$.\newline
We can now move to show the controls on the derivatives
of the mollified frozen density. It is not difficult to check that the arguments presented in the proofs of Proposition \ref{prop:Smoothing_effect}, Lemmas \ref{lemma:Control_p_M} and \ref{lemma:Control_P_N} can be applied again if we substitute there the L\'evy measure $\nu_{\tilde{S}^{\tau,\xi,t,s}}$ with the mollified one $\nu_{\tilde{S}^{\tau,\xi,t,s,\delta}}$. Indeed, taking into account the decomposition in \eqref{eq:decomposition_delta}, we notice that the L\'evy measure $\nu_{\tilde{S}^{\tau,\xi,t,s,\delta}}$ only considers an additional term ($\delta \nu_{\overline{Z}}$) that has the same $\alpha$-scaling nature considered before (but is however much less singular).\newline
To show instead that the estimates \eqref{eq:smoothing_effect_frozen_delta} are indeed uniform in the parameter $\delta$, it is sufficient to notice from \eqref{eq:decomposition_delta} that we have that
\[\nu_{\tilde{S}^{\tau,\xi,t,s,\delta}}(\mathcal{A}) \, \le \, \nu_{\tilde{S}^{\tau,\xi,t,s}}(\mathcal{A})+ \nu_{\overline{Z}}(\mathcal{A}), \quad \mathcal{A} \in \mathcal{B}(\R^N).\]
To conclude the proof of Proposition \ref{coroll:Smoothing_effect_delta}, it is then enough to take $\xi=y$, $\tau=s$ and to follow the same arguments introduced in the proof of Corollary \ref{coroll:Smoothing_effect}.

\subsection{Proof of the technical lemmas}
\label{SEC_TEC_LEMMA_APP}
\paragraph{Proof of Lemma \ref{lemma:bilip_control_flow} (Approximate Lipschitz condition of the flows)}
We start considering two measurable flows $\theta,\check{\theta}$ satisfying dynamics \eqref{eq:measurability_flow}. Recalling the decomposition $G(t,x)=A_tx +F(t,x)$, it follows immediately that:
\begin{equation}
\begin{split}
\label{eq:proof_bilip_control_flow1}
\mathbb{T}_{s-t}^{-1}(x-\theta_{t,s}(y))\, &= \, \mathbb{T}_{s-t}^{-1}\Bigl[\check{\theta}_{s,t}(x)-y-\int_t^s \Bigl(G(u,\check{\theta}_{u,t}(x))-G(u,\theta_{u,s}(y))\Bigr) \, du\Bigr] \\
&= \, \mathbb{T}_{s-t}^{-1}\bigl(\check{\theta}_{s,t}(x)-y\bigr)+ \mathcal{I}_{s,t}(x,y),
\end{split}
\end{equation}
where in the last step, we denoted
\[\mathcal{I}_{s,t}(x,y) \, = \, \mathbb{T}_{s-t}^{-1}\int_{t}^s \left[A_u\left(\theta_{u,s}(y)-\check{\theta}_{u,t}(x)\right)+\left(F(u,\theta_{u,s}(y))-F(u,\check{\theta}_{u,t}(x))\right)\right] \, du.\]
To conclude, we need to show the following bound for $\mathcal{I}_{s,t}(x,y)$:
\begin{equation}
\label{eq:proof_bilip_control_flow2}
|\mathcal{I}_{s,t}(x,y)| \, \le \,  C\left[1+(s-t)^{-1}\int_t^s |\mathbb{T}_{s-t}^{-1}(\check{\theta}_{u,t}(x)-\theta_{u,s}(y))|\,    du\right].
\end{equation}
Indeed, Control \eqref{eq:proof_bilip_control_flow2} together with \eqref{eq:proof_bilip_control_flow1} and the Gronwall lemma imply the right-hand side of Control \eqref{eq:bilip_control_flow}. The left-hand side one can be obtained analogously and we will not show it here.\newline
We start decomposing $\mathcal{I}_{s,t}$ into $\mathcal{I}^1_{s,t}+\mathcal{I}^2_{s,t}$, where we denote
\begin{align*}
\mathcal{I}^1_{s,t}(x,y) \, := \, \mathbb{T}_{s-t}^{-1}\int_{t}^s A_u\left(\theta_{u,s}(y)-\check{\theta}_{u,t}(x)\right) \, du; \\
\mathcal{I}^2_{s,t}(x,y) \, := \, \mathbb{T}_{s-t}^{-1}\int_{t}^s \left[F(u,\theta_{u,s}(y))-F(u,\check{\theta}_{u,t}(x))\right] \, du.
\end{align*}
The first remainder $\mathcal{I}^1_{s,t}$ can be controlled easily, exploiting the linearity of $z\to A_uz$. Indeed, for any $z$,$z'$ in $\R^N$ and any $u$ in $[s,t]$, we have that
\begin{equation}
\label{eq:proof_bilip_control_flow6}
\begin{split}
    \left|\mathbb{T}_{s-t}^{-1} A_u(z-z')\right| \, &\le \, \sum_{i=1}^n \sum_{j=(i-1)\vee 1}^n(s-t)^{-\frac{1+\alpha(i-1)}{\alpha}}|A^{i,j}_u|\ |(z-z')_j| \\
    &\le \, C(s-t)^{-1}|\mathbb{T}_{s-t}^{-1}(z-z')|.
\end{split}
\end{equation}
To control instead the second term $\mathcal{I}^2_{s,t}$, we will need to thoroughly exploit an appropriate smoothing method, due to the low regularity in space of the drift $F$. To overcome this problem, we are going to mollify the function $F$ in the following way. We start fixing a family  $\{\rho_i\colon i \in \llbracket 1,n\rrbracket\}$ of mollifiers on $\R^{D_i}$ where $D_i=N-\sum_{j=1}^{i-1}d_j$, i.e. for any $i$ in $\llbracket 1,n\rrbracket$, $\rho_i$ is a
compactly supported, non-negative, smooth function on $\R^{D_i}$ such that $\Vert \rho_i \Vert_{L^1}=1$, and a family $\{\delta_{ij}\colon i \le j\}$ of positive constants to be
chosen later. Then, the mollified version of the drift is defined by $F^\delta:=(F_1,F^{\delta}_2,\dots,F^{\delta}_n)$
where
\begin{equation}
\label{eq:multi_scale_mollification}
\begin{split}
    F^{\delta}_i(t,z) \, &:= \,  F_i \ast_x \rho^{\delta}_i(t,z) \\
    &:= \, \int_{\R^{D_i}}
 F_i(t, z_i-\omega_i,\dots,z_n-\omega_n)\frac{1}{\prod_{j=i}^n\delta_{ij}^{d_i}} \rho_i(\frac{\omega_i}{\delta_{ii}},\dots,\frac{\omega_n}{\delta_{in}}) \, d\omega.
\end{split}
\end{equation}
\textcolor{black}{Roughly speaking, we have mollified any component $F_i$ by convolution in space with a mollifier with multi-scaled dilations}.
Then, standard results on mollifier theory and our current assumptions on $F$ show us that the following controls hold
\begin{align}
\label{Proof:Controls_on_flows_mollifier} \vert  F_i(u,z) - F^{\delta}_i(u,z)\vert \, &\le \,C\sum_{j=i}^n\delta_{ij}^{\beta^j}, \\
\label{Proof:Controls_on_flows_mollifier1} | F^{\delta}_i(u,z) - F^{\delta}_i(u,z')| \, &\le \,C \sum_{j=i}^{n}\delta_{ij}^{\beta^j-1}|(z-z')_{j}\vert.
\end{align}
We can now pick $\delta_{ij}$ for any $i\le j$ in $\llbracket 2,n\rrbracket$ in order to have any contribution associated with the mollification
appearing in \eqref{Proof:Controls_on_flows_mollifier} at a good current scale time. Namely, we would like $\delta_{ij}$ to satisfy
\begin{equation}\label{proof:bilip_control}
    \left| \mathbb{T}^{-1}_{s-t}\left(F(u,z)-F^\delta(u,z)\right) \right| \,
\le \, C(s-t)^{-1},
\end{equation}
for any $u$ in $[t,s]$ and any $z$ in $\R^N$.
Using the mollifier controls \eqref{Proof:Controls_on_flows_mollifier}, it is enough to ask for
\begin{equation}
\sum_{i=2}^n(s-t)^{-\frac{1+\alpha(i-1)}{\alpha}}\sum_{j=i}^n\delta_{ij}^{\beta^j} \, \le \, C(s-t)^{-1}.\label{proof:SUM_COND_DELTA}
\end{equation}
This is true if we fix for example,
\begin{equation}\label{Proof:Controls_on_Flows_Choice_delta}
\delta_{ij} \, = \, (s-t)^{\frac{1+\alpha(i-2)}{\alpha\beta^j}} \quad \text{for $i\le j$ in $\llbracket 2,n\rrbracket$.}
\end{equation}

Next, we would like to show that, for our choice of the regularization parameter $\delta_{ij}$, the mollified drift $F^\delta$ satisfies an \emph{approximate} Lipschitz condition with a constant that, once the drift is integrated, does not yield any additional singularity. Namely, we want to derive the following control:
\begin{equation}
\label{approximateLippourF}
\left|\mathbb{T}_{s-t}^{-1}\left(F^\delta(u,z)-F^\delta(u,z')\right)\right|
\, \le \, C\left[(s-t)^{-\frac 1\alpha}+ (s-t)^{-1}|\mathbb{T}_{s-t}^{-1}(z-z')| \right] .
\end{equation}
To show it, we start noticing that $F_1 $ is H\"older continuous with H\"older index $\beta^1>0$. By Young inequality, it then yields that there exists a positive constant $C$ possibly depending on $\beta^1$ such that
$|z|^{\beta^1} \le C(1+|z|)$ for any $z$ in $\R^N$. It then follows from Equation \eqref{Proof:Controls_on_flows_mollifier1} that
\[\begin{split}
|\mathbb{T}_{s-t}^{-1}\bigl(F^\delta(u,z)-&F^\delta(u,z
')\bigr)| \\
&\le \,
C\Bigl[(s-t)^{-\frac{1}{\alpha}}(1+|(z-z')
|) + \sum_{i=2}^n\sum_{j=i}^n(s-t)^{-\frac{1+\alpha(i-1)}{\alpha}}\delta_{ij}^{\beta^j-1}|(z-z
')_j| \Bigr]\\
&\le\,  C\biggl[(s-t)^{-\frac{1}{\alpha}}+
|\mathbb{T}_{s-t}^{-1}(z-z')|\bigl(1+\sum_{i=2}^
n\sum_{j=i}^n \frac{(s-t)^{j-i}}{\delta_{ij}^{1-\beta^j}}
\bigr)\biggr].
\end{split}\]
Hence, Control \eqref{approximateLippourF} follows from the fact that, from our previous choice of $\delta_{ij}$, one gets
\begin{equation}\label{equilibri}
\frac{(s-t)^{j-i}}{\delta_{ij}^{1-\beta^j}} \, = \, (s-t)^{(j-i)-\frac{1+\alpha(i-2)}{\alpha\beta^j}(1-\beta^j)}\, \le \,  C (s-t)^{-1},
\end{equation}
\textcolor{black}{recalling that we assumed $s-t$ to be small enough and} since from the assumption \eqref{eq:thresholds_beta} on the indexes of H\"older continuity $\beta^j$ for $F$:
\[\beta^j>\frac{1+\alpha(i-2)}{1+\alpha(j-1)} \Leftrightarrow (j-i)-\frac{1+\alpha(i-2)}{\alpha\beta^j}(1-\beta^j) \, > \, -1.\]
We recall that the above inequality should precisely give the natural threshold, namely an exponent $\beta_i^j $ satisfying this condition. The current choice for $\beta^j $
is sufficient to ensure this bound holds for any $i\le j $ and is \textit{sharp} for $i=j$.
We can finally show the bound for the second remainder $\mathcal{I}^2_{s,t}(x,y)$ as given in \eqref{eq:proof_bilip_control_flow2}. It holds that:
\[\begin{split}
|\mathcal{I}^2_{s,t}(x,y)| \, &\le \, \int_{t}^{s}  \left|\mathbb{T}_{s-t}^{-1} \left(F(u,\theta_{u,s}(y))-F(u,\check{\theta}_{u,t}(x))\right)\right| \, du\\
&\le \, \int_t^s \left|\mathbb{T}_{s-t}^{-1}(F(u,\check{\theta}_{u,t}(x))-F^\delta(u,\check{\theta}_{u,t}(x)))\right|\, du \\
&\qquad\qquad\qquad + \int_t^s \left|\mathbb{T}_{s-t}^{-1}(F^\delta(u,\check{\theta}_{u,t}(x))-F^\delta(u,\theta_{u,s}(y)))\right| \, du \\
& \qquad\qquad\qquad\qquad\qquad + \int_t^s \left|\mathbb{T}_{s-t}^{-1}\left(F^\delta(u,\theta_{u,s}(y))-F(u,\theta_{u,s}(y))\right)\right| \, du\\
&=: \,  \mathcal{I}^{21}_{s,t}(x,y)+\mathcal{I}^{22}_{s,t}(x,y)+\mathcal{I}^{23}_{s,t}(x,y).
\end{split}\]
From Control \eqref{Proof:Controls_on_flows_mollifier} with our choice of $\delta_{ij}$,  we easily obtain \textcolor{black}{from Control \eqref{proof:bilip_control} that} there exists a positive constant $C:=C(T)$ such that
\begin{equation}
\label{eq:proof_bilip_control_flow5}
|\mathcal{I}^{21}_{s,t}(x,y)|+|\mathcal{I}^{23}_{s,t}(x,y)| \, \le \, C,
\end{equation}
for any $t\le s$ in $[0,T]$ and $x,y$ in $\R^N$. On the other hand, we exploit \eqref{approximateLippourF} to derive that
\[
|\mathcal{I}^{22}_{s,t}(x,y)| \, \le \, C\left[1+\int_t^s (s-t)^{-1}|\mathbb{T}_{s-t}^{-1}(\check{\theta}_{u,t}(x)-\theta_{u,s}(y))| \,    du\right]\]
for any $t\le s$ in $[0,T]$ and $x,y$ in $\R^N$.
To conclude, we finally derive \eqref{eq:proof_bilip_control_flow2} from the last inequality together with Controls \eqref{eq:proof_bilip_control_flow6}-\eqref{eq:proof_bilip_control_flow5}.

\paragraph{Proof of Lemma \ref{convergence_dirac} (Dirac convergence of frozen density).}
Fixed $(t,x)$ in $[0,T]\times \R^N$ and a bounded, continuous function $f\colon \R^N\to \R$, we want to show that the following limit
\[
\lim_{\epsilon \to 0}\left| \int_{\R^N} f(y) \tilde{p}^{t+\epsilon,y}(t,t+\epsilon,x,y)\, dy
-f(x) \right| \, = \, 0
\]
holds, uniformly in $t \in [0,T]$.\newline
We start rewriting the argument of the limit in the following way:
\begin{align}
\label{proof:control_deviation}
\int_{\R^N} f(y) \tilde{p}^{t+\epsilon,y}(t,t+\epsilon,\,&x,y)\, dy
-f(x) \\\notag
&= \, \int_{\R^N} f(y) \left[\tilde{p}^{t+\epsilon,y}(t,t+\epsilon,x,y)-\tilde{p}^{t,x}(t,t+\epsilon,x,y)\right] \, dy\\\notag
&\qquad\quad\qquad\qquad\qquad+ \int_{\R^N} f(y) \tilde{p}^{t,x}(t,t+\epsilon,x,y)\, dy-f(x).
\end{align}
By Proposition \ref{prop:Decomposition_Process_X}, we know that the second term in \eqref{proof:control_deviation} tends to zero, uniformly in $t$ in $[0,T]$ (scaling property of the upper bound for the density), when $\epsilon$ goes to zero. We can then focus on the first one. We start splitting the space $\R^N$ in the diagonal/off-diagonal regime associated with our anisotropic dynamics. Namely, we fix $\beta>0$ to be chosen later and we consider the following subsets:
\begin{align*}
   D_1 \, &:= \, \{y \in \R^N\colon \left|\mathbb{T}^{-1}_{\epsilon}(y-\theta_{t+\epsilon,t}(x))\right| \, \le \, \epsilon^{-\beta}\}; \\
    D_2 \, &:= \, \{y \in \R^N\colon \left|\mathbb{T}^{-1}_{\epsilon}(y-\theta_{t+\epsilon,t}(x))\right| \, > \, \epsilon^{-\beta}\},
\end{align*}
where $\mathbb{T}_\epsilon$ was defined in \eqref{eq:def_matrix_T}.
We can then decompose the first term in \eqref{proof:control_deviation} in the following way:
\begin{align}\notag
\Bigr{|}\int_{\R^N} f(y) \bigl[\tilde{p}^{t+\epsilon,y}(t,t+\epsilon,\,&x,y)-\tilde{p}^{t,x}(t,t+\epsilon,x,y)\bigr] \, dy\Bigr{|} \\
&\le \, \Vert f \Vert_\infty \int_{D_1}\left|\tilde{p}^{t+\epsilon,y}(t,t+\epsilon,x,y)-\tilde{p}^{t,x}(t,t+\epsilon,x,y)\right| dy \notag\\
&\qquad \qquad +\Vert f \Vert_\infty\int_{D_2}\left|\tilde{p}^{t+\epsilon,y}(t,t+\epsilon,x,y)-\tilde{p}^{t,x}(t,t+\epsilon,x,y)\right| dy \notag\\
&=:\, \Vert f \Vert_\infty\left(\mathcal{D}_1+\mathcal{D}_2\right)(t,t+\epsilon,x).\label{proof:control_deviation1}
\end{align}
We will follow different approaches to control the two terms $\mathcal{D}_1$, $\mathcal{D}_2$. In the off-diagonal regime $D_2$, the idea is to exploit tail estimates of the single densities while in the diagonal one $D_1$, a more thorough sensibility analysis between the spectral measures and the Fourier transform is needed. Let us consider first the off-diagonal term $\mathcal{D}_2$. We can write that
\begin{align*}
\mathcal{D}_2(t,t+\epsilon,x,y) \, &\le \, \int_{D_2}\left|\tilde{p}^{t+\epsilon,y}(t,t+\epsilon,x,y)\right| + \left|\tilde{p}^{t,x}(t,t+\epsilon,x,y)\right| \, dy\\
&\le \int_{D_2}\frac{1}{\det\mathbb{T}_{\epsilon}}\Big( {\bar p}(1,\mathbb{T}_{\epsilon}^{-1}(x-\theta_{t,t+\epsilon}(y))) + {\bar p}(1,\mathbb{T}_{\epsilon}^{-1}(\theta_{t+\epsilon,t}(x)-y))\Big)dy
\end{align*}
using Proposition \ref{prop:Smoothing_effect} together with Lemma \ref{lemma:identification_theta_m} for the last inequality.
From Lemma \ref{lemma:bilip_control_flow} (to use the \textit{approximate} Lipschitz property of the flows) and introducing 
\[\bar D_2 \, := \, \{y \in \R^N\colon \left|\mathbb{T}^{-1}_{\epsilon}(\theta_{t,t+\epsilon}(y)-x)\right| \, > \, \frac 12 \epsilon^{-\beta}\},\]
we thus deduce that for $\epsilon $ small enough we get:
\begin{multline*}
\mathcal{D}_2(t,t+\epsilon,x,y) \, \le \,  \int_{\bar D_2}\frac{1}{\det\mathbb{T}_{\epsilon}} {\bar p}(1,\mathbb{T}_{\epsilon}^{-1}(x-\theta_{t,t+\epsilon}(y)))\, dy\\
+\int_{D_2} \frac{1}{\det\mathbb{T}_{\epsilon}}{\bar p}(1,\mathbb{T}_{\epsilon}^{-1}(\theta_{t+\epsilon,t}(x)-y))\,dy.
\end{multline*}
Using now Equation \eqref{eq:smoothing_effect_frozen_y_GENERIC_FUNCTION} from Corollary \ref{coroll:Smoothing_effect}  for the first integral and the direct change of variable $z=\mathbb{T}^{-1}_\epsilon (y-\theta_{t+\epsilon,t}(x))$ for the second,  we can conclude that
\[\mathcal{D}_2(t,t+\epsilon,x) \,
\le \, C\int_{\R^N}\mathds{1}_{B^c(0,\frac 12\epsilon^{-\beta})}(z)(\check p+\overline{p})(1,z) \, dz,\]
where  $\check p $ is a density enjoying the same integrability properties as $\overline{p}$.
\newline
By dominated convergence theorem, it is easy to notice that $\mathcal{D}_2(t,t+\epsilon,x)$ tends to zero if $\epsilon$ goes to zero, \textcolor{black}{uniformly in the time variable $t$ in } $[0,T]$.\newline
We can now focus on the diagonal term $\mathcal{D}_1$ appearing in \eqref{proof:control_deviation1}. We start recalling from Equation \eqref{eq:definition_density_q} that the density $\tilde{p}^{\omega}(t,s,x,y)$ (for $\omega \in\{(t,x),(t+\epsilon,y)\}$) can be written as
\[
\tilde{p}^\omega(t,t+\epsilon,x,y) \, = \, \frac{\det \mathbb{M}^{-1}_\epsilon}{(2\pi)^N}\int_{\R^N}e^{\mathcal{F}_{t,t+\epsilon}(z,\omega)} \text{exp}\left(-i\langle \mathbb{M}^{-1}_\epsilon(y-\tilde{m}^{\omega}_{t+\epsilon,t}(x)),z\rangle\right) \, dz,
\]
where we have denoted:
\[\mathcal{F}_\epsilon(t,z,\omega) \, := \, \epsilon \int_{0}^{1}\int_{\R^d}\left[\cos\left(\langle z, \widehat{\mathcal{R}}_{v}B\tilde{\sigma}^{\omega}_{u(v)} p\rangle \right)-1\right] \,\nu(dp)dv,\]
with $u(v)=t+\epsilon v$ \textcolor{black}{(cf. notations in \eqref{proof:ref_notations} of Proposition \ref{prop:Decomposition_Process_X})} and $\Phi(p)$ the L\'evy symbol of the process $\{Z_t\}_{t\ge0}$.
We can now consider the two following terms
\begin{align*}
&\mathcal{P}_1(t,t+\epsilon,x,y) \, := \, \frac{\det \mathbb{M}^{-1}_\epsilon}{(2\pi)^N}\int_{\R^N}\left[e^{\mathcal{F}_\epsilon(t,z,t,x)}-e^{\mathcal{F}_\epsilon(t,z,t+\epsilon,y)}\right]e^{-i\langle \mathbb{M}^{-1}_\epsilon(y-\tilde{m}^{t,x}_{t+\epsilon,t}(x)),z\rangle} \, dz\\
&\mathcal{P}_2(t,t+\epsilon,x,y) \\
&\qquad\quad:= \,
\frac{\det \mathbb{M}^{-1}_\epsilon}{(2\pi)^N}\int_{\R^N}e^{\mathcal{F}_\epsilon(t,z,t+\epsilon,y)}\left[e^{-i\langle \mathbb{M}^{-1}_\epsilon(y-\tilde{m}^{t,x}_{t+\epsilon,t}(x)),z\rangle}-e^{-i\langle \mathbb{M}^{-1}_\epsilon(y-\tilde{m}^{t+\epsilon,y}_{t+\epsilon,t}(x)),z\rangle}\right] \, dz \notag
\end{align*}
and decompose $\mathcal{D}_1$ as follows:
\[\mathcal{D}_1 \, = \, \int_{D_1}\left|\mathcal{P}_1(t,t+\epsilon,x,y)\right|+\left|\mathcal{P}_2(t,t+\epsilon,x,y)\right| \, dy.\]
To control the first term $\mathcal{P}_1$, we can exploit a Taylor expansion. Indeed,
\begin{multline*}
   |\mathcal{P}_1(t,t+\epsilon,x,y)| \\
   \le \, \frac{C}{\det \mathbb{M}_\epsilon}\int_{\R^N}\int_0^1\left|\mathcal{F}_\epsilon(t,z,t+\epsilon,y)-\mathcal{F}_\epsilon(t,z,t,x)\right| e^{\lambda\mathcal{F}_\epsilon(t,z,t+\epsilon,y)+(1-\lambda)\mathcal{F}_\epsilon(t,z,t,x)} \, d\lambda dz. 
\end{multline*}
We then notice from \eqref{eq:control_Levy_symbol_S} that
\[\mathcal{F}_\epsilon(t,z,\omega) \, \le \, C\epsilon [1-|z|^\alpha],\]
and thus, we obtain that
\[e^{\lambda\mathcal{F}_\epsilon(t,z,t+\epsilon,y)+(1-\lambda)\mathcal{F}_\epsilon(t,z,t,x)} \, \le \, e^{C\epsilon(1-|z|^\alpha)},\]
for some constant $C$ independent from $\lambda$ in $[0,1]$.
\textcolor{black}{From our non-degenerate structure, any linear combination of the symbols remains homogeneous to a non-degenerate symbol.} Thus, we have that
\begin{equation}
\label{proof:control_deviation2}
    |\mathcal{P}_1(t,t+\epsilon,x,y)| \, \le \, \frac{C}{\det \mathbb{M}_\epsilon}\int_{\R^N}\left|\mathcal{F}_\epsilon(t,z,t+\epsilon,y)-\mathcal{F}_\epsilon(t,z,t,x)\right| e^{C\epsilon(1-|z|^\alpha)} \, dz.
\end{equation}
On the other hand,
we can decompose the difference in absolute value in the following way:
\begin{align}
\notag
|\mathcal{F}_\epsilon(t,z,t+\epsilon,y)-&\mathcal{F}_\epsilon(t,z,t,x)| \\ \notag
&\le \, \epsilon \int_{0}^{1}\Bigl{|}
\int_{\R^d}\left[\cos\left(\langle z, \widehat{\mathcal{R}}_{v}B\tilde{\sigma}^{t+\epsilon,y}_{u(v)} p\rangle\right)-\cos\left(\langle z,  \widehat{\mathcal{R}}_{v}B\tilde{\sigma}^{t,x}_{u(v)} p\rangle\right)\right] \,\nu(dp)\Bigr{|}dv \\ 
&\le \,  \epsilon \int_{0}^{1}\Bigl{|} \left(\Delta^{t,\epsilon,x,y}_s+\Delta^{t,\epsilon,x,y}_l\right)(v,z)\Bigl{|}\, dv,
\label{proof:A6}
\end{align}
where we denoted
\begin{align*}
\Delta^{t,\epsilon,x,y}_s(v,z)\, &= \, \int_{B(0,r_0)}\left[\cos\left(\langle z, \widehat{\mathcal{R}}_{v}B\tilde{\sigma}^{t+\epsilon,y}_{u(v)} p\rangle\right)-\cos\left(\langle z,  \widehat{\mathcal{R}}_{v}B\tilde{\sigma}^{t,x}_{u(v)} p\rangle\right)\right]
Q(p)\, \nu_\alpha(dy); \\
\Delta^{t,\epsilon,x,y}_l(v,z) \, &= \, \int_{B^c(0,r_0)}\left[\cos\left(\langle z, \widehat{\mathcal{R}}_{v}B\tilde{\sigma}^{t+\epsilon,y}_{u(v)} p\rangle\right)-\cos\left(\langle z,  \widehat{\mathcal{R}}_{v}B\tilde{\sigma}^{t,x}_{u(v)} p\rangle\right)\right]
Q(p)\, \nu_\alpha(dp),
\end{align*}
with $r_0$ defined in assumption [\textbf{ND}]. The term $\Delta^{t,\epsilon,x,y}_l$ involving the large jumps can be easily controlled using that $\sup_{p\in \R^d}Q(p)<\infty$:
\begin{align} \notag
|\Delta^{t,\epsilon,x,y}_l(v,z)| \, &\le \, \int_{B^c(0,r_0)}\left|\cos\left(\langle z, \widehat{\mathcal{R}}_{v}B\tilde{\sigma}^{t+\epsilon,y}_{u(v)} p\rangle \right)-\cos\left(\langle z,  \widehat{\mathcal{R}}_{v}B\tilde{\sigma}^{t,x}_{u(v)} p\rangle \right)\right|
\,\nu_\alpha(dp) \\
&\le \, C. \label{proof:A5}
\end{align}
To bound the term $\Delta^{t,\epsilon,x,y}_s$ associated with the small jumps, we want to exploit instead that $Q$ is Lipschitz continuous on $B(0,r_0)$. For this reason, we write that
\begin{align}
\notag
|\Delta^{t,\epsilon,x,y}_s&(v,z)| \\ \notag
&\le \, \left| \int_{B(0,r_0)}\left[\cos\left(\langle z, \widehat{\mathcal{R}}_{v}B\tilde{\sigma}^{t+\epsilon,y}_{u(v)} p\rangle \right)-\cos\left(\langle z,  \widehat{\mathcal{R}}_{v}B\tilde{\sigma}^{t,x}_{u(v)} p\rangle \right)\right][Q(p)-Q(0)]
\,\nu_\alpha(dp) \right|\\\notag
&\qquad\qquad\,\,+ \left| \int_{B(0,r_0)}\left[\cos\left(\langle z, \widehat{\mathcal{R}}_{v}B\tilde{\sigma}^{t+\epsilon,y}_{u(v)} p\rangle\right)-\cos\left(\langle z,  \widehat{\mathcal{R}}_{v}B\tilde{\sigma}^{t,x}_{u(v)} p\rangle \right)\right]Q(0)
\,\nu_\alpha(dp) \right| \\\label{proof:A1}
&=:\, \left(\Delta^{t,\epsilon,x,y}_{s,1}+\Delta^{t,\epsilon,x,y}_{s,2}\right)(v,z).
\end{align}
Since $Q$ and the cosine function are Lipschitz continuous in a neighborhood of $0$, we have that
\begin{equation}
\label{proof:A2}
\begin{split}
 \Delta^{t,\epsilon,x,y}_{s,1}(v,z)\, &\le \, C \int_{B(0,r_0)}|p||z|\left|\widehat{\mathcal{R}}_{v}B\tilde{\sigma}^{t+\epsilon,y}_{u(v)}p -\widehat{\mathcal{R}}_{v}B\tilde{\sigma}^{t,x}_{u(v)}p\right|
\,\nu_\alpha(dp) \\
&\le \, C \int_{B(0,r_0)}|p||z|\left|\sigma(u(v),\theta_{u(v),t+\epsilon}(y))p - \sigma(u(v),\theta_{u(v),t}(x))p\right|
\,\nu_\alpha(dp)\\
&\le \, C|z| \int_{B(0,r_0)}|p|^2 \,\nu_\alpha(dp) \, \le \, C|z|,
\end{split}
\end{equation}
where in the last step, we used that the diffusion coefficient $\sigma$ is bounded (cf. assumption [\textbf{UE}]).\newline
The control of the other term
$\Delta^{t,\epsilon,x,y}_{s,2}$ now follows from the classical characterization of the L\'evy symbol of a non-degenerate $\alpha$-stable process (see e.g. \cite{book:Sato99}). Indeed,
\[
\begin{split}
\Delta^{t,\epsilon,x,y}_{s,2}(v,z) \, &= \,\left| \int_{\R^d}\left[\cos\left(\langle z,\widehat{\mathcal{R}}_{v}B\tilde{\sigma}^{t+\epsilon,y}_{u(v)}p\rangle\right)-1\right]-\left[\cos\left(\langle z,\widehat{\mathcal{R}}_{v}B\tilde{\sigma}^{t,x}_{u(v)}p\rangle\right)-1\right] \, \nu(dp)\right| \\
&\le \, C\int_{\mathbb{S}^{d-1}}\left|\left|\langle z,\widehat{\mathcal{R}}_{v}B\tilde{\sigma}^{t+\epsilon,y}_{u(v)}s\rangle\right|^\alpha-\left|\langle z,\widehat{\mathcal{R}}_{v}B\tilde{\sigma}^{t,x}_{u(v)}s\rangle\right|^\alpha\right| \, \mu(ds).
\end{split}
\]
We now exploit the $\beta^1$-H\"older regularity in space of the diffusion coefficient $\sigma$ to show that
\begin{equation}
\label{proof:A3}
\begin{split}
    \Delta^{t,\epsilon,x,y}_{s,2}(v,z) \,&\le \, C
|z|^\alpha \left|\theta_{u(v),t+\epsilon}(y)-\theta_{u(v),t}(x)\right|^{\beta^1(\alpha\wedge 1)}  \\
&\le \, C |z|^\alpha \left[|y-\theta_{t+\epsilon,t}(x)|^{\beta^1}+\epsilon^{\beta^1}\right],
\end{split}
\end{equation}
where in the last step we used that $\alpha>1$ and the approximate Lipschitz property of the flow (cf. Lemma \ref{lemma:bilip_control_flow} \textcolor{black}{up to a normalization}, see also Lemma $1.1$ in \cite{Menozzi:Pesce:Zhang21}).\newline
We can now use Controls \eqref{proof:A2}-\eqref{proof:A3} in Equation \eqref{proof:A1} to show that
\begin{equation}
\label{proof:A4}
|\Delta^{t,\epsilon,x,y}_s(v,z)| \, \le \, C\left(|z|+    |z|^\alpha+ |y-\theta_{t+\epsilon,t}(x)|^{\beta^1}|z|^\alpha\right).
\end{equation}
Similarly, Controls \eqref{proof:A4}-\eqref{proof:A5} with Equation \eqref{proof:A6} allow us to conclude that
\begin{equation}
|\mathcal{F}_\epsilon(t,z,t+\epsilon,y)-\mathcal{F}_\epsilon(t,z,t,x)| \, \le \, C \epsilon \left(1+|z|+   \epsilon^{\beta^1} |z|^\alpha+ |y-\theta_{t+\epsilon,t}(x)|^{\beta^1}|z|^\alpha\right).
\end{equation}
We can now go back to Equation \eqref{proof:control_deviation2}. Changing variable and integrating over $z$, we find that
\[
\begin{split}
 |\mathcal{P}_1(t,t+\epsilon,x,y)| &\, \le \, \frac{C\epsilon}{\det \mathbb{M}_\epsilon}\int_{\R^N}\left(1+|z|+    \epsilon^{\beta^1}|z|^\alpha+ |y-\theta_{t+\epsilon,t}(x)|^{\beta^1}|z|^\alpha\right) e^{C\epsilon(1-|z|^\alpha)} \, dz \\
 &\le \, \frac{C}{\det \mathbb{T}_\epsilon}\int_{\R^N}\left(\epsilon+\epsilon^{\frac{\alpha-1}{\alpha}}|\tilde{z}|+  \epsilon^{\beta^1}  |\tilde{z}|^\alpha+ |y-\theta_{t+\epsilon,t}(x)|^{\beta^1}|\tilde{z}|^\alpha\right) e^{C(1-|\tilde{z}|^\alpha)} \, dz\\
&\le \, \frac{C}{\det \mathbb{T}_\epsilon}\left(\epsilon^{(1-\frac 1\alpha)\wedge \beta^1} +|y-\theta_{t+\epsilon,t}(x)|^{\beta^1}\right)
\end{split}
\]
where we recall that $\mathbb{T}_t=t^{1/\alpha}\mathbb{M}_t$.\newline
To conclude, we apply the change of variable $\tilde{y}=y-\theta_{t+\epsilon,t}(x)$:
\[
\begin{split}
\int_{D_1}|\mathcal{P}_1(t,t+\epsilon,x,y)| \, dy \, &\le \,\frac{C}{\det \mathbb{T}_\epsilon} \int_{D_1}   \left[|y-\theta_{t+\epsilon,t}(x)|^{\beta^1}+\epsilon^{(1-\frac 1\alpha)\wedge \beta^1}\right] \, dy \\
&= \, C\int_{|\tilde{y}|\le \epsilon^{-\beta}}   \left[|\mathbb{T}_\epsilon \tilde{y}|^{\beta^1}+\epsilon^{(1-\frac 1\alpha)\wedge \beta^1}\right] \, d\tilde{y} \\
&\le \, C[\epsilon^{\beta^1/\alpha-\beta(N+\beta^1)}+\epsilon^{(1-\frac 1\alpha)\wedge \beta^1-\beta N}].
\end{split}
\]
The above control then tends to zero letting $\epsilon$ go to zero, if we choose $\beta$ such that
\[0\,<\, \beta \,<\, \frac{\beta^1}{\alpha(N+\beta^1)}\wedge \frac{(1-\frac 1\alpha)\wedge \beta^1}{N}. \]
To control the second term $\mathcal{P}_2$, we use again Control \eqref{eq:control_Levy_symbol_S} and a Taylor expansion to write, similarly to above, that
\begin{align}\notag
|\mathcal{P}_2(t,&t+\epsilon,x,y) | \\ \notag
&\le \,
\frac{C}{\det \mathbb{M}_\epsilon}\int_{\R^N}e^{C\epsilon (1-|z|^\alpha)}\left|\langle \mathbb{M}^{-1}_\epsilon(y-\tilde{m}^{t,x}_{t+\epsilon,t}(x)),z\rangle-\langle \mathbb{M}^{-1}_\epsilon(y-\tilde{m}^{t+\epsilon,y}_{t+\epsilon,t}(x)),z\rangle \right| \, dz   \\
&\le \, \frac{C}{\det \mathbb{T}_\epsilon}\left|\mathbb{T}^{-1}_\epsilon\left(\theta_{t+\epsilon,t}(x)-\tilde{m}^{t+\epsilon,y}_{t+\epsilon,t}(x)\right)\right|,\label{proof:control_P2}
\end{align}
where in the last passage we used Lemma \ref{lemma:identification_theta_m}. To bound the above right-hand side, we now exploit Corollary \ref{coroll:control_error_flow} to show that
\[
\begin{split}
|\mathcal{P}_2(t,t+\epsilon,x,y) | \,
\le \, C \epsilon^{\frac{1}{\alpha} \wedge \zeta}\frac{1}{\det \mathbb{T}_\epsilon}\left(1+|\mathbb{T}_{\epsilon}^{-1}(\theta_{t+\epsilon,t}(x)-y)|\right).
\end{split}
\]
Similarly to above, we can then apply a change of variables:
\[
\int_{D_1}|\mathcal{P}_2(t,t+\epsilon,x,y) | \, dy \, \le \, C \epsilon^{\frac{1}{\alpha} \wedge \zeta}\int_{|z|\le\epsilon^{-\beta}}\left(1+|z| \right) \, dz.
\]
We can then notice again that the above control tends to zero letting $\epsilon$ goes to zero, if we choose $\beta$ small enough.

\paragraph{Proof of Lemma \ref{prop:convergence_LpLq}.}
As in the previous Lemma \ref{convergence_dirac}, we want to show the following limit:
\[
\lim_{\epsilon \to 0} \Vert
I_\epsilon f -f \Vert_{L^{p}_tL^q_x} \, = \, 0,
\]
for some $p\in (1,+\infty)$, $q\in (1,+\infty)$ and $f$ in $C_c^{1,2}([0,T)\times \R^N)$.
We start writing that
\[
\|I_\epsilon f-f\|^p_{L^{p}_tL^q_x} \, = \, \int_0^T \Vert I_\epsilon f(t,\cdot)-f(t,\cdot)\Vert_{L^q}^{p}\, dt.\]
We then notice that, up to a middle point-type argument, the indicator function in the definition \eqref{eq:def_f_epsilon} of $I_\epsilon f$ can be easily controlled. We can now write that
\begin{align}
\notag
\Vert I_\epsilon f(t,\cdot)-f(&t,\cdot)\Vert_{L^p}^{p}\, = \, \int_{\R^N} \left|\int_{\R^N}f(t+\epsilon,y) \tilde{p}^{t+\epsilon,y}(t,t+\epsilon,x,y) \, dy - f(t,x)\right|^{p}  dx \\ \notag
&\le C \Bigl(\int_{\R^N} \left|\int_{\R^N}f(t+\epsilon,y) \tilde{p}^{t+\epsilon,y}(t,t+\epsilon,x,y) \, dy - f(t+\epsilon,\theta_{t+\epsilon,t}(x))\right|^{p} dx\\ \notag
&\,\,\qquad\qquad\qquad\qquad\qquad\qquad\qquad + \int_{\R^N} |f(t+\epsilon,\theta_{t+\epsilon,t}(x))-f(t,x)|^p \, dx\Bigr)\\
&=:\, C\left(\mathcal{I}+\mathcal{I}'\right)(\epsilon,t).\label{proof:eq:LpLq_convergence}
\end{align}
Since $f$ is smooth  and with compact support in time and space, it follows immediately that $\mathcal{I}'(\epsilon,t)$ tends to zero if $\epsilon$ goes to zero, \textcolor{black}{thanks to the bounded convergence Theorem}.\newline
We can then focus on the first term $\mathcal{I}(\epsilon,t) $. We start splitting it in the following way:
\[
\begin{split}
\mathcal{I}(\epsilon,t) \, &\le \,  C\Bigl(\int_{\R^N}\left|\int_{\R^N}\left[f(t+\epsilon,y)-f(t+\epsilon,\theta_{t+\epsilon,t}(x))\right] \tilde{p}^{t+\epsilon,y}(t,t+\epsilon,x,y) \,  dy \right|^p dx\\
&\quad + \int_{\R^N}\left|f(t+\epsilon,\theta_{t+\epsilon,t}(x)) \int_{\R^N} \left[\tilde{p}^{t+\epsilon,y}(t,t+\epsilon,x,y)- \tilde{p}^{t,x}(t,t+\epsilon,x,y)\right]dy\right|^p dx\Bigr)\\
&=: \, C\left(\mathcal{I}_1+\mathcal{I}_2\right)(\epsilon,t),
\end{split}
\]
where we used that $\tilde{p}^{t,x}(t,s,x,y)$ is indeed a \textit{true} density with respect to $y$. The second term $\mathcal{I}_2(\epsilon,t)$ already appeared in the proof of Lemma \ref{convergence_dirac} (Dirac Convergence of frozen density) \textcolor{black}{(cf. term $\mathcal{D}_2$ in \eqref{proof:control_deviation1})} and a similar  analysis readily gives that $\mathcal{I}_2(\epsilon, t)\overset{\epsilon\rightarrow 0}{\longrightarrow}0 $.\newline
To control instead the first term $\mathcal{I}_1(\epsilon,t)$, we decompose the whole space $\R^N$ into $\Delta_1\cup \Delta_2$ given by
\begin{align*}
    \Delta_1 \, &:= \, \{x \in \R^N \colon |\theta_{t+\epsilon, t}(x)-\text{supp}[f(t+\epsilon,\cdot)]|\le 1\}; \\
    \Delta_2 \, &:= \,  \{x \in \R^N \colon |\theta_{t+\epsilon, t}(x)-\text{supp}[f(t+\epsilon,\cdot)]| > 1\}.
\end{align*}
Using Proposition \ref{prop:Smoothing_effect} with $(\tau,\xi)=(t+\epsilon,y)$, we write that
\[\begin{split}
\mathcal{I}_1(\epsilon,t) \, &\le \, \int_{\R^N}\left(\int_{\R^N}|f(t+\epsilon,y)-f(t+\epsilon,\theta_{t+\epsilon,t}(x))| \frac{\overline{p}\left(1,\mathbb{T}^{-1}_\epsilon(y-\tilde{m}^{t+\epsilon,y}_{t+\epsilon,t}(x))\right)}{\det \mathbb{T}_\epsilon} \,  dy  \right)^p dx \\
&\le \,\int_{\Delta_1}\left(\int_{\R^N}|f(t+\epsilon,y)-f(t+\epsilon,\theta_{t+\epsilon,t}(x))| \frac{\overline{p}\left(1,\mathbb{T}^{-1}_\epsilon(y-\tilde{m}^{t+\epsilon,y}_{t+\epsilon,t}(x))\right)}{\det \mathbb{T}_\epsilon} \,  dy  \right)^p dx \\
&\quad + \int_{\Delta_2}\left(\int_{\R^N}|f(t+\epsilon,y)-f(t+\epsilon,\theta_{t+\epsilon,t}(x))| \frac{\overline{p}\left(1,\mathbb{T}^{-1}_\epsilon(y-\tilde{m}^{t+\epsilon,y}_{t+\epsilon,t}(x))\right)}{\det \mathbb{T}_\epsilon} \,  dy  \right)^p dx\\
&=: \, \left(\mathcal{I}_{11}+\mathcal{I}_{12}\right)(\epsilon,t).
\end{split}
\]
To control $\mathcal{I}_{11}$, we start noticing that $f$ is H\"older continuous with a H\"older exponent $\gamma<\alpha$ in $(0,1]$, since it has a compact support. Moreover, $\Delta_1$ is a bounded set (uniformly in $\epsilon$). Then, from Lemma \ref{lemma:identification_theta_m} (cf.\ Equation \eqref{eq:identification_theta_m1}), Lemma \ref{lemma:bilip_control_flow} and Corollary \ref{coroll:Smoothing_effect},
\[
\begin{split}
\mathcal{I}_{11}(\epsilon,t) \, &\le \,C \int_{\Delta_1}\left(\int_{\R^N}|y-\theta_{t+\epsilon,t}(x)|^\gamma \frac{\overline{p}\left(1,\mathbb{T}^{-1}_\epsilon(y-\tilde{m}^{t+\epsilon,y}_{t+\epsilon,t}(x))\right)}{\det \mathbb{T}_\epsilon} \,  dy  \right)^p dx \\
&\le \,C \epsilon^{p\gamma/\alpha}\int_{\Delta_1}\left(\int_{\R^N}|\mathbb{T}^{-1}_\epsilon(y-\theta_{t+\epsilon,t}(x))|^\gamma \frac{\overline{p}\left(1,\mathbb{T}^{-1}_\epsilon(\theta_{t,t+\epsilon}(y)-x)\right)}{\det \mathbb{T}_\epsilon} \,  dy  \right)^p dx \\
&\le \, C \epsilon^{p\gamma/\alpha}\int_{\Delta_1}\left(\int_{\R^N}\left[|\mathbb{T}^{-1}_\epsilon(\theta_{t,t+\epsilon}(y)-x)|^\gamma+1\right] \frac{\overline{p}\left(1,\mathbb{T}^{-1}_\epsilon(\theta_{t,t+\epsilon}(y)-x)\right)}{\det \mathbb{T}_\epsilon} \,  dy  \right)^p dx \\
&\le \, C\epsilon^{p\gamma/\alpha}.
\end{split}
\]
To control instead $\mathcal{I}_{12}$ we firstly notice that if $x$ is in $\Delta_2$, then, $\theta_{t+\epsilon,t}(x)$ is not in the support of $f$. Thus,
\[\begin{split}
\mathcal{I}_{12}(\epsilon,t) \, &= \, \int_{\Delta_2}\left(\int_{\R^N}|f(t+\epsilon,y)-f(t+\epsilon,\theta_{t+\epsilon,t}(x))| \frac{\overline{p}\left(1,\mathbb{T}^{-1}_\epsilon(y-\tilde{m}^{t+\epsilon,y}_{t+\epsilon,t}(x))\right)}{\det \mathbb{T}_\epsilon} \,  dy  \right)^p dx\\
 &\le \, \int_{\Delta_2}\left(\int_{\text{supp}f}|f(t+\epsilon,y)| \frac{\overline{p}\left(1,\mathbb{T}^{-1}_\epsilon(\theta_{t,t+\epsilon}(y)-x)\right)}{\det \mathbb{T}_\epsilon} \,  dy  \right)^p dx\\
 &\le \,   \Vert f\Vert^p_\infty\int_{\Delta_2}\left(\int_{\text{supp}f} \frac{\overline{p}\left(1,\mathbb{T}^{-1}_\epsilon(\theta_{t,t+\epsilon}(y)-x)\right)}{\det \mathbb{T}_\epsilon} \,  dy\right)^{p-1+1} dx \\
  &\le \, C\int_{\text{supp}f}\int_{\Delta_2} \frac{\overline{p}\left(1,\mathbb{T}^{-1}_\epsilon(\theta_{t,t+\epsilon}(y)-x)\right)}{\det \mathbb{T}_\epsilon} \, dxdy,
\end{split}
\]
where in the last step we used that, from Corollary \ref{coroll:Smoothing_effect}
\[\left(\int_{\text{supp}f} \frac{\overline{p}\left(1,\mathbb{T}^{-1}_\epsilon(\theta_{t,t+\epsilon}(y)-x)\right)}{\det \mathbb{T}_\epsilon} \,  dy\right)^{p-1} \, \le \, \left(\int_{\R^N} \frac{\overline{p}\left(1,\mathbb{T}^{-1}_\epsilon(\theta_{t,t+\epsilon}(y)-x)\right)}{\det \mathbb{T}_\epsilon} \,  dy\right)^{p-1} \, \le  \, C_p.\]
We notice now that for any $y$ in $\text{supp}f$ and any $x$ in $\Delta_2$, we have that $|y-\theta_{t+\epsilon,t}(x))|\ge 1$. Exploiting Corollary \ref{coroll:Smoothing_effect} and Lemma \ref{lemma:bilip_control_flow}, we write that
\[\begin{split}
\mathcal{I}_{12}(\epsilon,t) \, &\le \,    \int_{\text{supp}f}\int_{\Delta_2} |y-\theta_{t+\epsilon,t}(x)| \frac{\overline{p}\left(1,\mathbb{T}^{-1}_\epsilon(\theta_{t,t+\epsilon}(y)-x)\right)}{\det \mathbb{T}_\epsilon} \, dxdy \\
&\le \, C\epsilon^{\frac{1}{\alpha}}\int_{\text{supp}f}\int_{\Delta_2}|\mathbb{T}^{-1}_\epsilon(y-\theta_{t+\epsilon,t}(x))| \frac{\overline{p}\left(1,\mathbb{T}^{-1}_\epsilon(\theta_{t,t+\epsilon}(y)-x)\right)}{\det \mathbb{T}_\epsilon} \, dxdy \\
&\le \, C\epsilon^{\frac{1}{\alpha}}\int_{\text{supp}f}\int_{\R^N}\left[|\mathbb{T}^{-1}_\epsilon(\theta_{t,t+\epsilon}(y)-x)|+1\right] \frac{\overline{p}\left(1,\mathbb{T}^{-1}_\epsilon(\theta_{t,t+\epsilon}(y)-x)\right)}{\det \mathbb{T}_\epsilon} \, dxdy \\
&\le \, C\epsilon^{\frac{1}{\alpha}}\int_{\text{supp}f}\int_{\R^N}\left[|z|+1\right] \overline{p}\left(1,z\right) \, dzdy \\
&\le \, C\epsilon^{\frac{1}{\alpha}}.
\end{split}
\]
Knowing the convergence of $\mathcal{I}(\epsilon,t)$ and $\mathcal{I}'(\epsilon,t)$ to zero, we can finally conclude the proof using the dominated convergence theorem in \eqref{proof:eq:LpLq_convergence}.

\subsection{Controls associated with the change of variable}

\subsubsection{Proof of Corollary \ref{coroll:Smoothing_effect}}
We first concentrate on the proof of Control \eqref{eq:smoothing_effect_frozen_y}. We start exploiting the decomposition of $\overline{p}(t,z)$ in terms of small and large jumps, as in \eqref{proof_decomposition_p_segnato}, to rewrite the left-hand side of Equation \eqref{eq:smoothing_effect_frozen_y} in the following way:
\begin{align*}
I(s,t,x)\, &:= \,
\int_{\R^N} \frac{|\mathbb{T}_{s-t}^{-1}(\theta_{t,s}(y)-x)|^\gamma}{\det\mathbb T_{s-t}}\bar p(1, \mathbb{ T}_{s-t}^{-1}(\theta_{t,s}(y)-x))\, dy  \\
&=\, \int_{\R^N} \frac{|\mathbb{T}_{s-t}^ {-1}(\theta_{t,s}(y)-x)|^\gamma}{\det \mathbb T_{s-t}}\int_{\R^N}p_{\overline{M}}(1, \mathbb{ T}_{s-t}^{-1}(\theta_{t,s}(y)-x)-w)\overline{P}_1(dw)dy,
\end{align*}
Then, the Fubini Theorem and the definition of $p_{\overline{M}}$ in \eqref{proof:control_pM_segnato} immediately imply that
\begin{align*}
I(s,t,x)\, &= \, \int_{\R^N} \int_{\R^N} \frac{|\mathbb{T}_{s-t}^ {-1}(\theta_{t,s}(y)-x)|^\gamma}{\det \mathbb T_{s-t}}p_{\overline{M}}(1, \mathbb{ T}_{s-t}^{-1}(\theta_{t,s}(y)-x)-w)\, dy\overline{P}_1(dw)\\
&\le \, C\int_{\R^N} \int_{\R^N}\det \mathbb T_{s-t}^{-1}\frac{ \left[|\mathbb{T}_{s-t}^ {-1}(\theta_{t,s}(y))-x-w|^\gamma+|w|^\gamma\right]}{\left[1+|\mathbb{ T}_{s-t}^{-1}(\theta_{t,s}(y)-x)-w|\right]^{N+2}}\, dy\overline{P}_1(dw).
\end{align*}
To conclude, it is now enough to show that for any $M> N+1$, there exists $C:=C(M)$ such that
\begin{equation}
\label{proof:control_smoothing2}
\int_{\R^N}  \frac{\det\mathbb{T}^{-1}_{s-t}}{\left[1+|\mathbb{T}_{s-t}^{-1}(\theta_{t,s}(y)-x) -w|\right]^M}\, dy\, \le \, C.
\end{equation}
Indeed, it would follow from Control \eqref{proof:control_smoothing2} that
\begin{align*}
    I(t,s,x) \, &\le \, C\int_{\R^N} \int_{\R^N} \frac{\det \mathbb T_{s-t}^{-1}}{\left[1+|\mathbb{ T}_{s-t}^{-1}(\theta_{t,s}(y)-x)-w|\right]^{N+2-\gamma}}\, dy\overline{P}_1(dw) \\
    &\qquad\qquad +\int_{\R^N} \int_{\R^N} \left[1+|w|^\gamma\right]\frac{\det \mathbb T_{s-t}^{-1}}{\left[1+|\mathbb{ T}_{s-t}^{-1}(\theta_{t,s}(y)-x)-w|\right]^{N+2}}\, dy\overline{P}_1(dw) \\
    &\le \,C\int_{\R^N} \left[1+|w|^\gamma\right]\overline{P}_1(dw) \, \le \, C.
\end{align*}
In order to show Control \eqref{proof:control_smoothing2}, we start noticing that it would be enough to apply the change of variable $\tilde{y}=\mathbb{T}_{s-t}^{-1}(x-\theta_{t,s}(y)) -w$ and then, to control the Jacobian matrix of the transformation. Unfortunately, our coefficients are not smooth enough in order to follow this kind of reasoning. Indeed, the drift $F$ is only H\"older continuous.\newline
As done already in in the proof of Lemma \ref{lemma:bilip_control_flow}, we firstly need to regularize $F$ through a multi-scale mollification procedure. Namely, we reintroduce the mollified drift $F^\delta:=(F^\delta_1,\dots,F^\delta_n)$ similarly to what we did in Equation \eqref{eq:multi_scale_mollification}. However we modify a bit the mollifying parameters and set
\begin{equation}\label{Proof:Controls_on_Flows_Choice_delta_LOC_LEMMA}
\delta_{ij} \, = \, \bar C(s-t)^{\frac{1+\alpha(j-2)}{\alpha\beta^j}} \quad \text{for $2\le i\le j\le n$,}
\end{equation}
for a constant $\bar C$ meant to be large enough. We also mollify the first component $ F_1$ at a macro scale, i.e. $\delta_{1j}=C_1$, with $C_1$ large enough as well.

In particular, this choice of parameters gives that the controls
\eqref{eq:proof_bilip_control_flow6}, \eqref{proof:bilip_control} and \eqref{approximateLippourF} hold again.
\newline
We can now define the mollified flow $\theta^\delta_{t,s}(y)$ associated with the drift $F^\delta$ given by
\begin{equation}\label{proof:definition_theta_delta}
\theta_{t,s}^\delta(y) \, =\, y-\int_t^s \left[A_u\theta_{u,s}^\delta+ F^\delta(u,\theta_{u,s}^\delta(y))\right] \,  du.
\end{equation}
Denoting now, for brevity, \[\Delta^\delta\theta_{u,s}(y) \, :=\, \theta_{u,s}(y)-\theta_{u,s}^\delta(y),\]
it is not difficult to check from the Gr\"onwall Lemma and Controls \eqref{eq:proof_bilip_control_flow6}, \eqref{proof:bilip_control} and \eqref{approximateLippourF} that
\begin{align} \label{SENSI_DELTA_THETA}
|\mathbb T_{s-t}^{-1}(&\theta_{t,s}(y)-\theta_{t,s}^\delta(y))| \, \le \, \left|\int_t^s \mathbb{T}_{s-t}^{-1}\left[A_u(\Delta^\delta\theta_{u,s}(y))+ F(u,\theta_{u,s}(y))-F^\delta(u,\theta_{u,s}^\delta(y))\right] \,  du\right| \notag\\
&\le \, \int_t^s \left|\mathbb{T}_{s-t}^{-1}A_u(\Delta^\delta\theta_{u,s}(y))\right|\, du+\int_t^s\left|\mathbb{T}_{s-t}^{-1}\left(F(u,\theta_{u,s}(y))-F^\delta(u,\theta_{u,s}(y))\right)\right| \,  du\notag\\
&\qquad \qquad\qquad \qquad\qquad \qquad\qquad + \int_t^s\left|\mathbb{T}_{s-t}^{-1}\left(F^\delta(u,\theta_{u,s}(y))-F^\delta(u,\theta_{u,s}^\delta(y))\right)\right| \,  du \notag\\
&\le C_0,
\end{align}
for some positive constant $C_0$. \newline
Exploiting now Control \eqref{SENSI_DELTA_THETA}, we firstly notice that for $C\ge 2C_0$,
\[
\begin{split}
C+|\mathbb{T}_{s-t}^{-1}(x-\theta_{t,s}(y)) -w|\, &\ge \, C+|\mathbb{T}_{s-t}^{-1}(x-\theta^\delta_{t,s}(y)) -w| - |\mathbb{T}_{s-t}^{-1}(\theta_{t,s}(y)-\theta_{t,s}^\delta(y))| \\
&\ge \, C_0+|\mathbb{T}_{s-t}^{-1}(x-\theta^\delta_{s,t}(y)) -w|
\end{split}\]
and we then use it to write that
\begin{align}
\int_{\R^N} \notag \frac{\det\mathbb{T}^{-1}_{s-t}}{\left(1+|\mathbb{T}_{s-t}^{-1}(x-\theta_{t,s}(y)) -w|\right)^M}\, dy\, &\le \,C \int_{\R^N}  \frac{\det\mathbb{T}^{-1}_{s-t}}{\left(1+|\mathbb{T}_{s-t}^{-1}(x-\theta^\delta_{t,s}(y)) -w|\right)^M}\, dy  \\
&= \, C\int_{\R^N}  \frac{1}{\left(1+|\tilde{y}|\right)^M}\frac{1}{\det  J_{t,s}^\delta(\tilde{y})}\, dy
\label{proof:sensibility1}
\end{align}
where in the last step we used the change of variables $\tilde{y}=\mathbb{T}_{s-t}^{-1}(x-\theta^\delta_{t,s}(y)) -w$ and denoted by $\textcolor{black}{J_{t,s}^\delta}(\tilde{y})$ the Jacobian matrix of $y\to \theta^\delta_{t,s}(y)$.\newline
It is now clear that the last term  in  \eqref{proof:sensibility1} is indeed controlled by a constant $C$, if we show the existence of a positive constant $c$, independent from $y$ in $\R^N$, $t<s$ in $[0,T]$ and $\delta$ , such that
\begin{equation}
\label{proof:sensibility3}
    |\det J_{t,s}^\delta(y)| \, \ge \, c >0.
\end{equation}
This is precisely the result provided by Lemma \ref{LEMMA_FOR_DET} below.
From the previous computations it is clear that \eqref{eq:smoothing_effect_frozen_y} holds.

Let us now turn to the proof of Control \eqref{eq:smoothing_effect_frozen_y_GENERIC_FUNCTION}. Following the previous approach, we can write
\begin{align*}
&\int_{\{|\mathbb{T}^{-1}_{s-t}(\theta_{t,s}(y)-x)|\ge K\}} \frac 1{\det \mathbb{T}_{s-t}}\overline{p}(1,\mathbb{T}^{-1}_{s-t}(\theta_{t,s}(y)-x)) \, dy \\
 &\le \, C \int_{\{|\mathbb{T}^{-1}_{s-t}(\theta_{t,s}^\delta(y)-x)|\ge K-|\mathbb{T}^{-1}_{s-t}(\Delta^\delta\theta_{u,s}(y))|\}}  \int_{\R^N}  \frac{\det\mathbb{T}^{-1}_{s-t}}{\left(1+|\mathbb{T}_{s-t}^{-1}(x-\theta^\delta_{s,t}(y)) -w|\right)^M}\overline{P}_1(dw) dy\\
 \le& \, C \int_{\{|\mathbb{T}^{-1}_{s-t}(\theta_{t,s}^\delta(y)-x)|\ge K-C_0\}}  \int_{\R^N}  \frac{\det\mathbb{T}^{-1}_{s-t}}{\left(1+|\mathbb{T}_{s-t}^{-1}(x-\theta^\delta_{s,t}(y)) -w|\right)^M}\overline{P}_1(dw) dy,
\end{align*}
exploiting also \eqref{SENSI_DELTA_THETA} for the last inequality. Using now the Fubini Theorem and the change of variables  $z=\mathbb{T}_{s-t}^{-1}(x-\theta^\delta_{s,t}(y))$, we derive from \eqref{proof:sensibility3} that
\begin{align*}
\int_{\{|{\mathbb{T}^{-1}_{s-t}(\theta_{t,s}(y)-x)}|\ge K\}} \frac 1{\det \mathbb{T}_{s-t}}&\overline{p}(1,\mathbb{T}^{-1}_{s-t}(\theta_{t,s}(y)-x)) \, dy \\
&\le \, C \int_{\R^N}   \int_{\{|z|\ge K-C_0\}} \frac{1}{\left(1+|z -w|\right)^M} \, dz \overline{P}_1(dw)\\
&=:\, C\int_{\{|z|\ge \frac K2|\}}\check p(1,z)\, dz,
\end{align*}
where $\check p $ is a density satisfying the same integrability properties as $\bar p $ assuming as well $K $ large enough.
Thus \eqref{eq:smoothing_effect_frozen_y_GENERIC_FUNCTION} holds and the proof of Corollary \ref{coroll:Smoothing_effect} is now complete.

\subsubsection{Jacobian of the mollified system}
This is a technical part dedicated to the proof of control \eqref{proof:sensibility3} appearing in the proof of key Corollary \ref{coroll:Smoothing_effect} which precisely gives the expected smoothing effect of the frozen density where the freezing parameters also correspond to the integration variable.
\begin{lemma}[Control of the determinant for the change of variable]\label{LEMMA_FOR_DET}
Let  $\theta^\delta_{t,s}(y)$ denote the mollified flow associated with the drift $F^\delta$ where the mollifying parameter $\delta$ has the form \eqref{Proof:Controls_on_Flows_Choice_delta_LOC_LEMMA}. Its dynamics writes:
\begin{equation*}
\theta_{t,s}^\delta(y) \, =\, y-\int_t^s \left[A_u\theta_{u,s}^\delta+ F^\delta(u,\theta_{u,s}^\delta(y))\right] \,  du.
\end{equation*}
Then, there exists a constants $c_0:=c_0(T)>0$ s.t., denoting for $0\le t\le s\le T  $ by $J_{t,s}^\delta (y)$ the  Jacobian matrix associated with the mapping $y\mapsto \theta_{s,t}^\delta (y) $
 $$ {\rm det}(J_{t,s}^{\delta}(y))\ge c_0 .$$
 Importantly, $c_0$ does not depend on $\delta$.
\end{lemma}
\begin{proof}
Let us first mention that the even though the coefficients $F^\delta$ are smooth, the above control is not direct because there is a subtle balance between the mollifying, matrix valued, parameter $\delta$ and the length of the considered time interval $[t,s] $. We recall indeed that the entries $\delta_{ij} $ given in \eqref{Proof:Controls_on_Flows_Choice_delta_LOC_LEMMA} do depend on $s-t$.

We also recall that, similarly to \eqref{Proof:Controls_on_flows_mollifier1}, it holds that
\begin{equation}
\label{proof:control_Derivative_F}
|D_{x_j}F_1^\delta(t,z)| \, \le \,  C (\delta_{1j})^{\beta^1-1},\ \forall 2\le i\le j\le n,\ |D_{x_j}F_i^\delta(t,z)| \, \le \,  C (\delta_{ij})^{\beta^j-1}.
\end{equation}
To prove the statement, we have thus to justify, somehow similarly to the control for the flows of Lemma \ref{lemma:bilip_control_flow}, that the explosive behavior of the Lipschitz moduli  is indeed well balanced by the time-integration.

Let us now start from the dynamics of $J^\delta (y)$ which writes:
\begin{align*}
J_{t,s}^\delta(y) = \ D_y\theta^\delta_{t,s}(y) &= \mathds{I} - \int_t^s \left[
\left(A_u +D_zF^\delta(u,z)|_{z=\theta^\delta_{u,s}(y)}\right)D_y\theta^\delta_{u,s}(y)\right] \, du\\
&= \mathds{I} - \int_t^s \left[
\left(A_u +D_zF^\delta(u,z)|_{z=\theta^\delta_{u,s}(y)}\right)J^\delta_{u,s}(y)\right] \, du.
\end{align*}
The above equation can be partially integrated using the resolvent $(R_{u,s})_{u\in [t,s]} $ associated with $A$, i.e. the $\R^N\otimes \R^N $ valued function satisfying
\begin{equation}\label{DYN_RES_STRUCT}
\frac{d}{du}R_{u,s}=A_uR_{u,s},\ R_{s,s}=I_{nd\times nd}.
\end{equation}
This yields:
\begin{equation}
\label{LIN_DYN_PARTIALLY_INTEGRATED}
J_{t,s}^\delta(y) =R_{t,s}-\int_t^s R_{t,u}D_zF^\delta(u,z)|_{z=\theta^\delta_{u,s}(y)}J^\delta_{u,s}(y)du.
\end{equation}

We actually have the following important structure property of the resolvent
$(R_{u,s})_{u\in [t,s]} $. There exists a non-degenerate family of  matrices $(\hat{R}_{\frac{u-t}{s-t}}^{t,s})_{u\in [t,s]}$, which  is bounded uniformly on $u\in [t,s] $ with constants depending on $T$ s.t.
\begin{equation}\label{scaling_relation}
R_{u,s}= \mathbb{T}_{s-t} \hat{R}_{\frac{u-t}{s-t}}^{t,s} (\mathbb{T}_{s-t})^{-1}.
\end{equation}
Indeed, setting for all $v\in [0,1],\ \hat R_v^{t,s}:=(\mathbb{T}_{s-t})^{-1}R_{t+v(s-t),s} \mathbb{T}_{s-t}$ and differentiating yields:
\begin{eqnarray*}
\partial_v \hat R_v^{t,s}&=&(s-t)(\mathbb{T}_{s-t})^{-1}A_{t+v(s-t)}R_{t+v(s-t),s} \mathbb{T}_{s-t}\\
&=&\biggl((s-t)(\mathbb{T}_{s-t})^{-1}A_{t+v(s-t)}\mathbb{T}_{s-t}\biggr)\hat R_v^{t,s}:=A_v^{t,s} \hat R_v^{t,s}.
\end{eqnarray*}
The identity \eqref{scaling_relation} then actually follows from the structure of the matrix $A_t$ (see assumption \textbf{[H]} and \eqref{eq:def_matrix_A}) which ensures that $(A_v^{t,s})_{v\in [0,1]} $ has bounded entries.

As a by-product of \eqref{scaling_relation}, we derive that there exists $C\ge 1 $ s.t. for all $(i,j)\in \llbracket 1,n \rrbracket$,
\begin{equation}\label{THE_BOUND_RES}
|(R_{t,u})_{ij}|\le C(\mathds{1}_{j\ge i}+(s-t)^{i-j}\mathds{1}_{i>j}).
\end{equation}
From \eqref{LIN_DYN_PARTIALLY_INTEGRATED} we thus derive
\begin{align}
|J_{t,s}^\delta(y)|\le& C+\int_t^s \sum_{i,j,k=1}^n \big|R_{t,u}D_zF^\delta(u,z)|_{z=\theta^\delta_{u,s}(y)}\big|_{ik} |J^\delta_{u,s}(y)|_{kj}du \notag\\
\le & C+\int_t^s \sum_{i,j,k=1}^n \big|R_{t,u} D_zF^\delta(u,z)|_{z=\theta^\delta_{u,s}(y)}\big|_{ik} |J^\delta_{u,s}(y)|_{kj}du.\label{SEMI_INT_2}
\end{align}
Remember now that $D_z F^\delta(u,z) $ is upper triangular. Then for fixed $(j,k)\in \llbracket 1,n \rrbracket^2 $, using \eqref{THE_BOUND_RES},
\begin{align}
\big|R_{t,u}D_zF^\delta(u,z)|_{z=\theta^\delta_{u,s}(y)}\big|_{ik}\le& \sum_{\ell=1}^k |R_{t,u}|_{i\ell}|DF^\delta_{\ell k}|_\infty\notag\\
\le &C \sum_{\ell=1}^k (\mathds{1}_{\ell\ge i}+(t-s)^{i-l}\mathds{1}_{\ell<i})|D_kF_\ell^{\delta}|_\infty.
\end{align}
It is now clearly seen that, for a fixed line index $i$ and $\ell\ge i $, there is no time regularity, contrarily to what happened with the control of the renormalized flows. Recall that if we had chosen
$F^\delta$ as in the proof of Lemma \ref{lemma:bilip_control_flow} then, for $\ell\ge 2 $ (recall that we regularize at macro scale $C_1$ for $F_1^\delta $) $|D_k F_\ell^\delta|_\infty \le C(\delta_{\ell k})^{-1+\beta_\ell^k}$. It is then clear that the $ \big(\delta_{lk}^{-1+\beta_\ell^k}\big)_{\ell \in \llbracket 2, k\rrbracket}$ must have the same order, which precisely prevents from the choice in \eqref{Proof:Controls_on_Flows_Choice_delta} which allows to consider minimal H\"older regularity exponents distinguishing the regularity with respect to the $k^{{\rm th}} $ variable in function of the level $\ell $ of the chain. We are here  led to consider $\beta_{\ell}^k =\beta_k^k =\beta^k$ (condition \eqref{Proof:Controls_on_Flows_Choice_delta_LOC_LEMMA}), imposing the strongest integrability threshold, associated with the diagonal perturbation at level $k$ all along the previous levels (up to the second one), which in principle lead to less singularity when the corresponding gradients are considered.

Such a phenomenon naturally appears when investigating the strong uniqueness of the SDE because of the Zvonkin approach, see e.g. \cite{Hao:Wu:Zhang19} for the Kinetic case deriving from our framework or \cite{Chaudru17} for the kinetic Brownian case. It was also the case, still for the Brownian kinetic case, in \cite{Chaudru18} where the parametrix approach freezing the initial coefficients was considered. The author had to impose the same regularity for the drift, in the degenerate variable, on the whole $F$. Hence, adapting the work \cite{Marino20} to derive pointwise bound of the gradients, which could have been another approach would have led to the same constraints. Here, we have slightly more freedom since we manage to have arbitrary smoothness indexes for the non-degenerate component of the drift.

We thus derive from \eqref{SEMI_INT_2} and for $\bar C,  C_1$ large enough there exists $c_0>0$ such that 
\[\left[ \sum_{k=2}^n\sum_{\ell=2}^k(\delta_{\ell k})^{-1+\beta^k}+\sum_{k=1}^n(\delta_{1k})^{-1+\beta_1^k}\right](s-t)\le c_0\]
meant to be small that, under the current assumptions, there exists $C\ge 1 $ s.t.
\begin{equation*}
|J_{t,s}^\delta(y)|\le C\exp(c_0),
\end{equation*}
and similarly, $\forall u\in [t,s]$,
\begin{equation}\label{CTR_BD_SUP}
|J_{u,s}^\delta(y)|\le C\exp(c_0).
\end{equation}
Rewriting:
\begin{align*}
J_{t,s}^\delta(y) =R_{t,s}\big(I-\int_t^s R_{s,u}D_zF^\delta(u,z)|_{z=\theta^\delta_{u,s}(y)}J^\delta_{u,s}(y)du\big),
\end{align*}
we derive from \eqref{THE_BOUND_RES}, \eqref{CTR_BD_SUP} that the matrix $\big(I-\int_t^s R_{s,u}D_zF^\delta(u,z)|_{z=\theta^\delta_{u,s}(y)}J^\delta_{u,s}(y)du\big)$
 has diagonal dominant and this gives, from the non degeneracy of $R$, the statement concerning the determinant.
\end{proof}

\setcounter{equation}{0}

\chapter{About the sharp
constants in Sobolev and Schauder estimates
for degenerate Kolmogorov operators}
\fancyhead[LE]{Chapter \thechapter. About the sharp constants}
\label{Chap:About_Sharp_Constant}

\paragraph{Abstract:}We consider a possibly degenerate Kolmogorov-Ornstein-Uhlenbeck operator of the form $L={\rm Tr}(BD^2)+\langle Az,D\rangle $, where  $A$, $B $ are $N\times N $ matrices, $z\in \R^N$, $N\ge 1 $, which satisfy the Kalman condition which is equivalent to the  hypoellipticity condition. We prove the following stability result:  the    Schauder and Sobolev estimates associated with the corresponding parabolic Cauchy problem remain valid, with the same constant, for the parabolic Cauchy problem associated with a second order perturbation of $L$, namely for $L+{\rm Tr}(S(t) D^2) $  where  $S(t)$ is a {\it non-negative definite}  $N\times N $ matrix depending continuously on $t \in [0,T]$. Our approach relies on the perturbative technique based on the Poisson process introduced in \cite{Krylov:Priola17}.

\section{Introduction}
\fancyhead[RO]{Section \thesection. Introduction}

Let us first consider the following parabolic Cauchy problem:
\begin{equation}
\label{ko_K}
\begin{cases}
   \partial_t u(t,x,y)\,  = \, \Delta_{x}u(t,x,y) + x\cdot \nabla_{y} u (t,x,y)   + f(t,x,y),\\
   u(0, x,y) \, = \, 0,
\end{cases}
\end{equation}
where $(t,x,y)$ is in $(0,+\infty)\times \R^{2d}$ \textcolor{black}{for some integer $d\ge 1$.}
The underlying differential operator
\[
L^{\text{K}}\, = \, \Delta_{x} +x\cdot \nabla_{y} =  \sum_{i=1}^{d}   \partial_{x_i x_i}^2 \, + \, \sum_{i=1}^{d} x_i  \partial_{y_i}
\]
is the so-called Kolmogorov operator whose fundamental solution
was derived in the seminal paper \cite{Kolmogorov34}. This particular operator was also mentioned  by H\"ormander as the starting point for his  theory of hypoelliptic operators \cite{Hormander67}.\newline
Let us  write $z = (x,y) \in  \R^{2d}$ and by $\partial_{z_j}$ and $\partial_{z_i z_j}^2$ we denote respectively \textcolor{black}{the} first and \textcolor{black}{the} second partial derivatives with $i,j = 1, \dots, 2d$.

We are interested in studying the influence of a second order perturbation on Equation \eqref{ko_K}.  Precisely, for a time-dependent \textcolor{black}{matrix $\{S(t)\colon t\ge 0\}$ in $\R^N\otimes \R^N$ such that} $t\mapsto S(t) $ is continuous and $S(t) $ is {\it symmetric} and \emph{non-negative definite} for any fixed $t$,  we consider the \textit{perturbed} Cauchy problem:
\begin{equation}
\label{ko_K_Pert}
\begin{cases}
\partial_t u_S(t, z) \, = \, L^{\text{K}}u_S(t,z) + \sum_{i,j=1}^{2d} S_{ij}(t)\,  \partial_{z_i z_j}^2  u_S(t,z)   + f(t,z)
\\
\, \quad \;\;\; \;\;\;\;\; \quad  =: \, L^{\text{K},S}\, u_S(t,z)+f(t,z), \\
u_S(0, z) =0, \;\; z \in \R^{2d}.
\end{cases}
\end{equation}
\textcolor{black}{In particular}, we \textcolor{black}{will} show that Sobolev  (and Schauder) estimates which hold for solutions $u$ \textcolor{black}{of the Cauchy Problem} \eqref{ko_K} are also true, with the same constants,
for solutions $u_S $ to \eqref{ko_K_Pert}. Clearly, \textcolor{black}{the operator} $L^{\text{K}, S}$ can be seen as a perturbation of $L^{\text{K}}$ involving second order partial derivatives with continuous  time-dependent coefficients.


For now, let us explain our main results in a special form for Equation \eqref{ko_K} in the case of $L^p$-estimates (or Sobolev estimates). \textcolor{black}{For a statement of our results in the whole generality}, we instead refer to Section 2.\newline
 For a fixed final time $T>0$ and a source $f$ in $C_0^\infty((0,T)\times \R^{2d})$,
it is known from the work of Bramanti \textit{et al.}  \cite{Bramanti:Cerrutti:Manfredini96}, \textcolor{black}{Theorem $3.1$,} that \textcolor{black}{Equation} \eqref{ko_K} admits a unique
classical bounded solution $u$ which
satisfies for \textcolor{black}{$p$ in $(1,+\infty)$}  the \textcolor{black}{following} estimates:
\begin{equation}\label{c_K}
 \| \Delta_{x} u \|_{L^p ((0,T) \times \R^{2d})} \le C_p\;  \| f \|_{L^p ((0,T) \times \R^{2d})} = C_p   \| \partial_t u - L^{\text{K}}  u \|_{L^p ((0,T) \times \R^{2d})}.
\end{equation}
 Note that  in this case $C_p = C_p(T, d)>0$.
We will actually manage to prove that the unique classical bounded  solution $u_S$ to \eqref{ko_K_Pert} satisfies the estimate
\begin{equation}\label{c_K_PERT}
 \| \Delta_{x} u_S \|_{L^p ((0,T) \times \R^{2d})} \le C_p\;  \| f \|_{L^p ((0,T) \times \R^{2d})} = C_p   \| \partial_t u - L^{\text{K},S}  u \|_{L^p ((0,T) \times \R^{2d})},
\end{equation}
with the \textbf{same} previous constant $C_p$ as in \eqref{c_K}. This result seems to be new   even in dimension $N=2$ and even if we only consider
$S(t) = S$, $t \in [0,T]$,  where $S$  is a $2 \times 2$ {\it symmetric non-negative definite} matrix.

For a uniformly elliptic second order perturbation  $S(t) = S$, $t \in [0,T]$,  where $S$  is \emph{positive definite},  we could also have appealed to  \cite{Bramanti:Cupini:Lanconelli:Priola10}
to derive \textcolor{black}{estimates like in} \eqref{c_K_PERT}. For related estimates in the uniformly elliptic case, see also Section $4$ in Metafune \textit{et al.} \cite{Metafune:Pruss:Rhandi:Schnaubelt02}.  However, note that
 from \cite{Bramanti:Cupini:Lanconelli:Priola10}  and \cite{Metafune:Pruss:Rhandi:Schnaubelt02}  we could  only deduce that
the constant $C_p$
depends on the ellipticity constant of the perturbation  (this is the first eigenvalue $\lambda_S $ of $S$ if $0< \lambda _S \le 1$) and   on the maximum eigenvalue of $S$ (on this respect,  see also
    \cite{Krylov02} and  \cite{Priola15_parabolic}).

The remarkable point in \eqref{c_K_PERT} is that the $L^p$-estimates are stable under second order perturbations, which can be  possibly degenerate. Namely,
the fact that $S(t)$ might be degenerate for some $t$ in $(0,T)$, or even in some non-empty sub-intervals of $(0,T)$, does not affect the estimates in \eqref{c_K_PERT}.

To prove \eqref{c_K_PERT}, we
combine the results  of \cite{Bramanti:Cerrutti:Manfredini96}
with a probabilistic \textit{perturbative} approach based on the Poisson process inspired by \cite{Krylov:Priola17}. 
There, it was established in particular that the $L^p$-estimates for
 non-degenerate parabolic heat equations with \textcolor{black}{space homogeneous}
 coefficients are
 valid with constants \textcolor{black}{that are} independent of the dimension.


\vskip 1mm

\begin{remark}
Importantly, the approach of \cite{Krylov:Priola17} turns out to be sufficiently robust to handle the estimates in the degenerate directions as well. We recall that the associated maximal $L^p$-regularity
was studied e.g. in \cite{Bouchut02}, \cite{Chen:Zhang19} or \cite{Huang:Menozzi:Priola19}. Fixed $p$ in $(1,+\infty)$, there exists $\tilde C_p>0$ such that for $f$ in $C_0^\infty((0,T)\times \R^{2d})$  the  unique classical bounded solution $u$ of \eqref{ko_K} verifies
\begin{equation}\label{LP_KOLMO_DEG}
  \| (\Delta_{y})^{\frac 13}u \|_{L^p ((0,T) \times \R^{2d})} \le \tilde C_p \;  \| f \|_{L^p ((0,T) \times \R^{2d})}=\tilde C_p   \| \partial_t u - L^{\text{K}}  u \|_{L^p ((0,T) \times \R^{2d})},
 \end{equation}
where $(\Delta_{y})^{\frac 13}$ denotes the fractional Laplacian with respect to the degenerate variables $y$ in $\R^{d}$. It turns out that this estimate is also stable for the previously described second order perturbation. Namely, for $u_S $ solving \eqref{ko_K_Pert},
\begin{equation}\label{LP_KOLMO_DEG_PERT}
 \| (\Delta_{y})^{\frac 13}u_S \|_{L^p ((0,T) \times \R^{2d})} \le \tilde C_p \;  \| f \|_{L^p ((0,T) \times \R^{2d})}=\tilde C_p   \| \partial_t u - L^{\text{K},S}  u \|_{L^p ((0,T) \times \R^{2d})},
 \end{equation}
where again $\tilde C_p $ is the same as in \eqref{LP_KOLMO_DEG}. \qed
\end{remark}

\begin{remark}
The same type of stability results  will also hold for the corresponding global Schauder estimates, first established in the framework of anisotropic H\"older spaces for the solution of \eqref{ko_K} by Lunardi \cite{Lunardi97} (see also \cite{Marino20, Marino21} and the references
 therein).  We refer to estimate \eqref{SCHAU_OU_PERT}.
\end{remark}

We point out that our results in Section $3$ could  be possibly obtained by using the general  theorems  of   Section $4$ in \cite{Krylov:Priola17}. This section in \cite{Krylov:Priola17} introduces a more  general probabilistic approach and  provides quite unexpected     regularity results. However    checking in our case all the assumptions given in  that  section   is quite involved. On the other hand,      we  provide    self-contained proofs   inspired by  Sections $2$ and $3$ of \cite{Krylov:Priola17}.

It remains a challenging open problem to have a purely analytic proof of the above regularity results.

 \paragraph{Lp-estimates for degenerate Ornstein-Uhlenbeck operators.}
Let us now describe \textcolor{black}{the more general framework we are going to consider here}.
\textcolor{black}{Fixed}  $\R^N = \R^{d} \times \R^{d'}$ where $d, d' $ are two {\it non-negative}   integers such that $d + d' =N$  and $d \ge 1$.
Let us introduce the non-negative, symmetric matrix \textcolor{black}{$B$ in $\R^N\otimes\R^N$ given by}
$$ B = \begin{pmatrix}
 B_0  & 0\\
 0 & 0
\end{pmatrix} ,$$ where $B_0$ is a symmetric, positive definite matrix in $\R^{d}\otimes \R^{d}$ such that
\[\nu  \sum_{i=1}^{d}  {\xi_i^2}
  \le \sum_{i,j =1}^{d}
  (B_0)_{ij} \xi_i \xi_j
  \le \frac{1}{\nu} \sum_{i=1}^{d}  {\xi_i^2},\]
for all $\xi \in \R^{d_0}$ and some $\nu >0$.

We will use, as underlying  \textit{proxy} operators, the family of degenerate Ornstein-Uhlenbeck generators of the form
\begin{equation}
\label{DEF_OU_OP_PROXY}
L^\text{ou} f(z) = 
{\rm Tr}(B D^2 f(z))
+ \langle A z , D f(z)\rangle,\;\;\;  \; z \in \R^N,
\end{equation}
for a matrix  $A$ in $\R^N\otimes \R^N$, where $\langle \cdot , \cdot \rangle$
denote\textcolor{black}{s} the usual inner product in $\R^N$.

\textcolor{black}{Moreover}, we assume  the Kalman condition:
\begin{trivlist}
\item[\textbf{[K]}] There exists a  non-negative integer $n$,
 such that
\begin{align} \label{kal}
{\rm Rank} [ B,AB,\cdots,
  A^{n-1}B
] =N,
 \end{align}
 where $ [B,AB,...,A^{n-1}B]$ is the $\R^N\otimes \R^{Nn}$ matrix whose 
 blocks
 are $B,AB,\cdots A^{n-1}B$.

From the non-degeneracy of $B_0$, the above condition amounts to say that the vectors
\begin{align} \label{kal_EQUIV}
 \{ e_1 , \ldots, e_{ d}, A e_1 ,
 \ldots, Ae_{ d}, \ldots, A^{n-1} e_1 , \ldots, A^{n-1}
 e_{ d}\} \;\;\; \mbox{generate} \;\; \R^N,
 \end{align}
 \end{trivlist}
where $\{e_i\colon i\in \{ 1,\cdots,d\}\}$ are the first $d$ vectors of  the canonical basis for $\R^N$.

Assumption \textbf{[K]} (which \textcolor{black}{also often} appears  in control theory; see e.g.\ \cite{book:Zabczyk95}) is equivalent to the  H\"ormander condition on the
commutators (cf.\ \cite{Hormander67})  ensuring the hypoellipticity of the operator $\partial_t-L^\text{ou}$. In particular, it implies the existence and \textcolor{black}{the} smoothness of a distributional solution for the following equation:
\begin{equation}
\label{eq:OU_initial:intro}
\begin{cases}
    \partial_t u(t,z) \, = \,  L^\text{ou} u(t,z) +f(t,z), &\mbox{ on }(0,T)\times \R^N;\\
   u(0,z)\, = \, 0, &\mbox{ on } \R^N,
\end{cases}
\end{equation}
 where $f$ is a function in $C_0^\infty((0,T)\times \R^N)$.

Similarly to  \cite{Krylov:Priola17}, we will prove below the existence and the uniqueness of bounded regular solutions to \eqref{eq:OU_initial:intro}
 assuming that the source $f$ belongs to $B_b\left(0,T;C^\infty_0(\R^N)\right)$, which contains $C_0^\infty((0,T)\times \R^N)$, \textcolor{black}{and that} \textcolor{black}{ can be roughly described as the family of functions which are bounded measurable in time and compactly supported in space, \textcolor{black}{uniformly} in time} (see Section \ref{SEC_NOT} for a precise definition). Equation  \eqref{eq:OU_initial:intro}
 will be \textcolor{black}{understood in an} integral form (cf.\ Equation \eqref{int2}).

 By Theorem $3$ in \cite{Bramanti:Cupini:Lanconelli:Priola10}  \textcolor{black}{and exploiting some explicit properties of the underlying heat kernel (see Section \ref{KH_OU} below), it can be derived that}
 for \textcolor{black}{any fixed $p$ in $(1,+\infty)$,} there exists $C_p = C_p(B_0,d,d', T)$  such that
 \begin{equation}\label{de2}
\| D^2_x u \|_{L^p ((0,T) \times \R^N)}\, \le \,  C_p  \| \partial_t u -  L^\text{ou} u    \|_{L^p ((0,T) \times \R^N)}\, = \, C_p  \| f    \|_{L^p ((0,T) \times \R^N)},
\end{equation}
\textcolor{black}{where} for any $(t,z)$ in $[0,T]\times\R^N$, $D_x^2u(t,z)$ \textcolor{black}{stands for the Hessian matrix in $\R^{d}\otimes\R^{d}$ with respect to the variable $x$}. We set
\[B_I \, :=\,  \begin{pmatrix}
 I_{d,d}  & 0_{d,d'}\\
 0_{d',d} & 0_{d',d'}
\end{pmatrix}\]
  and note, in particular, that \eqref{de2} can be rewritten in the following, equivalent way:
\begin{equation}\label{de3}
\begin{split}
\| B_I D^2 u \, B_I \|_{L^p ((0,T) \times \R^N)} \, &=
\, \|  D^2_x u \,  \|_{L^p ((0,T) \times \R^N)} \\
&\le \,  C_p \| \partial_t u -   L^\text{ou} u    \|_{L^p ((0,T) \times \R^N} \\
&= \, C_p  \| f    \|_{L^p ((0,T) \times \R^N)},
\end{split}
\end{equation}
where $D^2 u = D^2_z u$ represents instead the full Hessian matrix in $\R^N\otimes\R^N$ with respect to $z$.

\textcolor{black}{Fixed} a continuous mapping $t\mapsto S(t)$ \textcolor{black}{such that} $S(t)$ is a {\it symmetric} and \emph{non-negative definite} matrix in $\R^N\otimes\R^N$ for any $t\in [0,T]$, we consider again the following perturbation of $L^\text{ou}:$
\begin{equation} \label{ciao}
\begin{split}
L_t^{\text{ou},S} \phi(z) \, :&=\,  
{\rm Tr}(B D^2 \phi(z)) +
 {\rm Tr}(S(t) D^2 \phi(z))
+ \langle A z , D \phi(z)\rangle\\
&= \, L^\text{ou} \phi(z)+
{\rm Tr}(S(t) D^2 \phi(z)),
\end{split}
\end{equation}
where $z$ is in $\R^N$.
For the solution $u_s$ of the related Cauchy problem
\begin{equation}
\label{eq:OU_PERT0}
\begin{cases}
    \partial_t u_S(t,z) \, = \,  L_t^{\text{ou},S} u_S(t,z) +f(t,z), &\mbox{ on }(0,T)\times \R^N;\\
   u_S(0,z)\, = \, 0, &\mbox{ on } \R^N,
\end{cases}
\end{equation}
we will prove the following main theorem:

\begin{theorem} \label{d44}  Let us consider \eqref{eq:OU_PERT0} with $f \in B_b\left(0,T;C^\infty_0(\R^N)\right) $. \textcolor{black}{Then,} there exists a unique solution $u_S$ \textcolor{black}{of Cauchy Problem \eqref{eq:OU_PERT0}}
 which verifies,
with the same constant $C_p$ as in \eqref{de3},
\begin{equation}\label{de4}
\begin{split}
\|  D^2_x u_S \|_{L^p ((0,T) \times \R^N)}\, 
&= \, \| B_I D^2 u_S B_I \|_{L^p ((0,T) \times \R^N)} \\ 
&\le\,  C_p  \| \partial_t u_S -   {L}_t^{{\rm{ou}},S} u_S  
\|_{L^p ((0,T) \times \R^N)}\, = \, C_p   \| f
\|_{L^p ((0,T) \times \R^N)}.
\end{split}
\end{equation}
\end{theorem}

We point out  that for time-homogeneous non-negative definite matrices $S$,
 the  corresponding  elliptic $L^p$-estimates  as in formula (5) of \cite{Bramanti:Cupini:Lanconelli:Priola10}   (replacing $\mathcal A$ in \cite{Bramanti:Cupini:Lanconelli:Priola10} with  $ L^{\text{ou},S}: =$ $
{\rm Tr}(B D^2 \cdot ) +
 {\rm Tr}(S D^2 \cdot  )
+ \langle A z , D \cdot \rangle
)$ {\sl with a constant independent of $S$,}  could also be derived from \eqref{de4}  using an argument  given in \cite{Bramanti:Cupini:Lanconelli:Priola10}.

 For more information on the OU operator $L^{ou}$ we also  refer to  the recent work by Fornaro \textit{et al.} \cite{Fornaro:Metafune:Pallara:Schnaubelt21} about  full description of the spectrum of degenerate OU operators in $L^p$-spaces.\newline
Independently from the constant preservation, we also emphasize that the $L^p$ estimates in \eqref{de4} for the perturbed operator seem, to the best of our knowledge, to be new and have some interest by their own.\newline
Let us eventually mention that such stability results could turn out to be useful to investigate the well posedness of some related stochastic differential equations through the corresponding martingale problem.

We could  actually derive more general estimates, \textcolor{black}{possibly} involving the degenerate directions as well, \textcolor{black}{dependingly on} the structure of $A$. \textcolor{black}{Some results in that direction are gathered in Section \ref{ESTENSIONI}}.
Anyhow, to illustrate our approach we now briefly present the various steps to derive  \eqref{de4}.

\subsection{Strategy of the proof for Estimates (\ref{de4})}  \label{SEZ_STRAT}
Fixed a classical bounded solution $u$ to Cauchy Problem \eqref{eq:OU_initial:intro}, let us introduce $v(t,z) := u(t, e^{-tA} z)$.  This well-known transformation (cf.\ \cite{Daprato:Lunardi95})  precisely allows to get rid of the drift term in the PDE satisfied by  $v$.
Indeed, we have \textcolor{black}{that}  $u(t,z) = v(t, e^{tA} z)$
and since $ u$ solves \eqref{eq:OU_initial:intro}, it holds for any $(t,z)$ in $(0,T)\times \R^N$, that:
 \begin{equation}
 \begin{split}
 f {(t,z)}\, &= \,  \partial_t u(t,z) {\color{black} - L^{\text{ou}}}  u(t,z)\\
 &= \, v_t(t,e^{tA} z) + \langle D v
 (t, e^{tA} z) , A e^{tA} z\rangle  {\color{black} -} {\text Tr} \big( e^{tA} B
 e^{tA^*} D^2  v(t, e^{tA} z) \big) \\
& \qquad - \langle D v
 (t, e^{tA} z) , A e^{tA} z\rangle\\
&= \, v_t(t,e^{tA} z)    {\color{black} -} {\text Tr} \big( e^{tA} B
 e^{tA^*} D^2 v(t, e^{tA} z) \big).
 \end{split}
 \label{COMP_0H_VAR}
\end{equation}
\textcolor{black}{Denoting $\tilde{f}(t,z):= f(t,e^{-tA}z)$,}
It now follows that $v$ satisfies the PDE:
 \begin{equation} \label{ma}
 \begin{cases}
 \partial_tv(t, z)   \, = \,  {\text Tr} \big( e^{tA} B
 e^{tA^*} D^2 v(t, z) \big) + \tilde f(t,z) &\mbox{ on }(0,T)\times \R^N;\\
 v(0,z) =0 &\mbox{ on } \R^N.
 \end{cases}
\end{equation}
In terms of the function $v$, the known estimates in \eqref{de3} rewrites as:
\begin{equation} \label{1}
\| B_I e^{t A^*} D^2 v (t, e^{tA} \cdot ) \, e^{tA} B_I \|_{L^p((0,T)\times \R^N)} \le C_p  \|  \tilde f(t, e^{tA}  \cdot )  \|_{L^p((0,T)\times \R^N)},
\end{equation}
where we used the notation  $\| B_I e^{t A^*} D^2 v (t, e^{tA} \cdot ) \, e^{tA} B_I \|_{L^p((0,T)\times \R^N)} $  to stress the dependence on $t$
 instead of the more precise  formulation
\[{ \| B_I e^{ \cdot \,  A^*} D^2 v (\cdot , e^{ \cdot \, A} \cdot ) \, e^{ \cdot \, A} B_I \|_{L^p((0,T)\times \R^N)}}.\]

By changing variable in the integrals, Control \eqref{1}  is equivalent to
\begin{equation}\label{2}
\| B_I e^{tA^*} D^2 v (t,  \cdot ) \, e^{tA} B_I \|_{L^p((0,T)\times \R^N,{m})} \le C_p  \|  \tilde f  \|_{L^p((0,T)\times \R^N, {m})},
\end{equation}
where $L^p((0,T)\times \R^N,m)$ denotes the $L^p$ norms w.r.t. the measure \[m(dt,dx) \, :=\, {\rm det}(e^{-At}) dt dx.\]
Considering now the following, more general equation
\begin {equation} \label{d2}
\begin{cases}
 \partial_t w(t, z)    + {\text Tr} \big( e^{tA} B
 e^{tA^*} D^2 w(t, z) \big) + {\text Tr} \big( e^{tA} S(t)
 e^{tA^*} D^2 w(t, z) \big) \, = \, \tilde f(t,z);\\
 w(0,z)=0,
 \end{cases}
\end{equation}
we  can  establish \textcolor{black}{the well-posedness of} the Cauchy problem \eqref{d2}, \textcolor{black}{exploting}, for instance, probabilistic arguments \textcolor{black}{and} the underlying Gaussian process.

Now, the crucial  step consists in adapting  some arguments from \cite{Krylov:Priola17} based on the use of the Poisson process to derive that the same
 $L^p$-estimates in \eqref{2} still hold for $w$, independently \textcolor{black}{from} the non-negative definite, symmetric matrix $S(t)$. Precisely,
\begin{equation} \label{11}
\| B_I e^{tA^*} D^2 w (t, \cdot ) \, e^{tA} B_I \|_{L^p((0,T)\times \R^N, {m})} \le C_p \|  \tilde f(t,  \cdot )  \|_{L^p((0,T)\times \R^N, {m})},
\end{equation}
with the \emph{same} constant $C_p$ \textcolor{black}{appearing in \eqref{2}}.\newline
The last step then consists in  coming back  to the Ornstein-Uhlenbeck operators \textcolor{black}{framework}.  Namely, \textcolor{black}{we introduce} $\tilde u(t,z) := w(t, e^{tA} z)$ \textcolor{black}{which} solves, by definition, the following equation:
 \begin{equation*}
\begin{cases}
 \partial_t \tilde u(t,z) + L_t^{\text{ou}, S} \tilde u(t,z)\, =  \,f  (t,z),&\mbox{ on } (0,T)\times \R^N,\\
 \tilde u(0,z)\, = \, 0,&\mbox{ on } \R^N.
\end{cases}
\end{equation*}
Hence, it holds that $\tilde{u}=u_S$.
Noticing that 
\[D^2 w(t, \cdot)\, = \, D^2 [\tilde u(t, e^{-tA} \, \cdot )] \, =\, e^{-tA^*}D^2 \tilde u(t, e^{-tA} \cdot )e^{-tA},\]
we thus get from \eqref{11} that the following estimates hold:
\begin{equation}\label{de311}
\| B_I D^2 \tilde u \, B_I \|_{L^p ((0,T) \times \R^N)} \, \le \, 
C_p  \|  f \|_{L^p ((0,T) \times \R^N)} .
\end{equation}
Through the previous steps we have then constructed a solution $\tilde u $ \textcolor{black}{of Cauchy Problem} \eqref{eq:OU_PERT0} which indeed satisfies \textcolor{black}{the estimates in} \eqref{de4} with the same $C_p$, associated with the \textit{unperturbed} or \textit{proxy} operator.
 The maximum principle will eventually provide uniqueness for the solution $\tilde u $.

  \vskip 1mm
  \begin{remark} i) We point out that we could also consider more  general time-dependent Ornstein-Uhlenbeck operators like:
  \[M\phi(z) \, = \, {\rm Tr}(B(t) D^2 \phi(z))
+ \langle A z , D \phi(z)\rangle.\]
Arguing as before starting from $L^p$-estimates (or Schauder estimates) for $M$ we can derive the same $L^p$-estimates (or Schauder estimates) for a perturbation of $M$ like \eqref{ciao}.

ii) We could  extend the $L^p$-estimates (or the Schauder estimates) related to  $L^{{\rm ou}}$ to more general
 operators like
$$
L_t^{\text{ou},S} \phi(z) +  \langle b(t), D \phi(z) \rangle
$$
where $b: \R_+ \to \R^N$ is continuous. We can even add to $L_t$ a possibly degenerate   non-local perturbation  (cf. Section $7$ of \cite{Krylov:Priola17}). The $L^p$-estimates (or Schauder estimates) are still preserved with the same constant. For the sake of simplicity in the sequel we will only consider $b(t)=0$ and we will  not deal with  non-local perturbations of $L_t^{\text{ou},S}.$
  \end{remark}

\paragraph{Organisation of the Paper.}
The article is organised as follows. At the end of the current section, we first give some useful notations.
In Section \ref{SEC_DRIFT_LESS} we will then focus on driftless second order Cauchy problems associated with a non-negative, possibly degenerate, diffusion matrix. We will also consider its relation with the Ornstein-Uhlenbeck dynamics. We will establish through the probabilistic perturbation approach of \cite{Krylov:Priola17} that if some $L^p$-estimates hold for a particular diffusion matrix so does it, with the same associated constant, for a non-negative perturbation of the diffusion matrix (see Section $3$).
Finally by the argument of Section 1.1 we will obtain \eqref{de311}.
Stability results in anisotropic Sobolev space and Schauder estimates are given in Section $4$.

\subsection{Definition  of solution and useful  notations}\label{SEC_NOT}

Let us consider the following Cauchy problem:
\begin{equation}
\label{eq:Cauchy_problem_gen0}
\begin{cases}
\partial_t v(t,z) \, = \, \text{tr}\left(Q(t)D^2 v(t,z)\right) +
\langle  b(t,z), Dv(t,z) \rangle
+f(t,z), &\mbox{ on }(0,T)\times\R^N;\\
v(0,z) \, = \, 0, &\mbox{ on } \R^N;
\end{cases}
\end{equation}
where $Q\colon [0,T]\to \R^N\otimes \R^N$ is a {\it   continuous, symmetric, non-negative definite} matrix and $b: [0,T] \times \R^N \to \R^N$ is a {\it continuous} function such that  $|b(t,z)| \le K_T (1+ |z|)$, $(t,z) \in [0,T] \times \R^N$, for some constant $K_T>0$.

The function $f$ belongs to  $B_b\left(0,T;C^\infty_0(\R^N)\right)$, the space of all Borel bounded functions $\phi \colon [0,T]\times\R^N \to \R$ such
that $\phi(t,\cdot)$ is smooth and compactly supported for any $t$ in $[0,T]$, for any $k$ in $\N$ the $C^k(\R^N)$-norms of $\phi(t,\cdot)$ are bounded in time and the supports of the functions $\phi(t,\cdot)$ are contained in the same ball. Moreover, we require that, for any $z \in \R^N$, the mapping:
\begin{equation}\label{dd}
t \mapsto \phi(t,z)
\end{equation}
is a \emph{piece-wise continuous} function on $[0,T]$, i.e.\ it is continuous except \textcolor{black}{for} a finite number of points.

\begin{remark}
 Note that to perform the technique used in \cite{Krylov:Priola17} and based on the Poisson process we
need to consider equations like \eqref{eq:Cauchy_problem_gen0} with a source $f$ which is possibly discontinuous in time (cf.\ the proof in Section $2$ of \cite{Krylov:Priola17} and Section \ref{SEC_PERTURBATIVO} below).
\end{remark}

We interpret \textcolor{black}{Cauchy Problem \eqref{eq:Cauchy_problem_gen0}} in  an {\it integral} form:
\begin{equation}\label{int2}
v(t,z)=\int_0^t \Big(f(s,z)+{\rm Tr}(Q(s) D^2 v(s,z) ) +  \langle  b(s,z), Dv(s,z) \rangle\Big) ds.
\end{equation}
\textcolor{black}{In particular,} we say that a continuous and bounded function $v : [0, T] \times \R^N \to \R$ is a solution to Equation \eqref{eq:Cauchy_problem_gen0} if $v(t, \cdot)$ belongs to $C^{2}(\R^N)$, for any $t \in [0,T]$, and \eqref{int2} holds as well, for any $(t,z)$ in $[0,T]\times \R^N$.\newline
We finally note that, for any $z \in \R^N$, the function $t \mapsto v(t,z)$ is a $C^1$-piece-wise   function on $[0,T]$.\newline
By Theorem 4.1 in \cite{Krylov:Priola10} we deduce in a quite standard way that if a solution $v$ exists, then it is unique and the following maximum principle holds:
\begin{equation}\label{max}
 \sup_{(t,z)\in [0,T] \times \R^N}|v(t,z)|\le T\,  \sup_{(t,z)\in [0,T]\times \R^N}|f(t,z)|.
\end{equation}
About the proof of \eqref{max}, we only make some remarks.
 By considering $v$ and $-v$, we see that it is enough to prove
 that $v(t,z) \le T \| f\|_{\infty}$, for all $(t,z)\in [0,T] \times \R^N$. Moreover, setting $\tilde v = $ $
  v - t \| f\|_{\infty}$, we note that  $\tilde v$ verifies  \eqref{int2} with $f$ replaced by $f - \| f\|_{\infty} \le 0$.
 Finally, by considering the equation verified by  $e^{- t} \tilde v$, we can apply Theorem 4.1 in \cite{Krylov:Priola10}
   to obtain the result.

\section{Estimates for driftess second order operators and related perturbation}
\fancyhead[RO]{Section \thesection. Estimates for driftess second order operators and related perturbation}
\label{SEC_DRIFT_LESS}
Throughout this section, we consider the following Cauchy problem:
\begin{equation}
\label{eq:Cauchy_problem_gen}
\begin{cases}
\partial_t v(t,z) \, = \, \text{Tr}\left(Q(t)D^2 v(t,z)\right) +f(t,z), &\mbox{ on }(0,T)\times\R^N;\\
v(0,z) \, = \, 0, &\mbox{ on } \R^N,
\end{cases}
\end{equation}
 which can be seen as a special case of \eqref{eq:Cauchy_problem_gen0} when $b=0$. Moreover, we assume that $Q$ is not identically zero.

\subsection {Well-posedness }
\begin{prop}[Well-posedness in integral form for the driftless Cauchy problem]\label{PROP_WP_DRIFTLESS}
\textcolor{black}{Let $f$ be in $B_b\left(0,T;C^\infty_0(\R^N)\right)$. Then,} there exists a unique
 solution $v$ to Cauchy problem \eqref{eq:Cauchy_problem_gen} \textcolor{black}{in an integral sense}, i.e.\ it solves for  $(t,z)\in [0,T]\times \R^N$:
\begin{equation}
\label{INTEGRAL}
v(t,z)=\int_0^t \Big(f(s,z)+{\rm Tr}(Q(s) D^2 v(s,z)\Big)ds.
\end{equation}
We will denote in short $v=PDE(Q,f)$.

\end{prop}
\begin{proof}
By the maximum principle  (cf.\ Equation \eqref{max}) uniqueness  holds for Cauchy Problem \eqref{eq:Cauchy_problem_gen}. We can then focus on proving the existence of a solution.\newline
Let us introduce now \[v(t,z) \, := \, \int_0^t \mathbb E[f(s,z+I_{s,t})] \, ds\]
with the following notation: $I_{s,u}:=\sqrt 2\int_s^u Q(r)^{1/2}dW_r $, where $W$ is an $N$-dimensional Brownian motion on some probability space $(\Omega,\mathcal F,(\mathcal F_t)_{t\ge 0},\mathbb P) $ and $Q(r)^{1/2}$ stands for a square root of $Q(r)$, i.e. $Q(r)=Q(r)^{1/2}(Q(r)^{1/2})^* $.\newline
Applying the  It\^o formula in space to $f(s, z+I_{s,u} )_{u\in [s,t]}$, we get that
$$
\mathbb{E} f(s, z+  I_{s,t} )\, = \,f(s,z)+ \mathbb E\left[\int_s^{t} {\rm{Tr}}(Q(u) D^2 f(s,z+I_{s,u}))\,du\right].
$$
Hence,
$$
v(t,z)\, = \,\int_{0}^t \Big(f(s,z)+\mathbb E\left[\int_s^{t} {\rm{Tr}}(Q(u) D^2 f(s,z+I_{s,u}))\,du\right]\Big)ds,$$
from which it readily follows that
\[
\begin{split}
\partial_t v(t,z) \, &= \, f(t,z)+\int_0^t  \mathbb E\left[{\rm{Tr}}(Q(t) D^2f(s,z+I_{s,t}))\right]\, ds\\
&= \, f(t,z)+{\rm{Tr}}\Big(Q(t) D^2\int_0^t  \mathbb E\left[ f(s,z+I_{s,t})\right]\, ds\Big)\\
&= \, f(t,z)+{\rm{Tr}}\Big(Q(t) D^2v(t,z)\Big),
\end{split}
\]
\textcolor{black}{for almost every $t\in [0,T]$ and any $z \in\R^N $}.
\end{proof}

\subsection{Relation to the Ornstein-Uhlenbeck dynamics}\label{REL_OU}
If now in particular, $Q(t)$ has the particular form   $Q(t)= e^{tA} B
 e^{tA^*}$ (cf.\ Equation \eqref{ma}), we introduce
 $$
 u(t,z):=v(t,e^{tA}z),
 $$
  where $v$ is the solution to \eqref{INTEGRAL}   (see Proposition \ref{PROP_WP_DRIFTLESS}).
 Since we can differentiate with respect to $t$ the function $u(\cdot,z)$ for a.e. $t \in [0,T]$, we can perform
 computations similar to \eqref{COMP_0H_VAR} and get  that $ u (t,z)$  solves in integral form:
\begin{equation}
\label{eq:Cauchy_pro}
\begin{cases}
\partial_t u(t,z) \, = \,  L^\text{ou} \, u(t,z) +\bar f(t,z), &\mbox{ on }(0,T)\times\R^N;\\
u(0,z) \, = \, 0, &\mbox{ on } \R^N;
\end{cases}
\end{equation}
with $L^\text{ou}$ as in  \eqref{DEF_OU_OP_PROXY},
$\bar f(t,z)=f(t,e^{tA}z) $.  Precisely, for all $(t,z)\in [0,T]\times\R^N $,
\begin{equation}
\label{INTEGRAL_OU}
u(t,z)=\int_0^t \Big(\bar f(s,z)+ L^\text{ou} v(s,z)\Big) ds.
\end{equation}
Hence, $u$ is a solution to \eqref{eq:Cauchy_pro}.\newline
Let us also point out that the well-posedness of \eqref{eq:Cauchy_pro} could also have been obtained directly from Gaussian type calculations,
similar to those in the proof of Proposition \ref{PROP_WP_DRIFTLESS},
 introducing $u^{{\rm ou}}(t,z):=\int_0^t \mathbb E[\bar f(s,e^{(t-s)A}z+I_{s,t}^{{\rm ou}})] ds$ where $I_{s,u}^{{\rm ou}}:=\sqrt 2\int_s^u \textcolor{black}{e^{(u-v)A}} B dW_v  $.

\subsection{About the Lp-estimates in (\ref{de2}) }
\label{KH_OU}
{The aim of this section is to} fully justify the estimates in \eqref{de2}.  This is a consequence of the previous probabilistic representation and of Theorem $3$ in \cite{Bramanti:Cupini:Lanconelli:Priola10}. For $u$ solving \eqref{eq:OU_initial:intro}, it holds
that for all $(t,z)\in [0,T]\times \R^N $,
\begin{equation}\label{PROB_REP_OU}
u(t,z)\,= \, \int_0^t \mathbb E\left[ f(s,e^{A(t-s)}z+I_{s,t}^{{\rm ou}})\right] ds\, = \, \int_{0}^t\int_{\R^N } f(s,z') p^{{\rm ou}}(t-s,z,z') \, dz' ds,
\end{equation}
where for $v>0$, $p^{{\rm ou}}(v,z,\cdot) $ stands for the density at time $v$ of the stochastic process
\[
X^{{\rm ou}}_u \,:= \,e^{Au}z+\sqrt 2\int_0^u \textcolor{black}{e^{A(u-w)}} B dW_w\, = \, z+\int_0^u A X^{{\rm ou}}_w dw+ BW_u,\, u\ge 0.
\]
We recall from \cite{Lanconelli:Polidoro94} that assumption \textbf{[K]} is equivalent to the fact that there exist $n\in \N $ and positive integers $\{d_i\colon i\in 1,\cdots, n\}$ such that $d=d_1$, $\sum_{i=1}^n d_i=N$ and for all $i\in \{2,\cdots,n\}$ the matrixes
$${\mathscr A}^i\, := \, (A_{j,\ell})_{(j,\ell)\in \{\sum_{m=1}^{i-1}d_m+1,\cdots,\sum_{m=1}^{i} d_i\}\times \{\sum_{m=1}^{i-2}d_m +1,\cdots,\sum_{m=1}^{i-1} d_m\}},$$
with the natural notation $\sum_{m=1}^{0}=0$,
have rank $d_i$.  {The matrix $A$ writes:
\begin{equation}
\label{sotto}
A \, = \,
    \begin{pmatrix}
\ast   & \ast  & \dots  & \dots  & \ast   \\
{\mathscr A}^2  & \ast  & \ddots & \ddots  & \vdots   \\
0_{d_3,d}      & {\mathscr A}^3   & \ast  & \ddots & \vdots \\
\vdots &\ddots & \ddots& \ddots & \ast \\
0_{d_{n},d}      & \dots & 0_{d_{n}, d_{n-2}}     & {\mathscr A}^{n}    & \ast
    \end{pmatrix}.
\end{equation}
}
Following the proof of Lemma 5.5 in \cite{Delarue:Menozzi10}, where the case $ d_i=d$ for any $i$ in $\llbracket 1, n\rrbracket$ is addressed, it can be derived that there exists $C\ge 1$ s.t. for all $(v,z,z')\in (0,T]\times (\R^N)^2 $, 
\begin{equation}
\label{EST_DENS}
|D_{x}^2p^{{\rm ou}}(v,z,z')|\, \le \, \frac{C}{v^{\sum_{i=1}^{n} d_i(i-\frac{1}{2}) +1}}\exp\left(-C^{-1}v|\overline{\mathbb {M}}_v^{-1}(e^{Av}z-z')|^2 \right),
\end{equation}
where
$$\overline{\mathbb {M}}_v \, := \, \text{diag}(vI_{d\times d},v^2 I_{d_2\times d_2},\dots,v^{n}I_{d_{n}\times d_{n}}), \quad v\ge0,$$
reflects the various scales of the system. For a given function $f\in B_b\left(0,T;C^\infty_0(\R^N)\right)$, it is then clear from \eqref{PROB_REP_OU} and \eqref{EST_DENS} that for all $(t,z)\in (0,T]\times \R^N $:
\begin{equation}\label{REP_SING}
D_{x}^2 u(t,z)\, = \, {\rm p.v.}\int_{0}^t\int_{\R^N } f(s,z')D_x^2 p^{{\rm ou}}(t-s,z,z') dz' ds.
\end{equation}
It indeed suffices to observe that:
\[
\begin{split}
\Bigl|{\rm p.v.}&\int_{0}^t\int_{\R^N } f(s,z')D_x^2 p^{{\rm ou}}(t-s,z,z')\, dz' ds\Bigr|\\
&= \, \left|{\rm p.v.}\int_{0}^t\int_{\R^N } [f(s,z')-f(s,e^{A(t-s)}x)]D_x^2 p^{{\rm ou}}(t-s,z,z')\, dz' ds\right|\\
&\underset{\eqref{EST_DENS}}{\le} \, \sup_{s\in [0,T]}\|Df(s,\cdot)\|_\infty    \\ &\qquad\times \int_{0}^t\int_{\R^N }\frac{C}{(t-s)^{\sum_{i=1}^{n} d_i(i-\frac{1}{2}) +\frac{1}{2}}}\exp\left(-C^{-1}(t-s)|\mathbb T_{t-s}^{-1}(e^{A(t-s)}z-z')|^2 \right)\, dz'ds\\
&\le\,  C \sup_{s\in [0,T]}\|Df(s,\cdot)\|_\infty T^{\frac 12}.
\end{split}
\]
The estimates in \eqref{de2}
 now follows from the proof of Theorem 3 in \cite{Bramanti:Cupini:Lanconelli:Priola10}, starting from \eqref{REP_SING} instead of (16) therein. The strategy is clear. {It is necessary to}  introduce a cut-off function which separates the points $(s,z') $ which do not induce any singularity in \eqref{REP_SING} for the derivatives of the density, namely such that $t-s\ge c_0 $ or $|e^{A(t-s)}z-z'|\ge c_0 $, for some fixed constant $c_0>0$, from those who are close to the singularity.  For the non-singular part of the integral the expected $L^p$-control readily follows from \eqref{EST_DENS} and the Young inequality (see also Proposition 5 in \cite{Bramanti:Cupini:Lanconelli:Priola10}), whereas the derivation of the bound for the singular part requires some involved harmonic analysis, see Section 4  {on} the same reference.
  We can also refer to Theorem 11 and its proof in \cite{Priola15} for similar issues {linked with} the corresponding $L^p$-estimates for degenerate Ornstein-Uhlenbeck {operators} in an elliptic setting.

\subsection{The main result  for Equation (\ref{eq:Cauchy_problem_gen})}

{Let us fix $p$ in $(1,+\infty)$ and  assume that} there exist $R(t) \in \R^N \otimes \R^N$ depending continuously on $t \ge 0$ and a constant $C_p>0$, such that for any  $f$ in $B_b\left(0,T;C^\infty_0(\R^N)\right)$,
the unique solution $v = PDE(Q,f)$ to {equation} \eqref{eq:Cauchy_problem_gen} satisfies
\begin{equation}\label{s22}
 \| R(t)^* D^2 v  \, R(t) \,  \|_{L^p((0,T)\times \R^N, {\mathfrak m})} \le C_p  \|   f  \|_{L^p((0,T)\times \R^N,{\mathfrak m})},
\end{equation}
{for some absolutely continuous measure $\mathfrak m $ w.r.t. the Lebesgue measure on $ [0,T]\times \R^N$  such that  $\mathfrak m(dt,dx)=g(t)dtdx $ for some borel bounded function $g$  
(note that in \eqref{2} we have $R(t) =e^{tA} B^{1/2} $, $\mathfrak m(dt,dx)=g(t)dtdx={\rm det}(e^{-At}) dtdx $)}. 
\newline
We would like to {exhibit that a control} like \eqref{s22} {also} holds for the solution $w$ to {the following Cauchy Problem:}
\begin{equation}\label{w}
\begin{cases}
\partial_t w(t,z) \, = \, \text{tr}\left(Q(t)D^2 w(t,z)\right)+\text{tr}\left(Q'(t) D^2 w(t,z)\right)+f(t,z), &\mbox{ on }(0,T)\times \R^N;\\
w(0,z) \, = \, 0, &\mbox{ on }\R^N,
\end{cases}
\end{equation}
Namely we have to prove \textcolor{black}{the following result}.

\begin{theorem}\label{uno} Let us consider equations  \eqref{eq:Cauchy_problem_gen} and   \eqref{w}  where   $Q(t)$, $Q'(t)$ are two  continuous in time, non-negative definite matrices in $\R^N\otimes \R^N$ and   $f \in B_b\left(0,T;C^\infty_0(\R^N)\right)$. Assume that   estimate \eqref{s22} holds as explained above.

 Then the solution $w$ to \eqref{w} verifies
\begin{equation}\label{s224}
 \| R(t)^* D^2 w  \, R(t) \,  \|_{L^p((0,T)\times \R^N, {\mathfrak m})} \le C_p  \|   f  \|_{L^p((0,T)\times \R^N, {\mathfrak m})},
\end{equation}
$p \in (1, \infty)$ with the same constant $C_p$ as in \eqref{s22}.
\end{theorem}

From Theorem \ref{uno} using the argument of Section 1.1 we can easily derive Theorem \ref{d44}.

\section{Perturbation  argument for proving Theorem \ref{uno}}
\fancyhead[RO]{Section \thesection. A perturbation  argument for  proving Theorem \ref{uno}}
We aim here at applying the probabilistic perturbative approach considered in \cite{Krylov:Priola17}. The key  idea in that work was, for a well-posed PDE which enjoys some quantitative given estimates, to introduce a \textit{small} random perturbation in the source $f$  through a suitable  Poisson type process  and to investigate the properties of the associated PDE involving an unknown function $v$.
 After considering  a  {small} random perturbation of $v$,
we arrive at the useful integral formula \eqref{eq:1_OU}.  Taking the expectation, the contributions  associated with the jumps yield, for an appropriate intensity of the underlying Poisson process, a finite difference operator.
For the PDE satisfied by the expectation, involving the finite difference operator, {\it the  initial estimates are preserved}.
 Repeating the previous argument we can obtain a PDE  involving the composition of two  finite difference operators.\newline
 Compactness arguments  then allow to derive that, the initial estimates still hold at the limit with the  composition of two finite difference operators  replaced by the corresponding differential operator of order two. Iterating  this procedure we can obtain the result.

Below, we start recalling  basic properties of  Poisson type processes and   corresponding stochastic integrals, which  are needed for our approach.

\subsection{Poisson stochastic integrals}

We {briefly} recall here the very definition of {the stochastic}
 integral     {driven by a Poisson process}. 
We start reminding the construction of {such} processes. 
\newline
Given a probability space $(\Omega,\mathcal{F},\mathbb{P})$ to be fixed from this point further,  we start considering a sequence of independent real-valued random variables $\{\tau_m\}_{m \in \N}$ on $\Omega$ whose distribution is exponential of parameter $\lambda>0$:
\[\mathbb{P}(\tau_m>r) \,= \, e^{-r\lambda},  \quad r\ge 0.\]
We can then define the partial sums sequence $\{\sigma_n\}_{n \in \N}$  as follows:
\begin{align*}
    \sigma_0\, = \, 0; \;\;\;
    \sigma_m  =   \sum_{i=1}^m \tau_i, \quad n=1,2, \dots
\end{align*}
For any fixed $t\ge 0$, $\pi_t$ now denotes the number of consecutive sums of $\tau_i$ which lie on $[0,t]$, i.e.
\begin{equation}
\label{poisson}
  \pi_t=\sum_{n=0}^{\infty} \mathds{1}_{\sigma_m \le t},
\end{equation}
where $\mathds{1}_{\sigma_m \le t}$ represents the indicator function of the event $\{\sigma_m\le t\}$. The process $\{\pi_t\}_{t\ge 0}$ we have just constructed is usually known in the literature as a \emph{Poisson process} with intensity $\lambda$ (see, for instance, \cite{book:Protter04}).

Now, let $c : [0,T] \to \R^N$ be a continuous function.  We can define the  Poisson stochastic integral  as
\begin{gather} \label{s25}
b_t \, :=  \,\int_0^t c (s) d\pi_s \,= \, \sum_{\sigma_k \le t,\;\; k \ge 1}  {c(\sigma_k)}= \sum_{0 < s \le t}  {c(s)} (\pi_s - \pi_{s-})
\end{gather}
 $b_0=0$ (as usual $\pi_{s-}(\omega)$ denotes the left limit at $s$, for any $\omega$, $\mathbb{P}$-a.s.).  We now recall the following formula for the expectation of the stochastic integral:
\begin{equation}\label{see}
\mathbb{E} [\int_0^t c (s) d\pi_s]  = \lambda \int_0^t c (s) ds.
\end{equation}
 (cf.  Lemma 2.1  in \cite{Krylov:Priola17} for a direct proof; see also   Theorem 16 in \cite{book:Protter04} and Theorem 5.3 in \cite{Krylov:Priola17} for a more general formula involving stochastic integrals with predictable processes against the Poisson process).
We also recall  the following more general  result.

\vskip 1mm

\begin{lemma} \label{LEMME_POISSON}
Let $\{\pi_t\}_{t\ge0}$ be a Poisson Process of intensity $\lambda$ on a probability space $(\Omega,\mathcal{F},\mathbb{P}) $.
Let us consider a stochastic process $ (\xi_t )_{t \in [0,T]} $ with values in $\R$ which has cadlag  paths ($\mathbb{P}$-a.s.) and is ${\cal F}_t$-adapted where ${\cal F}_t$ is the augmented $\sigma$-algebra generated by the random variables $\pi_s$, $0 \le s \le t$.  Suppose that $\sup_{\omega \in \Omega, \; s \in [0,T]} |\xi_ s(\omega)| < \infty $. Then,
 \begin{equation}\label{rt}
\mathbb E\int_0^t \xi_{s-} \, d\pi_s=\lambda \int_0^t \mathbb E\xi_s ds.
\end{equation}
\end{lemma}

\subsection{Proof of 
Theorem \ref{uno}}
\label{SEC_PERTURBATIVO}
According to the notations appeared in Proposition \ref{PROP_WP_DRIFTLESS}, let  $v=PDE(Q ,f)$ and  ${w=PDE(Q+Q' ,f)}$ be  the unique solutions
 of \textcolor{black}{equations} \eqref{eq:Cauchy_problem_gen} and \eqref{w}, respectively.\newline
The proof of Theorem \ref{uno} will be obtained adapting the method developed in \cite{Krylov:Priola17} (see in particular Section $3$, therein).
Let $e_1$ be the first unit vector in $\R^N$. We define
$$
 X_t = \int_{0}^t \sqrt{Q'(r)} \, e_1 d \pi_r
$$
where $\sqrt{Q'(t)}$ is the unique $N \times N$ symmetric, non-negative definite square root of $Q'(t)$ and  $\{\pi_t\}_{t\ge0}$ is a Poisson Process of intensity $\lambda$ (cf. Equation \eqref{s25}). The parameter $\lambda$ will be chosen {appropriately} later on.\newline
Recall that the solution  $v$ to \eqref{eq:Cauchy_problem_gen} is given by
\begin{equation}\label{s55}
v(t,z) = \int_0^t ds\int_{\R^N} [f(s,z+ z') \mu_{\textcolor{black}{s,t}} (dz'),
\end{equation}
where $\mu_{\textcolor{black}{s,t}}$ is the Gaussian law of the stochastic integral $I_{s,t}:=\sqrt 2\int_s^t Q(v)^{1/2}dW_v$ (see the proof of Proposition \ref{PROP_WP_DRIFTLESS}).

Let us fix  $\epsilon>0$.  We notice   that the shifted source $f_\epsilon(t,z):=f(t,z-\epsilon X_t)$ (which also depends on $\omega$; we have omitted to write  such dependence on $\omega$) is again in $B_b\left(0,T;C^\infty_0(\R^N)\right)$.  \textcolor{black}{This is} the reason why we considered such a function space for the source. \textcolor{black}{It precisely allows to take into account the time discontinuities coming from the jumps of the Poisson process}.\newline
For any fixed $\omega$ in $\Omega$, Proposition \ref{PROP_WP_DRIFTLESS} readily gives that there exists a unique  solution $v_\epsilon=$  PDE$(Q,f(t,z-\epsilon X_t))$, depending also on $\epsilon$ and  $\omega$ as parameters, such that
\begin{align}
\label{proof:estimate1}
\sup_{(t,z)\in [0,T]\times \R^N}| v_\epsilon(t,z)| \, &\le \, T \sup_{(t,z)\in [0,T]\times \R^N}| f(t,z)|.
\end{align}
Moreover,  thanks to the invariance for translations of the $L^p$-norms, \textcolor{black}{it follows from \eqref{s22}} that
\begin{equation}\label{w11}
\| R(t)^* D^2 v_{\epsilon} R(t)  \|_{L^p((0,T)\times \R^N,\mathfrak m)}\, \le \, C_p  \|   f_\epsilon  \|_{L^p((0,T)\times \R^N,\textcolor{black}{\mathfrak m})}\, \le \,  C_p  \|   f  \|_{L^p((0,T)\times \R^N,\mathfrak m)}.
\end{equation}
 By \textcolor{black}{Equation} \eqref{s55}, we know that $v_{\epsilon} $ is given by
\begin{gather*}
v_{\epsilon}(t,z) = \int_0^t \int_{\R^N} [f(s,z - \epsilon X_s + z')\, \mu_{\textcolor{black}{s,t}} (dz')ds.
\end{gather*}
For each $z \in \R^N$, the stochastic process  $ (\textcolor{black}{v_{\epsilon}}(t,z))_{t \in [0,T]} $ has continuous paths ($\mathbb{P}$-a.s.) and is ${\cal F}_t$-adapted where
 ${\cal F}_t$ is the completed $\sigma$-algebra generated by the random variables $\pi_s$, $0 \le s \le t$.\newline
For  fixed $z \in \R^N$, let us introduce the process $(\textcolor{black}{v_\epsilon}(t,z+\epsilon X_t))_{t \in [0,T]}$ which is given by
$$
 v_{\epsilon}(t,z+\epsilon X_t) = \int_0^t \int_{\R^N} [f(s,z +\epsilon X_t - \epsilon X_s + z')\, \mu_{\textcolor{black}{s,t}} (dz')ds.
$$
\textcolor{black}{It is not difficut to check that it is} is  ${\cal F}_t$-adapted and \textcolor{black}{it has bounded and c\`adl\`ag paths}.\newline
Applying \eqref{INTEGRAL} on each interval $[\sigma_m,\sigma_{m+1}\wedge t ), m\in \{0,\cdots, \pi_t] $ on which $X_s$ is constant,  one then derives that:
\begin{equation}
\label{eq:1_OU}
v_\epsilon(t,z+\epsilon X_t) \, = \, \int_0^t\left(\text{tr}(Q(s)D^2_z v_\epsilon(s,z+\epsilon X_s))+ f(s,z)\right)\, ds +\int_0^t g_\epsilon(s,z) \,d\pi_s,
\end{equation}
where $g_\epsilon(s,z)=v_\epsilon(s,z+  \epsilon \sqrt{Q'(s)} \, e_1 + \epsilon X_{s-})-v_\epsilon(s,z+\epsilon X_{s-})$ is precisely the contribution associated with the jump  times. It is clear that $g_\epsilon(s,z)\neq 0$ if and only if $\pi_s$ has a jump at time $s$.
We then have by Lemma \ref{LEMME_POISSON}:
\begin{equation}
\label{proof:eq1}
    \mathbb{E}\int_0^tg_\epsilon(s,z)\, d\pi_s \, = \, \lambda \int_0^t\left(\bar v_\epsilon(s,z+\epsilon \sqrt{Q'(s)} \, e_1)-v_\epsilon(s,z)\right)\, ds,
\end{equation}
where $\bar v_\epsilon(s,z) = \mathbb{E}[v_\epsilon(s,z+\epsilon X_s)]$. Let \textcolor{black}{us denote}
$$
l(t) := \sqrt{Q'(t)} \, e_1.
$$
Taking the expectation on both sides of equation \eqref{eq:1_OU}, we find out that $\bar v_\epsilon$ is an integral solution of the following PDE:
\begin{equation}
\label{eq:2_OU}
    \partial_t \bar v_\epsilon(t,z) \, = \, \text{tr}(Q(t)D^2_z \bar v_\epsilon(t,z))+\lambda\left(\bar v_\epsilon(t,z+\epsilon l(t))-\bar v_\epsilon(t,z)\right)+f(t,z),
\end{equation}
with zero initial condition. Remark that uniqueness of bounded continuous solutions to \eqref{eq:2_OU} follows by the maximum principle,
arguing
as in the proof of Lemma 2.2 in \cite{Krylov:Priola17} (first one considers the case  $ \lambda T \le 1/4$ and then one iterates the procedure by steps of size $1/(4 \lambda)$). Moreover, by \eqref{w11} we obtain (using also the Jensen inequality and  the Fubini theorem) that
\[\begin{split}
\| R(t)^* D^2 \bar v_{\epsilon}  \, R(t)  \|_{L^p((0,T)\times \R^N,\textcolor{black}{\mathfrak m})}^p
\, &= \,  \int_{(0,T) \times \R^N} |R(t)^* D^2 \bar v_{\epsilon} (t,z) R(t)|^p \textcolor{black}{\mathfrak m(dt,dz)}
\\
&=\,  \int_0^T\int_{\R^N} |\mathbb{E} [R(t)^* D^2 v_{\epsilon} (t,z + \epsilon X_t) R(t)] \,|^p dz g(t)dt
\\
&\le   \int_0^T\int_{\R^N} \mathbb{E} [ | R(t)^* D^2 v_{\epsilon} (t,z + \epsilon X_t) R(t) \,|^p] dz g(t)dt
\\
&= \, \mathbb{E} \int_0^T\int_{\R^N}  | R(t)^* D^2 v_{\epsilon} (t,z + \epsilon X_t) R(t) \,|^p dz g(t)dt
\\
&= \, \textcolor{black}{\mathbb E}\int_0^T\int_{\R^N}  | R(t)^* D^2 v_{\epsilon} (t,\bar z ) R(t) \,|^p d\bar z g(t)dt
\\
&\le\,  C_p^p  \|   f  \|_{L^p((0,T)\times \R^N,\textcolor{black}{\mathfrak m})}^p,
\end{split}\]
using \eqref{w11} for the last inequality \textcolor{black}{($L^p$-estimate for the PDE with random source)}.
Choosing $\lambda = \epsilon^{-2}$ we have \textcolor{black}{
from \eqref{eq:2_OU} that
}
 \begin{equation}
\label{eq:2_OUe}
    \partial_t \bar v_\epsilon(t,z) \, = \, \text{tr}(Q(t)D^2_z\bar v_\epsilon(t,z))+\ \epsilon^{-2}\left(\bar v(t,z+\epsilon l(t))-\bar v_\epsilon(t,z)\right)+f(t,z),
\end{equation}
with zero initial condition and moreover
\begin{equation}\label{d3}
\| R(t)^* D^2 \bar v_{\epsilon}  R(t)  \|_{L^p((0,T)\times \R^N)}^p \, \le\,  C_p^p  \|   f  \|_{L^p((0,T)\times \R^N)}^p.
\end{equation}

Now, the idea is to apply again the same reasoning above to Equation \eqref{eq:2_OUe} with respect to $\bar v_\epsilon$ and $f(t,z+\epsilon X_t)$  again with  $\lambda=\epsilon^{-2}$. We obtain first a solution $p_{\epsilon}$ to \eqref{eq:2_OUe} corresponding to $f(t,z+\epsilon X_t)$ and we then \textcolor{black}{derive that}
$$
w_{\epsilon}(t,z) = \mathbb{E} [p_{\epsilon }(t , z - \epsilon X_t)]
$$
is the  unique bounded continuous (integral) solution $w_\epsilon$ of the following problem:
\begin{multline}
    \label{eq:3_OU}
\partial_t w_\epsilon(t,z) \,  = \, \text{tr}(Q(t)D^2 w_\epsilon(t,z))\\+\epsilon^{-2}\left[w_\epsilon(t,z+\epsilon l(t))-2w_\epsilon(t,z) +w_\epsilon(t,z-\epsilon l(t))\right] + f(t,z),
  \end{multline}
with initial condition $w_\epsilon(0,z) = 0$. Moreover,  the previous estimates  still hold  with $w_\epsilon$ instead of $v_\epsilon$, i.e.,
\begin{align}
\label{proof:estimate3_BIS}
\sup_{(t,z)\in[0,T]\times \R^N} |w_\epsilon(t,z)| \, &\le \, T\sup_{(t,z)\in[0,T]\times \R^N} |f(t,z)|;
\\
\label{proof:estimate4_BIS}
    \| R(t)^* D^2 w_{\epsilon}  \, R(t) \,  \|_{L^p((0,T)\times \R^N,\textcolor{black}{\mathfrak m})}  \, &\le\, C_p  \|   f  \|_{L^p((0,T)\times \R^N,\textcolor{black}{\mathfrak m})}.
\end{align}

We would like now to let $\epsilon$ goes to zero,
possibly passing to a subsequence $\epsilon_n \to 0$,
and prove that the associated limit $w$ solves
\begin{equation}
    \label{eq:3_OUf}
\begin{cases}
\partial_t w(t,z) \, = \, \text{tr}(Q(t)D^2 w(t,z))+ \langle D^2 w (t,z) \sqrt{Q'(t)}e_1,  \sqrt{Q'(t)}e_1 \rangle +f(t,z),\\
w(0,z) \, = \, 0
\end{cases}
\end{equation}
and  estimates \eqref{proof:estimate3_BIS} and \eqref{proof:estimate4_BIS} \textcolor{black}{hold} with $w_{\epsilon}$ replaced by $w$.\newline
To do so we will proceed by compactness. Namely, we are going to prove that the family of solutions $w_\epsilon$ solving \eqref{eq:3_OU}, indexed by the parameter $\epsilon$, is equi-Lipschitz on any compact subset of $[0,T]\times \R^N$ and the same holds for
  any derivative in space of $w_\epsilon$. Indeed, one can apply the finite difference operators with respect to $z$ at any order in \eqref{eq:3_OU}.
We recall that  for a \textit{smooth} function $\phi\colon \R^N\to \R$, the first finite difference $\delta_{h,i}\phi$, $i\in \{1,\cdots,N\} $ of step $h$ in the direction $e_i$ ($i^{\rm th}$ basis vector) is
given by
\[\delta_{h,i}\phi(z) \, = \, \frac{\phi(z+he_i)-\phi(z)}{h}, \quad z \in \R^N.\]For a given multi-index $\gamma \in \N^N$,  the $\gamma$-th order finite difference operator $\delta_{h,\gamma}$, is then defined, for any $h>0$, through composition. Namely,
\[\delta_{h,\gamma}\phi(z) \, = \, \delta^{\gamma_1}_{h,1}\delta^{\gamma_2}_{h,2}\dots \delta^{\gamma_N}_{h,N}\phi(z),\]
where $\delta^{\gamma_i}_{h,i}$ {denotes} the $\gamma_i$-th times composition of $\delta_{h,i}$ with itself.
\newline
Since any spatial derivative of $f$ belongs to $B_b\left(0,T;C^\infty_0(\R^N)\right) $, using \eqref{proof:estimate3_BIS} we deduce first that any finite difference of any order of $ w_\epsilon$ is bounded. Consequently, $ w_\epsilon$ is infinitely differentiable in space with bounded derivatives on $[0,T]\times \R^N $. Equation \eqref{eq:3_OU}, to be understood in its integral form similarly to \eqref{INTEGRAL}, then gives that those derivatives are themselves Lipschitz continuous in time. This precisely gives the equi-Lipschitz on any compact subset of $[0,T]\times \R^N$ of the family $w_\epsilon$ and any spatial derivative.
\newline
We can now apply the Arzel\`a-Ascoli theorem to $w_\epsilon$ showing the existence of a sub-sequence $\{w_{\epsilon_n}\}_{n \in \N}$ which converges uniformly on any compact set to a function $w\colon [0,T]\times \R^N\to \R$. Similarly, any derivative in space of $w_{\epsilon_n}$ tend to the respective derivatives of $w$, uniformly on the compact sets.\newline
Passing to the limit as $n\to \infty$ along the sequence $(\epsilon_n)_n $  in Equation \eqref{eq:3_OU}  (written in the  integral form), we can then conclude that $w$ solves \eqref{eq:3_OUf}.\newline
Moreover,  estimates \eqref{proof:estimate3_BIS} and \eqref{proof:estimate4_BIS} holds with $w_{\epsilon}$ replaced by $w$. Iterating the previous argument in $N$  steps we finally prove that the unique solution $w$ to
 \begin{equation}
    \label{eq:3_OUff}
\begin{cases}
\partial_t w(t,z) \, = \, \text{tr}(Q(t)D^2 w(t,z))+  \sum_{k=1}^N\langle D^2 w (t,z) \sqrt{Q'(t)}e_k,  \sqrt{Q'(t)}e_k \rangle +f(t,z),\\
w(0,z) \, = \, 0
\end{cases}
\end{equation}
verifies   estimates \eqref{proof:estimate3_BIS} and \eqref{proof:estimate4_BIS}  with $w_{\epsilon}$ replaced by $w$.
 The proof is complete.  \qed

\def\ciao{

%

\begin{prop}
\label{prop:link_ou_diffusive}
Let $f$ be in $B_b\left(0,T;C^\infty_0(\R^N)\right)$. Then, a function $u$ is the unique bounded, continuous solution of \textcolor{black}{the} Cauchy Problem \eqref{eq:def_ou_operator} if and only if the function $\tilde{u}$ given by $\tilde{u}(t,x):=u(t,e^{tA}x)$ is the unique bounded, continuous solution of the following:
\begin{equation}
\label{eq:PDE}
    \begin{cases}
    \partial(t) \tilde{u}(t,x) \, = \, \text{tr}\left(e^{-tA}Be^{-tA^\ast}D^2_x\tilde{u}(t,x)\right) +\tilde{f}(t,x), &(t,x)\in (0,T]\times \R^N;\\
    \tilde{u}(0,x) \, = \, 0, &x\in \R^N;
    \end{cases}
\end{equation}
where the function $\tilde{f}$ is given by $\tilde{f}(t,x)=f(t,e^{tA}x)$.
\end{prop}

The above result links the Cauchy Problem \eqref{eq:def_ou_operator}
involving Ornstein-Uhlenbeck operators with the diffusion one \eqref{eq:Cauchy_problem_gen}. In particular, it allows us to exploit Theorem \ref{thm:general} for the  Ornstein-Uhlenbeck operators.

From the assumed estimate on $u$ in \eqref{A_PRIO_EST}, we derive that for $\tilde u $ as in Proposition \ref{prop:link_ou_diffusive}:
\begin{align}\label{SHIFTED_A_PRIORI}
[u]_1:=[u(t,x)]_1=[\tilde u(t,e^{-t A}x)]_1=:\leftB \tilde u \rightB_1\le C[f]_2=[\tilde f(t,e^{-t A}x)]_2=:\leftB \tilde f\rightB_2,
\end{align}
with $\tilde{f}(t,x)=f(t,e^{tA}x) $ again as in Proposition \ref{prop:link_ou_diffusive}. The idea is now to apply Theorem \ref{thm:general}. To this end we first need to prove that $\leftB\cdot\rightB_1, \leftB\cdot\rightB_2 $ defined in \eqref{SHIFTED_A_PRIORI} are actually suitable semi-norms in the sense of Definition \ref{DEF_SUITABLE_NORM}. This property readily follows from the previous definition and the fact that $[\cdot]_1, [\cdot]_2 $ themselves are suitable seminorms.

It is also clear that, for $S$ as in the theorem, there exists a unique solution integral $u_S$ to \eqref{eq:OU_PERT} which is continuous and twice continuously differentiable in space. Similarly, we can associate to $ u_S$ the mapping $\tilde u_S(t,x)=u_S(t,e^{-At}x) $ which solves:
\begin{equation}
\label{eq:PDE_PERT}
    \begin{cases}
    \partial_t \tilde{u}_S(t,x) \, = \, \text{tr}\left(e^{-tA}(B+S(t))e^{-tA^\ast}D^2_x\tilde{u}_S(t,x)\right) +\tilde{f}(t,x), &\mbox{ on }[0,T]\times \R^N;\\
    \tilde{u}_S(0,x) \, = \, 0, &\mbox{ on }\R^N;
    \end{cases}.
\end{equation}
In short, with the notations of Section \ref{SEC_DRIFT_LESS}, we can write  $\tilde u=PDE(Q(t),\tilde f) $ with $Q(t)=e^{-tA}Be^{-tA^\ast} $ and $\tilde u_S=PDE(Q(t)+Q(t)',\tilde f) $ with $Q(t)'=e^{-tA}(B)e^{-tA^\ast} $. It is clear that $Q(t),Q(t)' $ are non-negative. From \eqref{SHIFTED_A_PRIORI}
and Theorem \ref{thm:general} we thus obtain:
\begin{align}\label{SHIFTED_A_PRIORI_FINAL}
\leftB \tilde u_S \rightB_1\le C\leftB \tilde f\rightB_2,
\end{align}
with the same previous constant $C$. From the definition of $(\leftB\cdot\rightB_i)_{i\in \{1,2\}} $ in \eqref{SHIFTED_A_PRIORI} we eventually derive:
$$[u_S]_1\le C[f]_1, $$
which concludes the proof.

\subsection{Concrete seminorms considered}\label{CONCRETE_SN}
In the final section of our article, we briefly present two natural examples of \emph{suitable} seminorms for which we can apply Theorem \ref{thm:OU_general} when the matrices $A,B$ satisfy the Kalman condition \textbf{[K]}.\newline
In order to introduce them, we firstly need to talk about the anisotropic nature of degenerate Ornstein-Uhlenbeck operators satisfying \textbf{[K]}. Intuitively, the appearance of this kind of phenomena is due essentially to the particular structure of the matrix $A$
that allows the smoothing effect of the second order operator $\text{tr}(BD^2_x)$, acting only on the first "component" given by $B_0$, to propagate into the system with lower and lower intensity. We recall indeed from \cite{Lanconelli:Polidoro94} that for degenerate Ornstein-Uhlenbeck operators, the Kalman condition is equivalent to the fact that $A,B$ must have the form:
\begin{equation}\label{eq:Lancon_Pol}
B \, = \,
    \begin{pmatrix}
        B_0 & 0_{d',N-d'}   \\
        0 _{N-d',d'}&   0_{N-d',N-d'}  \\
    \end{pmatrix}
\,\, \text{ and } \,\, A \, = \,
    \begin{pmatrix}
        \ast   & \ast  & \dots  & \dots  & \ast   \\
         A_2   & \ast  & \ddots & \ddots  & \vdots   \\
        0      & A_3   & \ast  & \ddots & \vdots \\
        \vdots &\ddots & \ddots& \ddots & \ast \\
        0      & \dots & 0     & A_n    & \ast
    \end{pmatrix}
\end{equation}
where $B_0$ is a non-degenerate matrix in $\R^{d'}\otimes \R^{d'}$ and $A_h$ are matrices in $\R^{d_h}\otimes \R^{d_{h-1}}$ with
$\text{rank}(A_h)=d_h$ for any $h$ in $\llbracket 2,n\rrbracket$. Moreover, $d'\ge d_2\ge \dots\ge d_n\ge 1$ and $d'+\dots+d_n=N$.

\paragraph{Anisotropic Geometry of Dynamics.} \textcolor{black}{Da S. a E.: forse da accorciare o sopprimere, la parte cinetica intendo.}
In order to illustrate how degenerate Ornstein-Uhlenbeck operators usually behave, we focus for the moment on the prototypical example we have already considered in the introduction, the Kolmogorov operator ${\Delta}_{\text{K}}$.\newline
Fixed $N=2d$ and $n=2$ we recall that $\Delta_{\text{K}}$ can be represented for any sufficiently regular function $\phi\colon \R^N\to \R$ as
\[\Delta_{\text{K}}\phi(x) \, = \, \Delta_{x_1}\phi(x)+ x_1\cdot D_{x_2}\phi(x) \quad \text{ on }\R^{2d},\]
where $(x_1,x_2)\in\R^{2d}$ and $\Delta_{x_1}$ is the standard Laplacian with respect to $x_1$. In our framework, it is associated with the matrices
\[A_{\text{K}} \, := \,
    \begin{pmatrix}
               0 & 0 \\
               I_{d\times d} & 0
    \end{pmatrix} \,\, \text{ and } \,\, B \, := \,
    \begin{pmatrix}
                I_{d\times d} & 0 \\
                 0 & 0
    \end{pmatrix}.
             \]
We then search for a dilation
\[\delta_\lambda\colon [0,\infty)\times \R^{2d} \to [0,\infty)\times \R^{2d}\]
which is invariant for the considered dynamics, i.e.\ a dilation that transforms solutions of the equation
\[\partial_tu(t,x)-\Delta_{\text{K}} u(t,x) \, = \, 0 \quad \text{ on }(0,\infty)\times\R^{2d}\]
into other solutions of the same equation. The idea of a dilation operator $\delta_\lambda$ that summarizes the multi-scaled behaviour of the dynamics was firstly introduced by Lanconelli and Polidoro in \cite{Lanconelli:Polidoro94}.  Since then, it has
become a "standard" tool in the analysis of degenerate equations (see e.g. \cite{Lunardi97}, \cite{Huang:Menozzi:Priola19} or  \cite{Marino20}, \cite{Marino21} in the fractional setting). \newline
Due to the particular sub-diagonal structure of $A_{\text{K}}$ and the natural scaling of $\Delta_{x_1}$ of order $2$, it makes sense to consider for any fixed $\lambda>0$, the
following
\[ \delta_\lambda(t,x_1,x_2)\, :=\,  (\lambda^2 t,\lambda x_1,\lambda^3 x_2).\]
It then holds that
\[\bigl(\partial_t -\Delta_{\text{K}}\bigr) u = 0 \, \Longrightarrow \bigl(\partial_t -\Delta_{\text{K}} \bigr)(u \circ
\delta_\lambda) = 0.\]
The previous reasoning suggests to introduce a parabolic distance $\mathbf{d}_P$ that is homogenous with respect to the dilation $\delta_\lambda$,
so that:
\[\mathbf{d}_P\bigl(\delta_\lambda(t,x);\delta_\lambda(s,x')\bigr) = \lambda \mathbf{d}_P\bigl((t,x);(s,x')\bigr).\]
Following the notations in \cite{Huang:Menozzi:Priola19}, we then introduce the distance $\mathbf{d}_P$ on $[0,\infty)\times \R^{2d}$  given by
\begin{equation}\label{Definition_distance_d_P}
\mathbf{d}_P\bigl((t,x),(s,x')\bigr)  \, := \, \vert s-t\vert^\frac{1}{2}+\mathbf{d}(x,x') \, := \,
\vert s-t\vert^\frac{1}{2} +|x_1-x_1'|+|x_2-x_2'|^{\frac{1}{3}}.
\end{equation}
As it can be seen, $\mathbf{d}_P$ is an extension of the standard parabolic distance adapted to respect the multi-scale nature
of the underlying operator. Indeed, the exponents appearing in \eqref{Definition_distance_d_P} are those which make each space component homogeneous to
the characteristic time scale $t^{1/2}$.\newline
From a more
probabilistic point of view, the exponents in Equation \eqref{Definition_distance_d_P} can also be related to the characteristic
time scales of the iterated integrals of a Brownian Motion. It can be easily seen from the example, noticing that the operator $\Delta_{\text{K}}$ corresponds to the generator of the Brownian Motion and its time integral.

Going back to the general framework with $A,B$ as in \eqref{eq:Lancon_Pol}, we can now construct the suitable distance $\mathbf{d}$ associated with the Cauchy Problem \eqref{eq:OU_initial} on $\R^N$.\newline
Recalling that $n$ denotes the smallest integer such that the Kalman rank condition [\textbf{K}] holds, we follow
\cite{Lunardi97} decomposing the space $\R^N$ with respect to the family of linear operators $B, AB,\dots,
A^{n-1}B$. \newline
We start defining the family $\{V_h\colon h\in \llbracket 1,n \rrbracket\}$ of subspaces of $\R^N$ through
\[
V_h \,  := \, \begin{cases}
            \text{Im} (B), & \mbox{if } h=1, \\
            \bigoplus_{k=1}^{h}\text{Im}(A^{k-1}B), & \mbox{otherwise}.
        \end{cases}\]
It is easy to notice that $V_h\neq V_k$ if $k\neq h$ and $V_1\subset V_2\subset\dots V_n=\R^N$. We can then construct iteratively the family
$\{E_h \colon h\in \llbracket 1,n \rrbracket\}$ of orthogonal projections from $\R^N$ as
\[E_h \,  := \,
        \begin{cases}
            \text{projection on } V_1, & \mbox{if } h=1; \\
            \text{projection on }(V_{h-1})^\perp \cap V_h, & \mbox{otherwise}.
        \end{cases}
\]
With a small abuse of notation, we will identify the
projection operators $E_h$ with the corresponding matrices
in $\R^N\otimes \R^N$. Actually, we have that $\dim E_h(\R^N)=d_h$ for any $h\ge 1$, where the integer $d_h$ already appeared in equation \eqref{eq:Lancon_Pol}.\newline
We can then define the distance $\mathbf{d}$ through the decomposition $\R^N=\bigoplus_{h=1}^nE_h(\R^N)$ as
\[\mathbf{d}(x,x') \,:= \, \sum_{h=1}^{n}\vert E_h(x-x')\vert^{\frac{1}{1+2(h-1)}}.\]
We highlight however that, technically speaking, $\mathbf{d}$ does not induce a norm on $[0,T]\times \R^N$ in the usual sense since it lacks of
linear homogeneity. \newline
The anisotropic distance $\mathbf{d}$ can be understood direction-wise: we firstly fix a "direction" $h$ in $\llbracket
1,n\rrbracket$ and then calculate the standard Euclidean distance on the associated subspace $E_h(\R^N)$, but scaled according to the
dilation of the system in that direction. We conclude summing the contributions associated with each component.\newline
Highlighted the particular anisotropic nature of the operator through the distance $\bm{d}$, we can now extend two standard functional spaces, the H\"older and the Sobolev ones, in order to reflect the intrinsic multi-scale behavior given by the dilation operator $\delta_\lambda$.
\textcolor{black}{Da S. a E.: veramente forse dovremmo parlare della parte infatti omogenea che d\`a gli esponenti.}

\paragraph{Anisotropic H\"older spaces.}
Fixed some $k$ in $\N_0:=\N\cup\{0\}$ and $\beta$ in $(0,1]$, we follow Krylov \cite{book:Krylov96} denoting the Zygmund-H\"older semi-norm for a function $\phi\colon \R^N\to \R$ as
\[[\phi]_{C^{k+\beta}} \, := \,
\begin{cases}
    \sup_{\vert \vartheta \vert= k}\sup_{x\neq y}\frac{\vert D^\vartheta\phi(x)-D^\vartheta\phi(y)\vert}{\vert x-y\vert^\beta} , & \mbox{if
    }\beta \neq 1; \\
    \sup_{\vert \vartheta \vert= k}\sup_{x\neq y}\frac{\bigl{\vert}D^\vartheta\phi(x)+D^\vartheta\phi(y)-2D^\vartheta\phi(\frac{x+y}{2})
    \bigr{\vert}}{\vert x-y \vert}, & \mbox{if } \beta =1.
                             \end{cases}\]
Consequently, The Zygmund-H\"older space $C^{k+\beta}_b(\R^N)$ is the family of functions $\phi\colon \R^N
\to\R$ such that $\phi$ and its derivatives up to order $k$ are continuous and the norm
\[\Vert \phi \Vert_{C^{k+\beta}_b} \,:=\, \sum_{i=1}^{k}\sup_{\vert\vartheta\vert = i}\Vert D^\vartheta\phi
\Vert_{\infty}+[\phi]_{C^{k+\beta}_b} \,\text{ is finite.}\]

We can define now the anisotropic Zygmund-H\"older spaces associated with the distance $\mathbf{d}$. Fixed $\gamma>0$, the space $C^{\gamma}_{b,d}(\R^N)$ is
the family of functions $\phi\colon \R^N\to \R$ such that for any $h$ in $\llbracket 1,n\rrbracket$ and any $x_0$ in $\R^N$, the function
\[z\in  E_h(\R^N)\,  \to \, \phi(x_0+z) \in \R \,\text{ belongs to }C^{\gamma/(1+2(h-1))}_b\left(E_h(\R^N)\right),\]
with a norm bounded by a constant independent from $x_0$. It is endowed with the norm
\begin{equation}\label{eq:def_anistotropic_norm}
\Vert\phi\Vert_{C^{\gamma}_{b,d}} \,:=\,\sum_{h=1}^{n}\sup_{x_0\in \R^N}\Vert\phi(x_0+\cdot)\Vert_{C^{\gamma/(1+2(h-1))}_b}.
\end{equation}
We highlight that it is possible to recover the expected joint regularity for the partial derivatives, when
they exist, as in the standard H\"older spaces. In such a case, they actually turn out to be H\"older continuous with respect to the
distance $\mathbf{d}$ with order one less than the function (See Lemma $2.1$ in \cite{Lunardi97} for more details).
\newline
Exploiting that for the "standard" Zygmund-H\"older norm $\Vert\cdot\Vert_{L^\infty((0,T),C^{\gamma}_{b})} $ is clearly suitable (in the sense of Definition \ref{DEF_SUITABLE_NORM}) and reasoning components by components, it is not difficult to conclude that also the anisotropic norm $\Vert\cdot\Vert_{L^\infty((0,T),C^{\gamma}_{b,d})}$ is suitable as well.\newline
Moreover, Lunardi in \cite{Lunardi97} showed that the following anisotropic Schauder estimates holds for the solution $u$ of Cauchy Problem \eqref{eq:OU_initial}:
\begin{equation}\label{SCHAU_OU}
\Vert u \Vert_{L^\infty((0,T),C^{2+\beta}_{b,d})} \, \le \, C \Vert f \Vert_{L^\infty((0,T),C^{\beta}_{b,d})},
\end{equation}
for some constant $C$ independent from $f$.
Thus, all the assumptions of Theorem \ref{thm:OU_general} are satisfied with $[\cdot]_1 
=\Vert\cdot\Vert_{L^\infty( (0,T),C^{2+\beta}_{b,d})}$ and $[ \cdot ]_2=\Vert\cdot\Vert_{L^\infty((0,T),C^{\beta}_{b,d})}$.

\paragraph{Anisotropic Sobolev spaces.}
Given $\alpha$ in $(0,1)$ and $h$ in $\llbracket 1, n \rrbracket$, we want to introduce the $\alpha$-fractional Laplacian $\Delta^\alpha_{x_h}$ along the "component" $x_h$.
To do it, we follow \cite{Huang:Menozzi:Priola19} by considering the orthogonal projection $p_h\colon \R^N\to \R^{d_h}$ such that $p_h(x)=x_h$ and denoting its adjoint by $B_h=p_h\colon \R^{d_h}\to \R^N$. Notice that $B=B_1B^\ast_1$ in a matricial form.\newline
We can now define the $\alpha$-fractional Laplacian $\Delta^\alpha_{x_h}$ as:
\[\Delta^{\alpha}_{x_h}\phi(x) \, := \, \text{p.v.}\int_{\R^{d_i}}\left[\phi(x+B_hz)-\phi(x)\right] \frac{dz}{|z|^{d_h+\alpha}}, \quad x \in \R^N,\]
for any sufficiently regular function $\phi \colon \R^N\to \R$.\newline
Fixed $p$ in $(1,+\infty)$, we recall that we have denoted by $L^p([0,T]\times \R^N)$ the standard $L^p$-space with respect to the Lebesgue measure.\newline
We can now define the anisotropic Sobolev space associated with the distance $\bm{d}$. For notational simplicity, let us denote \[\alpha_h \,:=\, \frac{1}{1+\alpha(h-1)}.\]
The space $\dot{W}^{2,p}_d([0,T]\times \R^N)$ is composed by all the functions $\varphi\colon [0,T]\times\R^N\to \R$ in $L^p([0,T]\times \R^N)$ such that for any $h$ in $\llbracket 1, n \rrbracket$,
\[\Delta^{\alpha_h}_{x_h}\varphi(t,x) \, := \,  \Delta^{\alpha_h}_{x_h}\varphi(t,\cdot)(x) \text{ belongs to } L^p([0,T]\times \R^N).\]
It is endowed with the natural norm
\[\Vert \varphi \Vert_{\dot{W}^{2,p}} \, = \, \sum_{h=1}^n\Vert \Delta^{\alpha_h}_{x_h}\varphi \Vert_{L^p}.\]
Similarly to the first example, we can reason componentwise in order to show that the anisotropic Sobolev norm $\Vert \cdot \Vert_{\dot{W}^{2,p}_d}$ is indeed a suitable seminorm in our definition.\newline
In \cite{Huang:Menozzi:Priola19}, see also \cite{Chen:Zhang19}  and \cite{Menozzi18} where time inhomogeneous coefficients are considered as well, it has been proven that when the strictly upper diagonal elements of $A$ in \eqref{eq:Lancon_Pol} are equal to zero then the following Sobolev estimates hold:
\begin{equation}\label{SOB_EST_GLOB}
\Vert u \Vert_{\dot{W}^{2,p}_d} \, \le \, C_p\Vert f \Vert_{L^p},
\end{equation}
where again $u$ is the unique bounded solution of Cauchy problem \eqref{eq:OU_initial}. In particular we get also the maximal smoothing effects w.r.t. the degenerate directions.
The specific structure assumed on $A$ is actually due
to the fact that for such matrices there is an underlying homogeneous space structure which makes easier to establish maximal regularity estimates (see e.g. \cite{book:Coifman:Weiss} in this general setting). Anyhow, in this case the assumptions of Theorem \ref{thm:OU_general} are satisfied with $[\cdot]_1
=\Vert\cdot\Vert_{\dot{W}^{2,p}_d}$ and $
[\cdot]_2=\Vert\cdot\Vert_{L^p}$.

For a general $A$ with the form in \eqref{eq:Lancon_Pol} such that $A$, $B$ satisfy \textbf{[K]}, we can refer to \cite{Bramanti:Cupini:Lanconelli:Priola10} (\textcolor{black}{Argomento da precisare per passare alla $f$, un taglio sfruttando la densit\`a esplicita e le sue propriet\`a}) to assert that the unique bounded solution of Cauchy problem \eqref{eq:OU_initial} actually satisfies
\begin{equation}\label{LP_OU_NON_DEG}
 \Vert D_{x_1}^2u \Vert_{L^p} \, \le \, C_p\Vert f \Vert_{L^p}.
\end{equation}
 We believe that the approach in \cite{Bramanti:Cupini:Lanconelli:Priola10} could extend to show that for \eqref{SOB_EST_GLOB} still hold in this general setting, but such estimates have not been, up to our best knowledge, proven yet.  In this case, to apply Theorem \ref{thm:OU_general} we consider $[\cdot]_1
=\Vert D_{x_1}^2\cdot\Vert_{L^{p}}$ and $
[\cdot]_2=\Vert\cdot\Vert_{L^p}$.

From the above information, namely \eqref{SCHAU_OU}, \eqref{SOB_EST_GLOB}, \eqref{LP_OU_NON_DEG}, and Theorem \ref{thm:OU_general} we obtain the following final result:

\begin{theorem}\label{CONCRETE_PERTURBATION_ESTIMATES} Let $A$, $B$ be two matrices of the form \eqref{eq:Lancon_Pol} such that the Kalman condition \textbf{[K]} holds.  Then, for $f$ in $B_b\left(0,T;C^\infty_0(\R^N)\right)$, the  unique bounded, continuous solution $u\colon [0,T]\times \R^N\rightarrow \R$ in integral form of the associated Cauchy Problem \eqref{eq:OU_initial} which is also smooth in space is such that:
\begin{trivlist}
\item[-] For all $\beta\in (0,1) $, there exists $C_\beta:=C(T,N,\beta)$ such that:
$$ \sup_{t\in[0,T]}\Vert u(t,\cdot)\Vert_{C^{2+\beta}_{b,d}} \, \le \, C_\beta \sup_{t \in (0,T)}\Vert f(t,\cdot) \Vert_{C^\beta_{b,d}}.$$
\item[-] For $p\in (1,+\infty) $, there exists $C_p:=C(T,N,p)$ such that:
$$\Vert D_{x_1}^2u \Vert_{L^p} \, \le \, C_p\Vert f \Vert_{L^p}.$$
If additionally $A_{ij}=0$ for $ j>i$ then
$$\Vert u \Vert_{\dot{W}^{2,p}_d} \, \le \,  C_p \Vert f \Vert_{L^p}.$$
\end{trivlist}
Let now $S:[0,T]\mapsto \R^N\otimes \R^N$ be a continuous matrix valued mapping s.t. for all $t\in [0,T] $, $S(t) $ is symmetric and non-negative. Then, the unique
bounded, continuous solution $u_S\colon [0,T]\times \R^N\rightarrow \R$ in integral form of the associated Cauchy Problem \eqref{eq:OU_PERT} which is also smooth in space also satisfies that:
\begin{trivlist}
\item[-] For all $\beta\in (0,1) $,
$$ \sup_{t\in[0,T]}\Vert u_S(t,\cdot)\Vert_{C^{2+\beta}_{b,d}} \, \le \, C_\beta \sup_{t \in (0,T)}\Vert f(t,\cdot) \Vert_{C^\beta_{b,d}}.$$
\item[-] For $p\in (1,+\infty) $,
$$\Vert D_{x_1}^2u_S \Vert_{L^p} \, \le \, C_p\Vert f \Vert_{L^p}.$$
If additionally $A_{ij}=0$ for $ j>i$ then
$$\Vert u_S \Vert_{\dot{W}^{2,p}_d} \, \le \,  C_p \Vert f \Vert_{L^p},$$
\end{trivlist}
for the \textbf{same previous constants} $C_\beta,C_p$.

\end{theorem}

\section{On some related applications}
\fancyhead[RO]{Section \thesection. On some related applications}
In order to illustrate one possible  application of the above estimates we will consider the following stochastic differential equation:
\begin{equation}\label{SDE_ATTRITO}
dX_t=\big(AX_t +F(t,X_t)\big) dt+\left ( \begin{array}{c}\sigma_1(t,X_t) dW_t^1\\
\sigma_2(t) dW_t^2\end{array}\right),
\end{equation}
where $(W_t^1,W_t^2) $ are $d$-dimensional Brownian motions defined on some probability space  $(\Omega,\mathcal F,\mathbb P) $, possibly dependent.  $X_t=(X_t^1,X_t^2)\in \R^{2d} $. We assume that:

\begin{trivlist}
\item[-]  The  matrices
$A$ and $B=\left ( \begin{array}{cc} I_{d,d}&  0_{d,d}\\ 0_{d,d}& 0_{d,d}\end{array}\right) \in \R^{2d}\otimes \R^{2d} $ satisfy the Kalman condition \textbf{[K]}.

\item[-\textbf{[F]}] The non linear drift $F=(F_1,F_2)$, is measurable in time and s.t.
\begin{itemize}
\item[-]$\exists K>0,\ \forall t\in [0,T],\ |b(t,0)|\le K $.
\item[-] $x=(x_1,x_2)\mapsto F_1(t,x) $  $\beta^1 $-H\"older continuous uniformly in $t\in [0,T]$.
\item[-] $x=(x_1,x_2)\mapsto F_2(t,x)=F_2(t,x_2)$ is uniformly $\beta^2>1/3$ H\"older continuous uniformly in $t\in [0,T]$.
\end{itemize}
\item[-\textbf{[UE]}]The diffusion coefficient $a_1=\sigma_1\sigma_1^*$ is uniformly elliptic and bounded, i.e.
there exists $\kappa\ge 1$ s.t. for all $t\in [0,T], x\in \R^{2d},\ \xi\in \R^d $:
\begin{equation}
\label{UE}
\kappa^{-1}|\xi|^2 \le \langle a_1(t,x)\xi,\xi\rangle \le \kappa |\xi|^2.
\end{equation}
\item[-\textbf{[L]}] Local condition for the diffusion. There exists  a measurable function $\varsigma:[0,T]\rightarrow {\mathcal S}_d $ (symmetric matrices of dimension $d$)
satisfying \textbf{[UE]} and such that
\begin{eqnarray}
\label{LOC_COND}
\varepsilon_{a_1}&:=&\sup_{0\le t \le T}\sup_{x\in \R^{2d} } |a_1(t,x)-\varsigma(t)|,
\end{eqnarray}
is \textit{small}. In particular, we do not assume any \textit{a priori} continuity of $a_1$.
\item[-\textbf{[C]}] The mapping $t\mapsto \sigma_2(t) $ is continuous.
\end{trivlist}
Importantly, we do not assume any non-degeneracy condition on the coefficient $\sigma_2$.

We have the following result.
\begin{theorem} \label{THM_SDE_ATTRITO}
Assume \textbf{[H]}, \textbf{[F]}, \textbf{[UE]}, \textbf{[L]}, \textbf{[C]} are in force. Then, the martingale problem associated with $(L_t)_{t\ge 0} $ where for $\varphi\in C_0^\infty(\R^{2d}) $, $(t,x)\in [0,T]\times \R^{2d} $,
\begin{align}
 L_t\varphi(x)=& \frac 12 {\rm Tr}\Big( a_1(t,x) D_{x_1}^2 \varphi(x)\Big)+\frac 12 {\rm Tr}\Big( a_2(t) D_{x_2}^2 \varphi(x)\Big)\notag\\
 &+\langle (Ax+F(t,x)),D\varphi(x)\rangle,\ a_2(t):=\sigma_2\sigma_2^*(t),
 \end{align}
is well posed. In particular, there exists a unique weak solution to equation \eqref{SDE_ATTRITO}. In particular, denoting by $\mathbb P$ the only solution of the martingale problem on $C([0,T],\R^{2d}) $, and denoting by $(X_t)_{t\ge 0}$ the associated canonical process, it holds that for given $p,q$ s.t. $4d/p+2/q<2 $, there exists $C_{p,q}\ge 1$ s.t.  for all $f\in L^q((0,T),L^p(\R^{2d}) $:
$$\mathbb E^{{\mathbb P}_{t,x}}[\int_t^T f(s,X_s)ds]\le C_{p,q}\|f\|_{L^q((0,T),L^p(\R^{2d}))} .$$
\end{theorem}
To prove Theorem \ref{THM_SDE_ATTRITO} we will proceed through a perturbative approach similarly to what was done in \cite{chau:meno:17}, \cite{Menozzi18}. To this end we first consider an underlying Gaussian proxy.
\subsection{Deterministic flow and Proxy Gaussian model}
 Introduce, for fixed $ T>0,\ y \in \R^{2d}$ and $t\in [0,T]$ the backward flow:
\begin{align}
\label{DET_SYST}
\overset{.}{\theta}_{t,T}(y)=&G(t,\theta_{t,T}(y)),\ \theta_{T,T}(y)=y\notag\\
G(t,z)=&Az+F(t,z),\ \in \R^{2d} .
\end{align}
\begin{remark}
We mention carefully that from  the Cauchy-Peano theorem, a solution to  \eqref{DET_SYST} exists. Indeed, the coefficients are continuous and have at most linear growth.
\end{remark}

\subsection{Linearized Multi-scale Gaussian Process and Associated Controls}\label{sec:freezing}
We now want to introduce the forward linearized flow around a solution of \eqref{DET_SYST}. Namely, we consider for $s\ge 0$ the deterministic ODE
\begin{equation}
\label{LIN_AROUND_BK_FLOW}
\overset{.}{\tilde{\phi}}_s =  A
\tilde{\phi}_s
+{ F}(s,\theta_{s,T}(y)).
\end{equation}

We consider now the Gaussian stochastic linearized dynamics $(\tilde X_s^{T,y})_{s\in [t,T]} $:
\begin{eqnarray}
&&\hspace*{-.5cm}d\tilde X_s^{T,y}=[ A\tilde X_s^{T,y}+F(s,\theta_{s,T}(y))]ds +\left(\begin{array}{c}\sigma_1(s,\theta_{s,T}(y)) dW_s^1\\
\sigma_2(s)dW_s^2,\end{array}\right),\nonumber\\
&& \hspace*{.75cm}\forall s\in  [t,T],\
 \tilde X_t^{T,y}=x. \label{FROZ}
 \end{eqnarray}

 From equations \eqref{LIN_AROUND_BK_FLOW}  we explicitly integrate \eqref{FROZ} to obtain for all $v\in [t,T] $:
 \begin{equation}
 \label{INTEGRATED}
\begin{split}
\tilde  X_v^{T,y}&=e^{A(v-t)}x+\int_t^v  e^{A(v-s)} F(s,\theta_{s,T}(y))ds\\
 & +\int_t^v e^{A(v-s)}\left(\begin{array}{c}\sigma_1(s,\theta_{s,T}(y)) dW_s^1\\
 \sigma_2(s)dW_s^2\end{array}\right).
 \end{split}
 \end{equation}
}

\section{Additional stability results}
\label{ESTENSIONI}

In this section we extend the previous approach to derive the stability with respect to a second order perturbation of the Ornstein-Uhlenbeck operator in \eqref{DEF_OU_OP_PROXY} under the Kalman condition \textbf{[K]}. Here, we consider also $L^p$-estimates involving the degenerate components of the operator and some associated Schauder estimates.

\subsection{Anisotropic Sobolev spaces and maximal Lp-regularity.}
With the notations of Section \ref{KH_OU} we write $z\in \R^N $ as $z=(x,y_2,\cdots, y_n) $ with $x\in\R^{d}$, $y_i\in \R^{d_i},\ i\in \{2,\cdots,n \}$, recalling also that $\sum_{i=2}^{n} d_i=d' $.\newline
Given $\beta $ in $(0,1)$ and $i$ in $\llbracket 2, n \rrbracket$, we want to introduce the $\beta$-fractional Laplacian $\Delta^\beta_{y_i}$ along the \textit{component} $y_i$.
To do so, we follow \cite{Huang:Menozzi:Priola19} by considering the orthogonal projection $p_i\colon \R^N\to \R^{d_i}$ such that $p_i(z)=p_i\big((x,y_2,\dots,y_n)\big)=y_i$ and denoting its adjoint by $E_i\colon \R^{d_i}\to \R^N$.
We can now define the $\beta$-fractional Laplacian $\Delta^\beta_{y_i}$ as:
\[\Delta^{\beta}_{y_i}\phi(z) \, := \, \text{p.v.}\int_{\R^{ d_i}}\left[\phi(z+E_i w)-\phi(z)\right] \frac{dw}{|w|^{d_i+\textcolor{black}{2}\beta}}, \quad z \in \R^N,\]
for any sufficiently regular function $\phi \colon \R^N\to \R$.\newline
Fixed $p$ in $(1,+\infty)$, we recall that we have denoted by $L^p((0,T)\times \R^N)$ the standard $L^p$-space with respect to the Lebesgue measure.\newline
We can now define the appropriate anisotropic Sobolev space to state our results.
For notational simplicity, let us denote
\begin{equation}\label{INDEXES}
\alpha_{i} \,:=\, \frac{1}{2i-1}.
\end{equation}
Set now $\alpha:=(\alpha_1,\cdots,\alpha_k) \in \R^{k}$. The \textit{homogeneous} space $\dot{W}^{2,p}_\alpha([0,T]\times \R^N)$ is composed by all the functions $\varphi\colon [0,T]\times\R^N\to \R$ in $L^p([0,T]\times \R^N)$ such that $(t,z)\in [0,T]\times\R^N\mapsto \Delta_x \varphi(t,z) \in L^p([0,T]\times \R^N) $, where $\Delta_x \varphi$ is intended  in distributional sense,   and for any $i$ in $\llbracket 2, n \rrbracket$, $\Delta^{\alpha_i}_{y_i}\varphi(t,z)$ is well defined for almost every $(t,z)$ and 
\[\Delta^{\alpha_i}_{y_i}\varphi(t,z) \, := \,  \Delta^{\alpha_i}_{y_i}\varphi(t,\cdot)(z)   \text{ belongs to } L^p([0,T]\times \R^N).\]
It is endowed with the natural \textit{semi}-norm $\Vert \varphi \Vert_{\dot{W}_\alpha^{2,p}}$  where 
 \begin{equation}\label{SEMI_NORM_SOB}
 \Vert \varphi \Vert_{\dot{W}_\alpha^{2,p}}^p \, = \, \|\Delta_x \varphi\|_{L^p}^p +\sum_{i=2}^n\Vert \Delta^{\alpha_i}_{y_i}\varphi \Vert_{L^p}^p.
\end{equation} 
The thresholds in \eqref{INDEXES} might seem awkward at first sight. They actually correspond to the indexes needed to get stability of the harmonic functions associated with the principal part of \eqref{DEF_OU_OP_PROXY}, that is considering $A_0$ consisting in the sub-diagonal part of $A$ only (i.e.\ considering \eqref{sotto} when the diagonal elements and the  strictly upper diagonal elements  are equal  to zero), along  an associated dilation operator. Namely, setting  \begin{equation}
\label{DEF_OU_OP_PROXY_0}
{L}_0^\text{ou} f(z) = 
{\rm Tr}(B D^2 f(z))
+ \langle A_0 z , D f(z)\rangle,\;\;\;  \; z \in \R^N,
\end{equation}
so that $A_0,B $ satisfy \textbf{[K]}, if it holds that
\[ (\partial_t -{L}_0^\text{ou})u(t,z)=0\]
then for all $ \lambda>0$, we have that
\[(\partial_t -L_0^\text{ou})u\big(\delta_\lambda(t,z) \big)\, = \, 0\] 
where the dilation operator
$$\delta_\lambda(t,z)\, = \, (\lambda^{1/2}t , \lambda x, \lambda^{1/3} y_1,\cdots, y_k^{1/(2n-1)})$$
precisely exhibits the exponents in \eqref{INDEXES} for the degenerate components.\newline 
In \cite{Huang:Menozzi:Priola19}, see also \cite{Chen:Zhang19}  and \cite{Menozzi18} where time inhomogeneous coefficients are considered as well, it has been proven that if $A$, $B$ satisfy \textbf{[K]} and the diagonal and the strictly upper diagonal elements of $A$ in \eqref{DEF_OU_OP_PROXY} are equal to zero (i.e.\ when $A=A_0$), then the following Sobolev estimates hold:
\begin{equation}\label{SOB_EST_GLOB}
\Vert u \Vert_{\dot{W}^{2,p}_\alpha} \, \le \, C_p\Vert f \Vert_{L^p},
\end{equation}
with $C_p:=C_p(\nu,A,d,d')$, where again $u$ is the unique bounded solution  to the corresponding Cauchy problem \eqref{eq:OU_initial:intro}. In particular we get also the maximal smoothing effects with respect to the degenerate directions. Note that the solution $u$ to \eqref{ko_K}
 verifies \eqref{SOB_EST_GLOB}.
The specific structure assumed on $A$ is actually due
to the fact that for such matrices there is an underlying homogeneous space structure which makes easier to establish maximal regularity estimates (see e.g. \cite{book:Coifman:Weiss} in this general setting).

If $A$, $B$ satisfy the Kalman condition \textbf{[K]} with a general $A$ as in \eqref{sotto}, having non zero strictly upper diagonal entries (non zero entries in the diagonal should  not create difficulties), we believe that the approach in \cite{Bramanti:Cupini:Lanconelli:Priola10} could extend to show that Estimates \eqref{SOB_EST_GLOB} still hold in this general setting. However, such estimates have not been, up to our best knowledge, proven yet.
 
\paragraph {\bf $L^p$-estimates for  the degenerate directions of Ornstein-Uhlenbeck operators.}
Setting, as in Section \ref{SEZ_STRAT}, $u(t,z) = v(t, e^{tA} z)$
and since $ u$ solves \eqref{eq:OU_initial:intro} we have that $v$ in turn solves \eqref{ma}. From  the previous computations and setting 
\[B_I\, := \, \begin{pmatrix}
 I_{d,d}  & 0_{d,d'}\\
 0_{d',d} & 0_{d',d'}
\end{pmatrix},\]
and considering $A$ as in \cite{Huang:Menozzi:Priola19}, with the diagonal and strictly upper diagonal elements of $A$ equal to zero in \eqref{sotto},
we derive that
\begin{align*}
\|D_x^2 u  \|_{L^p((0,T)\times \R^N)}\, &=\, \|  
B_Ie^{t A^*} D^2 v (t, e^{tA} \cdot ) e^{tA} B_I
\|_{L^p((0,T)\times \R^N)} \\
&\le\,  C_p  \|  \tilde f(t, e^{tA}  \cdot )  \|_{L^p((0,T)\times \R^N)}.
\end{align*}
On the other hand, for all $i\in \{2,\cdots,n\} $ and with $\alpha_i $ as in \eqref{INDEXES},
\begin{align*}
 \|\Delta_{y_i}^{\alpha_i} u  &\|_{L^p((0,T)\times \R^N)}^p\, =\,\int_0^T\int_{\R^N} \Big|{\rm p.v.}\int_{\R^{ d_i}}[u(t,z+E_iw)-u(t,z)] \frac{dw}{|w|^{ d_i+2\alpha_i}}\Big|^p\, dzdt\\
&\qquad\qquad= \, \int_0^T \int_{\R^N} \Big|{\rm p.v.}\int_{\R^{ d_i}}[v(t,e^{tA}(z+E_iw))-v(t,e^{tA}z)] \frac{dw}{|w|^{ d_i+2\alpha_i}}\Big|^p\, dzdt\\
&\qquad\qquad=: \,\|\Delta^{\alpha_i,i,A}v\|_{L^p((0,T)\times \R^N)}^p,
\end{align*}
using that ${\rm Tr} (A)=0$.
Hence, setting 
\[
\begin{split}
\|\Delta^{\alpha_0,0,A}v\|_{L^p((0,T)\times \R^N)}^p\, &:= \, \| {\rm Tr} \big(B_Ie^{ t A^*} D^2 v (t, e^{tA} \cdot ) e^{tA}B_I^{ 1/2}\big)  \|_{L^p((0,T)\times \R^N)}^p\\
&= \,  \| {\rm Tr} \big(B_Ie^{ t A^*} D^2 v e^{tA}B_I^{ 1/2}\big)  \|_{L^p((0,T)\times \R^N)}^p,
\end{split}\] 
we get from Definition \eqref{SEMI_NORM_SOB} that the Estimates \eqref{SOB_EST_GLOB} rewrite in term of $v$ as:
\begin{equation} 
\label{EST_SOB_V}
\Vert v \Vert_{\dot{W}^{2,p,A}_\alpha}^p\, := \, \sum_{i=1}^{n} \|\Delta^{\alpha_i,i,A}v\|_{L^p((0,T)\times \R^N)}^p\le \tilde C_p \|f\|_{L^p((0,T)\times \R^N)}^p,
\end{equation}
with $\tilde C_p = C_p^p.$ We now want to prove that for $w$ solving \eqref{d2}, namely
\begin {equation*}
\begin{cases}
 \partial_t w(t, z)    + {\text Tr} \big( e^{tA} B
 e^{tA^*} D^2 w(t, z) \big) + {\text Tr} \big( e^{tA} S(t)
 e^{tA^*} D^2 w(t, z) \big)\,  = \, \tilde f(t,z),\\
 w(0,z)\, = \, 0,
 \end{cases}
\end{equation*}
it also holds that
\begin{equation}
\label{EST_SOB_W}
\Vert w\Vert_{\dot{W}^{2,p,A}_\alpha}^p\, :=\, \sum_{i=1}^n \|\Delta^{\alpha_i,i,A}w\|_{L^p((0,T)\times \R^N)}^p \, \le \,  \tilde C_p \|f\|_{L^p((0,T)\times \R^N)}^p,
\end{equation}
with the same constants $\tilde C_p$ as  in \eqref{EST_SOB_V}. This can be done through the previous perturbative approach of Section \ref{SEC_PERTURBATIVO} employed to prove Theorem \ref{uno}, which actually gives the expected control for the  the second order derivatives contribution of the semi-norm $\Vert \cdot\Vert_{\dot{W}^{2,p,A}_\alpha} $.\newline
For the other contributions and with the notations of Section \ref{SEC_PERTURBATIVO}, with $Q'(s)=e^{sA}S(s)e^{sA^*}$ and  with $m $ which is the Lebesgue measure on $[0,T] \times \R^N$ (indeed in the present  case $g(t)= {\rm det}(e^{-At}) =1$, for all $t$), we would get that
\[\begin{split}
\Vert \bar v_\epsilon\Vert_{\dot{W}^{2,p,A}_\alpha}^p\, &= \, \sum_{i=1}^n\|   \Delta^{\alpha_i,i,A}\bar v_{\epsilon}   \,  \|_{L^p((0,T)\times \R^N)}^p\\
&= \,  \sum_{i=1}^n \int_0^T \int_{\R^N} |\Delta^{\alpha_i,i,A} \bar v_{\epsilon} (t,z) |^p \, dz dt
\\ 
&= \, \sum_{i=1}^n \int_0^T \int_{\R^N} \left| \mathbb{E} [\Delta^{\alpha_i,i,A}  v_{\epsilon} (t,z + \epsilon X_t) ] \,\right|^p \, dz dt
\\
&\le 
\, \sum_{i=1}^n   \int_0^T \int_{\R^N} \mathbb{E} \left[ | \Delta^{\alpha_i,i,A}  v_{\epsilon} (t,z + \epsilon X_t)  \,|^p\right] \, dz dt
\\
&= \,  \sum_{i=1}^n\mathbb{E} \left[\int_0^T \int_{\R^N} \,    | \Delta^{\alpha_i,i,A}  v_{\epsilon}  (t,z + \epsilon X_t)  \,|^p\, dz dt\right]
\\
&= \,  \sum_{i=1}^n\textcolor{black}{\mathbb E}\left[\int_0^T \int_{\R^N}  |  \Delta^{\alpha_i,i,A}  v_{\epsilon} (t,\bar z )  \,|^p \, d\bar z dt\right]
\\
&\le \,  \tilde C_p  \|   f  \|_{L^p((0,T)\times \R^N)}^p,
\end{split}
\]
{using for the last inequality that $v_\epsilon$ also satisfies \eqref{EST_SOB_V} (similarly to what had been established in \eqref{w11})}. The same previous compactness argument then yields \eqref{EST_SOB_W}.\newline
Setting eventually $\tilde u(t,z) := w(t, e^{tA} z)$, which is the unique integral solution (smooth in space) of
 \begin{equation*}
\begin{cases}
 \partial_t \tilde u(t,z) +  L_t^{\text{ou}, S} \tilde u(t,z) = f  (t,z), &\mbox{ on } (0,T)\times \R^N,\\
 \tilde u(0,z)=0, &\mbox{ on }  \R^N,
\end{cases}
\end{equation*}
where  $ L_t^{\text{ou}, S}$ introduced in \eqref{ciao}  is the Ornstein-Uhlenbeck operator  perturbed at second order, we derive that
\begin{equation}\label{SOB_EST_GLOB_PERT}
\Vert \tilde u \Vert_{\dot{W}^{2,p}_\alpha} \, \le \, C_p\Vert f \Vert_{L^p},
\end{equation}
with $C_p$ as in \eqref{SOB_EST_GLOB}. We have thus extended the results of Theorem \ref{d44} for the anisotropic Sobolev semi-norm in \eqref{SEMI_NORM_SOB}. The estimate \eqref{SOB_EST_GLOB} is stable for a continuous, non-negative definite, second order perturbation of the underlying degenerate Ornstein-Uhlenbeck operator.

\subsection{Anisotropic Schauder estimates}
Following Krylov \cite{book:Krylov96}, for  some fixed $\ell$ in $\N_0:=\N\cup\{0\}$ and $\beta$ in $(0,1]$, we introduce for a function $\phi\colon \R^N\to \R$ the Zygmund-H\"older semi-norm  as
\[[\phi]_{C^{\ell+\beta}} \, := \,
\begin{cases}
    \sup_{\vert \vartheta \vert= \ell}\sup_{x\neq y}\frac{\vert D^\vartheta\phi(x)-D^\vartheta\phi(y)\vert}{\vert x-y\vert^\beta} , & \mbox{if
    }\beta \neq 1; \\
    \sup_{\vert \vartheta \vert= \ell}\sup_{x\neq y}\frac{\bigl{\vert}D^\vartheta\phi(x)+D^\vartheta\phi(y)-2D^\vartheta\phi(\frac{x+y}{2})
    \bigr{\vert}}{\vert x-y \vert}, & \mbox{if } \beta =1
                             \end{cases}\]
(we are using usual multi-indices $\vartheta$ for the partial derivatives). Consequently, the Zygmund-H\"older space $C^{\ell+\beta}_b(\R^N)$ is the family of bounded functions $\phi\colon \R^N
\to\R$ such that $\phi$ and its derivatives up to order $\ell$ are continuous and the norm
\[\Vert \phi \Vert_{C^{\ell+\beta}_b} \,:=\,  \sum_{i=0}^{\ell}\sup_{\vert\vartheta\vert = i}\Vert D^\vartheta\phi
\Vert_{\infty}+[\phi]_{C^{\ell+\beta}} \,\text{ is finite.}\]

We can now define the anisotropic Zygmund-H\"older spaces associated with the current setting and which again reflect the various scales already introduced in \eqref{INDEXES}. Fixed $\gamma\in (0,3)$, the space $C^{\gamma}_{b,d}(\R^N)$ is
the family of functions $\phi\colon \R^N\to \R$ such that for any $i$ in $\llbracket 1,n\rrbracket$ and any $z_0$ in $\R^N$, the function
\[w\in  \R^{ d_i}\,  \to \, \phi(z_0+E_i(w)) \in \R \,\text{ belongs to }C^{\gamma/(2i-1)}_b\left(\R^{ d_i}\right),\]
with a norm bounded by a constant independent from $z_0$. In the above expression, we recall that
 the $\{E_i\colon i\in \{2,\cdots,n\}\}$ have been defined in the previous paragraph, $d_1=d$ and $E_1$ is the embedding matrix from $\R^{d} $ into $\R^N $.
 It is endowed with the norm
\begin{equation}\label{eq:def_anistotropic_norm}
\Vert\phi\Vert_{C^{\gamma}_{b,d}} \,:=\,\sup_{z_0\in \R^N}\Vert\phi\big(z_0+E_0(\cdot)\big)\Vert_{C^{\gamma}_b (\R^{d})  }+\sum_{i=1}^{k}\sup_{z_0\in \R^N}[\phi\big(z_0+E_i(\cdot)\big)]_{C^{\gamma/(1+2i)} (\R^{d_i} )   }.
\end{equation}
We denote by $C^{\gamma}_{b,d}$ this function space because the regularity exponents reflect again the multi-scale features of the system. The norm could equivalently be defined through the corresponding spatial parabolic distance $d$ defined as follows. For all $z=(x,y)$, $z'=(x',y')$ in $\R^N=\R^d\times\R^{d'}$:
$$d(z,z'):=|x-x'|+\sum_{i=2}^{n}\vert y_i-y_i'\vert^{\frac{1}{2i-1}},$$
where the exponents are again those who appeared in \eqref{INDEXES}.
\newline
Let as before  $f$ be in $B_b\left(0,T;C^\infty_0(\R^N)\right)$.
Under \textbf{[K]}, by the results of Lunardi \cite{Lunardi97}, it follows that the unique bounded solution of the Cauchy Problem \eqref{eq:OU_initial:intro} (written in integral form) verifies the following anisotropic Schauder estimates:
\begin{equation}\label{SCHAU_OU}
\Vert u \Vert_{L^\infty((0,T),C^{2+\beta}_{b,d})} \, \le \, \textcolor{black}{C_\beta} \Vert f \Vert_{L^\infty((0,T),C^{\beta}_{b,d})},
\end{equation}
for some constant ${C_\beta}$ independent from $f$, i.e.,
\[
\sup_{0 \le t \le T}\Vert u (t, \cdot )\Vert_{C^{2+\beta}_{b,d}} \, \le \, \textcolor{black}{C_\beta}  \sup_{0 \le t \le T }\Vert f (t, \cdot )\Vert_{C^{\beta}_{b,d}},
\]
We again set as in the previous paragraph $u(t,z) = v(t, e^{tA} z)$
and since $ u$ solves \eqref{eq:OU_initial:intro} we have that $v$ in turn solves \eqref{ma}. Write:
\begin{equation}
    \begin{split}
    \Vert u \Vert_{L^\infty((0,T),C^{2+\beta}_{b,d})}\, &= \, \Vert v(t,e^{tA}\cdot) \Vert_{L^\infty((0,T),C^{2+\beta}_{b,d})}\, = : \, \Vert v \Vert_{L^\infty((0,T),C^{2+\beta}_{b,d,A})}\\
&\le \,\textcolor{black}{C_\beta} \Vert f \Vert_{L^\infty((0,T),C^{\beta}_{b,d})} 
 \,\le\,  {C_\beta}\Vert \tilde f(t,e^{tA}\cdot) \Vert_{L^\infty((0,T),C^{\beta}_{b,d})}\\
 &= \,  {C_\beta}\Vert \tilde f\Vert_{L^\infty((0,T),C^{\beta}_{b,d,A})},
    \end{split}
    \label{SC_EST_V}    
\end{equation}
denoting $\tilde{f}(t,z):= f(t,e^{-tA}z)$.
We again want to prove, as in Section $1.1$ ,
that for $w$ solving \eqref{d2},
\begin{align}\label{STIMA_W_PER_SC}
\Vert w \Vert_{L^\infty((0,T),C^{2+\beta}_{b,d,A})}\le \textcolor{black}{C_\beta}\Vert \tilde f\Vert_{L^\infty((0,T),C^{\beta}_{b,d,A})}
\end{align}
with the same  constant $\textcolor{black}{C_\beta}$ as in \eqref{SC_EST_V}. We proceed \textcolor{black}{one more time} through the previous perturbative approach of Section \ref{SEC_PERTURBATIVO}.
With the notations employed therein, we deduce that there exists a unique  solution $v_\epsilon=$  PDE$(Q,\tilde f(t,z-\epsilon X_t))$, depending also on $\epsilon$ and  $\omega$ as parameters  such that
\begin{align}
\label{proof:estimate1_SC}
\sup_{(t,z)\in [0,T]\times \R^N}| v_\epsilon(t,z)| \, &\le \, T \sup_{(t,z)\in [0,T]\times \R^N}|\tilde  f(t,z)|.
\end{align}
By the translation invariance of the H\"older-norms, using also that $X_t=e^{tA}e^{-tA}X_t$, it is not difficult to prove that, for any $\omega$, $\mathbb{P}$-a.s., 
\begin{equation}\label{wqs}
\|  \tilde f  \|_{L^\infty((0,T),C^{\beta}_{b,d,A})} = \|  \tilde f(\cdot , \cdot \, -\epsilon X_{\cdot}) \|_{L^\infty((0,T),C^{\beta}_{b,d,A})} .
\end{equation}
Thus it also holds from \eqref{SC_EST_V}
\begin{equation}\label{w11_BIS}
\|  v_{\epsilon}    \|_{L^\infty((0,T),C^{2+\beta}_{b,d,A})} \le \textcolor{black}{C_\beta}  \|  \tilde f  \|_{L^\infty((0,T),C^{\beta}_{b,d,A})}.
\end{equation}
Recalling now that $\bar v_\epsilon(s,z) = \mathbb{E}[v_\epsilon(s,z+\epsilon X_s)]$ 
is   an integral solution of 
\begin{equation*}
    \partial_t \bar v_\epsilon(t,z) \, = \, \text{tr}(Q(t)D^2_z \bar v_\epsilon(t,z))+\lambda\left(\bar v_\epsilon(t,z+\epsilon l(t))-\bar v_\epsilon(t,z)\right)+ 
    \tilde f(t,z),
\end{equation*}
with zero initial condition, we write for any $i$ in $\{2,\cdots, n\}$, $w,w'$ in $\R^{d_i}$ e $ (t,z_0)$ in $[0,T]\times \R^N$, that
\begin{align*}
|\bar{v}_\epsilon(t,&e^{At}(z_0+E_i(w))- \bar{v}_\epsilon(t,e^{At}(z_0+E_i(w'))|\\
&\le \, \mathbb{E}\left[ |v_\epsilon(t,e^{At}(z_0+E_i(w))+\epsilon  e^{At}e^{-At} X_t)- v_\epsilon(t,e^{At}(z_0+E_i(w'))+\epsilon e^{At}e^{-At}  X_t)|\right]  
\\
&\le \, \mathbb{E}\left[ [v_\epsilon(t,e^{At}(z_0+E_i(\cdot)))]_{C^{\frac{2+\beta}{2i-1}}} \right]|w-w'|^{\frac{2+\beta}{2i-1}}.
\end{align*}
 Hence,
\[
\left[\bar v_\epsilon(t,e^{-At}(z_0+E_i(\cdot)))\right]_{C^{\frac{2+\beta}{2i-1}}}\, \le \, \mathbb{E}\left[ [v_\epsilon(t,e^{-At}(z_0+E_i(\cdot)))]_{C^{\frac{2+\beta}{2i-1}}} \right].
\]
We would get similarly, that
\[
\left[D_x^2\bar v_\epsilon(t,e^{At}(z_0+E_1(\cdot)))\right]_{C^{\beta}}\, \le \,\mathbb{E}\left[ [D_x^2v_\epsilon(t,e^{At}(z_0+E_1(\cdot)))]_{C^{\beta}} \right],
\]
and for all $k\in \{0,1,2\} $,
\[
\left|D_x^k\bar v_\epsilon(t,e^{At}(z_0+E_1(\cdot)))\right|_\infty \, \le \,  \mathbb{E}\left[|D_x^kv_\epsilon(t,e^{At}(z_0+E_1(\cdot)))|_\infty\right].
\]
Summing all those contributions, we thus derive from \eqref{eq:def_anistotropic_norm}, \eqref{SC_EST_V} that: 
\begin{equation}
\|\bar v_\epsilon\|_{L^\infty((0,T),C^{2+\beta}_{b,d,A})}\, \le \,  \sup_{0 \le t \le T}\mathbb{E}\left[\|v_\epsilon (t, \cdot)\|_{C^{2+\beta}_{b,d,A})}\right] \,  \le\, \textcolor{black}{C_\beta}  \|  \tilde f  \|_{L^\infty((0,T),C^{\beta}_{b,d,A})},
\end{equation}
using \eqref{w11_BIS} for the last inequality.  Now, continuing   as in Section 3.2, using also  a compactness argument, one would derive
 that \eqref{STIMA_W_PER_SC} indeed holds.

Going  backwards, setting  $\tilde u(t,z) := w(t, e^{tA} z)$,   we find that  $\tilde u  $  is the unique (integral) solution $u_S$ to \eqref{eq:OU_PERT0}:
 we finally derive that
\begin{equation}\label{SCHAU_OU_PERT}
\Vert  u_S \Vert_{L^\infty((0,T),C^{2+\beta}_{b,d})} \, \le \, \textcolor{black}{C_\beta} \Vert f \Vert_{L^\infty((0,T),C^{\beta}_{b,d})},
\end{equation}
where $C_{\beta}$ is the same constant as in \eqref{SCHAU_OU}. Estimate  \eqref{SCHAU_OU_PERT} provides the extension of Theorem \ref{d44} for the anisotropic Schauder estimates.

\begin{remark}
Let us mention that for the perturbative method to work, roughly speaking, few properties were actually needed on the underlying norm. Namely, we used the translation invariance and some kind of commutation between the norm (or a function of the norm in the $L^p$-case) and  expectation. Hence, this approach could  possibly be applied to a much wider class of estimates in other function spaces (like e.g. Besov spaces). This will concern further research.
 \end{remark}

\pagebreak

\fancyhead[LE]{Bibliography}
\fancyhead[RO]{}
\bibliography{bibli}

\bibliographystyle{alpha}
\addcontentsline{toc}{chapter}{Bibliography}


\newgeometry{top=1.5cm, bottom=1.25cm, left=2cm, right=2cm}
\fontfamily{rm}\selectfont

\lhead{}
\rhead{}
\rfoot{}
\cfoot{}
\lfoot{}
\thispagestyle{empty}
\noindent
\thispagestyle{empty}
\includegraphics[height=2.45cm]{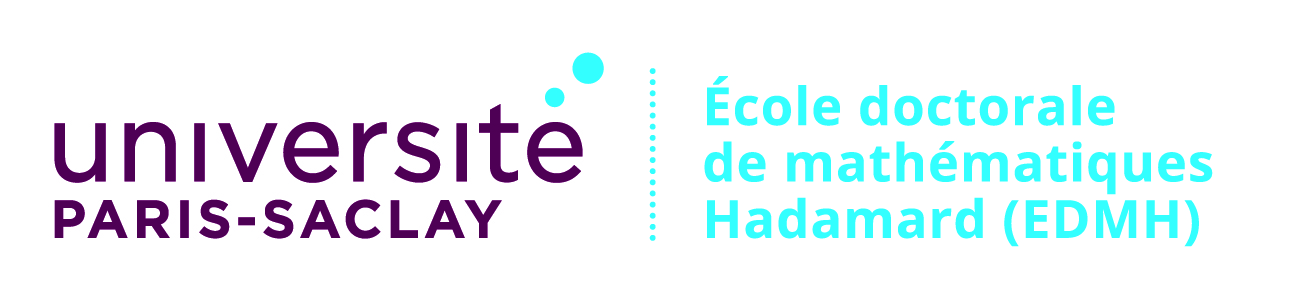}
\vspace{1cm}
\begin{mdframed}[linecolor=Prune,linewidth=1]
\vspace{-.25cm}
\paragraph*{Titre:} 
Régularisation faible par un bruit de Lévy dégénéré et applications
\begin{small}
\vspace{-.25cm}
\paragraph*{Mots clés:} 
EDOs mal posées, Solutions faibles via régularisations, \'Equations dégénérées de Kolmogorov, Estimées de Schauder pour opérateurs intégro-différentiels, Méthodes de la parametrix, Processus de Lévy
\vspace{-.5cm}
\paragraph*{Résumé:}
Après une introduction générale sur le phénomène de régularisation par le bruit dans le cadre dégénéré, la première partie de cette thèse est consacrée à l'obtention d'estimées de Schauder, un outil analytique utile pour établir le caractère bien posé des équations différentielles stochastiques (EDS), pour deux différents classes d'équations de Kolmogorov sous une condition de type Hörmander faible, dont les coefficients appartiennent à des espaces de Holder anisotropes appropriés avec des multi-indices de régularité.  
La première classe considère un système non linéaire dirigé par un opérateur $\alpha$-stable symétrique agissant uniquement sur certaines composantes. Notre méthode de preuve repose sur une approche perturbative basée sur des développements parametrix progressifs par des formules de type Duhamel. En raison des faibles propriétés de régularisation données par le cadre dégénéré, nous exploitons également certains contrôles sur les normes de Besov, afin de traiter la perturbation non-linéaire.
Dans le prolongement de la première, nous présentons également les estimées de Schauder pour un opérateur dégénéré d'Ornstein-Uhlenbeck associé à une classe plus large d'opérateurs de type $\alpha$-stable, comme l'opérateur stable relativiste ou de Lamperti. La preuve de ce résultat repose plutôt sur une analyse précise du comportement du semi-groupe de Markov correspondante entre les espaces de Hölder anisotropes et quelques techniques d'interpolation.
En exploitant une approche parametrix rétrograde, la deuxième partie de cette thèse cherche à établir le caractère bien-posé au sens faible d'une chaîne dégénérée de EDS dirigées par la même classe de processus de type $\alpha$-stable, sous des hypothèses de régularité de Hölder minimale sur les coefficients. Comme corollaire de notre méthode, nous présentons également des estimations de type Krylov d'intérêt indépendant pour le processus canonique sous-jacent. Enfin, nous soulignons à travers des contre-exemples appropriés qu'il existe bien un seuil (presque) optimal sur les exposants de régularité assurant le caractère faiblement bien posé pour l'EDS.
En lien avec quelques applications mécaniques pour des dynamiques cinétique avec frottement, nous concluons en étudiant la stabilité des perturbations du second ordre pour  des opérateurs de Kolmogorov dégénérés en normes Lp et Hölder.
\end{small}
\end{mdframed}
\newpage

\thispagestyle{empty}

\begin{mdframed}[linecolor=Prune,linewidth=1]
\vspace{-.25cm}
\paragraph*{Title:} 
Weak regularization by degenerate Lévy noise and its applications

\begin{small}
\vspace{-.25cm}
\paragraph*{Keywords:} 

Ill posed ODEs, Weak solutions through noise regularization, Kolmogorov degenerate equations, Schauder estimates for integro-differential operators, Parametrix Methods, Lévy processes

\vspace{-.5cm}
\paragraph*{Abstract:} 

After a general introduction about the regularization by noise phenomenon in the degenerate setting, the first part of this thesis focuses at establishing the Schauder estimates, a useful analytical tool to prove also the well-posedness of stochastic differential equations (SDEs), for two different classes of Kolmogorov equations under a weak Hörmander-like condition, whose coefficients lie in suitable anisotropic Hölder spaces with multi-indices of regularity.
The first class considers a nonlinear system controlled by a symmetric $\alpha$-stable operator acting only on some components. Our method of proof relies on a perturbative approach based on forward parametrix expansions through Duhamel-type formulas. Due to the low regularizing properties given by the degenerate setting, we also exploit some controls on Besov norms, in order to deal with the non-linear perturbation.
As an extension of the first one, we also present Schauder estimates associated with a degenerate Ornstein-Uhlenbeck operator driven by a larger class of $\alpha$-stable-like operators, like the relativistic or the Lamperti stable one. The proof of this result relies instead on a precise analysis of the behaviour of the associated Markov semigroup between anisotropic Hölder spaces and some interpolation techniques.
Exploiting a backward parametrix approach, the second part of this thesis aims at establishing the well-posedness in a weak sense of a degenerate chain of SDEs driven by the same class of $\alpha$-stable-like processes, under the assumptions of the minimal Hölder regularity on the coefficients.  As a by-product of our method, we also present
Krylov-type estimates of independent interest for the associated canonical process. 
Finally, we emphasize through suitable counter-examples that there exists indeed an (almost) sharp threshold on the regularity exponents ensuring the weak well-posedness for the SDE.
In connection with some mechanical applications for kinetic dynamics with friction, we conclude by investigating the stability of second-order perturbations for degenerate Kolmogorov operators in Lp and Hölder norms.
\end{small}
\end{mdframed}

\vspace{3cm} 

\fontfamily{fvs}\fontseries{m}\selectfont
\begin{tabular}{p{14cm}r}
\multirow{3}{16cm}[+0mm]{{\color{Prune} Université Paris-Saclay\\
Espace Technologique / Immeuble Discovery\\
Route de l’Orme aux Merisiers RD 128 / 91190 Saint-Aubin, France}} & \multirow{3}{2.19cm}[+9mm]{
}\\
\end{tabular}

\end{document}